\documentclass[a4paper,12pt,oneside,reqno]{amsart}

\usepackage{amssymb,amsfonts,amsmath,amsthm,amscd, bbm}
\usepackage{graphicx, color}
\usepackage{bm}
\usepackage{comment}
\usepackage{cases}
\usepackage[pdf]{pstricks}
\usepackage{algorithm, algpseudocode}

\numberwithin{equation}{section}  

\newtheorem{teor}{Theorem}%[section]

\usepackage[utf8]{inputenc}
\usepackage{placeins}
\usepackage{float}
\usepackage{subcaption}
\newcommand{\EE}{{\sf{E}}}

\newcommand{\beq}{\begin{equation}}
\newcommand{\eeq}{\end{equation}}
\newcommand{\bea}{\begin{aligned}}
\newcommand{\eea}{\end{aligned}}

\newcommand{\PP}{\mathbb P}
\newcommand{\R}{\mathbb R}
\newcommand{\N}{\mathbb N}
\newcommand{\E}{\mathbb E}
\newcommand{\tr}{{\rm tr} \,}

\newcommand{\bfm}{{\boldsymbol m}}

\newcommand{\be}{\beta}
\newcommand{\vare}{\varepsilon}

\begin{document}

\title{TAP equations are repulsive.}
%\author[S. Gufler]{Stephan Gufler}
%\author[J.L. Igelbrink]{Jan Lukas Igelbrink}
%\author[N. Kistler]{Nicola Kistler}
%\address[S. Gufler, N. Kistler]{ J.W. Goethe-Universit\"at Frankfurt, Germany.}
%\address[J.L. Igelbrink]{Institut f\"ur Mathematik, Johannes Gutenberg-Universit\"at Mainz, Germany.}
%\email[S. Gufler]{gufler@math.uni-frankfurt.de}
%\email[J.L. Igelbrink]{jigelbri@uni-mainz.de}
%\email[N. Kistler]{kistler@math.uni-frankfurt.de}
 \author[S. Gufler]{Stephan Gufler}
 \address{Stephan Gufler \\ J.W. Goethe-Universit\"at Frankfurt, Germany.}
 \email{gufler@math.uni-frankfurt.de}

 \author[Jan Lukas Igelbrink]{Jan Lukas Igelbrink}
 \address{Jan Lukas Igelbrink  \\ Institut f\"ur Mathematik, 
 Johannes Gutenberg-Universit\"at Mainz, Germany.}
 \email{jigelbri@uni-mainz.de}

 \author[N. Kistler]{Nicola Kistler}
 \address{Nicola Kistler, J.W. Goethe-Universit\"at Frankfurt, Germany.}
 \email{kistler@math.uni-frankfurt.de}

\thanks{This work has been supported by a DFG research grant, contract number 2337/1-1. We acknowledge the \emph{Allianz f\"ur Hochleistungsrechnen Rheinland-Pfalz} for granting us access to the High Performance Computing \textsc{Elwetritsch}, on which our numerical simulations have been performed. It is furthermore a pleasure to thank Yan V. Fyodorov, Giorgio Parisi, Timm Plefka, Federico Ricci-Tersenghi and Marius A. Schmidt for enlightening conversations.}

%\subjclass[2000]{60J80, 60G70, 82B44} \keywords{Mean Field Spin Glasses, Large Deviations, Gibbs-Boltzmann and Parisi Variational Principles} 
\date{\today}
\begin{abstract}
We show that for low enough temperatures, but still above the AT line, the Jacobian of the TAP equations for the SK model has a macroscopic fraction of eigenvalues outside the unit interval. This provides a simple explanation for the numerical instability of the fixed points, which thus occurs already in high temperature. The insight leads to some algorithmic considerations on the low temperature regime, also briefly discussed.
\end{abstract}
\vspace*{-1.5cm}\maketitle
\vspace*{-0.5cm}\section{Introduction} 
Let $N \in \mathbb N, \beta>0, h\in\mathbb R$, and consider independent standard Gaussians $\boldsymbol g = (g_{ij})_{1\le i<j \le N}$ issued on some probability space $(\Omega, \mathcal F, \PP)$. The TAP equations \cite[after Thouless, Anderson and Palmer]{tap} for the spin magnetizations $\boldsymbol m = (m_i)_{i=1}^N \in [-1,1]^N$ in the SK-model \cite[after Sherrington and Kirkpatrick]{sk} at inverse temperature $\be$ and external field $h$ read
\begin{equation}\label{e:TAP}
	m_i=\tanh\left(h + \frac{\beta}{\sqrt N} \sum_{j:j\neq i} g_{ij} m_j -\beta^2(1-q)m_i \right),\quad i=1,\ldots, N.
\end{equation}
The $\be^2$-term is the (limiting) Onsager correction: it involves the {\it high temperature} order parameter
 $q$ which is the (for $\beta\le 1$ or $h\ne 0$ unique, see \cite[Proposition 1.3.8]{Talagrand}) solution of the fixed point equation
\begin{equation}\label{e:q}
	q=\mathbb \EE \tanh^2\left(h+\beta\sqrt{q}Z\right),
\end{equation}
$Z$ being a standard Gaussian, and $\EE$ its expectation. 

The concept of {\it high temperature} is related to the AT-line \cite[after de Almeida and Thouless]{at}, i.e. the $(\be, h)$-region satisfying 
\begin{equation}\label{e:AT}
	\beta^2\EE\frac{1}{\cosh^4(h+\beta\sqrt{q}Z)} = 1.
\end{equation}
For $(\beta, h)$ where the l.h.s. is strictly less than unity, the system is allegedly in the replica symmetric phase (high temperature), see e.g. Adhikari {\it et.al.} \cite{acsy} for a recent thorough discussion of this issue, and Chen \cite {chen} for evidence supporting the conjecture.

Due to the fixed point nature of the TAP equations \eqref{e:q}, one is perhaps tempted to solve them numerically via classical Banach iterations, i.e. to approximate solutions via   
\beq\label{e:1step}
m_i^{(k+1)}=\tanh\left(h + \frac{\beta}{\sqrt N} \sum_{j:j\neq i} g_{ij} m_j^{(k)} -\beta^2(1-q)m_i^{(k)} \right),\quad i=1,\ldots, N,
\eeq
As it turns out, these plain iterations are astonishingly {\it in}-efficient in finding stable fixed points\footnote{this observation has been made, at times anecdotally, ever since: it is already present in M\'ezard, Parisi and Virasoro \cite[Section II.4]{mpv}.} insofar the mean squared error between two iterates, to wit:
\beq \label{mse}
\text{MSE}\left( \boldsymbol m^{(k+1)}, \boldsymbol m^{(k)}\right) \equiv \frac{1}{N} \sum_{i=1}^N \left( m_i^{(k+1)} - m_i^{(k)} \right)^2, 
\eeq
often remains large, even for very large $k's$. Even more surprising, this issue is not restricted to the low temperature phase, cfr. Figure \ref{figure1} below. 
\begin{figure}[h]  
\includegraphics[width=8cm, height = 5cm, trim={0.45cm 0.4cm 0 0},clip]{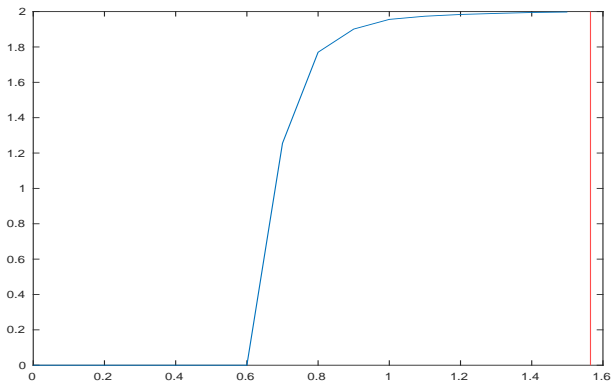}
\caption{MSE between (the last) two iterates as a function of $\beta$, to fixed magnetic field $h=0.5$, 
a single realization of the disorder, and system size $N=25$. We have run Banach algorithm $k= 1000$ times, randomly initialized (uniformly chosen $\boldsymbol m^{(0)}$). Visibly, iterations stabilize for small $\beta$, but  diverge beyond the threshold $\approx 0.6$, way below the AT-line (red).}
\label{figure1}
\end{figure}
Bolthausen \cite{b1} bypasses this problem by means of a modification of the Banach algorithm (recalled below) which converges up to the AT-line. Notwithstanding, the origin of the phenomenon captured by Figure \ref{figure1} has not yet been, to our knowledge, identified. It is the purpose of this work to fill this gap. Precisely, we show in Theorem \ref{t:rep} below that classical Banach iterates become unstable for the simplest reason: for large enough $\beta$, but still below the AT-line, TAP equations become {\it repulsive}. This should be contrasted with the classical counterpart of the SK-model: it is well known (and a simple fact) that {\it relevant} solutions of the fixed point equations for the Curie-Weiss model are attractive at {\it any} temperature, with the {\it irrelevant} solution even becoming repulsive in low temperature.
\section{Main result}
\subsection{Iterative procedure for the magnetizations}
Bolthausen~\cite{b1} constructs magnetizations $\bfm^{(k)} =\left(m^{(k)}_i\right)_{i\leq N}$ for given disorder $(g_{ij})_{1\le i<j\le N}$ through an iterative procedure in $k=1, 2, \dots$. These magnetizations approximate fixed points of the TAP equations in the limit $N\to\infty$ followed by $k\to\infty$. The algorithm uses the initial values $\bfm^{(0)}=\boldsymbol 0$, $\bfm^{(1)}=\sqrt{q}\boldsymbol{1}$. The iteration step reads
\begin{equation}\label{eq:twostep}
	m^{(k+1)}_i=\tanh\left( h+\frac{\beta}{\sqrt N} \sum_{j:\: j\ne i} g_{ij} m^{(k)}_j - \beta^2 (1-q) m^{(k-1)}_i \right),\end{equation}
for $i\leq N$, and $k \in \N$. Remark in particular that contrary to the classical Banach algorithm, the above scheme invokes a time delay in the Onsager correction: we will thus refer to \eqref{eq:twostep} as {\it Two Steps Banach algorithm}, \textsc{2SteB} for short. By \cite[Proposition 2.5]{b1}, the quantity $q_N:=N^{-1}\sum_{i=1}^{N}(m^{(k)}_i)^2$ converges to $q$ in probability and in expectation, as $N\to\infty$ followed by $k\to\infty$. Moreover, by \cite[Theorem 2.1]{b1},
\begin{equation}\label{e:convit}
	\lim_{k,k'\to\infty}\limsup_{N\to\infty}\E\left[ \frac{1}{N} \sum_{i=1}^N \left(m^{(k)}_i - m^{(k')}_i \right)^2 \right]=0,
\end{equation}
provided $(\beta,h)$ is below the AT-line, i.e.\ if the left-hand side of~\eqref{e:AT} is less than unity.
\subsection{Spectrum of the Jacobian}
Denote by $F_i(\boldsymbol m^{(k)})$ the right-hand side of~\eqref{e:1step}, and consider the Jacobi matrix $J^{(k)}$ with entries
\begin{equation}
	J^{(k)}_{ij} := \frac{\partial F_i}{\partial m_j} (\bfm^{(k)}) = \begin{cases}
		\frac{\beta}{\sqrt{N}} g_{ij} \left(1-F_i\left(\bfm^{(k)}\right)^2\right): i\ne j\\
		-\beta^2(1-q) \left(1-F_i\left(\bfm^{(k)}\right)^2\right): i=j,
	\end{cases}
\end{equation}
omitting the obvious $N$-dependence to lighten notations.

The Jacobi matrix $J^{(k)}$ is not symmetric, but we claim that it is nonetheless diagonalizable, and that all eigenvalues are real. To see this we write the Jacobian as a product
\begin{equation} \label{repr}
	J^{(k)} = \hat J {\rm diag}\left( 1-F_\cdot\left(\bfm^{(k)}\right)^2 \right),
\end{equation}
where
\begin{equation} \label{e:jacobi-product}
	\hat J_{ij}  := \begin{cases}
		\frac{\beta}{\sqrt{N}} g_{ij} & i < j, \\
		\frac{\beta}{\sqrt{N}} g_{ji} & j < i, \\
		-\beta^2(1-q) & i= j.
	\end{cases}
\end{equation}
The first matrix on the r.h.s. of \eqref{repr} is symmetric, whereas the second is negative definite: it thus readily follows from \cite[Theorem 2]{DH65} that $J^{(k)}$ itself is diagonalizable, and that all its eigenvalues are real, settling our claim. 

Denoting the ordered sequence of the real eigenvalues of a matrix $M$ by $\lambda_1(M)\ge \ldots\ge \lambda_N(M)$, we then consider the empirical spectral measure of the Jacobian
\begin{equation}
	\mu^{(k)} :=\frac{1}{N}\sum_{i=1}^N\delta_{\lambda_i(J^{(k)})}.
\end{equation}
Our main result\footnote{See also \cite{cfm} for a treatment similar in spirit to our considerations, albeit with radically different tools, for $\mathbb Z_2$-synchronization.} states that in a region below the AT-line, $\mu^{(k)}$ has mass outside the unit interval. 
\begin{teor}\label{t:rep}
	 For all $\beta\in(\sqrt 2-1,1)$ and $\epsilon >0$, there exists $h_0>0$ such that the following holds true: for all $h\in[0,h_0]$, there exists $k_0\in\N$, and for all $k\ge k_0$, there exists $N_0\in\N$, such that $\PP\left( \mu^{(k)}(-\infty,-1)>0 \right) \geq 1- \epsilon$, for all $N\ge N_0$. 
\end{teor}
TAP equations thus become repulsive due a macroscopic fraction of low-lying ($< -1 $) eigenvalues. Before giving a proof of this statement, we remark that the measures $\mu^{(k)}$ converge to a limiting measure $\mu$ which can be stated as a free multiplicative convolution
	\begin{equation}
		\label{e:freeconvJ}
		\mu=\lambda_{\beta,q}\boxtimes\nu,
	\end{equation}
	where $\lambda_{\beta, q}$ is the law of $\beta X - \beta^2(1-q)$ for $X$ distributed according to the 	standard semicircular law with density $x \ni [-2,2] \mapsto (2\pi)^{-1}\sqrt{(4-x^2}$, and $\nu$ is the law of $1- \tanh^2(h + \beta \sqrt{q} Z)$ for $Z$ standard Gaussian. To sketch a proof of this claim, we use the decomposition~\eqref{repr} and make the following observations: {\it i)} The empirical spectral distribution of the first factor $\hat J$ weakly-converges almost surely to the scaled/shifted semicircular law $\lambda_{\beta, q}$; {\it ii)} The empirical spectral distribution of the second factor $1- F_\cdot(\bfm^{(k)})^2$ can similarly be shown to converge to $\nu$. These observations, together with the (asymptotic) independence of the two factors which follows from \cite{b2}, and finally \cite[Theorem 5.4.2]{AGZ} then yield the representation \eqref{e:freeconvJ}. 

The free convolution can also be evaluated more explicitly using Voiculescu's S-transform via inversion of moment generating functions, see e.g. \cite[Chapter~5.3]{AGZ}. We believe this approach allows to remove the \emph{small-$h$ condition} in Theorem \ref{t:rep}. The ensuing analysis is however both long and (tediously) technical. As the outcome arguably adds little to the main observation of this work, we refrain from pursuing this route here. 
%\end{rem}
\begin{proof}[Proof of Theorem~\ref{t:rep}]
	In a first step we consider a simplified version of the Jacobian, namely the matrix $\hat J$ from~\eqref{e:jacobi-product} in place of $J^{(k)}$. We write $\hat J=W-\beta^2(1-q)I$, where $W$ is a Wigner matrix. As a consequence of Wigner's theorem, see e.g.~\cite[Theorem~2.4.2]{tao}, the empirical spectral measure $\hat\mu$ associated with $\hat J$ converges a.s.\ with respect to the vague topology to
	the law of $\beta X-\beta^2(1-q)$, where $X$ has the standard semicircular density. We note that this limit law has mass in any right vicinity of $-2\beta-\beta^2(1-q)$.
	
	Next we show that $\mu^{(k)}$ and $\hat\mu$ converge to the same limit as $N\to\infty$ followed by $k\to\infty$, and finally $h\to 0$. To this aim,  let $R^{(k)}(z)=(J^{(k)}-zI)^{-1}$ and $\hat R(z)=(\hat J-zI)^{-1}$ denote the resolvents of $J^{(k)}$ and $\hat J$, respectively. It suffices to show that the Stieltjes transforms $N^{-1}\tr \hat R(z)$ and $N^{-1}\tr R^{(k)}(z)$ of $\hat\mu$ and $\mu^{(k)}$, respectively, converge to the same limit in probability as $N\to\infty$ followed by $k\to\infty$ and $h\to 0$, pointwise for all $z\in\mathbb C$ with $\Im z>0$ (see e.g.~\cite[Section~2.4.3]{tao} for properties of the Stieltjes transform).
	By the resolvent identity,
	\begin{equation}
	\label{e:p:rep:Etr}
	N^{-1}\tr R^{(k)}(z)-N^{-1}\tr \hat R(z)=N^{-1}\tr R^{(k)}(z) \left( \hat J - J^{(k)} \right) \hat R(z).
	\end{equation}
	The $p$-Schatten norm of an $N\times N$ matrix $M$ whose eigenvalues are all real is defined by
	\begin{equation}
	\|M\|_p:=\left(\sum_{i=1}^N |\lambda_i(M)|^p\right)^{1/p} \text{ for } p\in[ 1,\infty),\quad
	\|M\|_\infty:= \max\left\{\left|\lambda_i(M)\right|: i=1,\ldots,N \right\}
	\end{equation}
	which satisfies $\|M\|_1\ge |\tr M|$ and the H\"older inequality. Hence, the expression in~\eqref{e:p:rep:Etr} is bounded in absolute value by
	\begin{equation}
	\label{e:p:trep:Holder}
	\|R^{(k)}(z)\|_\infty \|\hat R(z)\|_\infty \|\hat J\|_\infty N^{-1}\|{\rm diag } F_\cdot(\bfm^{(k)})^2\|_1,
	\end{equation}
	where we evaluated $\hat J-J^{(k)}$ using~\eqref{repr}. Each of the first two terms in~\eqref{e:p:trep:Holder} is bounded by $1/|\Im z|$, which follows from the definition of the resolvent. The third term $\| \hat J\|_\infty$ converges in probability to $2\beta+\beta^2(1-q)$ as $N\to\infty$ by Wigner's theorem in conjunction with edge scaling for Wigner matrices, see e.g. \cite[Chapter 3]{AGZ}, and \cite{AEKS20} for a recent reference.
	From the definitions of $F_i(\boldsymbol m^{(k)})$ and $m^{(k+1)}_i$, and as the function $x \mapsto \tanh^2(x)$ is 2-Lipschitz continuous, we obtain
	\begin{multline}
	\left| F_i(\boldsymbol m^{(k)})^2 - {m^{(k+1)}_i}^2 \right| =
	\bigg| \tanh^2\left(h+\frac{\beta}{\sqrt N} \sum_{j:j\ne i}  g_{ij} m^{(k)}_j -\beta^2(1-q)m^{(k)}_i\right)\\
	-\tanh^2\left(h+\frac{\beta}{\sqrt N} \sum_{j:j\ne i}  g_{ij} m^{(k)}_j -\beta^2(1-q)m^{(k-1)}_i\right)\bigg|
	\le 2\beta^2(1-q)\left| m^{(k)}_i - m^{(k-1)}_i \right|.
	\end{multline}
	Hence, by definition of the $1$-Schatten norm,
	\begin{equation}\label{e:p-TAP-Schatten}
		N^{-1}\|{\rm diag } F_\cdot(\bfm^{(k)})^2\|_1\le
		N^{-1}\sum_{i=1}^N \left( m^{(k)}_i \right)^2  +2 N^{-1}\sum_{i=1}^N \left| 1\cdot\left(m^{(k)}_i - m^{(k-1)}_i \right)\right|.
	\end{equation}
	The second term on the r.h.s.\ is bounded by
	\begin{equation}
		2N^{-1/2}\left[\sum_{i=1}^N \left( m^{(k)}_i - m^{(k-1)}_i \right)^2 \right]^{1/2}
	\end{equation}
	by the Cauchy-Schwarz inequality and thus converges to $0$ in probability as $N\to\infty$ followed by $k\to\infty$ by~\eqref{e:convit}.
	The first term on the r.h.s.\ of~\eqref{e:p-TAP-Schatten} equals $q_N$ and thus converges to $q$ in probability as $N\to\infty$ followed by $k\to\infty$.
	From~\eqref{e:q}, we obtain
	 \begin{equation}
	 	q=\mathbb \EE \tanh^2\left(h+\beta\sqrt{q}Z\right)\le h^2 + \beta^2 q\,,
	 \end{equation}
 	and
 	\begin{equation}
 		0\le q \le \frac{h^2}{1-\beta^2}\,,
 	\end{equation}
 	for $\beta\in(0,1)$, hence $h\to 0$ implies $q\to 0$.
	
	From the above, it follows that $\mu^{(k)}$ and $\hat\mu$ converge in probability to the same vague limit $\mu$ as $N\to\infty$ followed by $k\to\infty$ and $h\to 0$. For $h=0$ and $\beta\in(\sqrt{2}-1,1)$, we have $\mu(-\infty,-1)>0$ a.s.\ as a consequence of the first part of the proof. The assertion now follows from the vague convergence in probability of $\mu^{(k)}$ to $\mu$. 
\end{proof}
\section{Yet another algorithm, and simulations}%\FloatBarrier
The above considerations pertain to high temperature, but naturally lead to insights into the low temperature regime, where finding solutions of the TAP equations is notoriously hard. Indeed, Theorem \ref{t:rep}  suggests that the numerical instability of classical iteration schemes is ``nothing structural", but heavily depends on the way TAP equations are written. Under this light, the following iteration scheme, which we refer to as {\it $\vare$-\textsc{Banach}}, is quite natural: for $\varepsilon \in \R$  an additional free parameter, it reads
\begin{equation}\label{ebanach}
	m^{(k+1)}_i= \varepsilon m^{(k)}+ (1-\varepsilon) \tanh\left( h+\frac{\beta}{\sqrt N} \sum_{j:\: j\ne i} g_{ij} m^{(k)}_j - \beta^2 (1-q_N) m^{(k)}_i \right),
      \end{equation}
for $i=1 \dots N, k \in \N$. The idea\footnote{This scheme is well-known in the numerical literature: it has been implemented e.g. by Aspelmeier et. al \cite{am} to probe marginal stability of TAP solutions ``at the edge of chaos".} behind the ``$\varepsilon$-splitting" is of course to mitigate the impact of large negative Jacobian-eigenvalues.  We emphasize that:
\begin{itemize}
\item Here and henceforth we consider the finite-$N$ Onsager reaction term\footnote{Analogously for \textsc{2SteB}: the approximation \eqref{e:q} is of course wrong in low temperature.}  
\[ 
q_N := \frac{1}{N} \sum_{i=1}^N {m_j^{(k)}}^2.
\] 
\item Iterations are classical: unlike \eqref{eq:twostep}, no time-delay appears in Onsager's correction.
\end{itemize}

In this section we present a summary of numerical simulations run on the High Performance Computer \textsc{Elwetritsch}; additional material may be found in the supplementary material. 

Anticipating, \textsc{2SteB} is largely outperformed, in low temperature, by $\vare$-\textsc{Banach}: this is due to the fact that the former is extremely sensitive to the initialization, to the point of becoming eventually hopeless at finding stable fixed points; this limitation is not shared by the latter, in virtue (also) of the additional ``degree of freedom" $\varepsilon$ which may be tweaked at our discretion. 
\subsection{Calibrating $\vare$-\textsc{Banach}: high temperature} Given our limited theoretical understanding, effective choices of the $\vare$ parameter can only be found empirically. Here we show that a proper calibration does lead to an algorithmic performance  of $\vare$-\textsc{Banach} which is, in high temperature, fully comparable to \textsc{2SteB}. To do so, we appeal to the {\it TAP free energy}, TAP FE for short, which we recall is given by 
\beq \bea
& N f_{\text{TAP}}(\boldsymbol m) \equiv \frac{\be}{\sqrt{N}} \sum_{i < j}^N g_{ij} m_i m_j + h \sum_{i\leq N} m_i + \frac{\be^2}{2N} \sum_{i < j}^N (1-m_i^2)(1-m_j^2) - \sum_{i\leq N} I(m_i), 
\eea \eeq
where $\boldsymbol m \in [-1,1]^N$, and $I(x) \equiv \frac{1+x}{2} \log (1+x) + \frac{1-x}{2} \log(1-x)$.

\begin{figure}[h]
\includegraphics[height=4cm, trim={0.4cm 0.4cm 0 0},clip]{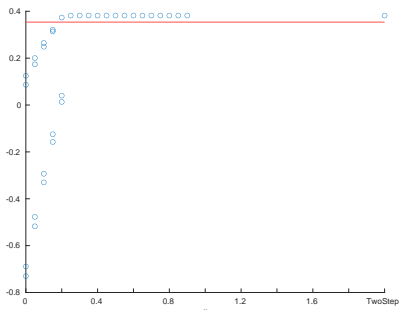}
\caption{TAP FE (y-axis) vs. $\vare$ (x-axis) for $\beta = 1$ and $h=0.5$ (way below AT-line). 
The red line is the FE computed from Parisi replica symmetric formula. 
We have run $\vare$-\textsc{Banach} for 1000 uniformly distributed start values: blue dots correspond to solutions which lie in the hypercube after $k= 1000$ iterations. The rightmost dot corresponds to the TAP FE of the (approximate) fixed point found by \textsc{2SteB}. For $\vare \approx 0.5$, say, the performance of $\vare$-\textsc{Banach} is thus fully equivalent to the \textsc{2SteB}. }
\label{calibrating}
\end{figure}

Critical points of the TAP FE are solutions of the TAP equations: in our simulations we have thus  computed the TAP FE of the fixed points found by both \textsc{2SteB}  and $\vare$-\textsc{Banach}. This is no simple task, for multiple reasons. First, the choice of the system's size is a priori not clear: {\it we have chosen this to be 
$N = 25$ throughout}. This might look at first sight unreasonably small, but numerical evidence rejects the objection. As a matter of fact, we show below that the (way) more delicate issue pertains to the number of initial values where iterations must be started for algorithms to find anything reasonable {\it at all}. 
\begin{itemize}
\item[{\sf C1.}] Even more challenging is the number of iterations, i.e.\ to decide whether these have stabilized. In this regard, our simulations (see also the supplementary material) suggest that the MSE \eqref{mse} is way too inaccurate a tool: we have thus discarded it in favor of the {\it maximum absolute error} between iterates, i.e. 
\beq
\text{MAE}\left( \boldsymbol m^{(k+1)}, \boldsymbol m^{(k)} \right) \equiv \max_{i=1}^N \left| m_i^{(k+1)} - m_i^{(k)} \right|\,,
\eeq 
thereby stopping iterations as soon as $\text{MAE} \leq 10^{-7}$.
\item[{\sf C2.}] Criteria for the validity of the TAP-Plefka expansions (which lead to the TAP FE) are vastly unknown. The only criterion which seems to be unanimously accepted is the one by Plefka 
\cite{plefka}. We thus require that (approximate) fixed points also satisfy Plefka's condition, i.e.
\beq \label{plefka}
\text{Plefka}\left( \boldsymbol m^{(k)} \right) \equiv \frac{\be^2}{N} \sum_{i=1}^N \left(1 - {m_i^{(k)}}^2 \right)^2 \leq 1. 
\eeq
\end{itemize}
Yet another key test for an approximate solution to be physically relevant is the associated TAP FE, which in turn should coincide\footnote{The low temperature SK-model allegedly requires $\infty$-many replica symmetry breakings (RSB). We use the Parisi FE obtained from a 2RSB approximation as the numerical error is mostly irrelevant \cite{mpv}.} with the Parisi FE \cite{pcs}. Figure \ref{calibrating} summarizes our findings concerning the calibration of $\vare$-\textsc{Banach} to yield TAP solutions which fulfill all above criteria. 
\subsection{\textsc{2SteB} vs. $\vare$-\textsc{Banach}: low temperature} Henceforth the physical parameters are $h=0.5$ and $\beta = 3$ (well above the AT-line). The size of the system is again $N= 25$. 

In low temperature, the issue of initialization becomes salient. Since Plefka criterion must be satisfied, this amounts to starting the algorithm {\it close to the facets of the hypercube}: by this we mean that we will consider $\boldsymbol m^{(0)}$ drawn uniformly in the subset of the hypercube where coordinates satisfy $m_\cdot^{(0)} \in [0.98, 1]\cup[-1, -0.98]$. 

A second issue is the $\vare$-calibration: the somewhat counterintuitive upshot is summarized in Figure \ref{fig:posorneg} below. Anticipating, one evinces that calibrating $\vare$-\textsc{Banach} with {\it positive} $\vare$ leads to a wealth of approximate fixed points  with unphysical TAP FE. We do not have any compelling explanation to this riddle, but it is tempting to believe that it finds its origin in the following: Plefka's criterion  \eqref{plefka}, which is known to be a necessary condition for convergence of the TAP-Plefka expansion, is possibly  one of (many?) still unknown criteria. Finally, $\vare$-\textsc{Banach} to negative $\vare$ possibly yields iterates that {\it temporarily} exit the hypercube: curiously, this appears to be an efficient way to overcome the repulsive nature of the fixed points which, as manifested by Figure \ref{rep_low_temp} below, becomes even more dramatic in low temperature. 
\begin{figure}%[h]
  \centering
   \includegraphics[height=4cm, width = 6cm, trim={0cm 0 0 0},clip]{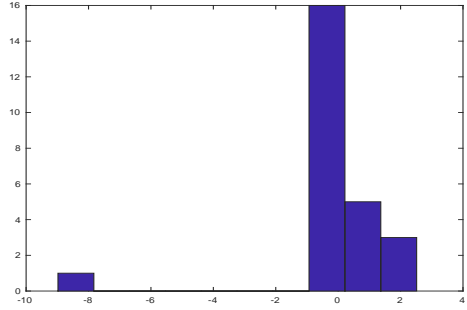}
  \caption{Eigenvalues of the Jacobian evaluated in one solution fulfilling {\sf C1-2)}. The Plefka value is $0.073206039645813 \ll 1$.} \label{rep_low_temp}
\end{figure}
\begin{figure}%[h]
     \centering
\begin{subfigure}{0.27\textwidth} \includegraphics[height=4cm, trim={0.45cm 0.4cm 0 0},clip]{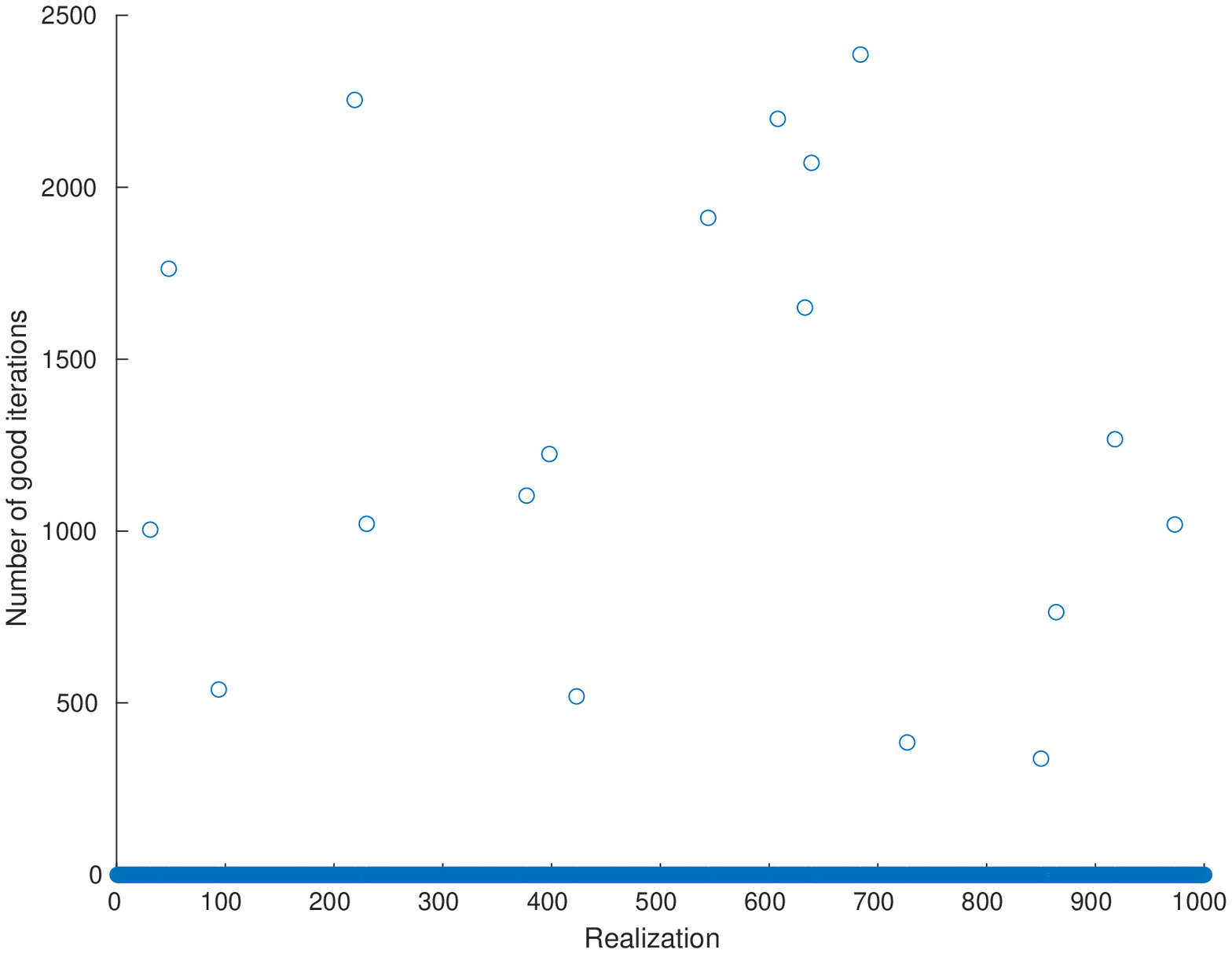}
    \caption{Number of iterations (y-axis) fulfilling criteria {\sf C1-2}) for negative $\vare =-0.705 - j \cdot 0.001, j =0,\ldots, 20$ vs. disorder realisations (x-axis). }
  \end{subfigure}\hfill
  \begin{subfigure}{0.27\textwidth} \includegraphics[height=4cm, trim={0.4cm 0.4cm 0 0},clip]{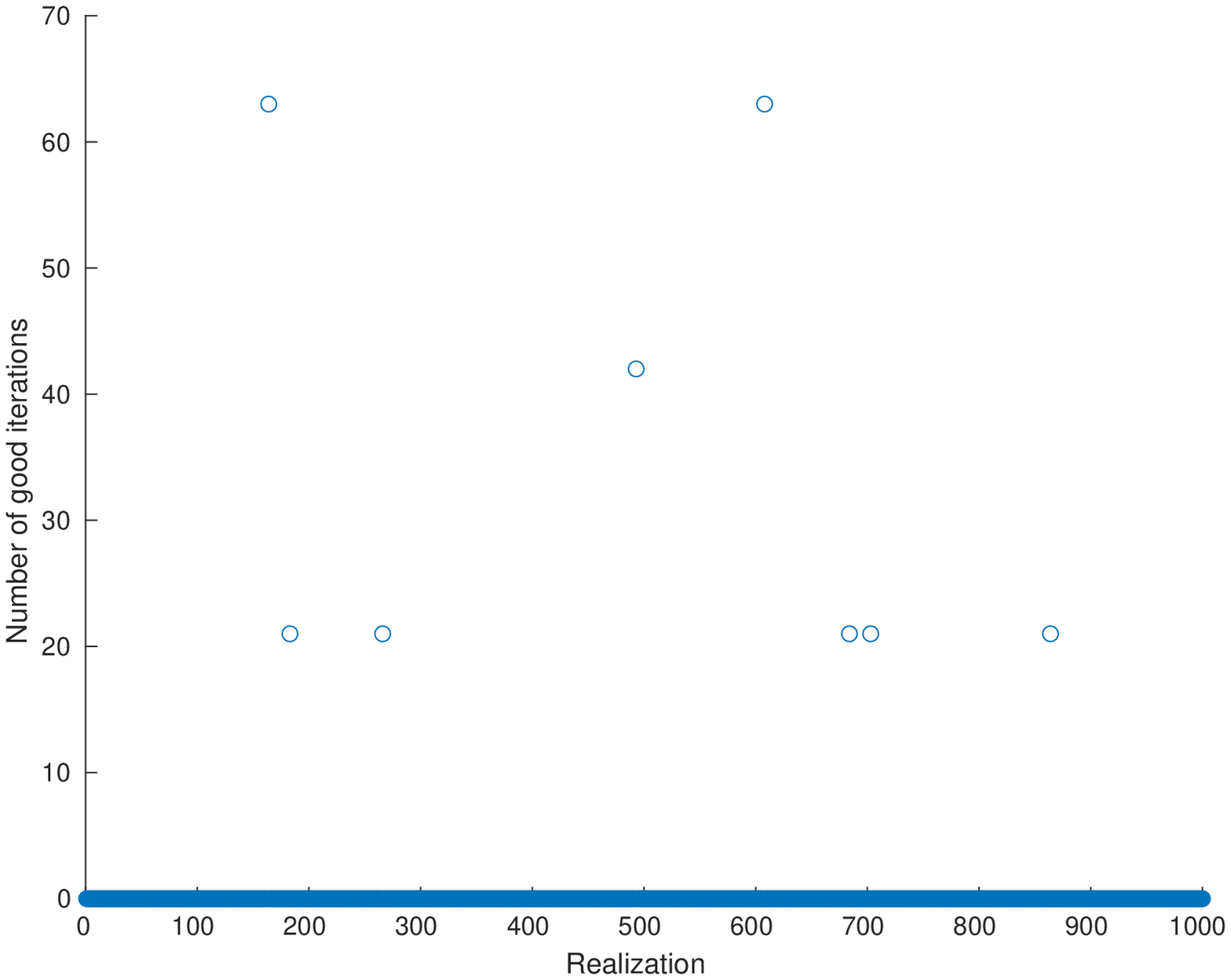}
    \caption{Number of iterations (y-axis) fulfilling criteria {\sf C1-2)} for positive $\varepsilon=0.705 + j \cdot 0.001, j =0,\ldots, 20$ vs. disorder realisations (x-axis).}
  \end{subfigure}%\vspace{1cm}
  \hfill
   \begin{subfigure}{0.27\textwidth} \includegraphics[height=4cm, trim={0cm 0.4cm 0 0},clip]{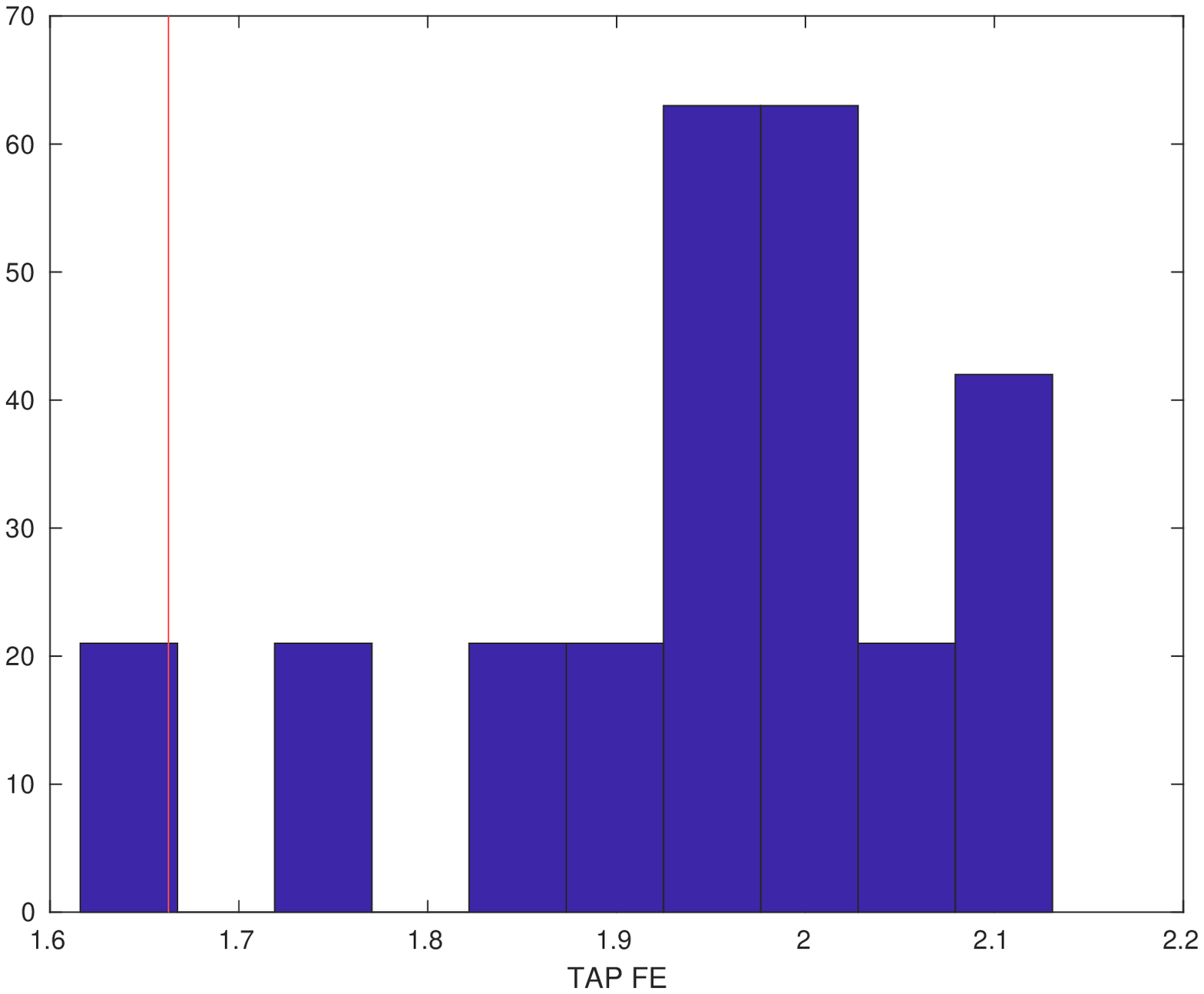}
    \caption{Histogram of TAP FE satisfying {\sf C1-2)} with positive $\varepsilon=0.705 + j \cdot 0.001, j=0,\ldots, 20$. The red line corresponds to the 2RSB Parisi FE.}
  \end{subfigure}
    \caption{$\vare$-\textsc{Banach}: 500 start values close to the facets of the hypercube, 1000 disorder realizations.  1000 iterations for each start value and $\varepsilon$. Since the Parisi FE is an upper-bound \cite{guerra}, fixed points to the right of the red line in (C) must be rejected.}
     \label{fig:posorneg}
  \end{figure} 
\begin{figure}%[h]
  \centering
  \begin{subfigure}{0.4\textwidth}
    \includegraphics[height=4cm, trim={0.4cm 0.45cm 0 0},clip]{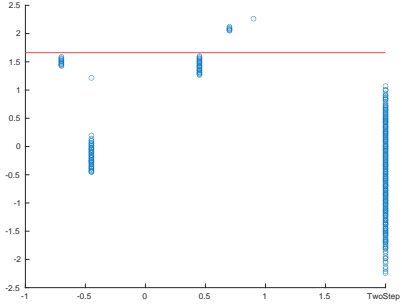}
    \caption{TAP FE of all iterations landing inside the hypercube. The rightmost line of blue dots corresponds to \textsc{2SteB}. Red is the Parisi FE (within the 2RSB approximation).} 
  \end{subfigure}
\hfil  
  %\hfill
   % \begin{subfigure}{0.4\textwidth}
   % \includegraphics[height=5cm]{preprintrsynccopy/final_simulations/starting_in_uniform_selected_corners_of_the_cube/banach_epsilon_and_twostep_l2norm_tap_fe/2tap_plefka_fullfilled.eps}
   % \caption{TAP free energy of all iterations ending with a \emph{good} solution.}
  %\end{subfigure}\hfill
 \begin{subfigure}{0.4\textwidth} \includegraphics[height=4cm, trim={0.4cm 0.45cm 0 0},clip]{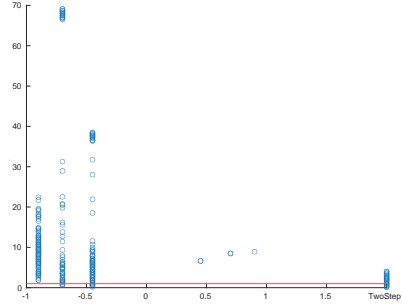}
    \caption{Plefka-values (y-axis) of all iterations (x-axis): the large values ($\text{Plefka}\left( \boldsymbol m^{(k)} \right) \gg 10$) correspond to (unviable) iterates falling out of the hypercube. The red line is at $y=1$.}
  \end{subfigure}
    \caption{$\vare$-\textsc{Banach} vs. \textsc{2SteB}: 1000 uniformly chosen initializations, 1000 iterations each. 
    One clearly evinces that \textsc{2SteB}  comes, in low temperature, to a stall: the algorithm doesn't even come close to a reasonable TAP FE (Figure A). In case of $\vare$-\textsc{Banach}, two choices lead to reasonable TAP free energies, but caution is needed: positive $\vare =0.5$ yields solutions violating Plefka's criterion (Figure B). Only negative $\vare \approx - 0.7$ lead to viable solutions: a numerical analysis of these is presented in Figure \ref{fig:beta3_1}.}\label{fig:sensitivity1}
  \end{figure}
\begin{figure}%[h]
  \centering
  \begin{subfigure}{0.27\textwidth}
    \includegraphics[height=4cm, trim={0.4cm 0.45cm 0 0},clip]{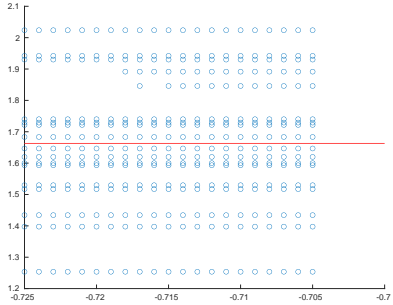}
    \caption{TAP FE (y-axis) vs. $\vare$ (x-axis).}
  \end{subfigure}\hfill
    \begin{subfigure}{0.27\textwidth}
    \includegraphics[height=4cm, trim={0.4cm 0.4cm 0 0},clip]{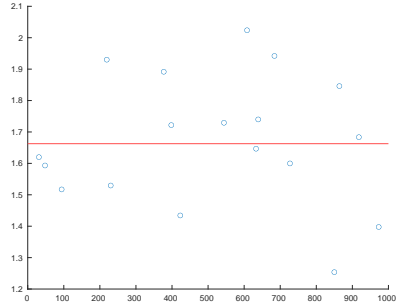}
   \caption{TAP FE (y-axis) of the realizations (x-axis).}
 \end{subfigure}%\vspace{1cm}
 \hfill
 \begin{subfigure}{0.27\textwidth} \includegraphics[height=4cm,trim={0cm 0.4cm 0 0},clip]{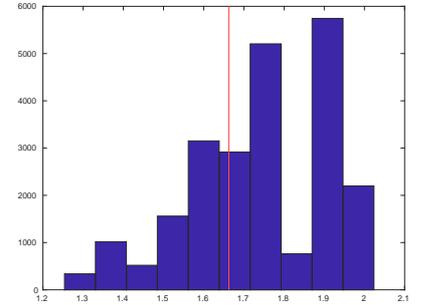}
    \caption{TAP FE Histogram. ${}$ ${}$${}$${}$}
  \end{subfigure}
    \caption{$\vare$-\textsc{Banach}: TAP FE of all solutions satisfying  {\sf C1-2)}, 500 initializations close to the facets of the hypercube, 1000 disorder realizations, $\varepsilon = -0.705 - 0.001 \cdot j$ for $j=0,\ldots, 20$. Red is the 2RSB Parisi FE.}\label{fig:beta3_1}
\end{figure}
%\FloatBarrier
\subsection{One, or more solutions?} 
The Parisi theory predicts, in low temperature, a large number of TAP solutions per realization. As can be evinced from Figure \ref{fig:onlyone} below,  numerical evidence for this is hard to come by: all solutions differ only after the 7th digit. 
\begin{figure}%[h]
    \centering
  \includegraphics[height=4cm, width = 7cm, trim={0.4cm 0.4cm 0 0},clip]{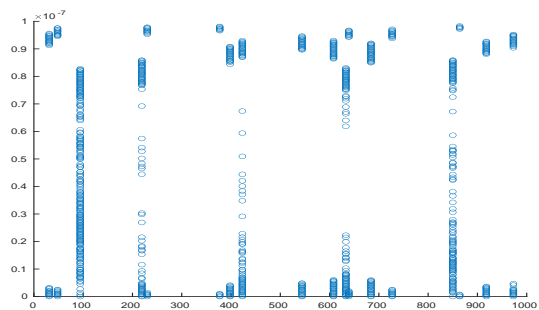}
    \caption{$\ell _\infty$-Distance (MAE) to a reference solution, {\sf C1-2)} being satisfied, found by $\vare$-\textsc{Banach} in 1000 realizations, each with 500 start values close to the facets of the hypercube, and $\varepsilon = -0.705 - 0.001 \cdot j, j=0,\ldots, 20$. }\label{fig:onlyone}
 \end{figure}
The crux of the matter seems to lie in the following: upon closer inspection, one finds that even after 1000 iterations, solutions (barely, but) still move. We have thus performed a more accurate numerical study, thereby recursively validating a fixed point as {\it new} whenever its $\ell_\infty$-distance to all previously validated is greater than $10^{-7}$. In order to facilitate the search for viable solutions we have also increased the accuracy of the $\vare$-mesh: the outcome is summarized in Figure \ref{fig:notonlyone} below. 
  \begin{figure}%[h] 
  \centering
  \begin{subfigure}{0.4\textwidth} \includegraphics[height=4cm, trim={0.4cm 0.4cm 0 0},clip]{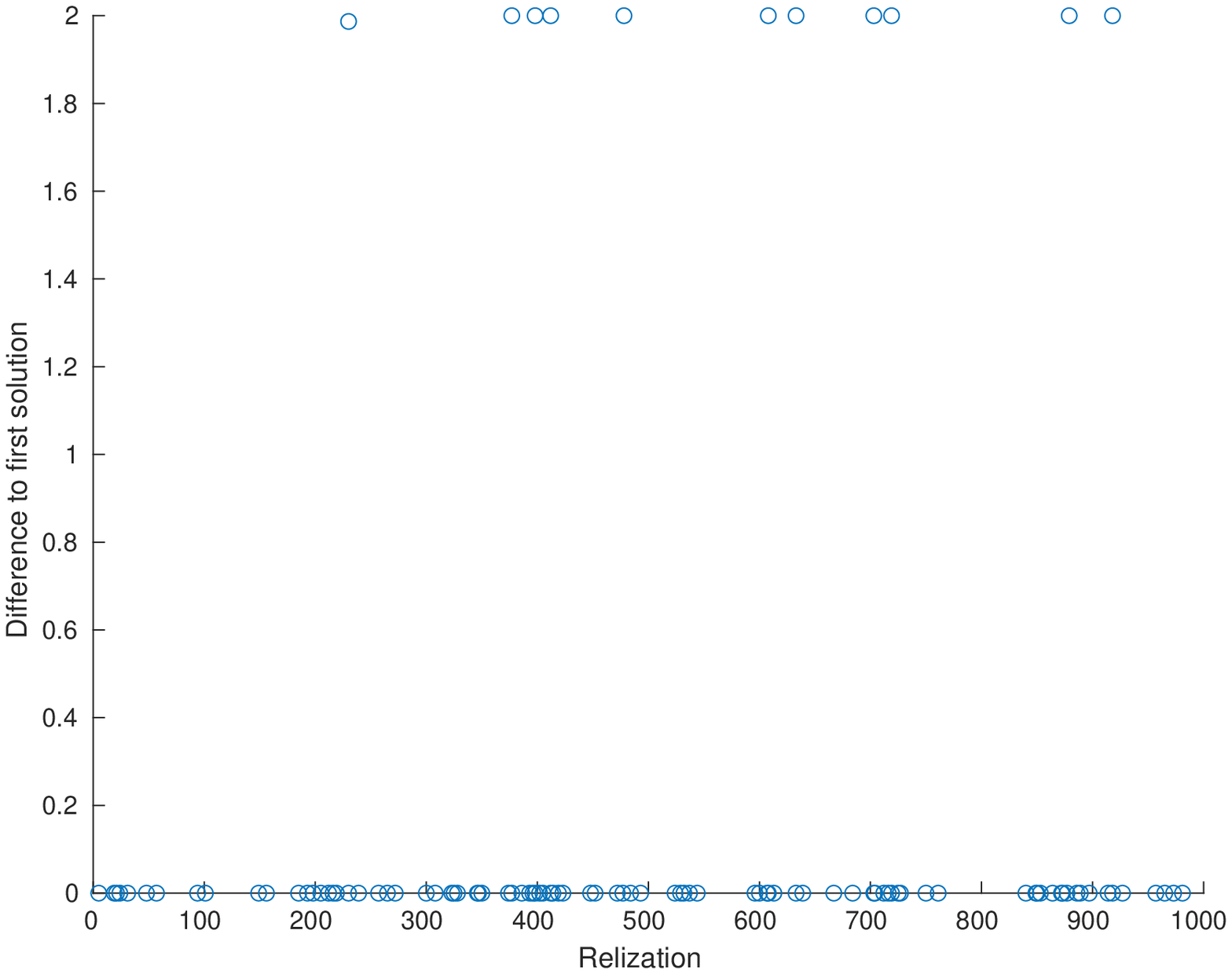}
    \caption{As in Fig. \ref{fig:onlyone}, but finer $\vare$-mesh.}
  \end{subfigure}
  \hfil
  \begin{subfigure}{0.4\textwidth} \includegraphics[height=3.8cm, width = 6cm, trim={0.4cm 0.4cm 0 0},clip]{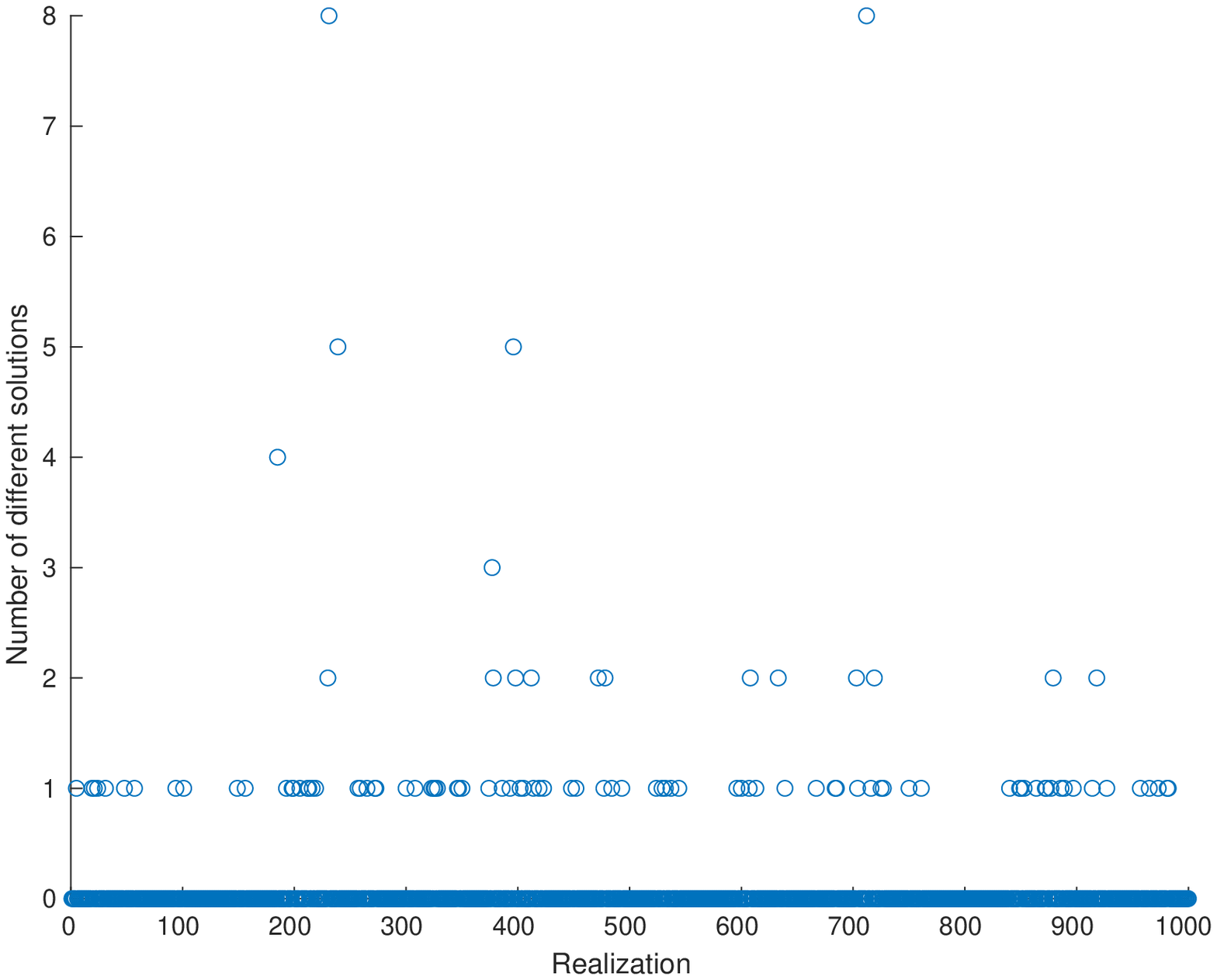}
    \caption{ \# solutions per realization.}
  \end{subfigure} 
    \caption{ Iterations satisfying {\sf C1-2)}, and $\varepsilon = -0.505 + 0.001 \cdot j, j \leq 250$.}\label{fig:notonlyone}
  \end{figure}
\begin{figure}%[h]
  \centering
  \begin{subfigure}{0.27\textwidth}
    \includegraphics[height=4cm, trim={0.4cm 0.45cm 0 0},clip]{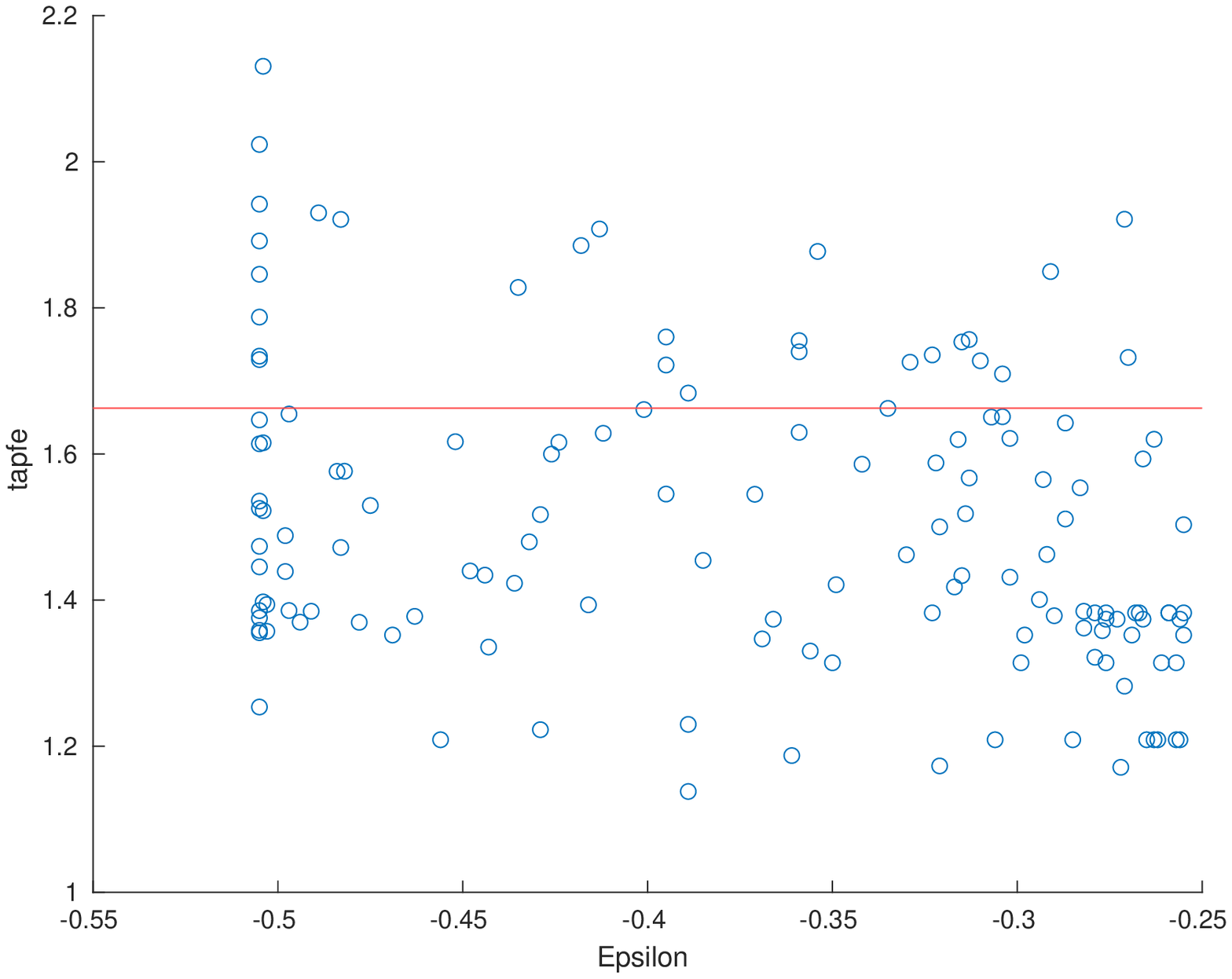}
    \caption{For each $\varepsilon$ all TAP free energy values of all different solutions found.}
  \end{subfigure}\hfill
    \begin{subfigure}{0.27\textwidth}
    \includegraphics[height=4cm, trim={0.4cm 0.4cm 0 0},clip]{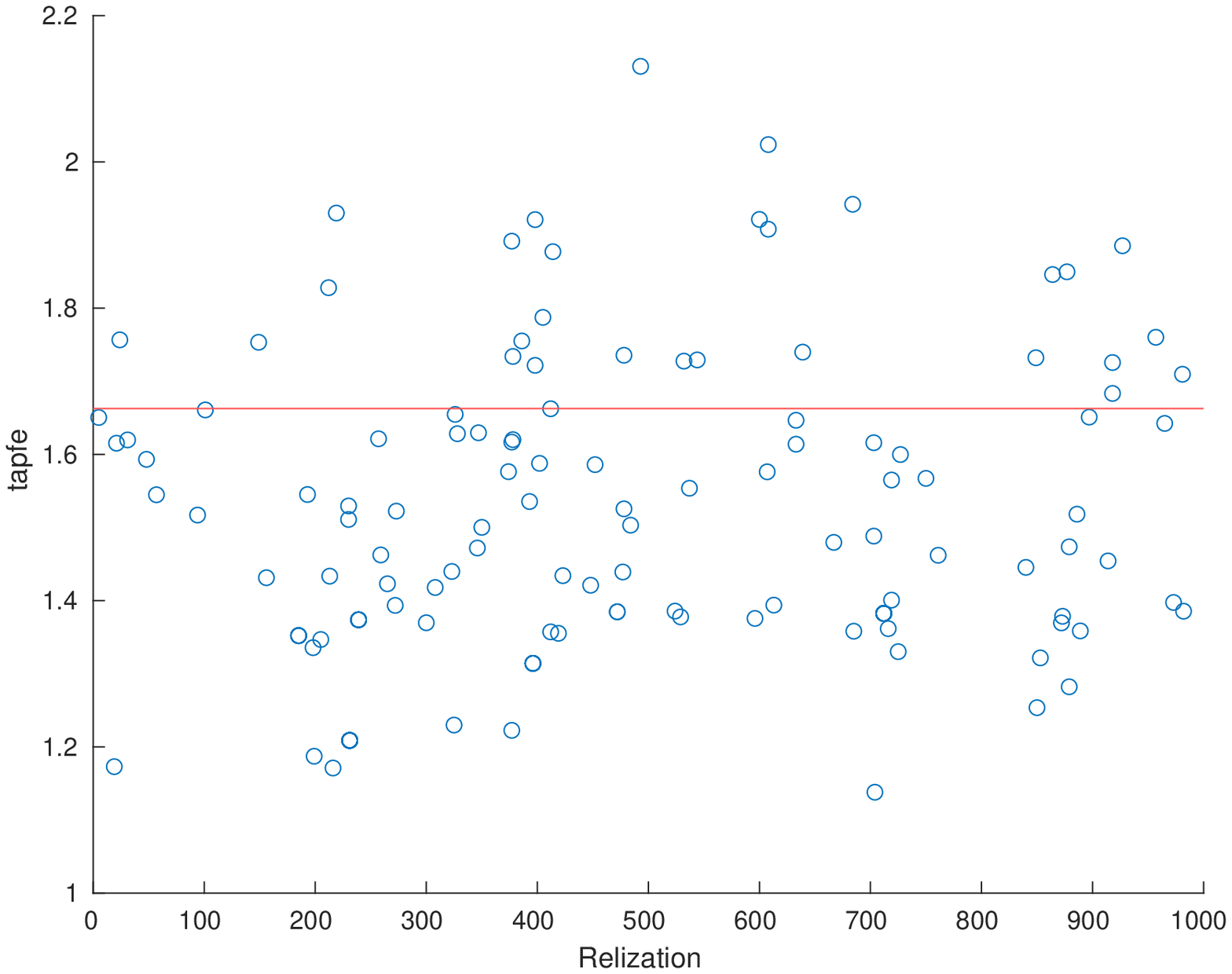}
   \caption{TAP FE of all realizations and all different solutions.}
 \end{subfigure}%\vspace{1cm}
 \hfill
 \begin{subfigure}{0.27\textwidth} \includegraphics[height=4cm,trim={0cm 0.4cm 0 0},clip]{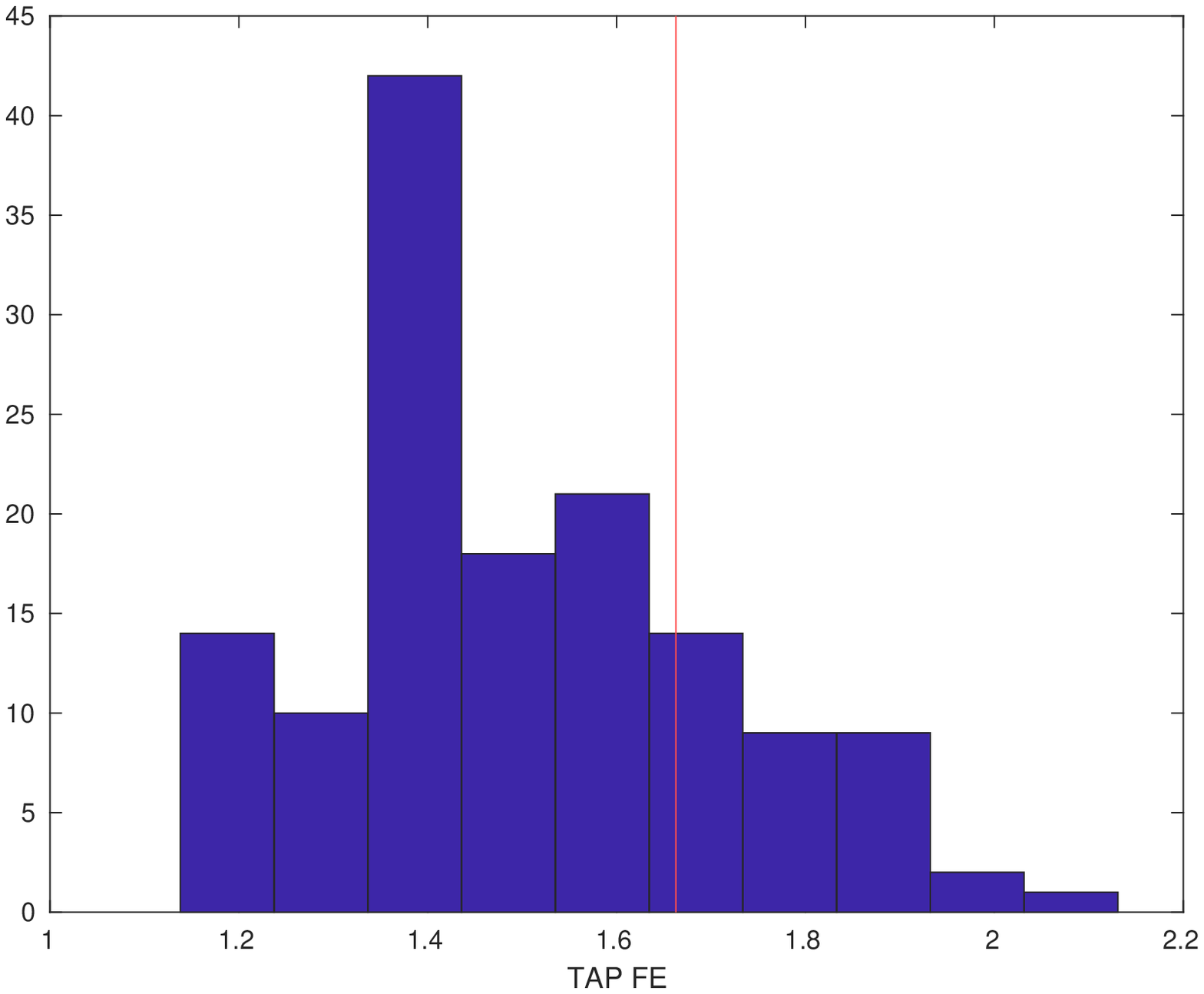}
    \caption{Histogram of the TAP free energy values of all different solutions found.}
  \end{subfigure}
    \caption{TAP FE as in Fig. \ref{fig:beta3_1}, but with $\varepsilon = -0.505 + 0.001 \cdot j, j \leq 250$.} \label{fig:last}
\end{figure}
All fixed points appearing in Figure \ref{fig:last} satisfy {\sf C1-2)}: one evinces that $\vare$-\textsc{Banach} has indeed found multiple solutions. However, many of these still yield TAP FE larger than the Parisi FE and must therefore be discarded. We believe this to be yet another instance of the aforementioned issue of unknown convergence criteria for the TAP-Plefka expansion.
% \FloatBarrier

\end{document}

\title{Supplementary material}

\author[S. Gufler]{Stephan Gufler}
\address{Stephan Gufler \\ J.W. Goethe-Universität Frankfurt, Germany.}
\email{gufler@math.uni-frankfurt.de}

\author[Jan Lukas Igelbrink]{Jan Lukas Igelbrink}
\address{Jan Lukas Igelbrink  \\ Institut f\"ur Mathematik, 
Johannes Gutenberg-Universit\"at Mainz, Germany.}
\email{jigelbri@uni-mainz.de}

\author[N. Kistler]{Nicola Kistler}
\address{Nicola Kistler, J.W. Goethe-Universit\"at Frankfurt, Germany.}
\email{kistler@math.uni-frankfurt.de}
\maketitle
\noindent
\tableofcontents
\newpage
\section{Parisi FE by replica symmetric formula}
\begin{figure}[H]
  \centering
  
  \includegraphics[width=\textwidth ,trim={0.45cm 0.4cm 0 0},clip]{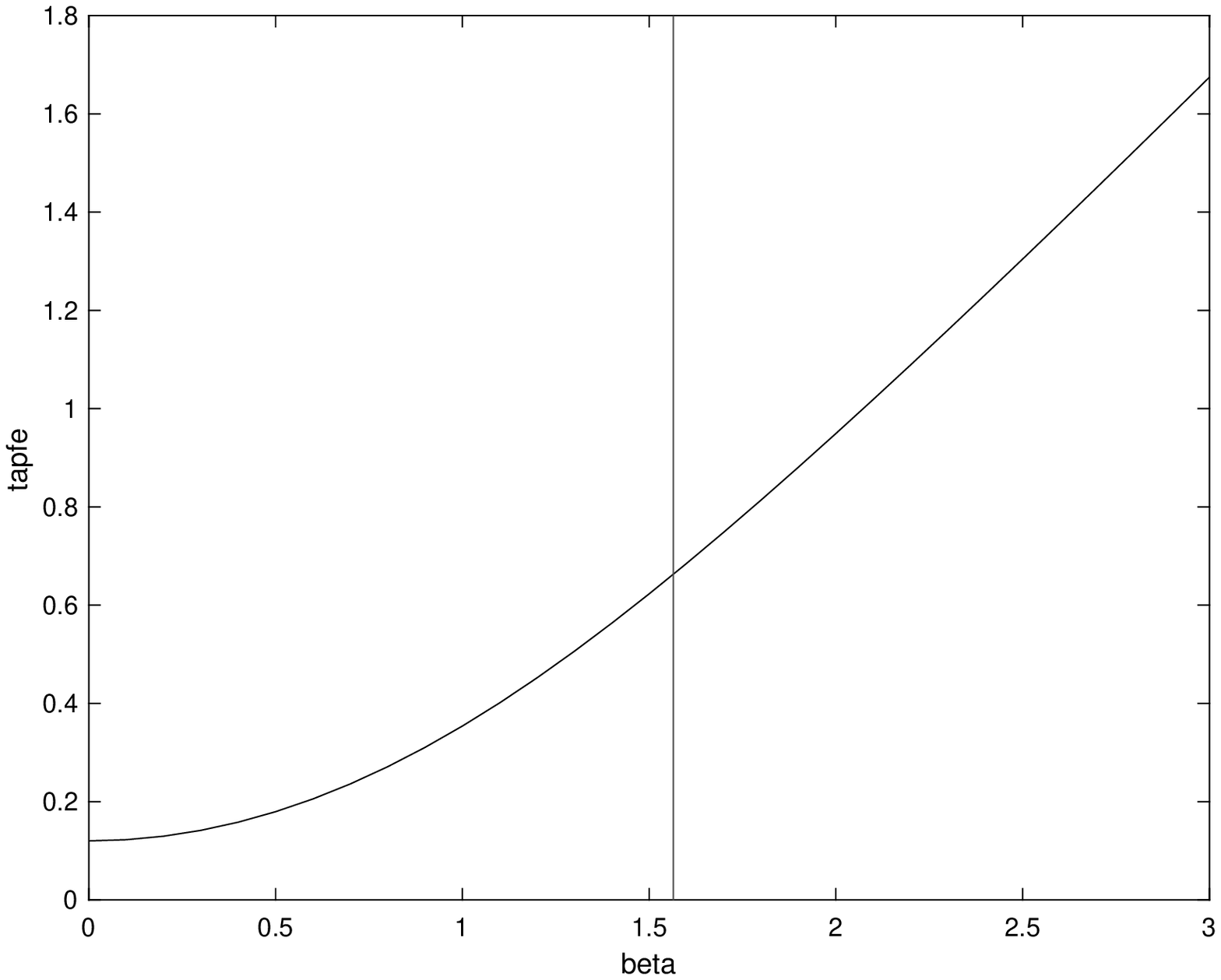}
  \caption{{\tiny} \space  TAP free energy by RS for $\beta \in [0, 3]$ and $h=0.5$}
\end{figure}
\FloatBarrier\clearpage
\section{Measuring the stability of our algorithms}
We call
\begin{equation}
  %\frac{1}{n}||m_{n}-m_{n-1}||_2^2
  \text{MSE}_N\left( \boldsymbol m^{(k+1)}, \boldsymbol m^{(k)}\right) \equiv \frac{1}{N} \sum_{i=1}^N \left( m_i^{(k+1)} - m_i^{(k)} \right)^2
\end{equation}
the \emph{mean squared error} between the two iterates $ \boldsymbol m^{(k+1)}$ and  $\boldsymbol m^{(k)}$. And
\begin{equation}
  % \max\left\{ |m^{(i)}_n-m^{(i)}_{n-1}| : i=1,\ldots, N \right\}
  \text{MAE}_N\left( \boldsymbol m^{(k+1)}, \boldsymbol m^{(k)} \right) \equiv \max_{i=1}^N \left| m_i^{(k+1)} - m_i^{(k)} \right|
\end{equation}
the {\it maximum absolute error}  between the two iterates $ \boldsymbol m^{(k+1)}$ and  $\boldsymbol m^{(k)}$.
\FloatBarrier\clearpage
\section{How we chose the start values}
We chose the start values uniformly out of the hypercube or uniformly out of all points inside the hypercube whose coordinates are in  $[0.98, 1]\cup[-1, -0.98]$. The second method will be called as \emph{chosen uniformly out of the corners of the hypercube}.
\FloatBarrier\clearpage
\section{Standard Banach algorithm}
\subsection{Using the {\it maximum absolute error}}
\subsubsection{Uniformly on the hypercube distributed start values}
\leavevmode
\begin{figure}[H]
  \centering
  \includegraphics[width=0.5\textwidth,  trim={0.45cm 0.4cm 0 0},clip]{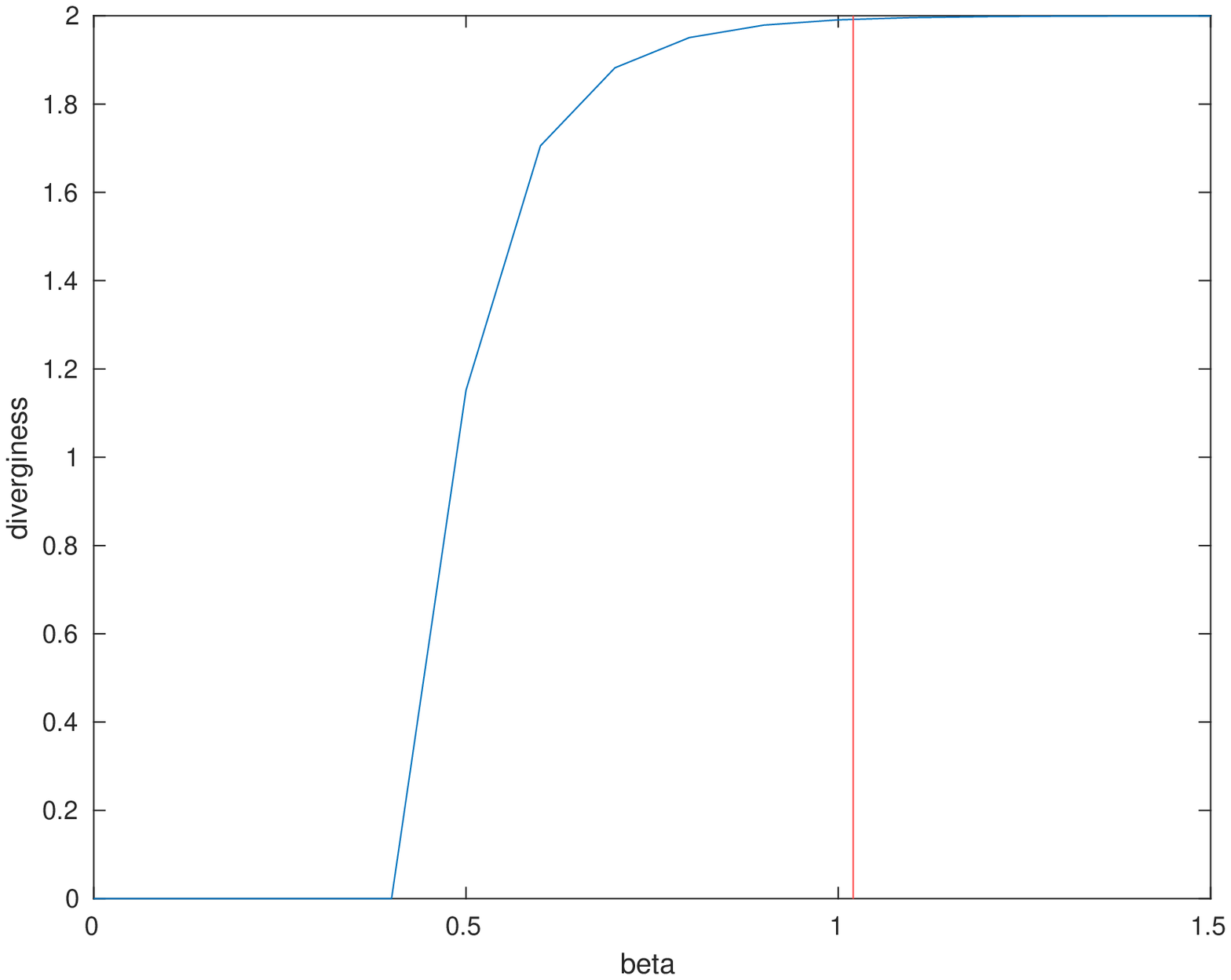}
  \caption{{\tiny} \space  1000 iterations of the standard Banach algorithm for $N=10$ and $h=0$. For $\beta\in [0, 1.5]$ we analyze the stability of this algorithm and see how quickly it takes off. To measure if the iterations stabilized we use the {\italic maximum absolute error} between the last two iterates. The vertical line marks the AT line.}
\end{figure}
\begin{figure}[H]
  \centering
  \includegraphics[width=0.5\textwidth,  trim={0.45cm 0.4cm 0 0},clip]{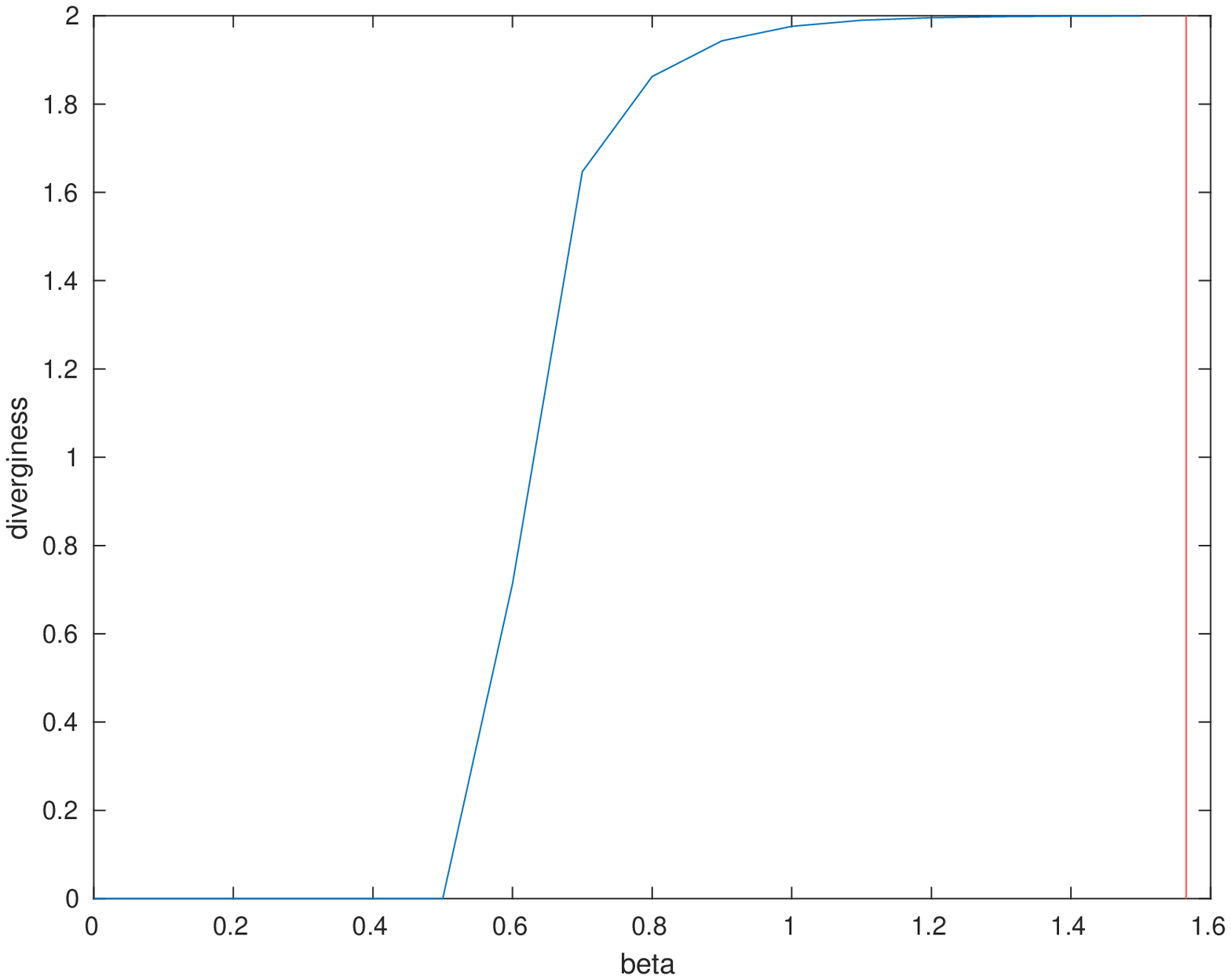}
  \caption{{\tiny} \space  1000 iterations of the standard Banach algorithm for $N=10$ and $h=0.5$. For $\beta\in [0, 1.5]$ we analyze the stability of this algorithm and see how quickly it takes off. To measure if the iterations stabilized we use the {\italic maximum absolute error} between the last two iterates. The vertical line marks the AT line.}
\end{figure}
\begin{figure}[H]
  \centering
  \includegraphics[width=0.5\textwidth,  trim={0.45cm 0.4cm 0 0},clip]{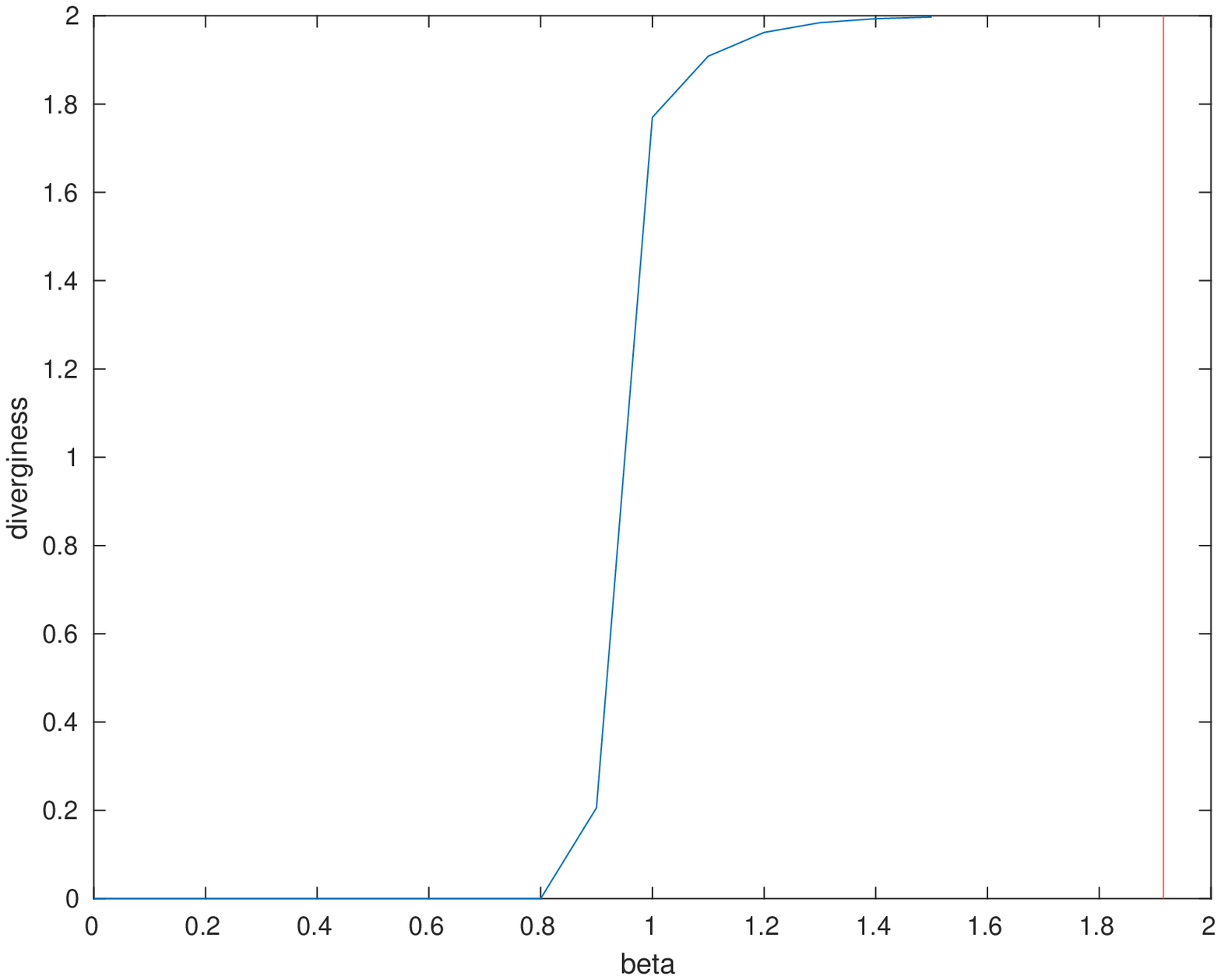}
  \caption{{\tiny} \space  1000 iterations of the standard Banach algorithm for $N=10$ and $h=1$. For $\beta\in [0, 1.5]$ we analyze the stability of this algorithm and see how quickly it takes off. To measure if the iterations stabilized we use the {\italic maximum absolute error} between the last two iterates. The vertical line marks the AT line.}
\end{figure}
\begin{figure}[H]
  \centering
  \includegraphics[width=0.5\textwidth,  trim={0.45cm 0.4cm 0 0},clip]{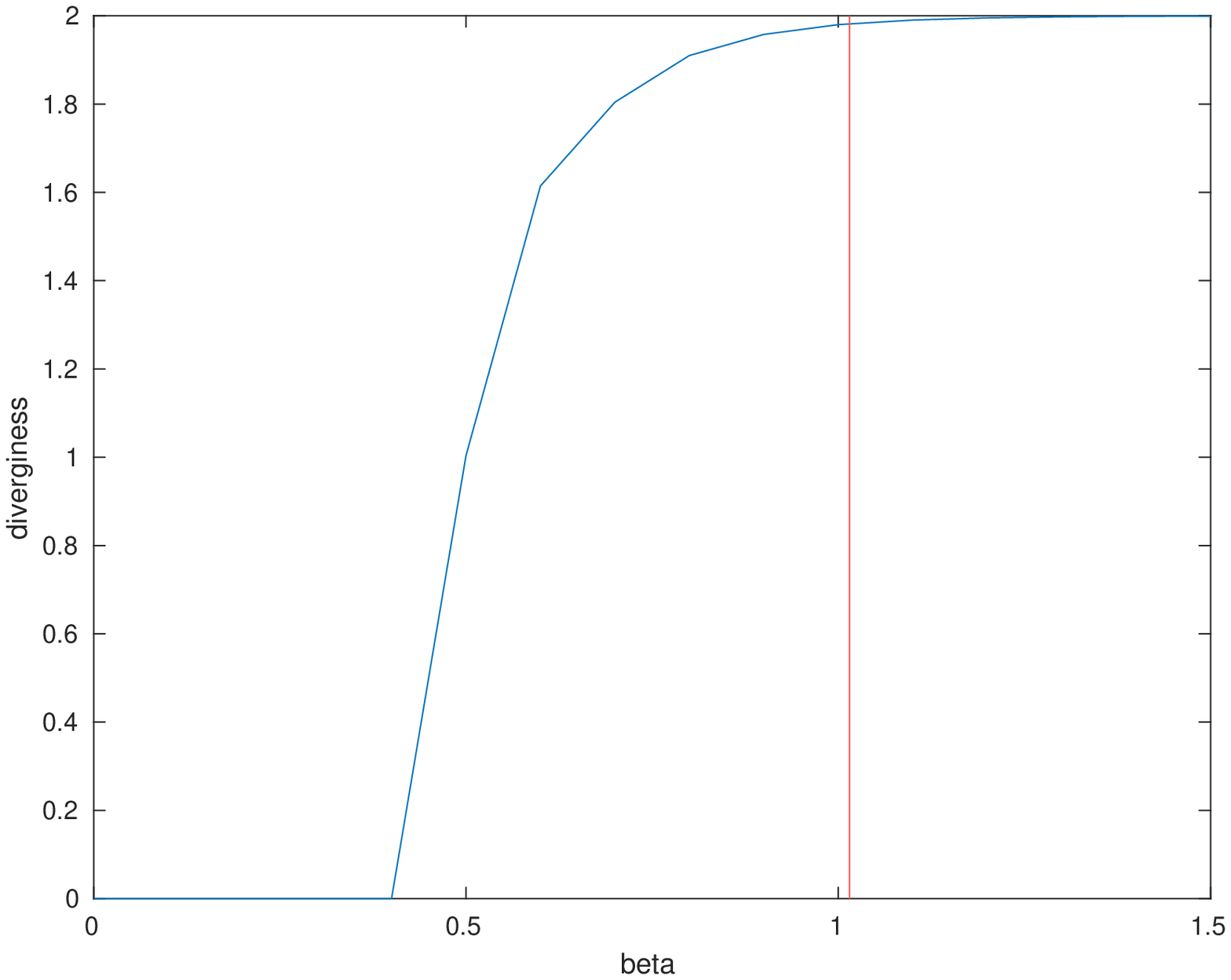}
  \caption{{\tiny} \space  1000 iterations of the standard Banach algorithm for $N=25$ and $h=0$. For $\beta\in [0, 1.5]$ we analyze the stability of this algorithm and see how quickly it takes off. To measure if the iterations stabilized we use the {\italic maximum absolute error} between the last two iterates. The vertical line marks the AT line.}
\end{figure}
\begin{figure}[H]
  \centering
  \includegraphics[width=0.5\textwidth,  trim={0.45cm 0.4cm 0 0},clip]{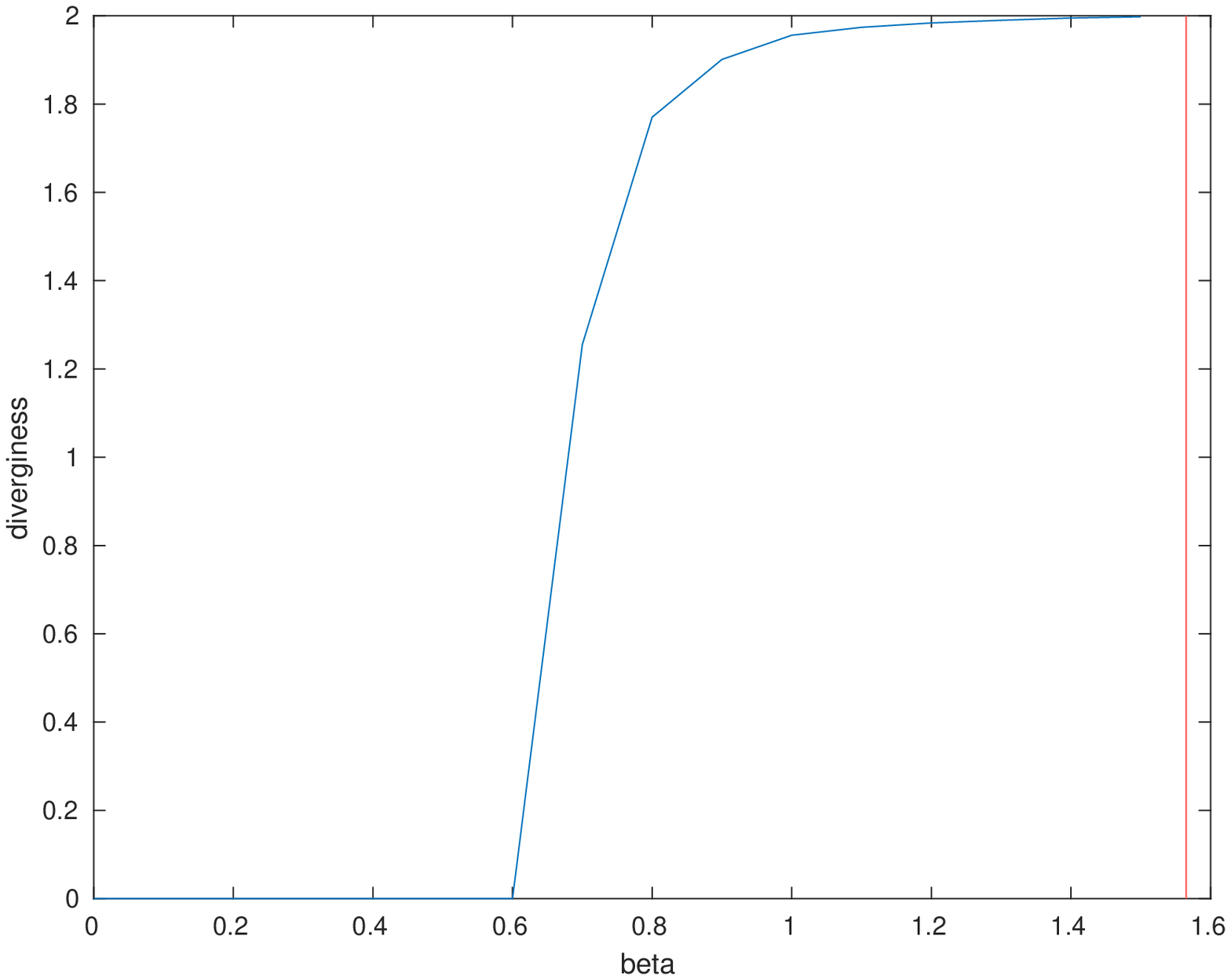}
  \caption{{\tiny} \space  1000 iterations of the standard Banach algorithm for $N=25$ and $h=0.5$. For $\beta\in [0, 1.5]$ we analyze the stability of this algorithm and see how quickly it takes off. To measure if the iterations stabilized we use the {\italic maximum absolute error} between the last two iterates. The vertical line marks the AT line.}
\end{figure}
\begin{figure}[H]
  \centering
  \includegraphics[width=0.5\textwidth,  trim={0.45cm 0.4cm 0 0},clip]{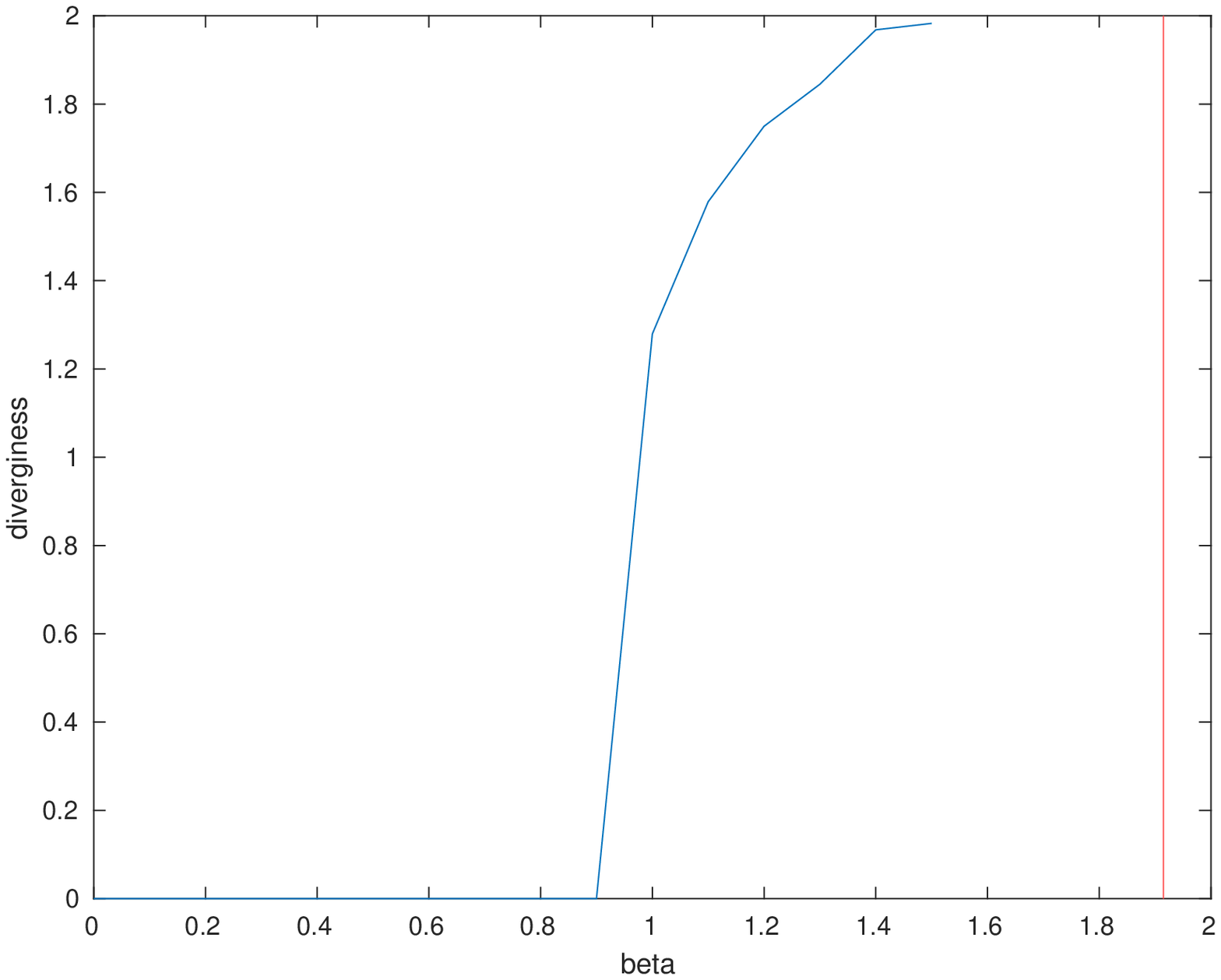}
  \caption{{\tiny} \space  1000 iterations of the standard Banach algorithm for $N=25$ and $h=1$. For $\beta\in [0, 1.5]$ we analyze the stability of this algorithm and see how quickly it takes off. To measure if the iterations stabilized we use the {\italic maximum absolute error} between the last two iterates. The vertical line marks the AT line.}
\end{figure}
\begin{figure}[H]
  \centering
  \includegraphics[width=0.5\textwidth,  trim={0.45cm 0.4cm 0 0},clip]{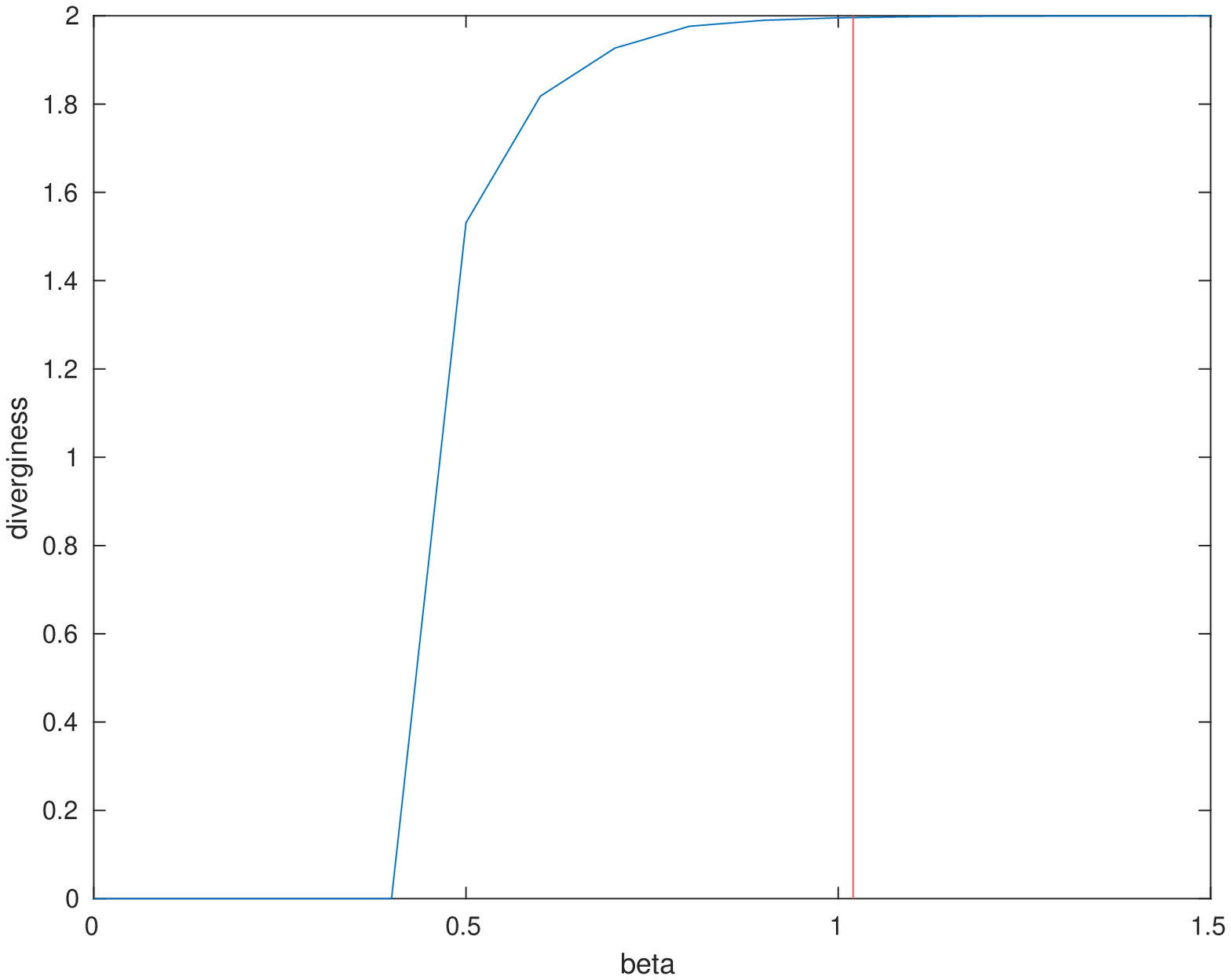}
  \caption{{\tiny} \space  1000 iterations of the standard Banach algorithm for $N=100$ and $h=0$. For $\beta\in [0, 1.5]$ we analyze the stability of this algorithm and see how quickly it takes off. To measure if the iterations stabilized we use the {\italic maximum absolute error} between the last two iterates. The vertical line marks the AT line.}
\end{figure}
\begin{figure}[H]
  \centering
  \includegraphics[width=0.5\textwidth,  trim={0.45cm 0.4cm 0 0},clip]{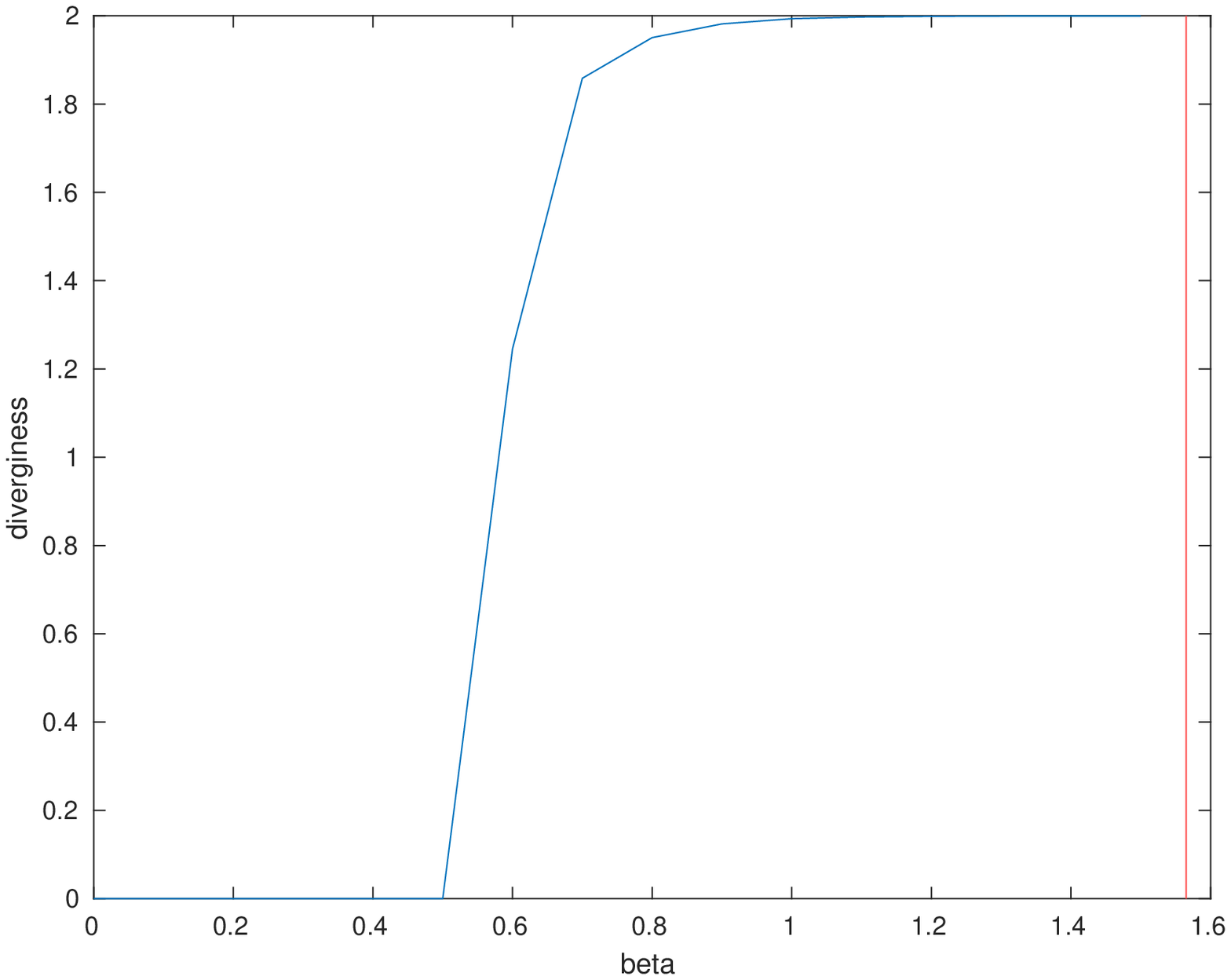}
  \caption{{\tiny} \space  1000 iterations of the standard Banach algorithm for $N=100$ and $h=0.5$. For $\beta\in [0, 1.5]$ we analyze the stability of this algorithm and see how quickly it takes off. To measure if the iterations stabilized we use the {\italic maximum absolute error} between the last two iterates. The vertical line marks the AT line.}
\end{figure}
\begin{figure}[H]
  \centering
  \includegraphics[width=0.5\textwidth,  trim={0.45cm 0.4cm 0 0},clip]{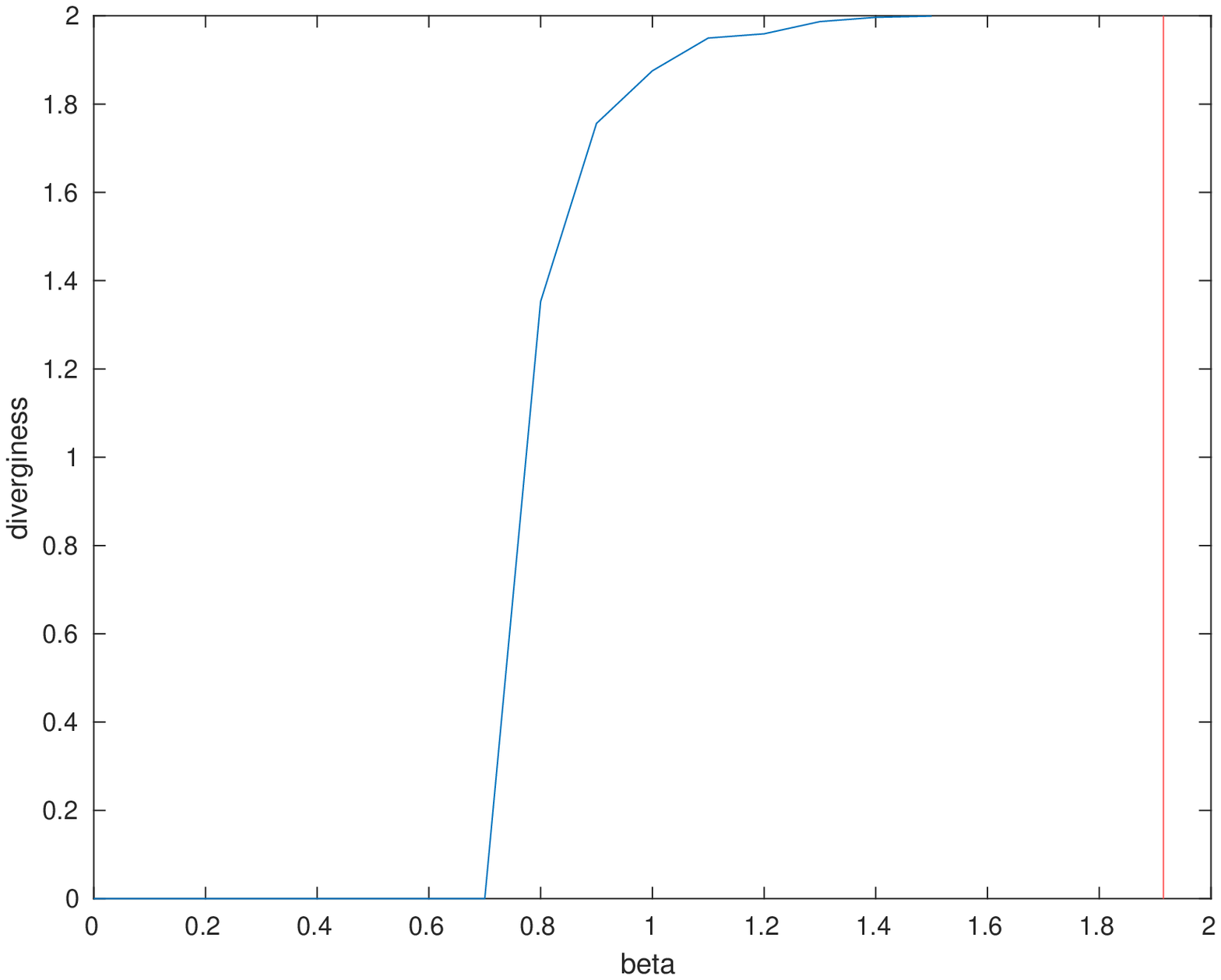}
  \caption{{\tiny} \space  1000 iterations of the standard Banach algorithm for $N=100$ and $h=1$. For $\beta\in [0, 1.5]$ we analyze the stability of this algorithm and see how quickly it takes off. To measure if the iterations stabilized we use the {\italic maximum absolute error} between the last two iterates. The vertical line marks the AT line.}
\end{figure}
\begin{figure}[H]
  \centering
  \includegraphics[width=0.5\textwidth,  trim={0.45cm 0.4cm 0 0},clip]{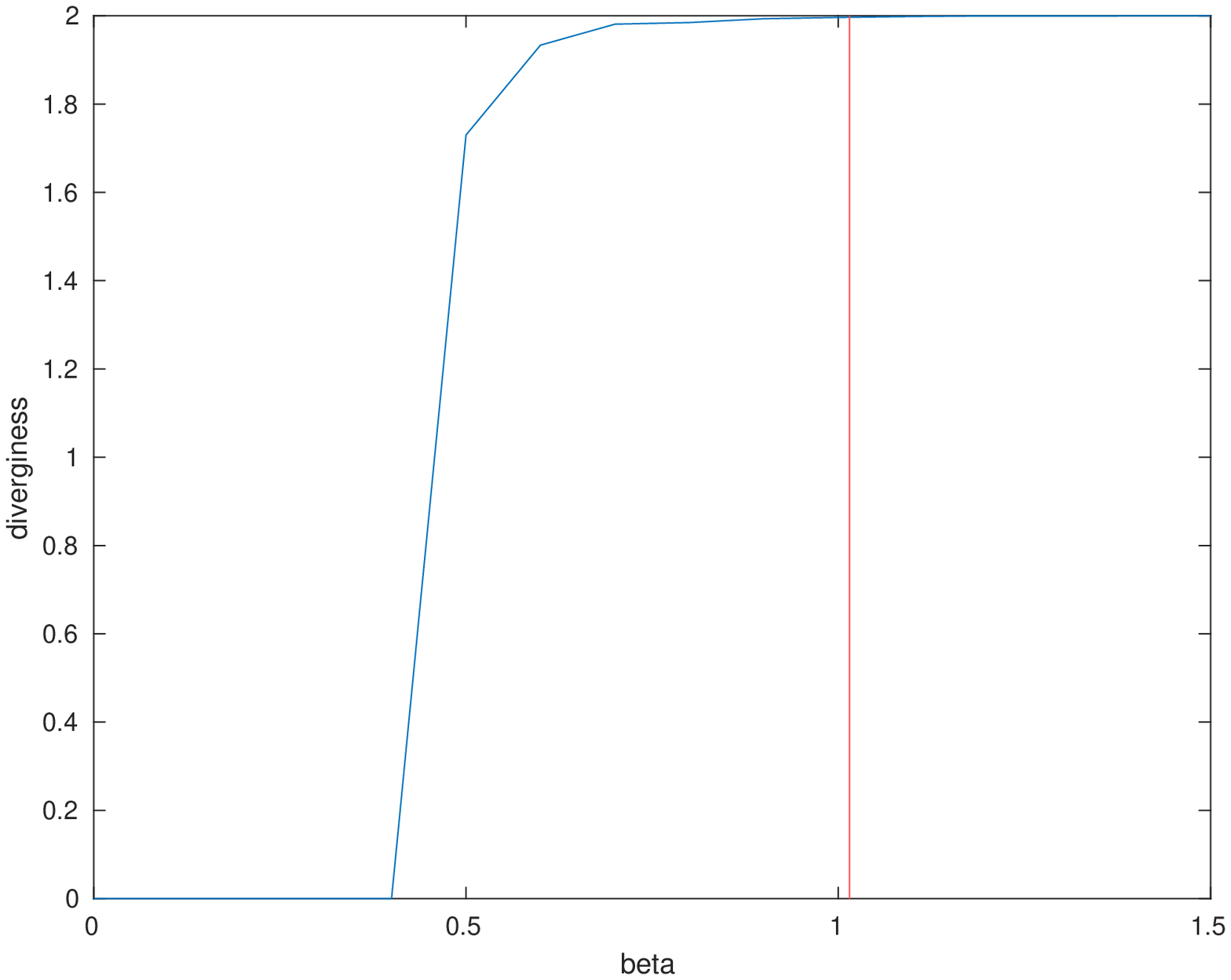}
  \caption{{\tiny} \space  1000 iterations of the standard Banach algorithm for $N=500$ and $h=0$. For $\beta\in [0, 1.5]$ we analyze the stability of this algorithm and see how quickly it takes off. To measure if the iterations stabilized we use the {\italic maximum absolute error} between the last two iterates. The vertical line marks the AT line.}
\end{figure}
\begin{figure}[H]
  \centering
  \includegraphics[width=0.5\textwidth,  trim={0.45cm 0.4cm 0 0},clip]{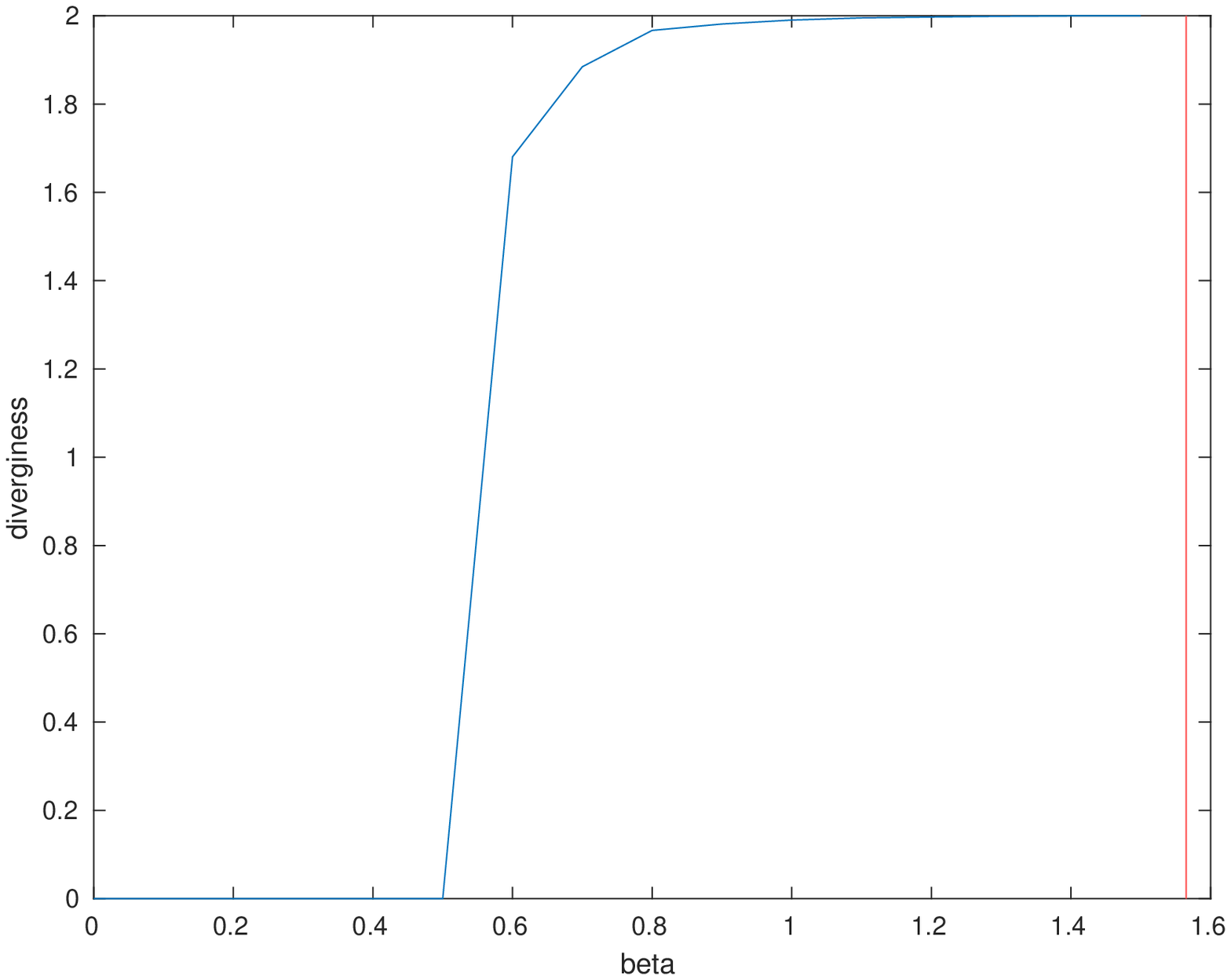}
  \caption{{\tiny} \space  1000 iterations of the standard Banach algorithm for $N=500$ and $h=0.5$. For $\beta\in [0, 1.5]$ we analyze the stability of this algorithm and see how quickly it takes off. To measure if the iterations stabilized we use the {\italic maximum absolute error} between the last two iterates. The vertical line marks the AT line.}
\end{figure}
\begin{figure}[H]
  \centering
  \includegraphics[width=0.5\textwidth,  trim={0.45cm 0.4cm 0 0},clip]{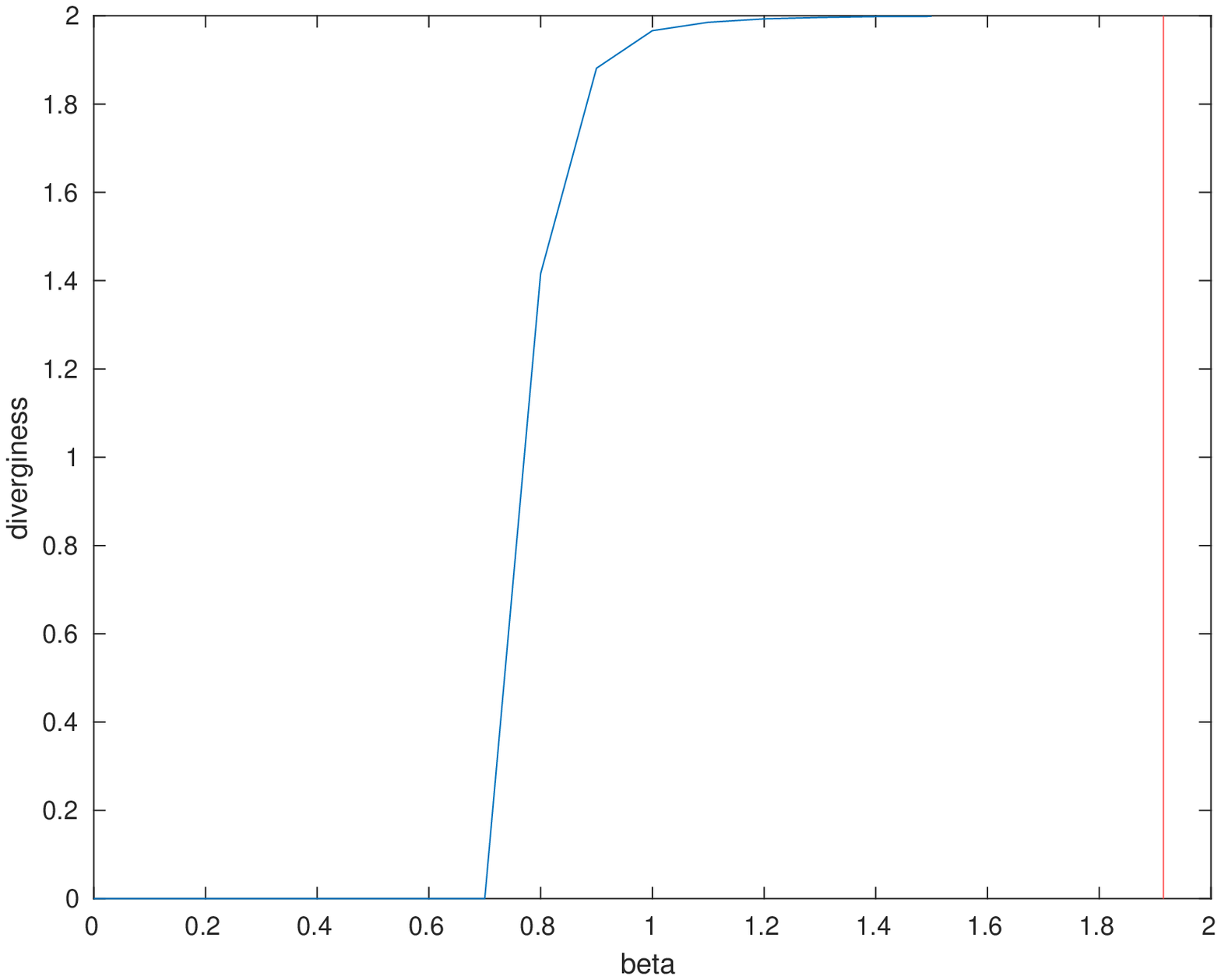}
  \caption{{\tiny} \space  1000 iterations of the standard Banach algorithm for $N=500$ and $h=1$. For $\beta\in [0, 1.5]$ we analyze the stability of this algorithm and see how quickly it takes off. To measure if the iterations stabilized we use the {\italic maximum absolute error} between the last two iterates. The vertical line marks the AT line.}
\end{figure}
\begin{figure}[H]
  \centering
  \includegraphics[width=0.5\textwidth,  trim={0.45cm 0.4cm 0 0},clip]{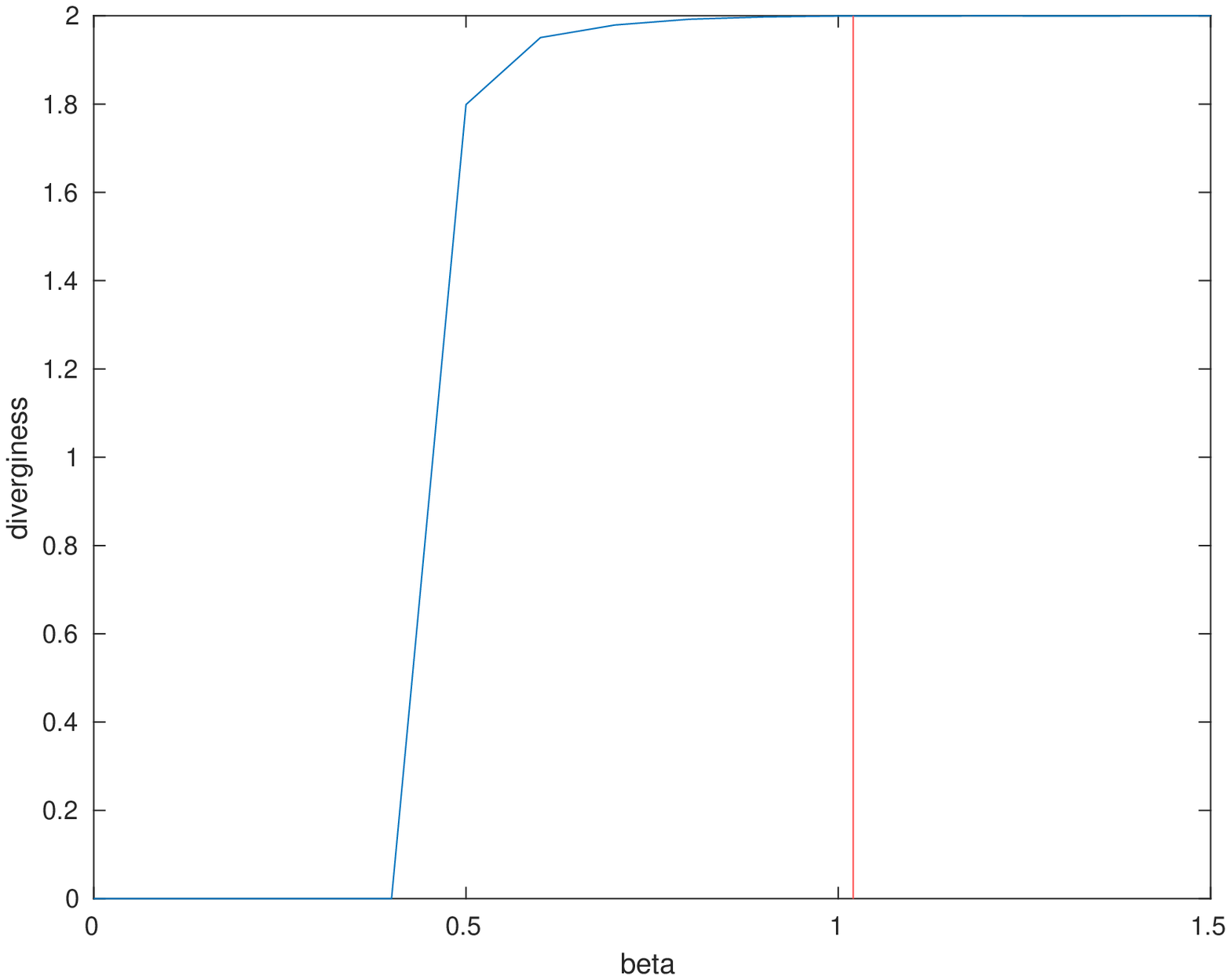}
  \caption{{\tiny} \space  1000 iterations of the standard Banach algorithm for $N=1000$ and $h=0$. For $\beta\in [0, 1.5]$ we analyze the stability of this algorithm and see how quickly it takes off. To measure if the iterations stabilized we use the {\italic maximum absolute error} between the last two iterates. The vertical line marks the AT line.}
\end{figure}
\begin{figure}[H]
  \centering
  \includegraphics[width=0.5\textwidth,  trim={0.45cm 0.4cm 0 0},clip]{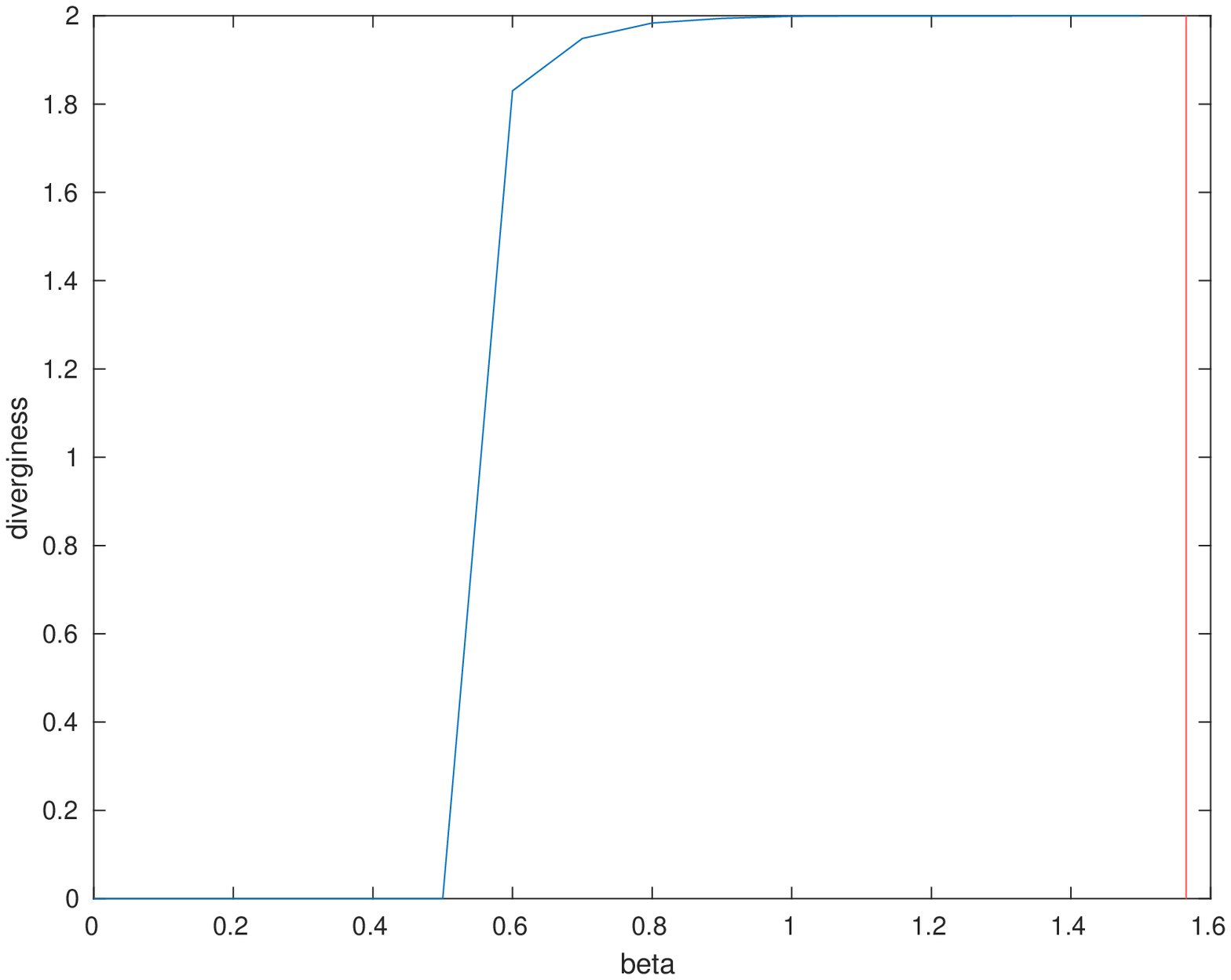}
  \caption{{\tiny} \space  1000 iterations of the standard Banach algorithm for $N=1000$ and $h=0.5$. For $\beta\in [0, 1.5]$ we analyze the stability of this algorithm and see how quickly it takes off. To measure if the iterations stabilized we use the {\italic maximum absolute error} between the last two iterates. The vertical line marks the AT line.}
\end{figure}
\begin{figure}[H]
  \centering
  \includegraphics[width=0.5\textwidth,  trim={0.45cm 0.4cm 0 0},clip]{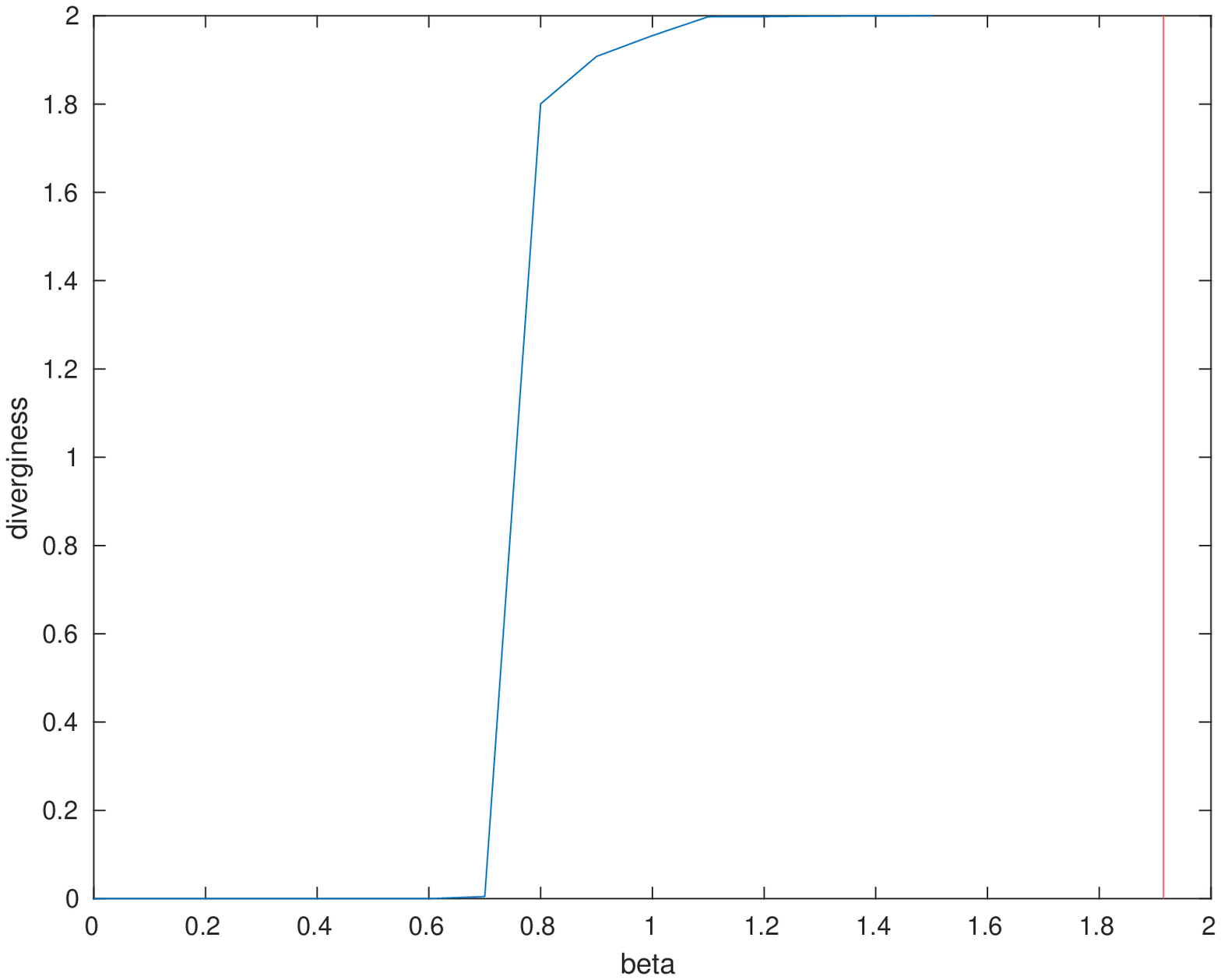}
  \caption{{\tiny} \space  1000 iterations of the standard Banach algorithm for $N=1000$ and $h=1$. For $\beta\in [0, 1.5]$ we analyze the stability of this algorithm and see how quickly it takes off. To measure if the iterations stabilized we use the {\italic maximum absolute error} between the last two iterates. The vertical line marks the AT line.}
\end{figure}\FloatBarrier\clearpage
\subsubsection{Uniformly in the corners of the hypercube  distributed start values}
\leavevmode
\begin{figure}[H]
  \centering
  \includegraphics[width=0.5\textwidth,  trim={0.45cm 0.4cm 0 0},clip]{rsynccopy/final_simulations/starting_uniform_in_the_hypercube/maxdiverginess_banach_explodes_at_line/1Banach.eps}
  \caption{{\tiny} \space  1000 iterations of the standard Banach algorithm for $N=10$ and $h=0$. For $\beta\in [0, 1.5]$ we analyze the stability of this algorithm and see how quickly it takes off. To measure if the iterations stabilized we use the {\italic maximum absolute error} between the last two iterates. The vertical line marks the AT line.}
\end{figure}
\begin{figure}[H]
  \centering
  \includegraphics[width=0.5\textwidth,  trim={0 0 0 0},clip]{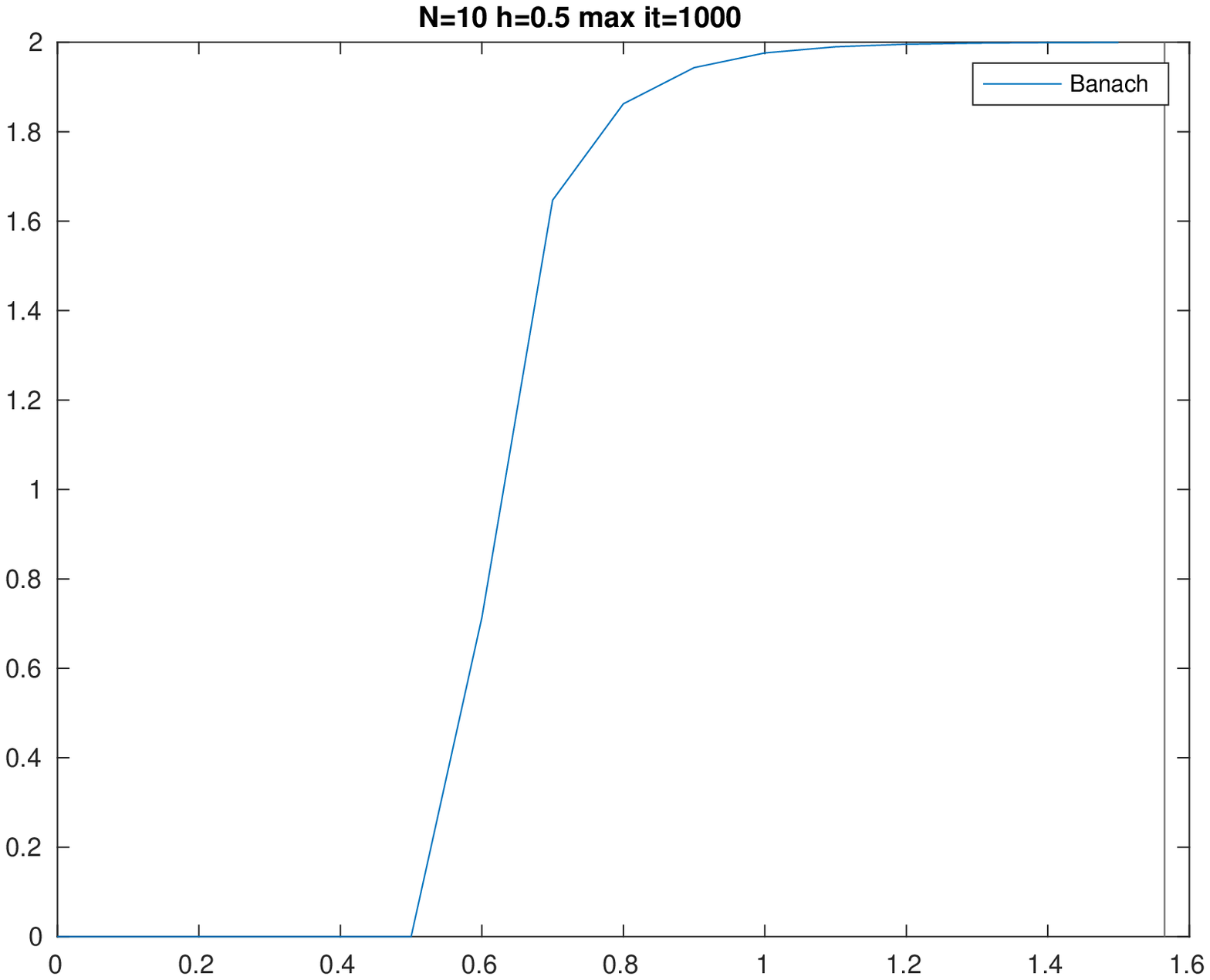}
  \caption{{\tiny} \space  1000 iterations of the standard Banach algorithm for $N=10$ and $h=0.5$. For $\beta\in [0, 1.5]$ we analyze the stability of this algorithm and see how quickly it takes off. To measure if the iterations stabilized we use the {\italic maximum absolute error} between the last two iterates. The vertical line marks the AT line.}
\end{figure}
\begin{figure}[H]
  \centering
  \includegraphics[width=0.5\textwidth,  trim={0 0 0 0},clip]{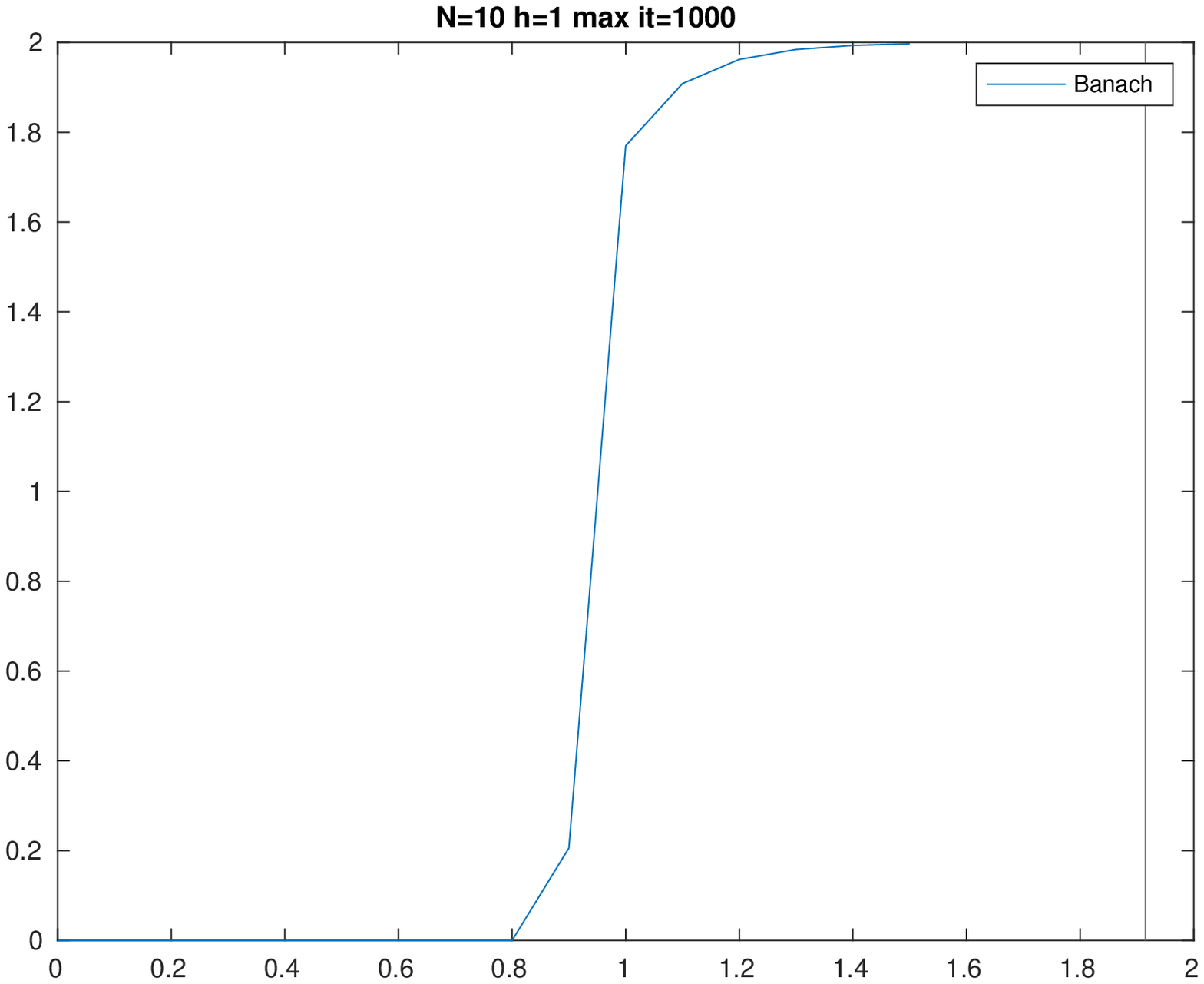}
  \caption{{\tiny} \space  1000 iterations of the standard Banach algorithm for $N=10$ and $h=1$. For $\beta\in [0, 1.5]$ we analyze the stability of this algorithm and see how quickly it takes off. To measure if the iterations stabilized we use the {\italic maximum absolute error} between the last two iterates. The vertical line marks the AT line.}
\end{figure}
\begin{figure}[H]
  \centering
  \includegraphics[width=0.5\textwidth,  trim={0 0 0 0},clip]{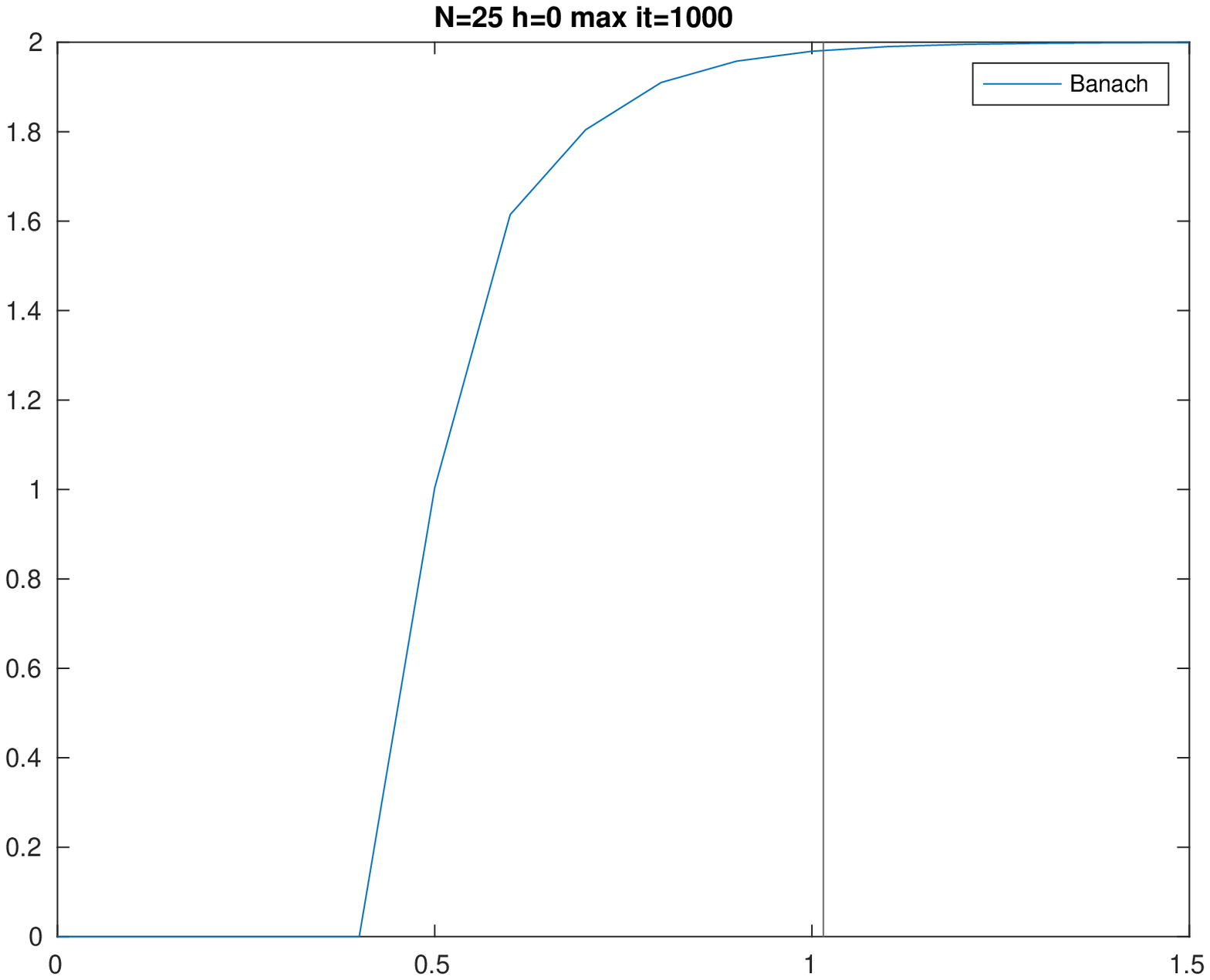}
  \caption{{\tiny} \space  1000 iterations of the standard Banach algorithm for $N=25$ and $h=0$. For $\beta\in [0, 1.5]$ we analyze the stability of this algorithm and see how quickly it takes off. To measure if the iterations stabilized we use the {\italic maximum absolute error} between the last two iterates. The vertical line marks the AT line.}
\end{figure}
\begin{figure}[H]
  \centering
  \includegraphics[width=0.5\textwidth,  trim={0 0 0 0},clip]{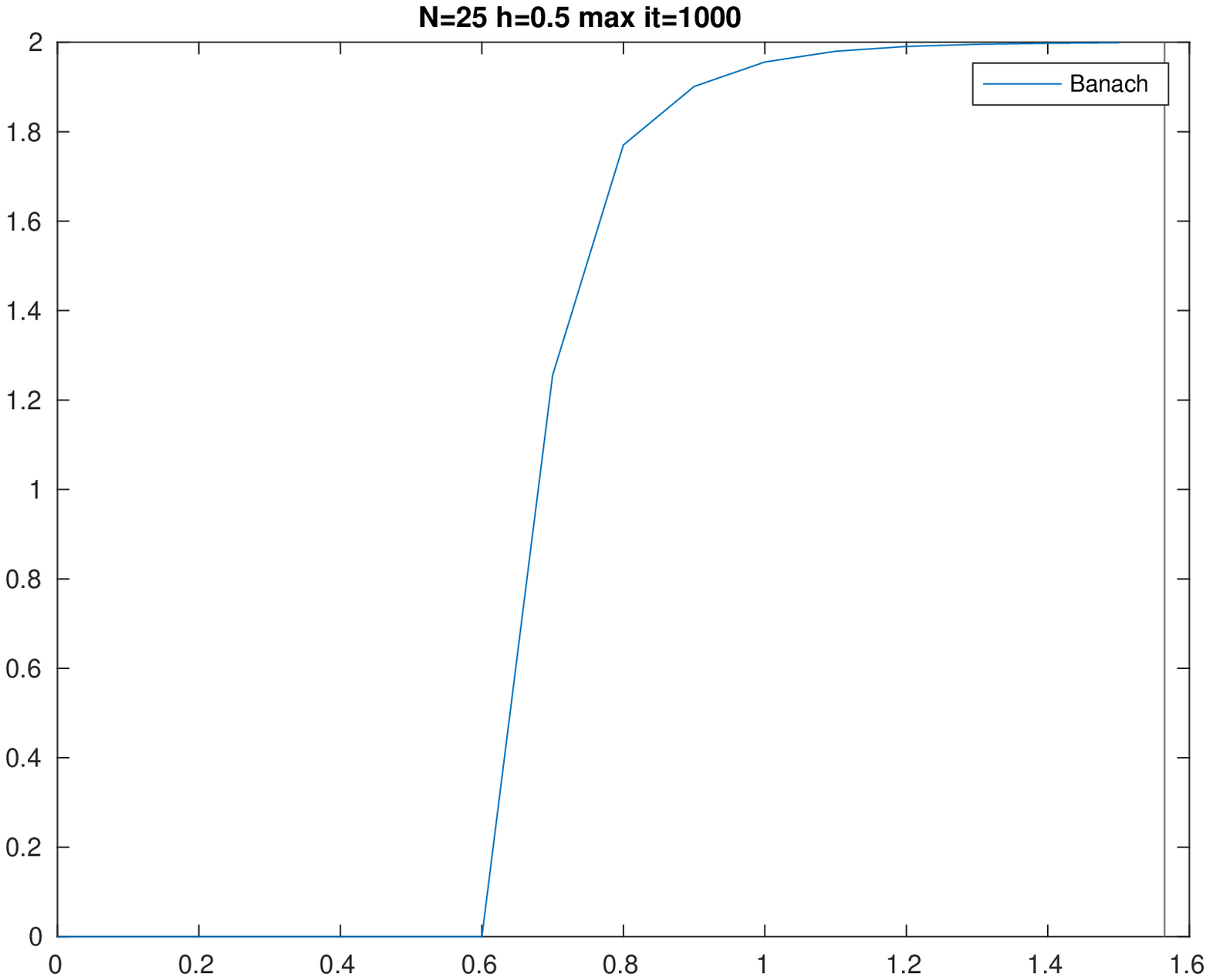}
  \caption{{\tiny} \space  1000 iterations of the standard Banach algorithm for $N=25$ and $h=0.5$. For $\beta\in [0, 1.5]$ we analyze the stability of this algorithm and see how quickly it takes off. To measure if the iterations stabilized we use the {\italic maximum absolute error} between the last two iterates. The vertical line marks the AT line.}
\end{figure}
\begin{figure}[H]
  \centering
  \includegraphics[width=0.5\textwidth,  trim={0 0 0 0},clip]{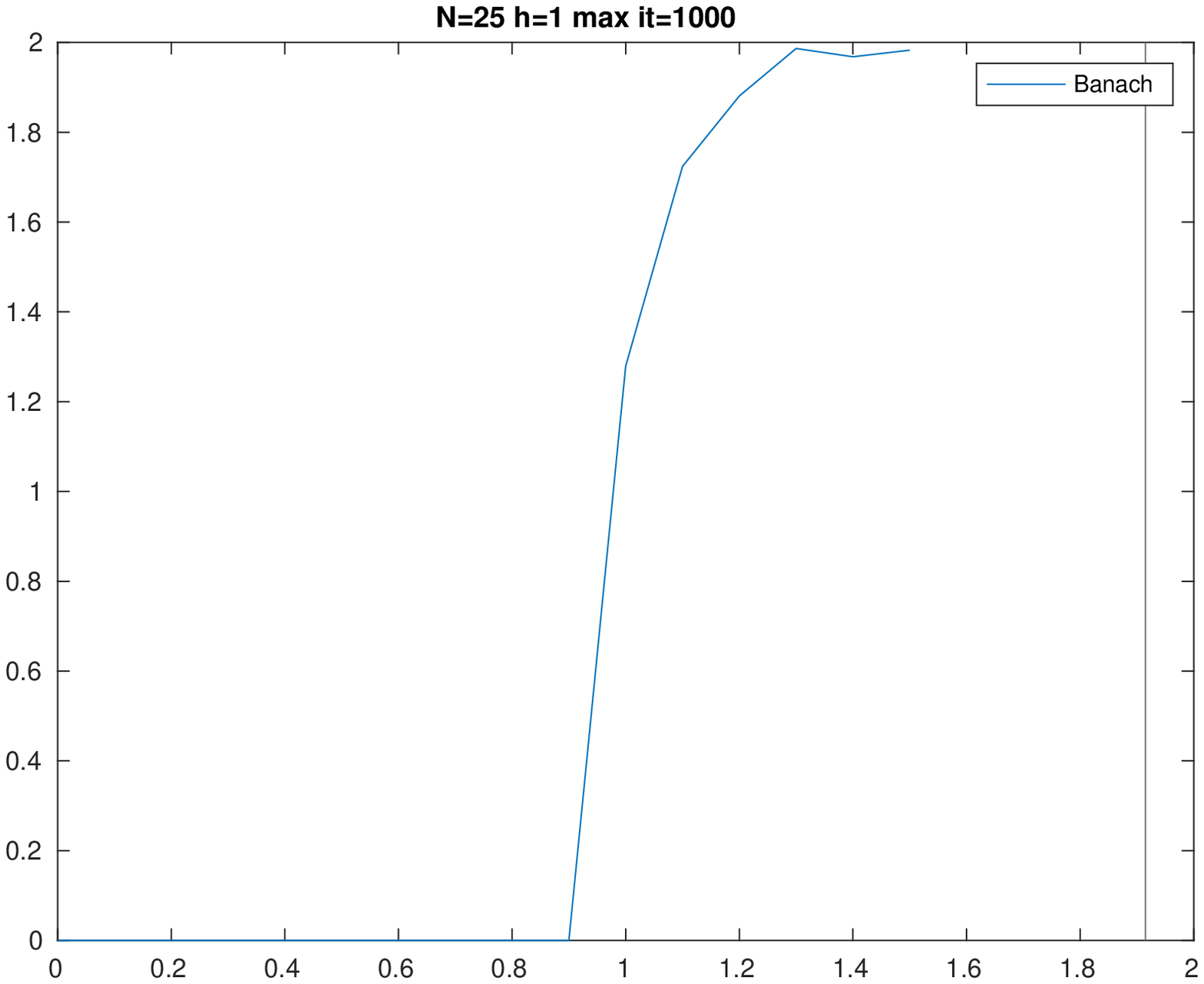}
  \caption{{\tiny} \space  1000 iterations of the standard Banach algorithm for $N=25$ and $h=1$. For $\beta\in [0, 1.5]$ we analyze the stability of this algorithm and see how quickly it takes off. To measure if the iterations stabilized we use the {\italic maximum absolute error} between the last two iterates. The vertical line marks the AT line.}
\end{figure}
\begin{figure}[H]
  \centering
  \includegraphics[width=0.5\textwidth,  trim={0 0 0 0},clip]{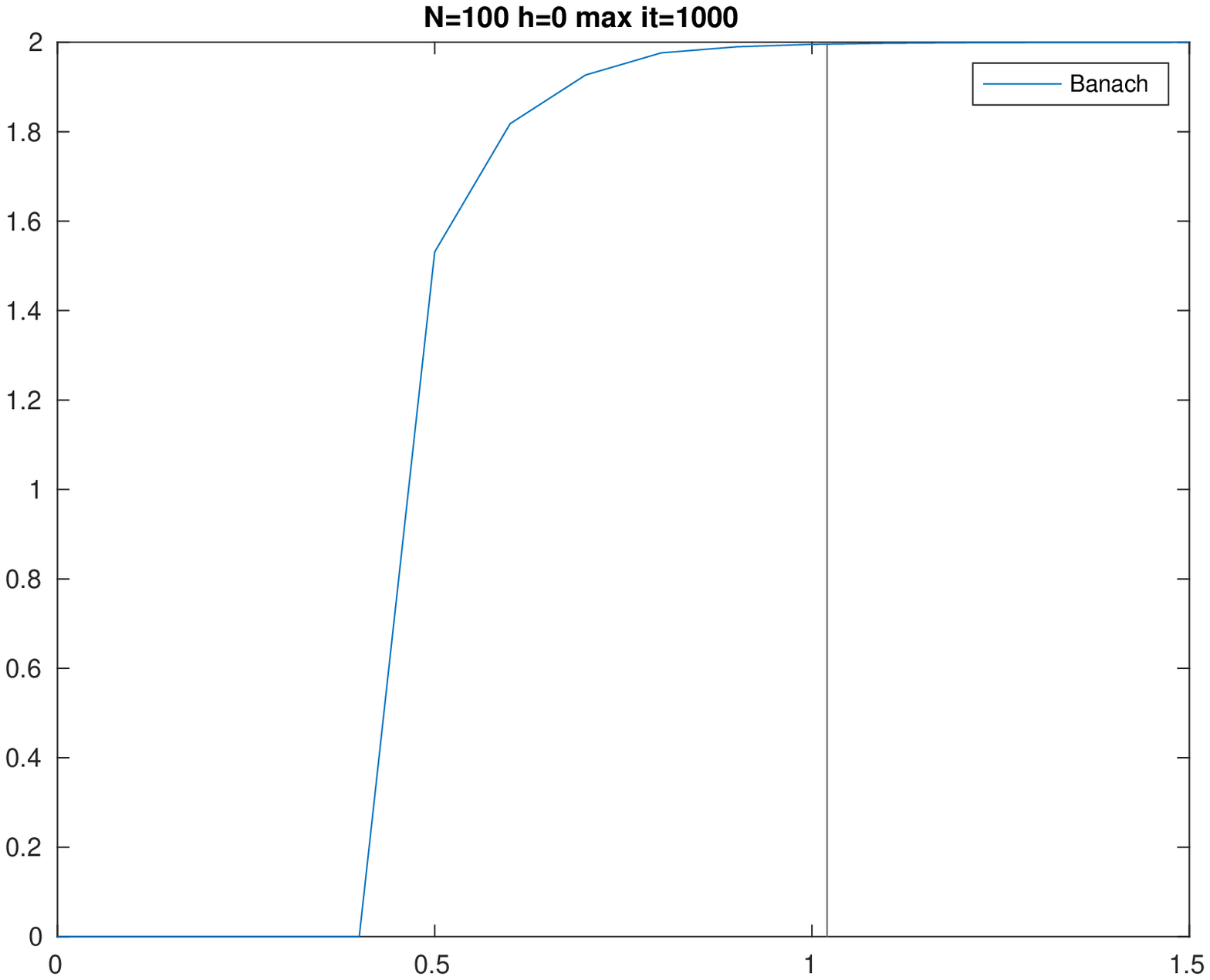}
  \caption{{\tiny} \space  1000 iterations of the standard Banach algorithm for $N=100$ and $h=0$. For $\beta\in [0, 1.5]$ we analyze the stability of this algorithm and see how quickly it takes off. To measure if the iterations stabilized we use the {\italic maximum absolute error} between the last two iterates. The vertical line marks the AT line.}
\end{figure}
\begin{figure}[H]
  \centering
  \includegraphics[width=0.5\textwidth,  trim={0 0 0 0},clip]{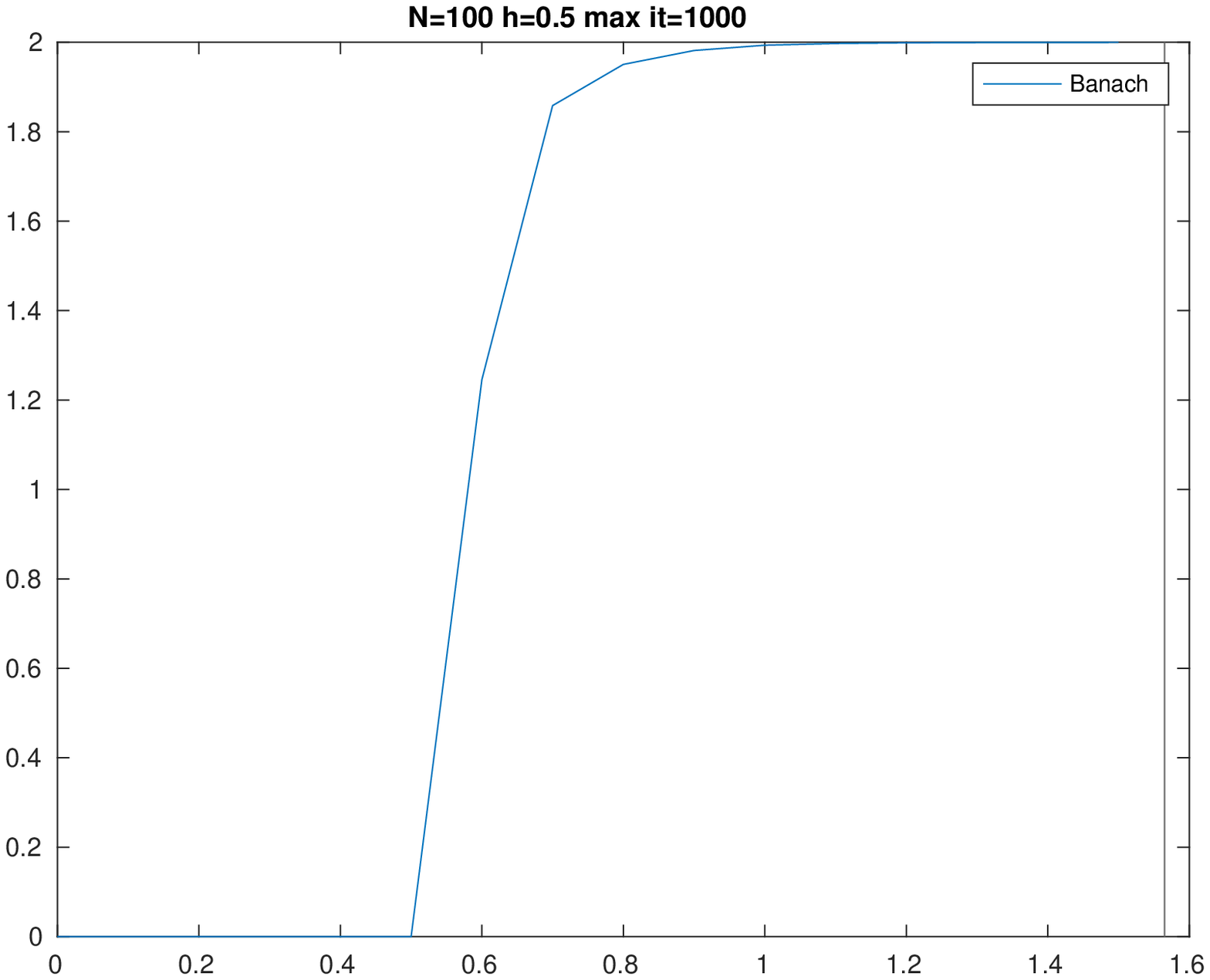}
  \caption{{\tiny} \space  1000 iterations of the standard Banach algorithm for $N=100$ and $h=0.5$. For $\beta\in [0, 1.5]$ we analyze the stability of this algorithm and see how quickly it takes off. To measure if the iterations stabilized we use the {\italic maximum absolute error} between the last two iterates. The vertical line marks the AT line.}
\end{figure}
\begin{figure}[H]
  \centering
  \includegraphics[width=0.5\textwidth,  trim={0 0 0 0},clip]{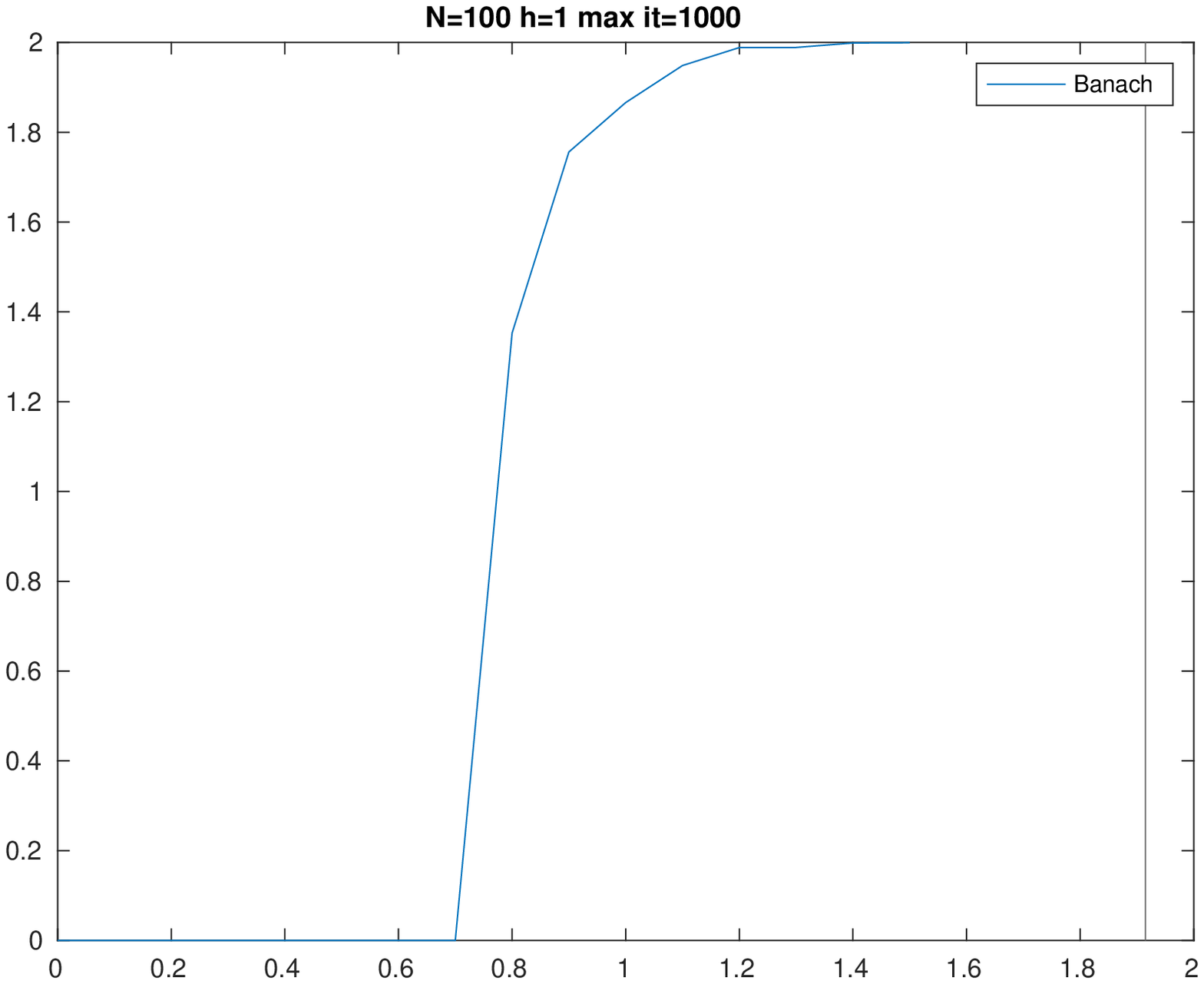}
  \caption{{\tiny} \space  1000 iterations of the standard Banach algorithm for $N=100$ and $h=1$. For $\beta\in [0, 1.5]$ we analyze the stability of this algorithm and see how quickly it takes off. To measure if the iterations stabilized we use the {\italic maximum absolute error} between the last two iterates. The vertical line marks the AT line.}
\end{figure}
\begin{figure}[H]
  \centering
  \includegraphics[width=0.5\textwidth,  trim={0 0 0 0},clip]{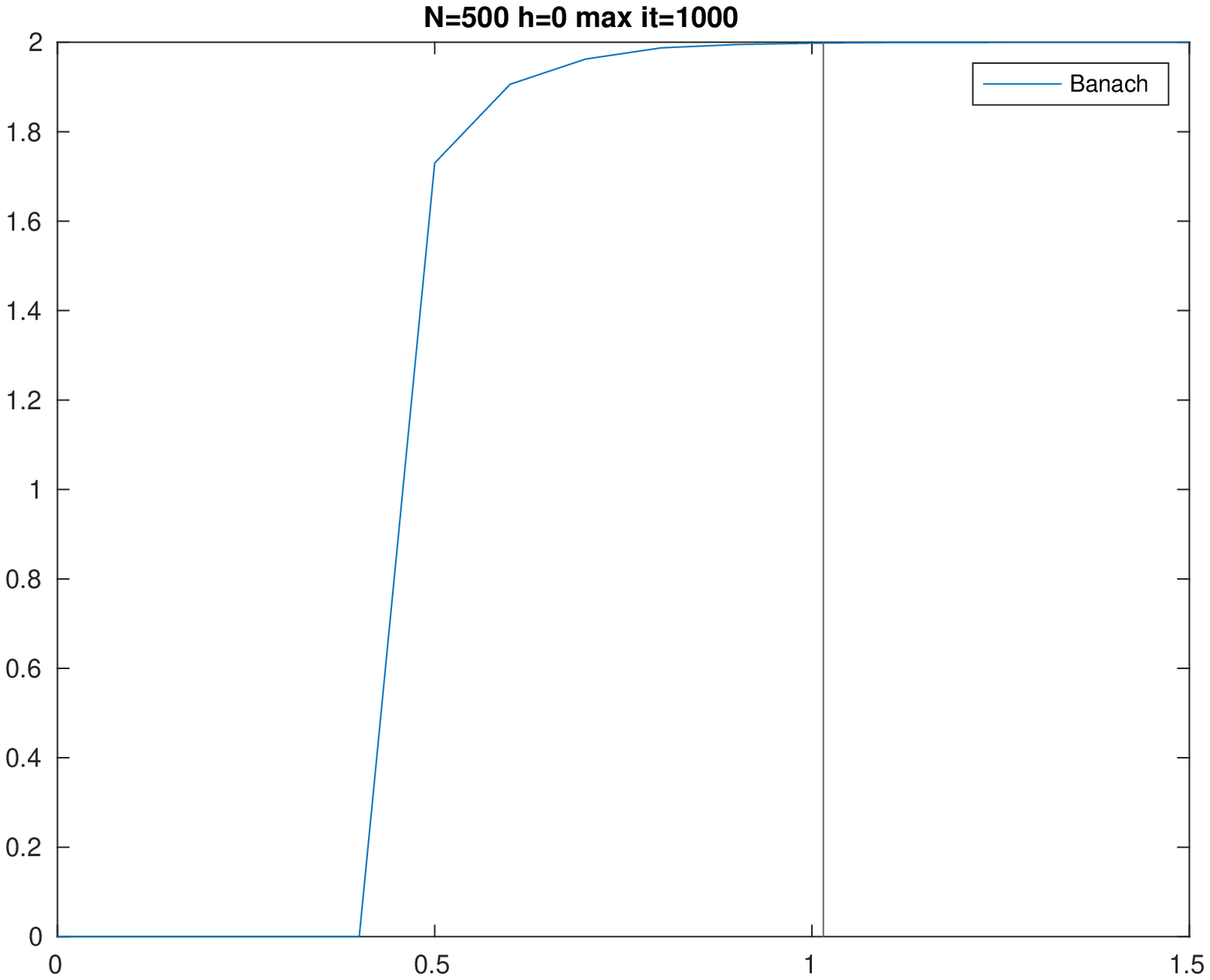}
  \caption{{\tiny} \space  1000 iterations of the standard Banach algorithm for $N=500$ and $h=0$. For $\beta\in [0, 1.5]$ we analyze the stability of this algorithm and see how quickly it takes off. To measure if the iterations stabilized we use the {\italic maximum absolute error} between the last two iterates. The vertical line marks the AT line.}
\end{figure}
\begin{figure}[H]
  \centering
  \includegraphics[width=0.5\textwidth,  trim={0 0 0 0},clip]{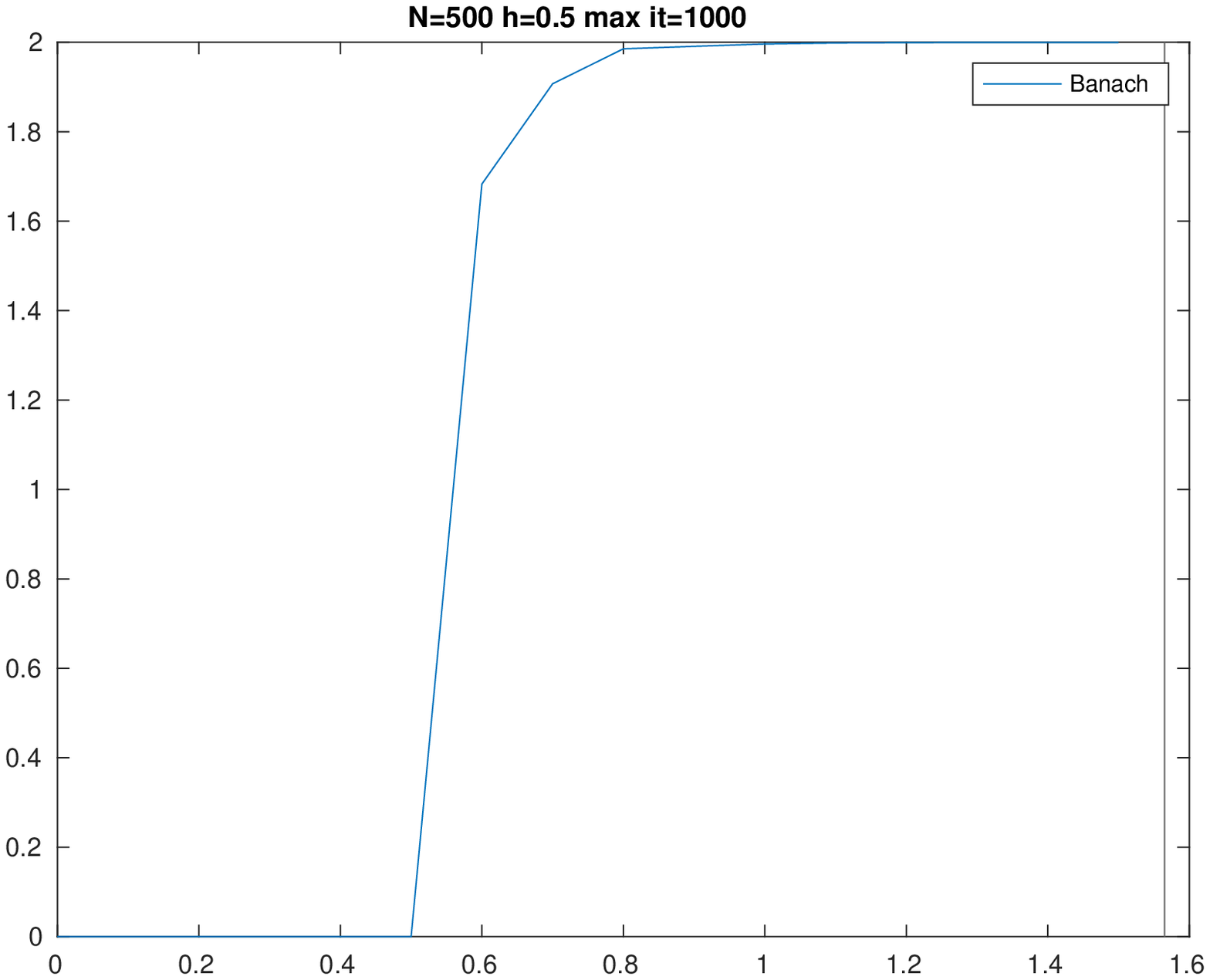}
  \caption{{\tiny} \space  1000 iterations of the standard Banach algorithm for $N=500$ and $h=0.5$. For $\beta\in [0, 1.5]$ we analyze the stability of this algorithm and see how quickly it takes off. To measure if the iterations stabilized we use the {\italic maximum absolute error} between the last two iterates. The vertical line marks the AT line.}
\end{figure}
\begin{figure}[H]
  \centering
  \includegraphics[width=0.5\textwidth,  trim={0 0 0 0},clip]{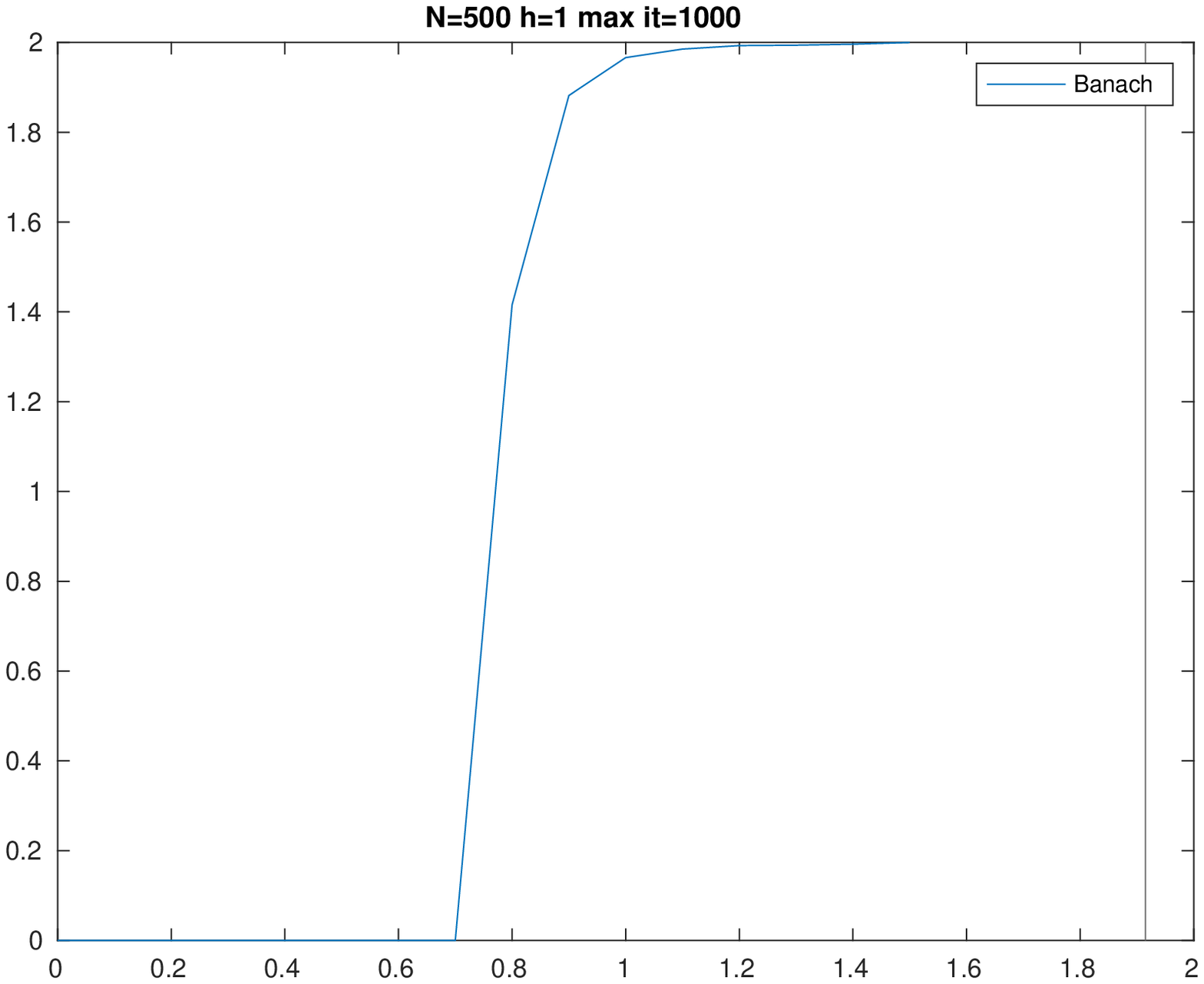}
  \caption{{\tiny} \space  1000 iterations of the standard Banach algorithm for $N=500$ and $h=1$. For $\beta\in [0, 1.5]$ we analyze the stability of this algorithm and see how quickly it takes off. To measure if the iterations stabilized we use the {\italic maximum absolute error} between the last two iterates. The vertical line marks the AT line.}
\end{figure}
\begin{figure}[H]
  \centering
  \includegraphics[width=0.5\textwidth,  trim={0 0 0 0},clip]{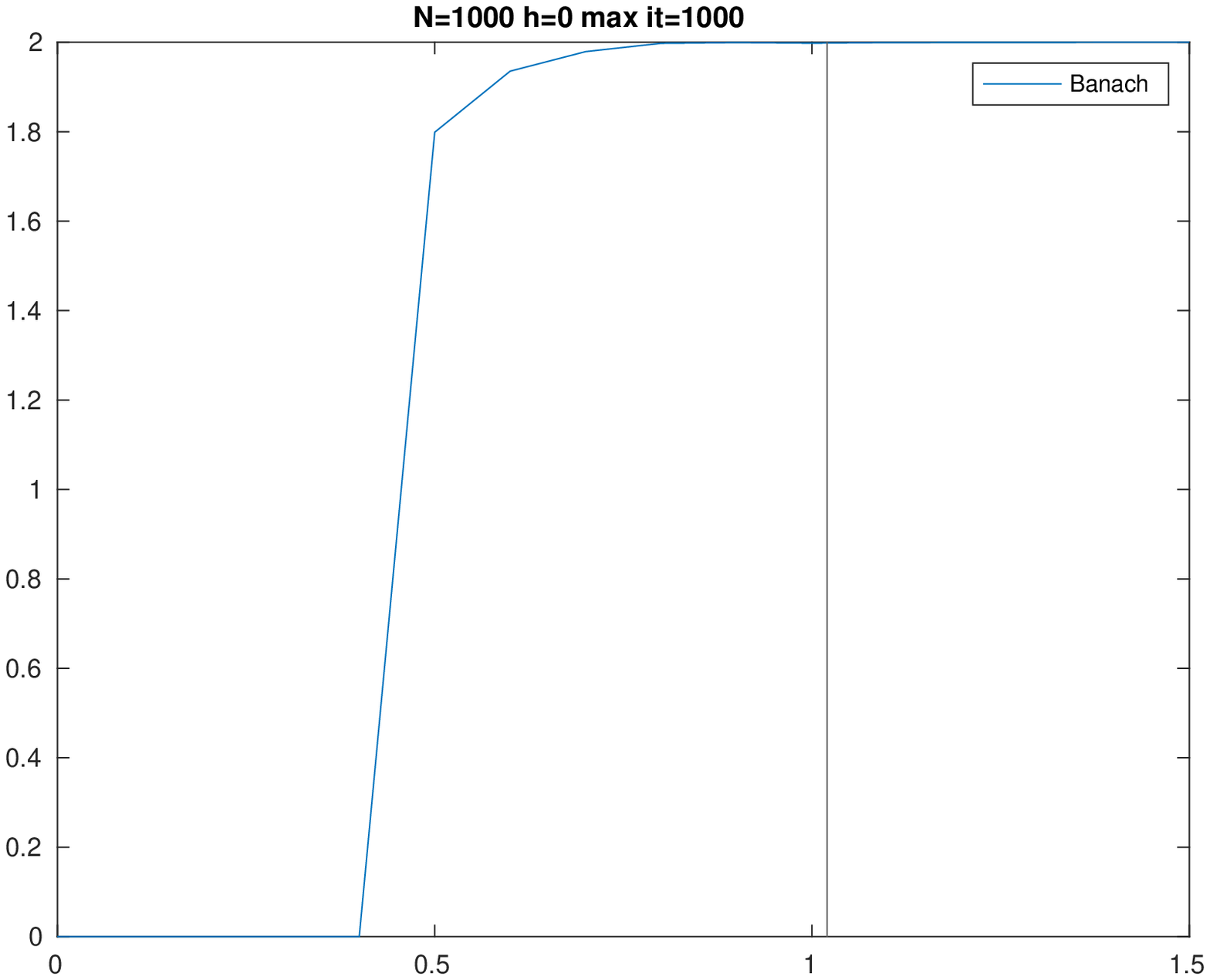}
  \caption{{\tiny} \space  1000 iterations of the standard Banach algorithm for $N=1000$ and $h=0$. For $\beta\in [0, 1.5]$ we analyze the stability of this algorithm and see how quickly it takes off. To measure if the iterations stabilized we use the {\italic maximum absolute error} between the last two iterates. The vertical line marks the AT line.}
\end{figure}
\begin{figure}[H]
  \centering
  \includegraphics[width=0.5\textwidth,  trim={0 0 0 0},clip]{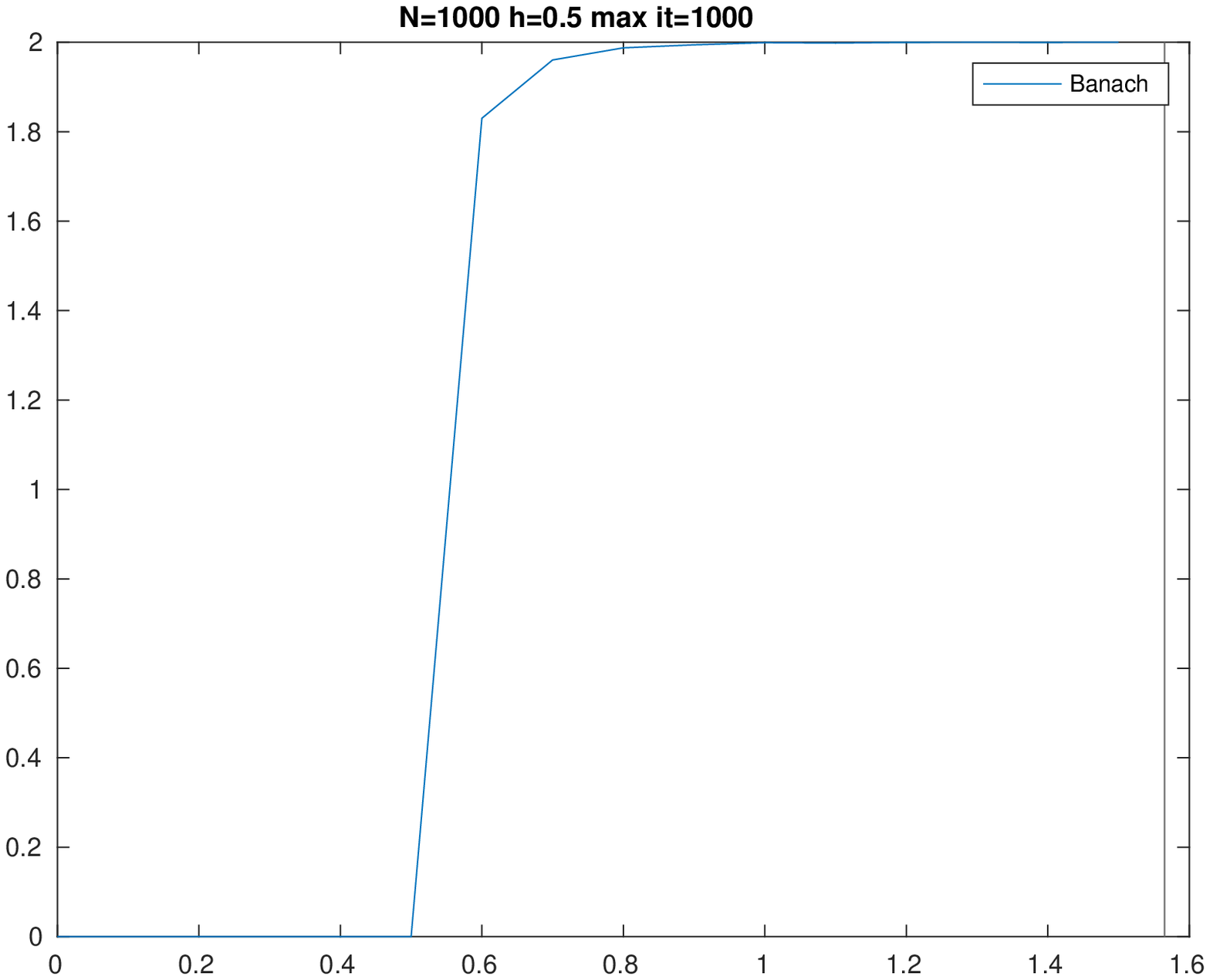}
  \caption{{\tiny} \space  1000 iterations of the standard Banach algorithm for $N=1000$ and $h=0.5$. For $\beta\in [0, 1.5]$ we analyze the stability of this algorithm and see how quickly it takes off. To measure if the iterations stabilized we use the {\italic maximum absolute error} between the last two iterates. The vertical line marks the AT line.}
\end{figure}
\begin{figure}[H]
  \centering
  \includegraphics[width=0.5\textwidth,  trim={0 0 0 0},clip]{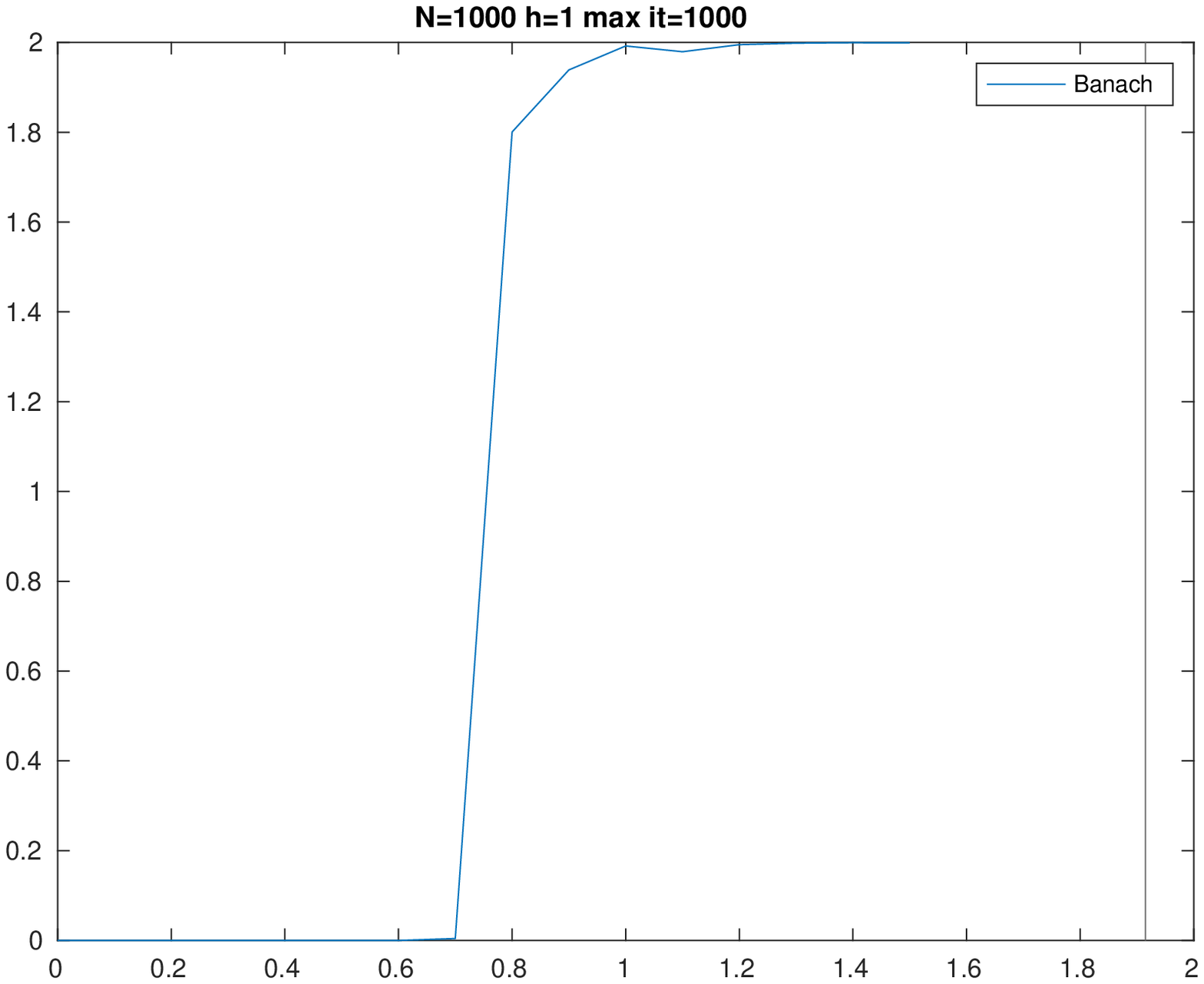}
  \caption{{\tiny} \space  1000 iterations of the standard Banach algorithm for $N=1000$ and $h=1$. For $\beta\in [0, 1.5]$ we analyze the stability of this algorithm and see how quickly it takes off. To measure if the iterations stabilized we use the {\italic maximum absolute error} between the last two iterates. The vertical line marks the AT line.}
\end{figure}\FloatBarrier\clearpage
\subsection{Using the {\it mean squared error}}
\subsubsection{Uniformly in the hypercube distributed start values}
\leavevmode
\begin{figure}[H]
  \centering
  \includegraphics[width=0.5\textwidth,  trim={0.45cm 0.4cm 0 0},clip]{rsynccopy/final_simulations/starting_uniform_in_the_hypercube/maxdiverginess_banach_explodes_at_line/1Banach.eps}
  \caption{{\tiny} \space  1000 iterations of the standard Banach algorithm for $N=10$ and $h=0$. For $\beta\in [0, 1.5]$ we analyze the stability of this algorithm and see how quickly it takes off. To measure if the iterations stabilized we use the {\italic mean squared error} between the last two iterates. The vertical line marks the AT line.}
\end{figure}
\begin{figure}[H]
  \centering
  \includegraphics[width=0.5\textwidth,  trim={0.45cm 0.4cm 0 0},clip]{rsynccopy/final_simulations/starting_uniform_in_the_hypercube/maxdiverginess_banach_explodes_at_line/2Banach.eps}
  \caption{{\tiny} \space  1000 iterations of the standard Banach algorithm for $N=10$ and $h=0.5$. For $\beta\in [0, 1.5]$ we analyze the stability of this algorithm and see how quickly it takes off. To measure if the iterations stabilized we use the {\italic mean squared error} between the last two iterates. The vertical line marks the AT line.}
\end{figure}
\begin{figure}[H]
  \centering
  \includegraphics[width=0.5\textwidth,  trim={0 0 0 0},clip]{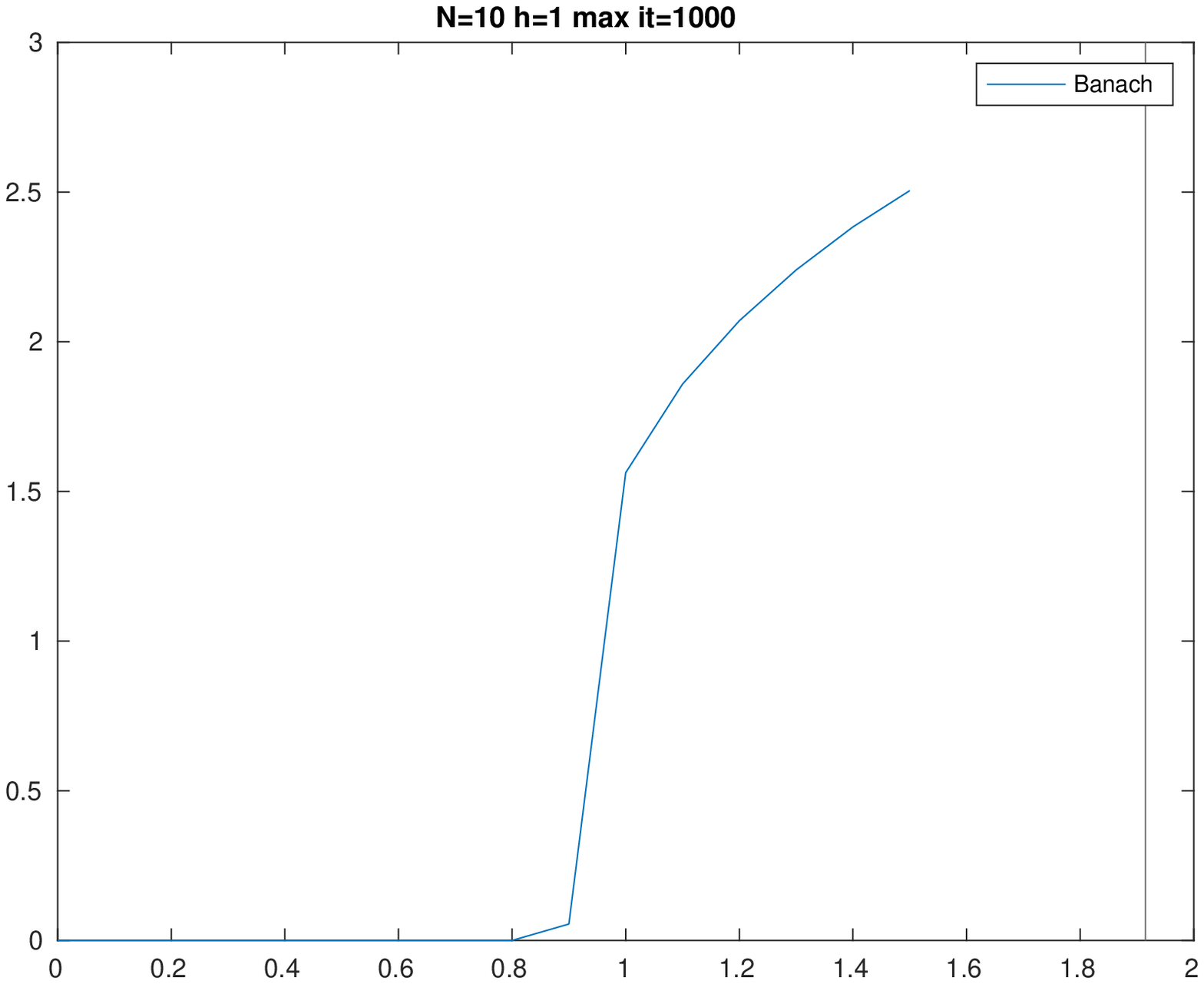}
  \caption{{\tiny} \space  1000 iterations of the standard Banach algorithm for $N=10$ and $h=1$. For $\beta\in [0, 1.5]$ we analyze the stability of this algorithm and see how quickly it takes off. To measure if the iterations stabilized we use the {\italic mean squared error} between the last two iterates. The vertical line marks the AT line.}
\end{figure}
\begin{figure}[H]
  \centering
  \includegraphics[width=0.5\textwidth,  trim={0 0 0 0},clip]{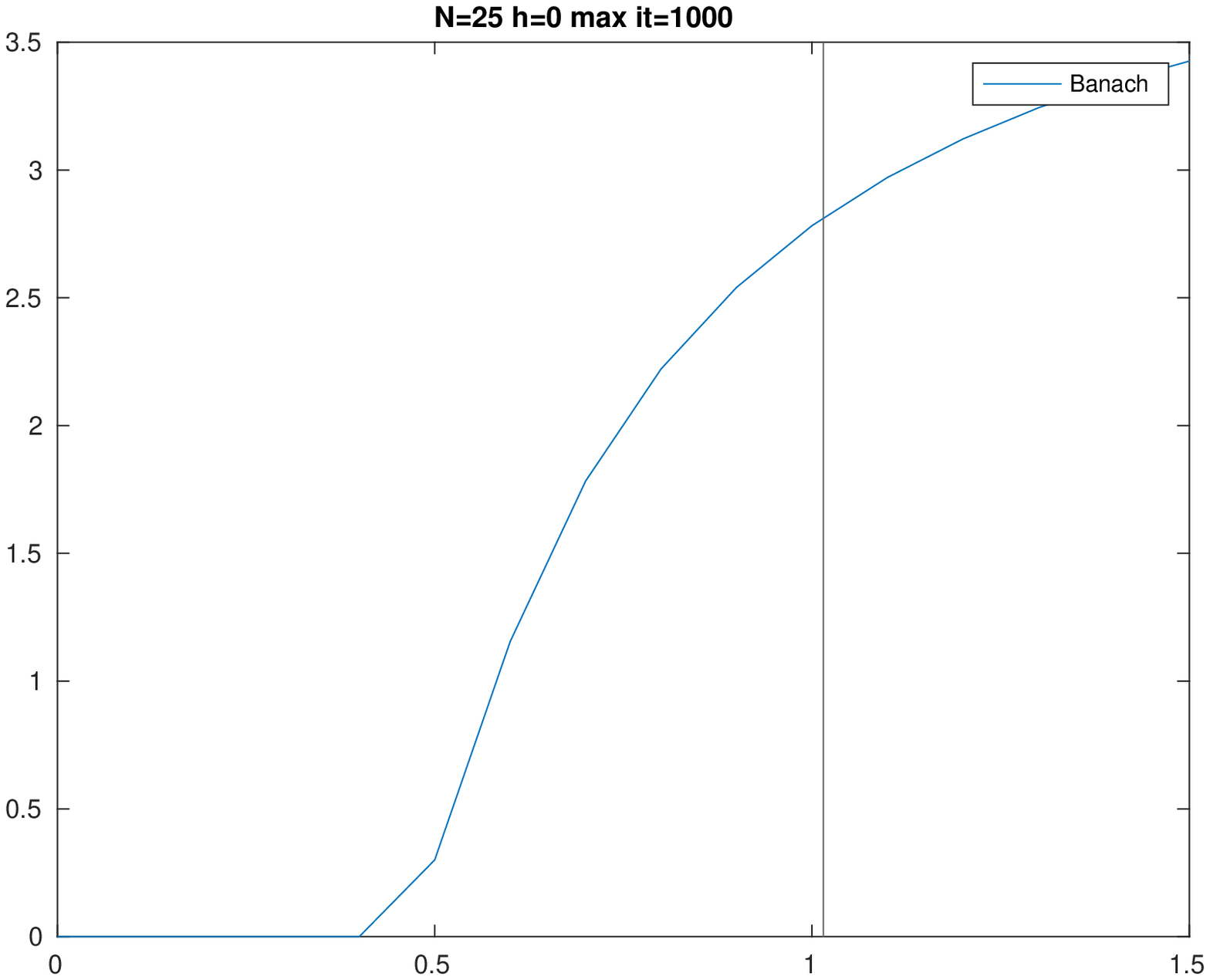}
  \caption{{\tiny} \space  1000 iterations of the standard Banach algorithm for $N=25$ and $h=0$. For $\beta\in [0, 1.5]$ we analyze the stability of this algorithm and see how quickly it takes off. To measure if the iterations stabilized we use the {\italic mean squared error} between the last two iterates. The vertical line marks the AT line.}
\end{figure}
\begin{figure}[H]
  \centering
  \includegraphics[width=0.5\textwidth,  trim={0 0 0 0},clip]{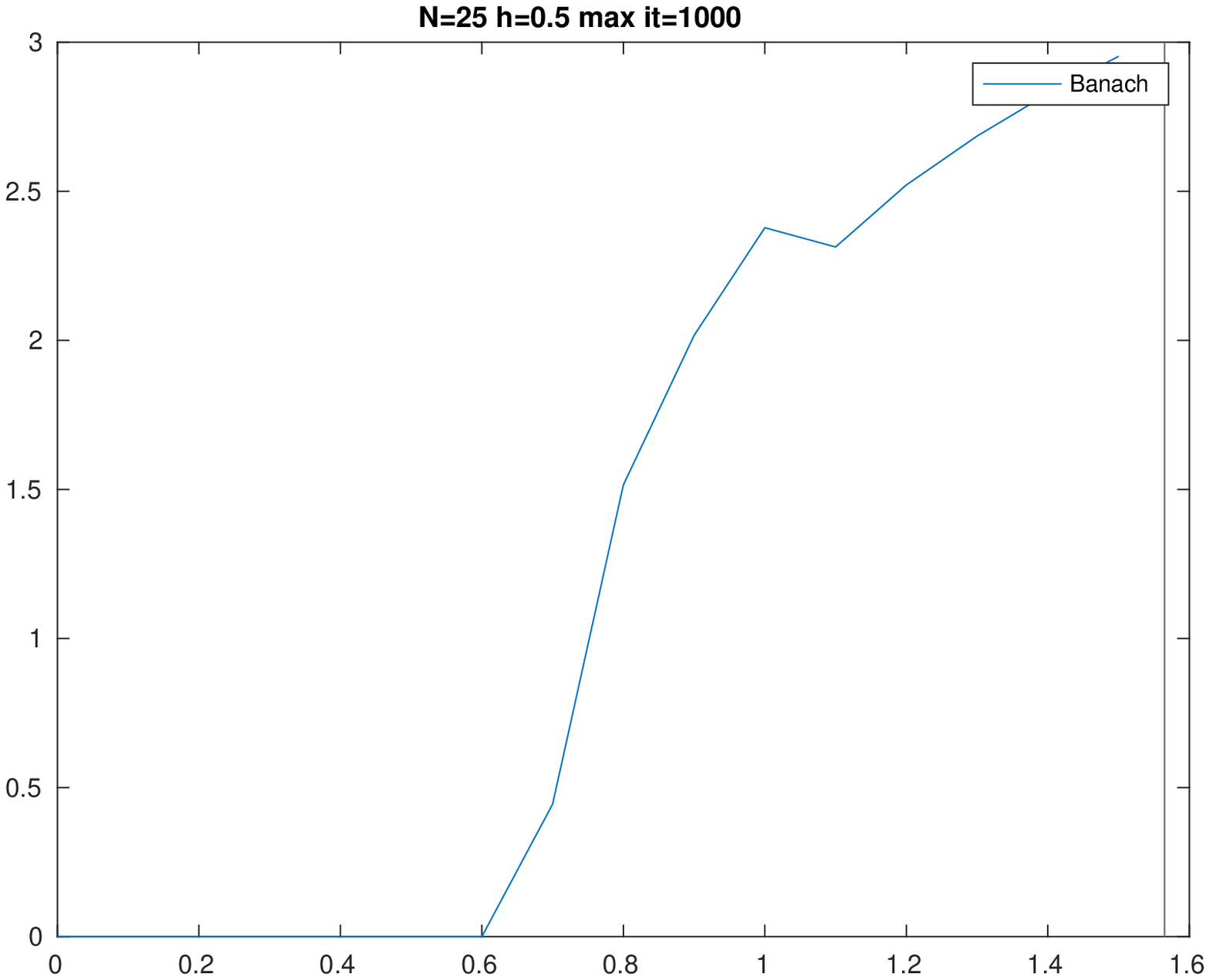}
  \caption{{\tiny} \space  1000 iterations of the standard Banach algorithm for $N=25$ and $h=0.5$. For $\beta\in [0, 1.5]$ we analyze the stability of this algorithm and see how quickly it takes off. To measure if the iterations stabilized we use the {\italic mean squared error} between the last two iterates. The vertical line marks the AT line.}
\end{figure}
\begin{figure}[H]
  \centering
  \includegraphics[width=0.5\textwidth,  trim={0 0 0 0},clip]{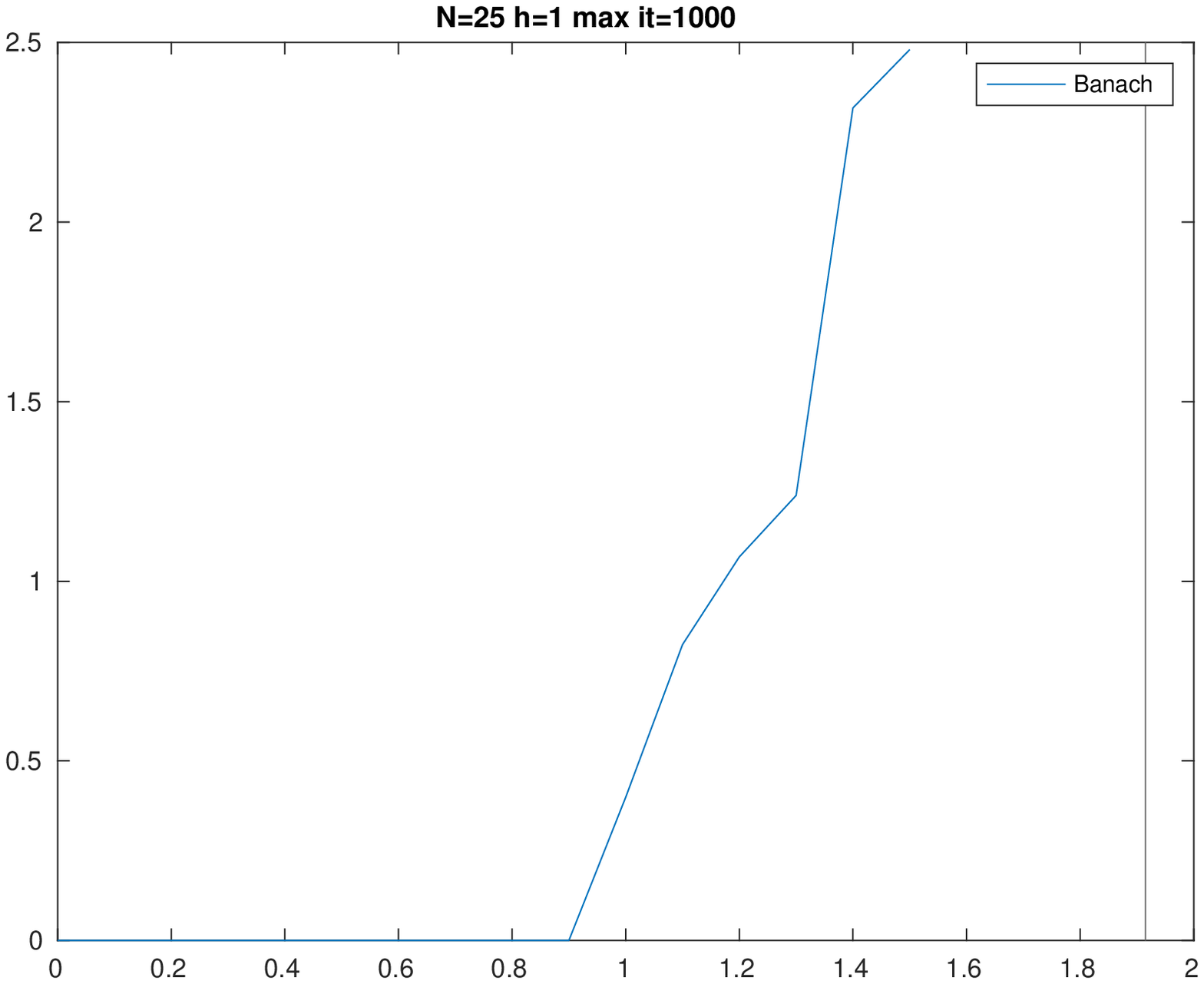}
  \caption{{\tiny} \space  1000 iterations of the standard Banach algorithm for $N=25$ and $h=1$. For $\beta\in [0, 1.5]$ we analyze the stability of this algorithm and see how quickly it takes off. To measure if the iterations stabilized we use the {\italic mean squared error} between the last two iterates. The vertical line marks the AT line.}
\end{figure}
\begin{figure}[H]
  \centering
  \includegraphics[width=0.5\textwidth,  trim={0 0 0 0},clip]{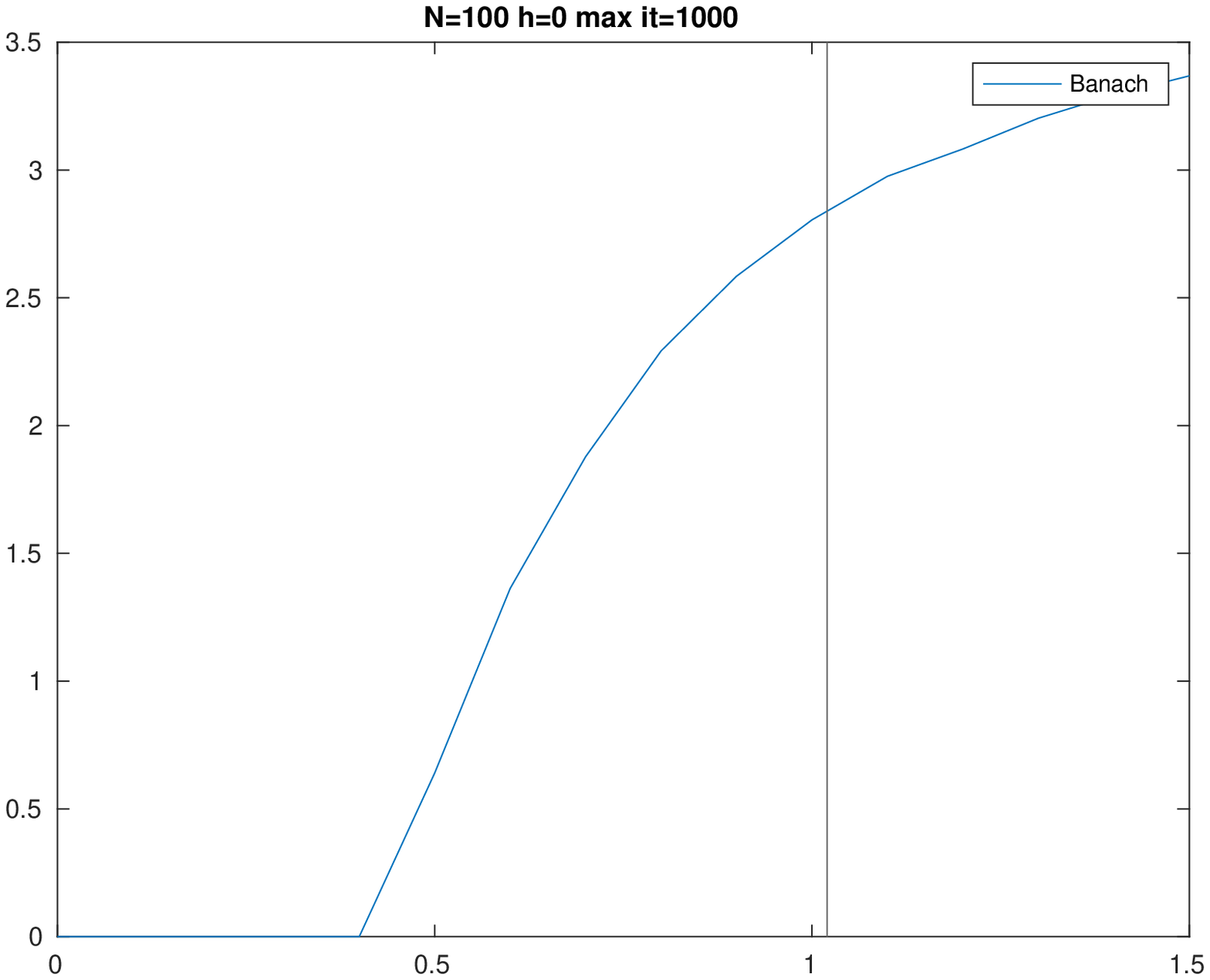}
  \caption{{\tiny} \space  1000 iterations of the standard Banach algorithm for $N=100$ and $h=0$. For $\beta\in [0, 1.5]$ we analyze the stability of this algorithm and see how quickly it takes off. To measure if the iterations stabilized we use the {\italic mean squared error} between the last two iterates. The vertical line marks the AT line.}
\end{figure}
\begin{figure}[H]
  \centering
  \includegraphics[width=0.5\textwidth,  trim={0 0 0 0},clip]{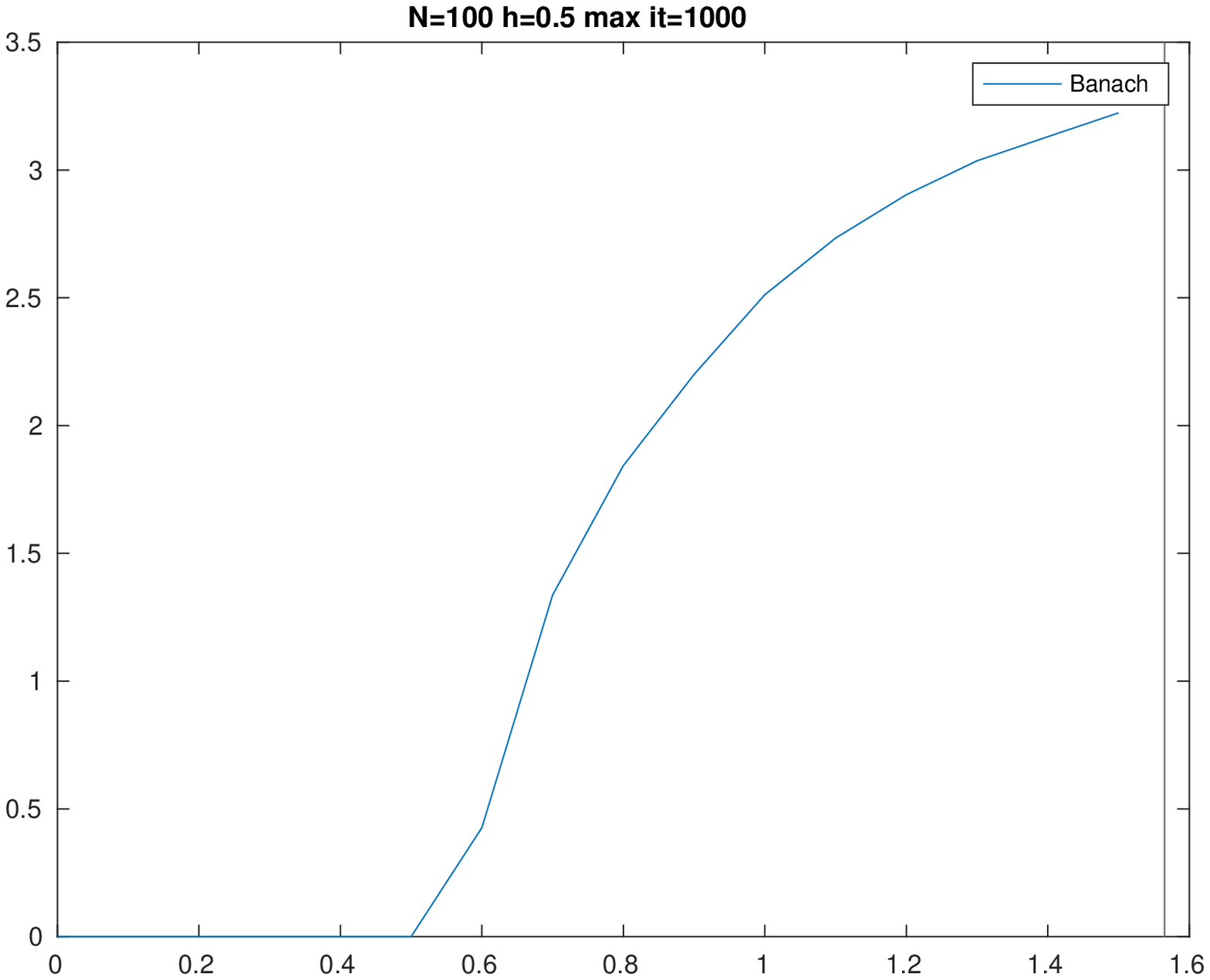}
  \caption{{\tiny} \space  1000 iterations of the standard Banach algorithm for $N=100$ and $h=0.5$. For $\beta\in [0, 1.5]$ we analyze the stability of this algorithm and see how quickly it takes off. To measure if the iterations stabilized we use the {\italic mean squared error} between the last two iterates. The vertical line marks the AT line.}
\end{figure}
\begin{figure}[H]
  \centering
  \includegraphics[width=0.5\textwidth,  trim={0 0 0 0},clip]{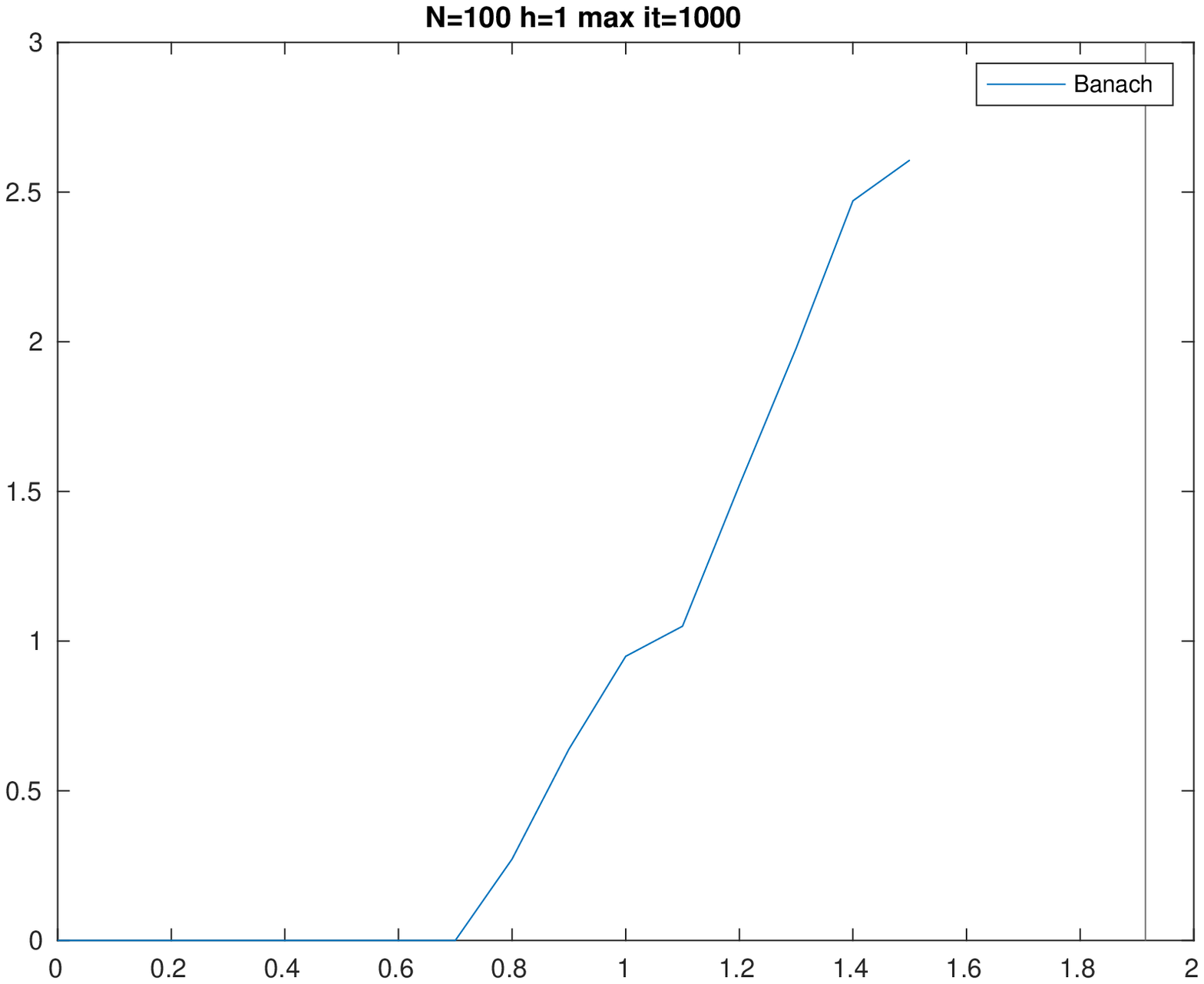}
  \caption{{\tiny} \space  1000 iterations of the standard Banach algorithm for $N=100$ and $h=1$. For $\beta\in [0, 1.5]$ we analyze the stability of this algorithm and see how quickly it takes off. To measure if the iterations stabilized we use the {\italic mean squared error} between the last two iterates. The vertical line marks the AT line.}
\end{figure}
\begin{figure}[H]
  \centering
  \includegraphics[width=0.5\textwidth,  trim={0 0 0 0},clip]{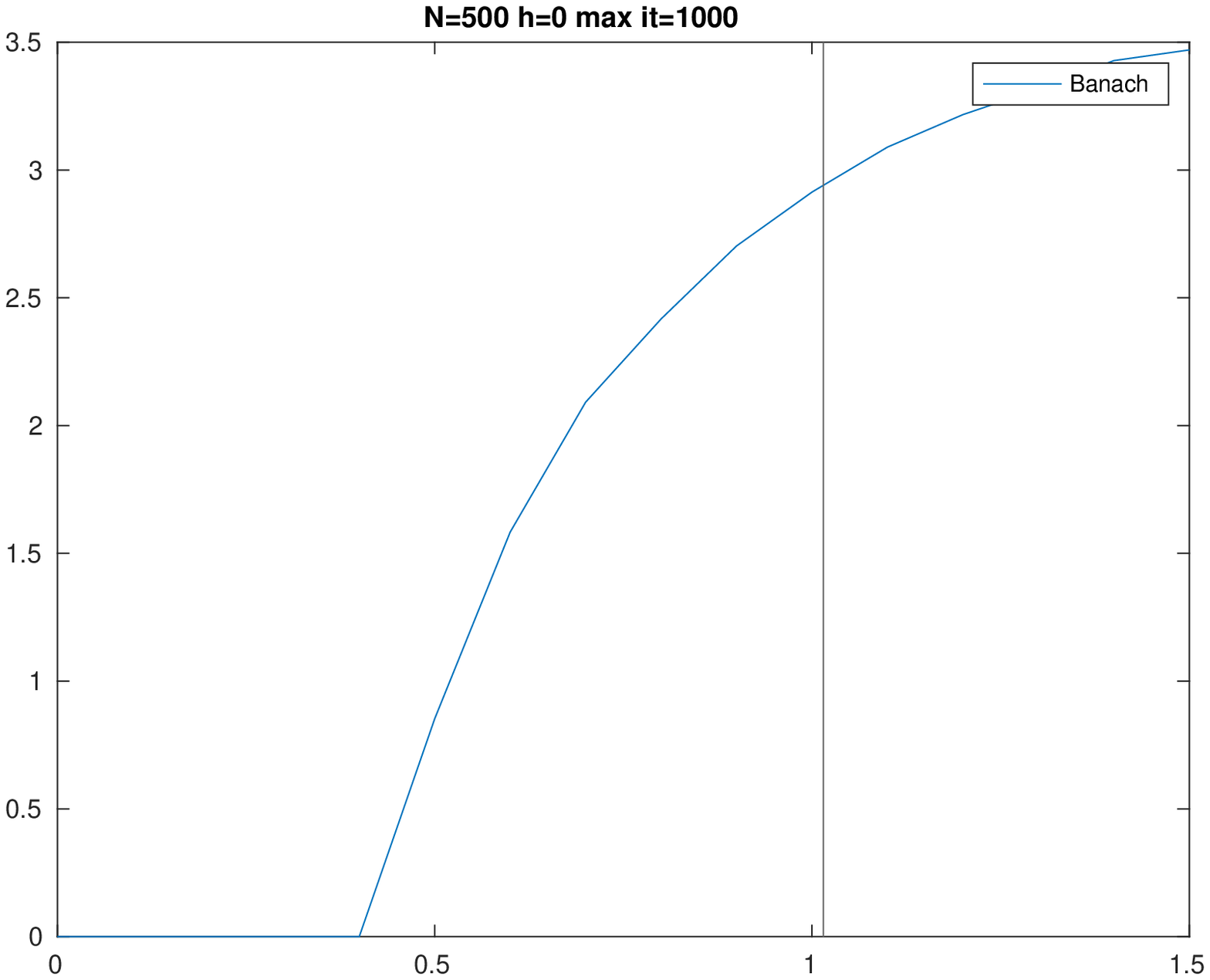}
  \caption{{\tiny} \space  1000 iterations of the standard Banach algorithm for $N=500$ and $h=0$. For $\beta\in [0, 1.5]$ we analyze the stability of this algorithm and see how quickly it takes off. To measure if the iterations stabilized we use the {\italic mean squared error} between the last two iterates. The vertical line marks the AT line.}
\end{figure}
\begin{figure}[H]
  \centering
  \includegraphics[width=0.5\textwidth,  trim={0 0 0 0},clip]{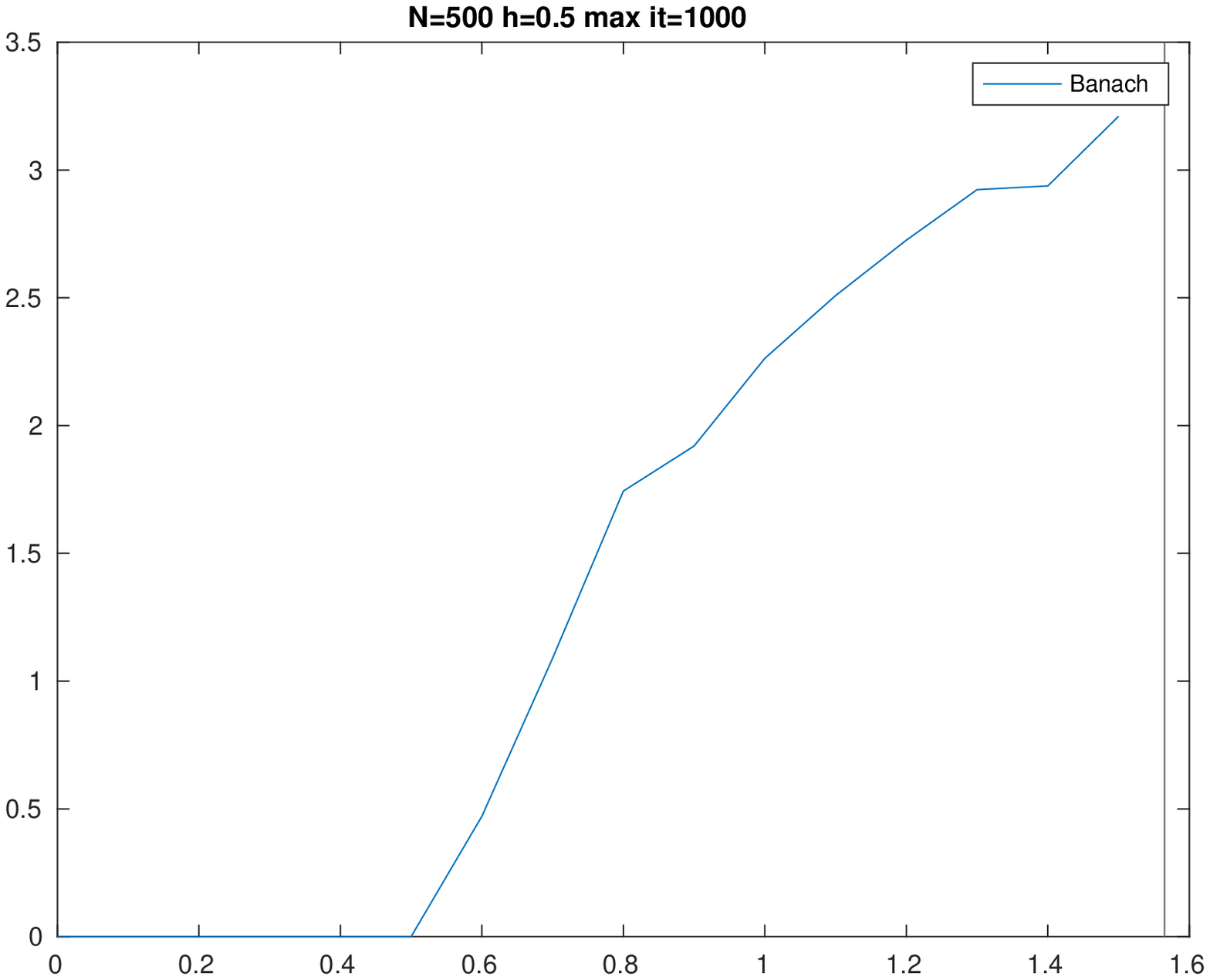}
  \caption{{\tiny} \space  1000 iterations of the standard Banach algorithm for $N=500$ and $h=0.5$. For $\beta\in [0, 1.5]$ we analyze the stability of this algorithm and see how quickly it takes off. To measure if the iterations stabilized we use the {\italic mean squared error} between the last two iterates. The vertical line marks the AT line.}
\end{figure}
\begin{figure}[H]
  \centering
  \includegraphics[width=0.5\textwidth,  trim={0 0 0 0},clip]{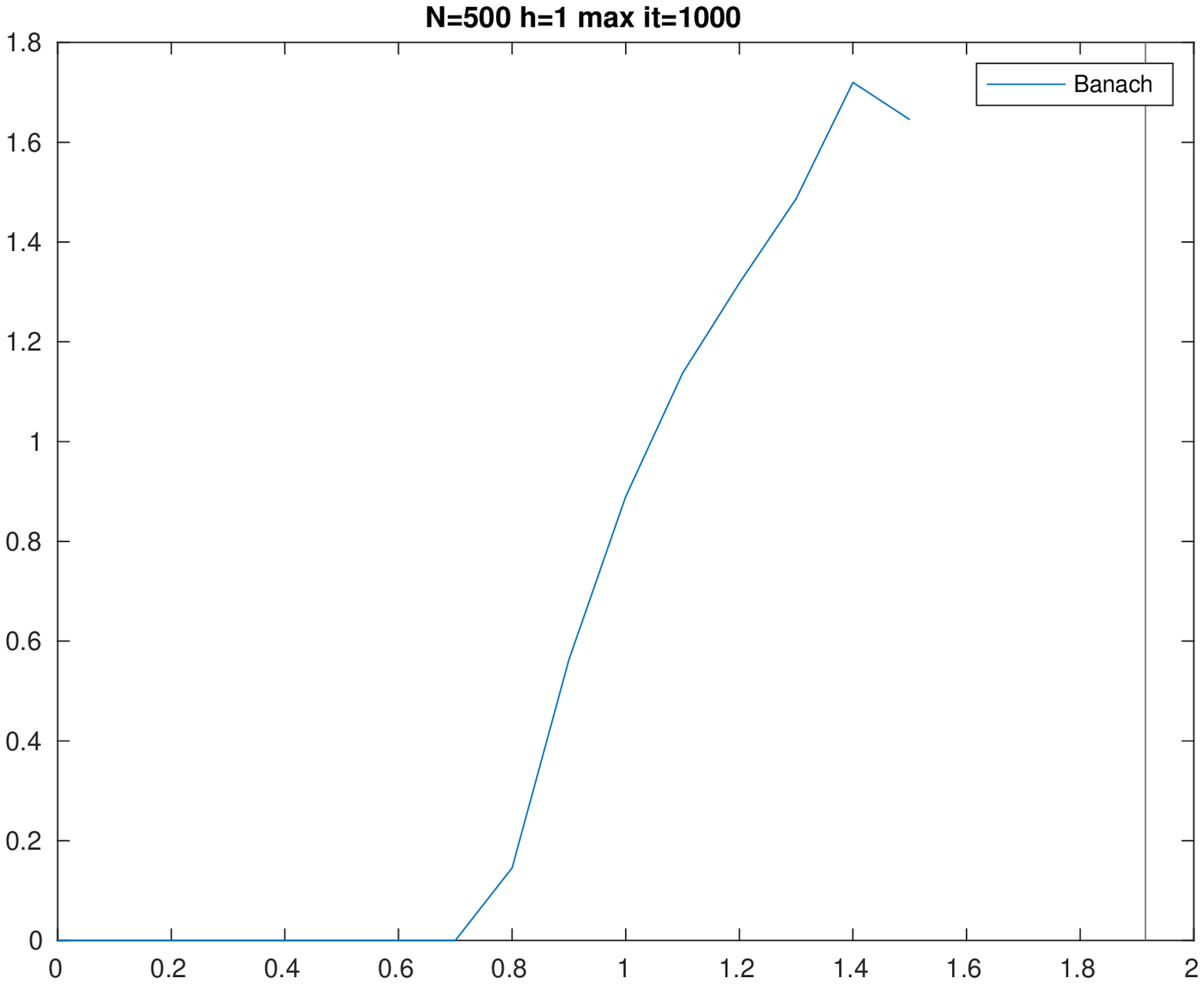}
  \caption{{\tiny} \space  1000 iterations of the standard Banach algorithm for $N=500$ and $h=1$. For $\beta\in [0, 1.5]$ we analyze the stability of this algorithm and see how quickly it takes off. To measure if the iterations stabilized we use the {\italic mean squared error} between the last two iterates. The vertical line marks the AT line.}
\end{figure}
\begin{figure}[H]
  \centering
  \includegraphics[width=0.5\textwidth,  trim={0 0 0 0},clip]{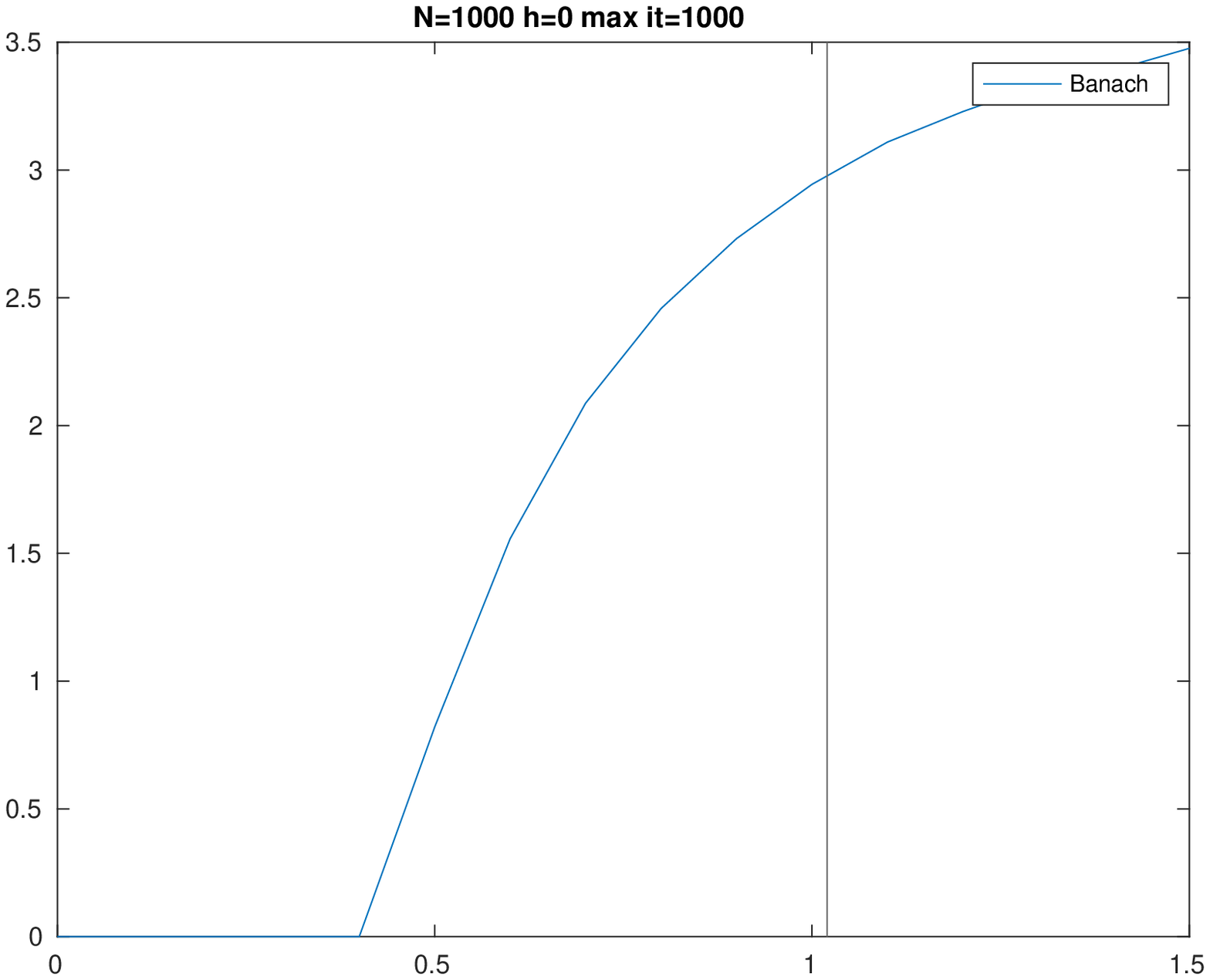}
  \caption{{\tiny} \space  1000 iterations of the standard Banach algorithm for $N=1000$ and $h=0$. For $\beta\in [0, 1.5]$ we analyze the stability of this algorithm and see how quickly it takes off. To measure if the iterations stabilized we use the {\italic mean squared error} between the last two iterates. The vertical line marks the AT line.}
\end{figure}
\begin{figure}[H]
  \centering
  \includegraphics[width=0.5\textwidth,  trim={0 0 0 0},clip]{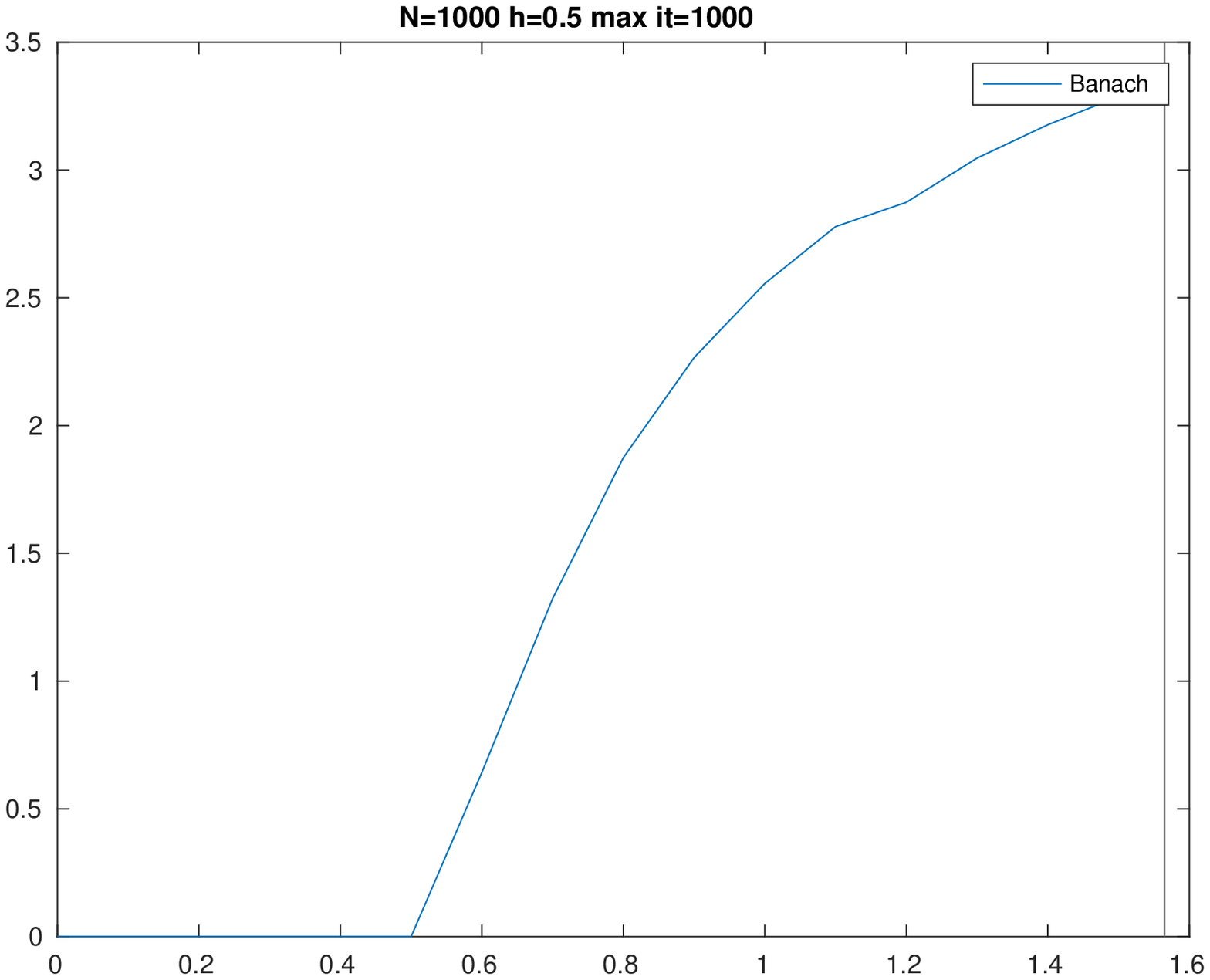}
  \caption{{\tiny} \space  1000 iterations of the standard Banach algorithm for $N=1000$ and $h=0.5$. For $\beta\in [0, 1.5]$ we analyze the stability of this algorithm and see how quickly it takes off. To measure if the iterations stabilized we use the {\italic mean squared error} between the last two iterates. The vertical line marks the AT line.}
\end{figure}
\begin{figure}[H]
  \centering
  \includegraphics[width=0.5\textwidth,  trim={0 0 0 0},clip]{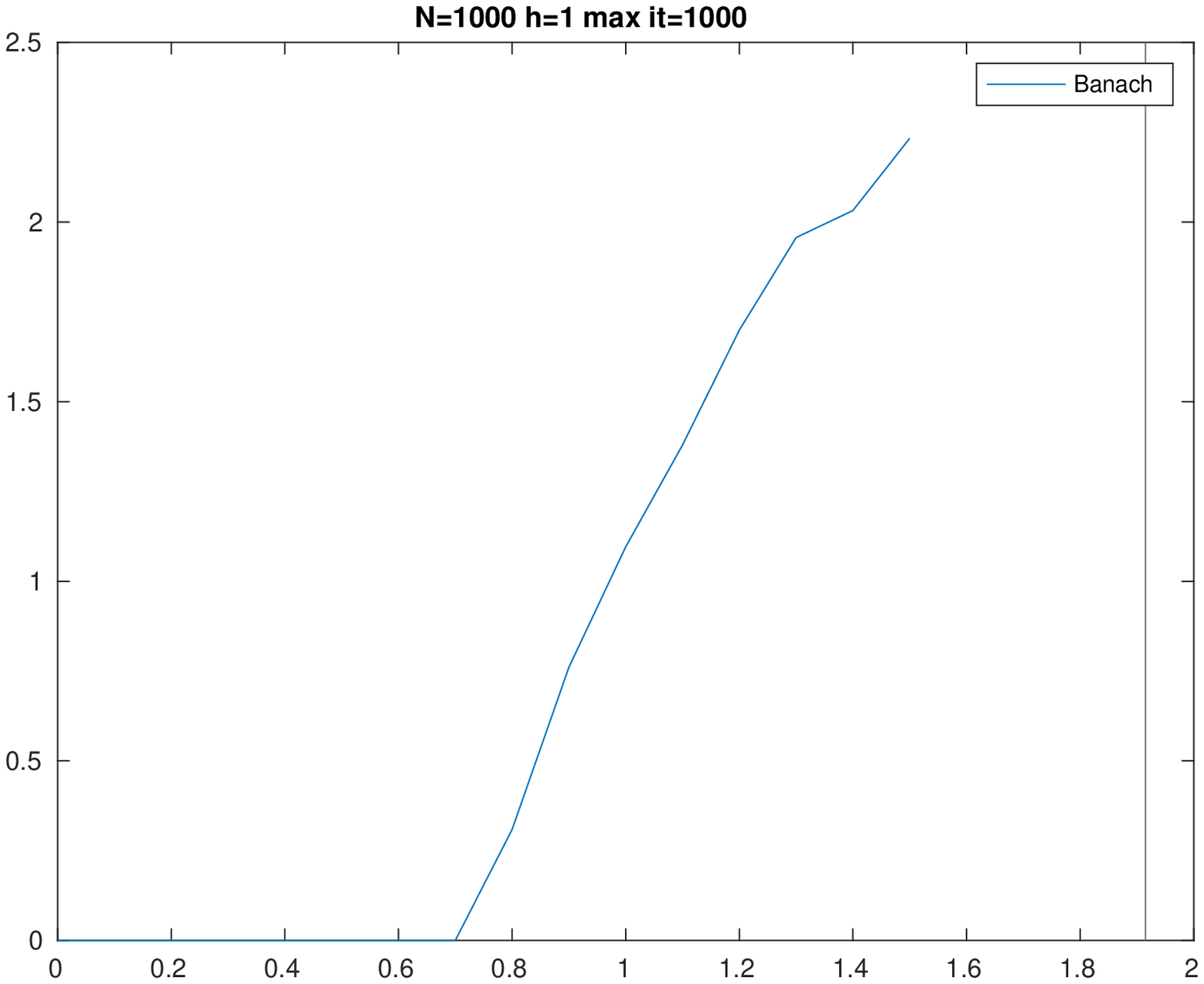}
  \caption{{\tiny} \space  1000 iterations of the standard Banach algorithm for $N=1000$ and $h=1$. For $\beta\in [0, 1.5]$ we analyze the stability of this algorithm and see how quickly it takes off. To measure if the iterations stabilized we use the {\italic mean squared error} between the last two iterates. The vertical line marks the AT line.}
\end{figure}\FloatBarrier\clearpage
\subsubsection{Uniformly in the corners of the hypercube distributed start values}
\leavevmode
\begin{figure}[H]
  \centering
  \includegraphics[width=0.5\textwidth,  trim={0 0 0 0},clip]{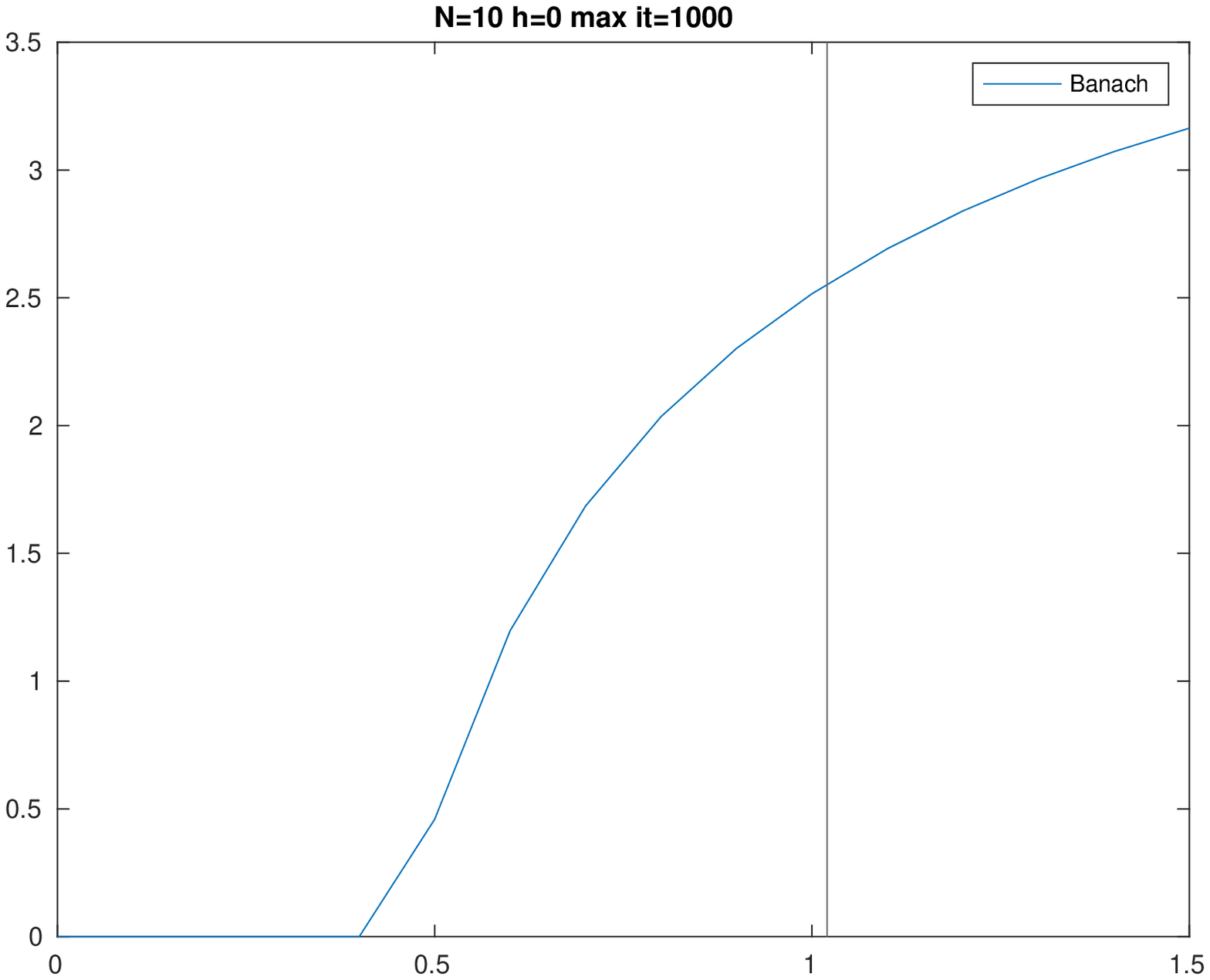}
  \caption{{\tiny} \space  1000 iterations of the standard Banach algorithm for $N=10$ and $h=0$. For $\beta\in [0, 1.5]$ we analyze the stability of this algorithm and see how quickly it takes off. To measure if the iterations stabilized we use the {\italic mean squared error} between the last two iterates. The vertical line marks the AT line.}
\end{figure}
\begin{figure}[H]
  \centering
  \includegraphics[width=0.5\textwidth,  trim={0 0 0 0},clip]{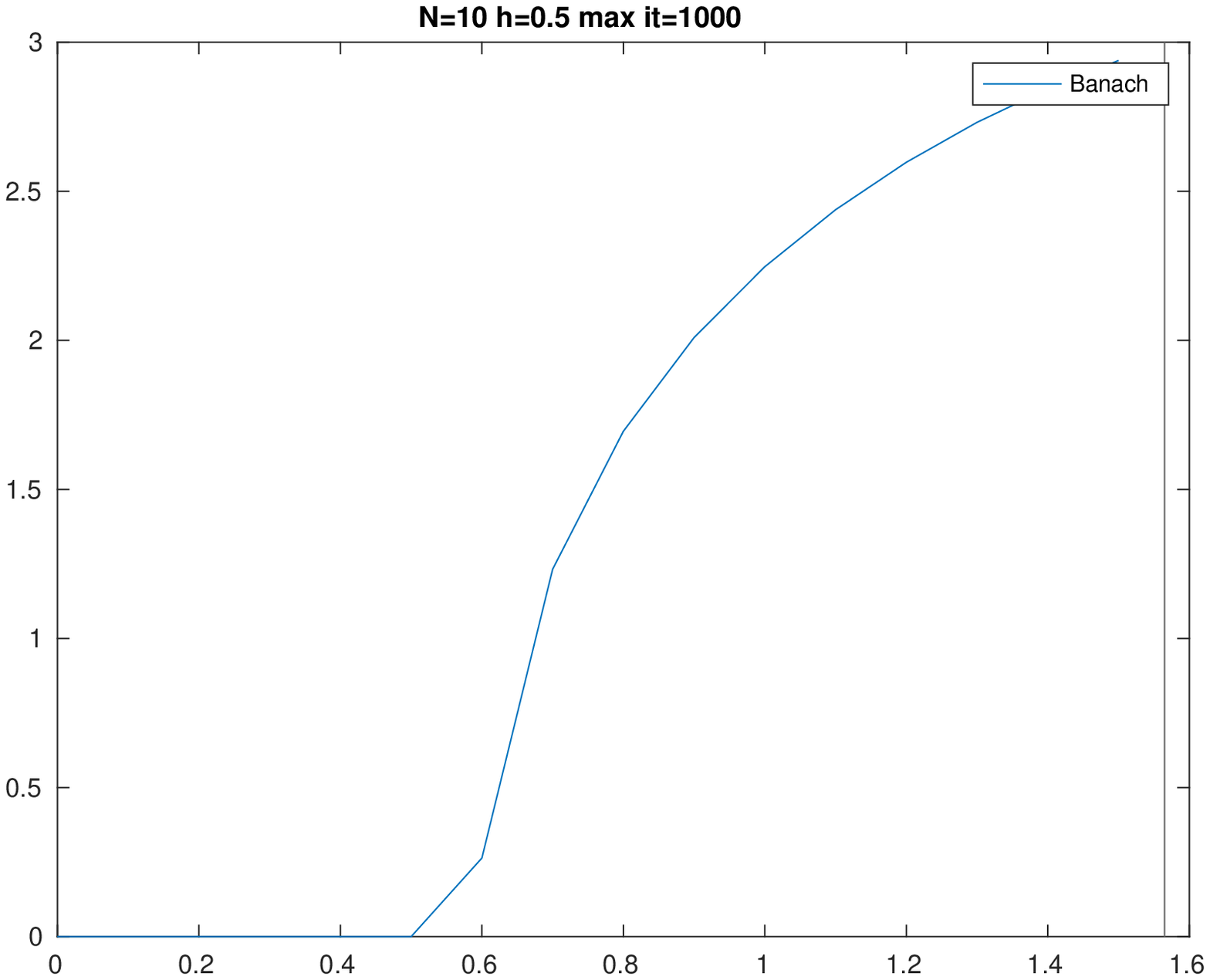}
  \caption{{\tiny} \space  1000 iterations of the standard Banach algorithm for $N=10$ and $h=0.5$. For $\beta\in [0, 1.5]$ we analyze the stability of this algorithm and see how quickly it takes off. To measure if the iterations stabilized we use the {\italic mean squared error} between the last two iterates. The vertical line marks the AT line.}
\end{figure}
\begin{figure}[H]
  \centering
  \includegraphics[width=0.5\textwidth,  trim={0 0 0 0},clip]{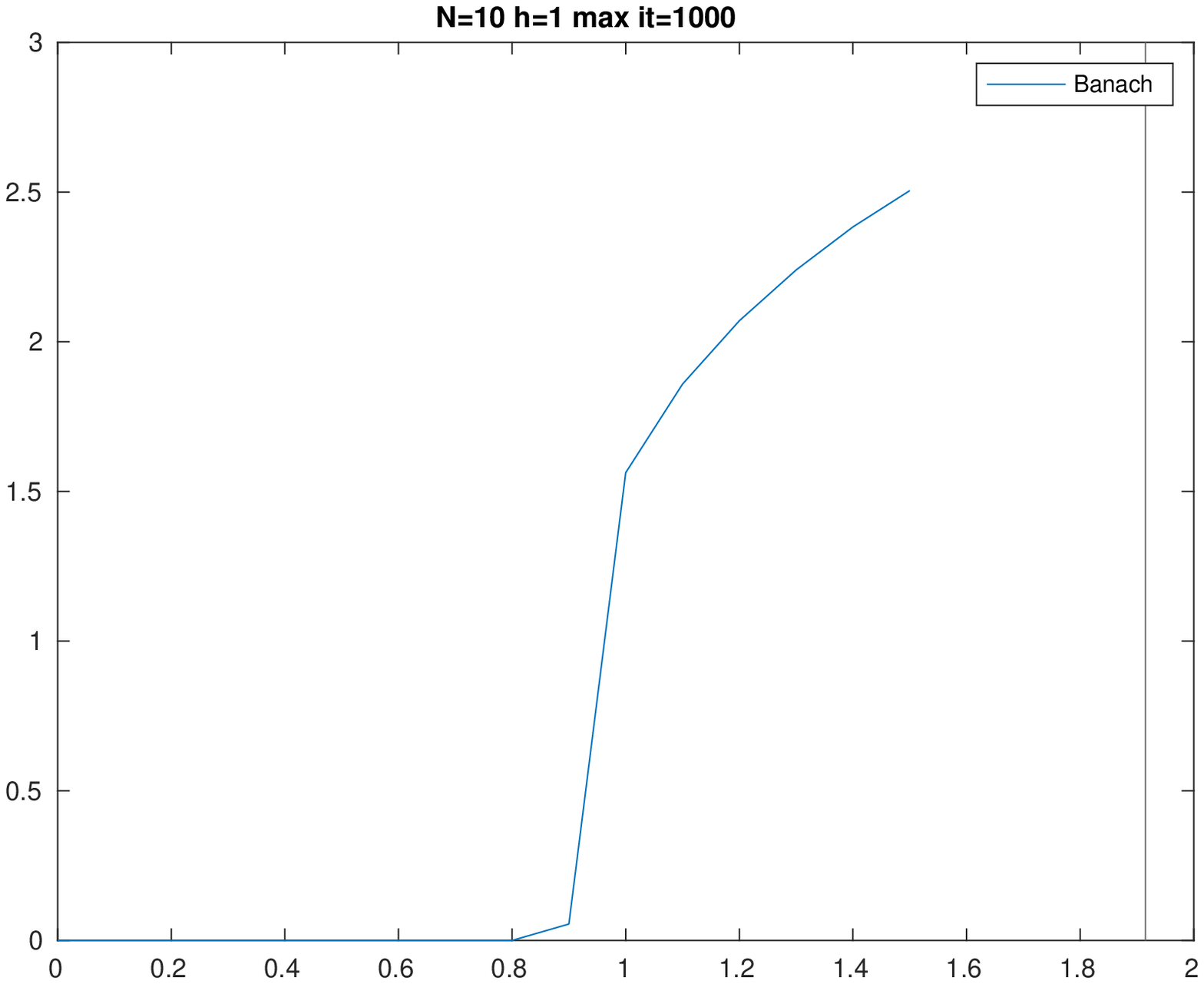}
  \caption{{\tiny} \space  1000 iterations of the standard Banach algorithm for $N=10$ and $h=1$. For $\beta\in [0, 1.5]$ we analyze the stability of this algorithm and see how quickly it takes off. To measure if the iterations stabilized we use the {\italic mean squared error} between the last two iterates. The vertical line marks the AT line.}
\end{figure}
\begin{figure}[H]
  \centering
  \includegraphics[width=0.5\textwidth,  trim={0 0 0 0},clip]{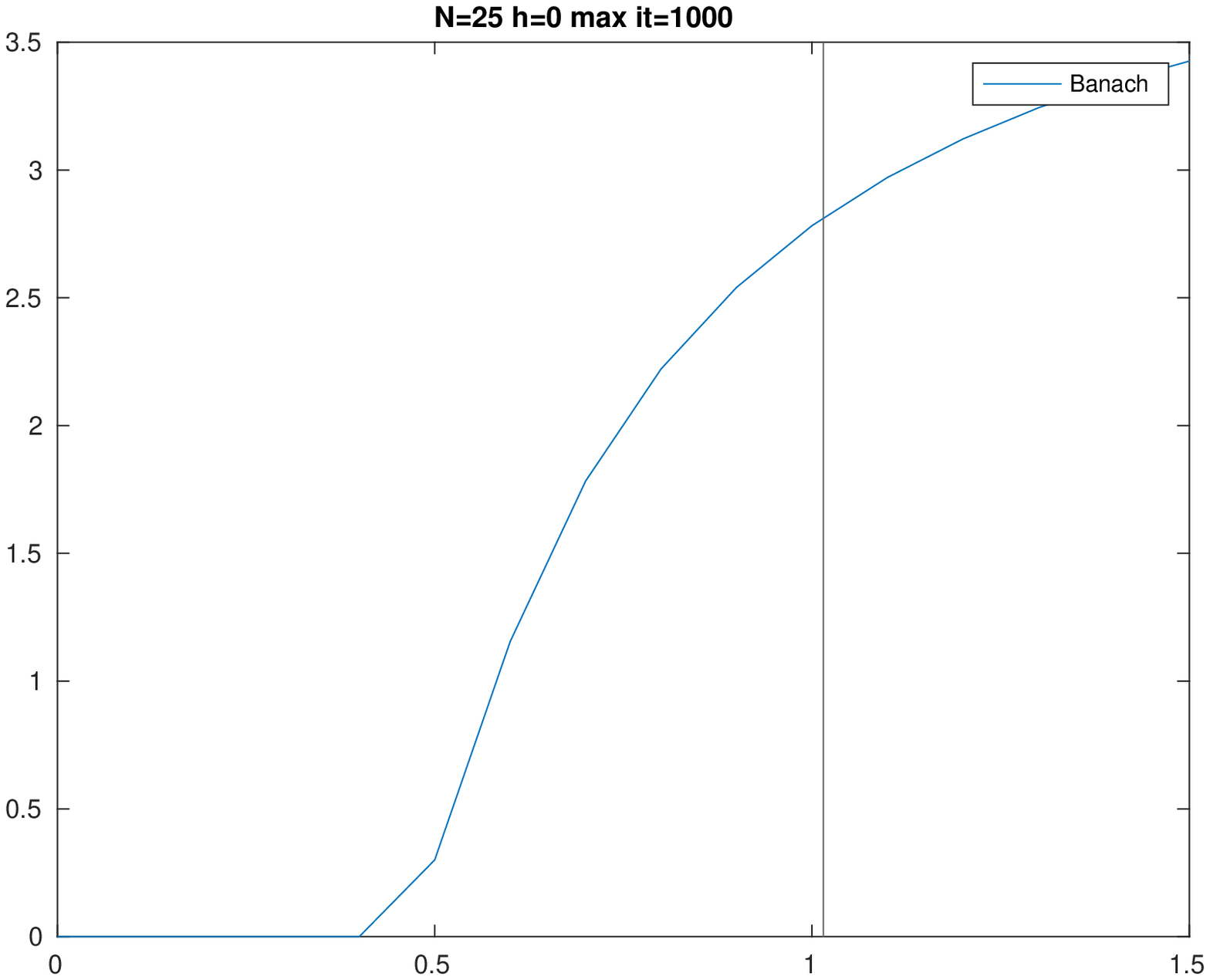}
  \caption{{\tiny} \space  1000 iterations of the standard Banach algorithm for $N=25$ and $h=0$. For $\beta\in [0, 1.5]$ we analyze the stability of this algorithm and see how quickly it takes off. To measure if the iterations stabilized we use the {\italic mean squared error} between the last two iterates. The vertical line marks the AT line.}
\end{figure}
\begin{figure}[H]
  \centering
  \includegraphics[width=0.5\textwidth,  trim={0 0 0 0},clip]{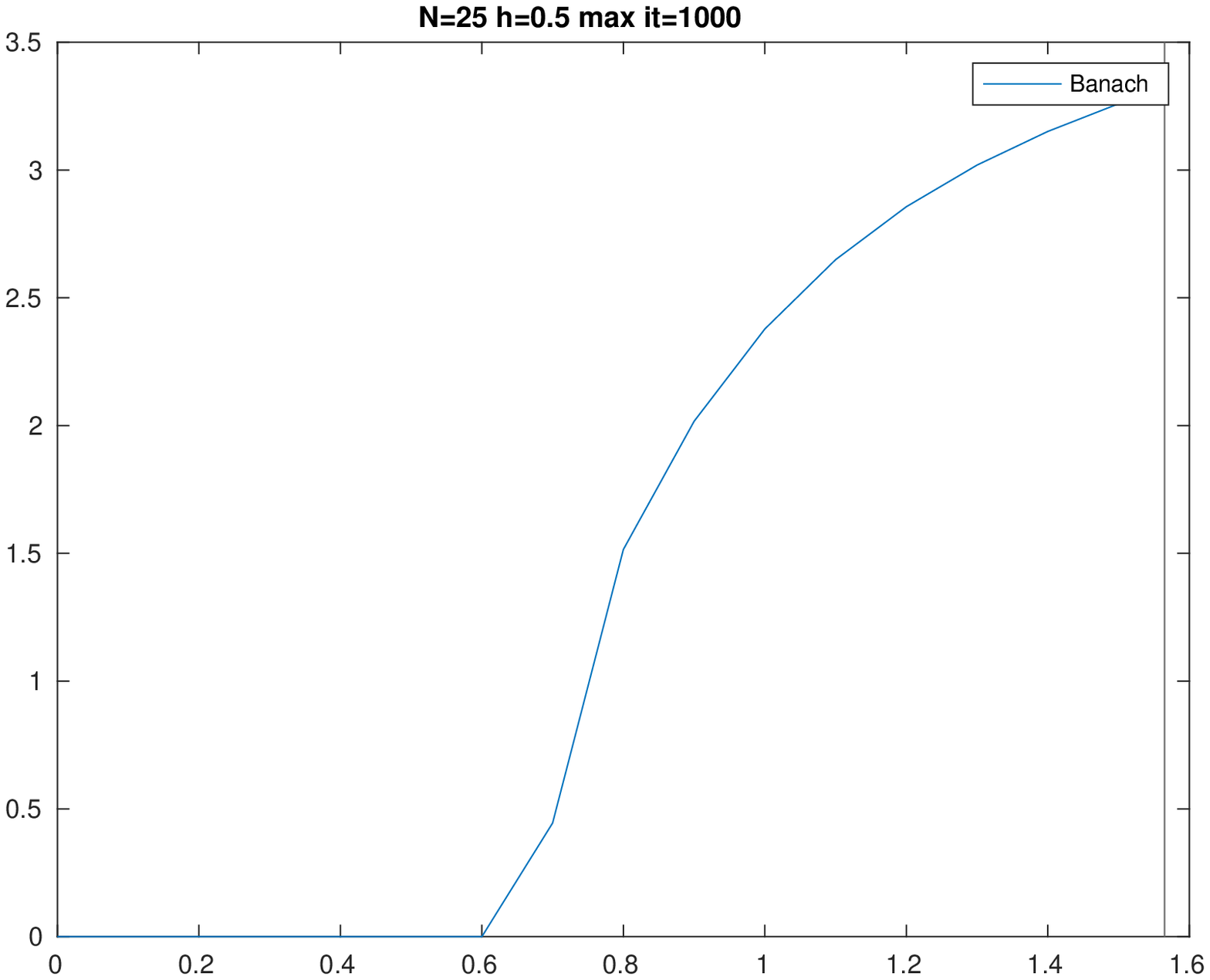}
  \caption{{\tiny} \space  1000 iterations of the standard Banach algorithm for $N=25$ and $h=0.5$. For $\beta\in [0, 1.5]$ we analyze the stability of this algorithm and see how quickly it takes off. To measure if the iterations stabilized we use the {\italic mean squared error} between the last two iterates. The vertical line marks the AT line.}
\end{figure}
\begin{figure}[H]
  \centering
  \includegraphics[width=0.5\textwidth,  trim={0 0 0 0},clip]{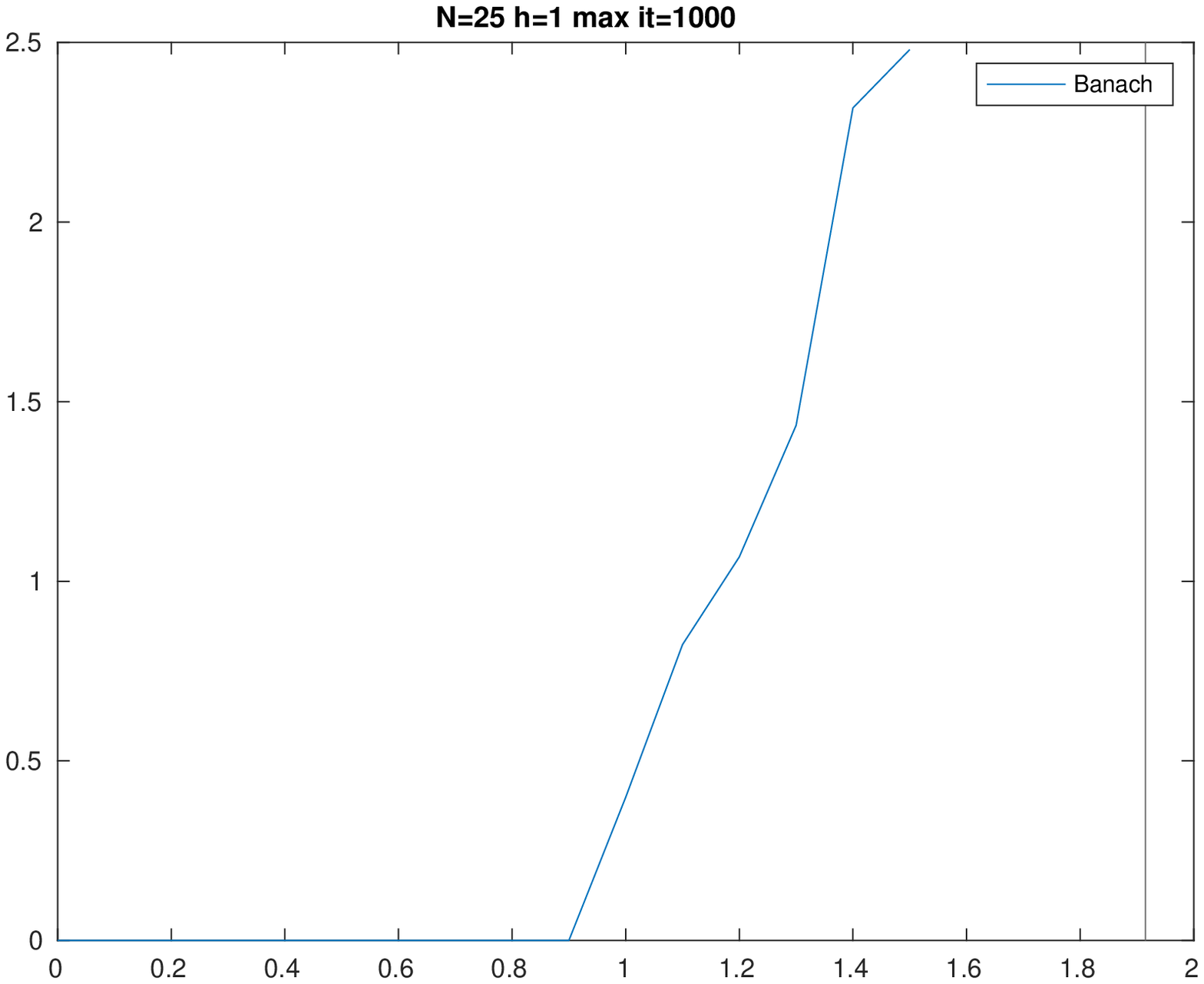}
  \caption{{\tiny} \space  1000 iterations of the standard Banach algorithm for $N=25$ and $h=1$. For $\beta\in [0, 1.5]$ we analyze the stability of this algorithm and see how quickly it takes off. To measure if the iterations stabilized we use the {\italic mean squared error} between the last two iterates. The vertical line marks the AT line.}
\end{figure}
\begin{figure}[H]
  \centering
  \includegraphics[width=0.5\textwidth,  trim={0 0 0 0},clip]{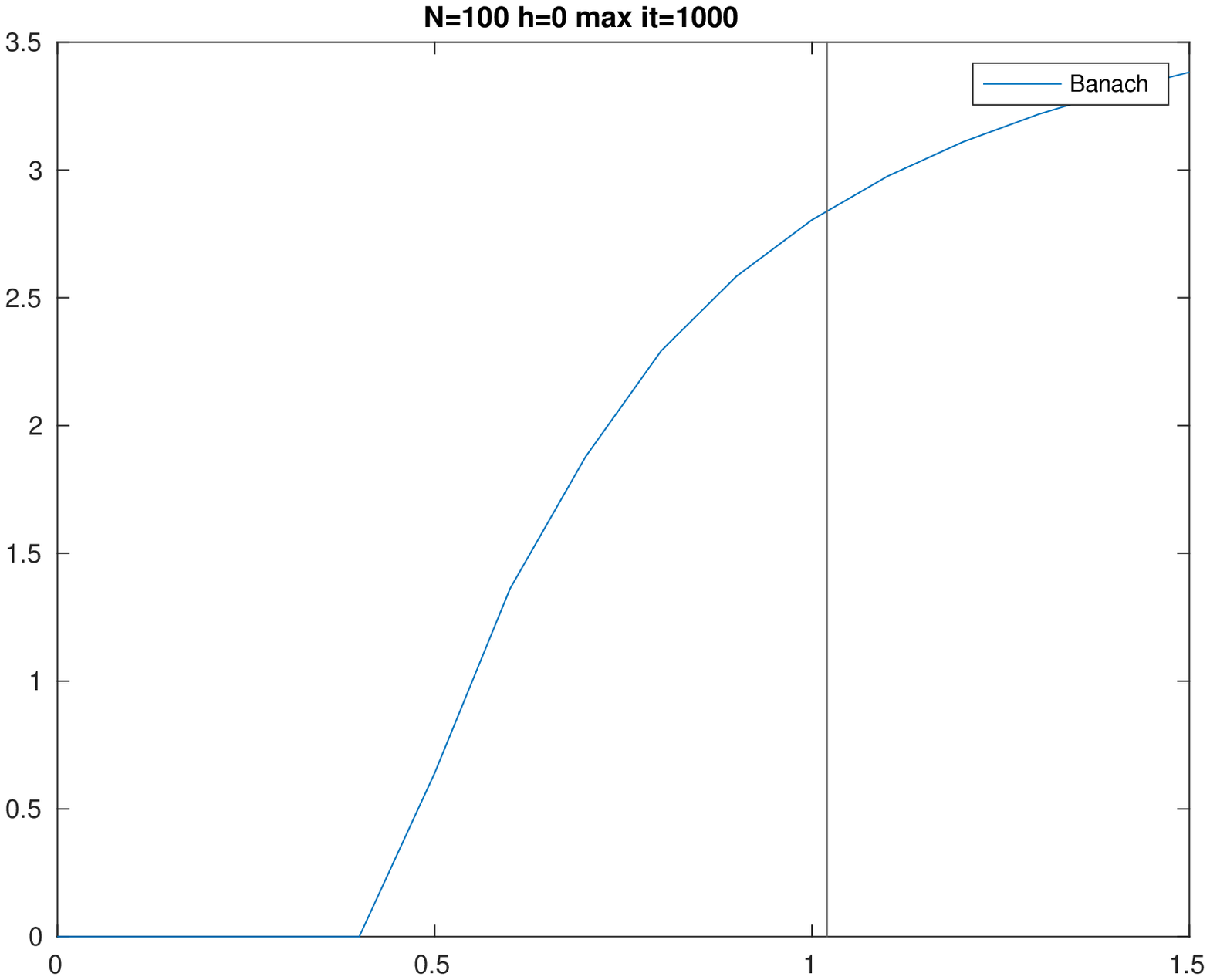}
  \caption{{\tiny} \space  1000 iterations of the standard Banach algorithm for $N=100$ and $h=0$. For $\beta\in [0, 1.5]$ we analyze the stability of this algorithm and see how quickly it takes off. To measure if the iterations stabilized we use the {\italic mean squared error} between the last two iterates. The vertical line marks the AT line.}
\end{figure}
\begin{figure}[H]
  \centering
  \includegraphics[width=0.5\textwidth,  trim={0 0 0 0},clip]{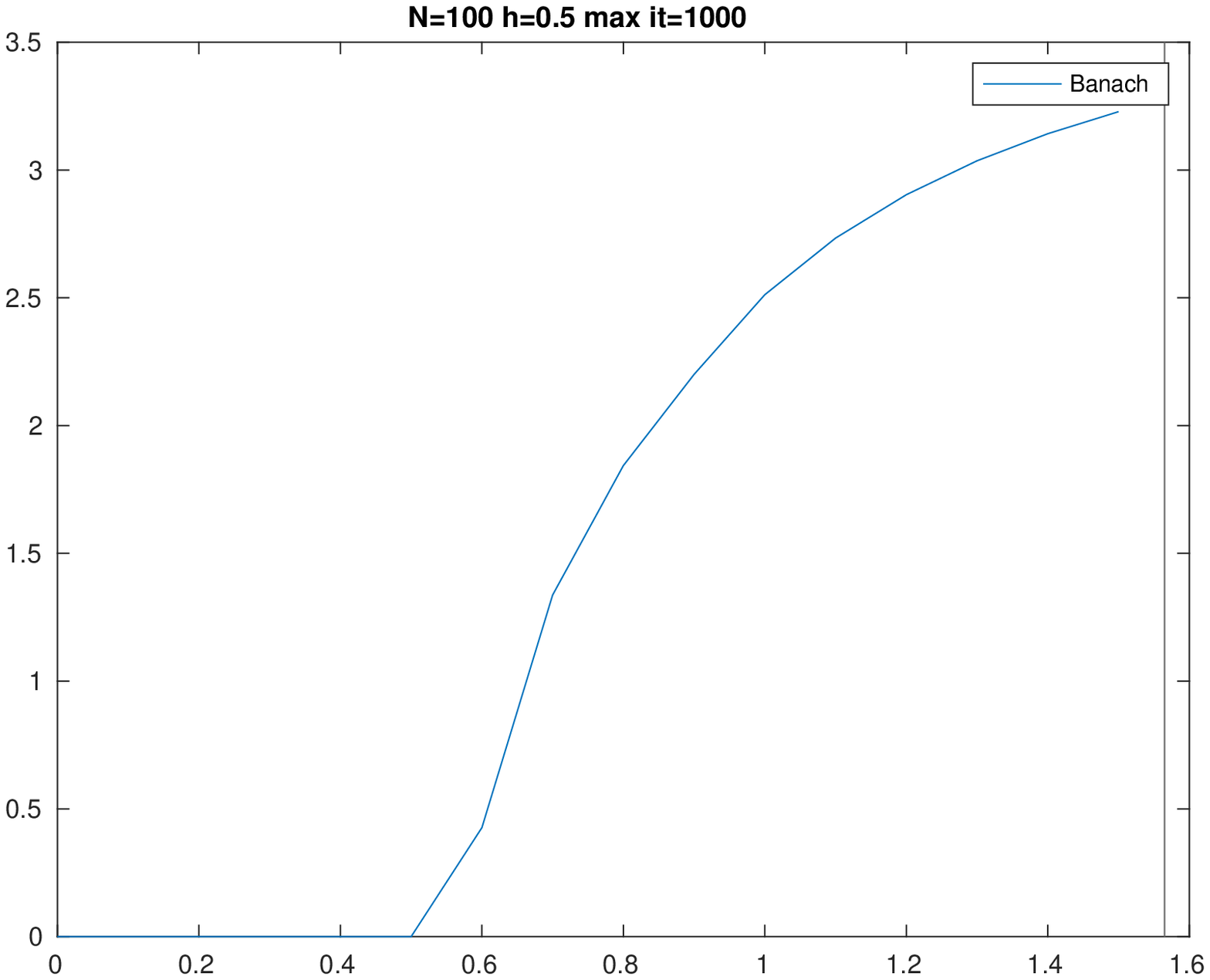}
  \caption{{\tiny} \space  1000 iterations of the standard Banach algorithm for $N=100$ and $h=0.5$. For $\beta\in [0, 1.5]$ we analyze the stability of this algorithm and see how quickly it takes off. To measure if the iterations stabilized we use the {\italic mean squared error} between the last two iterates. The vertical line marks the AT line.}
\end{figure}
\begin{figure}[H]
  \centering
  \includegraphics[width=0.5\textwidth,  trim={0 0 0 0},clip]{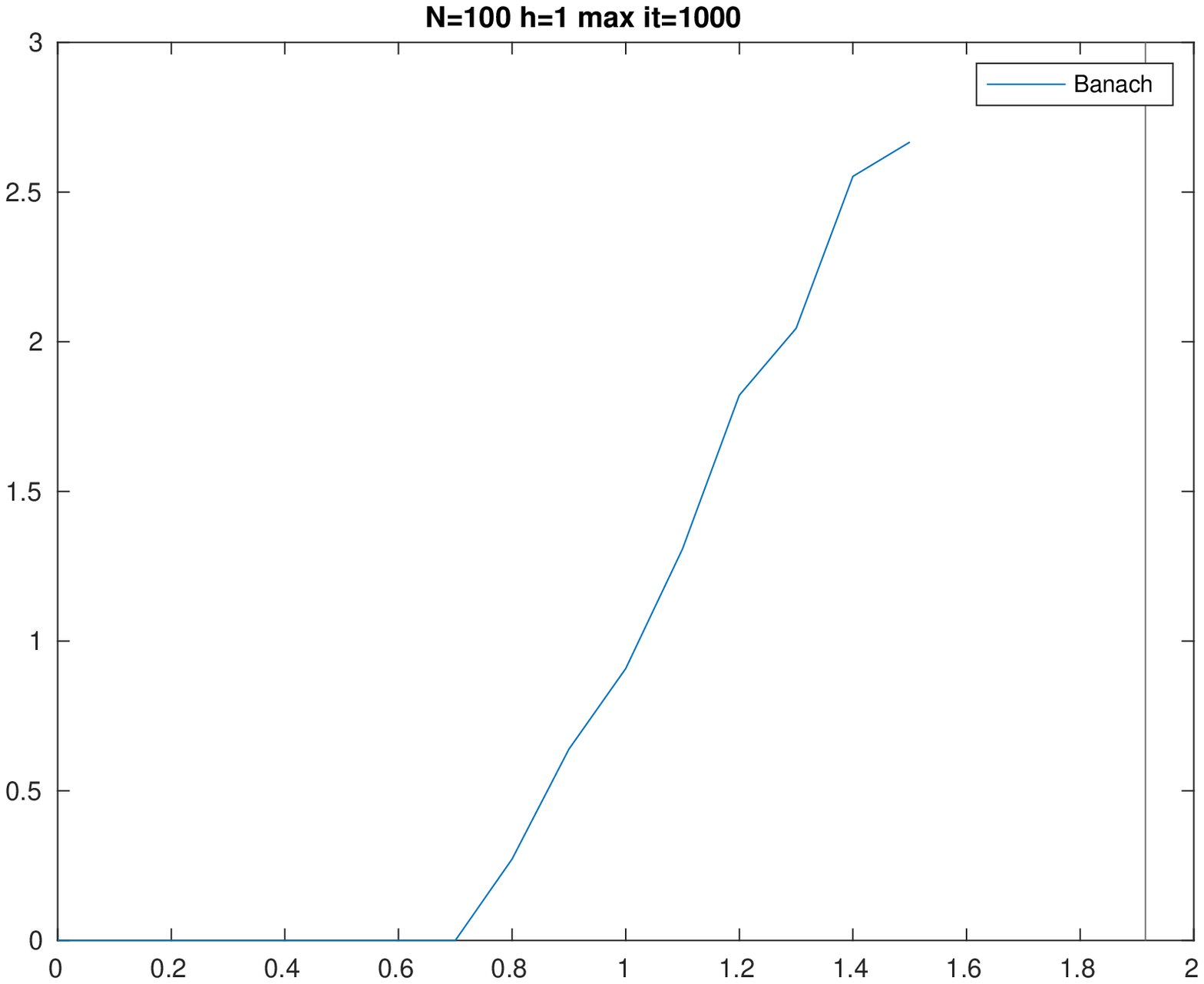}
  \caption{{\tiny} \space  1000 iterations of the standard Banach algorithm for $N=100$ and $h=1$. For $\beta\in [0, 1.5]$ we analyze the stability of this algorithm and see how quickly it takes off. To measure if the iterations stabilized we use the {\italic mean squared error} between the last two iterates. The vertical line marks the AT line.}
\end{figure}
\begin{figure}[H]
  \centering
  \includegraphics[width=0.5\textwidth,  trim={0 0 0 0},clip]{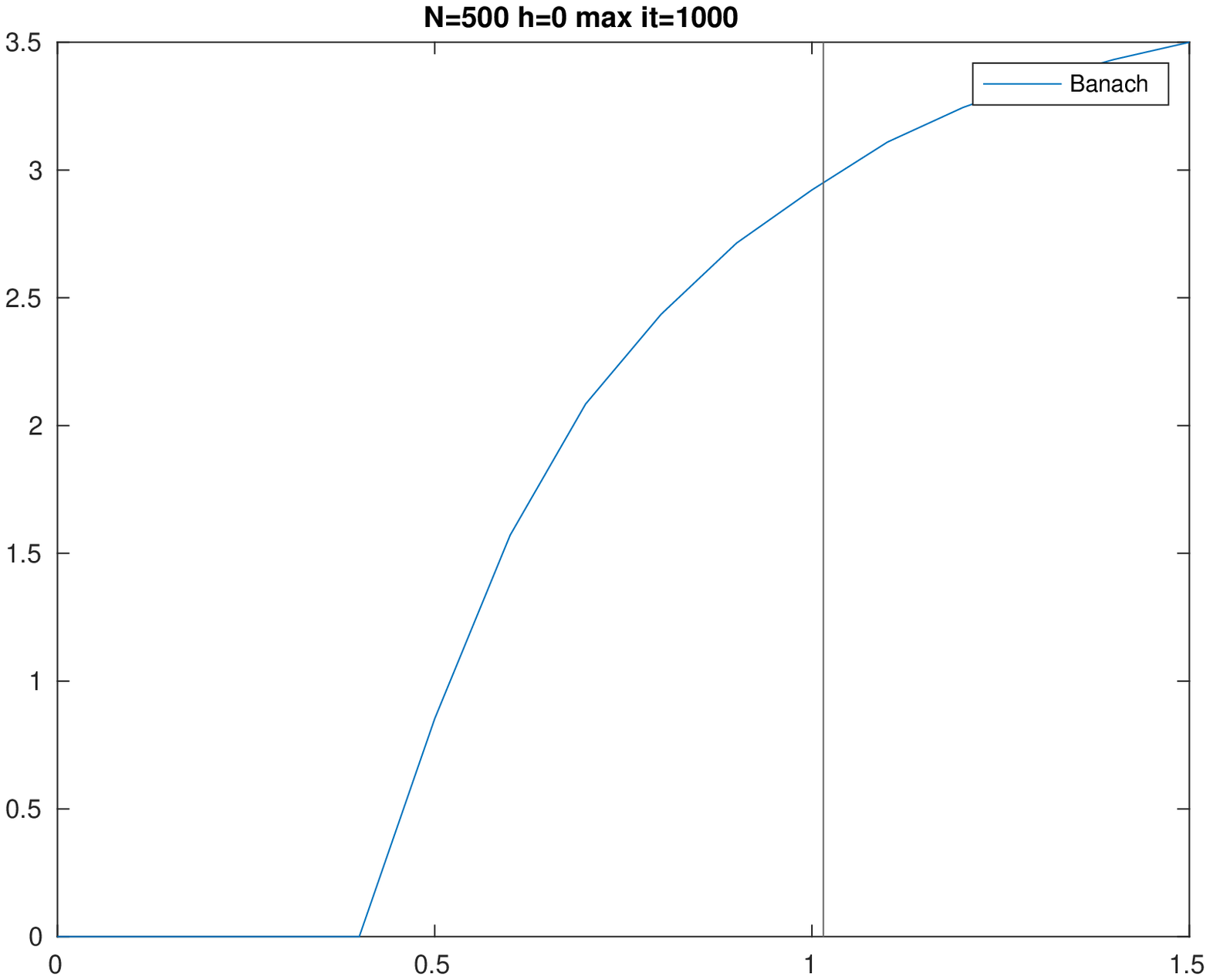}
  \caption{{\tiny} \space  1000 iterations of the standard Banach algorithm for $N=500$ and $h=0$. For $\beta\in [0, 1.5]$ we analyze the stability of this algorithm and see how quickly it takes off. To measure if the iterations stabilized we use the {\italic mean squared error} between the last two iterates. The vertical line marks the AT line.}
\end{figure}
\begin{figure}[H]
  \centering
  \includegraphics[width=0.5\textwidth,  trim={0 0 0 0},clip]{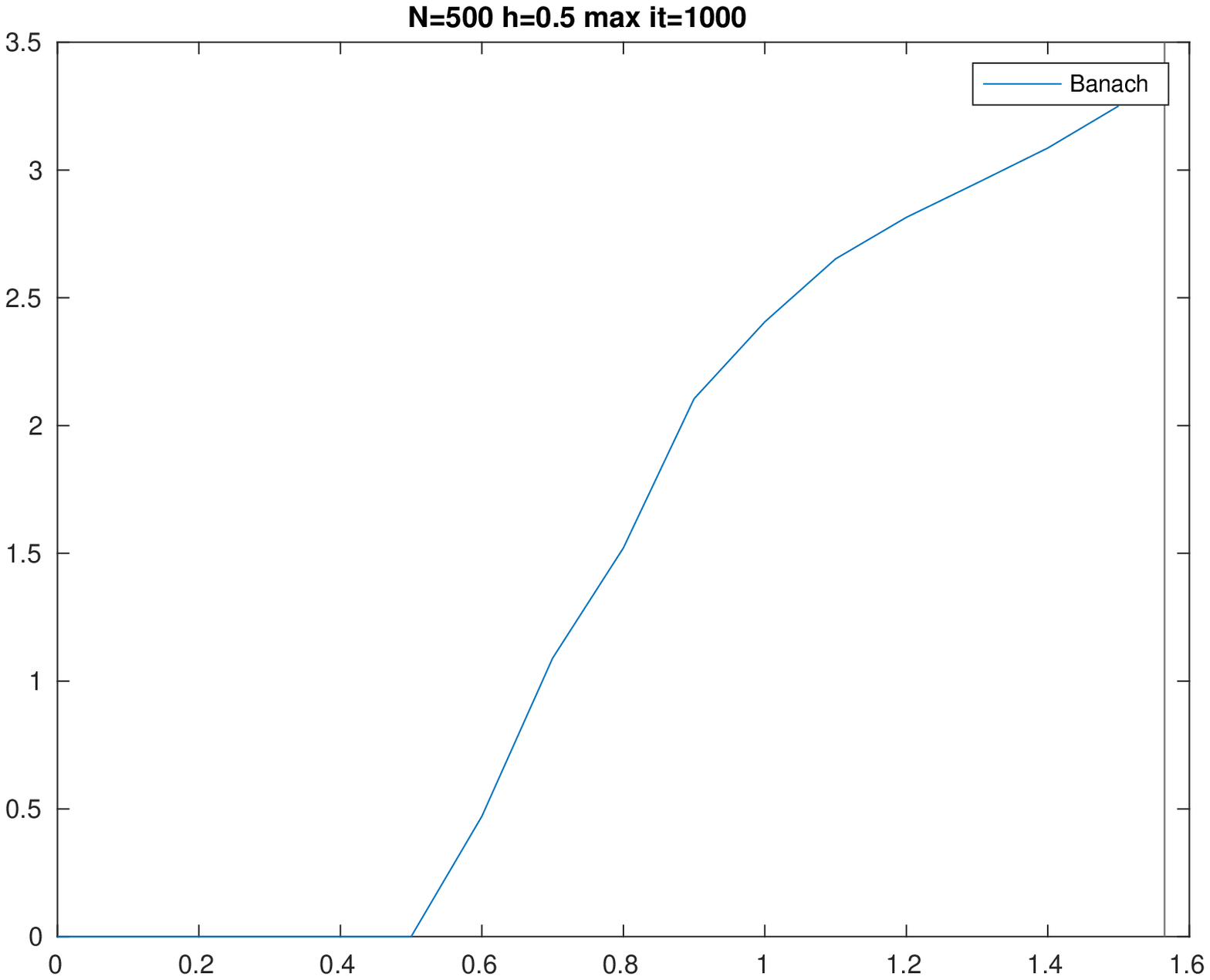}
  \caption{{\tiny} \space  1000 iterations of the standard Banach algorithm for $N=500$ and $h=0.5$. For $\beta\in [0, 1.5]$ we analyze the stability of this algorithm and see how quickly it takes off. To measure if the iterations stabilized we use the {\italic mean squared error} between the last two iterates. The vertical line marks the AT line.}
\end{figure}
\begin{figure}[H]
  \centering
  \includegraphics[width=0.5\textwidth,  trim={0 0 0 0},clip]{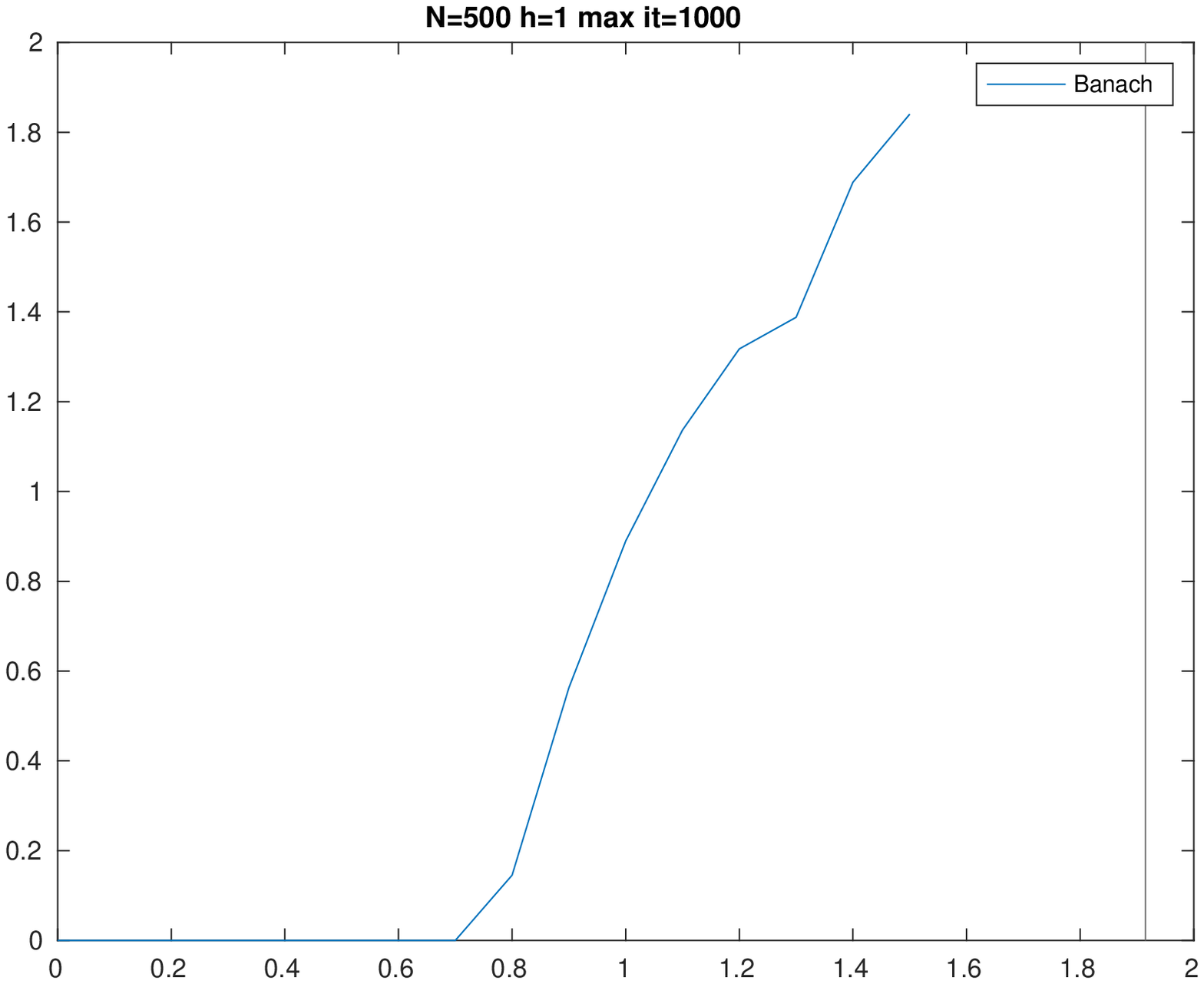}
  \caption{{\tiny} \space  1000 iterations of the standard Banach algorithm for $N=500$ and $h=1$. For $\beta\in [0, 1.5]$ we analyze the stability of this algorithm and see how quickly it takes off. To measure if the iterations stabilized we use the {\italic mean squared error} between the last two iterates. The vertical line marks the AT line.}
\end{figure}
\begin{figure}[H]
  \centering
  \includegraphics[width=0.5\textwidth,  trim={0 0 0 0},clip]{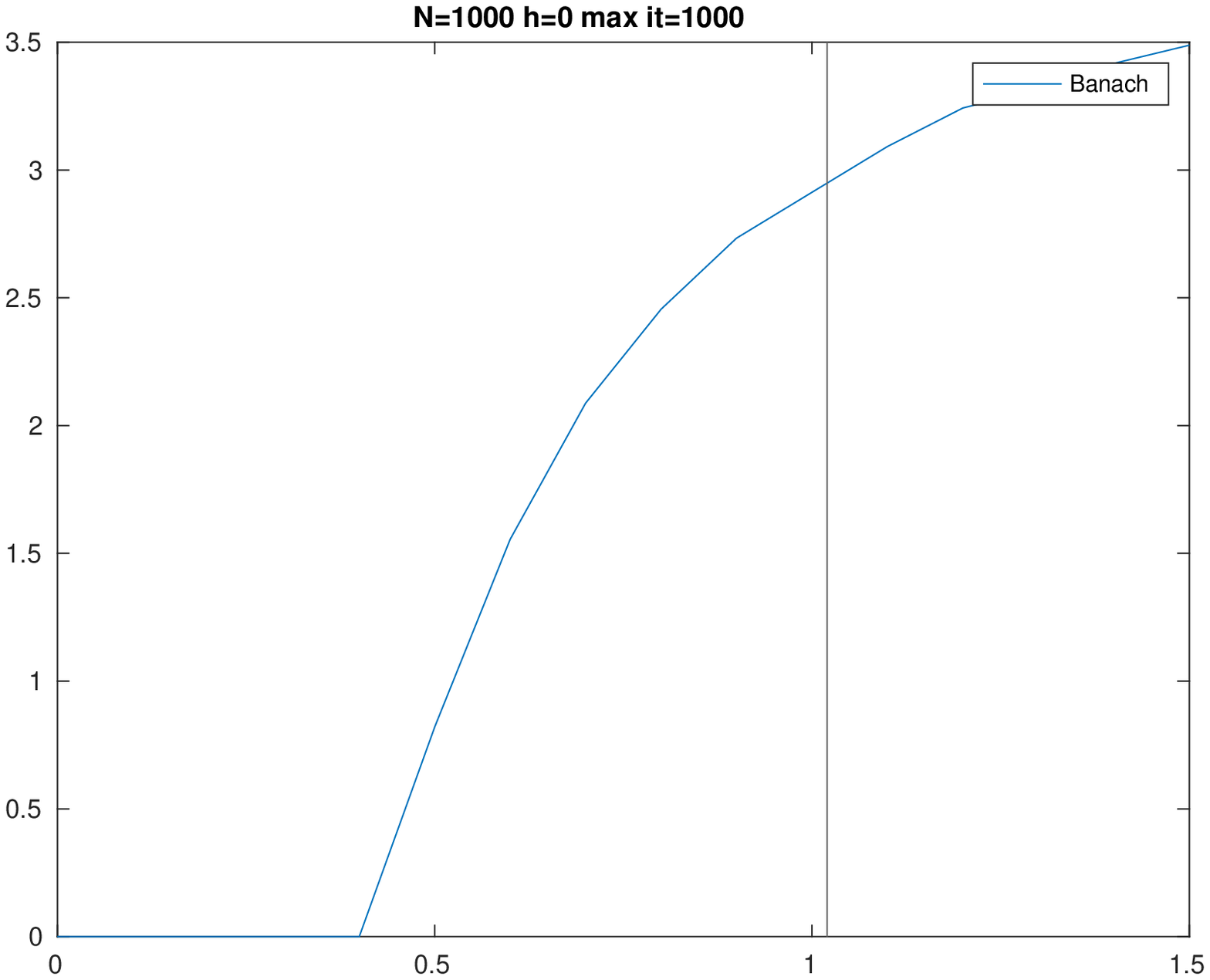}
  \caption{{\tiny} \space  1000 iterations of the standard Banach algorithm for $N=1000$ and $h=0$. For $\beta\in [0, 1.5]$ we analyze the stability of this algorithm and see how quickly it takes off. To measure if the iterations stabilized we use the {\italic mean squared error} between the last two iterates. The vertical line marks the AT line.}
\end{figure}
\begin{figure}[H]
  \centering
  \includegraphics[width=0.5\textwidth,  trim={0 0 0 0},clip]{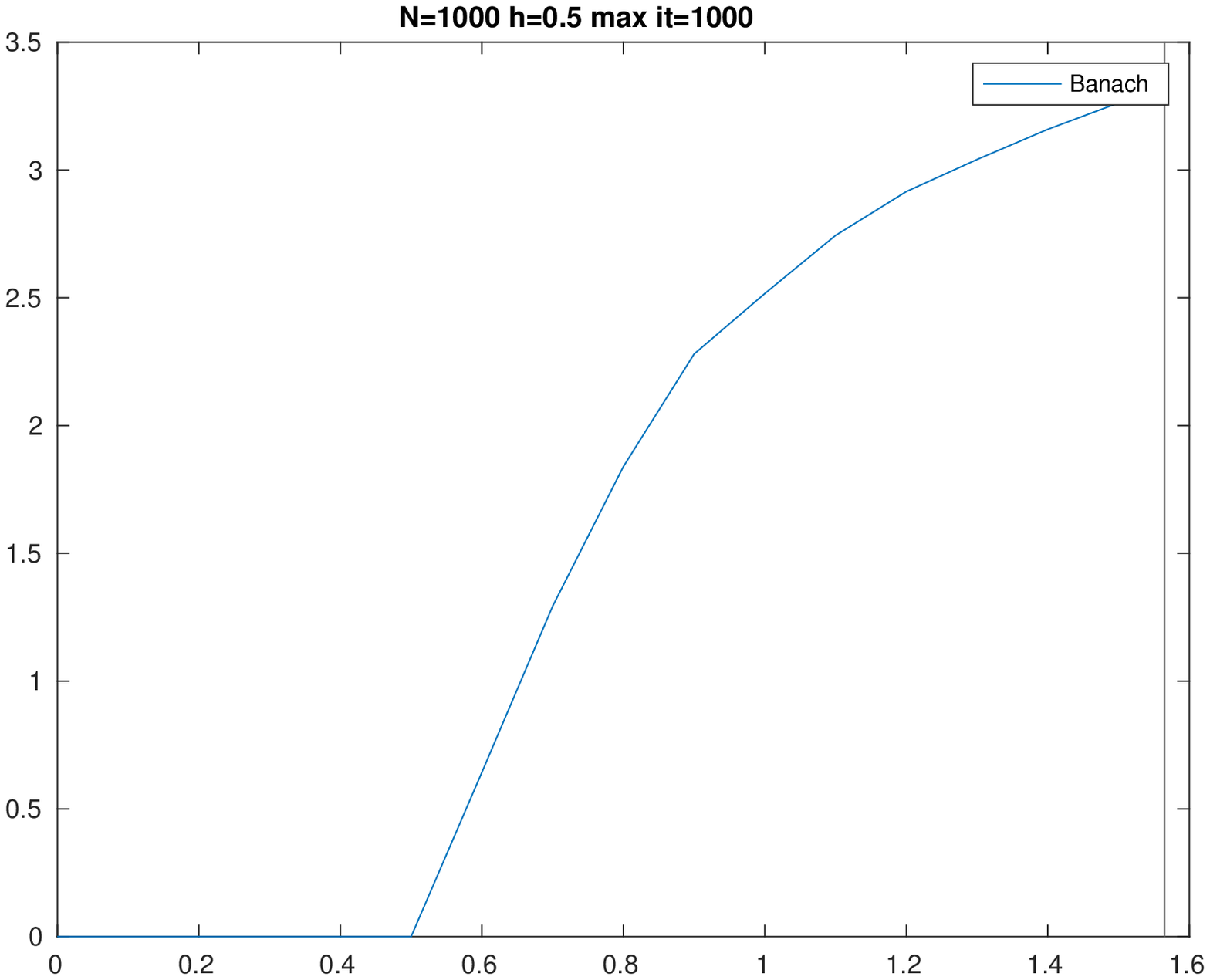}
  \caption{{\tiny} \space  1000 iterations of the standard Banach algorithm for $N=1000$ and $h=0.5$. For $\beta\in [0, 1.5]$ we analyze the stability of this algorithm and see how quickly it takes off. To measure if the iterations stabilized we use the {\italic mean squared error} between the last two iterates. The vertical line marks the AT line.}
\end{figure}
\begin{figure}[H]
  \centering
  \includegraphics[width=0.5\textwidth,  trim={0 0 0 0},clip]{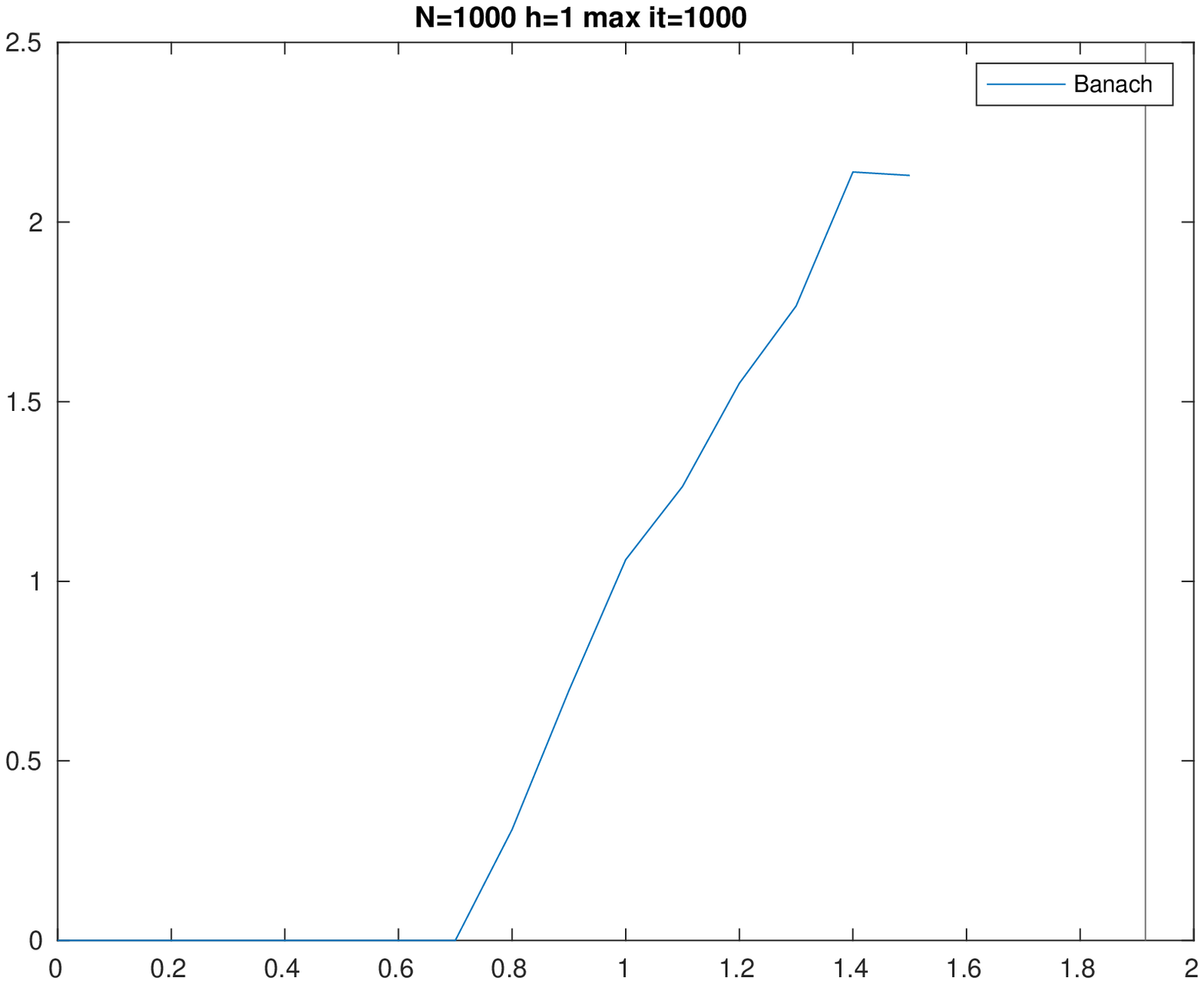}
  \caption{{\tiny} \space  1000 iterations of the standard Banach algorithm for $N=1000$ and $h=1$. For $\beta\in [0, 1.5]$ we analyze the stability of this algorithm and see how quickly it takes off. To measure if the iterations stabilized we use the {\italic mean squared error} between the last two iterates. The vertical line marks the AT line.}
\end{figure}
\FloatBarrier\clearpage
\subsection{Short summary}
We see that no matter how we choose the start values the standard Banach-algorithm does not converge (no matter which of the two measures for the stability we use!) for $\beta$ large enough where \emph{large enough} is way before the AT-line. We see that for $h=0, h=0.5, h=1$ and $N=25, 100, 500, 1000$. 
\FloatBarrier\clearpage
\section{$\varepsilon$-Banach and Two-Step-Banach in high temperature for $\beta=1$ and $h=0.5$}
\subsection{Uniformly in the corners of the hypercube distributed start values}\leavevmode
\begin{figure}[H]
  \centering
  \includegraphics[width=0.5\textwidth,  trim={0.45cm 0.4cm 0 0.4cm},clip]{rsynccopy/final_simulations/starting_in_uniform_selected_corners_of_the_cube/banach_epsilon_and_twostep_beta_1/1div.eps}
  \caption{{\tiny} \space  1000 iterations of the $\varepsilon$-Banach algorithm for $N=10, \beta=1$ and $h=0.5$. To measure if the iterations stabilized we use the {\italic mean squared error} between the last two iterates. We simulated 1000 uniformly in the corners of the hypercube distributed start values for $\varepsilon = -0.9, -0.85,\ldots, 0.0$. The {\italic mean squared error} seems to depend highly on the choice of $\varepsilon$. It's obvious that the {\italic mean squared error} will be smaller for $\epsilon\rightarrow 1$ and this does not mean that solutions found with those $\varepsilon$ are better ones. }
\end{figure}
\begin{figure}[H]
  \centering
  \includegraphics[width=0.5\textwidth,  trim={0.45cm 0.4cm 0 0.4cm},clip]{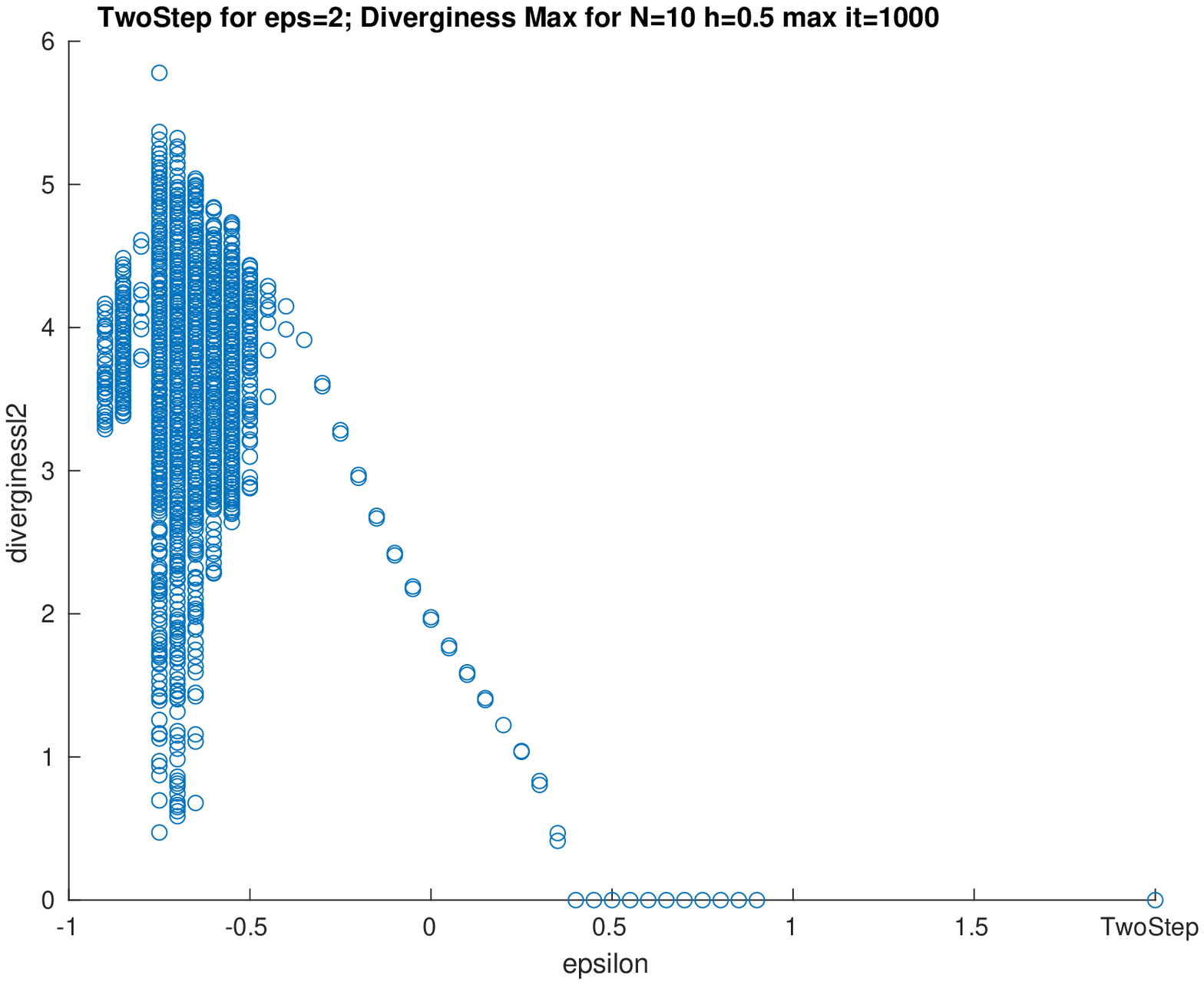}
  \caption{{\tiny} \space  1000 iterations of the $\varepsilon$-Banach algorithm for $N=10, \beta=1$ and $h=0.5$.To measure if the iterations stabilized we use the {\italic maximum absolute error} between the last two iterates. We simulated 1000 uniformly in the corners of the hypercube distributed start values for $\varepsilon = -0.9, -0.85,\ldots, 0.0$. The {\italic maximum absolute error} seems to depend highly on the choice of $\varepsilon$. It's obvious that the {\italic maximum absolute error} will be smaller for $\epsilon\rightarrow 1$ and this does not mean that solutions found with those $\varepsilon$ are better ones. }
\end{figure}
\begin{figure}[H]
  \centering
  % \psfrag{2}{\tiny Two Step}
  \includegraphics[width=0.5\textwidth,  trim={0.45cm 0.4cm 0 0.4cm},clip]{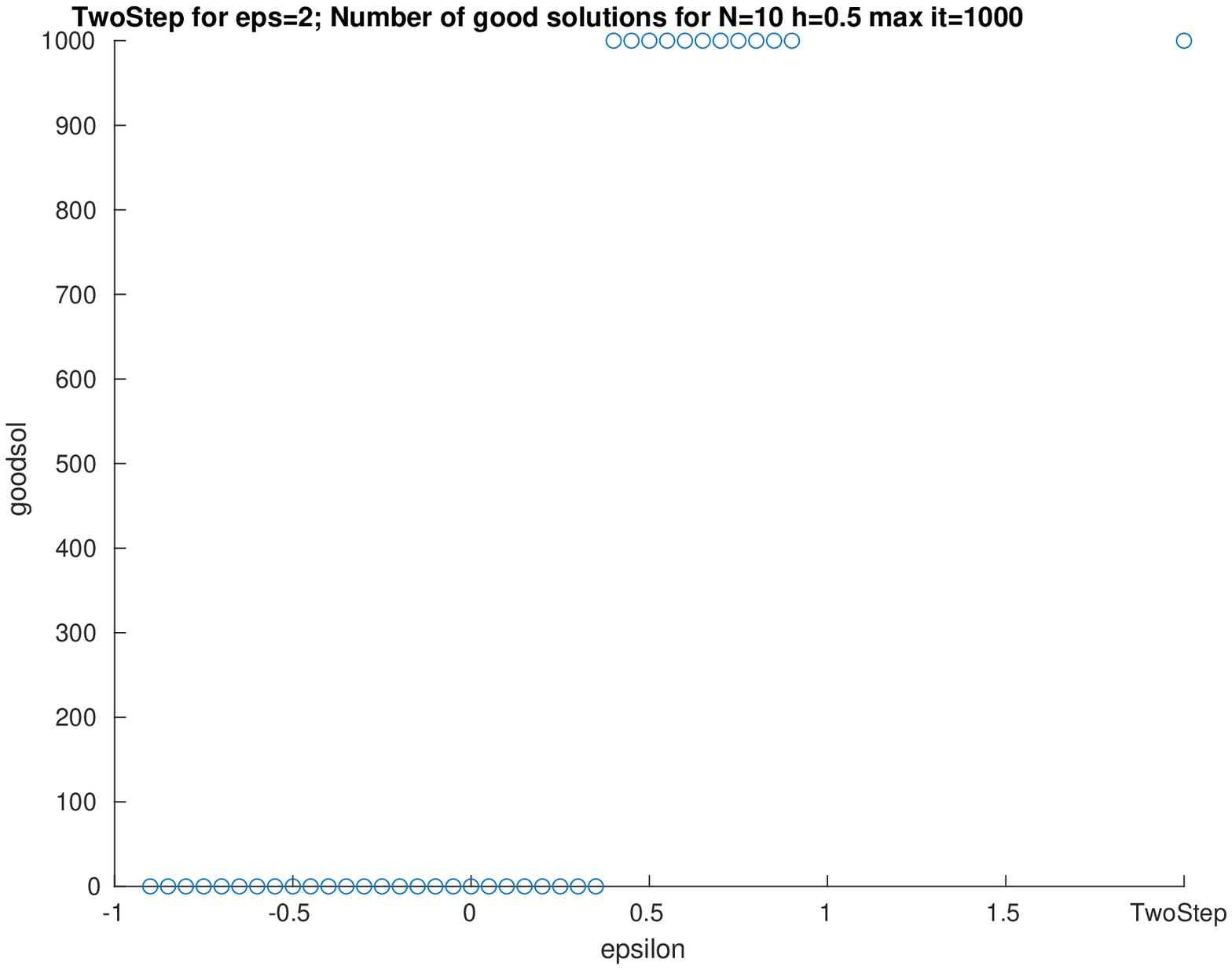}
  \caption{{\tiny} \space  1000 iterations of the $\varepsilon$-Banach algorithm for $N=10, \beta=1$ and $h=0.5$. To measure if the iterations stabilized we use the {\italic mean squared error} between the last two iterates. We simulated 1000 uniformly in the corners of the hypercube distributed start values for $\varepsilon = -0.9, -0.85,\ldots, 0.0$. We say that a found fixed point is a good one if the found TAP free energy is between  $0.35383069483591\pm 0.05$ ($0.35383069483591$ is the RS Parisi FE), the {\italic mean squared error} is smaller than 0.0001 and the Plefka condition is fulfilled.  }
\end{figure}
\begin{figure}[H]
  \centering
  % \psfrag{2}{\tiny Two Step}
  \includegraphics[width=0.5\textwidth,  trim={0.45cm 0.4cm 0 0.4cm},clip]{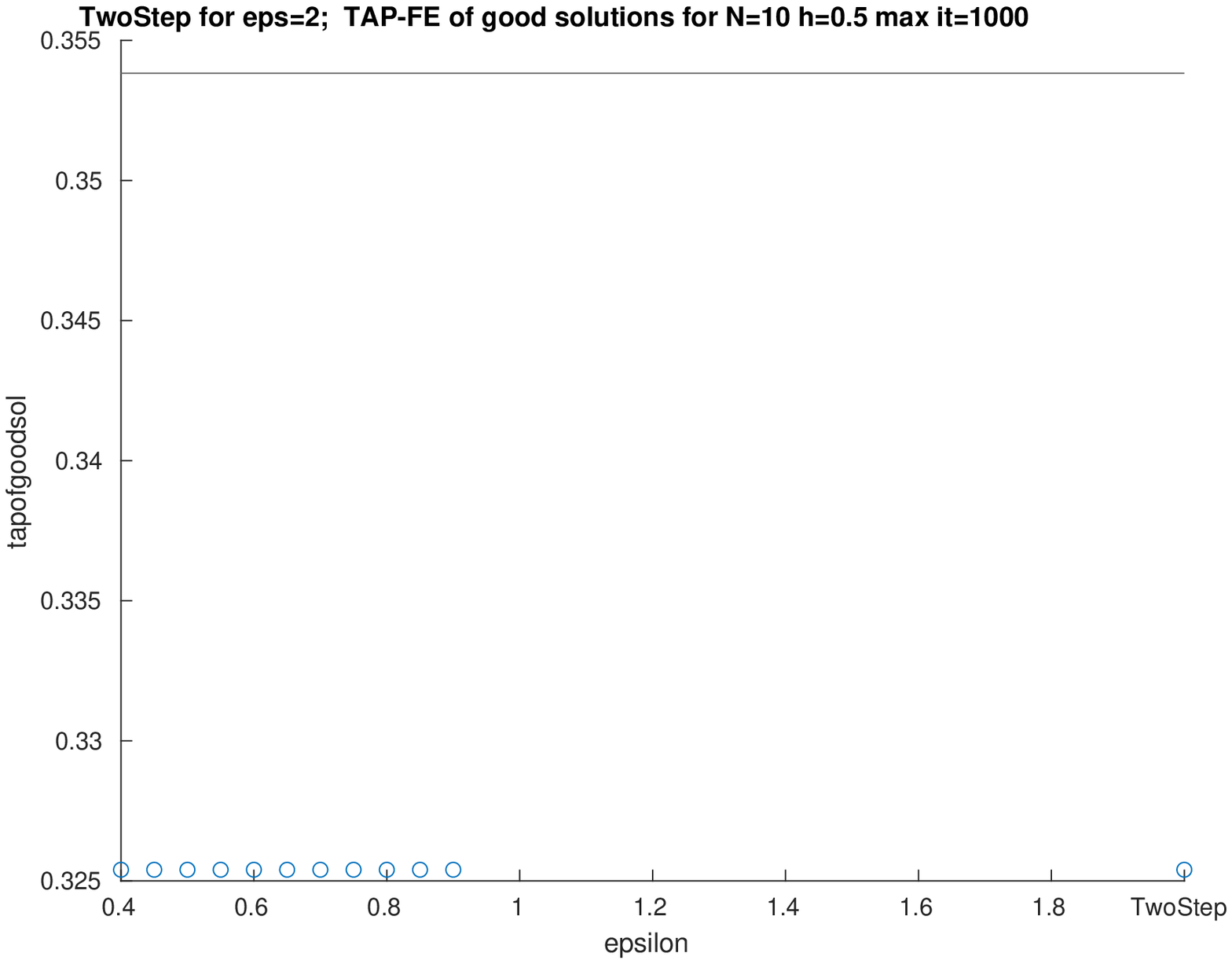}
  \caption{{\tiny} \space  1000 iterations of the $\varepsilon$-Banach algorithm for $N=10, \beta=1$ and $h=0.5$. We simulated 1000 uniformly in the corners of the hypercube distributed start values for $\varepsilon = -0.9, -0.85,\ldots, 0.0$. We say that a found fixed point is a good one if the found TAP free energy is between  $0.35383069483591\pm 0.05$ ($0.35383069483591$ is the RS Parisi FE), the {\italic mean squared error} is smaller than 0.0001 and the Plefka condition is fulfilled. We plotted all found values of the TAP free energy values of good iterations and excluded all iterations which landed outside the hypercube. The vertical line marks the Parisi free energy value found with the RS formulae.}
\end{figure}

\begin{figure}[H]
  \centering
  % \psfrag{2}{\tiny Two Step}
  \includegraphics[width=0.5\textwidth,  trim={0.45cm 0.4cm 0 0.4cm},clip]{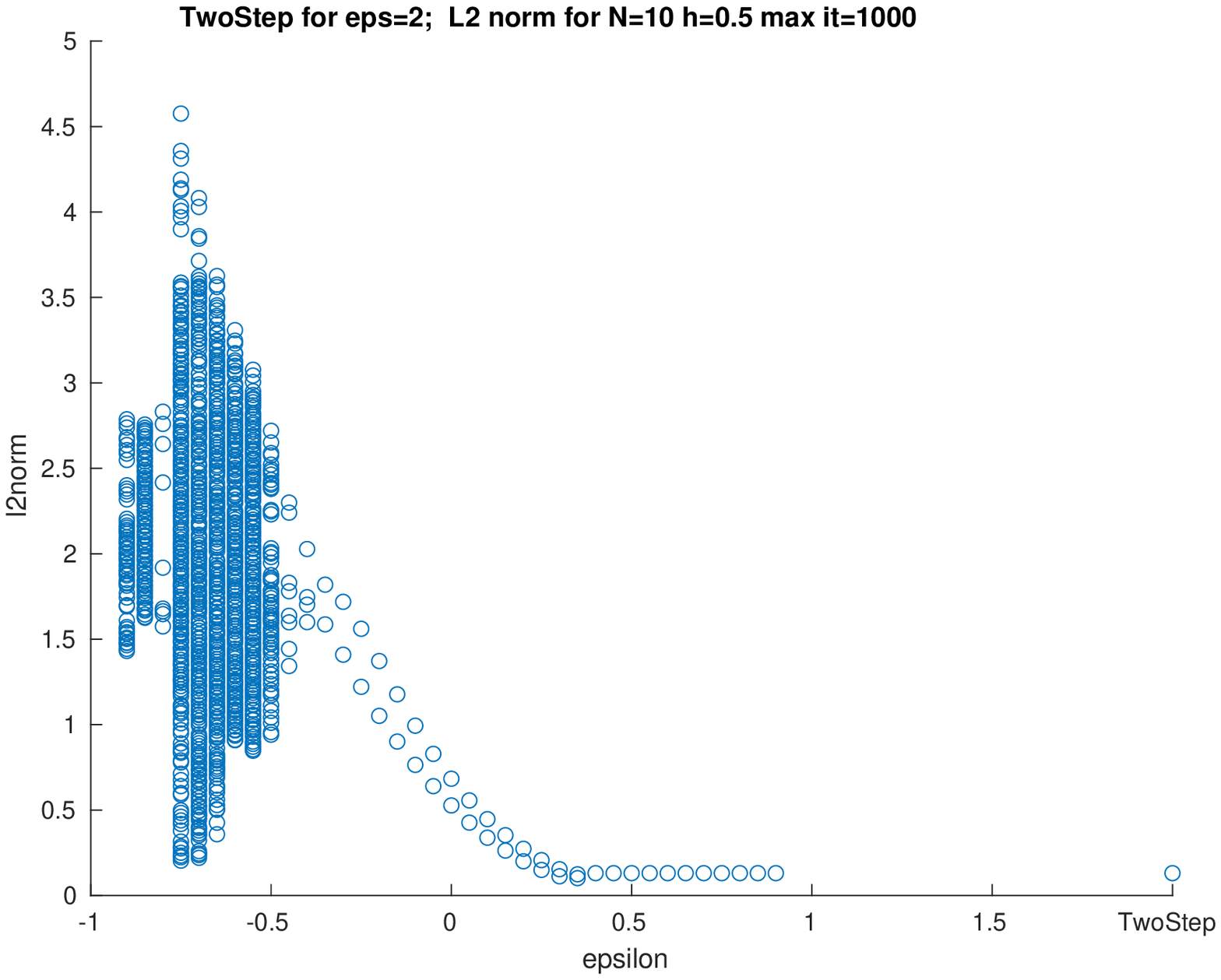}
  \caption{{\tiny} \space  1000 iterations of the $\varepsilon$-Banach algorithm for $N=10, \beta=1$ and $h=0.5$. We simulated 1000 uniformly  in the corners of the hypercube distributed start values for $\varepsilon = -0.9, -0.85,\ldots, 0.0$. We plotted all found {\it mean squared errors} between the last two iterates for found fixed points.   }
\end{figure}
\begin{figure}[H]
  \centering
  % \psfrag{2}{\tiny Two Step}
  \includegraphics[width=0.5\textwidth,  trim={0.45cm 0.4cm 0 0.4cm},clip]{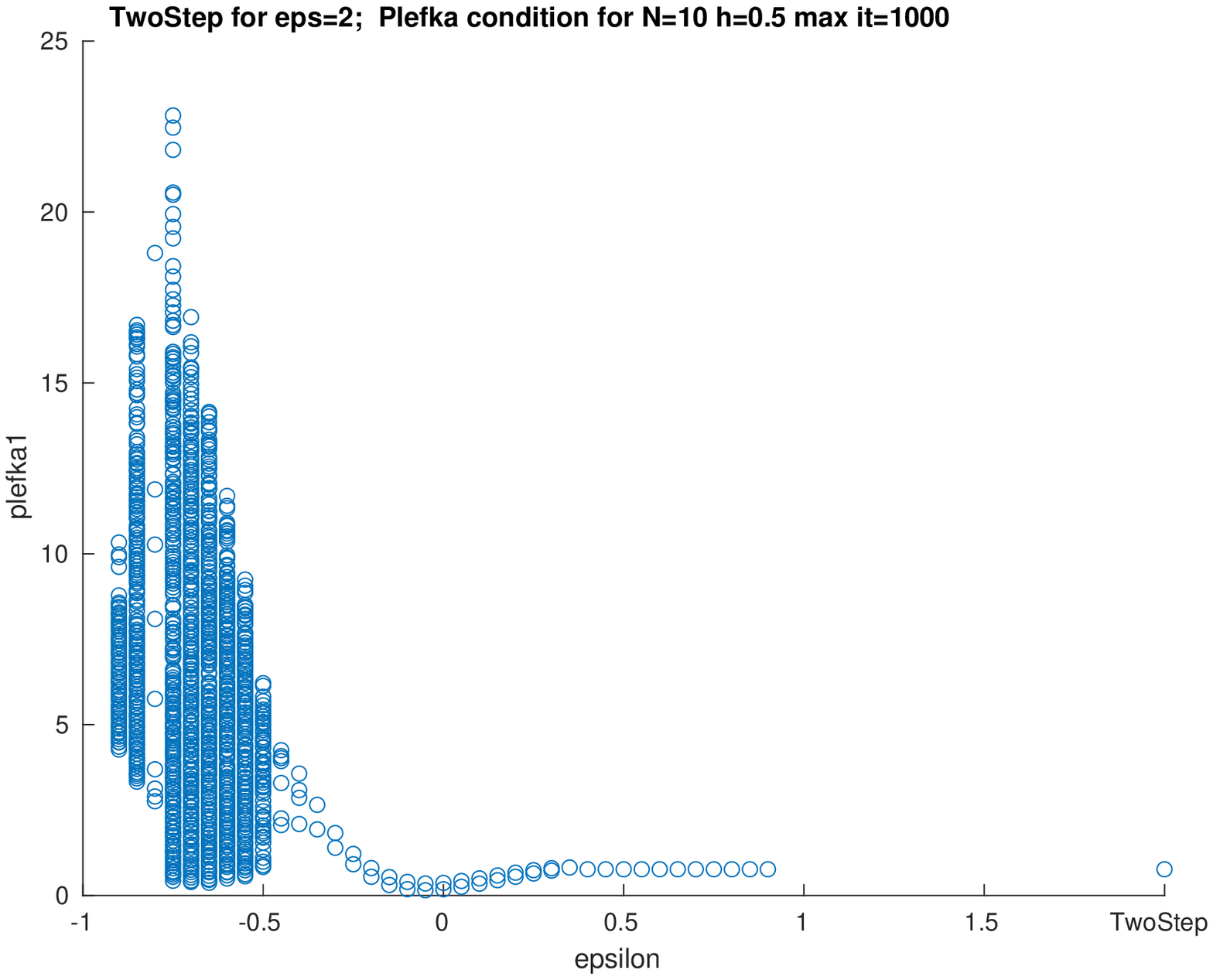}
  \caption{{\tiny} \space  1000 iterations of the $\varepsilon$-Banach algorithm for $N=10, \beta=1$ and $h=0.5$. We simulated 1000 uniformly in the corners of the hypercube distributed start values for $\varepsilon = -0.9, -0.85,\ldots, 0.0$. We plotted all found Plefka values of found fixed points.}
\end{figure}
\begin{figure}[H]
  \centering
  % \psfrag{2}{\tiny Two Step}
  \includegraphics[width=0.5\textwidth,  trim={0.45cm 0.4cm 0 0.4cm},clip]{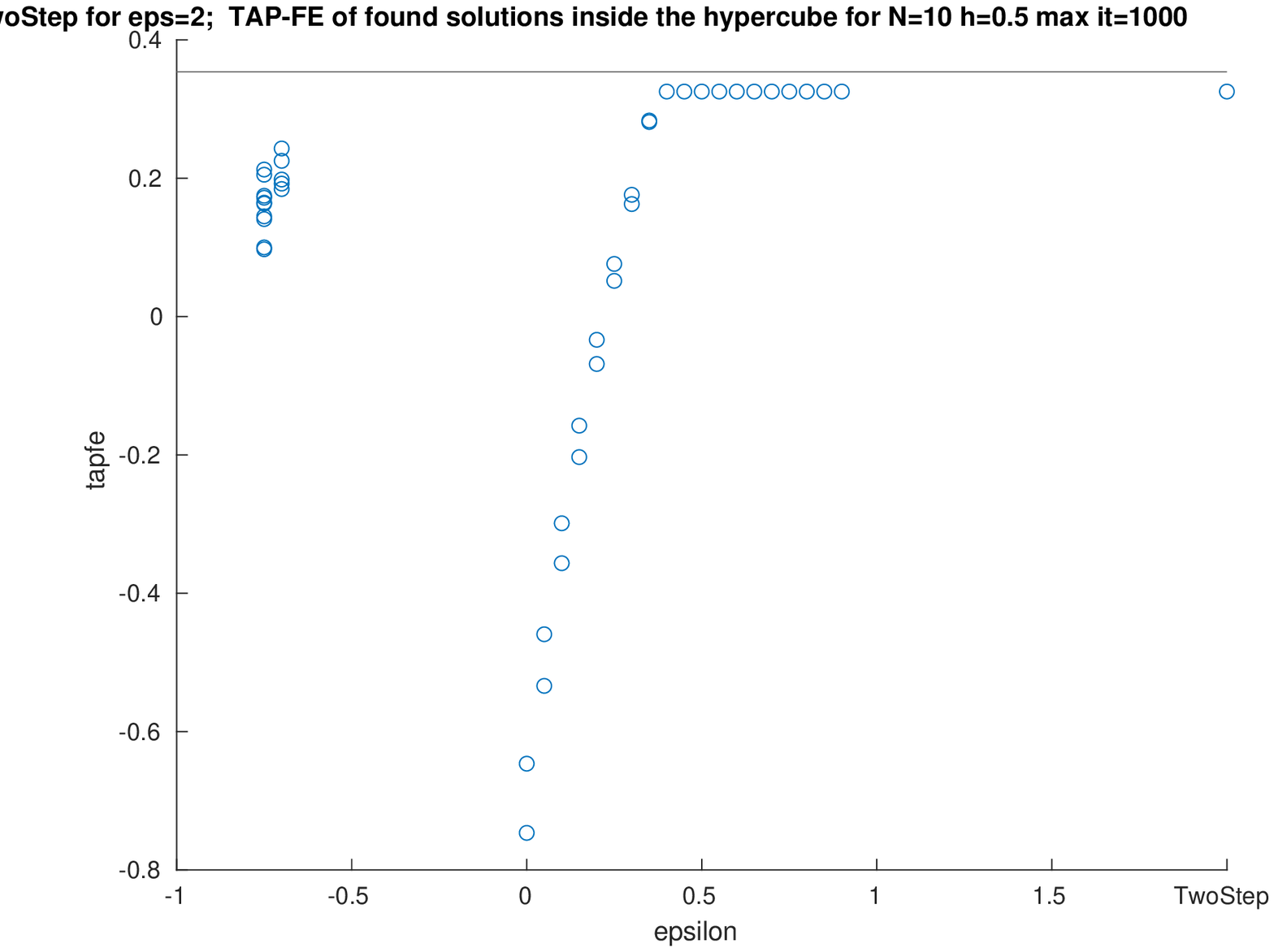}
  \caption{{\tiny} \space  1000 iterations of the $\varepsilon$-Banach algorithm for $N=10, \beta=1$ and $h=0.5$. We simulated 1000 uniformly in the corners of the hypercube distributed start values for $\varepsilon = -0.9, -0.85,\ldots, 0.0$. We plotted all found values of the TAP free energy and excluded all iterations which landed outside the hypercube. The vertical line marks the Parisi free energy value found with the RS formulae.}
\end{figure}

\begin{figure}[H]
  \centering
  \includegraphics[width=0.5\textwidth,  trim={0.45cm 0.4cm 0 0.4cm},clip]{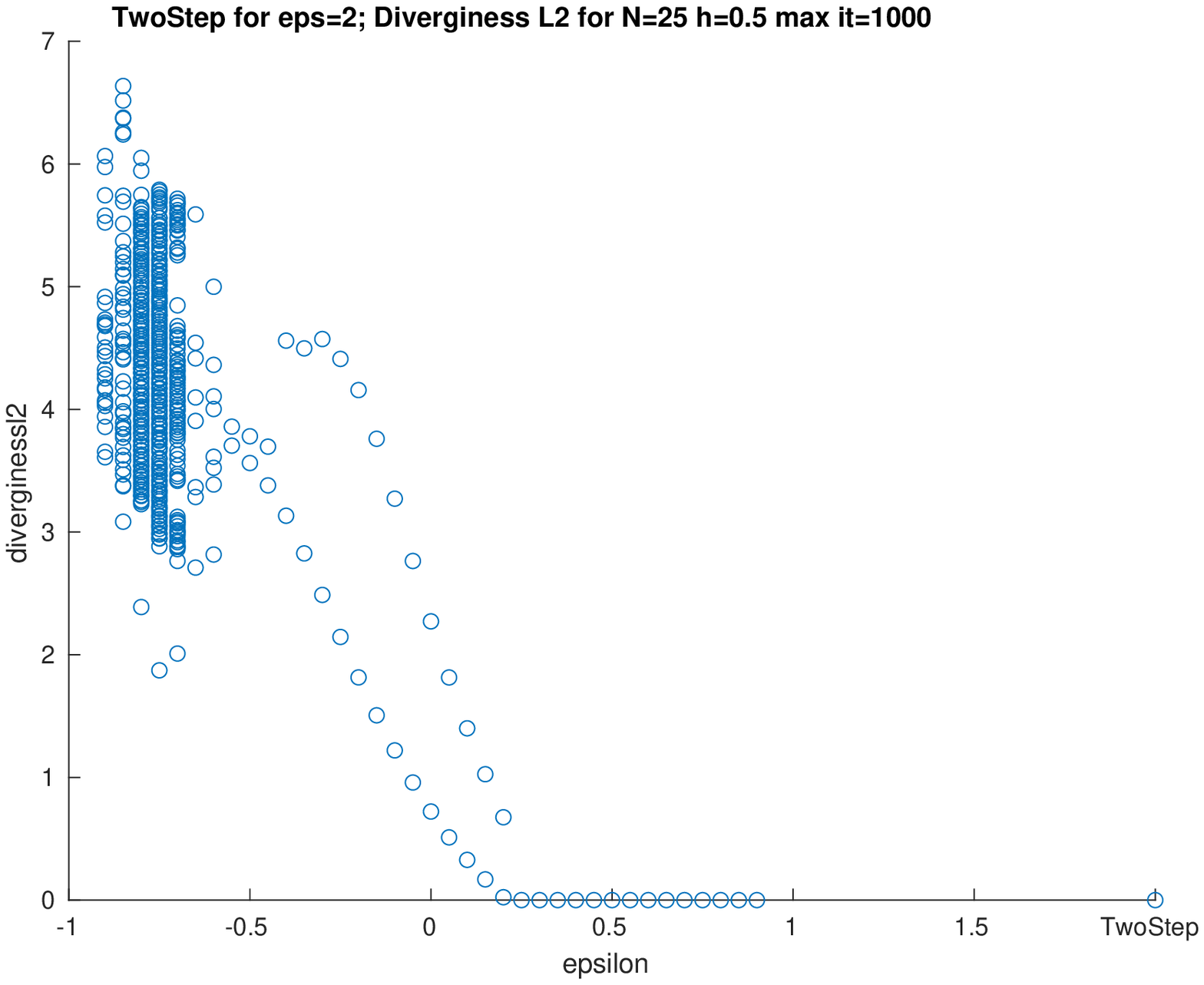}
  \caption{{\tiny} \space  1000 iterations of the $\varepsilon$-Banach algorithm for $N=25, \beta=1$ and $h=0.5$. To measure if the iterations stabilized we use the {\italic mean squared error} between the last two iterates. We simulated 1000 uniformly in the corners of the hypercube distributed start values for $\varepsilon = -0.9, -0.85,\ldots, 0.0$. The {\italic mean squared error} seems to depend highly on the choice of $\varepsilon$. It's obvious that the {\italic mean squared error} will be smaller for $\epsilon\rightarrow 1$ and this does not mean that solutions found with those $\varepsilon$ are better ones. }
\end{figure}
\begin{figure}[H]
  \centering
  \includegraphics[width=0.5\textwidth,  trim={0.45cm 0.4cm 0 0.4cm},clip]{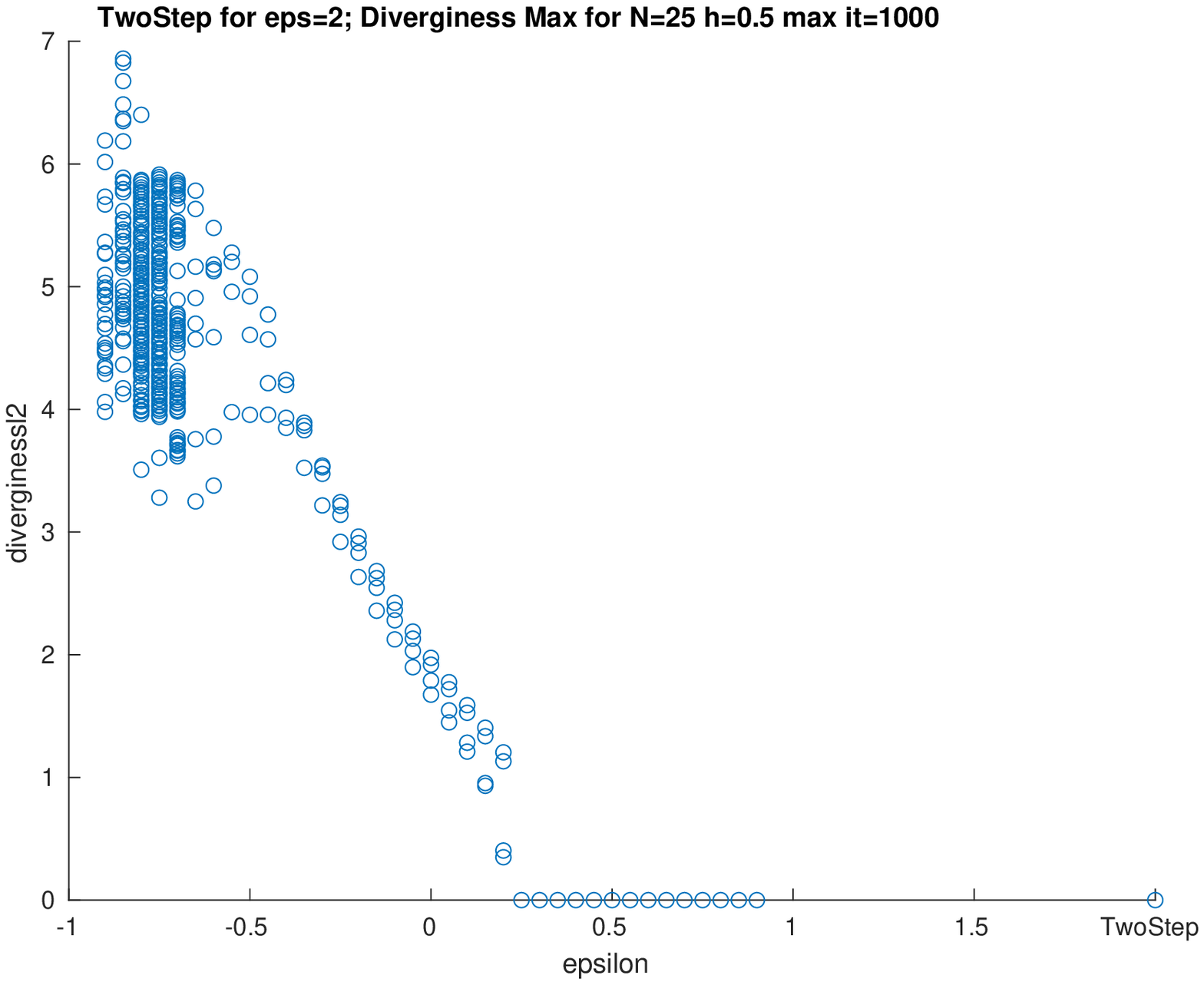}
  \caption{{\tiny} \space  1000 iterations of the $\varepsilon$-Banach algorithm for $N=25, \beta=1$ and $h=0.5$.To measure if the iterations stabilized we use the {\italic maximum absolute error} between the last two iterates. We simulated 1000 uniformly in the corners of the hypercube distributed start values for $\varepsilon = -0.9, -0.85,\ldots, 0.0$. The {\italic maximum absolute error} seems to depend highly on the choice of $\varepsilon$. It's obvious that the {\italic maximum absolute error} will be smaller for $\epsilon\rightarrow 1$ and this does not mean that solutions found with those $\varepsilon$ are better ones. }
\end{figure}
\begin{figure}[H]
  \centering
  % \psfrag{2}{\tiny Two Step}
  \includegraphics[width=0.5\textwidth,  trim={0.45cm 0.4cm 0 0.4cm},clip]{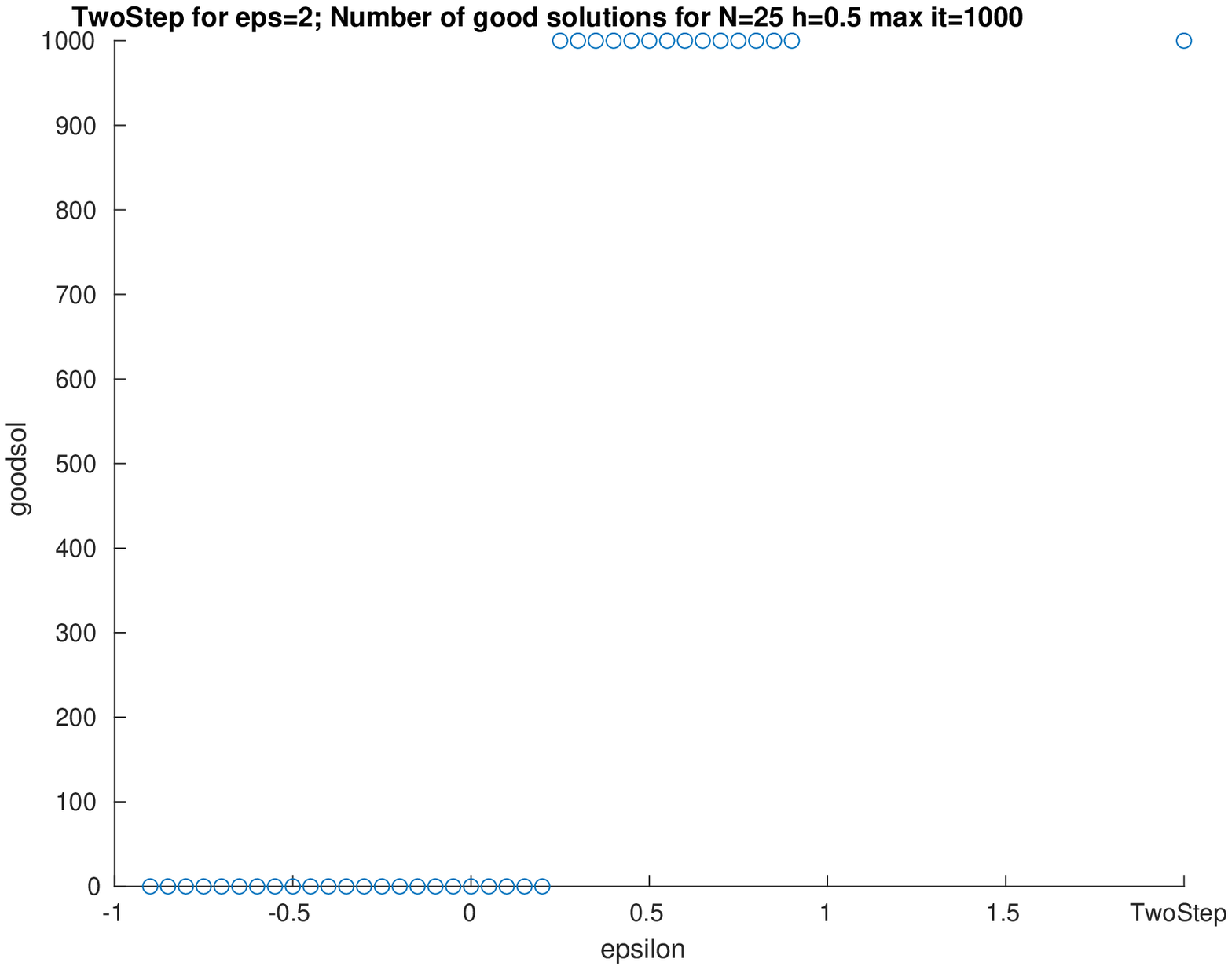}
  \caption{{\tiny} \space  1000 iterations of the $\varepsilon$-Banach algorithm for $N=25, \beta=1$ and $h=0.5$. To measure if the iterations stabilized we use the {\italic mean squared error} between the last two iterates. We simulated 1000 uniformly in the corners of the hypercube distributed start values for $\varepsilon = -0.9, -0.85,\ldots, 0.0$. We say that a found fixed point is a good one if the found TAP free energy is between  $0.35383069483591\pm 0.05$ ($0.35383069483591$ is the RS Parisi FE), the {\italic mean squared error} is smaller than 0.0001 and the Plefka condition is fulfilled.  }
\end{figure}
\begin{figure}[H]
  \centering
  % \psfrag{2}{\tiny Two Step}
  \includegraphics[width=0.5\textwidth,  trim={0.45cm 0.4cm 0 0.4cm},clip]{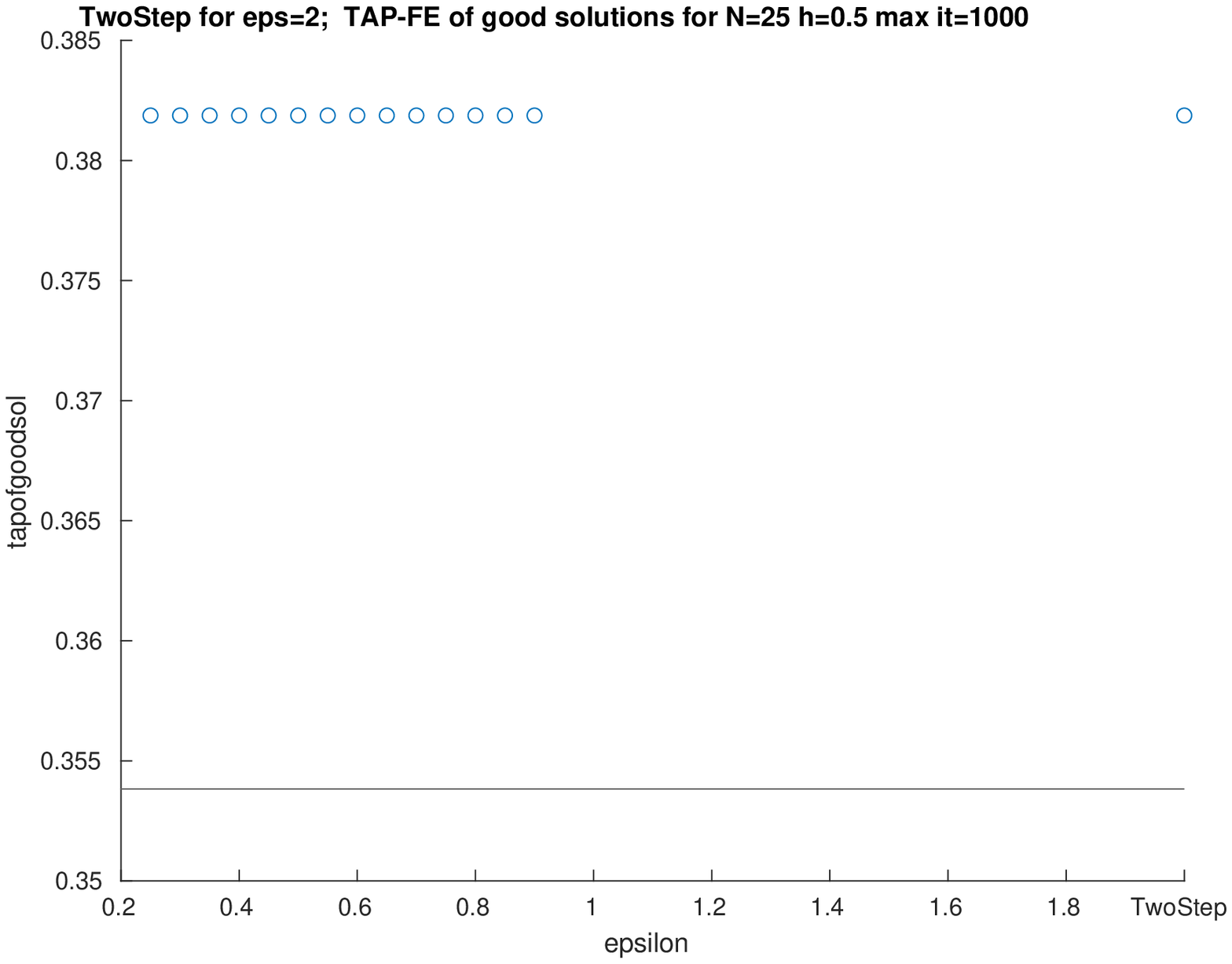}
  \caption{{\tiny} \space  1000 iterations of the $\varepsilon$-Banach algorithm for $N=25, \beta=1$ and $h=0.5$. We simulated 1000 uniformly in the corners of the hypercube distributed start values for $\varepsilon = -0.9, -0.85,\ldots, 0.0$. We say that a found fixed point is a good one if the found TAP free energy is between  $0.35383069483591\pm 0.05$ ($0.35383069483591$ is the RS Parisi FE), the {\italic mean squared error} is smaller than 0.0001 and the Plefka condition is fulfilled. We plotted all found values of the TAP free energy values of good iterations and excluded all iterations which landed outside the hypercube. The vertical line marks the Parisi free energy value found with the RS formulae.}
\end{figure}

\begin{figure}[H]
  \centering
  % \psfrag{2}{\tiny Two Step}
  \includegraphics[width=0.5\textwidth,  trim={0.45cm 0.4cm 0 0.4cm},clip]{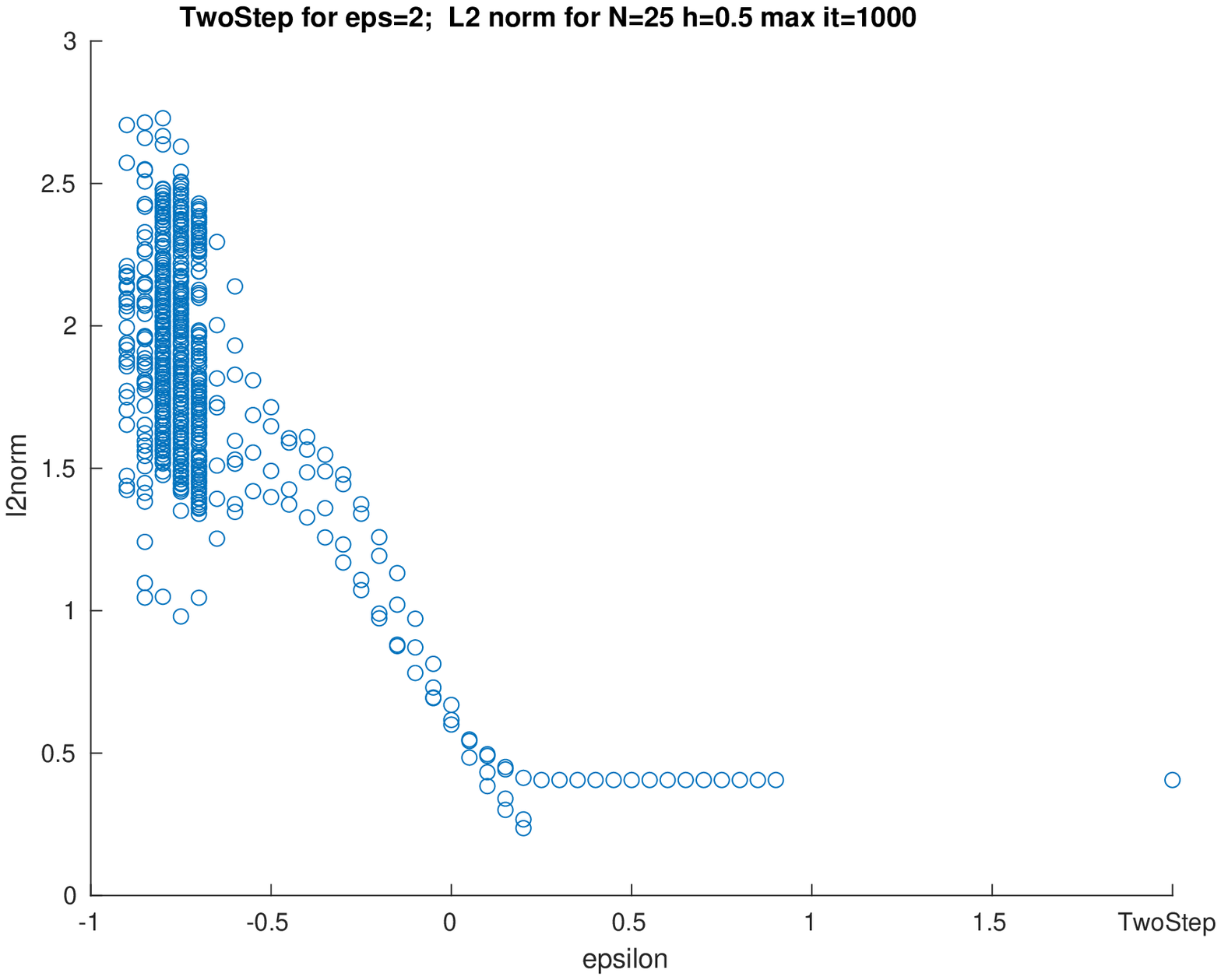}
  \caption{{\tiny} \space  1000 iterations of the $\varepsilon$-Banach algorithm for $N=25, \beta=1$ and $h=0.5$. We simulated 1000 uniformly in the corners of the hypercube distributed start values for $\varepsilon = -0.9, -0.85,\ldots, 0.0$. We plotted all found {\it mean squared errors} between the last two iterates for found fixed points.   }
\end{figure}
\begin{figure}[H]
  \centering
  % \psfrag{2}{\tiny Two Step}
  \includegraphics[width=0.5\textwidth,  trim={0.45cm 0.4cm 0 0.4cm},clip]{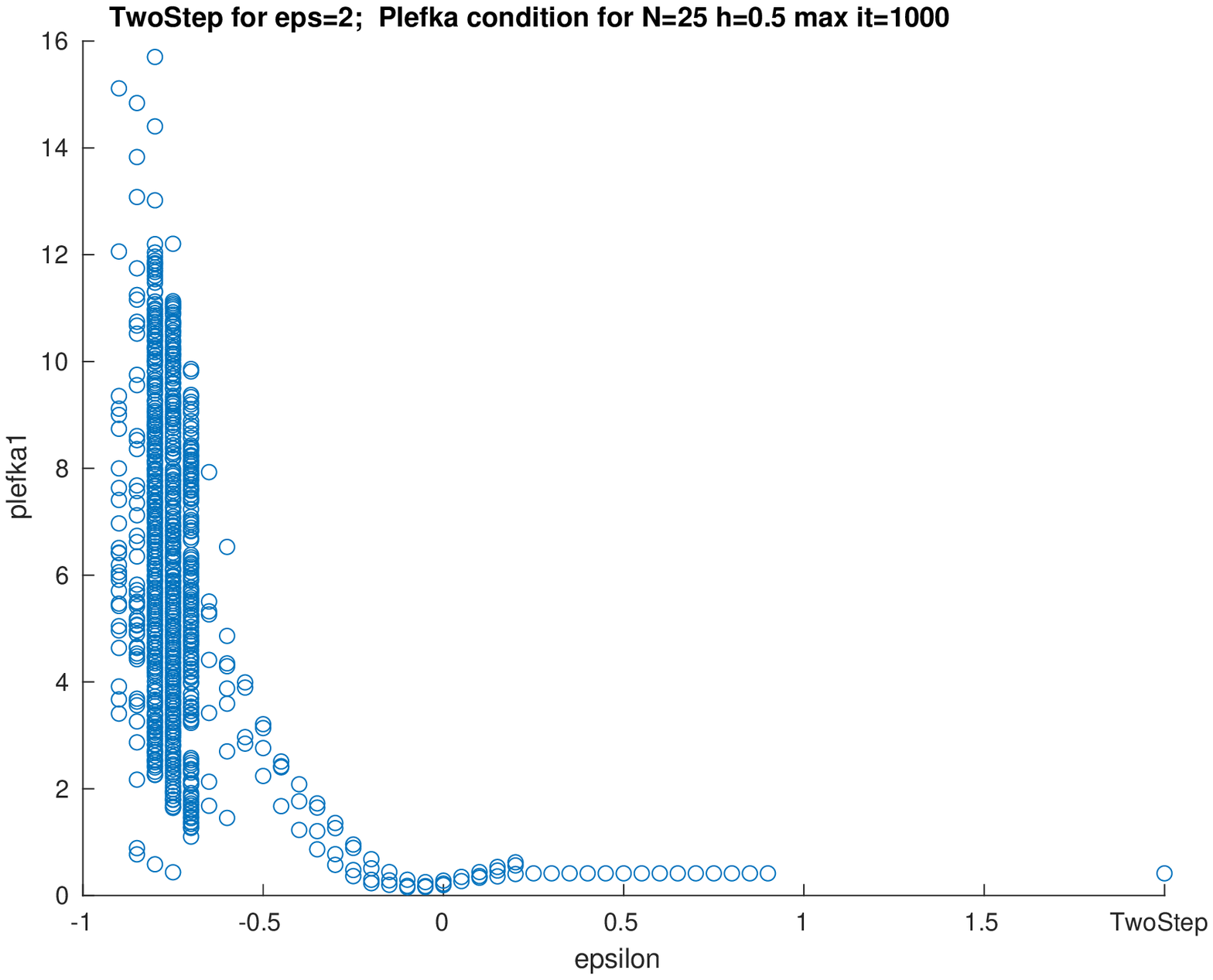}
  \caption{{\tiny} \space  1000 iterations of the $\varepsilon$-Banach algorithm for $N=25, \beta=1$ and $h=0.5$. We simulated 1000 uniformly in the corners of the hypercube distributed start values for $\varepsilon = -0.9, -0.85,\ldots, 0.0$. We plotted all found Plefka values of found fixed points.}
\end{figure}
\begin{figure}[H]
  \centering
  % \psfrag{2}{\tiny Two Step}
  \includegraphics[width=0.5\textwidth,  trim={0 0.4cm 0 0.4cm},clip]{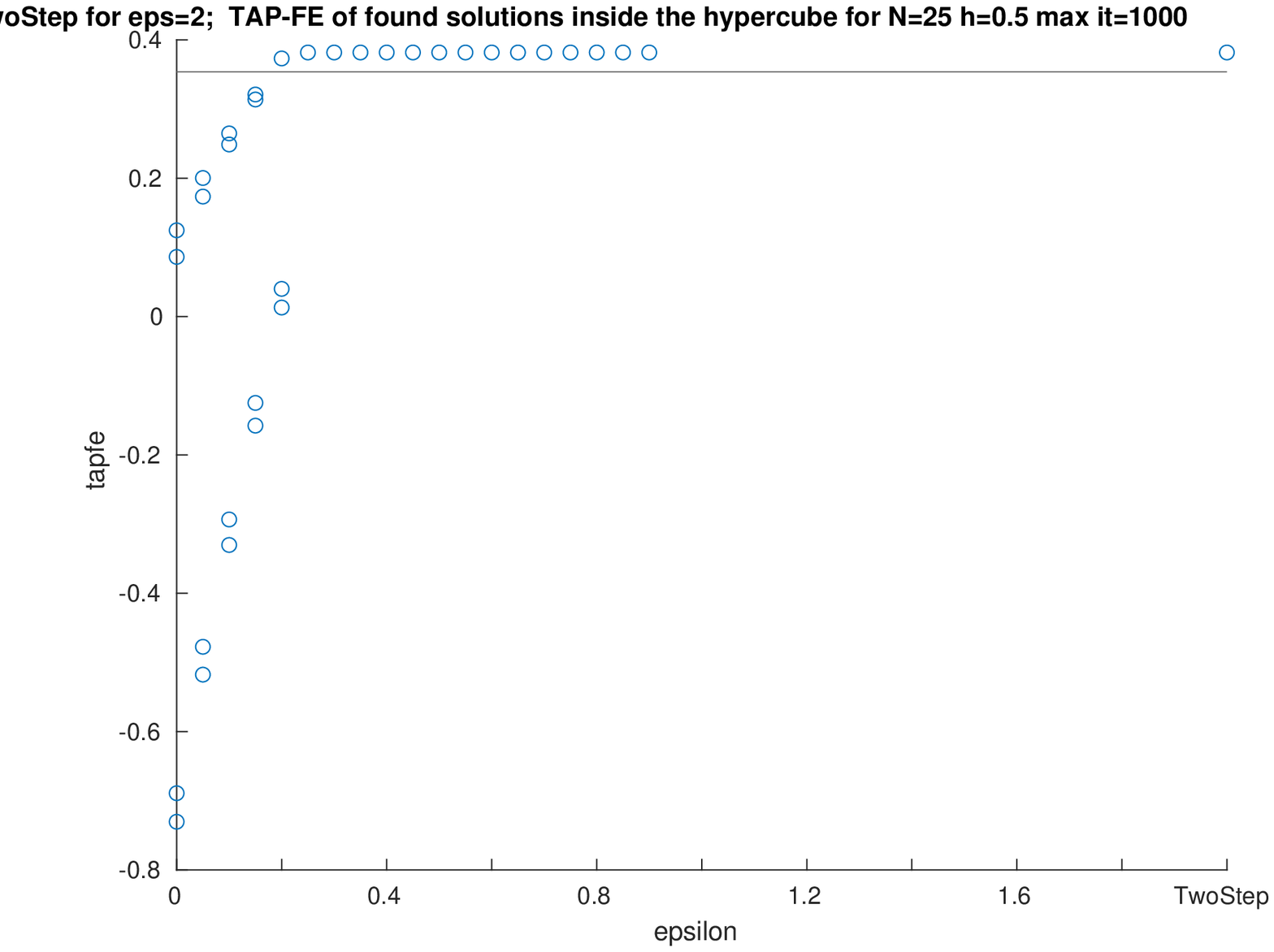}
  \caption{{\tiny} \space  1000 iterations of the $\varepsilon$-Banach algorithm for $N=25, \beta=1$ and $h=0.5$. We simulated 1000 uniformly in the corners of the hypercube distributed start values for $\varepsilon = -0.9, -0.85,\ldots, 0.0$. We plotted all found values of the TAP free energy and excluded all iterations which landed outside the hypercube. The vertical line marks the Parisi free energy value found with the RS formulae.}
\end{figure}

\FloatBarrier\clearpage
\subsection{Uniformly in the hypercube distributed start values}
\leavevmode

\begin{figure}[H]
  \centering
  \includegraphics[width=0.5\textwidth,  trim={0cm 0cm 0 0.4cm},clip]{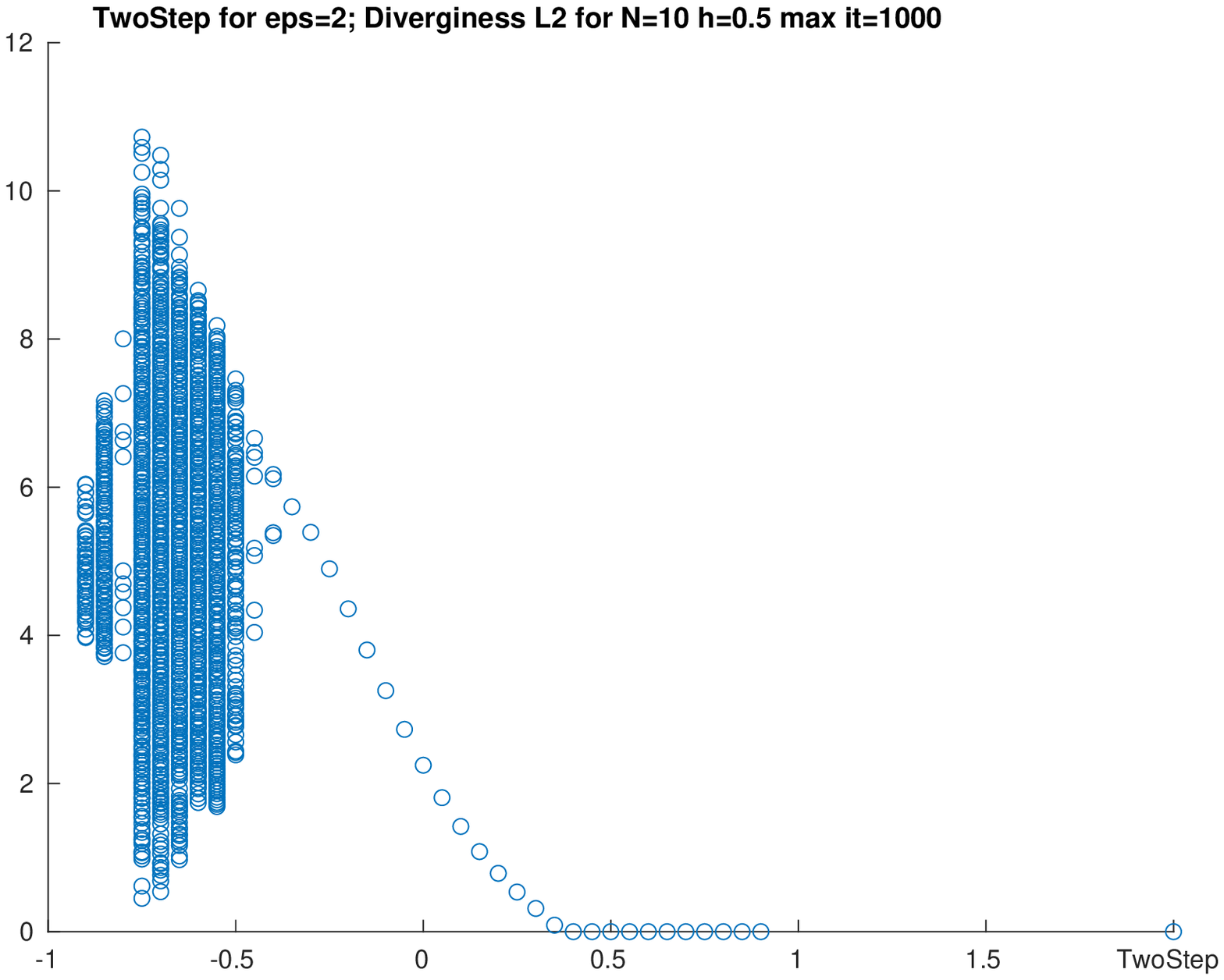}
  \caption{{\tiny} \space  1000 iterations of the $\varepsilon$-Banach algorithm for $N=10, \beta=1$ and $h=0.5$. To measure if the iterations stabilized we use the {\italic mean squared error} between the last two iterates. We simulated 1000 uniformly in the hypercube distributed start values for $\varepsilon = -0.9, -0.85,\ldots, 0.0$. The {\italic mean squared error} seems to depend highly on the choice of $\varepsilon$. It's obvious that the {\italic mean squared error} will be smaller for $\epsilon\rightarrow 1$ and this does not mean that solutions found with those $\varepsilon$ are better ones. }
\end{figure}
\begin{figure}[H]
  \centering
  \includegraphics[width=0.5\textwidth,  trim={0 0 0 0.4cm},clip]{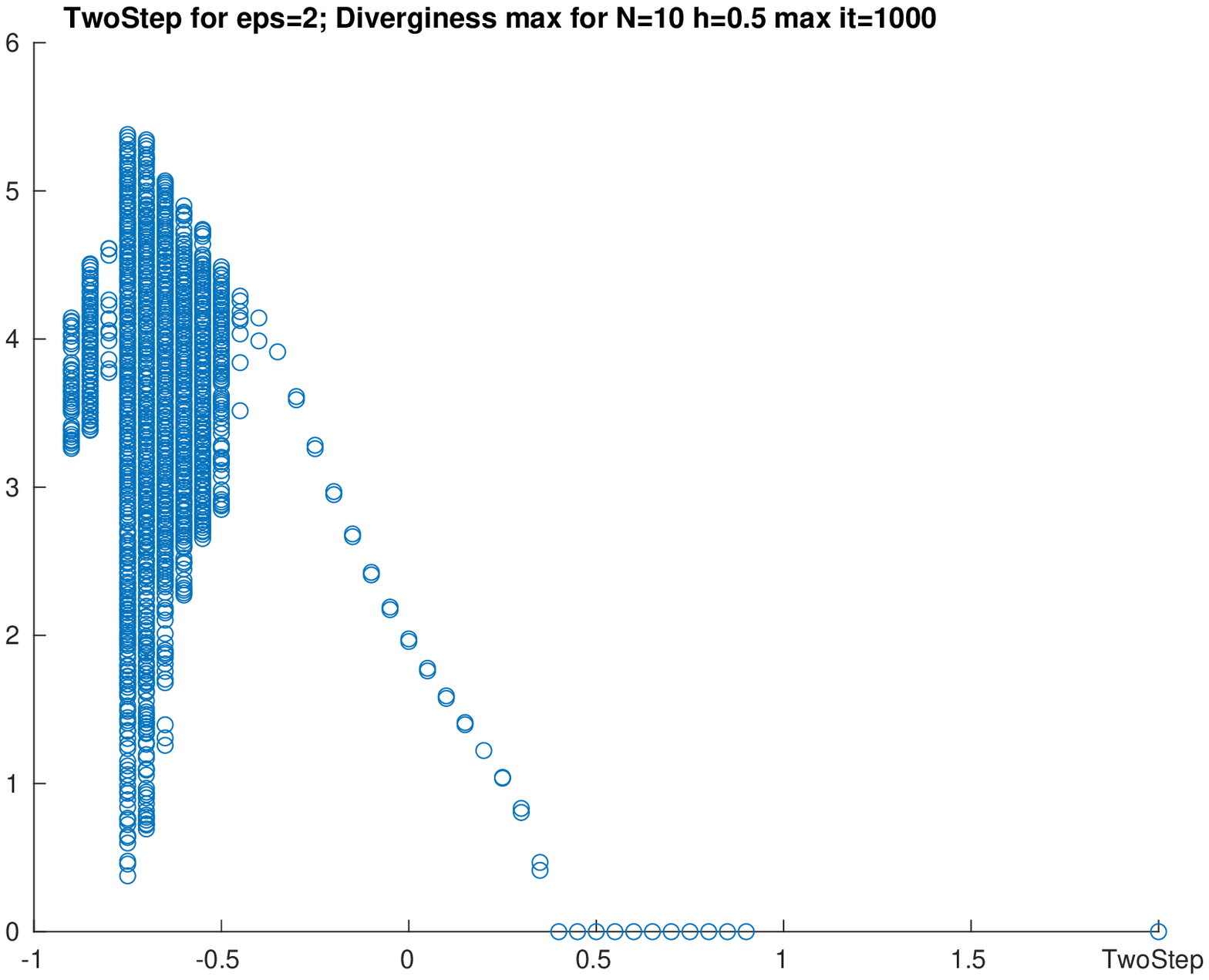}
  \caption{{\tiny} \space  1000 iterations of the $\varepsilon$-Banach algorithm for $N=10, \beta=1$ and $h=0.5$.To measure if the iterations stabilized we use the {\italic maximum absolute error} between the last two iterates. We simulated 1000 uniformly in the hypercube distributed start values for $\varepsilon = -0.9, -0.85,\ldots, 0.0$. The  {\italic maximum absolute error} seems to depend highly on the choice of $\varepsilon$. It's obvious that the {\italic maximum absolute error} will be smaller for $\epsilon\rightarrow 1$ and this does not mean that solutions found with those $\varepsilon$ are better ones. }
\end{figure}
\begin{figure}[H]
  \centering
  % \psfrag{2}{\tiny Two Step}
  \includegraphics[width=0.5\textwidth,  trim={0 0 0 0.4cm},clip]{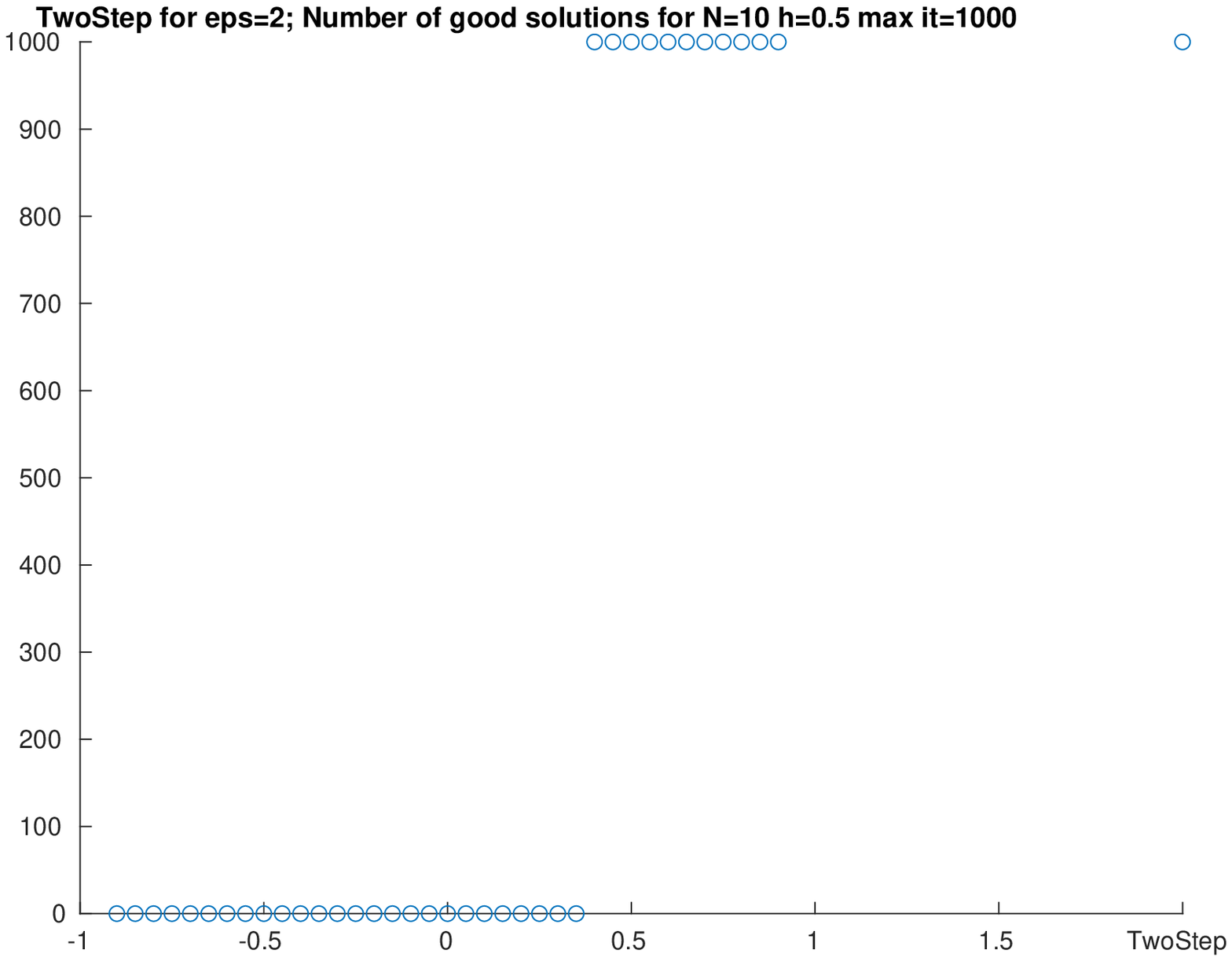}
  \caption{{\tiny} \space  1000 iterations of the $\varepsilon$-Banach algorithm for $N=10, \beta=1$ and $h=0.5$. To measure if the iterations stabilized we use the {\italic mean squared error} between the last two iterates. We simulated 1000 uniformly in the hypercube distributed start values for $\varepsilon = -0.9, -0.85,\ldots, 0.0$. We say that a found fixed point is a good one if the found TAP free energy is between  $0.35383069483591\pm 0.05$ ($0.35383069483591$ is the RS Parisi FE), the {\italic mean squared error} is smaller than 0.0001 and the Plefka condition is fulfilled.  }
\end{figure}
\begin{figure}[H]
  \centering
  % \psfrag{2}{\tiny Two Step}
  \includegraphics[width=0.5\textwidth,  trim={0 0 0 0},clip]{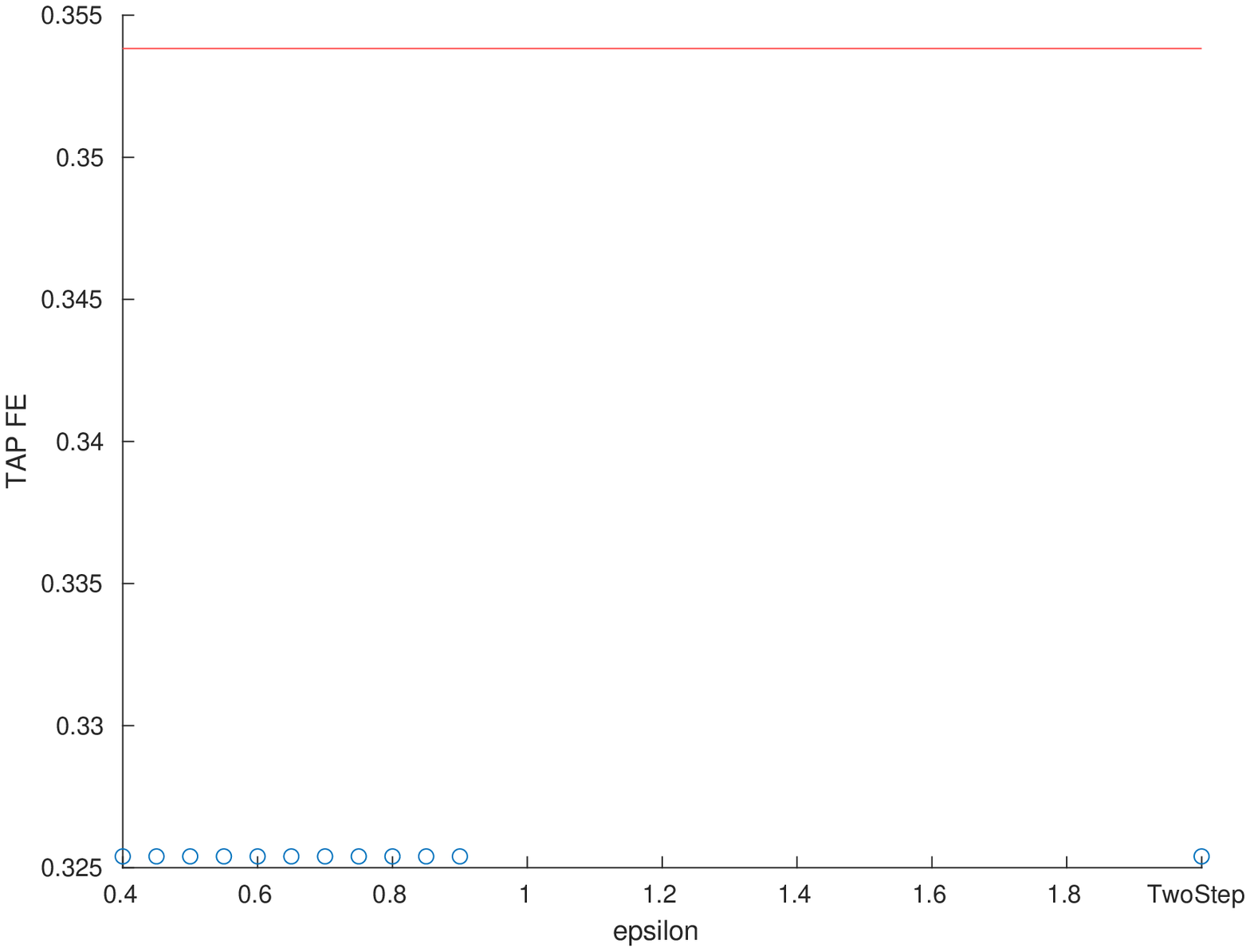}
  \caption{{\tiny} \space  1000 iterations of the $\varepsilon$-Banach algorithm for $N=10, \beta=1$ and $h=0.5$. We simulated 1000 uniformly in the hypercube distributed start values for $\varepsilon = -0.9, -0.85,\ldots, 0.0$. We say that a found fixed point is a good one if the found TAP free energy is between  $0.35383069483591\pm 0.05$ ($0.35383069483591$ is the RS Parisi FE), the {\italic mean squared error} is smaller than 0.0001 and the Plefka condition is fulfilled. We plotted all found values of the TAP free energy values of good iterations and excluded all iterations which landed outside the hypercube. The vertical line marks the Parisi free energy value found with the RS formulae.}
\end{figure}

\begin{figure}[H]
  \centering
  % \psfrag{2}{\tiny Two Step}
  \includegraphics[width=0.5\textwidth,  trim={0 0 0 0.4cm},clip]{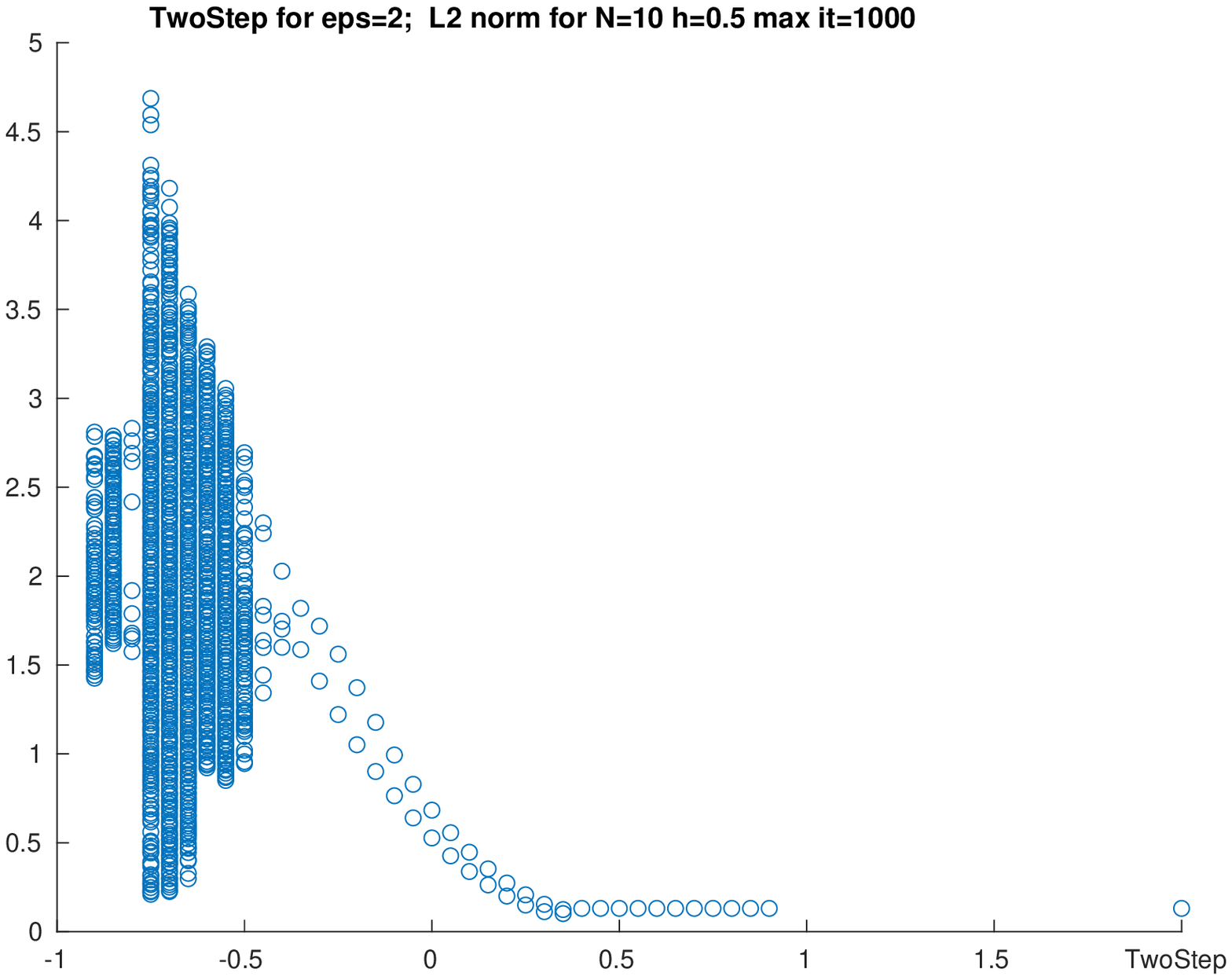}
  \caption{{\tiny} \space  1000 iterations of the $\varepsilon$-Banach algorithm for $N=10, \beta=1$ and $h=0.5$. We simulated 1000 uniformly in the hypercube distributed start values for $\varepsilon = -0.9, -0.85,\ldots, 0.0$. We plotted all found {\it mean squared errors} of found fixed points.   }
\end{figure}
\begin{figure}[H]
  \centering
  % \psfrag{2}{\tiny Two Step}
  \includegraphics[width=0.5\textwidth,  trim={0 0 0 0.4cm},clip]{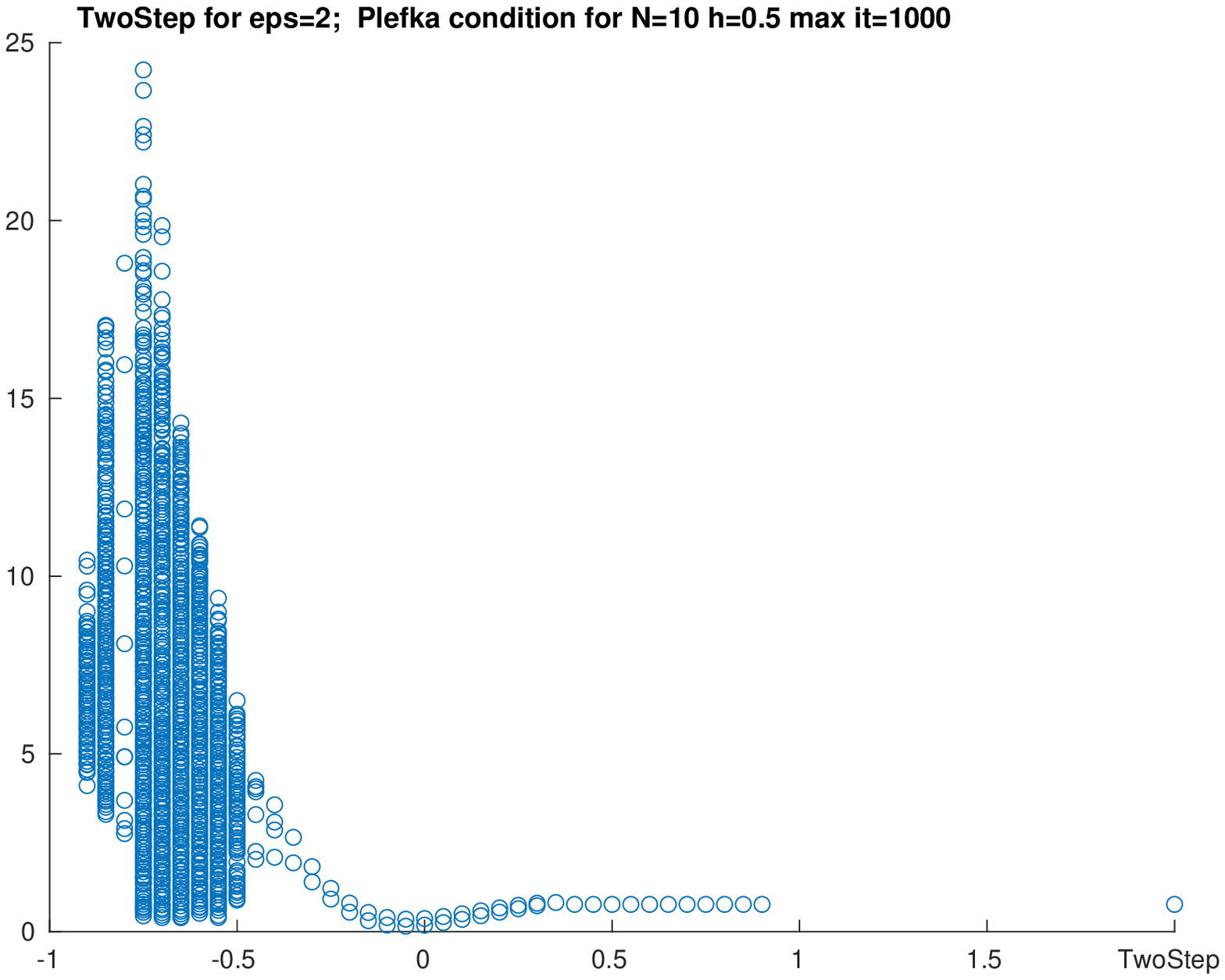}
  \caption{{\tiny} \space  1000 iterations of the $\varepsilon$-Banach algorithm for $N=10, \beta=1$ and $h=0.5$. We simulated 1000 uniformly in the hypercube distributed start values for $\varepsilon = -0.9, -0.85,\ldots, 0.0$. We plotted all found Plefka values of found fixed points.}
\end{figure}
\begin{figure}[H]
  \centering
  % \psfrag{2}{\tiny Two Step}
  \includegraphics[width=0.5\textwidth,  trim={0.45cm 0.4cm 0 0},clip]{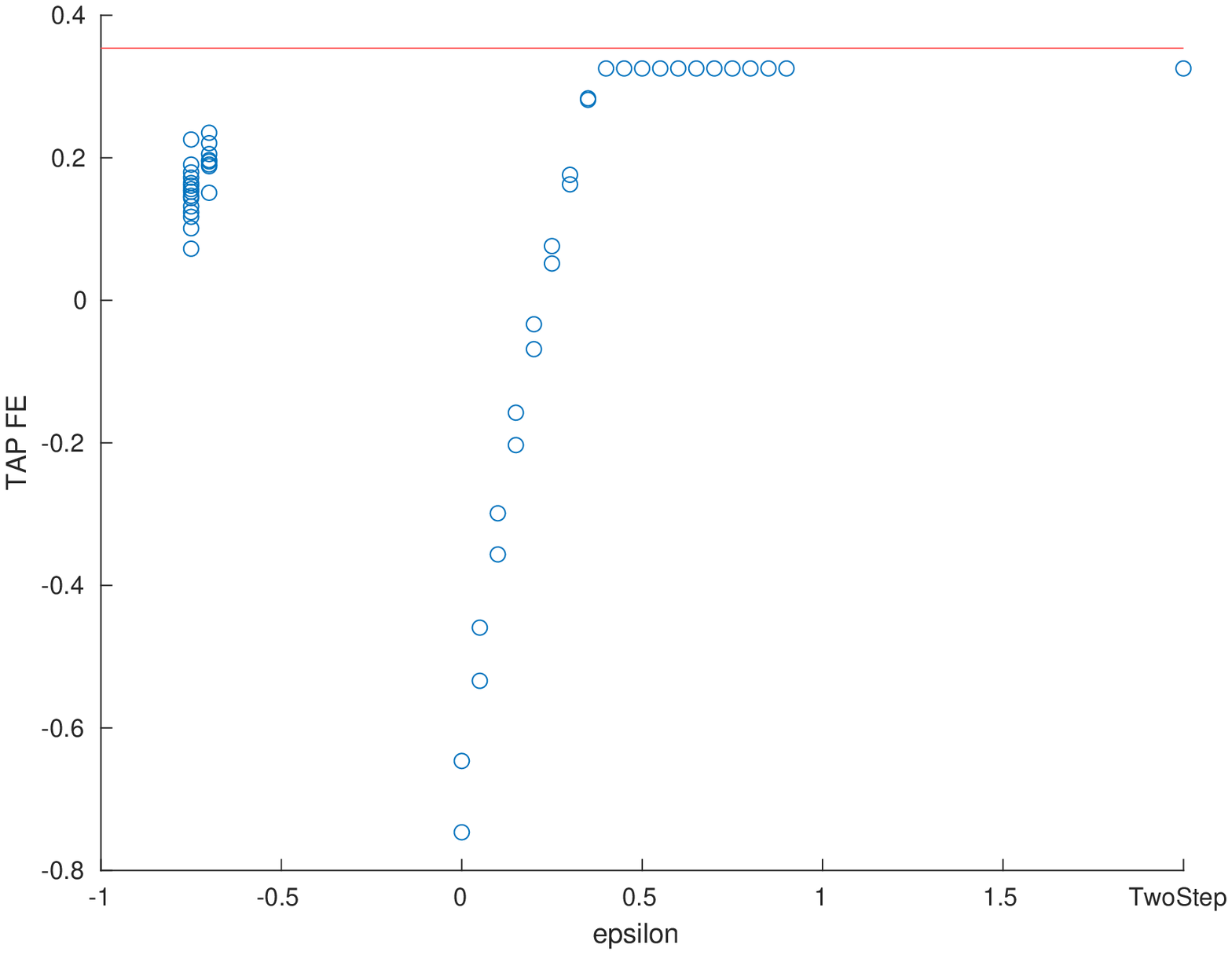}
  \caption{{\tiny} \space  1000 iterations of the $\varepsilon$-Banach algorithm for $N=10, \beta=1$ and $h=0.5$. We simulated 1000 uniformly in the hypercube distributed start values for $\varepsilon = -0.9, -0.85,\ldots, 0.0$. We plotted all found values of the TAP free energy and excluded all iterations which landed outside the hypercube. The vertical line marks the Parisi free energy value found with the RS formulae.}
\end{figure}

\begin{figure}[H]
  \centering
  \includegraphics[width=0.5\textwidth,  trim={0 0 0 0.4cm},clip]{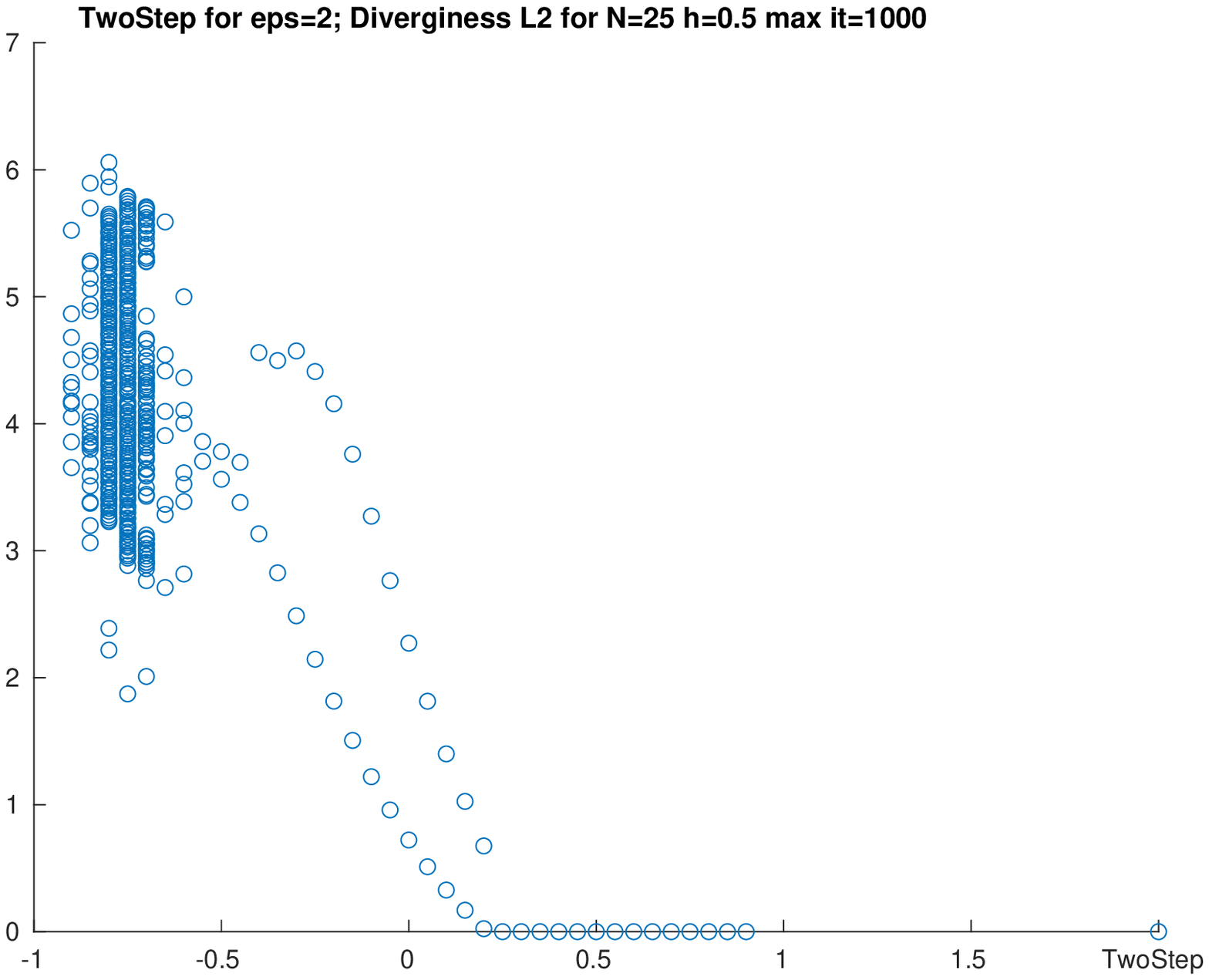}
  \caption{{\tiny} \space  1000 iterations of the $\varepsilon$-Banach algorithm for $N=25, \beta=1$ and $h=0.5$. To measure if the iterations stabilized we use the {\italic mean squared error} between the last two iterates. We simulated 1000 uniformly in the hypercube distributed start values for $\varepsilon = -0.9, -0.85,\ldots, 0.0$. The {\italic mean squared error} seems to depend highly on the choice of $\varepsilon$. It's obvious that the {\italic mean squared error} will be smaller for $\epsilon\rightarrow 1$ and this does not mean that solutions found with those $\varepsilon$ are better ones. }
\end{figure}
\begin{figure}[H]
  \centering
  \includegraphics[width=0.5\textwidth,  trim={0 0 0 0.4cm},clip]{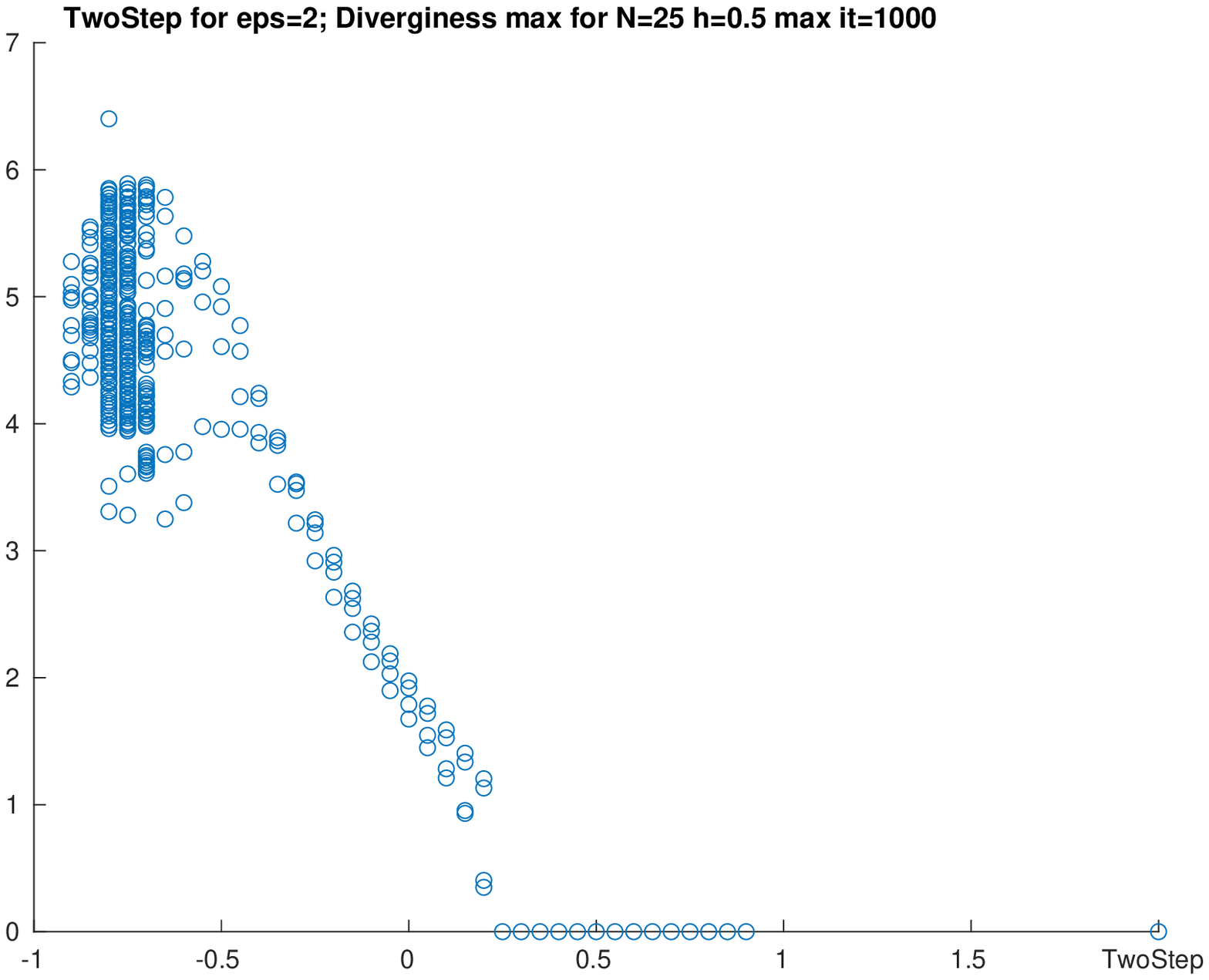}
  \caption{{\tiny} \space  1000 iterations of the $\varepsilon$-Banach algorithm for $N=25, \beta=1$ and $h=0.5$. To measure if the iterations stabilized we use the {\italic maximum absolute error} between the last two iterates. We simulated 1000 uniformly in the hypercube distributed start values for $\varepsilon = -0.9, -0.85,\ldots, 0.0$. The {\italic maximum absolute error} seems to depend highly on the choice of $\varepsilon$. It's obvious that the {\italic maximum absolute error} will be smaller for $\epsilon\rightarrow 1$ and this does not mean that solutions found with those $\varepsilon$ are better ones. }
\end{figure}
\begin{figure}[H]
  \centering
  % \psfrag{2}{\tiny Two Step}
  \includegraphics[width=0.5\textwidth,  trim={0 0 0 0.4cm},clip]{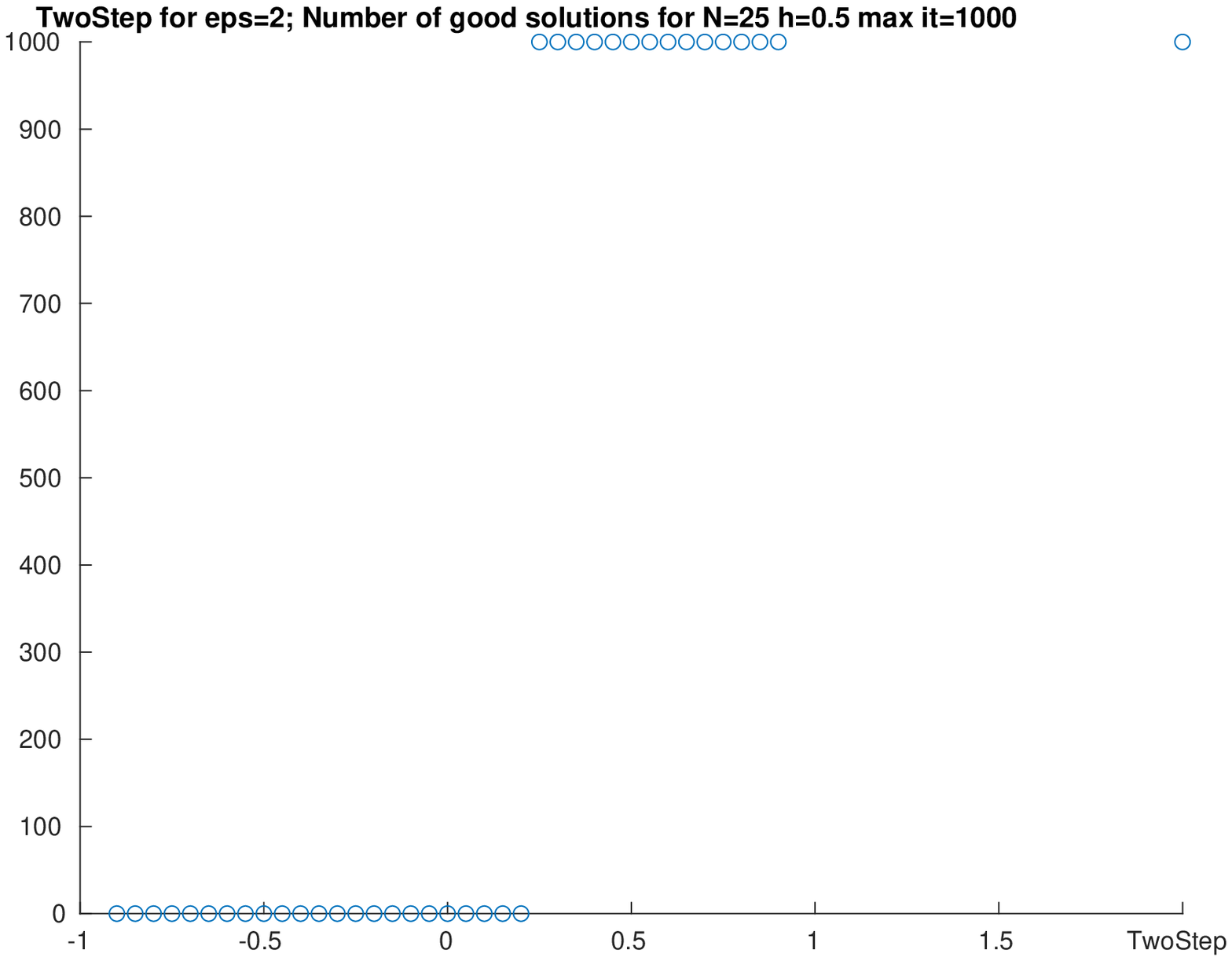}
  \caption{{\tiny} \space  1000 iterations of the $\varepsilon$-Banach algorithm for $N=25, \beta=1$ and $h=0.5$. To measure if the iterations stabilized we use the {\italic mean squared error} between the last two iterates. We simulated 1000 uniformly in the hypercube distributed start values for $\varepsilon = -0.9, -0.85,\ldots, 0.0$. We say that a found fixed point is a good one if the found TAP free energy is between  $0.35383069483591\pm 0.05$ ($0.35383069483591$ is the RS Parisi FE), the {\italic mean squared error} is smaller than 0.0001 and the Plefka condition is fulfilled.  }
\end{figure}
\begin{figure}[H]
  \centering
  % \psfrag{2}{\tiny Two Step}
  \includegraphics[width=0.5\textwidth,  trim={0.45cm 0.4cm 0 0},clip]{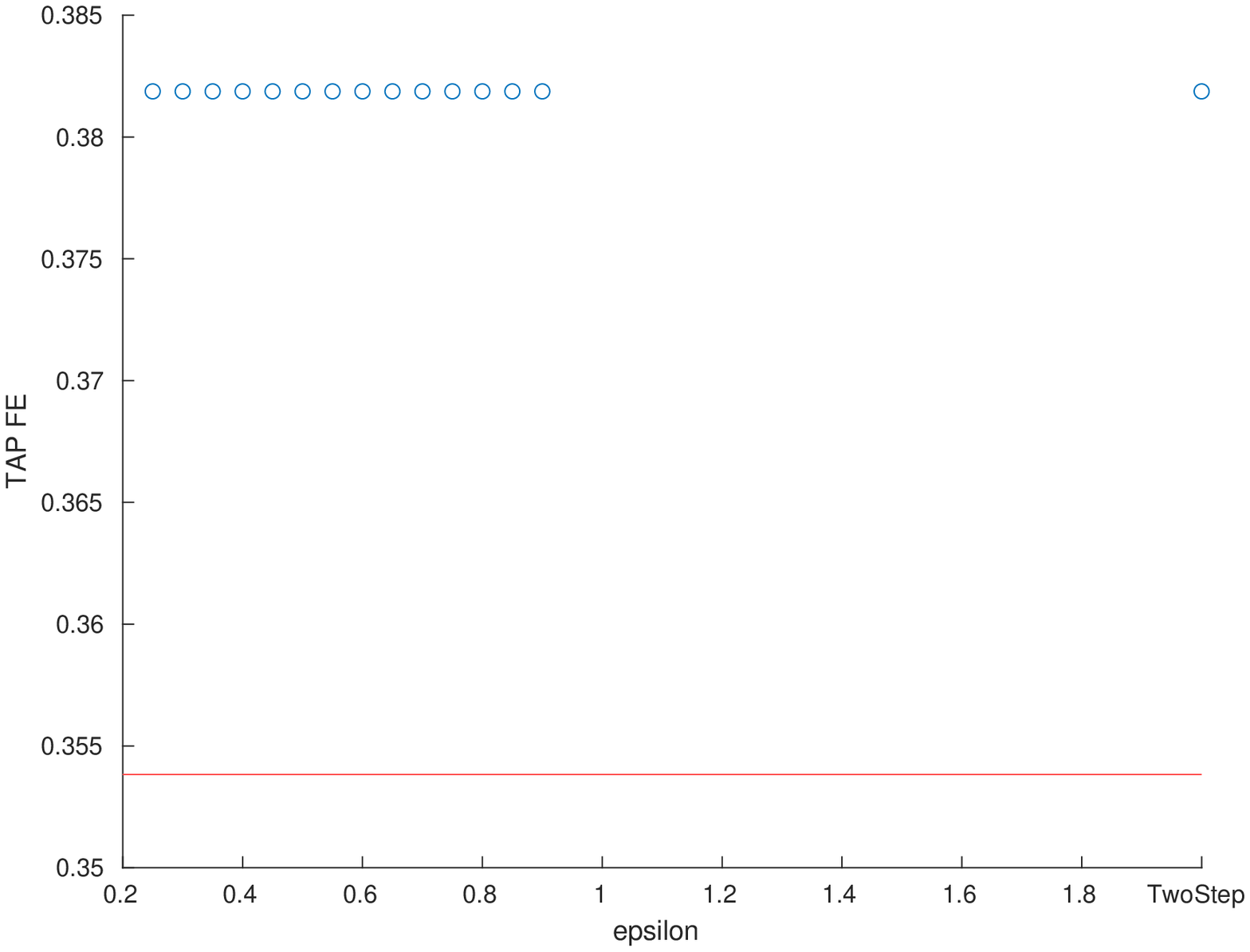}
  \caption{{\tiny} \space  1000 iterations of the $\varepsilon$-Banach algorithm for $N=25, \beta=1$ and $h=0.5$. We simulated 1000 uniformly in the hypercube distributed start values for $\varepsilon = -0.9, -0.85,\ldots, 0.0$. We say that a found fixed point is a good one if the found TAP free energy is between  $0.35383069483591\pm 0.05$ ($0.35383069483591$ is the RS Parisi FE), the  {\italic mean squared error} is smaller than 0.0001 and the Plefka condition is fulfilled. We plotted all found values of the TAP free energy values of good iterations and excluded all iterations which landed outside the hypercube. The vertical line marks the Parisi free energy value found with the RS formulae.}
\end{figure}

\begin{figure}[H]
  \centering
  % \psfrag{2}{\tiny Two Step}
  \includegraphics[width=0.5\textwidth,  trim={0 0 0 0.4cm},clip]{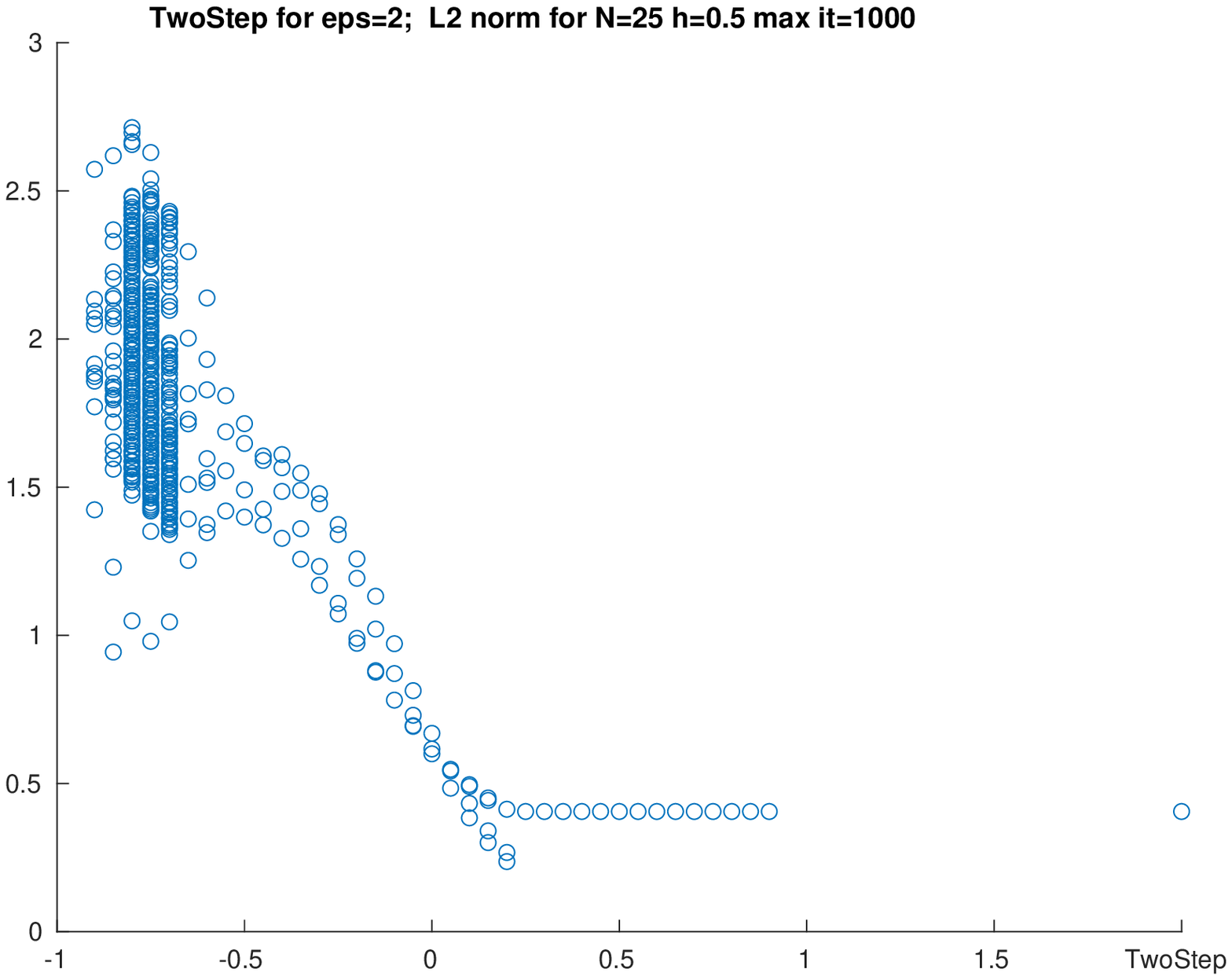}
  \caption{{\tiny} \space  1000 iterations of the $\varepsilon$-Banach algorithm for $N=25, \beta=1$ and $h=0.5$. We simulated 1000 uniformly in the hypercube distributed start values for $\varepsilon = -0.9, -0.85,\ldots, 0.0$. We plotted all found {\it mean squared errors} between the last two iterates for found fixed points.   }
\end{figure}
\begin{figure}[H]
  \centering
  % \psfrag{2}{\tiny Two Step}
  \includegraphics[width=0.5\textwidth,  trim={0 0 0 0.4cm},clip]{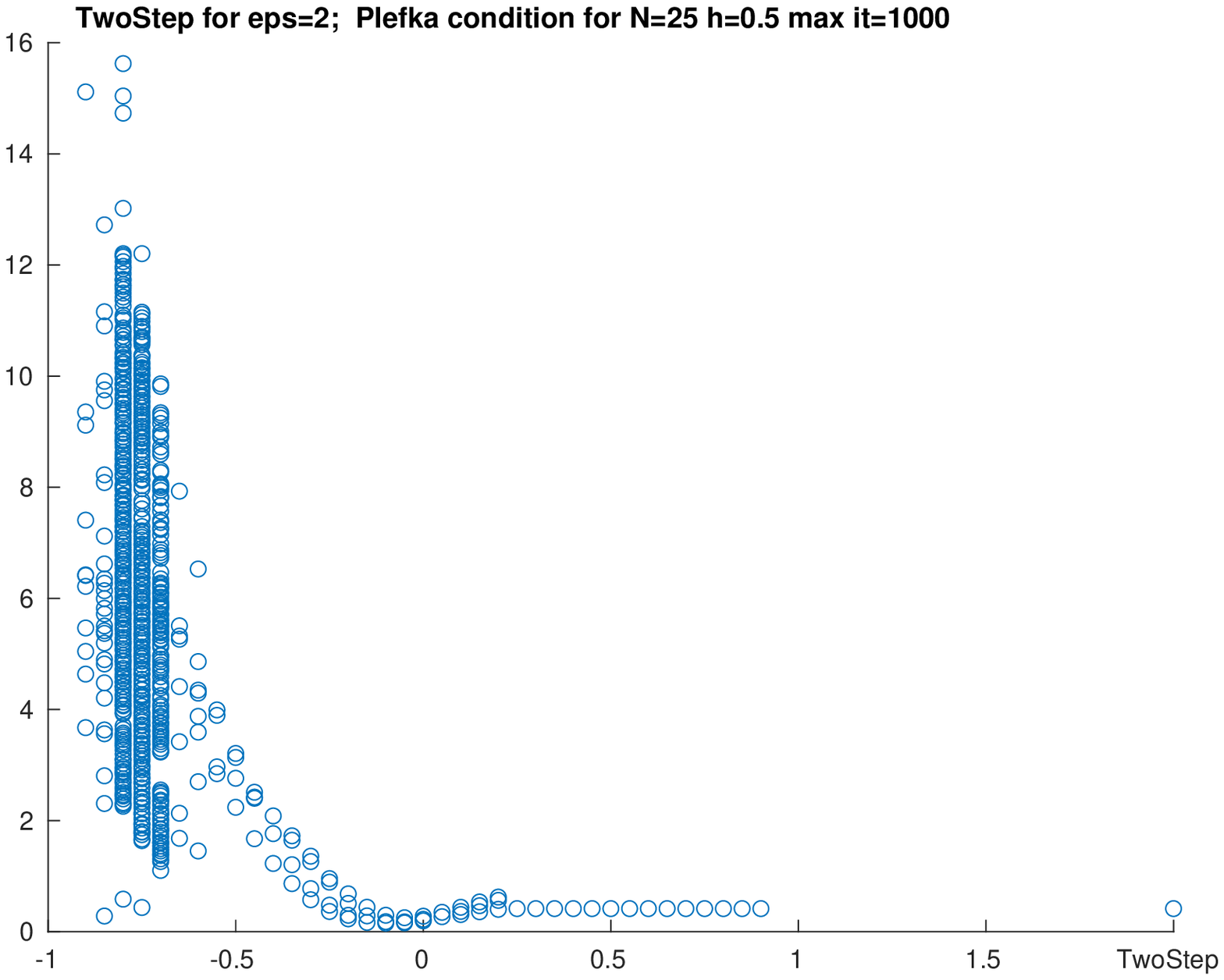}
  \caption{{\tiny} \space  1000 iterations of the $\varepsilon$-Banach algorithm for $N=25, \beta=1$ and $h=0.5$. We simulated 1000 uniformly in the hypercube distributed start values for $\varepsilon = -0.9, -0.85,\ldots, 0.0$. We plotted all found Plefka values of found fixed points.}
\end{figure}
\begin{figure}[H]
  \centering
  % \psfrag{2}{\tiny Two Step}
  \includegraphics[width=0.5\textwidth,  trim={0.45cm 0.4cm 0 0},clip]{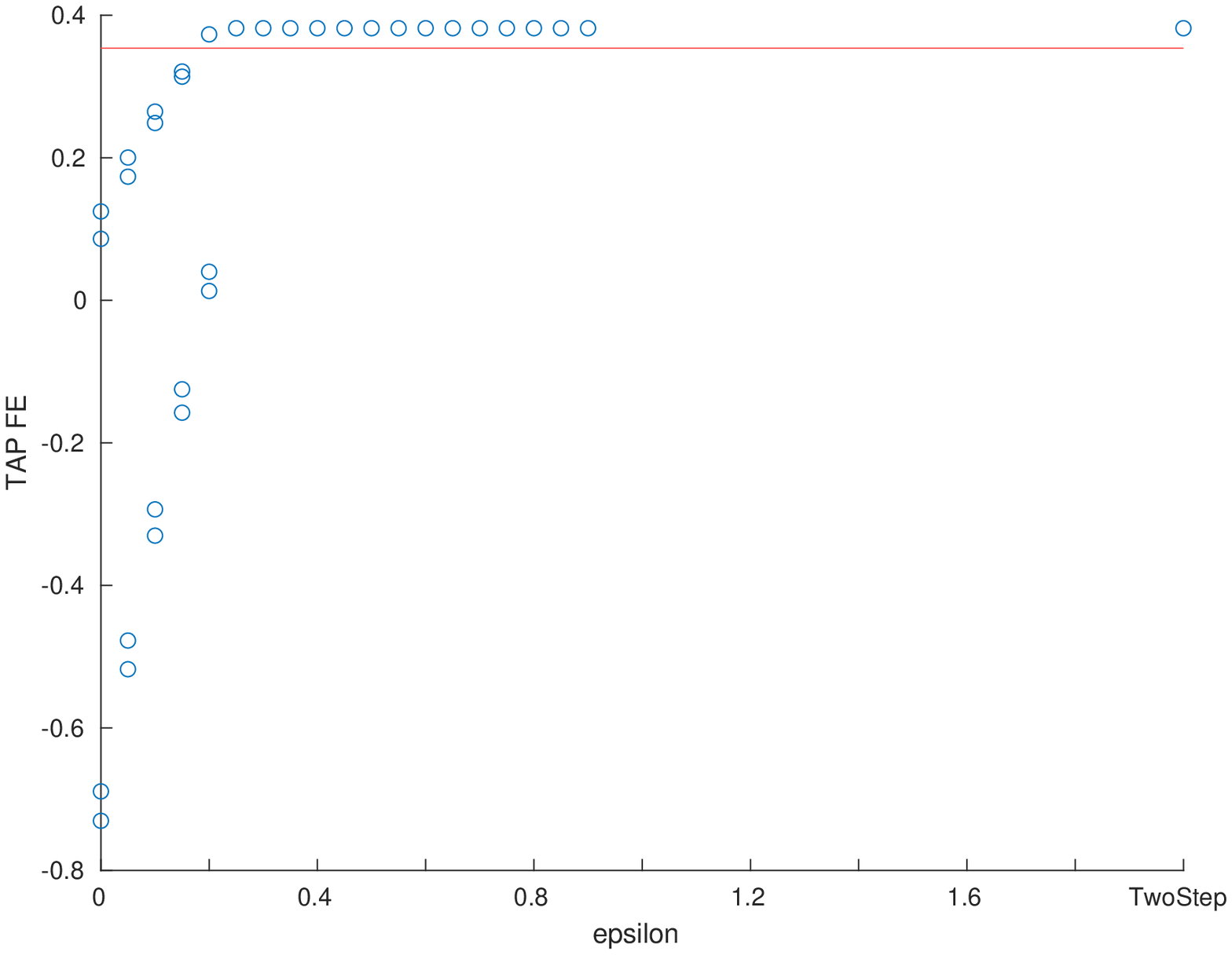}
  \caption{{\tiny} \space  1000 iterations of the $\varepsilon$-Banach algorithm for $N=25, \beta=1$ and $h=0.5$. We simulated 1000 uniformly in the hypercube distributed start values for $\varepsilon = -0.9, -0.85,\ldots, 0.0$. We plotted all found values of the TAP free energy and excluded all iterations which landed outside the hypercube. The vertical line marks the Parisi free energy value found with the RS formula.}
\end{figure}\FloatBarrier\clearpage
\subsection{Short summary}
In the high temperature regime we see that Two-Step-Banach is not sensitive with respect to the start values. $\varepsilon$-Banach is very sensitive with respect to the start values for $\varepsilon$ negative enough. In the high temperature regime there seems to be a positive lower bound on $\varepsilon$ such that we find \emph{good} solutions. It is very interesting that we always find a TAP free energy different from the RS Parisi FE.\FloatBarrier\clearpage
\section{Low temperature regime: $\varepsilon$-Banach and Two-Step}
\subsection{Uniformly in the hypercube distributed start values for  $\varepsilon$-Banach}
\subsubsection{Using the {\it mean squared error}}\FloatBarrier\leavevmode
\begin{figure}
  \centering
  % \psfrag{2}{\tiny Two Step}
  \includegraphics[width=0.5\textwidth,  trim={0 0 0 0},clip]{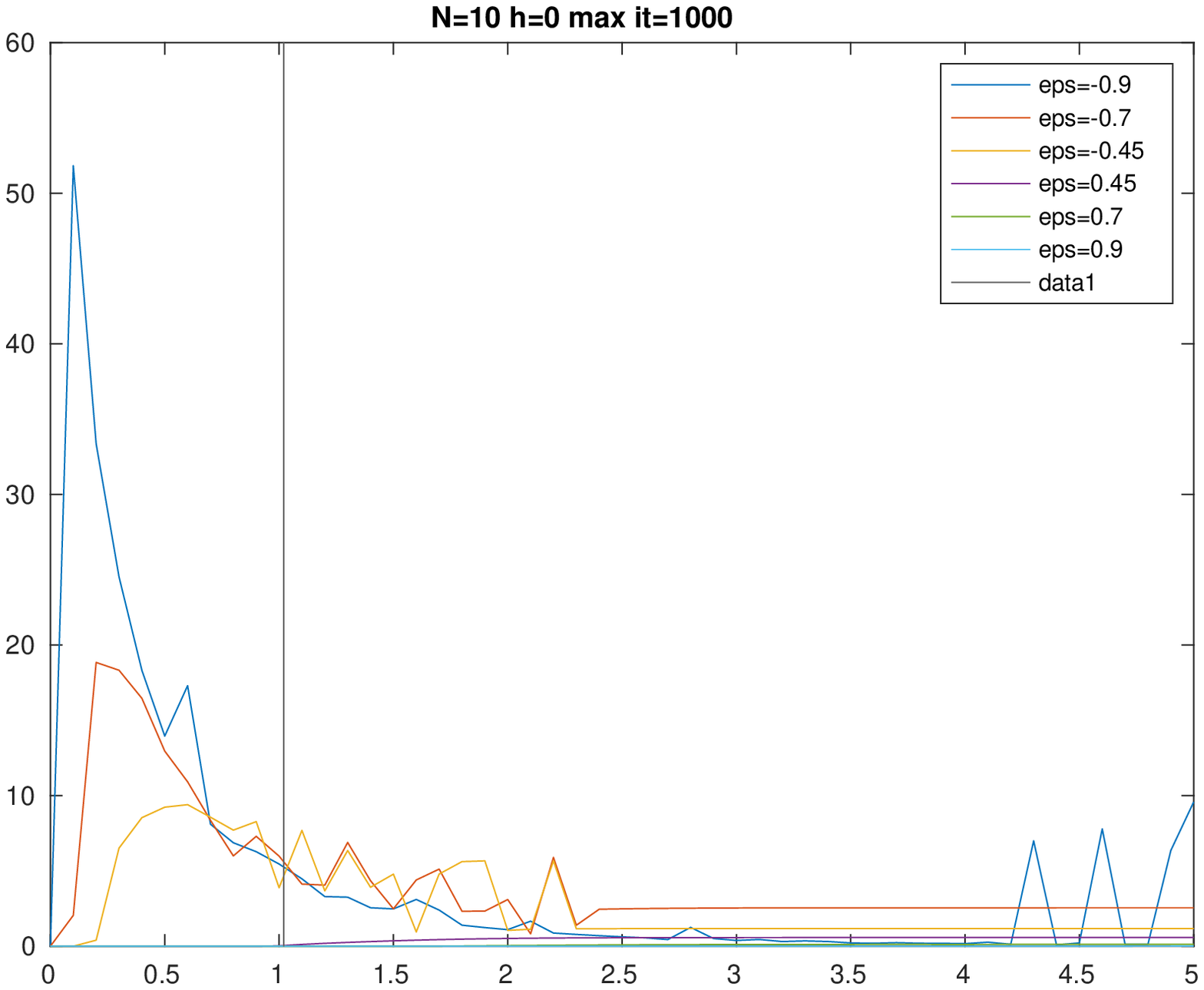}
  \caption{{\tiny} \space  1000 iterations of the $\varepsilon$-Banach algorithm for $N=10, \beta\in [0, 5]$ and $h=0$. We simulated one uniformly in the hypercube distributed start value for $\varepsilon = -0.9, -0.7, -0.45, 0.45, 0.7, 0.9$. We plotted the {\italic mean squared error} for the given parameters and used the same start value for all $\beta$.}
\end{figure}
\begin{figure}[H]
  \centering
  % \psfrag{2}{\tiny Two Step}
  \includegraphics[width=0.5\textwidth,  trim={0 0 0 0},clip]{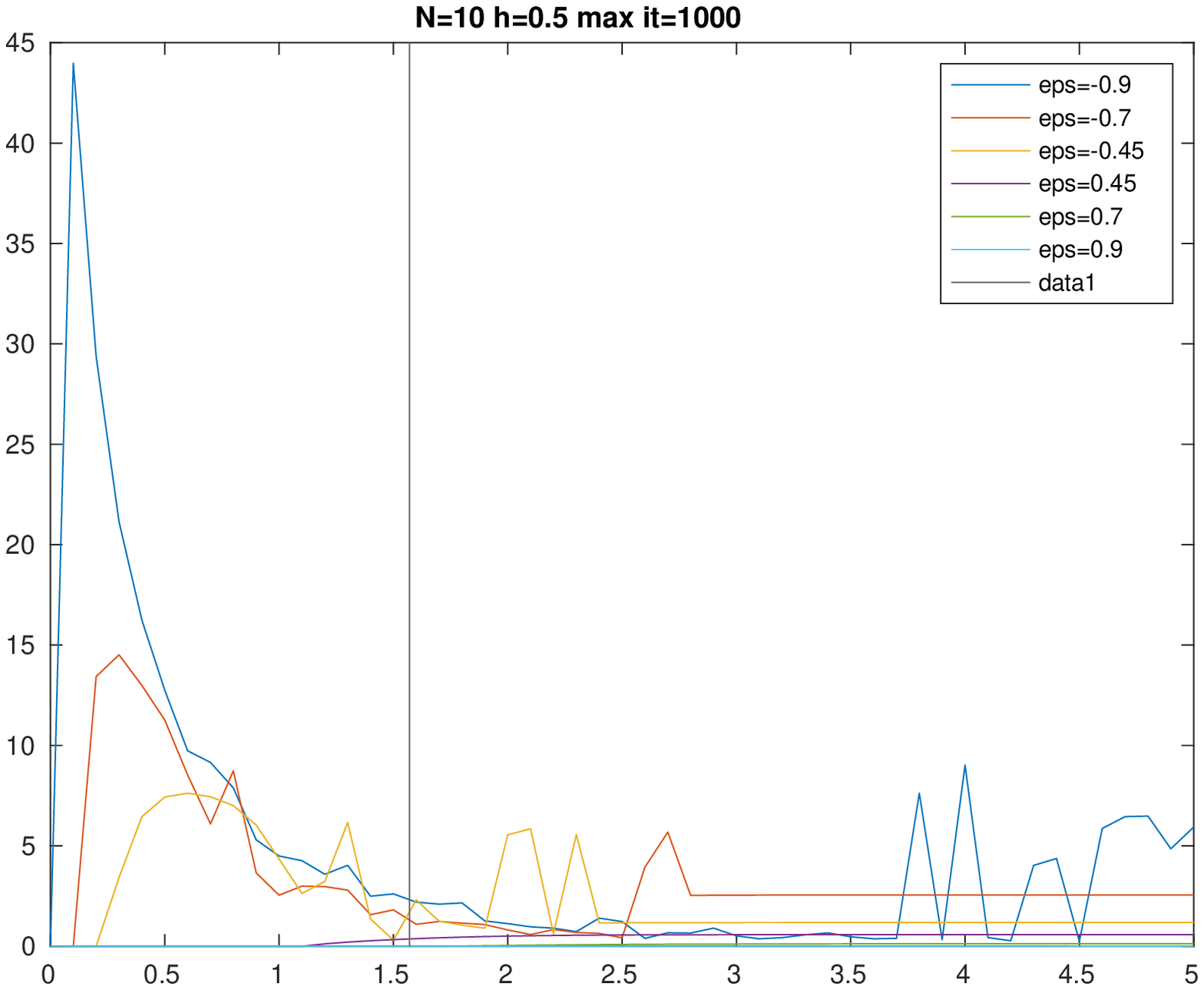}
  \caption{{\tiny} \space  1000 iterations of the $\varepsilon$-Banach algorithm for $N=10, \beta\in [0, 5]$ and $h=0.5$. We simulated one uniformly in the hypercube distributed start value for $\varepsilon = -0.9, -0.7, -0.45, 0.45, 0.7, 0.9$. We plotted the  {\italic mean squared error} for the given parameters and used the same start value for all $\beta$.}
\end{figure}
\begin{figure}[H]
  \centering
  % \psfrag{2}{\tiny Two Step}
  \includegraphics[width=0.5\textwidth,  trim={0 0 0 0},clip]{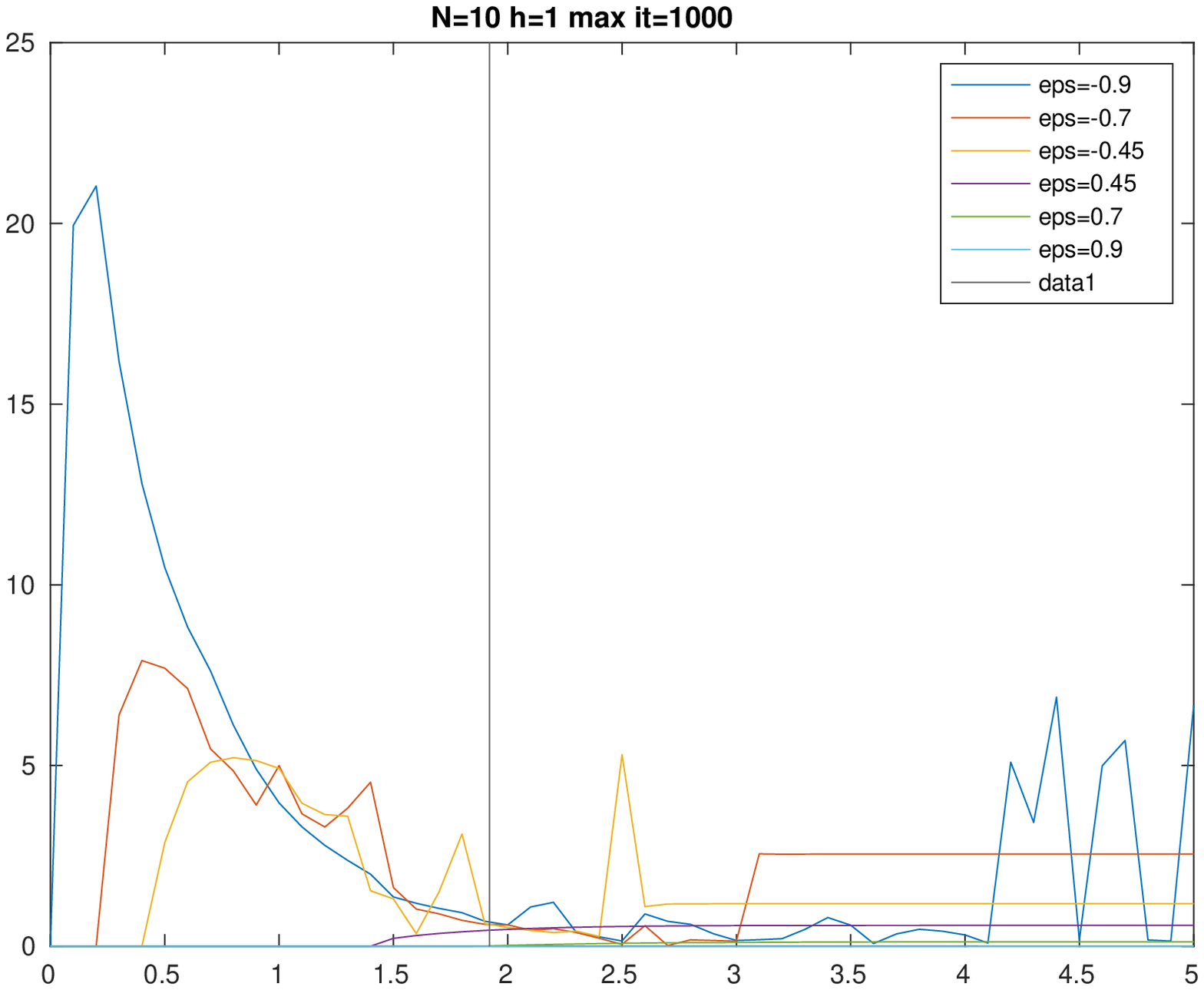}
  \caption{{\tiny} \space  1000 iterations of the $\varepsilon$-Banach algorithm for $N=10, \beta\in [0, 5]$ and $h=1$. We simulated one uniformly in the hypercube distributed start value for $\varepsilon = -0.9, -0.7, -0.45, 0.45, 0.7, 0.9$. We plotted the  {\italic mean squared error} for the given parameters and used the same start value for all $\beta$.}
\end{figure}
%%%%%% 
%%%%% 
\begin{figure}[H]
  \centering
  % \psfrag{2}{\tiny Two Step}
  \includegraphics[width=0.5\textwidth,  trim={0 0 0 0},clip]{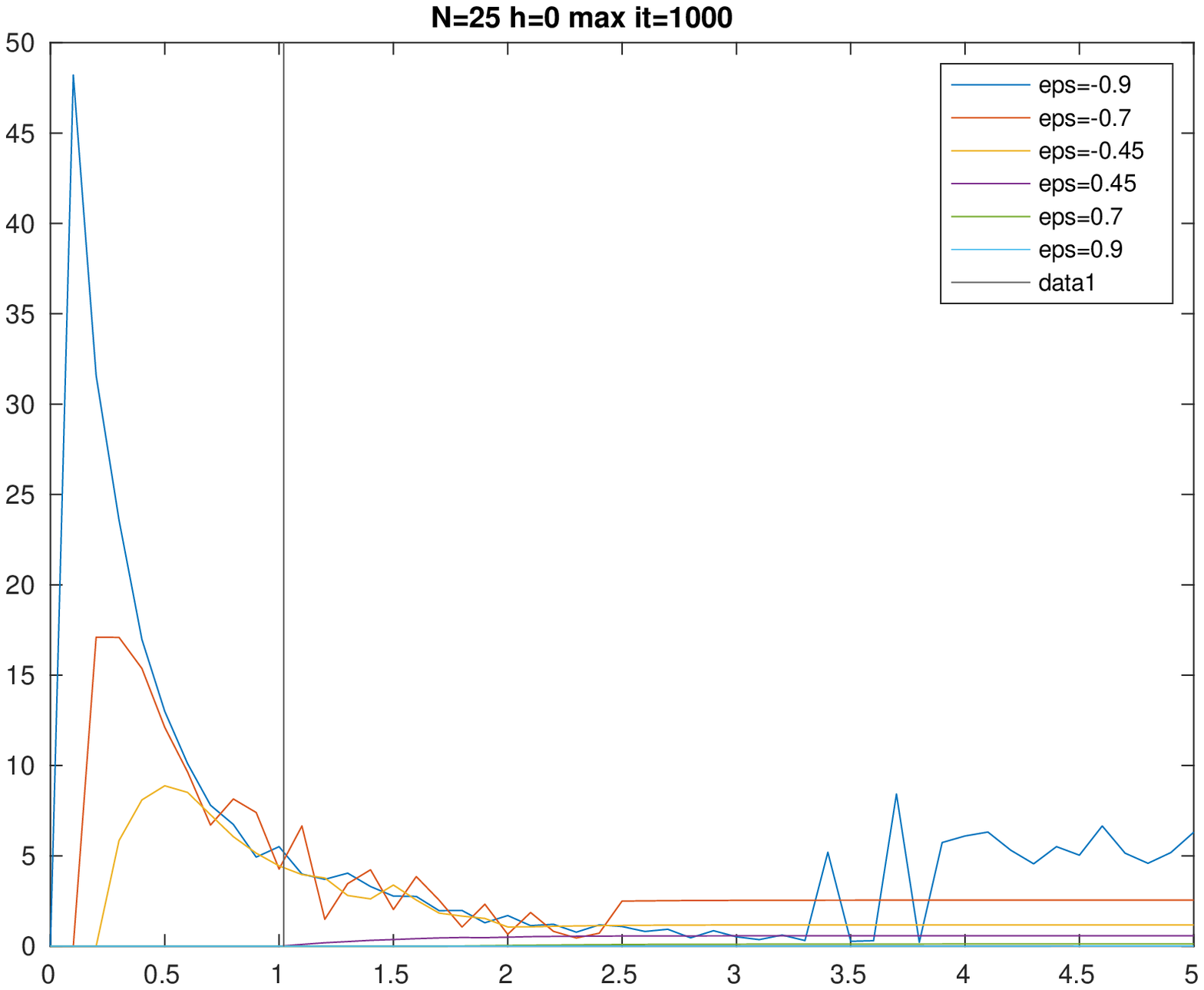}
  \caption{{\tiny} \space  1000 iterations of the $\varepsilon$-Banach algorithm for $N=25, \beta\in [0, 5]$ and $h=0$. We simulated one uniformly in the hypercube distributed start value for $\varepsilon = -0.9, -0.7, -0.45, 0.45, 0.7, 0.9$. We plotted the  {\italic mean squared error} for the given parameters and used the same start value for all $\beta$.}
\end{figure}
\begin{figure}[H]
  \centering
  % \psfrag{2}{\tiny Two Step}
  \includegraphics[width=0.5\textwidth,  trim={0 0 0 0},clip]{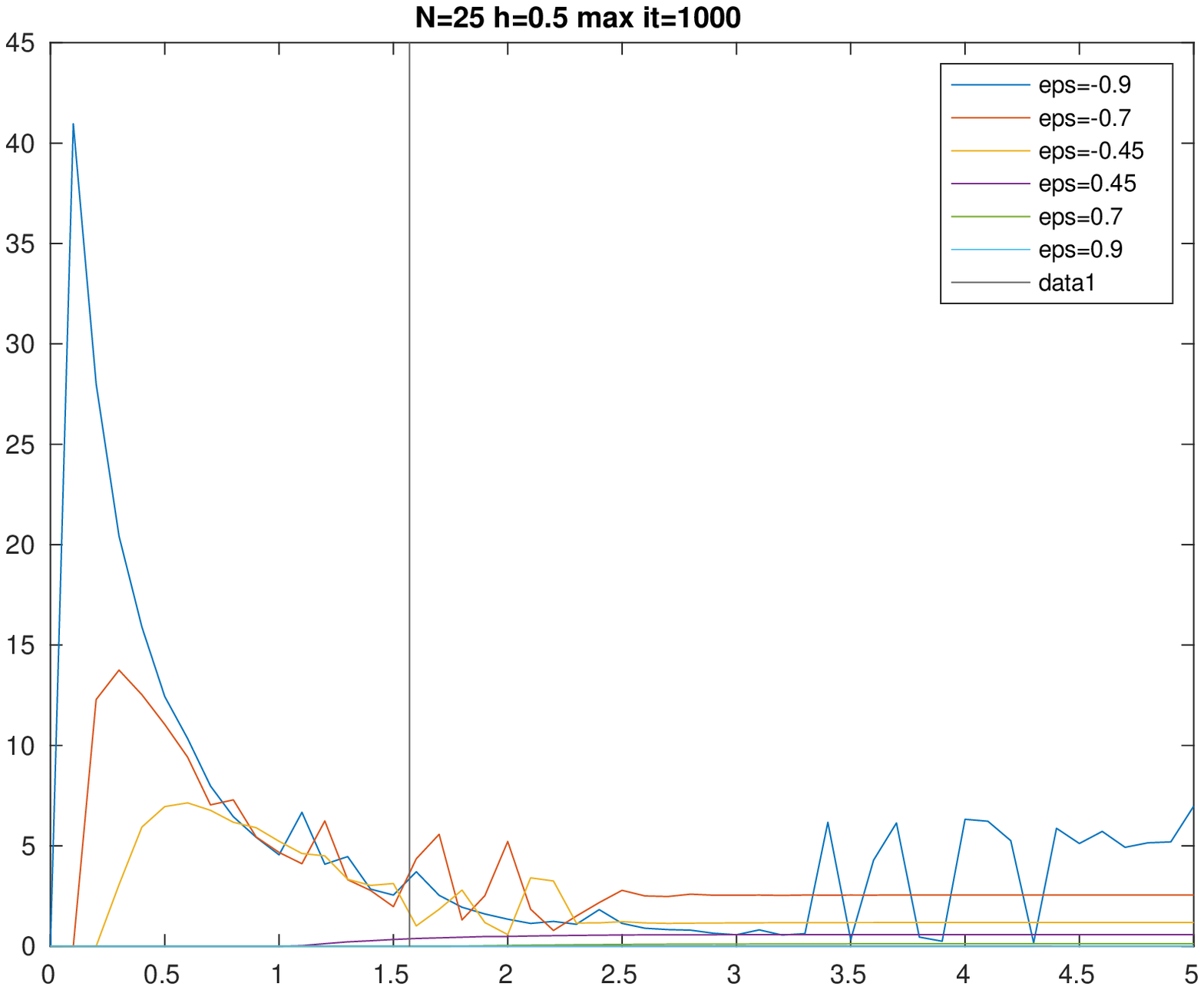}
  \caption{{\tiny} \space  1000 iterations of the $\varepsilon$-Banach algorithm for $N=25, \beta\in [0, 5]$ and $h=0.5$. We simulated one uniformly in the hypercube distributed start value for $\varepsilon = -0.9, -0.7, -0.45, 0.45, 0.7, 0.9$. We plotted the  {\italic mean squared error} for the given parameters and used the same start value for all $\beta$.}
\end{figure}
\begin{figure}[H]
  \centering
  % \psfrag{2}{\tiny Two Step}
  \includegraphics[width=0.5\textwidth,  trim={0 0 0 0},clip]{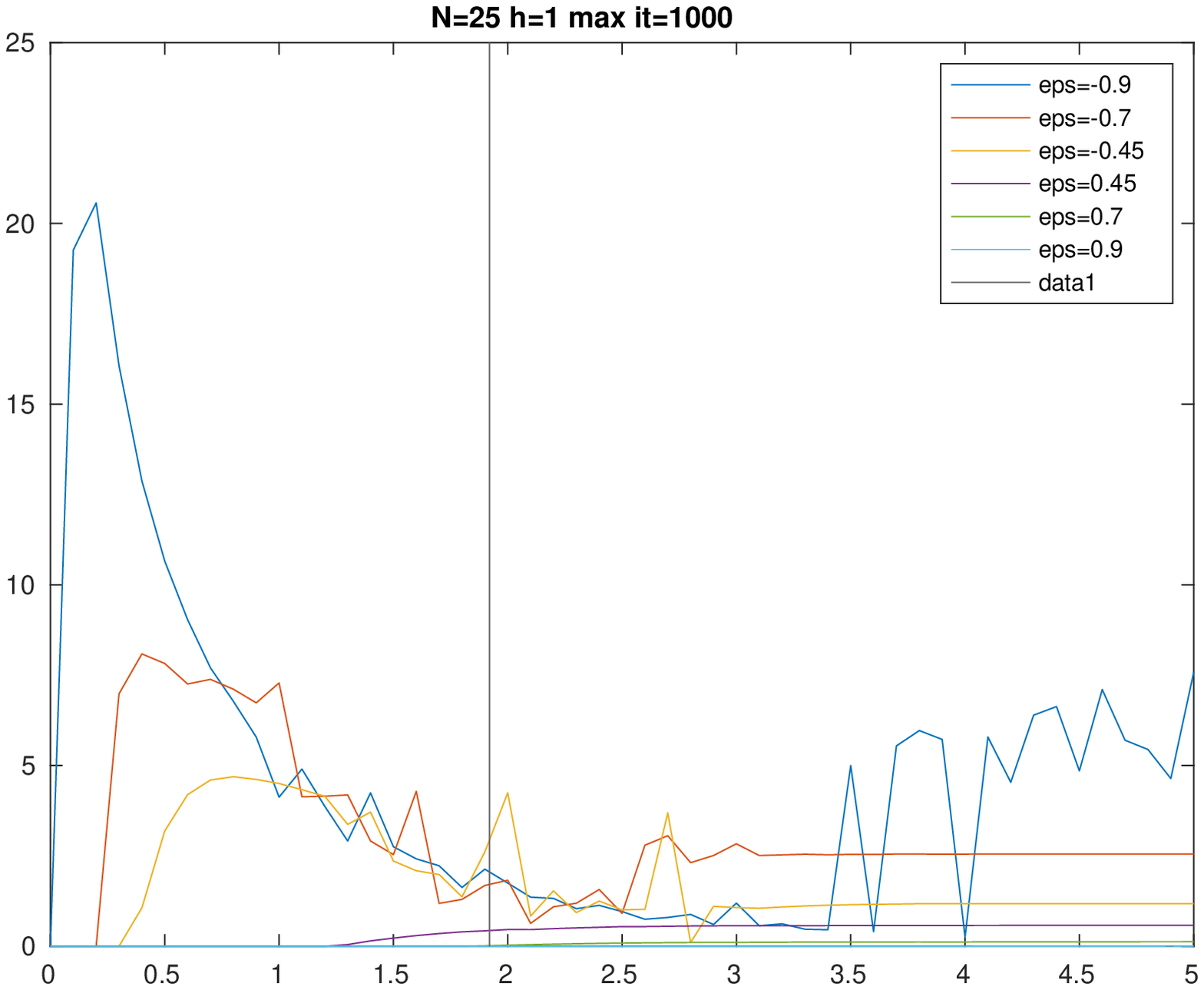}
  \caption{{\tiny} \space  1000 iterations of the $\varepsilon$-Banach algorithm for $N=25, \beta\in [0, 5]$ and $h=1$. We simulated one uniformly in the hypercube distributed start value for $\varepsilon = -0.9, -0.7, -0.45, 0.45, 0.7, 0.9$. We plotted the  {\italic mean squared error} for the given parameters and used the same start value for all $\beta$.}
\end{figure}
%%%%%% 
%%%%% 
\begin{figure}[H]
  \centering
  % \psfrag{2}{\tiny Two Step}
  \includegraphics[width=0.5\textwidth,  trim={0 0 0 0},clip]{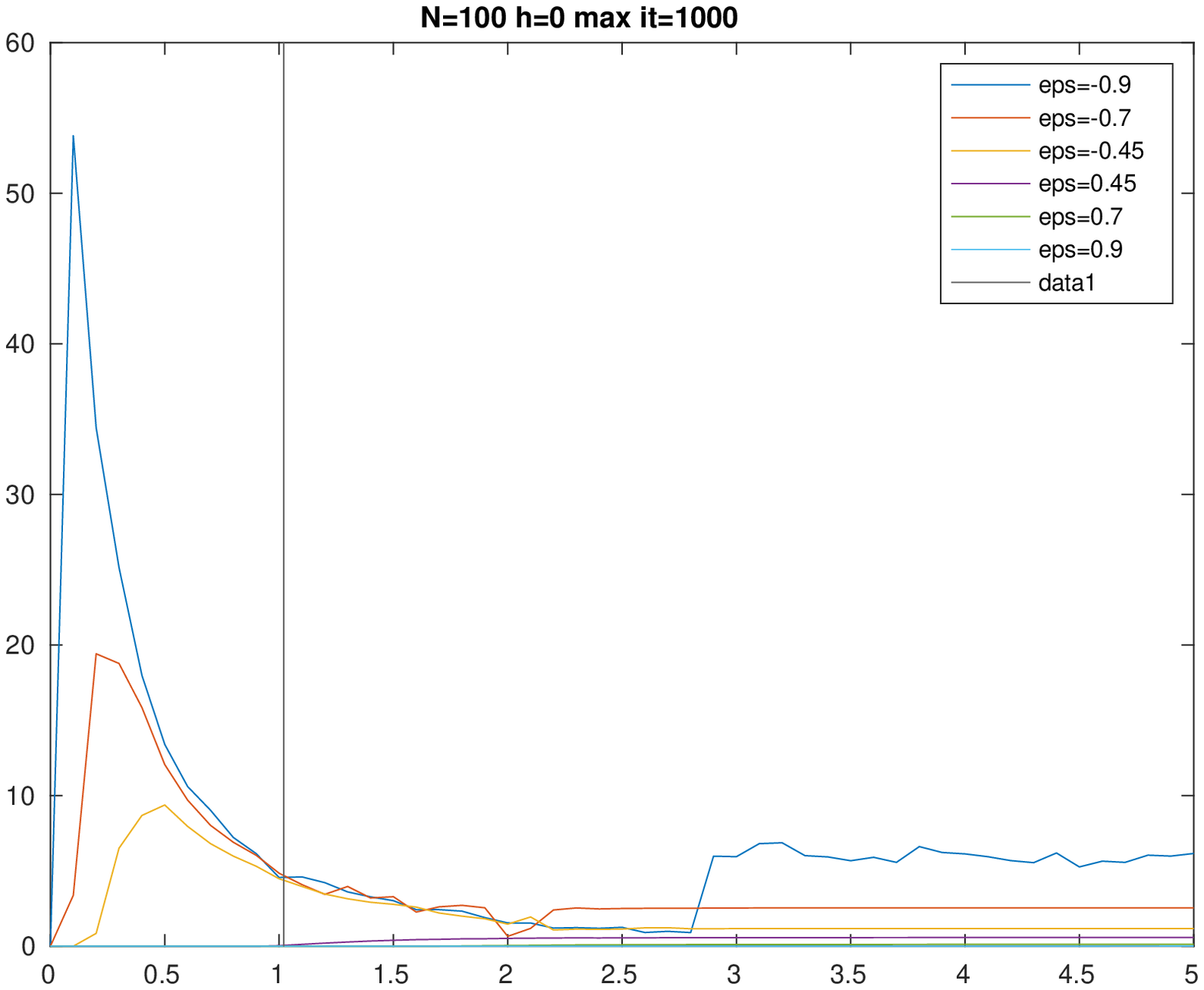}
  \caption{{\tiny} \space  1000 iterations of the $\varepsilon$-Banach algorithm for $N=100, \beta\in [0, 5]$ and $h=0$. We simulated one uniformly in the hypercube distributed start value for $\varepsilon = -0.9, -0.7, -0.45, 0.45, 0.7, 0.9$. We plotted the  {\italic mean squared error} for the given parameters and used the same start value for all $\beta$.}
\end{figure}
\begin{figure}[H]
  \centering
  % \psfrag{2}{\tiny Two Step}
  \includegraphics[width=0.5\textwidth,  trim={0 0 0 0},clip]{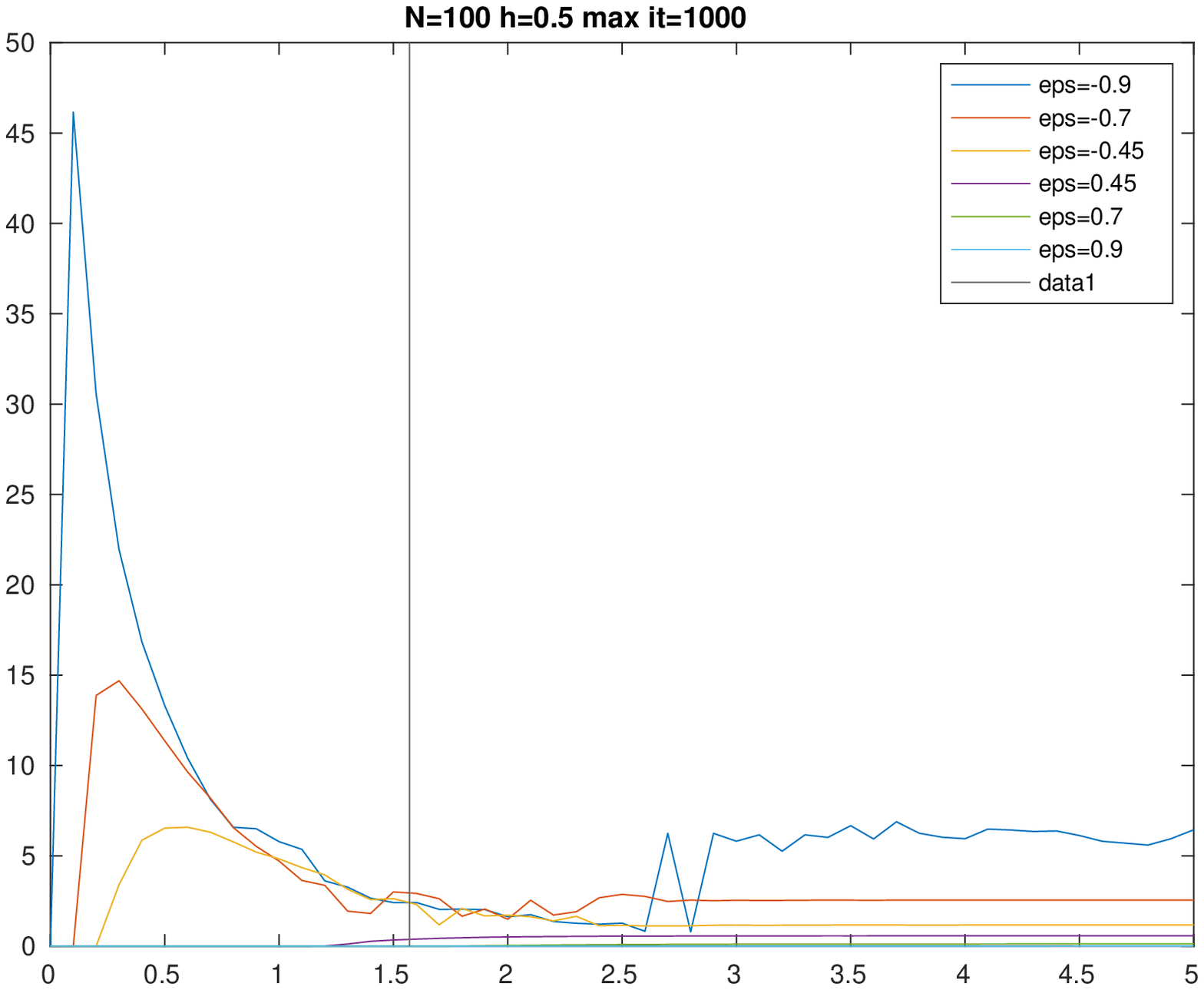}
  \caption{{\tiny} \space  1000 iterations of the $\varepsilon$-Banach algorithm for $N=100, \beta\in [0, 5]$ and $h=0.5$. We simulated one uniformly in the hypercube distributed start value for $\varepsilon = -0.9, -0.7, -0.45, 0.45, 0.7, 0.9$. We plotted the  {\italic mean squared error} for the given parameters and used the same start value for all $\beta$.}
\end{figure}
\begin{figure}[H]
  \centering
  % \psfrag{2}{\tiny Two Step}
  \includegraphics[width=0.5\textwidth,  trim={0 0 0 0},clip]{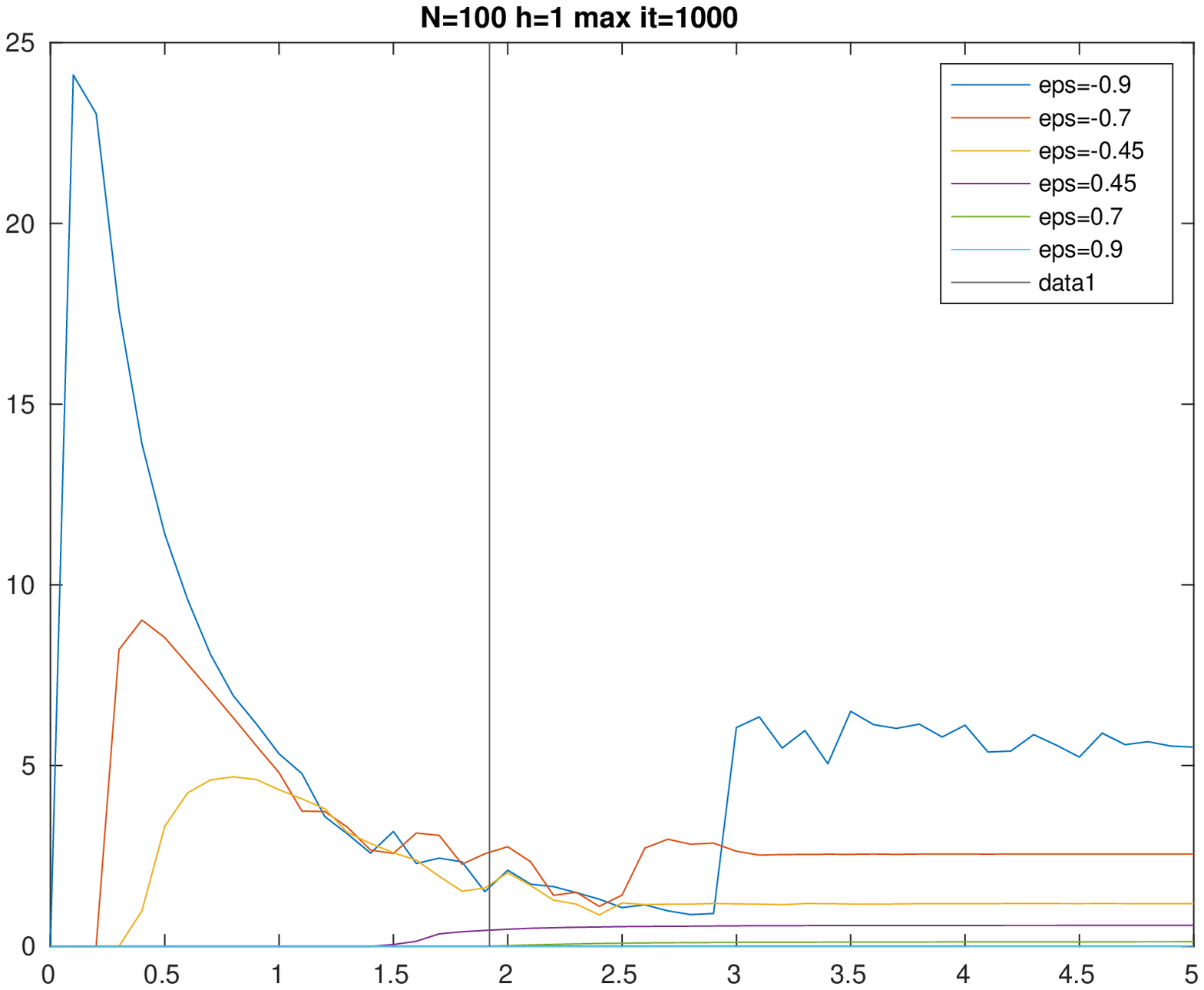}
  \caption{{\tiny} \space  1000 iterations of the $\varepsilon$-Banach algorithm for $N=100, \beta\in [0, 5]$ and $h=1$. We simulated one uniformly in the hypercube distributed start value for $\varepsilon = -0.9, -0.7, -0.45, 0.45, 0.7, 0.9$. We plotted the  {\italic mean squared error} for the given parameters and used the same start value for all $\beta$.}
\end{figure}
%%%%%% 
%%%%% 
\begin{figure}[H]
  \centering
  % \psfrag{2}{\tiny Two Step}
  \includegraphics[width=0.5\textwidth,  trim={0 0 0 0},clip]{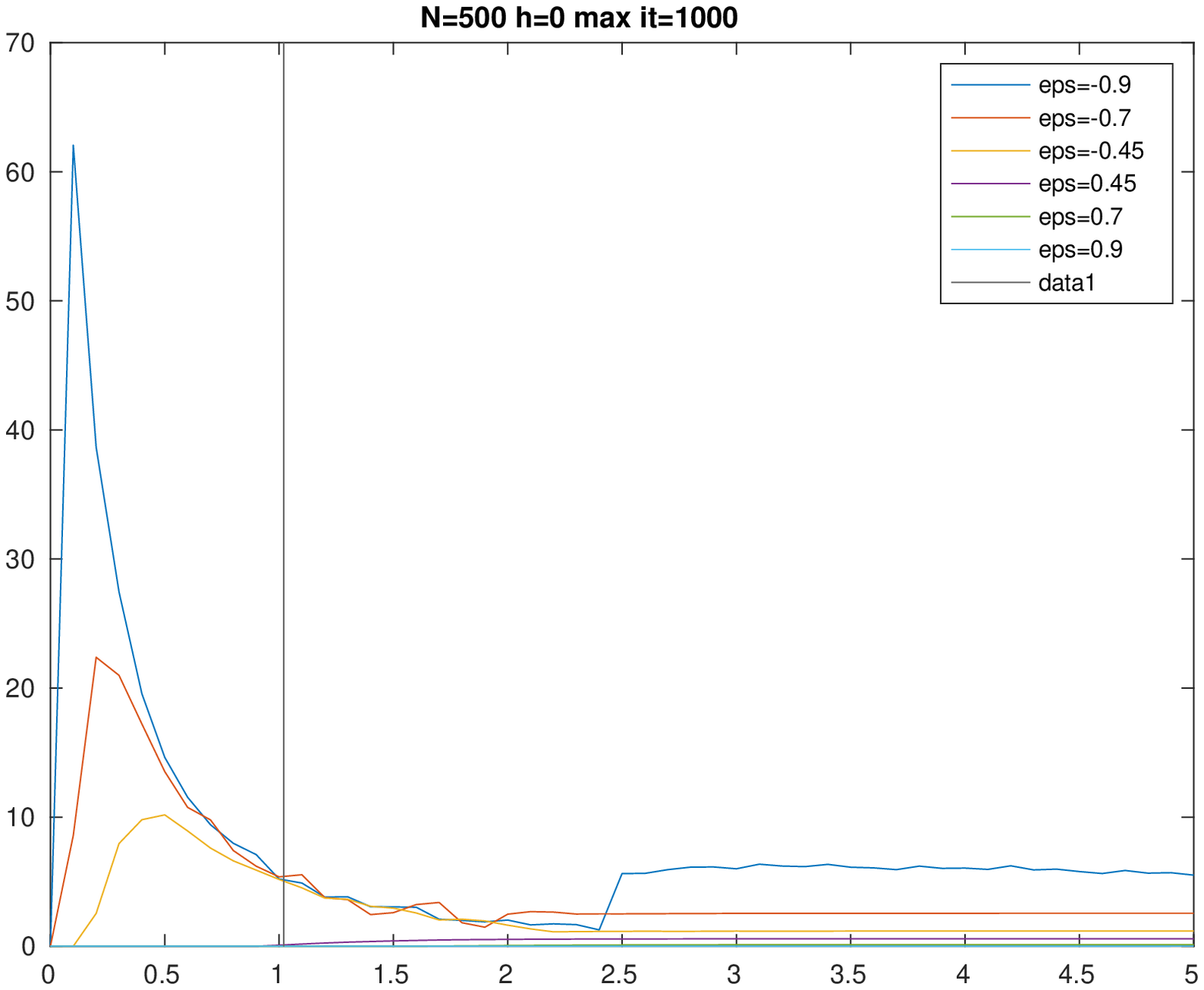}
  \caption{{\tiny} \space  1000 iterations of the $\varepsilon$-Banach algorithm for $N=500, \beta\in [0, 5]$ and $h=0$. We simulated one uniformly in the hypercube distributed start value for $\varepsilon = -0.9, -0.7, -0.45, 0.45, 0.7, 0.9$. We plotted the  {\italic mean squared error} for the given parameters and used the same start value for all $\beta$.}
\end{figure}
\begin{figure}[H]
  \centering
  % \psfrag{2}{\tiny Two Step}
  \includegraphics[width=0.5\textwidth,  trim={0 0 0 0},clip]{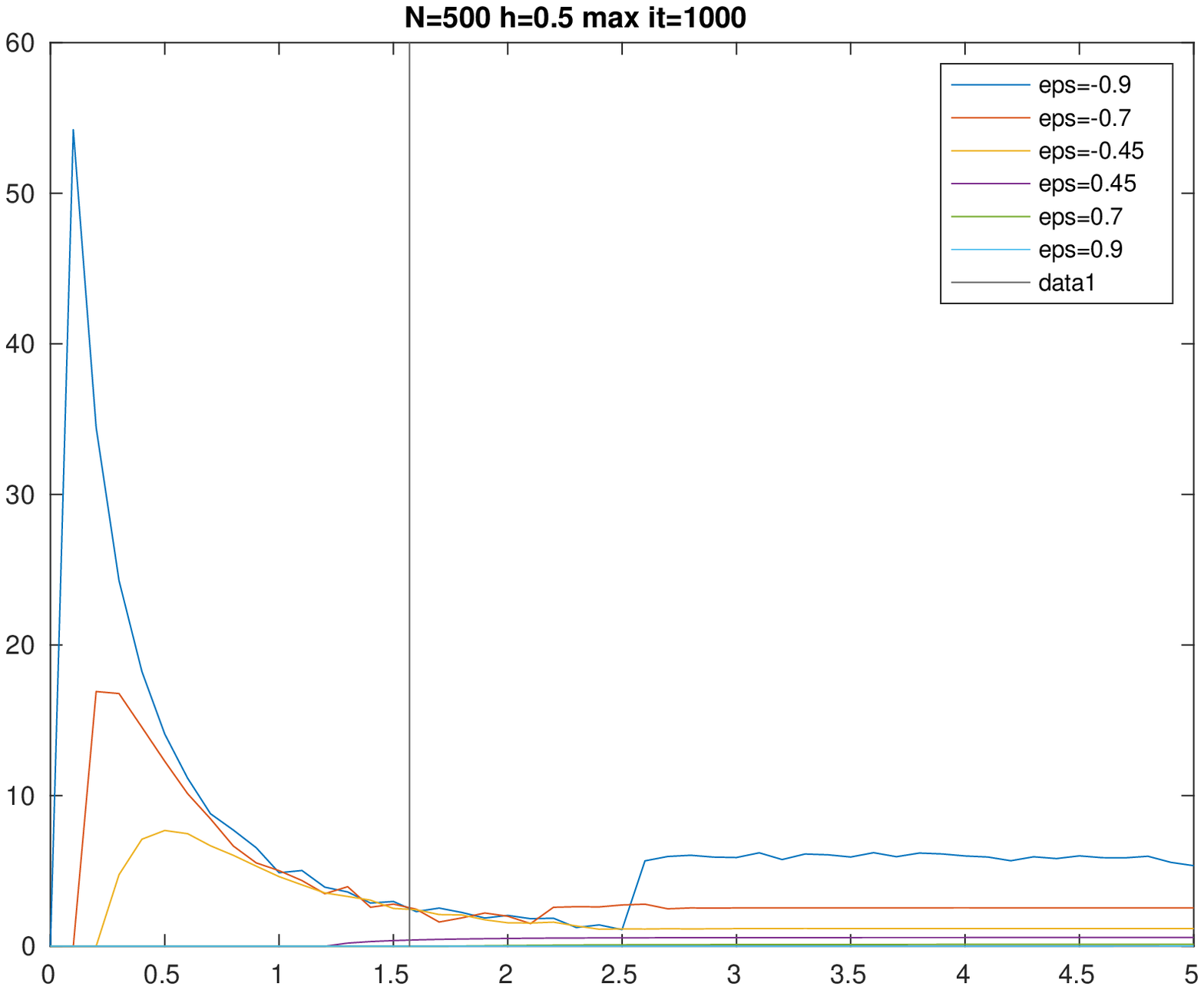}
  \caption{{\tiny} \space  1000 iterations of the $\varepsilon$-Banach algorithm for $N=500, \beta\in [0, 5]$ and $h=0.5$. We simulated one uniformly in the hypercube distributed start value for $\varepsilon = -0.9, -0.7, -0.45, 0.45, 0.7, 0.9$. We plotted the  {\italic mean squared error} for the given parameters and used the same start value for all $\beta$.}
\end{figure}
\begin{figure}[H]
  \centering
  % \psfrag{2}{\tiny Two Step}
  \includegraphics[width=0.5\textwidth,  trim={0 0 0 0},clip]{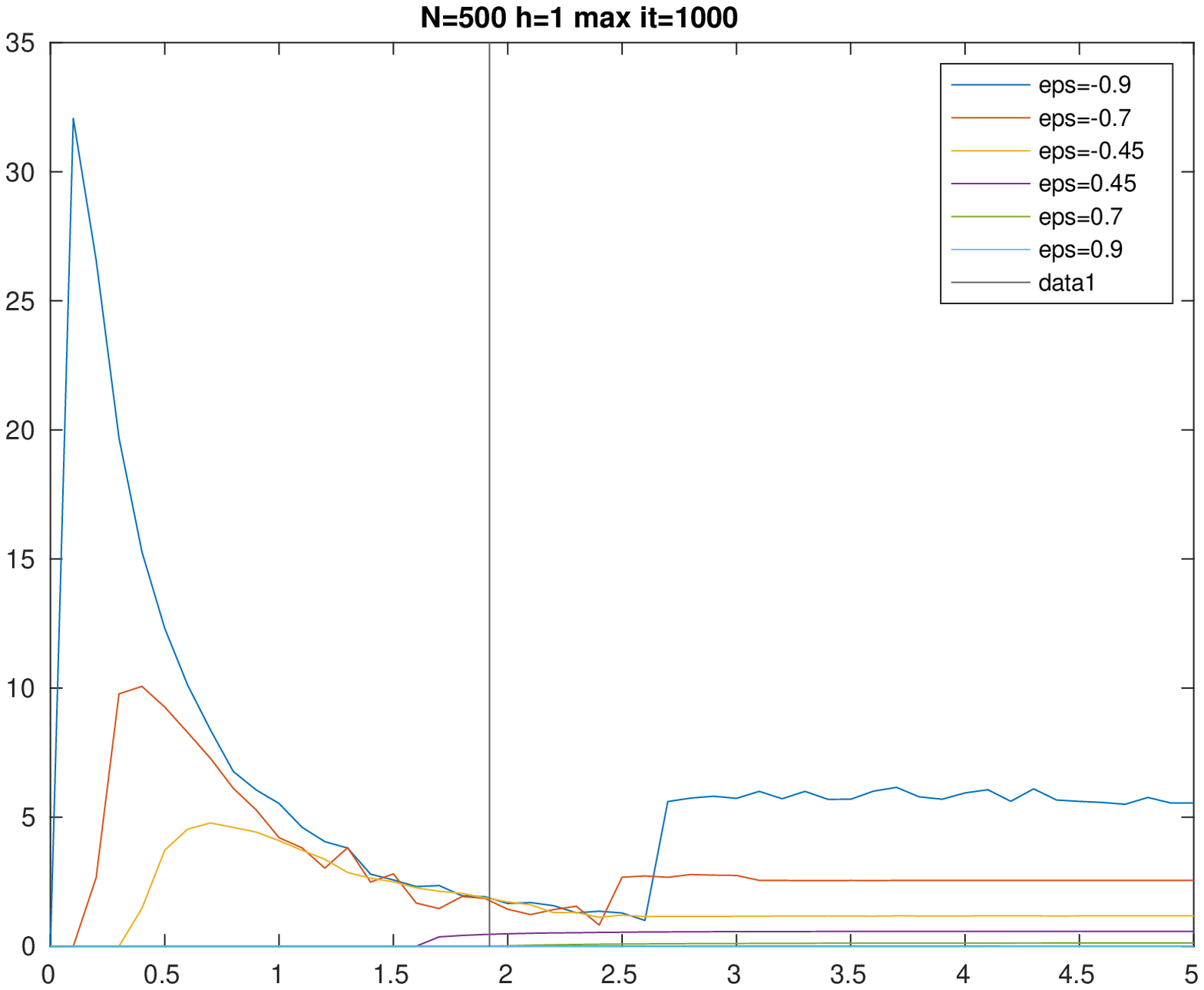}
  \caption{{\tiny} \space  1000 iterations of the $\varepsilon$-Banach algorithm for $N=500, \beta\in [0, 5]$ and $h=1$. We simulated one uniformly in the hypercube distributed start value for $\varepsilon = -0.9, -0.7, -0.45, 0.45, 0.7, 0.9$. We plotted the  {\italic mean squared error} for the given parameters and used the same start value for all $\beta$.}
\end{figure}
%%%%%% 
%%%%% 
\begin{figure}[H]
  \centering
  % \psfrag{2}{\tiny Two Step}
  \includegraphics[width=0.5\textwidth,  trim={0 0 0 0},clip]{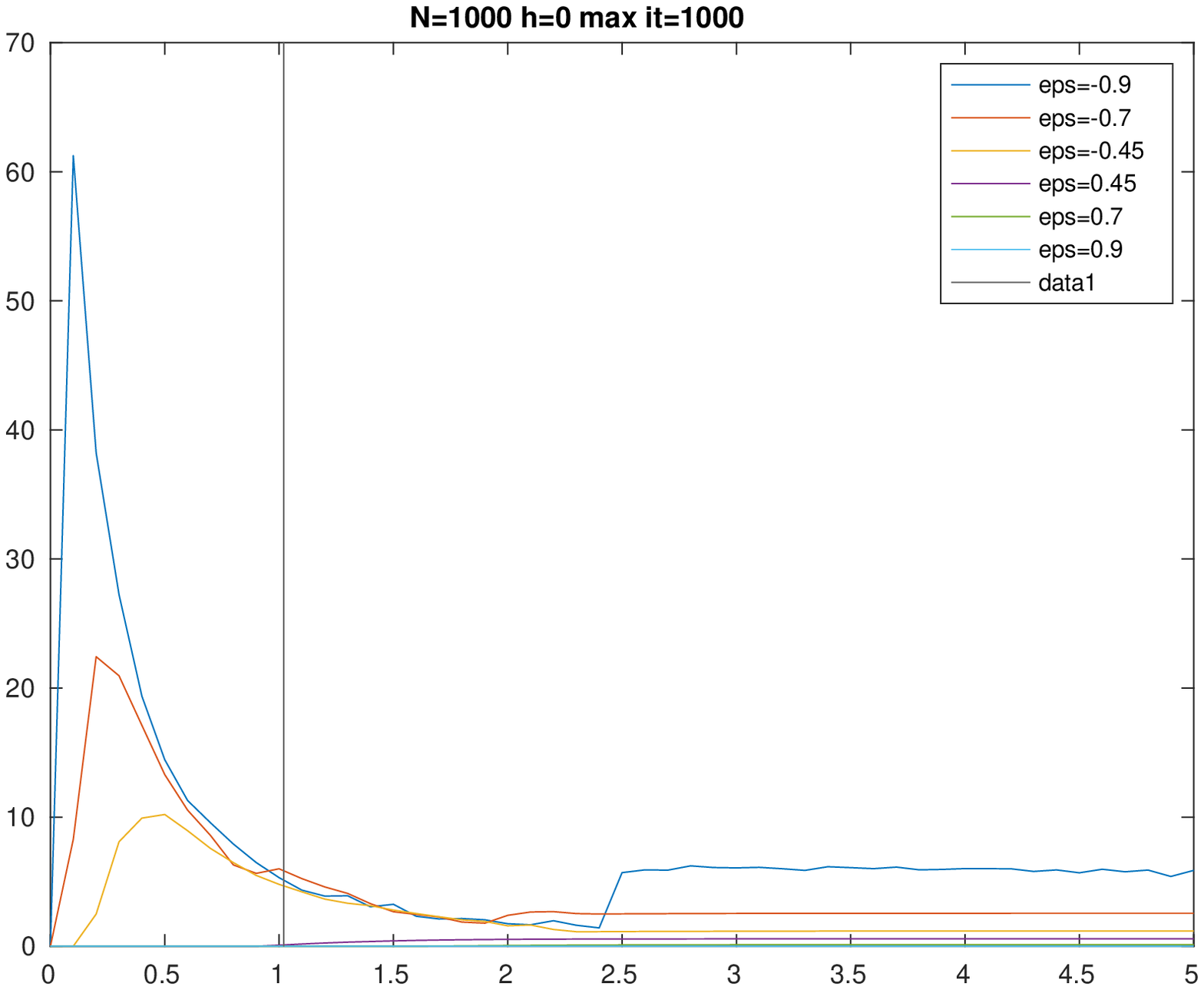}
  \caption{{\tiny} \space  1000 iterations of the $\varepsilon$-Banach algorithm for $N=1000, \beta\in [0, 5]$ and $h=0$. We simulated one uniformly in the hypercube distributed start value for $\varepsilon = -0.9, -0.7, -0.45, 0.45, 0.7, 0.9$. We plotted the {\italic mean squared error} for the given parameters and used the same start value for all $\beta$.}
\end{figure}
\begin{figure}[H]
  \centering
  % \psfrag{2}{\tiny Two Step}
  \includegraphics[width=0.5\textwidth,  trim={0 0 0 0},clip]{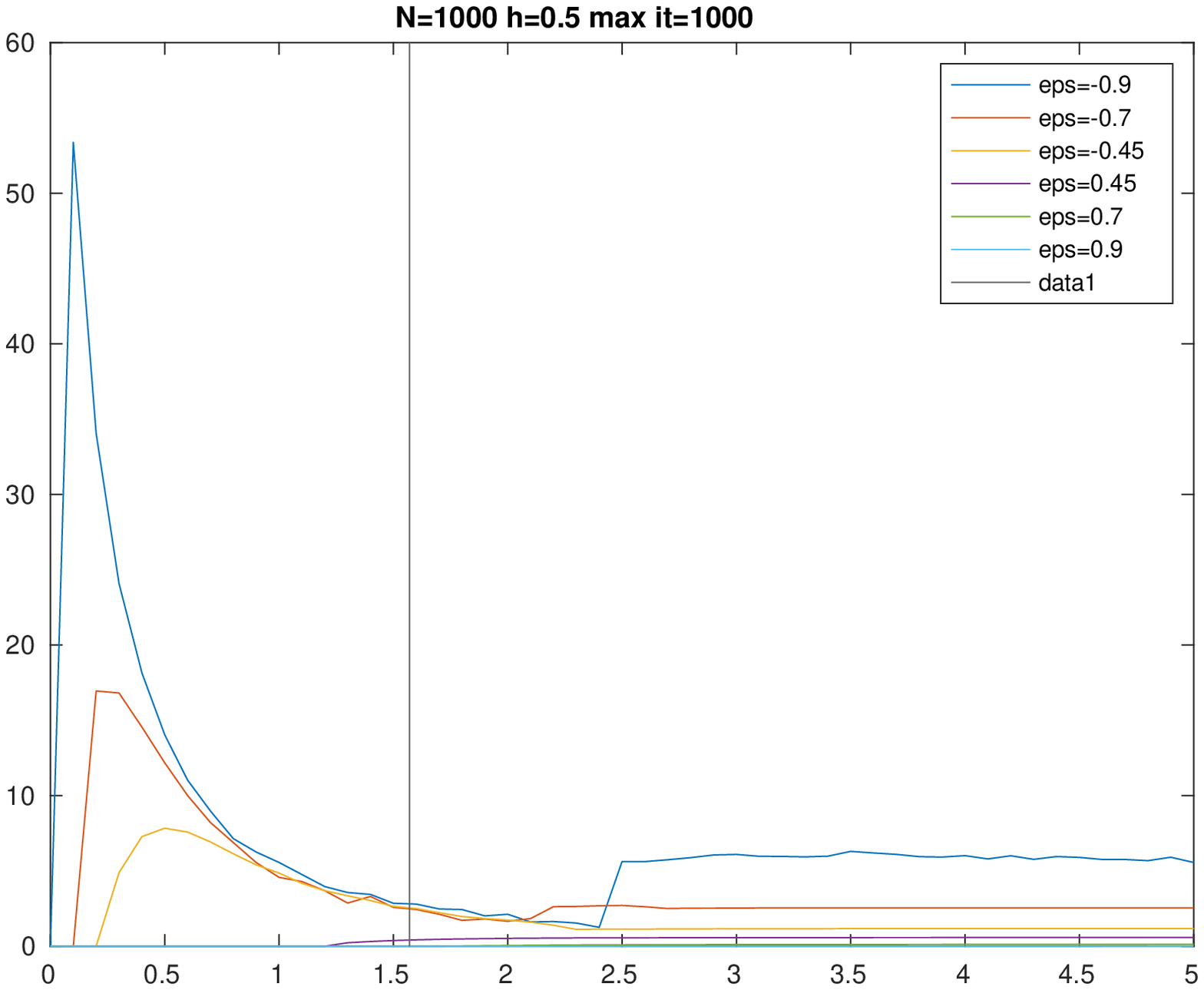}
  \caption{{\tiny} \space  1000 iterations of the $\varepsilon$-Banach algorithm for $N=1000, \beta\in [0, 5]$ and $h=0.5$. We simulated one uniformly in the hypercube distributed start value for $\varepsilon = -0.9, -0.7, -0.45, 0.45, 0.7, 0.9$. We plotted the {\italic mean squared error} for the given parameters and used the same start value for all $\beta$.}
\end{figure}
\begin{figure}[H]
  \centering
  % \psfrag{2}{\tiny Two Step}
  \includegraphics[width=0.5\textwidth,  trim={0 0 0 0},clip]{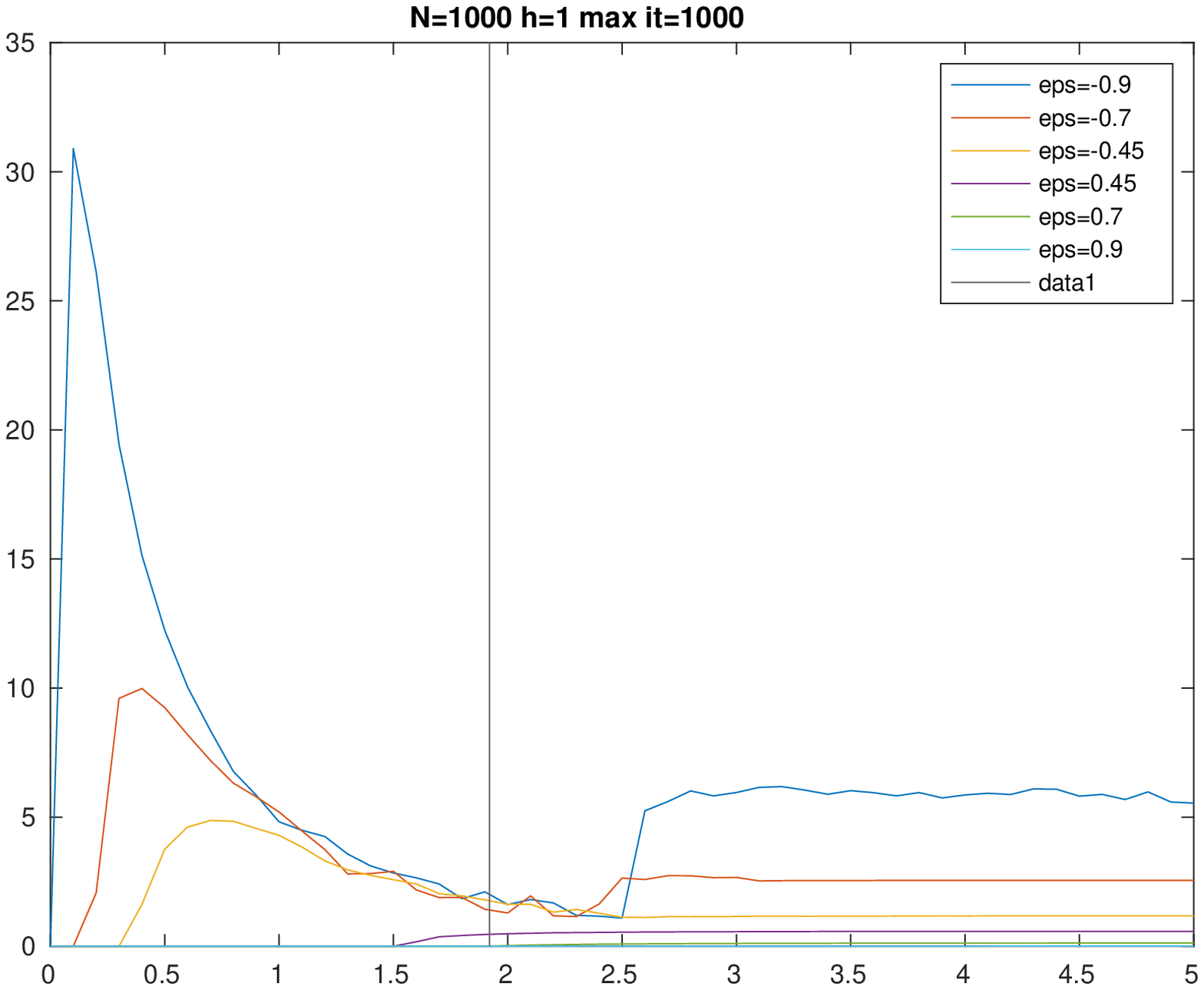}
  \caption{{\tiny} \space  1000 iterations of the $\varepsilon$-Banach algorithm for $N=1000, \beta\in [0, 5]$ and $h=1$. We simulated one uniformly in the hypercube distributed start value for $\varepsilon = -0.9, -0.7, -0.45, 0.45, 0.7, 0.9$. We plotted the {\italic mean squared error} for the given parameters and used the same start value for all $\beta$.}
\end{figure}
\FloatBarrier\clearpage
\subsubsection{Using the {\it maximum absolute error}}\leavevmode
\begin{figure}[H]
  \centering
  % \psfrag{2}{\tiny Two Step}
  \includegraphics[width=0.5\textwidth,  trim={0 0 0 0},clip]{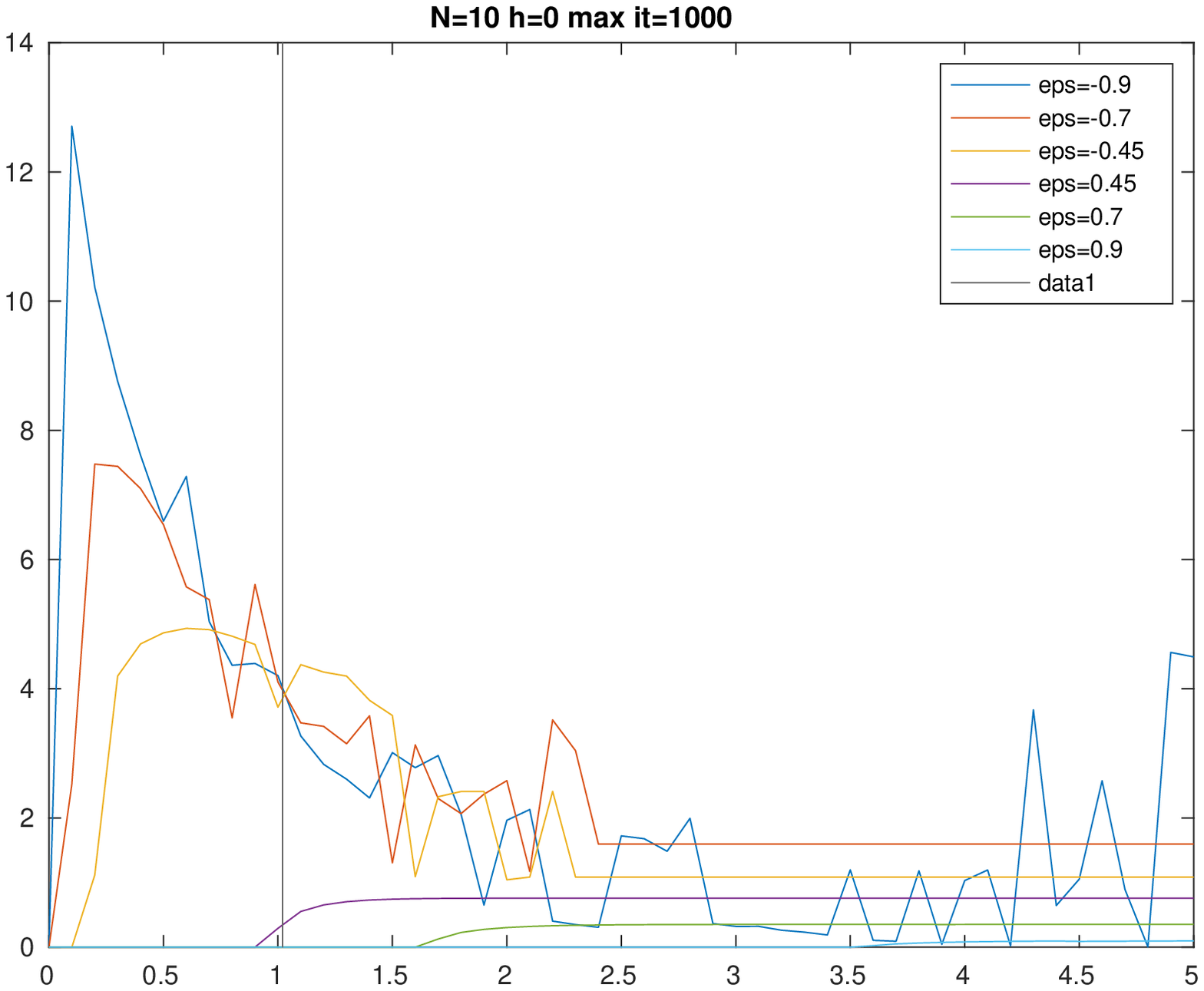}
  \caption{{\tiny} \space  1000 iterations of the $\varepsilon$-Banach algorithm for $N=10, \beta\in [0, 5]$ and $h=0$. We simulated one uniformly in the hypercube distributed start value for $\varepsilon = -0.9, -0.7, -0.45, 0.45, 0.7, 0.9$. We plotted the {\it maximum absolute error} for the given parameters and used the same start value for all $\beta$.}
\end{figure}
\begin{figure}[H]
  \centering
  % \psfrag{2}{\tiny Two Step}
  \includegraphics[width=0.5\textwidth,  trim={0 0 0 0},clip]{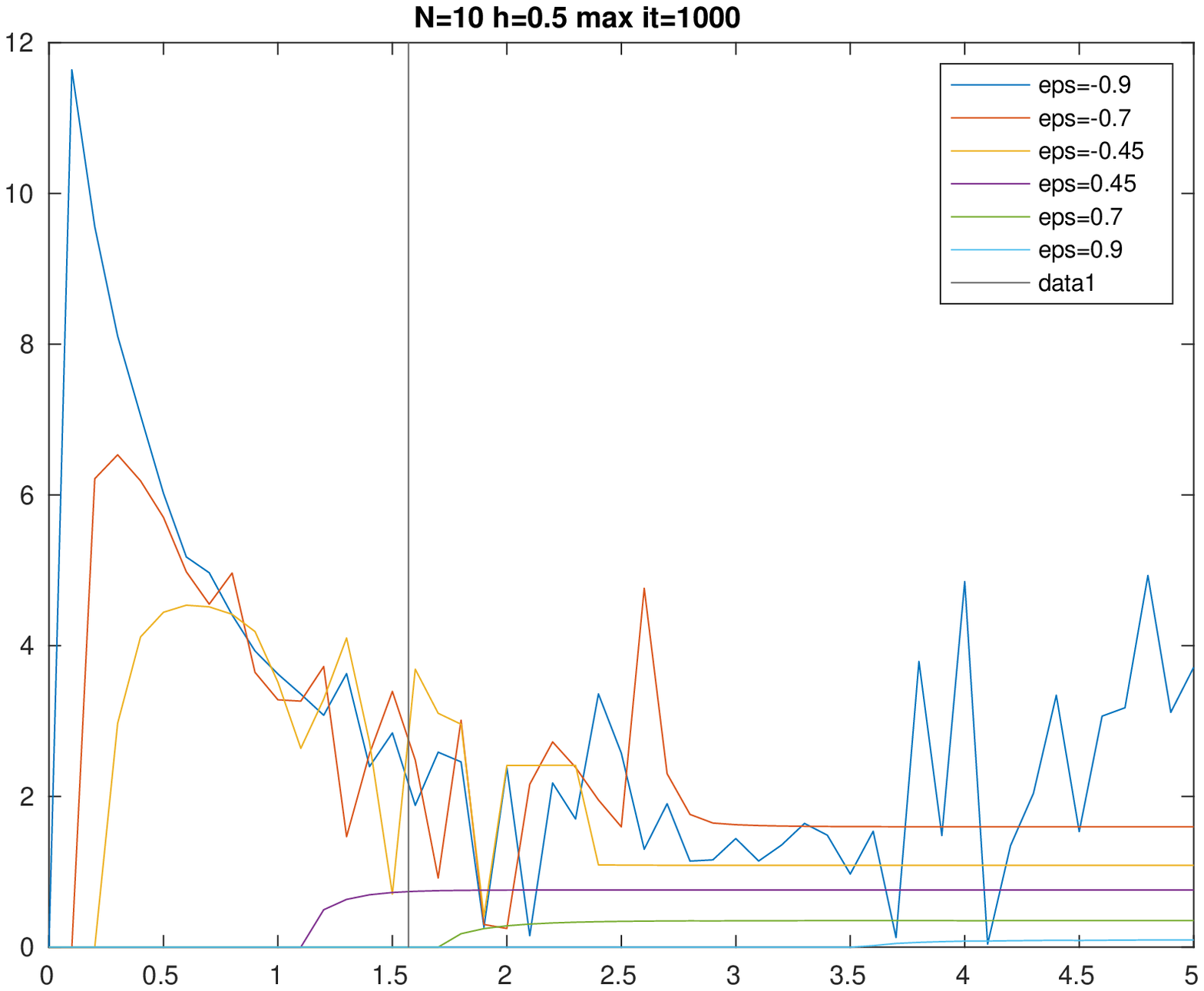}
  \caption{{\tiny} \space  1000 iterations of the $\varepsilon$-Banach algorithm for $N=10, \beta\in [0, 5]$ and $h=0.5$. We simulated one uniformly in the hypercube distributed start value for $\varepsilon = -0.9, -0.7, -0.45, 0.45, 0.7, 0.9$. We plotted the {\it maximum absolute error} for the given parameters and used the same start value for all $\beta$.}
\end{figure}
\begin{figure}[H]
  \centering
  % \psfrag{2}{\tiny Two Step}
  \includegraphics[width=0.5\textwidth,  trim={0 0 0 0},clip]{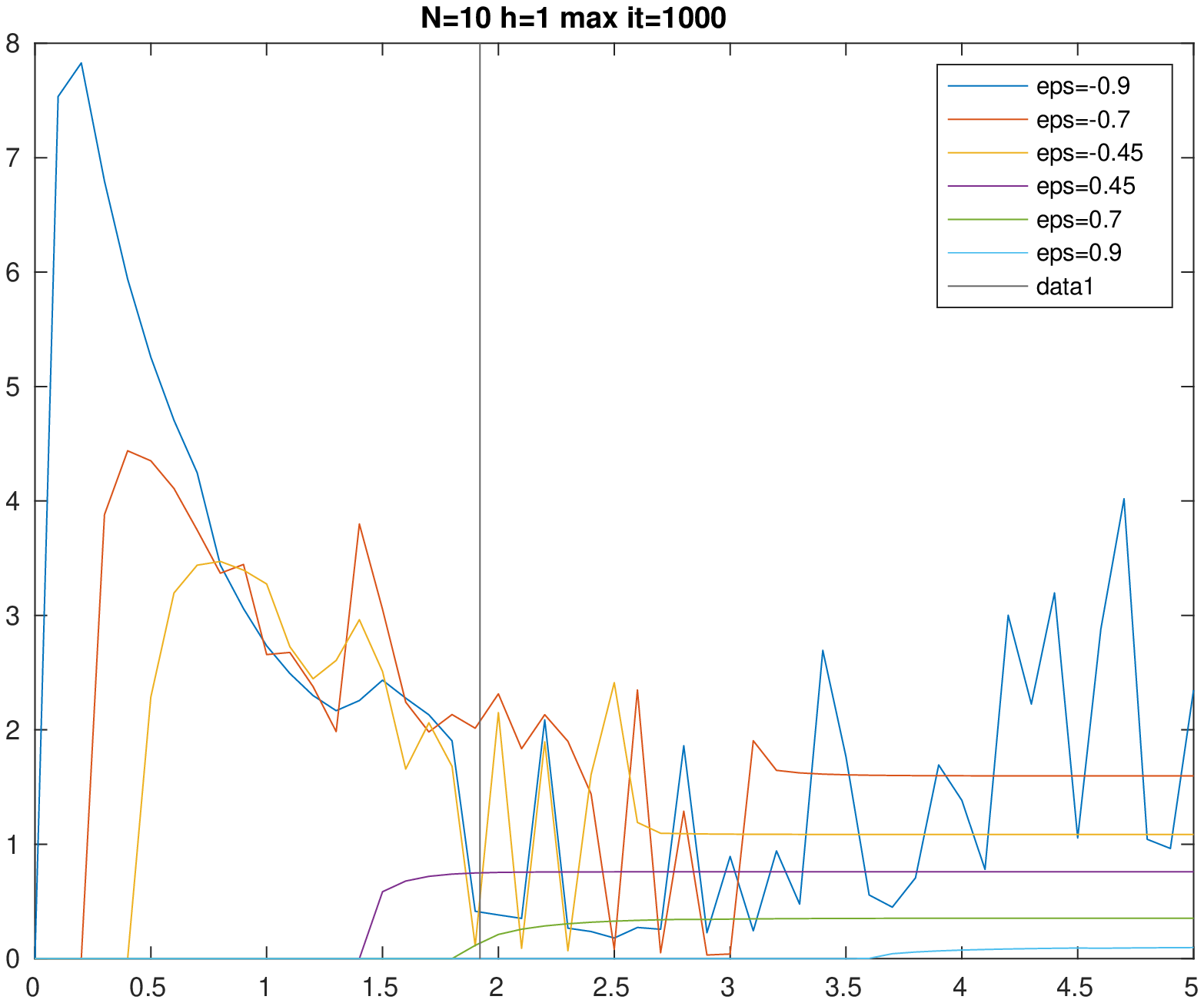}
  \caption{{\tiny} \space  1000 iterations of the $\varepsilon$-Banach algorithm for $N=10, \beta\in [0, 5]$ and $h=1$. We simulated one uniformly in the hypercube distributed start value for $\varepsilon = -0.9, -0.7, -0.45, 0.45, 0.7, 0.9$. We plotted the {\it maximum absolute error} for the given parameters and used the same start value for all $\beta$.}
\end{figure}
%%%%%% 
%%%%% 
\begin{figure}[H]
  \centering
  % \psfrag{2}{\tiny Two Step}
  \includegraphics[width=0.5\textwidth,  trim={0 0 0 0},clip]{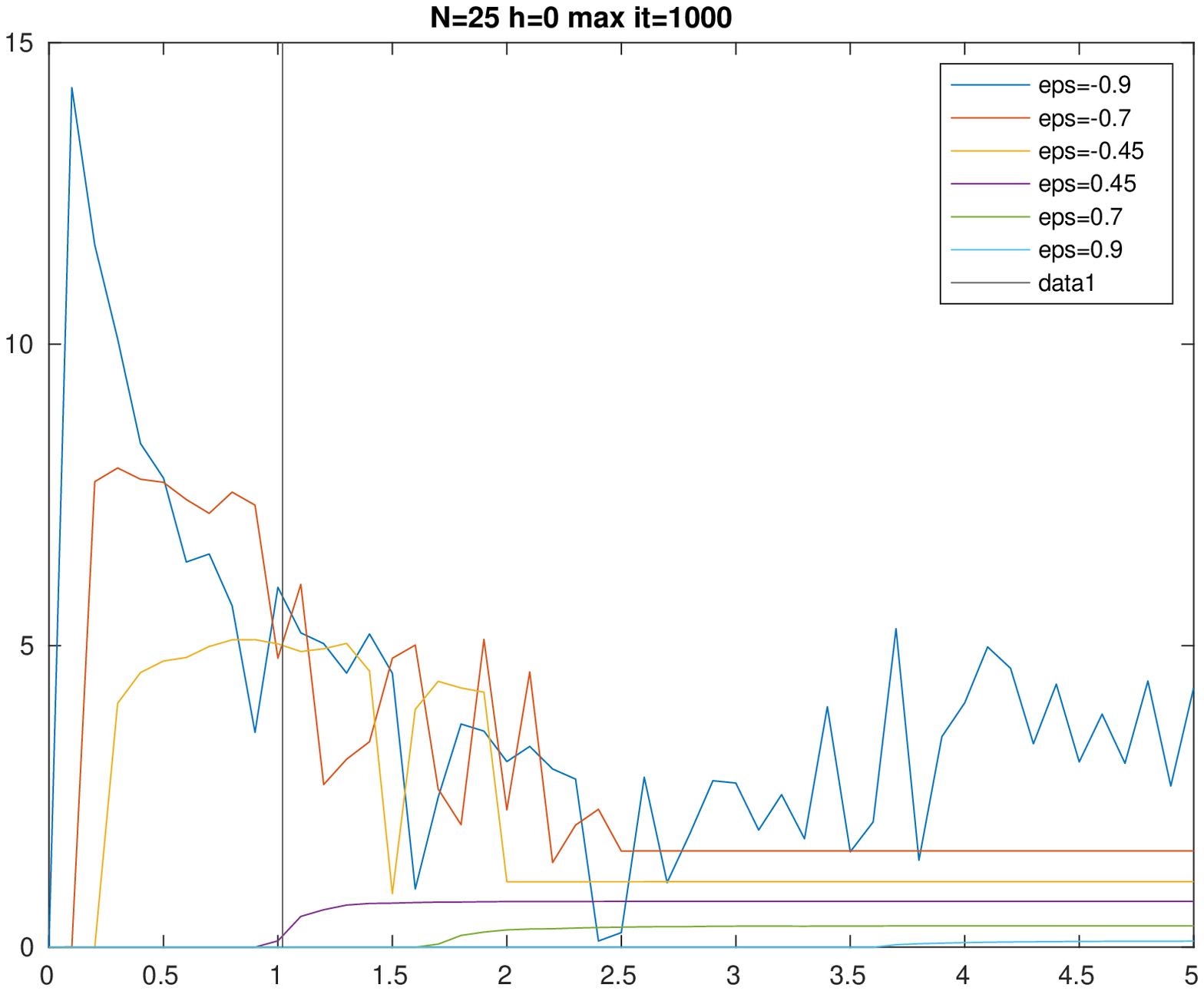}
  \caption{{\tiny} \space  1000 iterations of the $\varepsilon$-Banach algorithm for $N=25, \beta\in [0, 5]$ and $h=0$. We simulated one uniformly in the hypercube distributed start value for $\varepsilon = -0.9, -0.7, -0.45, 0.45, 0.7, 0.9$. We plotted the {\it maximum absolute error} for the given parameters and used the same start value for all $\beta$.}
\end{figure}
\begin{figure}[H]
  \centering
  % \psfrag{2}{\tiny Two Step}
  \includegraphics[width=0.5\textwidth,  trim={0 0 0 0},clip]{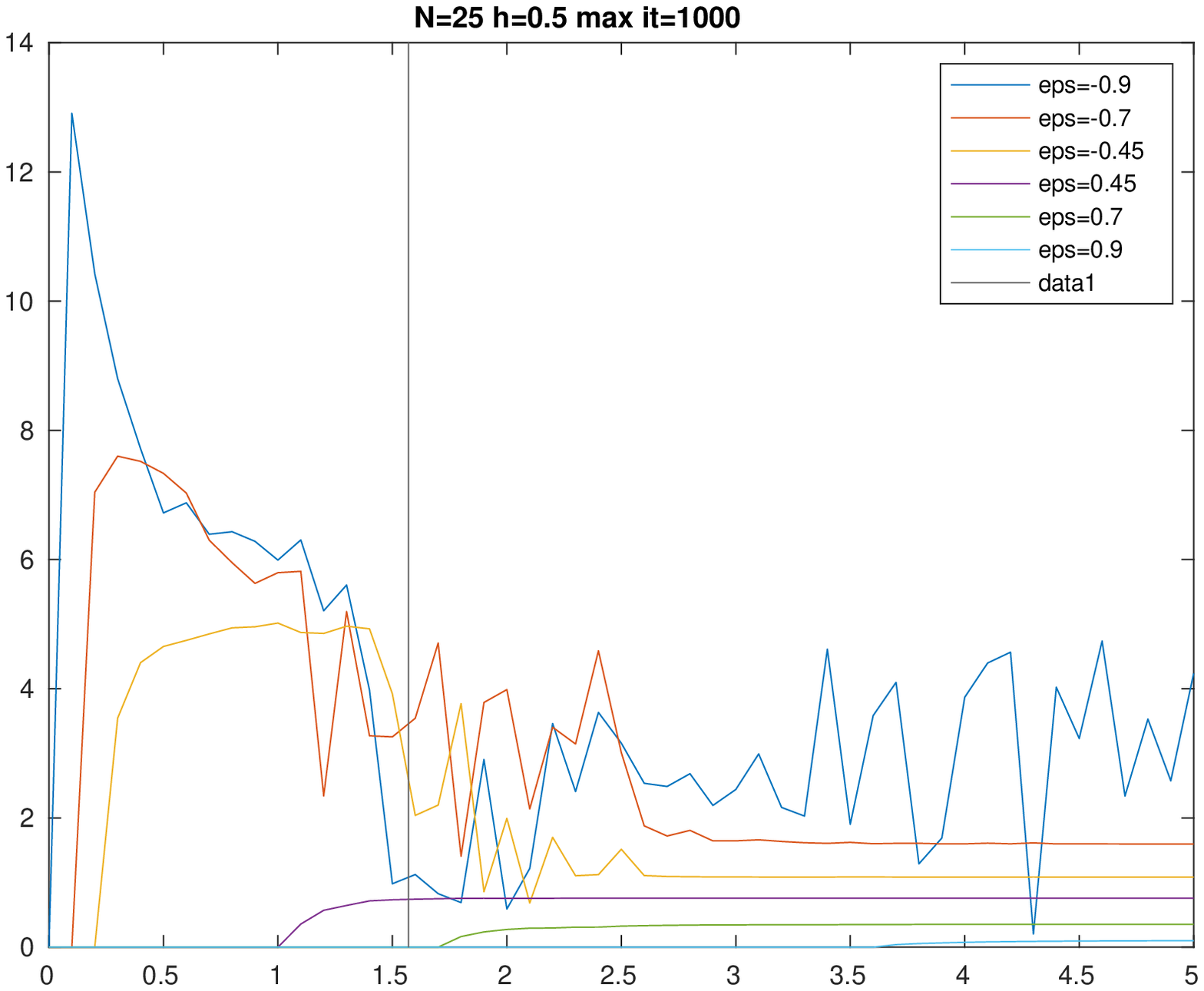}
  \caption{{\tiny} \space  1000 iterations of the $\varepsilon$-Banach algorithm for $N=25, \beta\in [0, 5]$ and $h=0.5$. We simulated one uniformly in the hypercube distributed start value for $\varepsilon = -0.9, -0.7, -0.45, 0.45, 0.7, 0.9$. We plotted the {\it maximum absolute error} for the given parameters and used the same start value for all $\beta$.}
\end{figure}
\begin{figure}[H]
  \centering
  % \psfrag{2}{\tiny Two Step}
  \includegraphics[width=0.5\textwidth,  trim={0 0 0 0},clip]{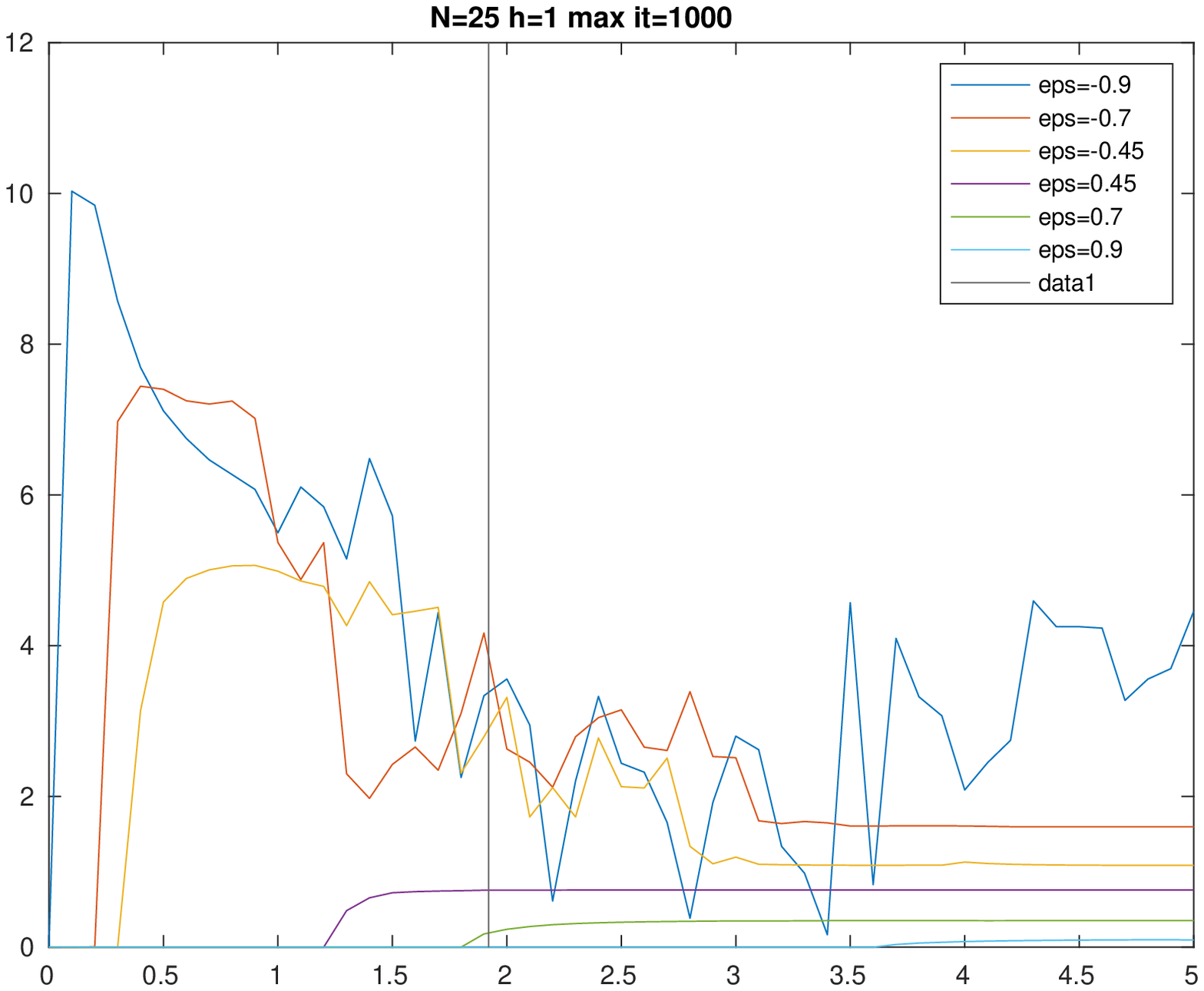}
  \caption{{\tiny} \space  1000 iterations of the $\varepsilon$-Banach algorithm for $N=25, \beta\in [0, 5]$ and $h=1$. We simulated one uniformly in the hypercube distributed start value for $\varepsilon = -0.9, -0.7, -0.45, 0.45, 0.7, 0.9$. We plotted the {\it maximum absolute error} for the given parameters and used the same start value for all $\beta$.}
\end{figure}
%%%%%% 
%%%%% 
\begin{figure}[H]
  \centering
  % \psfrag{2}{\tiny Two Step}
  \includegraphics[width=0.5\textwidth,  trim={0 0 0 0},clip]{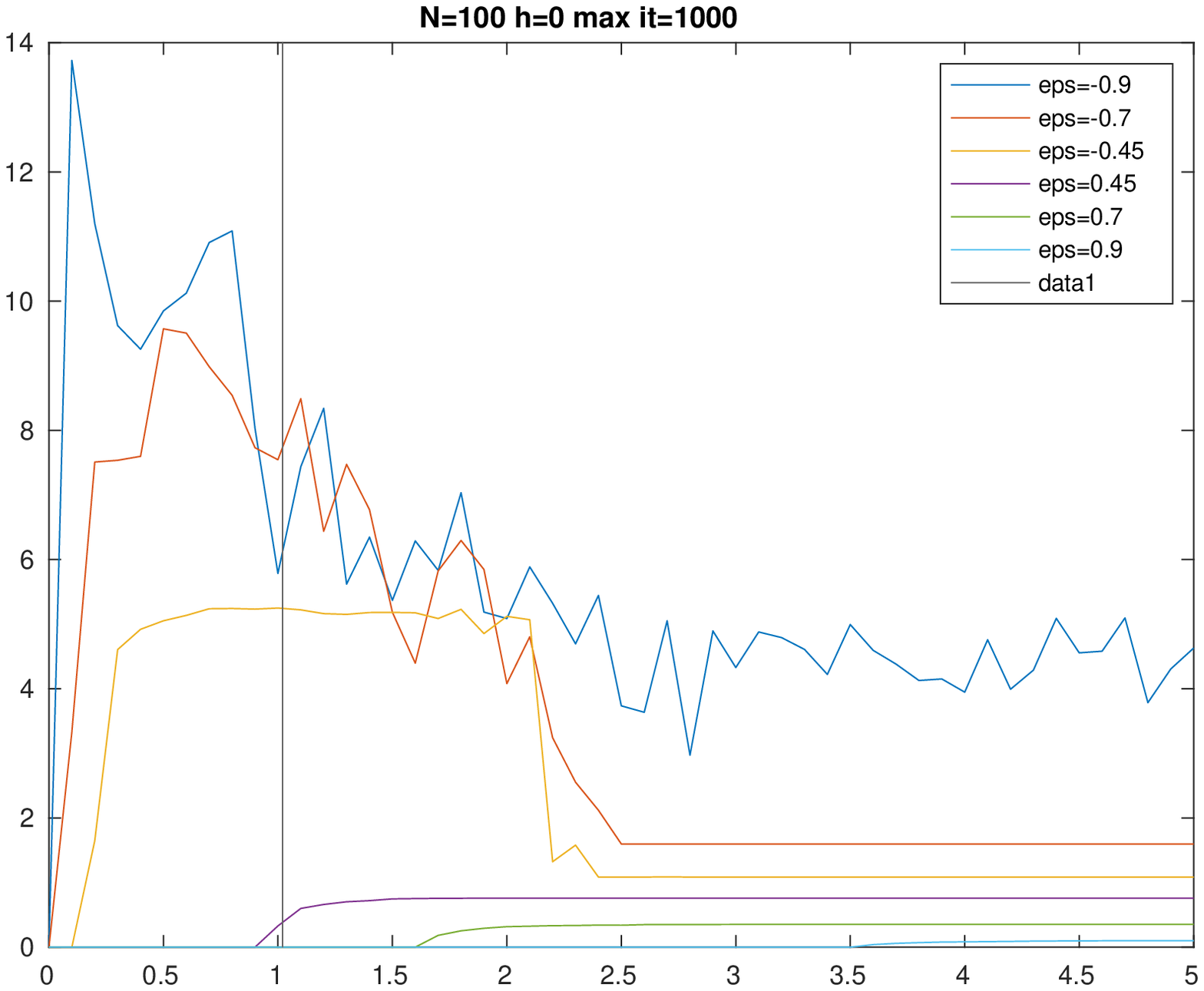}
  \caption{{\tiny} \space  1000 iterations of the $\varepsilon$-Banach algorithm for $N=100, \beta\in [0, 5]$ and $h=0$. We simulated one uniformly in the hypercube distributed start value for $\varepsilon = -0.9, -0.7, -0.45, 0.45, 0.7, 0.9$. We plotted the {\it maximum absolute error} for the given parameters and used the same start value for all $\beta$.}
\end{figure}
\begin{figure}[H]
  \centering
  % \psfrag{2}{\tiny Two Step}
  \includegraphics[width=0.5\textwidth,  trim={0 0 0 0},clip]{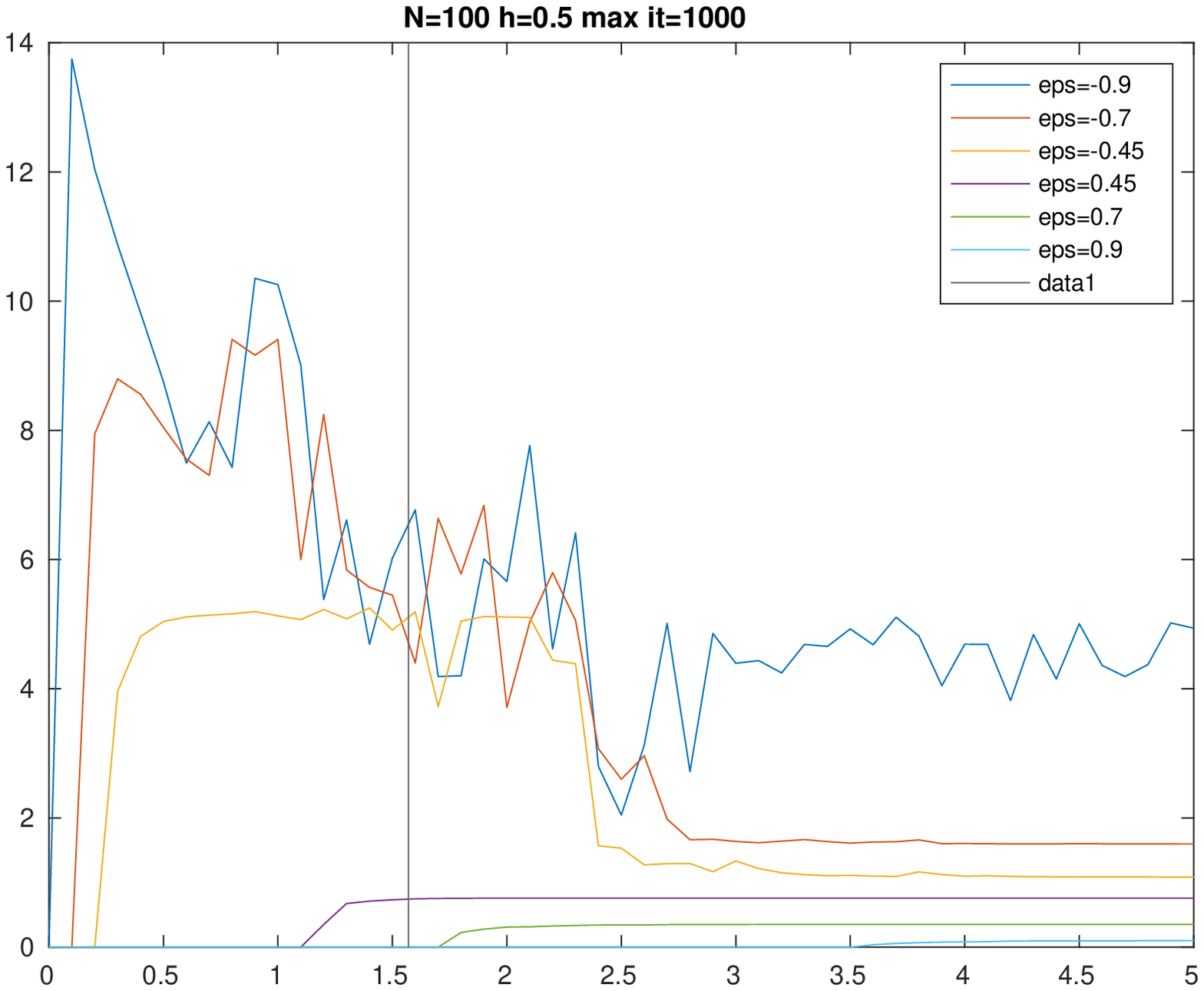}
  \caption{{\tiny} \space  1000 iterations of the $\varepsilon$-Banach algorithm for $N=100, \beta\in [0, 5]$ and $h=0.5$. We simulated one uniformly in the hypercube distributed start value for $\varepsilon = -0.9, -0.7, -0.45, 0.45, 0.7, 0.9$. We plotted the {\it maximum absolute error} for the given parameters and used the same start value for all $\beta$.}
\end{figure}
\begin{figure}[H]
  \centering
  % \psfrag{2}{\tiny Two Step}
  \includegraphics[width=0.5\textwidth,  trim={0 0 0 0},clip]{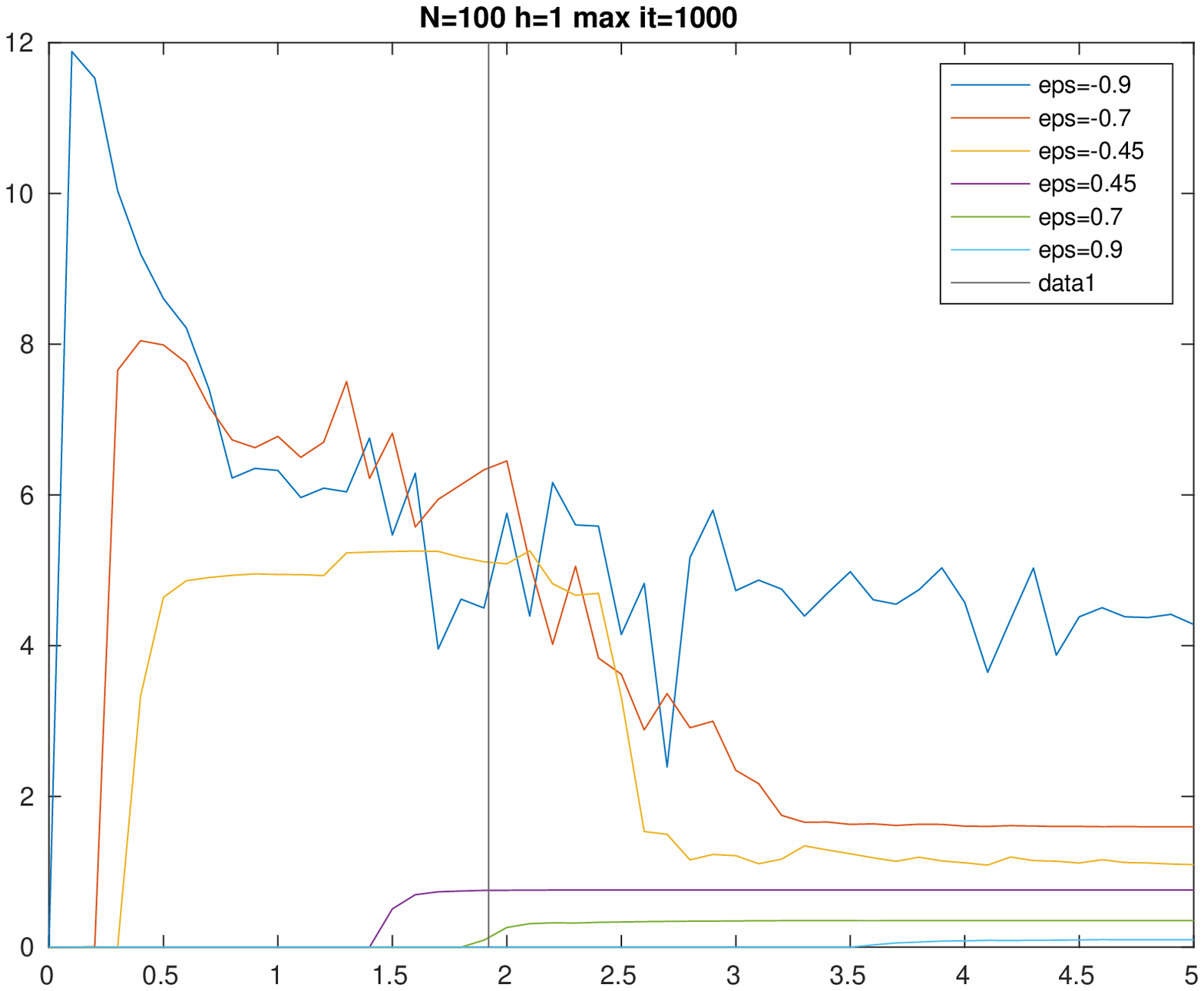}
  \caption{{\tiny} \space  1000 iterations of the $\varepsilon$-Banach algorithm for $N=100, \beta\in [0, 5]$ and $h=1$. We simulated one uniformly in the hypercube distributed start value for $\varepsilon = -0.9, -0.7, -0.45, 0.45, 0.7, 0.9$. We plotted the {\it maximum absolute error} for the given parameters and used the same start value for all $\beta$.}
\end{figure}
%%%%%% 
%%%%% 
\begin{figure}[H]
  \centering
  % \psfrag{2}{\tiny Two Step}
  \includegraphics[width=0.5\textwidth,  trim={0 0 0 0},clip]{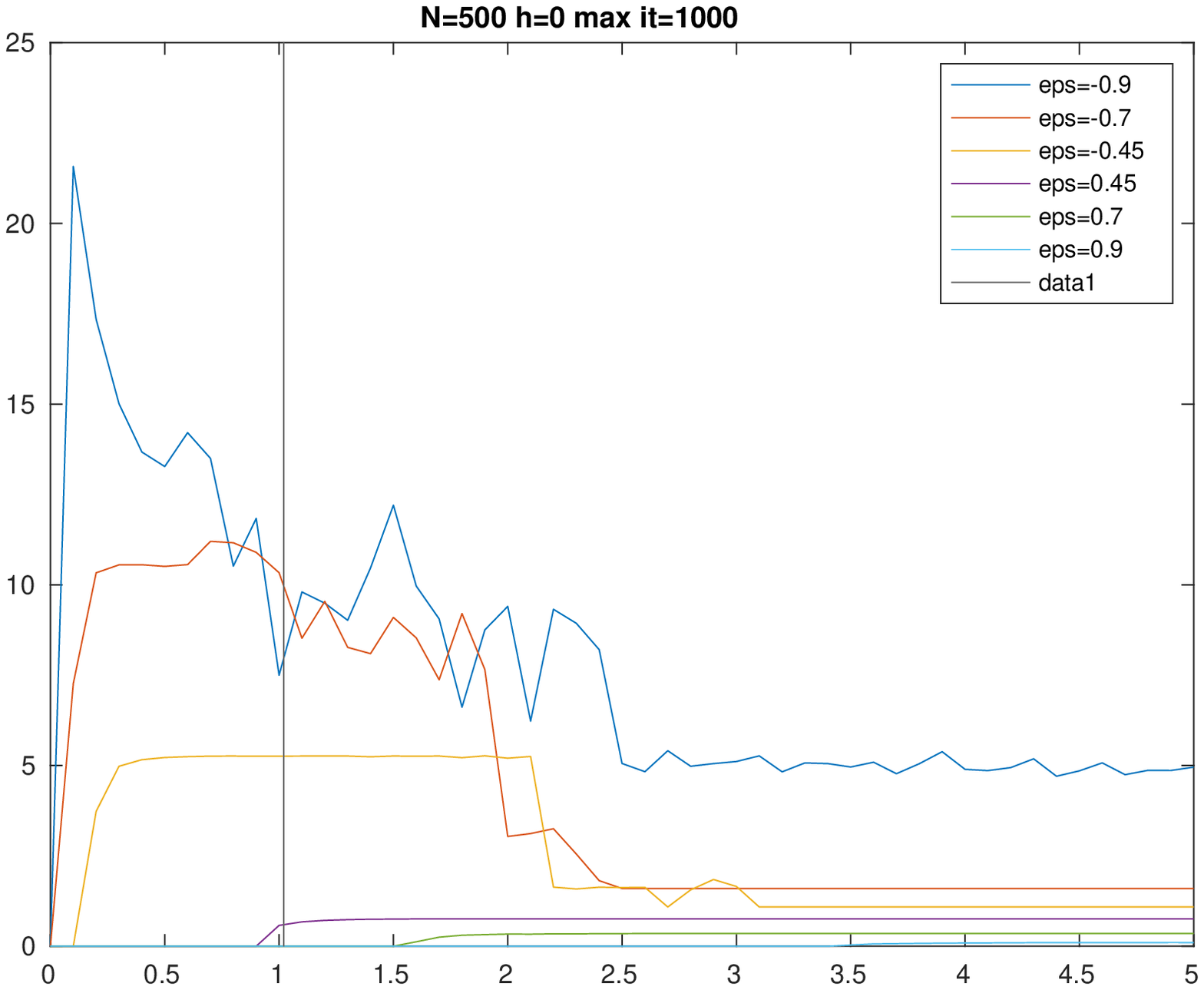}
  \caption{{\tiny} \space  1000 iterations of the $\varepsilon$-Banach algorithm for $N=500, \beta\in [0, 5]$ and $h=0$. We simulated one uniformly in the hypercube distributed start value for $\varepsilon = -0.9, -0.7, -0.45, 0.45, 0.7, 0.9$. We plotted the {\it maximum absolute error} for the given parameters and used the same start value for all $\beta$.}
\end{figure}
\begin{figure}[H]
  \centering
  % \psfrag{2}{\tiny Two Step}
  \includegraphics[width=0.5\textwidth,  trim={0 0 0 0},clip]{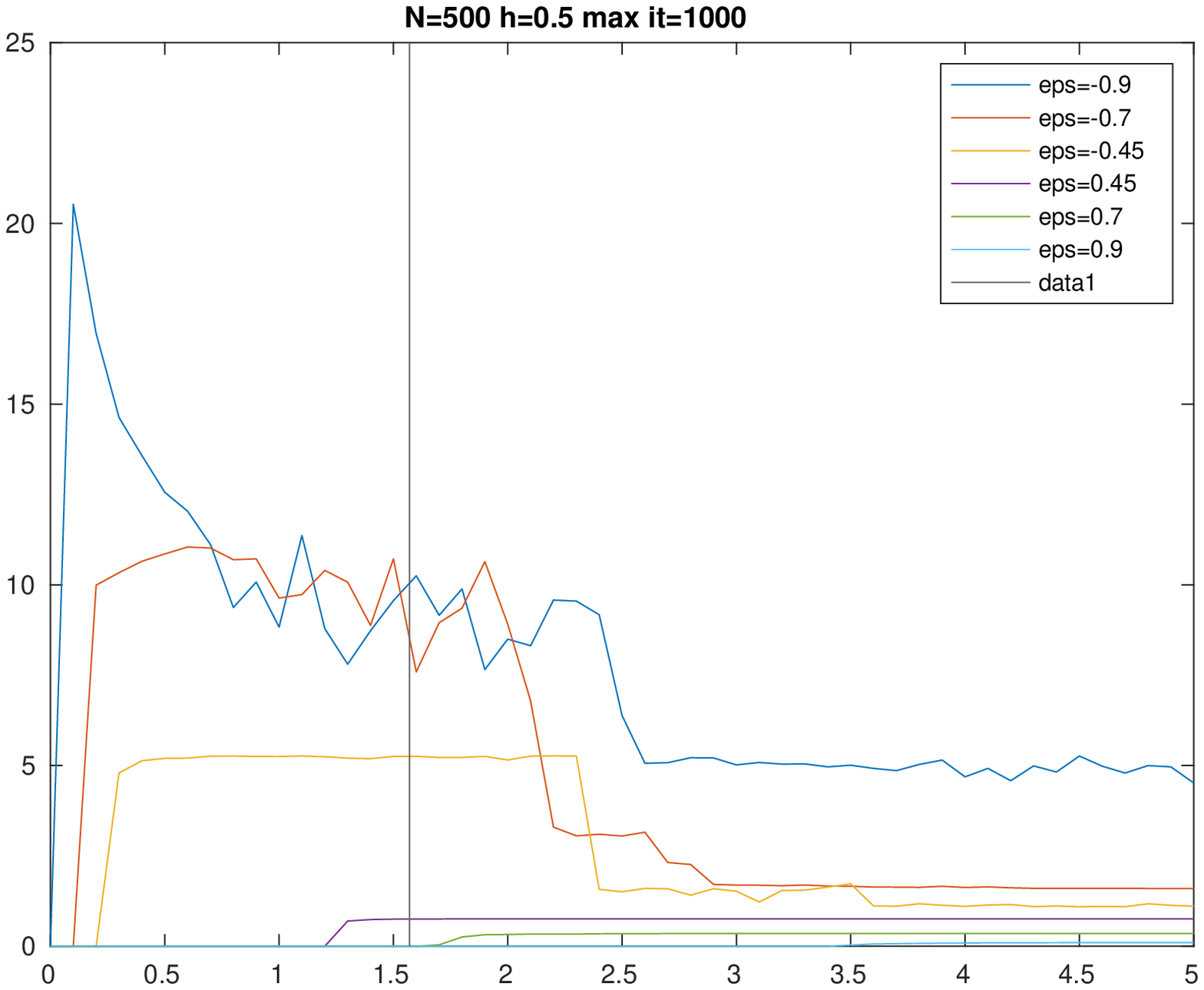}
  \caption{{\tiny} \space  1000 iterations of the $\varepsilon$-Banach algorithm for $N=500, \beta\in [0, 5]$ and $h=0.5$. We simulated one uniformly in the hypercube distributed start value for $\varepsilon = -0.9, -0.7, -0.45, 0.45, 0.7, 0.9$. We plotted the {\it maximum absolute error} for the given parameters and used the same start value for all $\beta$.}
\end{figure}
\begin{figure}[H]
  \centering
  % \psfrag{2}{\tiny Two Step}
  \includegraphics[width=0.5\textwidth,  trim={0 0 0 0},clip]{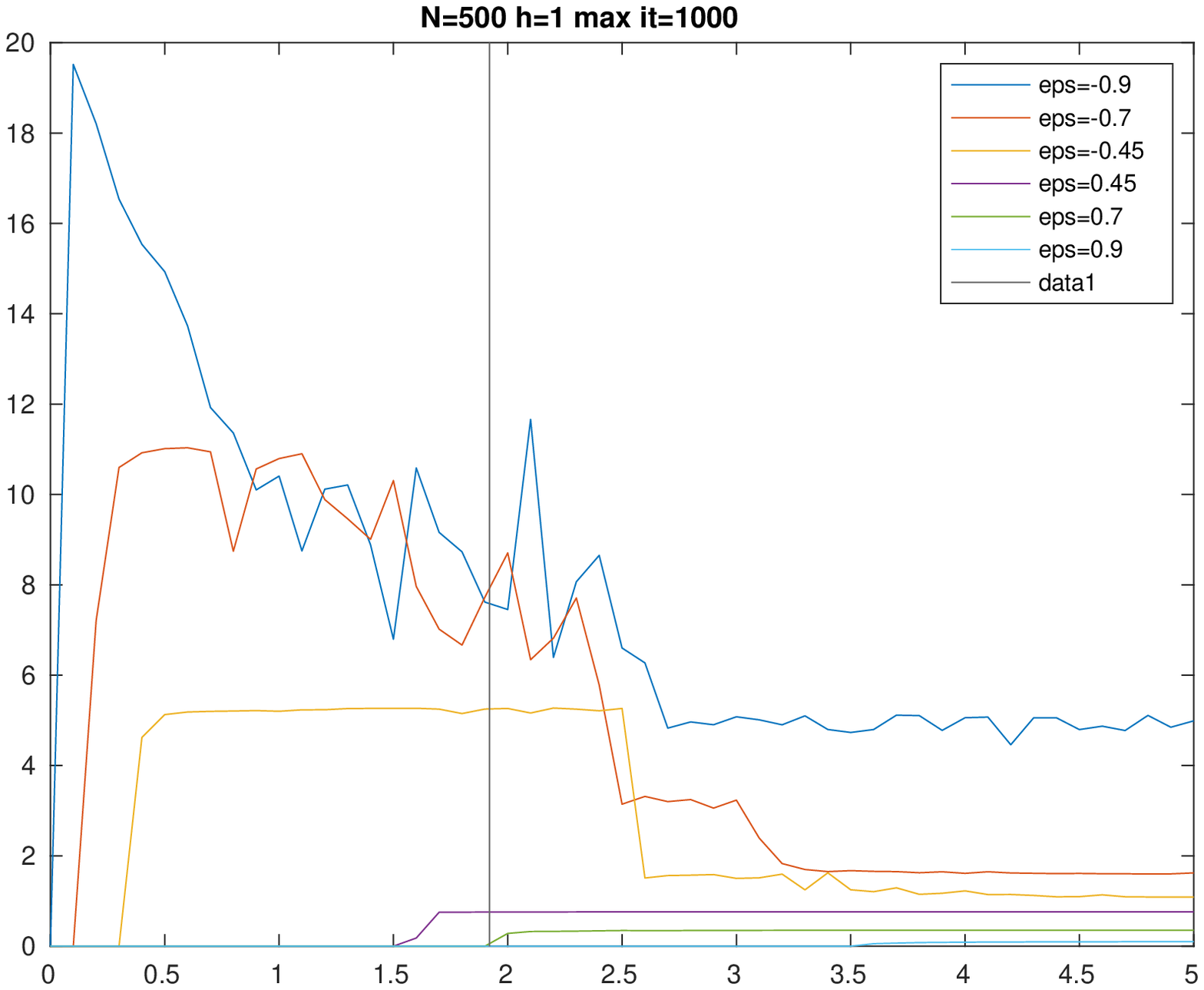}
  \caption{{\tiny} \space  1000 iterations of the $\varepsilon$-Banach algorithm for $N=500, \beta\in [0, 5]$ and $h=1$. We simulated one uniformly in the hypercube distributed start value for $\varepsilon = -0.9, -0.7, -0.45, 0.45, 0.7, 0.9$. We plotted the {\it maximum absolute error} for the given parameters and used the same start value for all $\beta$.}
\end{figure}
%%%%%% 
%%%%% 
\begin{figure}[H]
  \centering
  % \psfrag{2}{\tiny Two Step}
  \includegraphics[width=0.5\textwidth,  trim={0 0 0 0},clip]{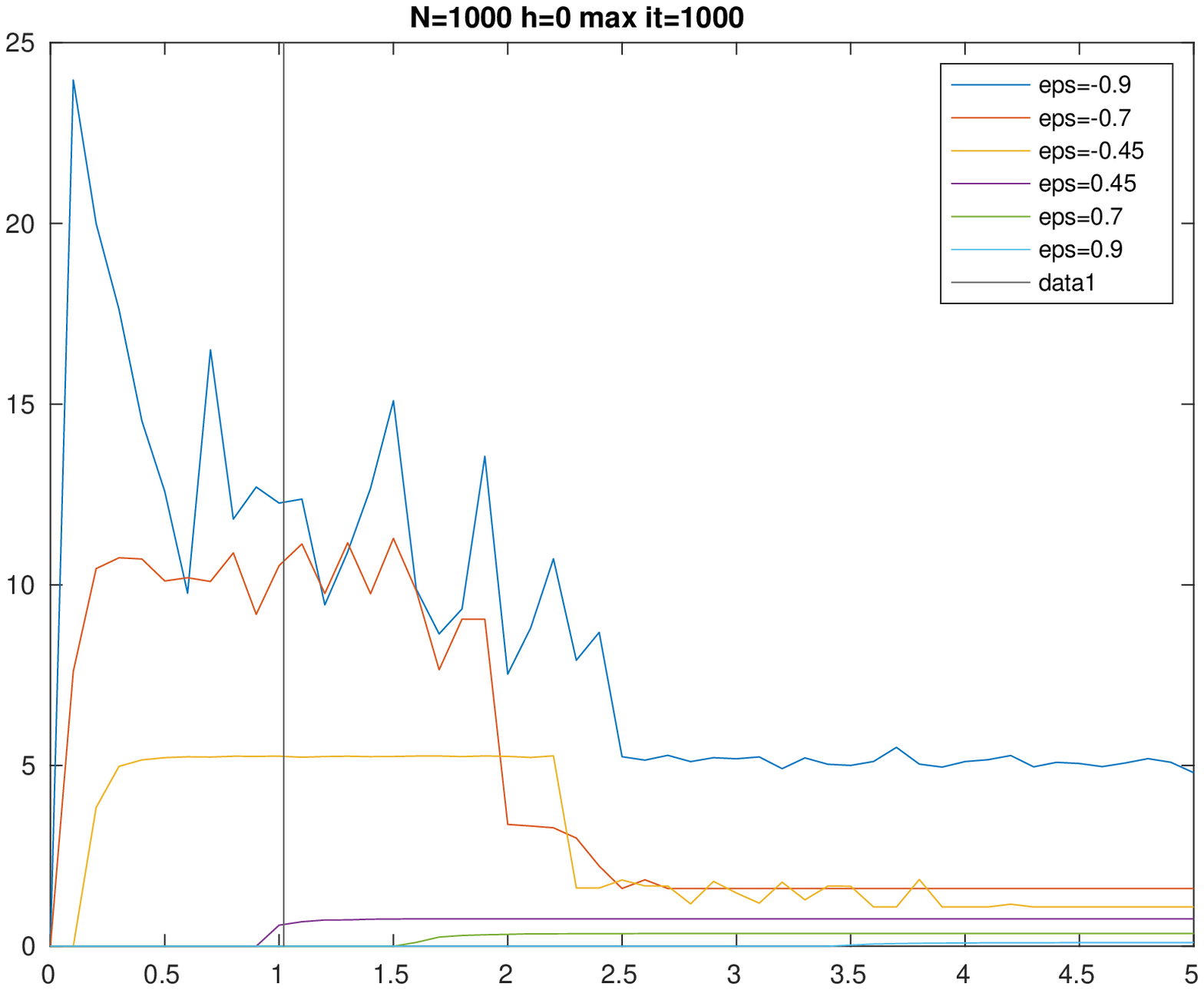}
  \caption{{\tiny} \space  1000 iterations of the $\varepsilon$-Banach algorithm for $N=1000, \beta\in [0, 5]$ and $h=0$. We simulated one uniformly in the hypercube distributed start value for $\varepsilon = -0.9, -0.7, -0.45, 0.45, 0.7, 0.9$. We plotted the {\it maximum absolute error} for the given parameters and used the same start value for all $\beta$.}
\end{figure}
\begin{figure}[H]
  \centering
  % \psfrag{2}{\tiny Two Step}
  \includegraphics[width=0.5\textwidth,  trim={0 0 0 0},clip]{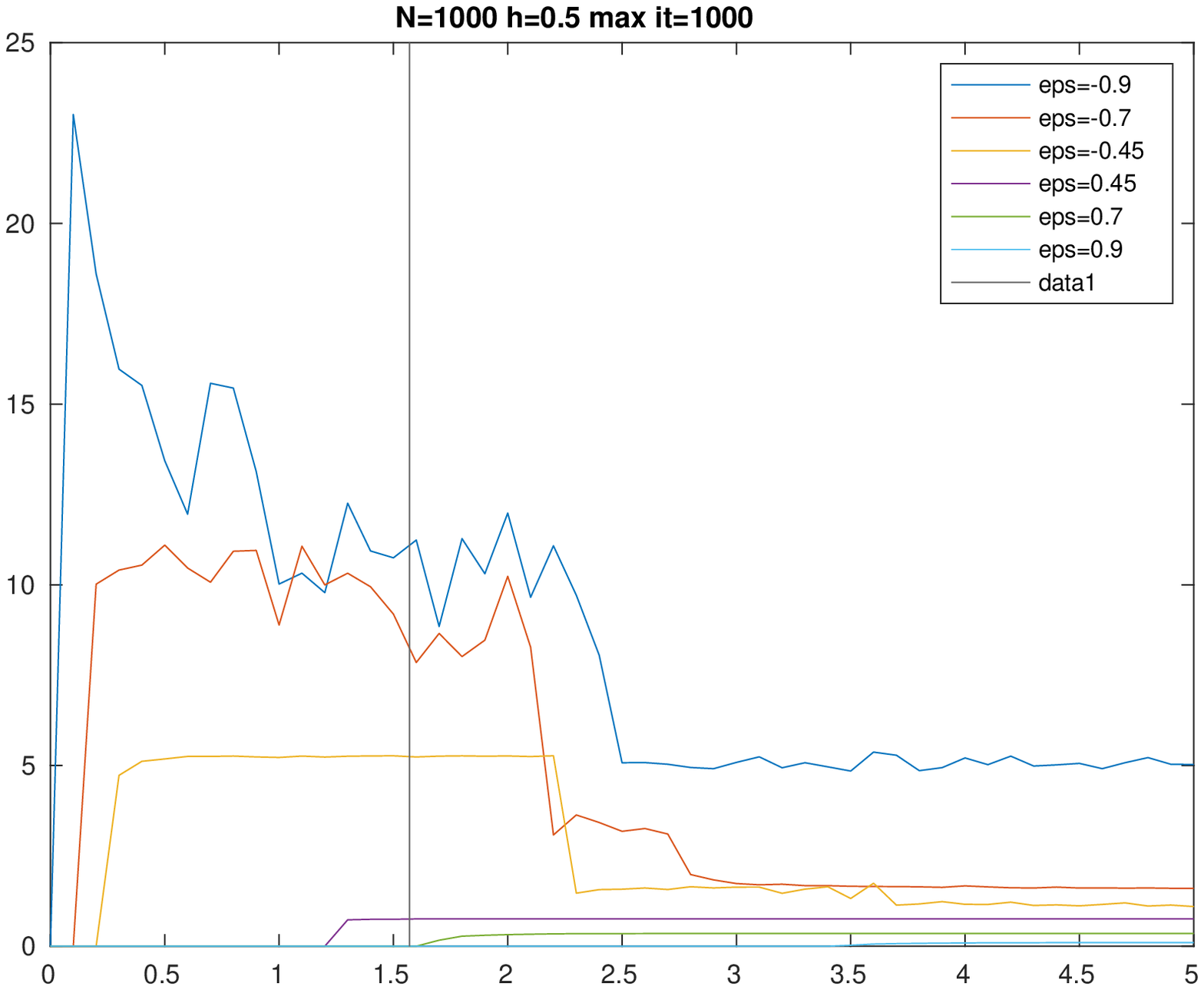}
  \caption{{\tiny} \space  1000 iterations of the $\varepsilon$-Banach algorithm for $N=1000, \beta\in [0, 5]$ and $h=0.5$. We simulated one uniformly in the hypercube distributed start value for $\varepsilon = -0.9, -0.7, -0.45, 0.45, 0.7, 0.9$. We plotted the {\it maximum absolute error} for the given parameters and used the same start value for all $\beta$.}
\end{figure}
\begin{figure}[H]
  \centering
  % \psfrag{2}{\tiny Two Step}
  \includegraphics[width=0.5\textwidth,  trim={0 0 0 0},clip]{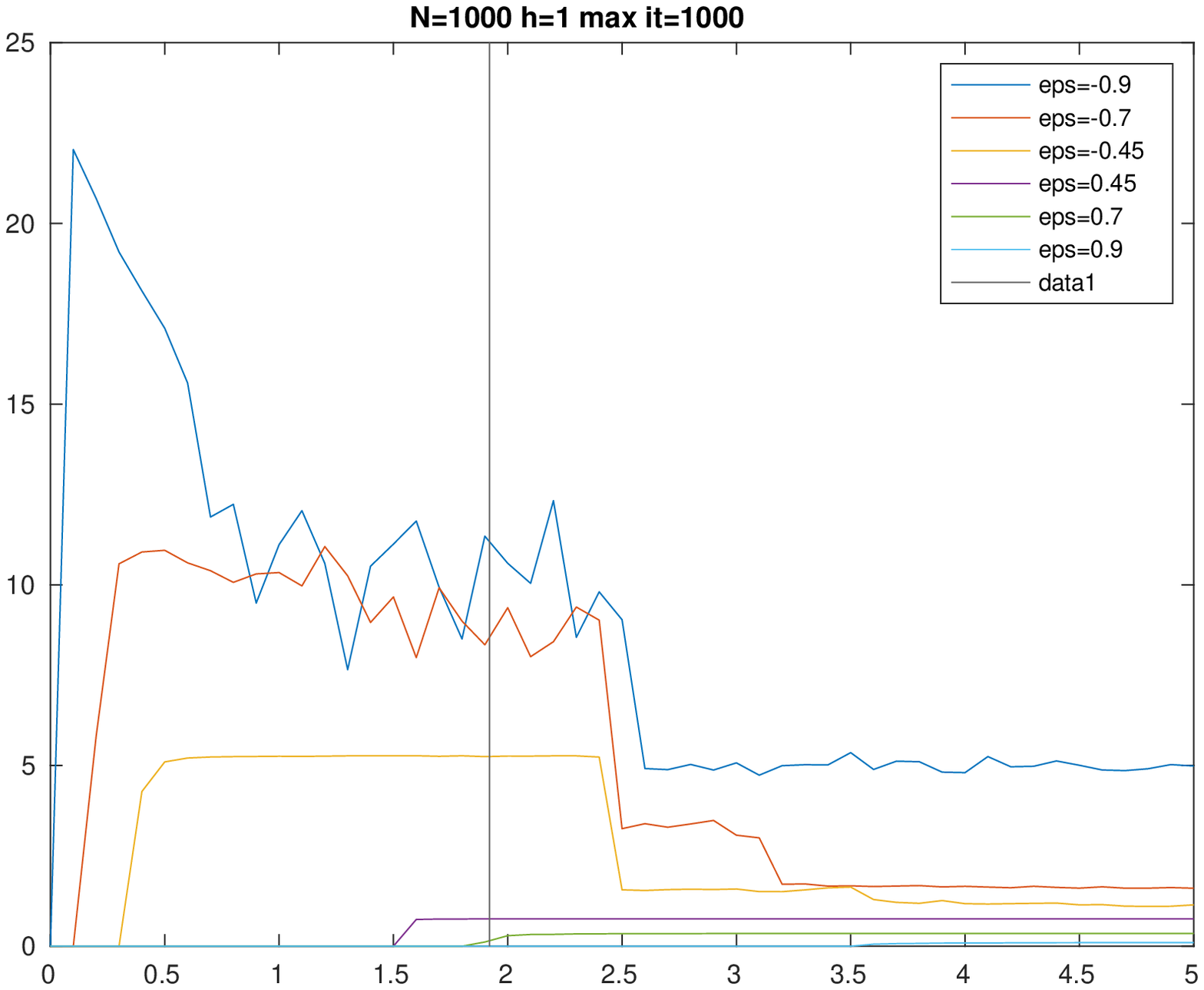}
  \caption{{\tiny} \space  1000 iterations of the $\varepsilon$-Banach algorithm for $N=1000, \beta\in [0, 5]$ and $h=1$. We simulated one uniformly in the hypercube distributed start value for $\varepsilon = -0.9, -0.7, -0.45, 0.45, 0.7, 0.9$. We plotted the {\it maximum absolute error} for the given parameters and used the same start value for all $\beta$.}
\end{figure}
\FloatBarrier\clearpage
\subsection{Uniformly in the corners of the hypercube distributed start values for  $\varepsilon$-Banach}
\subsubsection{{\italic mean squared error}}\FloatBarrier\leavevmode
\begin{figure}[H]
  \centering
  % \psfrag{2}{\tiny Two Step}
  \includegraphics[width=0.5\textwidth,  trim={0 0 0 0},clip]{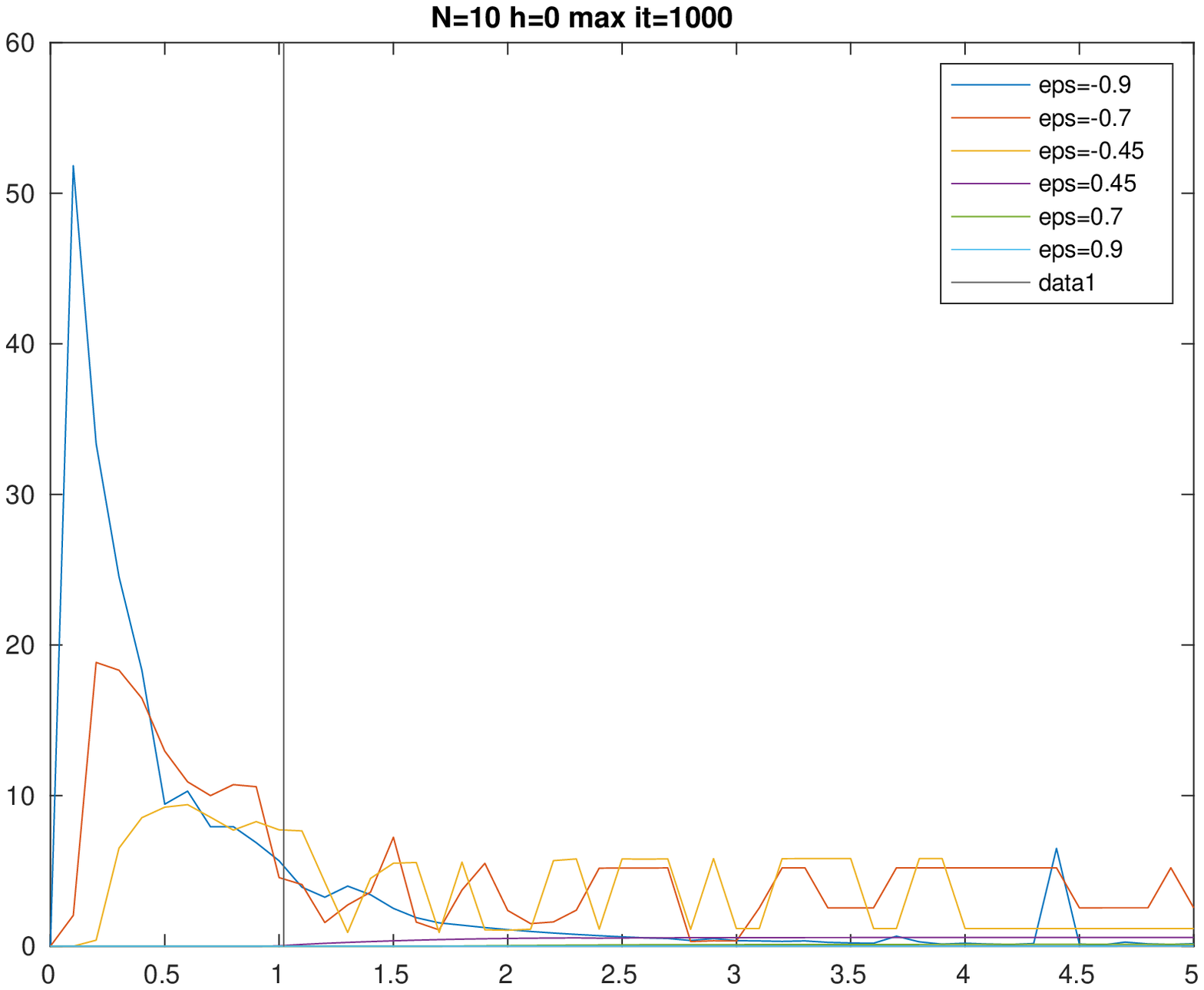}
  \caption{{\tiny} \space  1000 iterations of the $\varepsilon$-Banach algorithm for $N=10, \beta\in [0, 5]$ and $h=0$. We simulated one uniformly in the corners of the hypercube distributed start value for $\varepsilon = -0.9, -0.7, -0.45, 0.45, 0.7, 0.9$. We plotted the {\italic mean squared error} for the given parameters and used the same start value for all $\beta$.}
\end{figure}
\begin{figure}[H]
  \centering
  % \psfrag{2}{\tiny Two Step}
  \includegraphics[width=0.5\textwidth,  trim={0 0 0 0},clip]{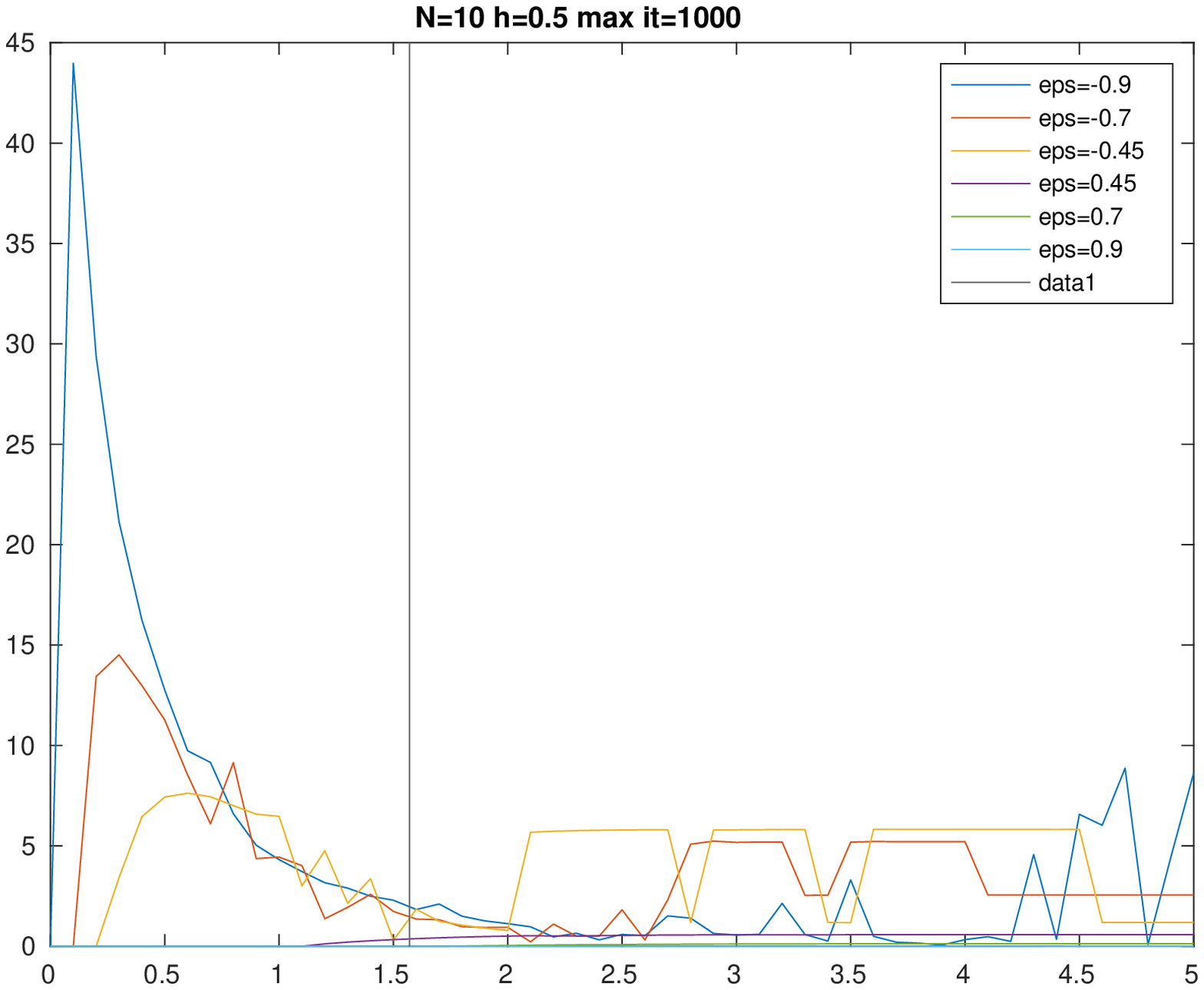}
  \caption{{\tiny} \space  1000 iterations of the $\varepsilon$-Banach algorithm for $N=10, \beta\in [0, 5]$ and $h=0.5$. We simulated one uniformly in the corners of the hypercube distributed start value for $\varepsilon = -0.9, -0.7, -0.45, 0.45, 0.7, 0.9$. We plotted the {\italic mean squared error} for the given parameters and used the same start value for all $\beta$.}
\end{figure}
\begin{figure}[H]
  \centering
  % \psfrag{2}{\tiny Two Step}
  \includegraphics[width=0.5\textwidth,  trim={0 0 0 0},clip]{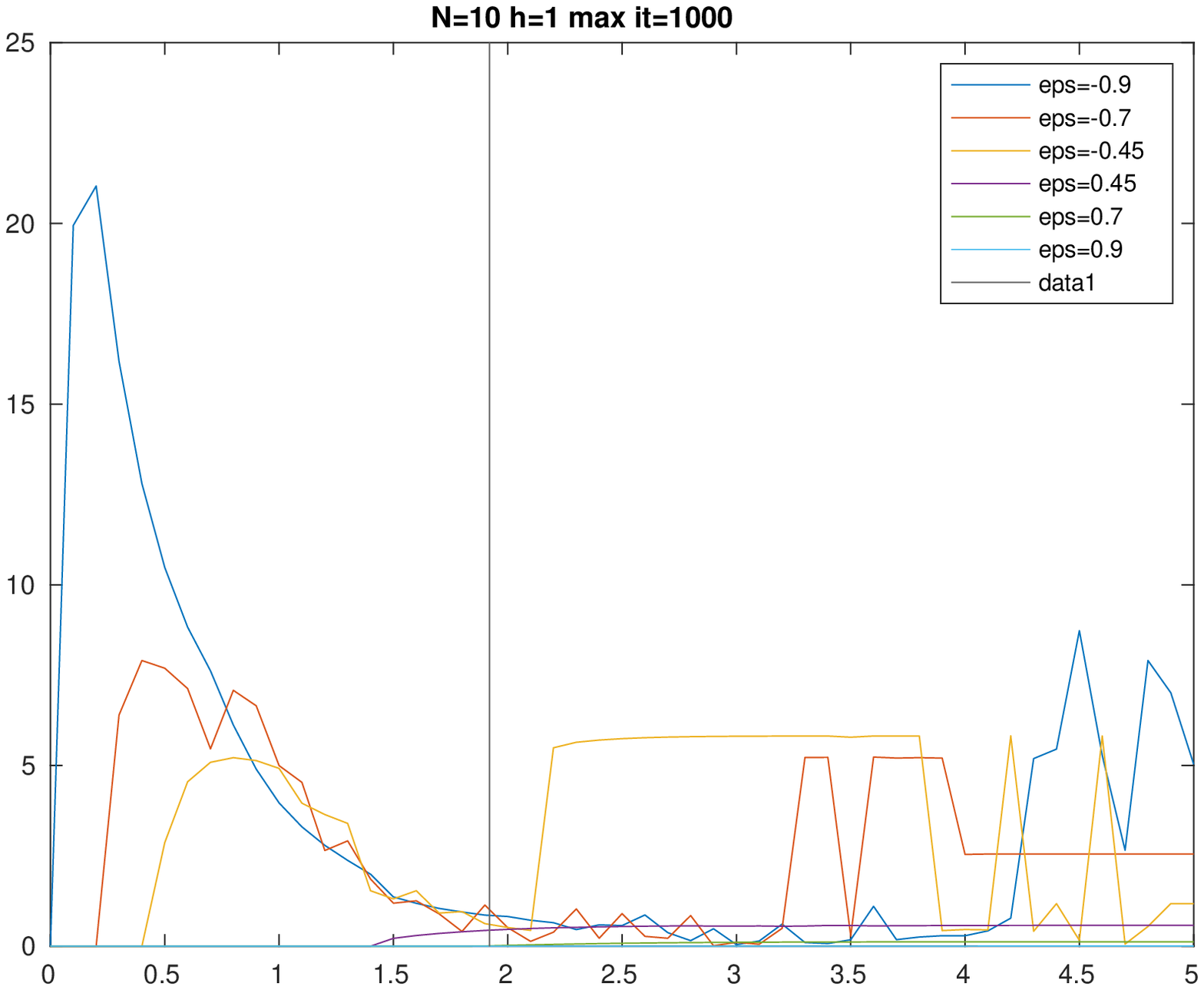}
  \caption{{\tiny} \space  1000 iterations of the $\varepsilon$-Banach algorithm for $N=10, \beta\in [0, 5]$ and $h=1$. We simulated one uniformly in the corners of the hypercube distributed start value for $\varepsilon = -0.9, -0.7, -0.45, 0.45, 0.7, 0.9$. We plotted the {\italic mean squared error} for the given parameters and used the same start value for all $\beta$.}
\end{figure}
%%%%%% 
%%%%% 
\begin{figure}[H]
  \centering
  % \psfrag{2}{\tiny Two Step}
  \includegraphics[width=0.5\textwidth,  trim={0 0 0 0},clip]{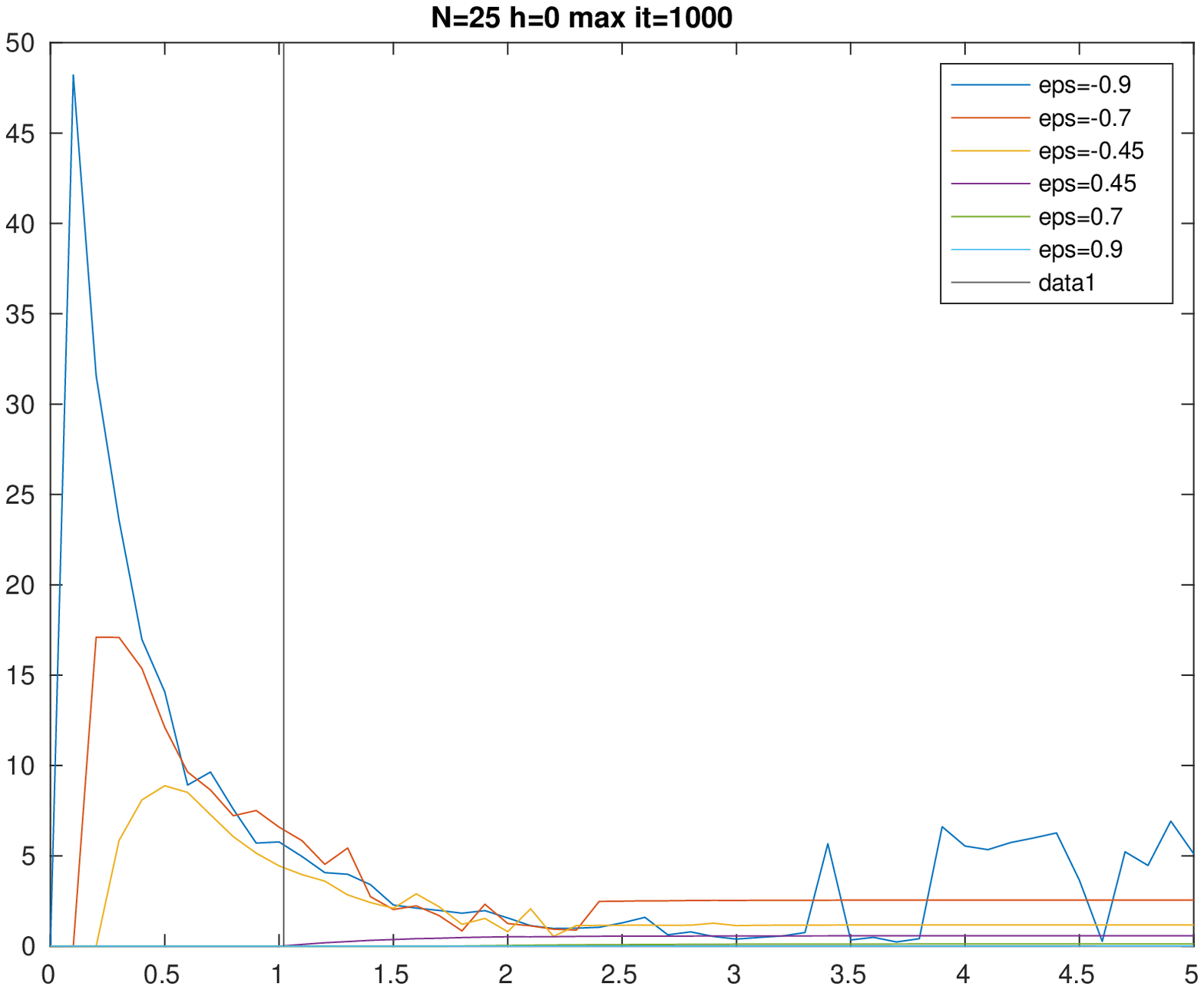}
  \caption{{\tiny} \space  1000 iterations of the $\varepsilon$-Banach algorithm for $N=25, \beta\in [0, 5]$ and $h=0$. We simulated one uniformly in the corners of the hypercube distributed start value for $\varepsilon = -0.9, -0.7, -0.45, 0.45, 0.7, 0.9$. We plotted the {\italic mean squared error} for the given parameters and used the same start value for all $\beta$.}
\end{figure}
\begin{figure}[H]
  \centering
  % \psfrag{2}{\tiny Two Step}
  \includegraphics[width=0.5\textwidth,  trim={0 0 0 0},clip]{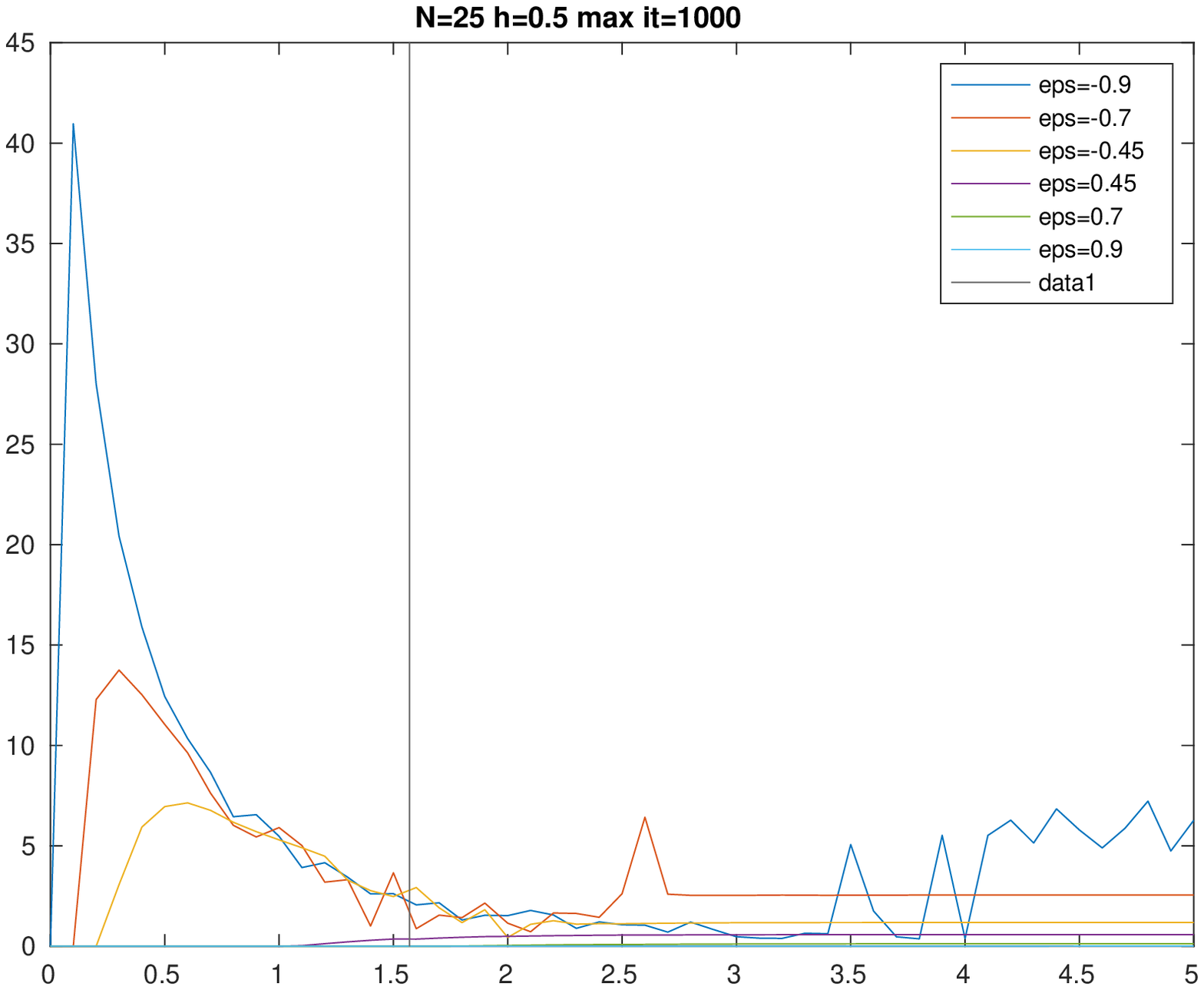}
  \caption{{\tiny} \space  1000 iterations of the $\varepsilon$-Banach algorithm for $N=25, \beta\in [0, 5]$ and $h=0.5$. We simulated one uniformly in the corners of the hypercube distributed start value for $\varepsilon = -0.9, -0.7, -0.45, 0.45, 0.7, 0.9$. We plotted the {\italic mean squared error} for the given parameters and used the same start value for all $\beta$.}
\end{figure}
\begin{figure}[H]
  \centering
  % \psfrag{2}{\tiny Two Step}
  \includegraphics[width=0.5\textwidth,  trim={0 0 0 0},clip]{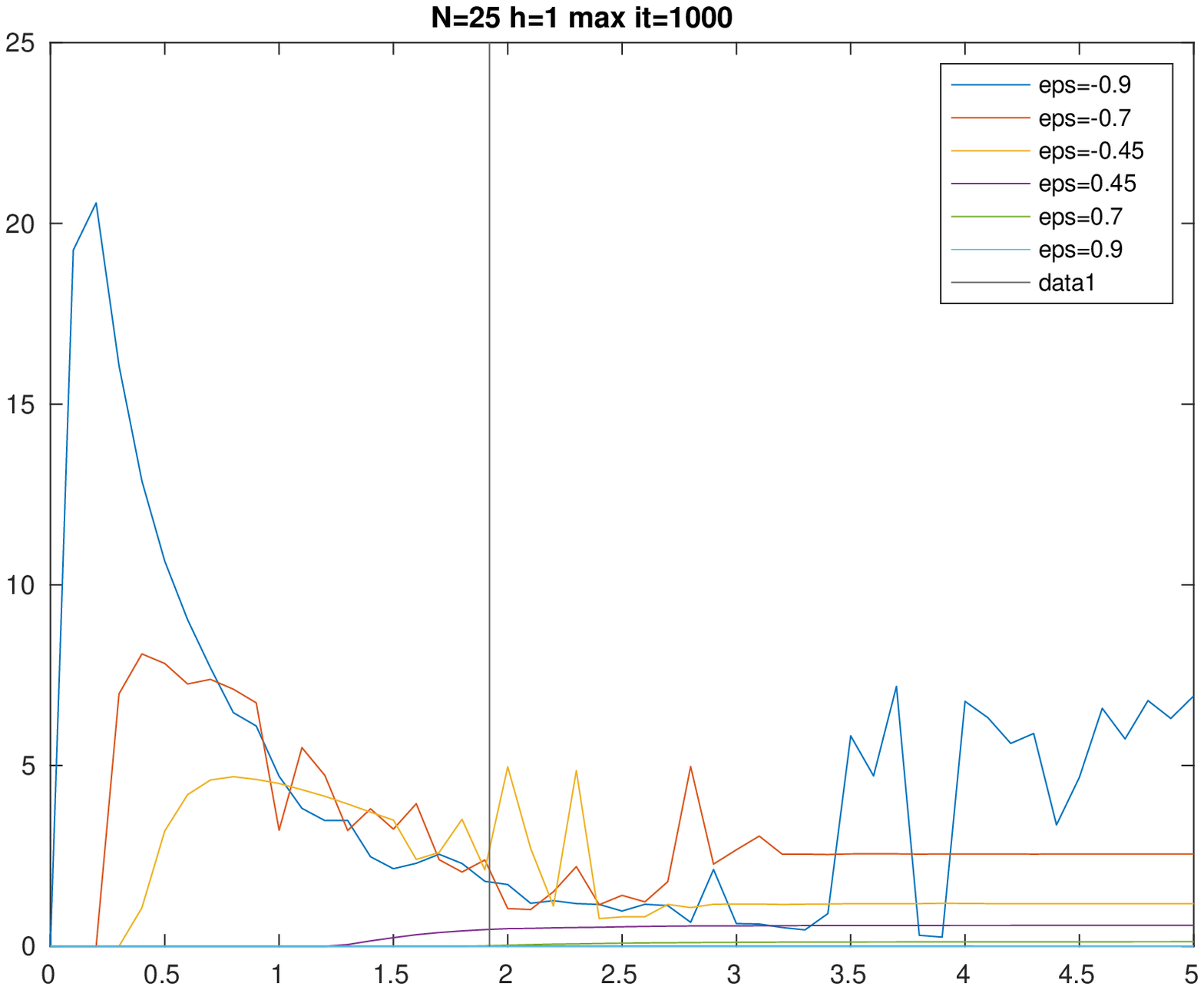}
  \caption{{\tiny} \space  1000 iterations of the $\varepsilon$-Banach algorithm for $N=25, \beta\in [0, 5]$ and $h=1$. We simulated one uniformly in the corners of the hypercube distributed start value for $\varepsilon = -0.9, -0.7, -0.45, 0.45, 0.7, 0.9$. We plotted the {\italic mean squared error} for the given parameters and used the same start value for all $\beta$.}
\end{figure}
%%%%%% 
%%%%% 
\begin{figure}[H]
  \centering
  % \psfrag{2}{\tiny Two Step}
  \includegraphics[width=0.5\textwidth,  trim={0 0 0 0},clip]{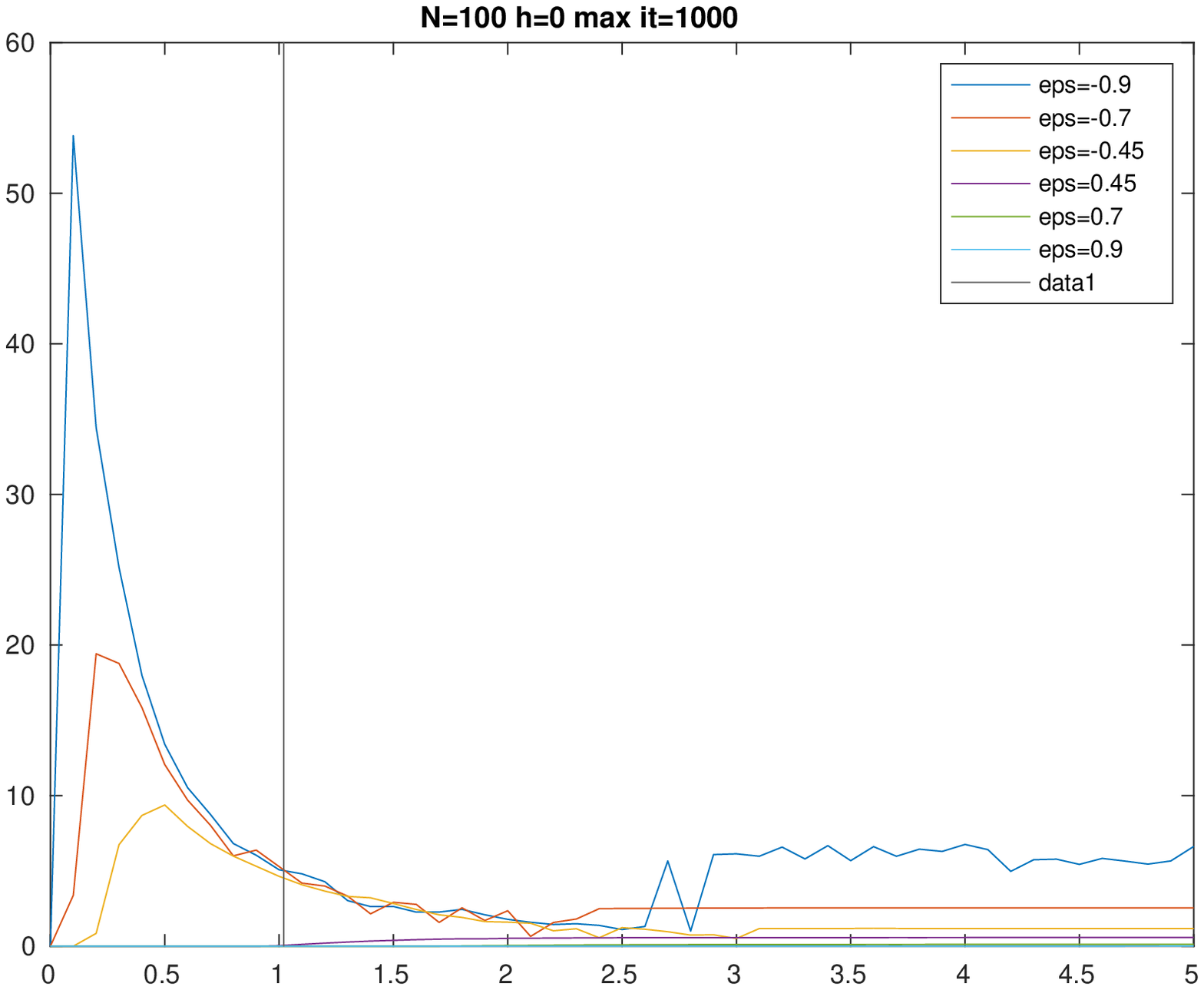}
  \caption{{\tiny} \space  1000 iterations of the $\varepsilon$-Banach algorithm for $N=100, \beta\in [0, 5]$ and $h=0$. We simulated one uniformly in the corners of the hypercube distributed start value for $\varepsilon = -0.9, -0.7, -0.45, 0.45, 0.7, 0.9$. We plotted the {\italic mean squared error} for the given parameters and used the same start value for all $\beta$.}
\end{figure}
\begin{figure}[H]
  \centering
  % \psfrag{2}{\tiny Two Step}
  \includegraphics[width=0.5\textwidth,  trim={0 0 0 0},clip]{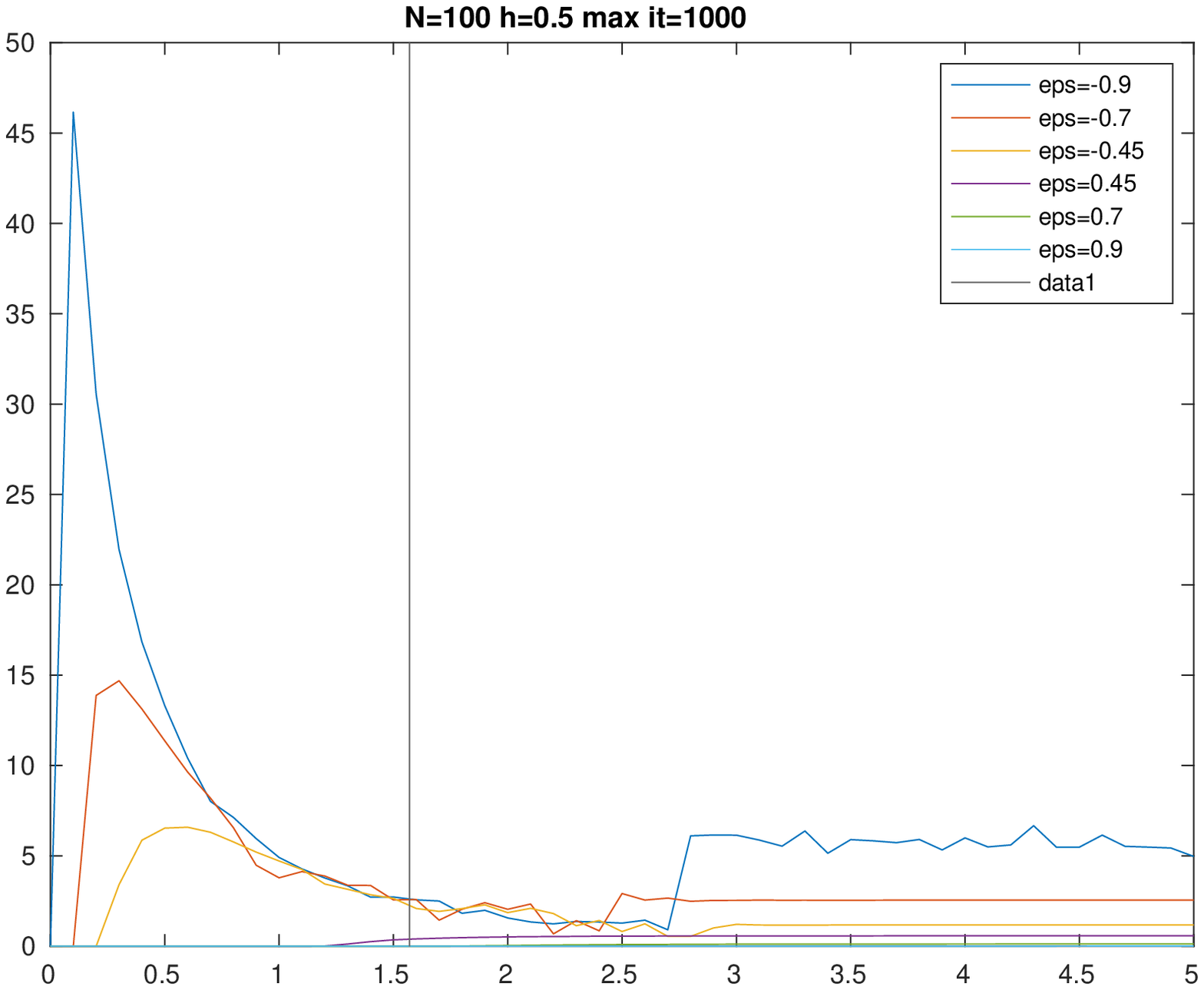}
  \caption{{\tiny} \space  1000 iterations of the $\varepsilon$-Banach algorithm for $N=100, \beta\in [0, 5]$ and $h=0.5$. We simulated one uniformly in the corners of the hypercube distributed start value for $\varepsilon = -0.9, -0.7, -0.45, 0.45, 0.7, 0.9$. We plotted the {\italic mean squared error} for the given parameters and used the same start value for all $\beta$.}
\end{figure}
\begin{figure}[H]
  \centering
  % \psfrag{2}{\tiny Two Step}
  \includegraphics[width=0.5\textwidth,  trim={0 0 0 0},clip]{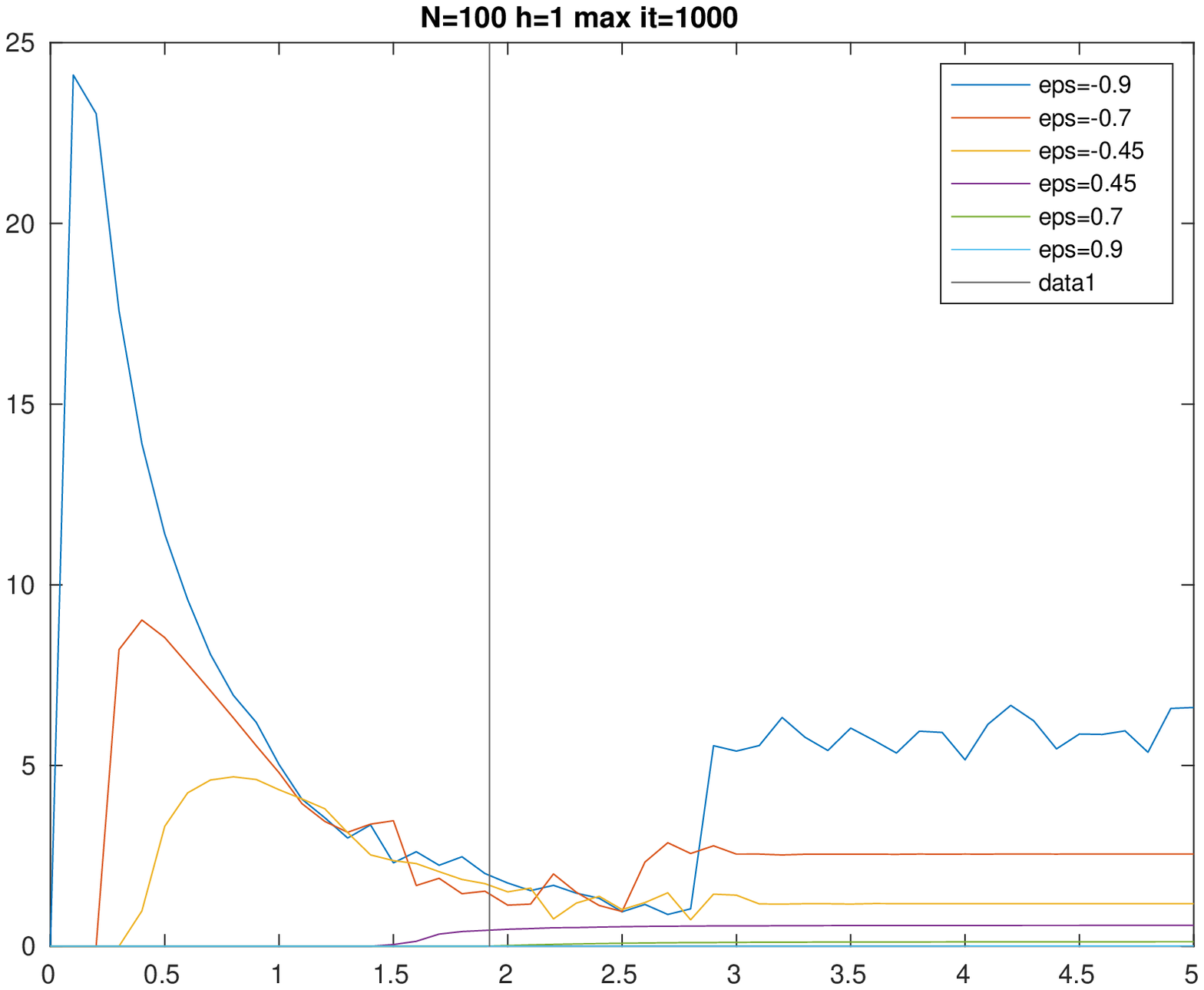}
  \caption{{\tiny} \space  1000 iterations of the $\varepsilon$-Banach algorithm for $N=100, \beta\in [0, 5]$ and $h=1$. We simulated one uniformly in the corners of the hypercube distributed start value for $\varepsilon = -0.9, -0.7, -0.45, 0.45, 0.7, 0.9$. We plotted the {\italic mean squared error} for the given parameters and used the same start value for all $\beta$.}
\end{figure}
%%%%%% 
%%%%% 
\begin{figure}[H]
  \centering
  % \psfrag{2}{\tiny Two Step}
  \includegraphics[width=0.5\textwidth,  trim={0 0 0 0},clip]{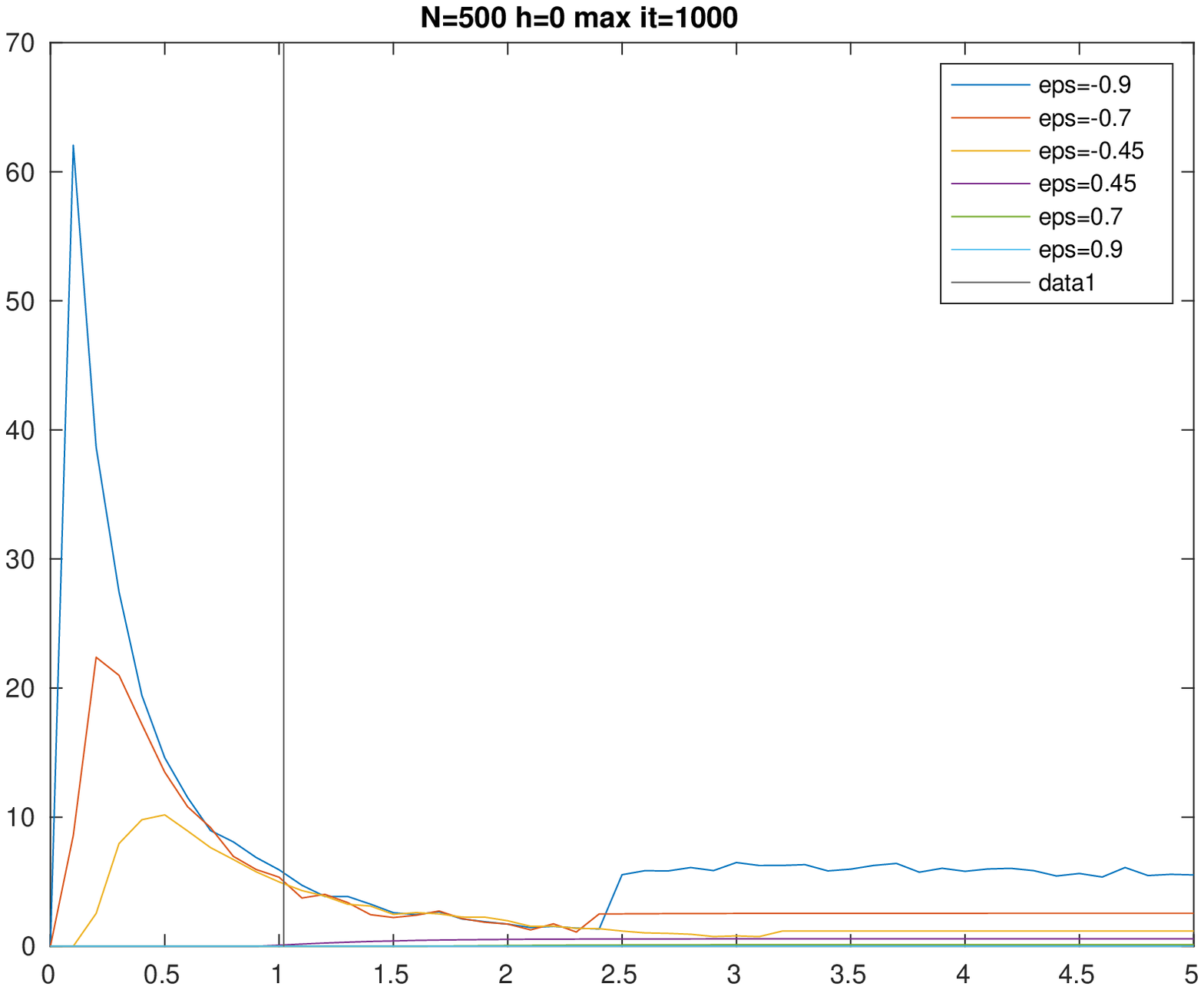}
  \caption{{\tiny} \space  1000 iterations of the $\varepsilon$-Banach algorithm for $N=500, \beta\in [0, 5]$ and $h=0$. We simulated one uniformly in the corners of the hypercube distributed start value for $\varepsilon = -0.9, -0.7, -0.45, 0.45, 0.7, 0.9$. We plotted the {\italic mean squared error} for the given parameters and used the same start value for all $\beta$.}
\end{figure}
\begin{figure}[H]
  \centering
  % \psfrag{2}{\tiny Two Step}
  \includegraphics[width=0.5\textwidth,  trim={0 0 0 0},clip]{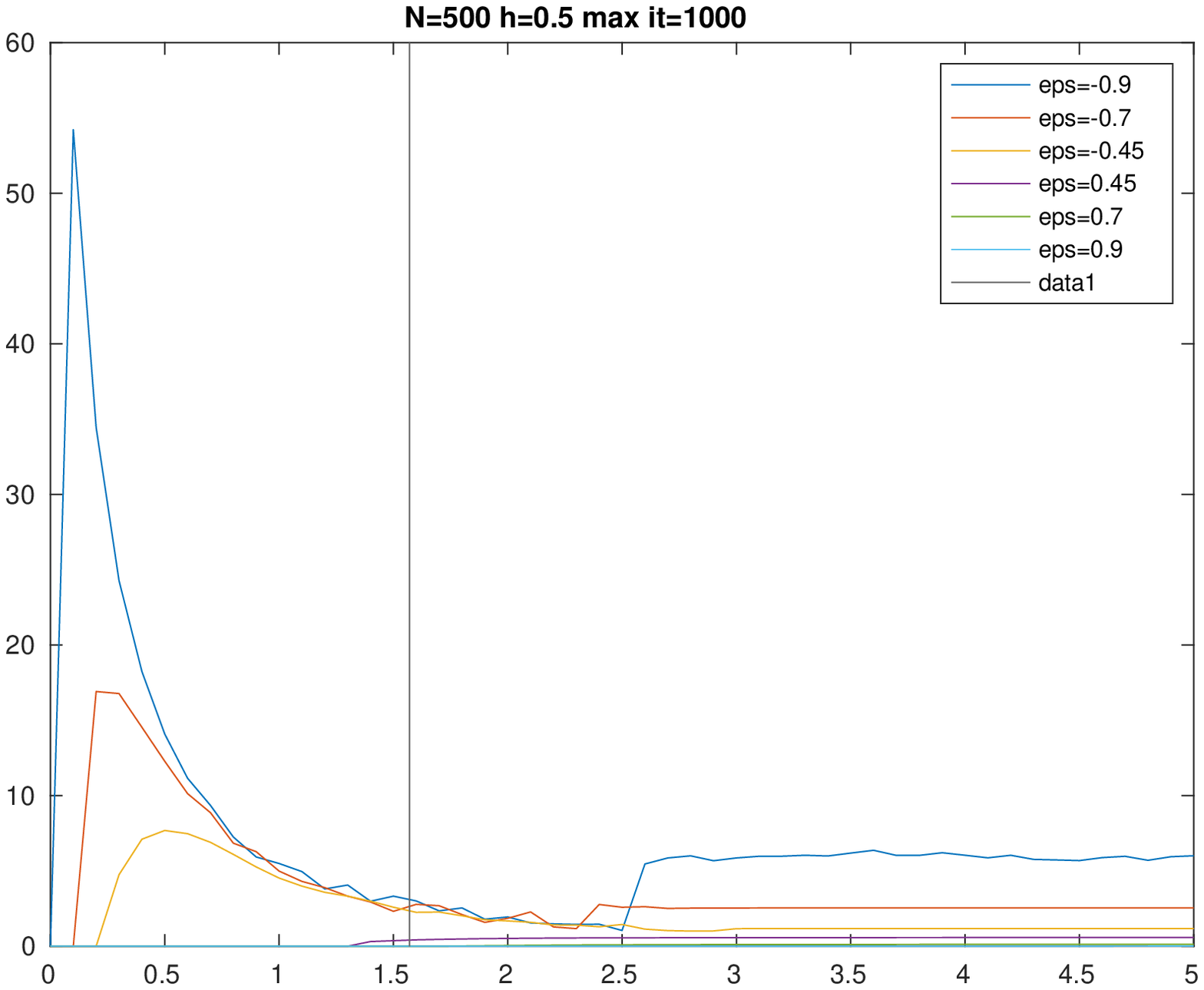}
  \caption{{\tiny} \space  1000 iterations of the $\varepsilon$-Banach algorithm for $N=500, \beta\in [0, 5]$ and $h=0.5$. We simulated one uniformly in the corners of the hypercube distributed start value for $\varepsilon = -0.9, -0.7, -0.45, 0.45, 0.7, 0.9$. We plotted the {\italic mean squared error} for the given parameters and used the same start value for all $\beta$.}
\end{figure}
\begin{figure}[H]
  \centering
  % \psfrag{2}{\tiny Two Step}
  \includegraphics[width=0.5\textwidth,  trim={0 0 0 0},clip]{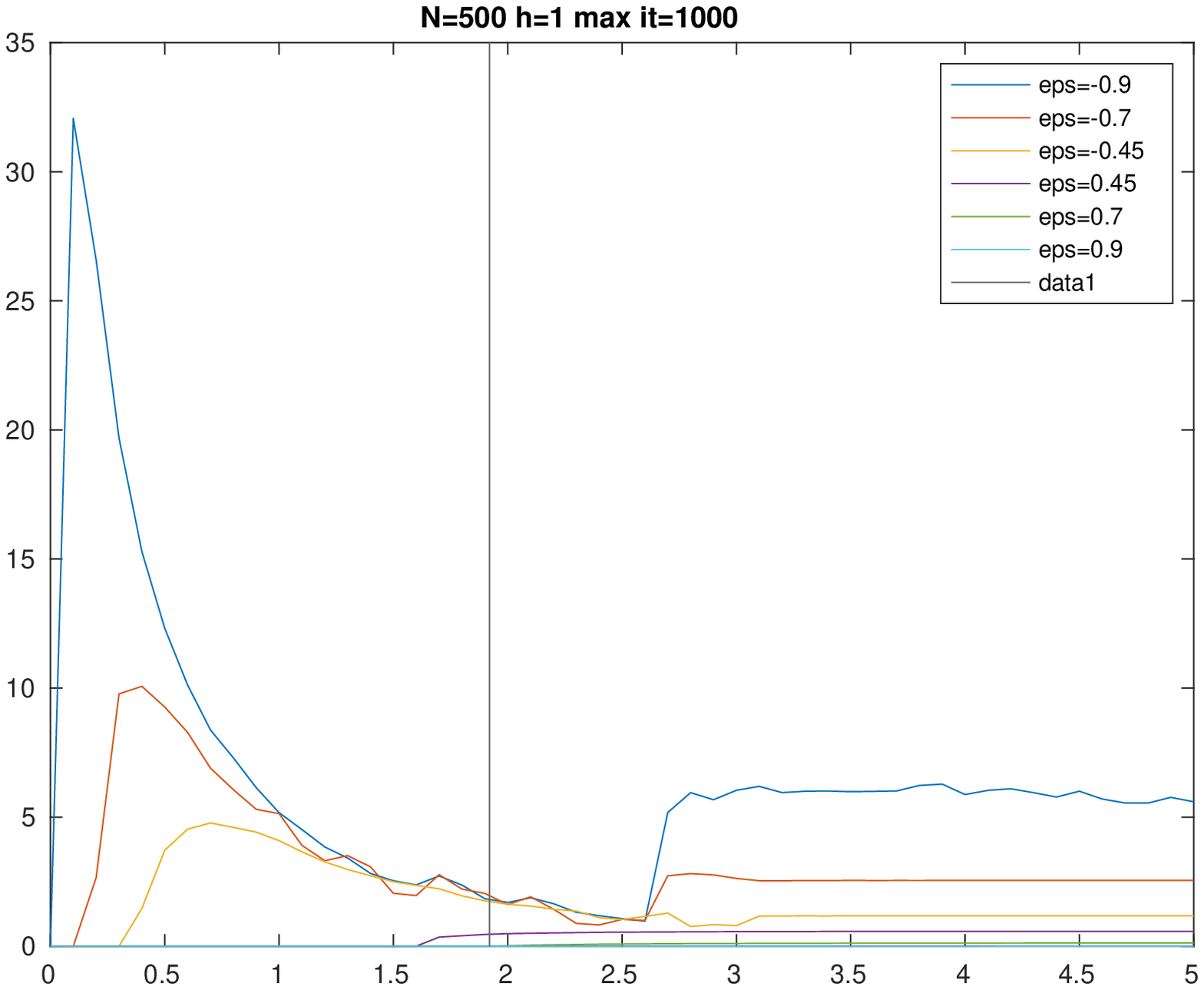}
  \caption{{\tiny} \space  1000 iterations of the $\varepsilon$-Banach algorithm for $N=500, \beta\in [0, 5]$ and $h=1$. We simulated one uniformly in the corners of the hypercube distributed start value for $\varepsilon = -0.9, -0.7, -0.45, 0.45, 0.7, 0.9$. We plotted the {\italic mean squared error} for the given parameters and used the same start value for all $\beta$.}
\end{figure}
%%%%%% 
%%%%% 
\begin{figure}[H]
  \centering
  % \psfrag{2}{\tiny Two Step}
  \includegraphics[width=0.5\textwidth,  trim={0 0 0 0},clip]{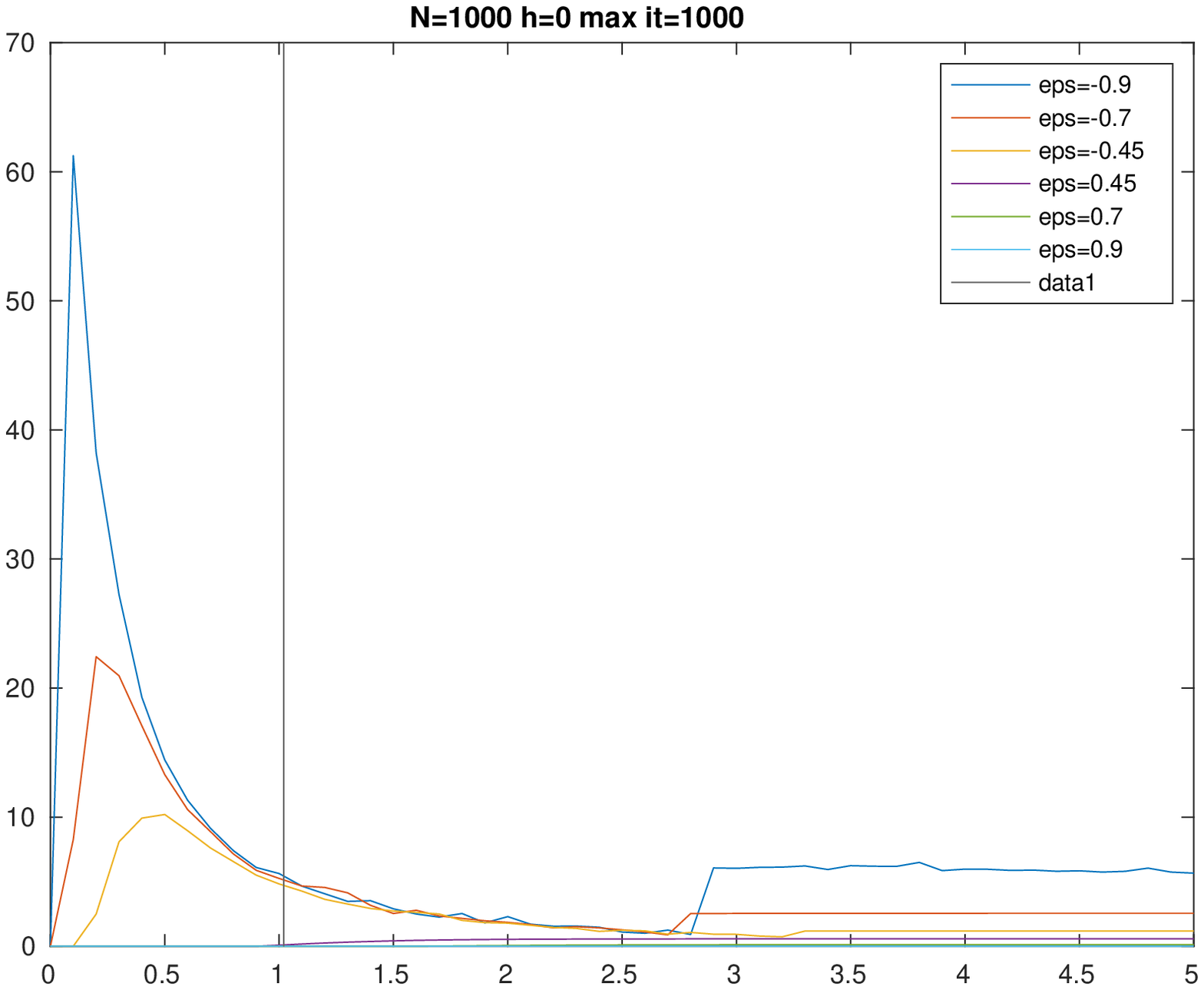}
  \caption{{\tiny} \space  1000 iterations of the $\varepsilon$-Banach algorithm for $N=1000, \beta\in [0, 5]$ and $h=0$. We simulated one uniformly in the corners of the hypercube distributed start value for $\varepsilon = -0.9, -0.7, -0.45, 0.45, 0.7, 0.9$. We plotted the {\italic mean squared error} for the given parameters and used the same start value for all $\beta$.}
\end{figure}
\begin{figure}[H]
  \centering
  % \psfrag{2}{\tiny Two Step}
  \includegraphics[width=0.5\textwidth,  trim={0 0 0 0},clip]{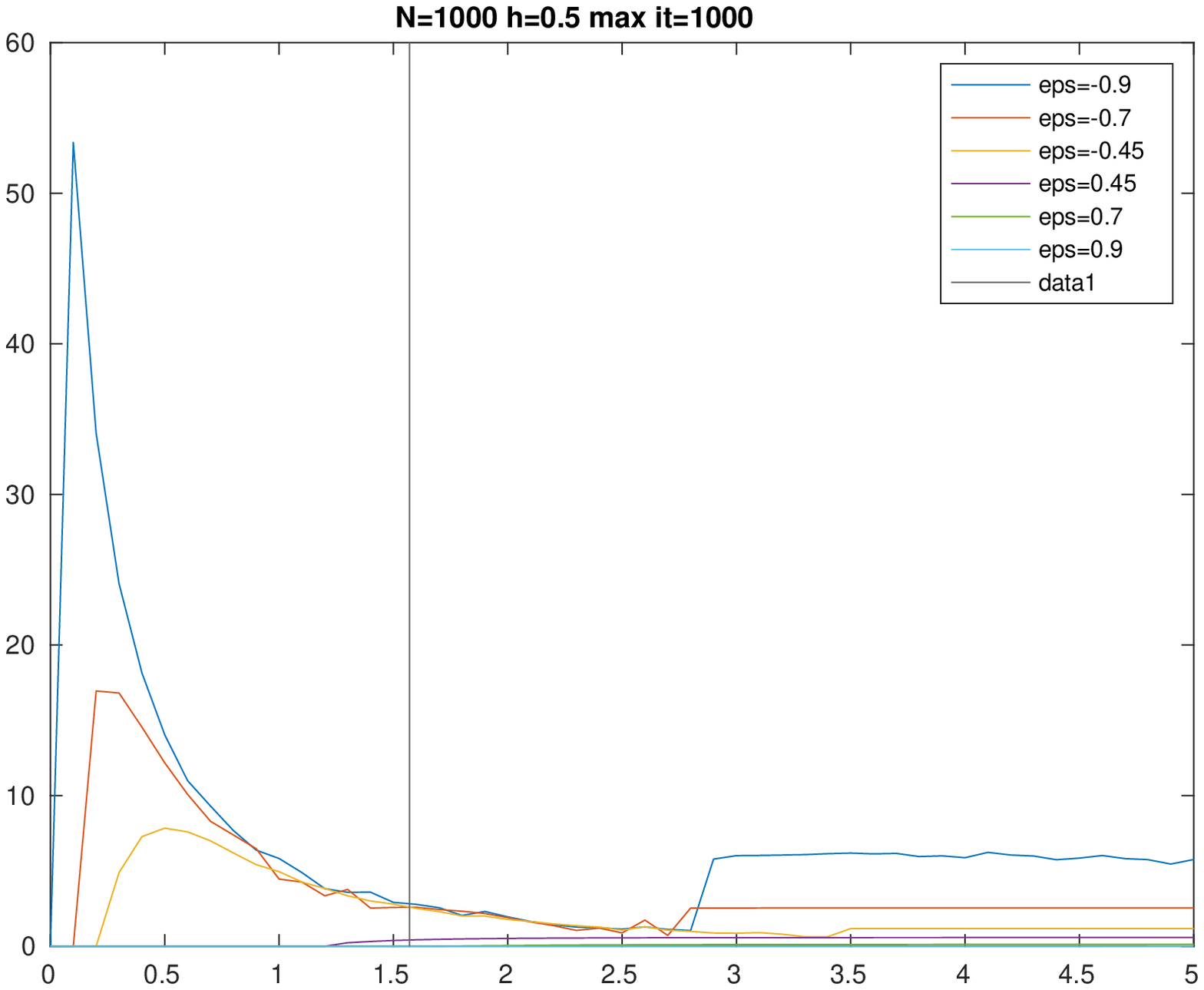}
  \caption{{\tiny} \space  1000 iterations of the $\varepsilon$-Banach algorithm for $N=1000, \beta\in [0, 5]$ and $h=0.5$. We simulated one uniformly in the corners of the hypercube distributed start value for $\varepsilon = -0.9, -0.7, -0.45, 0.45, 0.7, 0.9$. We plotted the {\italic mean squared error} for the given parameters and used the same start value for all $\beta$.}
\end{figure}
\begin{figure}[H]
  \centering
  % \psfrag{2}{\tiny Two Step}
  \includegraphics[width=0.5\textwidth,  trim={0 0 0 0},clip]{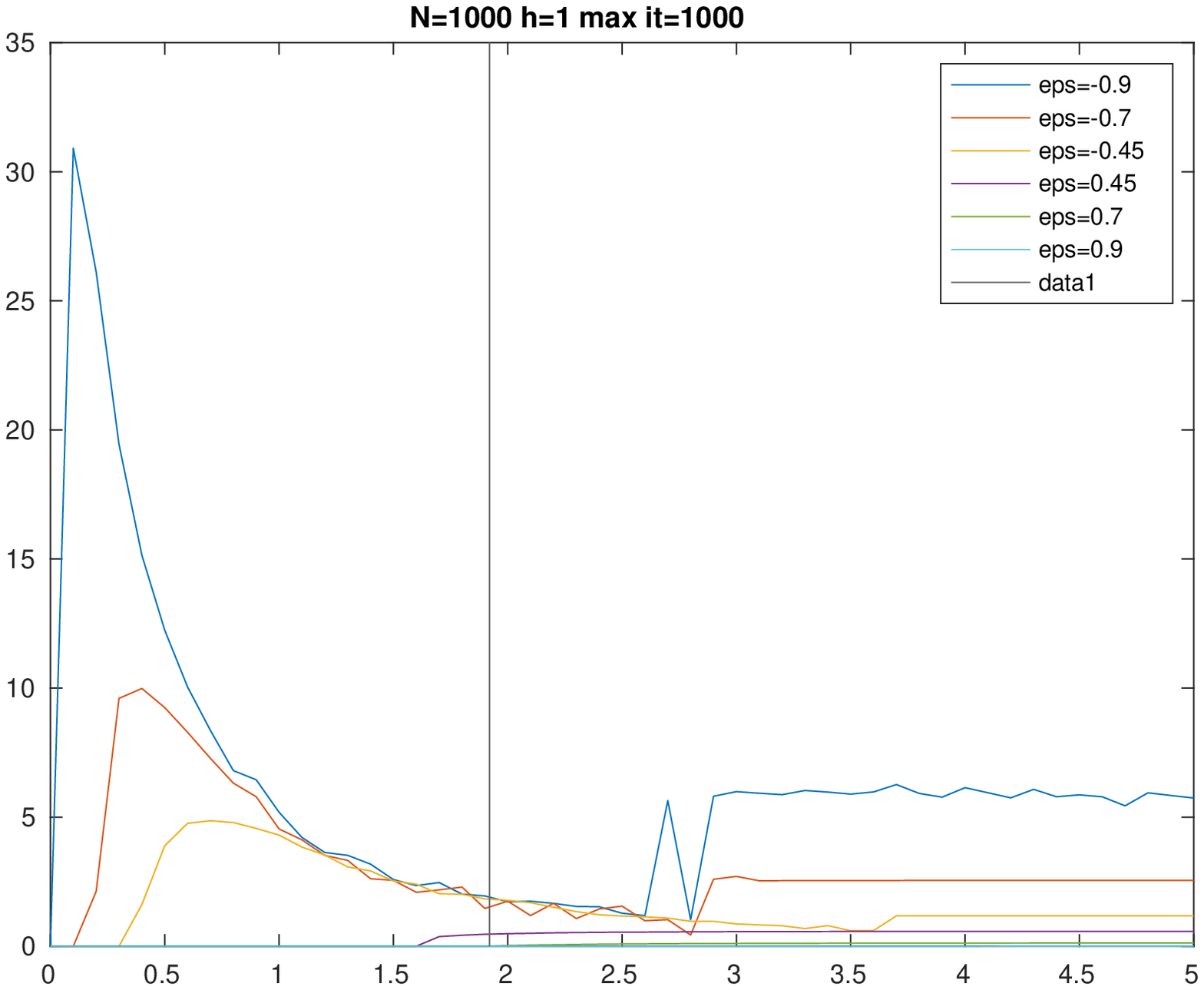}
  \caption{{\tiny} \space  1000 iterations of the $\varepsilon$-Banach algorithm for $N=1000, \beta\in [0, 5]$ and $h=1$. We simulated one uniformly in the corners of the hypercube distributed start value for $\varepsilon = -0.9, -0.7, -0.45, 0.45, 0.7, 0.9$. We plotted the {\italic mean squared error} for the given parameters and used the same start value for all $\beta$.}
\end{figure}
\FloatBarrier\clearpage
\subsubsection{Using the {\it maximum absolute error}}\leavevmode
\begin{figure}[H]
  \centering
  % \psfrag{2}{\tiny Two Step}
  \includegraphics[width=0.5\textwidth,  trim={0 0 0 0},clip]{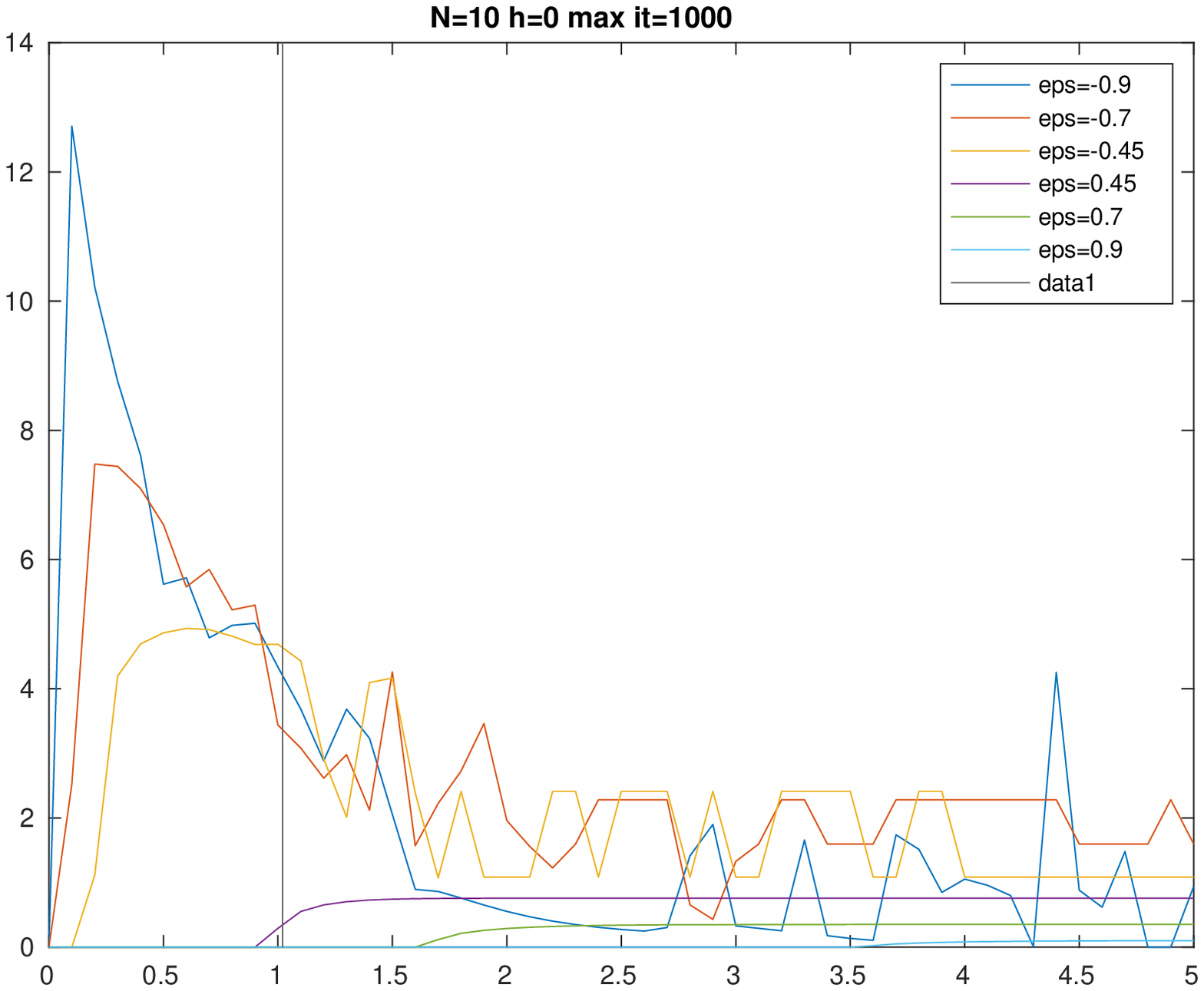}
  \caption{{\tiny} \space  1000 iterations of the $\varepsilon$-Banach algorithm for $N=10, \beta\in [0, 5]$ and $h=0$. We simulated one uniformly in the corners of the hypercube distributed start value for $\varepsilon = -0.9, -0.7, -0.45, 0.45, 0.7, 0.9$. We plotted the {\it maximum absolute error} for the given parameters and used the same start value for all $\beta$.}
\end{figure}
\begin{figure}[H]
  \centering
  % \psfrag{2}{\tiny Two Step}
  \includegraphics[width=0.5\textwidth,  trim={0 0 0 0},clip]{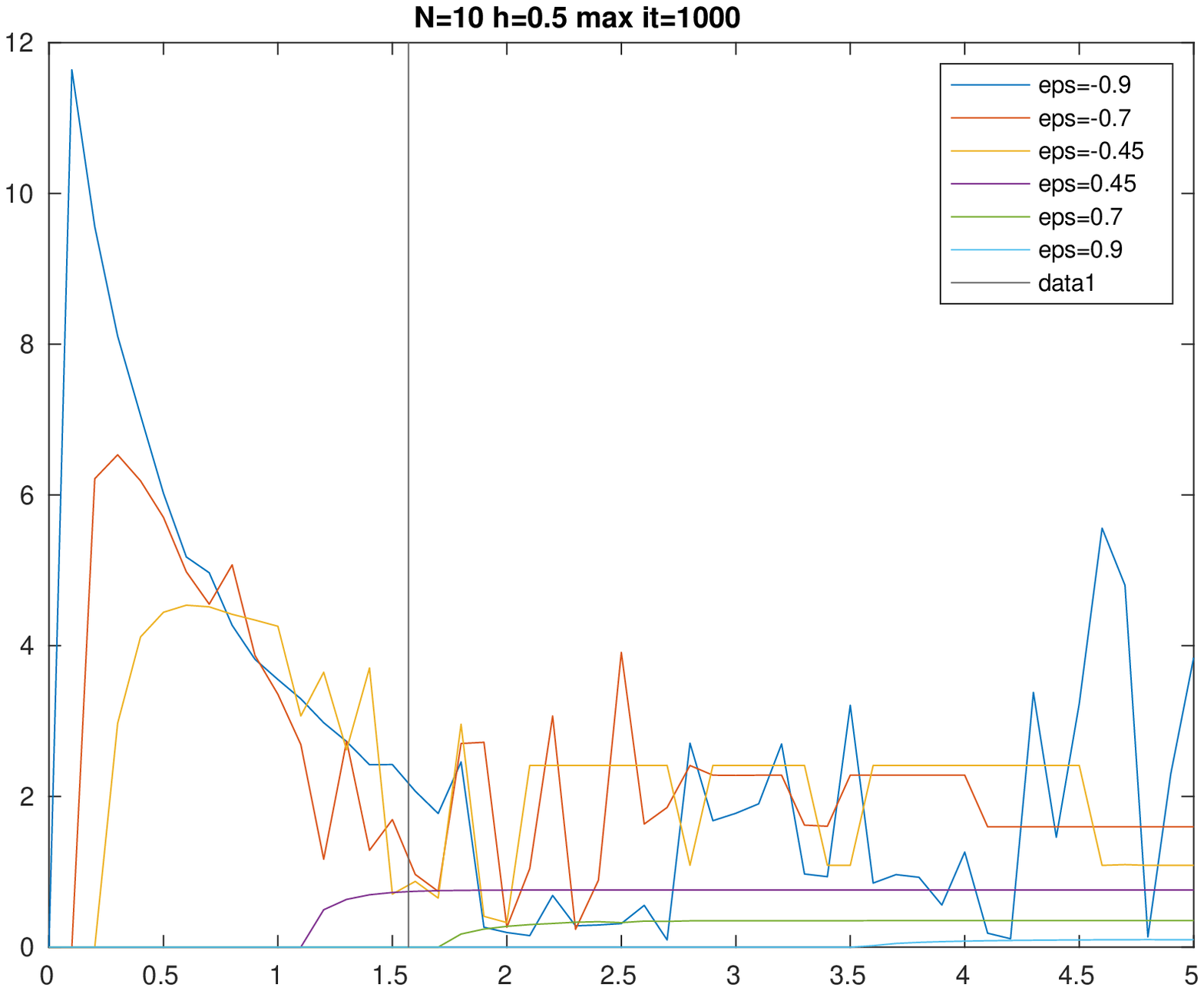}
  \caption{{\tiny} \space  1000 iterations of the $\varepsilon$-Banach algorithm for $N=10, \beta\in [0, 5]$ and $h=0.5$. We simulated one uniformly in the corners of the hypercube distributed start value for $\varepsilon = -0.9, -0.7, -0.45, 0.45, 0.7, 0.9$. We plotted the {\it maximum absolute error} for the given parameters and used the same start value for all $\beta$.}
\end{figure}
\begin{figure}[H]
  \centering
  % \psfrag{2}{\tiny Two Step}
  \includegraphics[width=0.5\textwidth,  trim={0 0 0 0},clip]{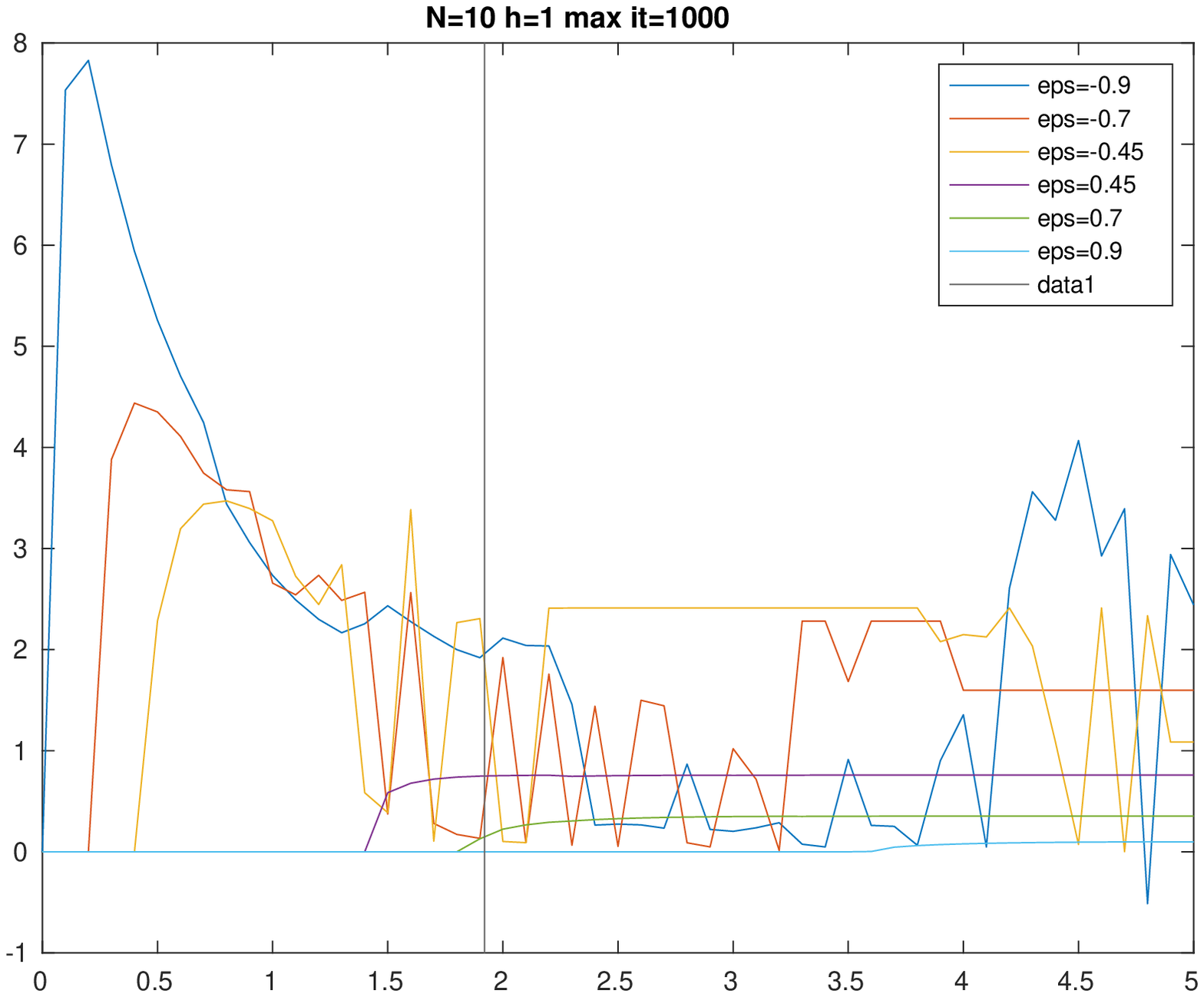}
  \caption{{\tiny} \space  1000 iterations of the $\varepsilon$-Banach algorithm for $N=10, \beta\in [0, 5]$ and $h=1$. We simulated one uniformly in the corners of the hypercube distributed start value for $\varepsilon = -0.9, -0.7, -0.45, 0.45, 0.7, 0.9$. We plotted the {\it maximum absolute error} for the given parameters and used the same start value for all $\beta$.}
\end{figure}
%%%%%% 
%%%%% 
\begin{figure}[H]
  \centering
  % \psfrag{2}{\tiny Two Step}
  \includegraphics[width=0.5\textwidth,  trim={0 0 0 0},clip]{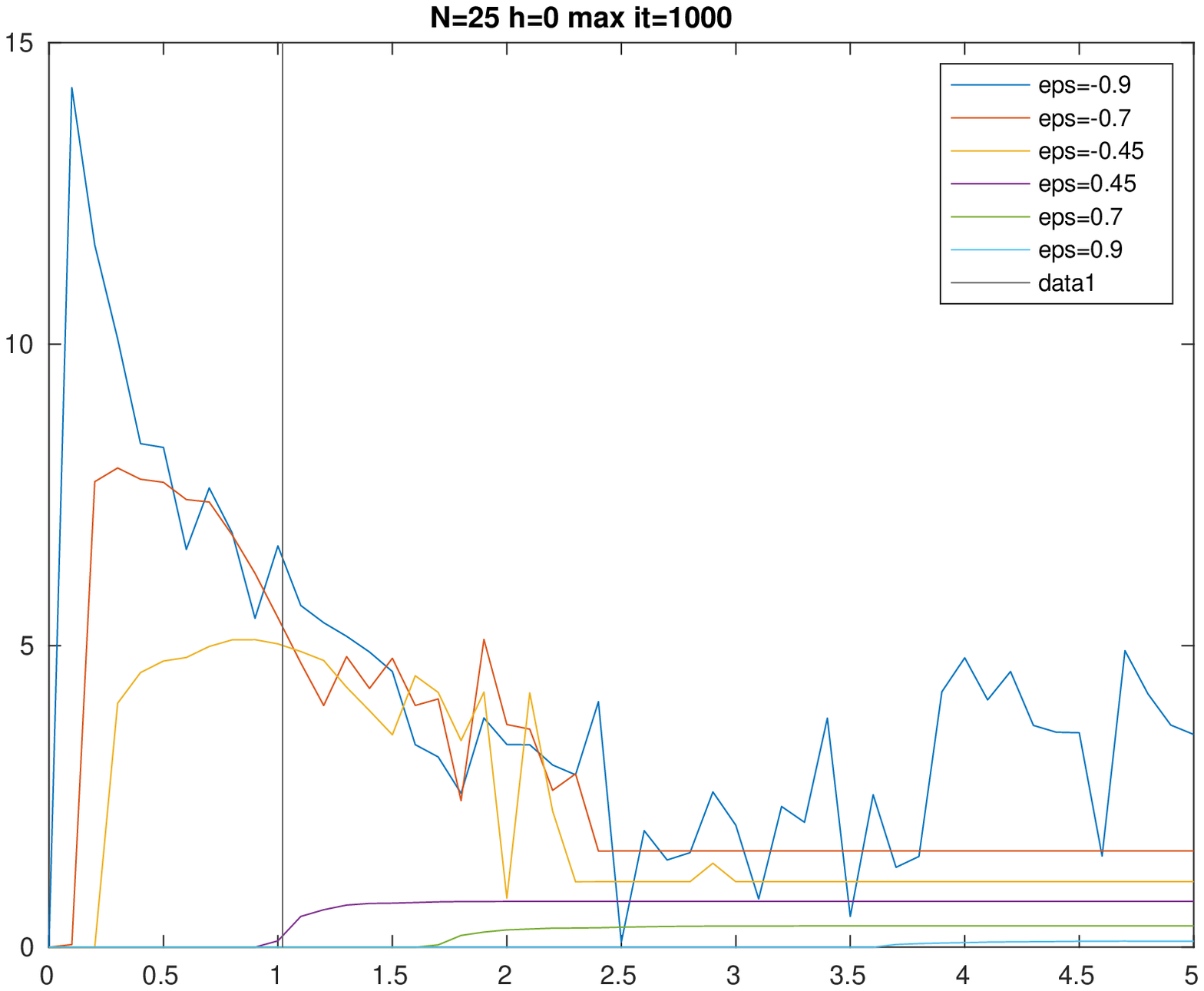}
  \caption{{\tiny} \space  1000 iterations of the $\varepsilon$-Banach algorithm for $N=25, \beta\in [0, 5]$ and $h=0$. We simulated one uniformly in the corners of the hypercube distributed start value for $\varepsilon = -0.9, -0.7, -0.45, 0.45, 0.7, 0.9$. We plotted the {\it maximum absolute error} for the given parameters and used the same start value for all $\beta$.}
\end{figure}
\begin{figure}[H]
  \centering
  % \psfrag{2}{\tiny Two Step}
  \includegraphics[width=0.5\textwidth,  trim={0 0 0 0},clip]{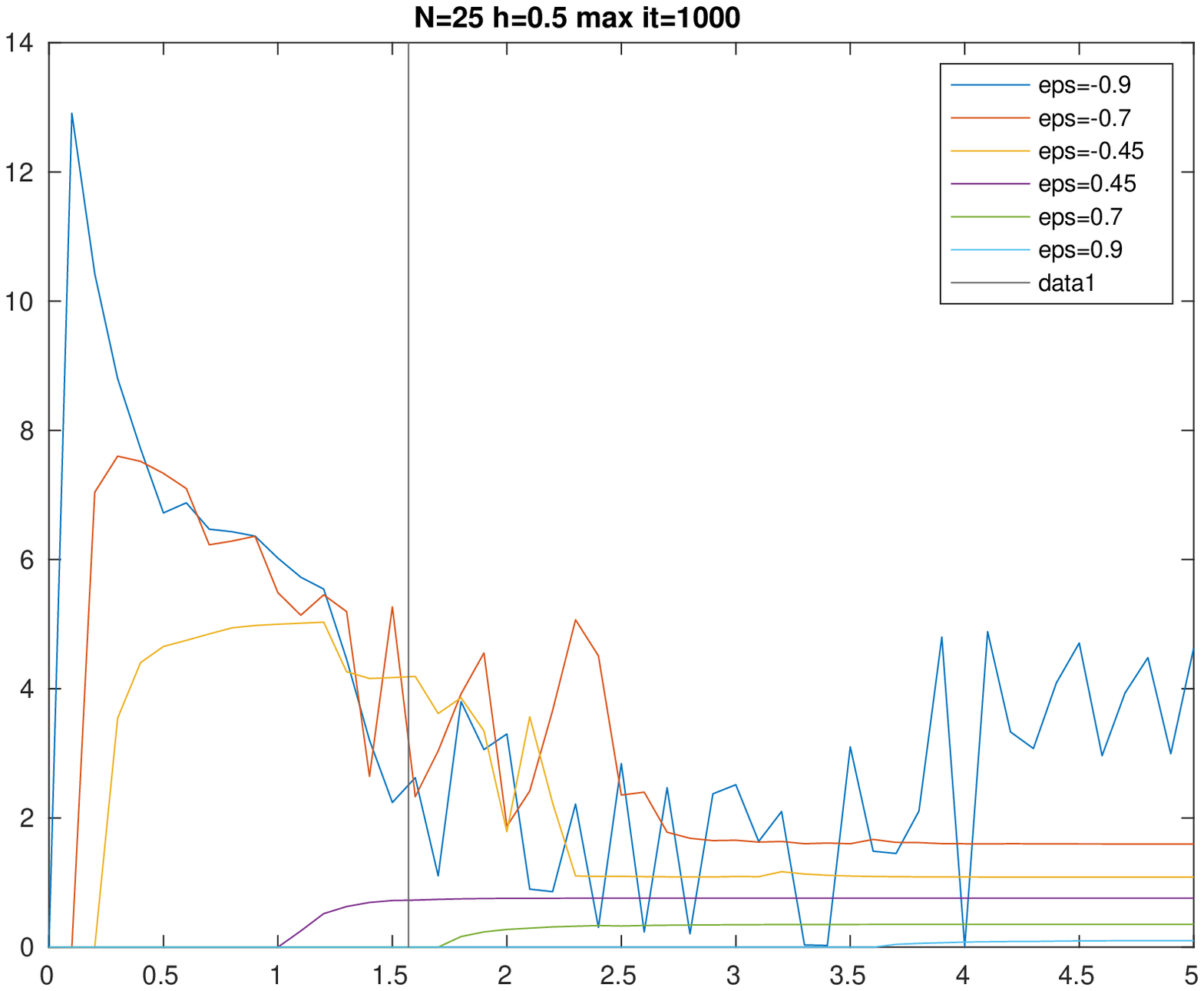}
  \caption{{\tiny} \space  1000 iterations of the $\varepsilon$-Banach algorithm for $N=25, \beta\in [0, 5]$ and $h=0.5$. We simulated one uniformly in the corners of the hypercube distributed start value for $\varepsilon = -0.9, -0.7, -0.45, 0.45, 0.7, 0.9$. We plotted the {\it maximum absolute error} for the given parameters and used the same start value for all $\beta$.}
\end{figure}
\begin{figure}[H]
  \centering
  % \psfrag{2}{\tiny Two Step}
  \includegraphics[width=0.5\textwidth,  trim={0 0 0 0},clip]{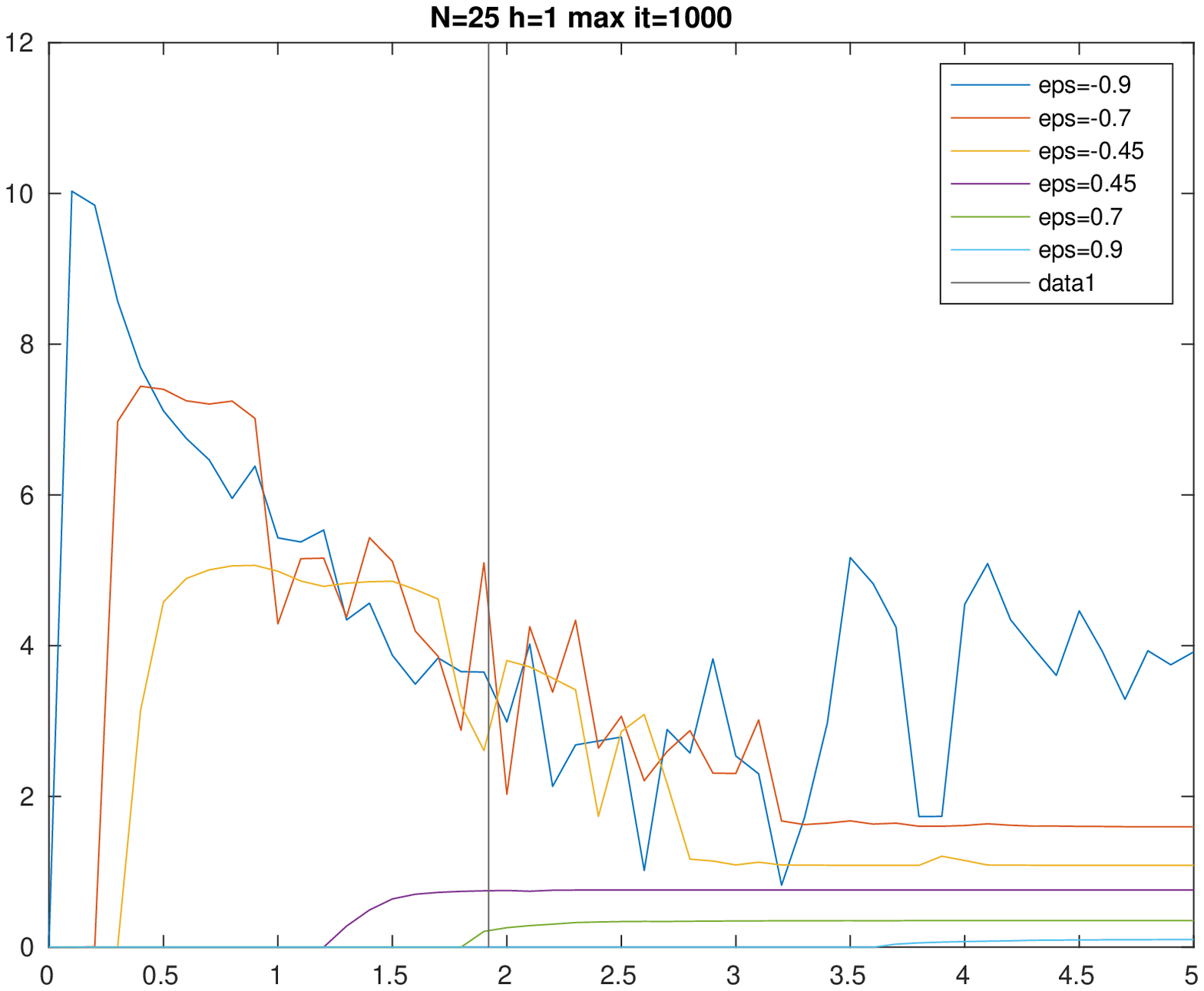}
  \caption{{\tiny} \space  1000 iterations of the $\varepsilon$-Banach algorithm for $N=25, \beta\in [0, 5]$ and $h=1$. We simulated one uniformly in the corners of the hypercube distributed start value for $\varepsilon = -0.9, -0.7, -0.45, 0.45, 0.7, 0.9$. We plotted the {\it maximum absolute error} for the given parameters and used the same start value for all $\beta$.}
\end{figure}
%%%%%% 
%%%%% 
\begin{figure}[H]
  \centering
  % \psfrag{2}{\tiny Two Step}
  \includegraphics[width=0.5\textwidth,  trim={0 0 0 0},clip]{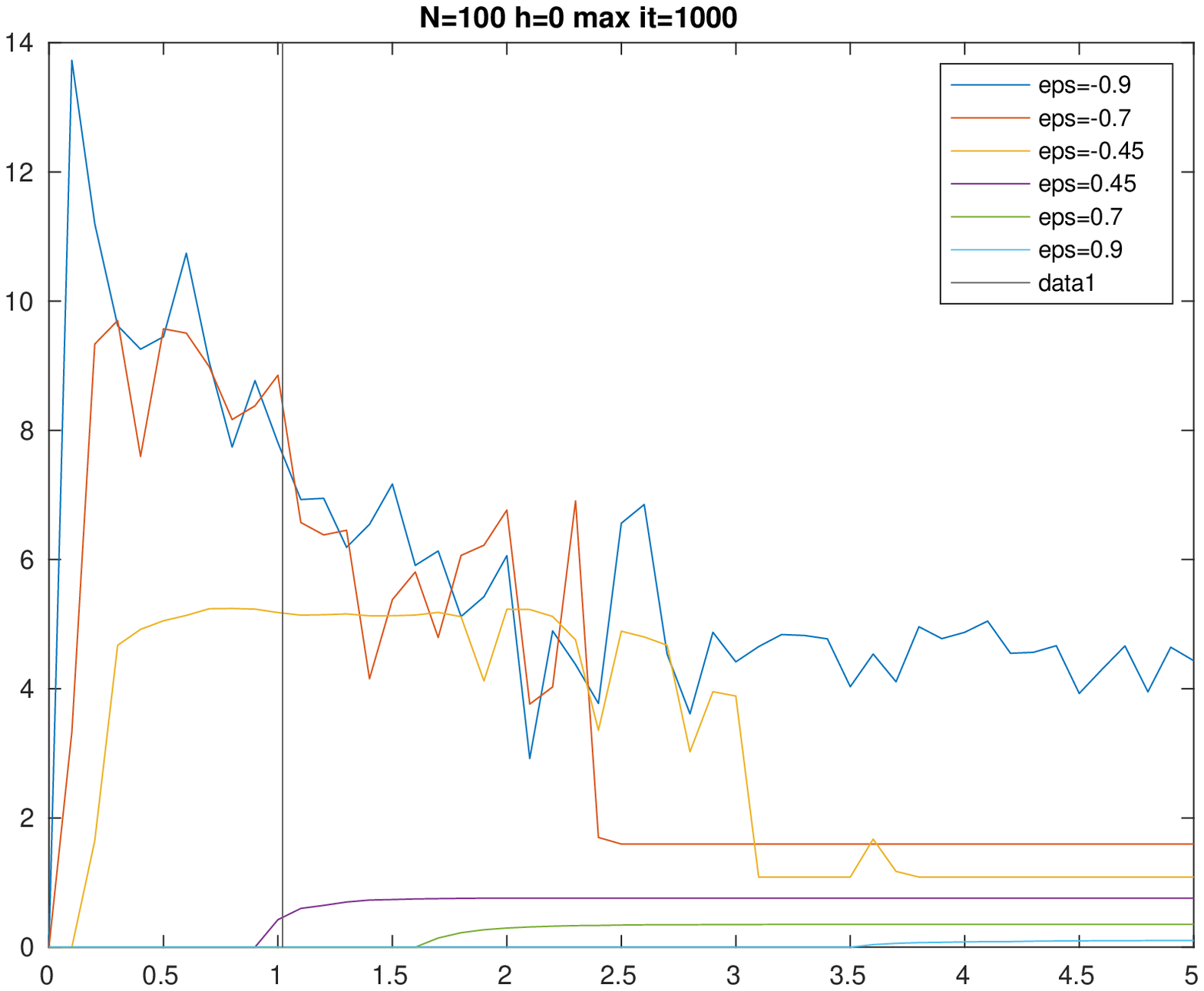}
  \caption{{\tiny} \space  1000 iterations of the $\varepsilon$-Banach algorithm for $N=100, \beta\in [0, 5]$ and $h=0$. We simulated one uniformly in the corners of the hypercube distributed start value for $\varepsilon = -0.9, -0.7, -0.45, 0.45, 0.7, 0.9$. We plotted the {\it maximum absolute error} for the given parameters and used the same start value for all $\beta$.}
\end{figure}
\begin{figure}[H]
  \centering
  % \psfrag{2}{\tiny Two Step}
  \includegraphics[width=0.5\textwidth,  trim={0 0 0 0},clip]{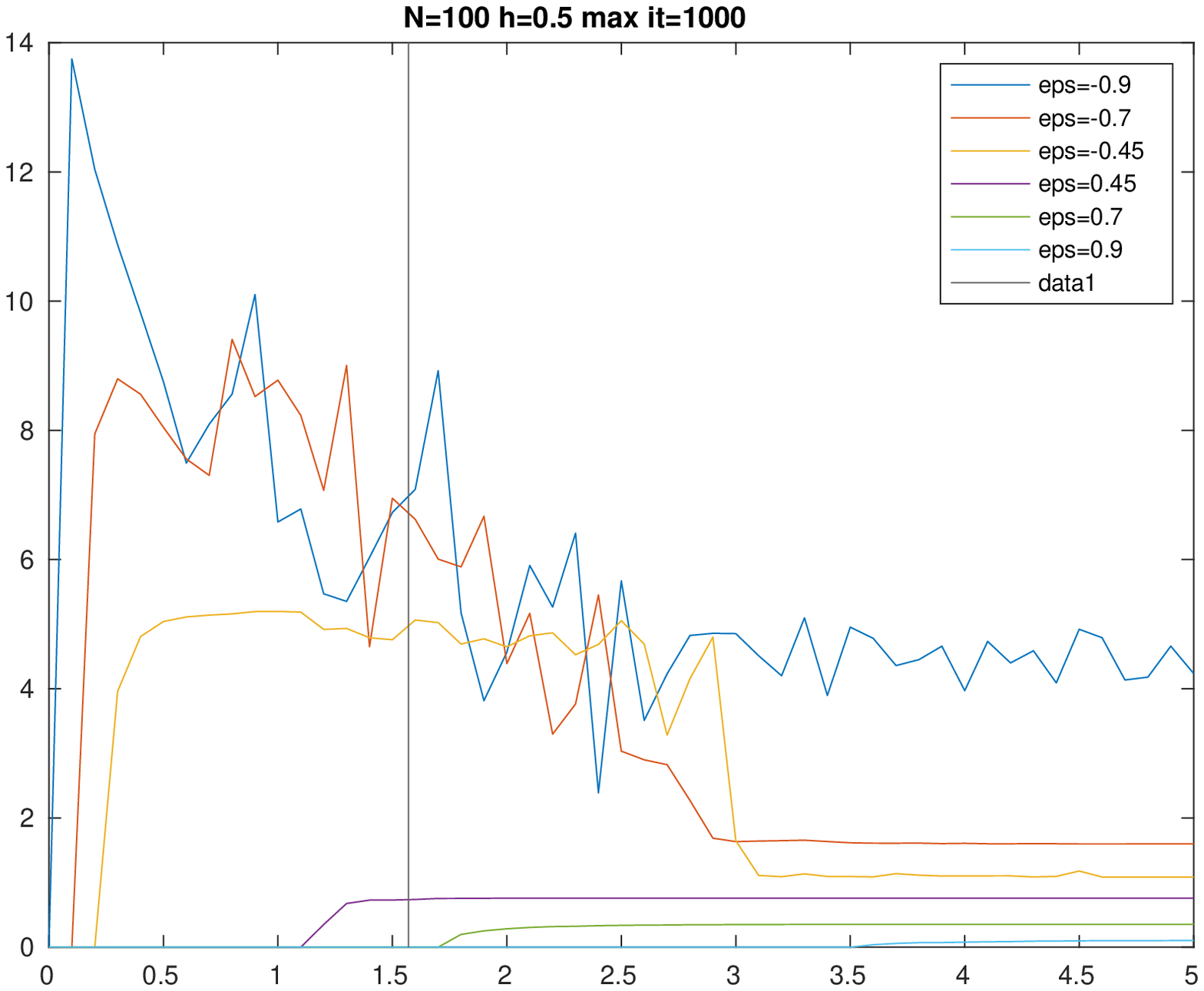}
  \caption{{\tiny} \space  1000 iterations of the $\varepsilon$-Banach algorithm for $N=100, \beta\in [0, 5]$ and $h=0.5$. We simulated one uniformly in the corners of the hypercube distributed start value for $\varepsilon = -0.9, -0.7, -0.45, 0.45, 0.7, 0.9$. We plotted the {\it maximum absolute error} for the given parameters and used the same start value for all $\beta$.}
\end{figure}
\begin{figure}[H]
  \centering
  % \psfrag{2}{\tiny Two Step}
  \includegraphics[width=0.5\textwidth,  trim={0 0 0 0},clip]{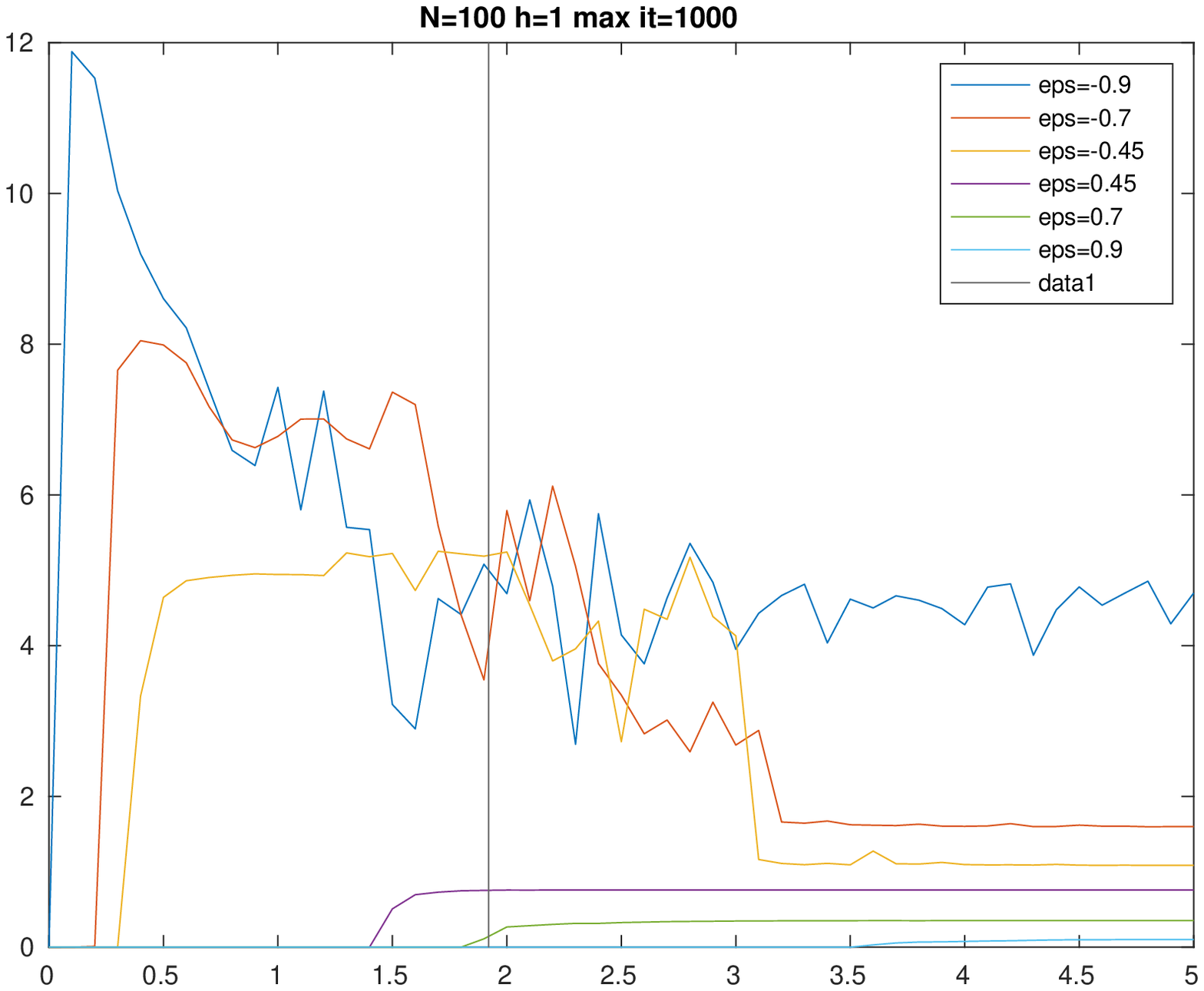}
  \caption{{\tiny} \space  1000 iterations of the $\varepsilon$-Banach algorithm for $N=100, \beta\in [0, 5]$ and $h=1$. We simulated one uniformly in the corners of the hypercube distributed start value for $\varepsilon = -0.9, -0.7, -0.45, 0.45, 0.7, 0.9$. We plotted the {\it maximum absolute error} for the given parameters and used the same start value for all $\beta$.}
\end{figure}
%%%%%% 
%%%%% 
\begin{figure}[H]
  \centering
  % \psfrag{2}{\tiny Two Step}
  \includegraphics[width=0.5\textwidth,  trim={0 0 0 0},clip]{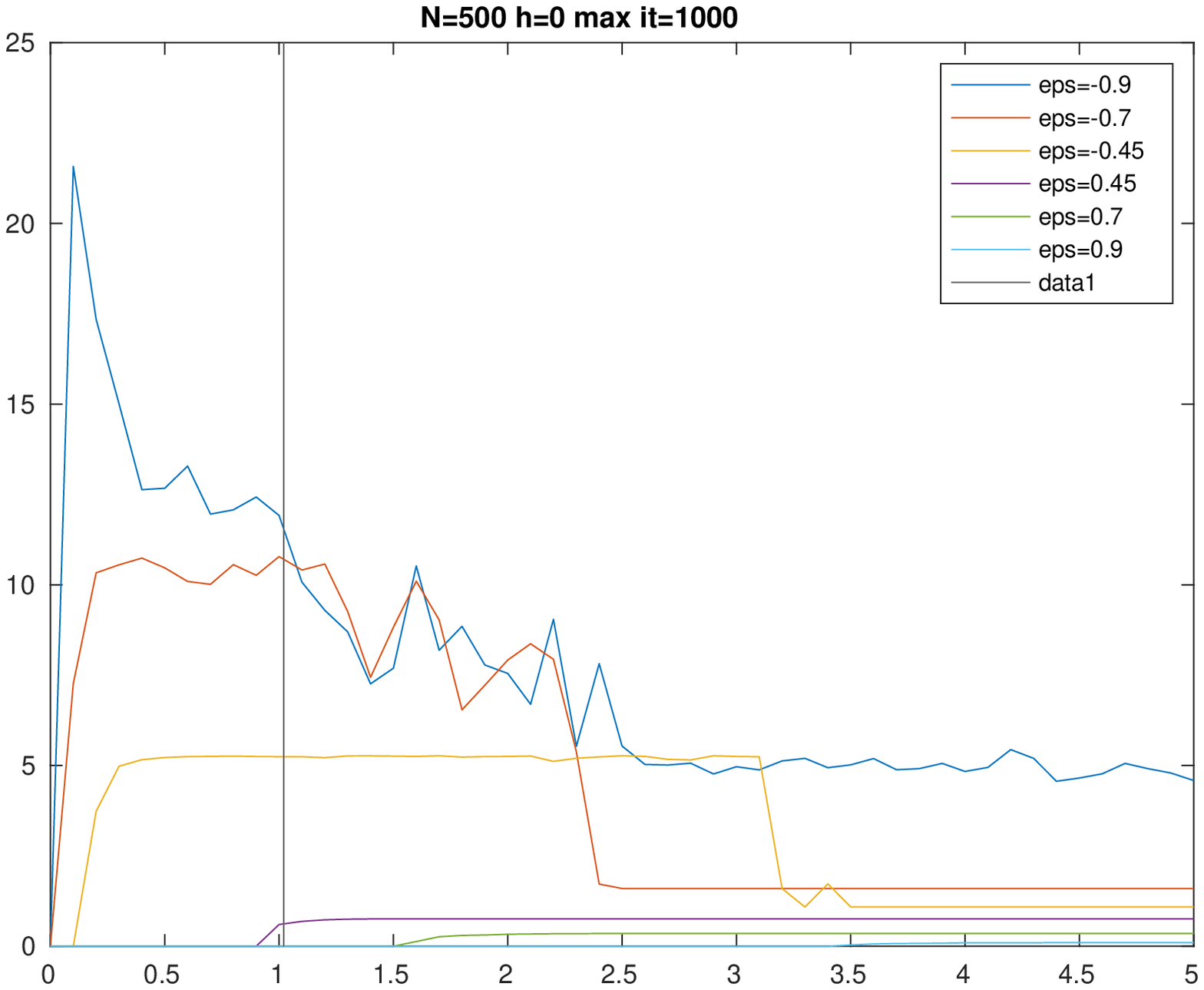}
  \caption{{\tiny} \space  1000 iterations of the $\varepsilon$-Banach algorithm for $N=500, \beta\in [0, 5]$ and $h=0$. We simulated one uniformly in the corners of the hypercube distributed start value for $\varepsilon = -0.9, -0.7, -0.45, 0.45, 0.7, 0.9$. We plotted the {\it maximum absolute error} for the given parameters and used the same start value for all $\beta$.}
\end{figure}
\begin{figure}[H]
  \centering
  % \psfrag{2}{\tiny Two Step}
  \includegraphics[width=0.5\textwidth,  trim={0 0 0 0},clip]{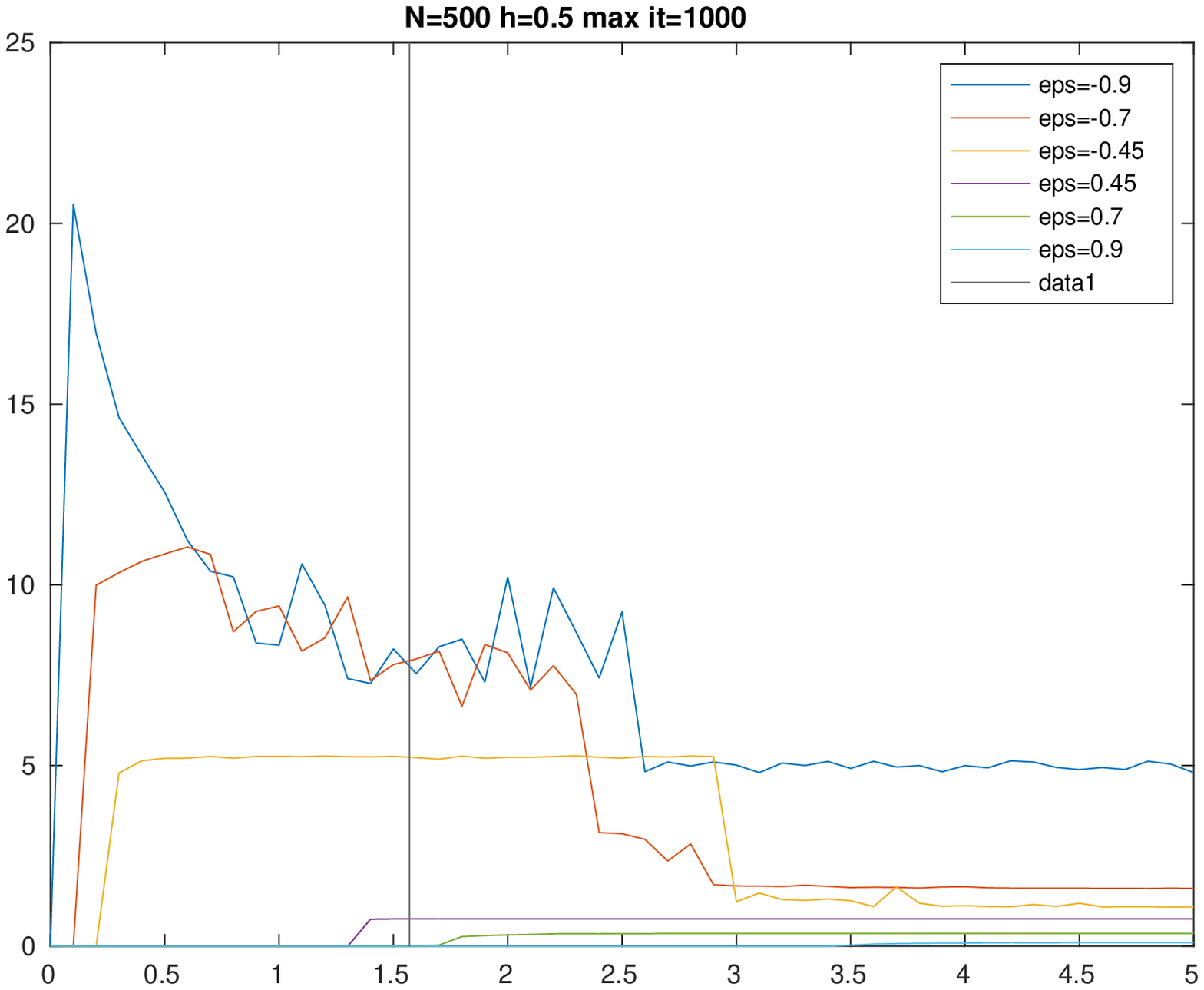}
  \caption{{\tiny} \space  1000 iterations of the $\varepsilon$-Banach algorithm for $N=500, \beta\in [0, 5]$ and $h=0.5$. We simulated one uniformly in the corners of the hypercube distributed start value for $\varepsilon = -0.9, -0.7, -0.45, 0.45, 0.7, 0.9$. We plotted the {\it maximum absolute error} for the given parameters and used the same start value for all $\beta$.}
\end{figure}
\begin{figure}[H]
  \centering
  % \psfrag{2}{\tiny Two Step}
  \includegraphics[width=0.5\textwidth,  trim={0 0 0 0},clip]{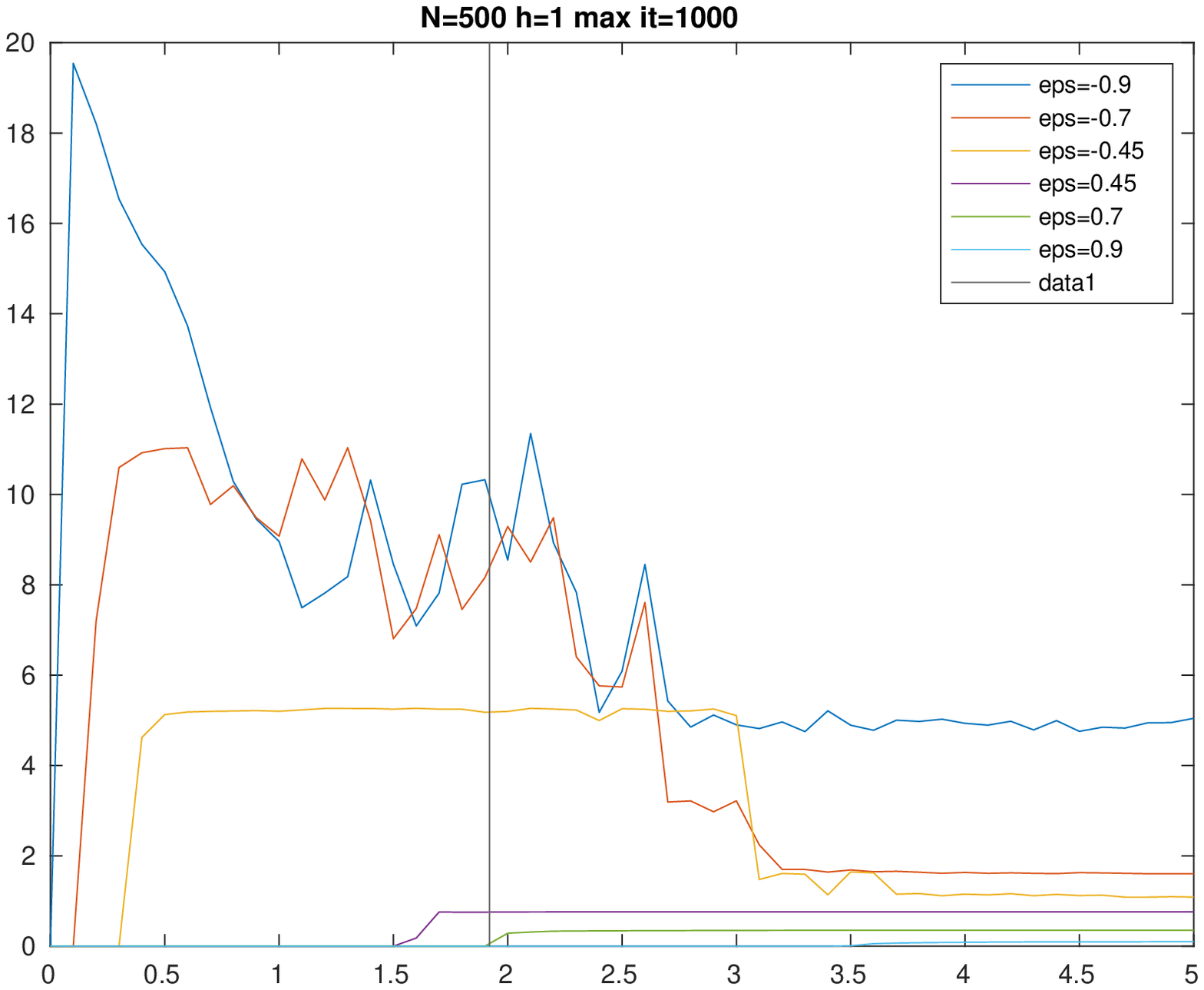}
  \caption{{\tiny} \space  1000 iterations of the $\varepsilon$-Banach algorithm for $N=500, \beta\in [0, 5]$ and $h=1$. We simulated one uniformly in the corners of the hypercube distributed start value for $\varepsilon = -0.9, -0.7, -0.45, 0.45, 0.7, 0.9$. We plotted the {\it maximum absolute error} for the given parameters and used the same start value for all $\beta$.}
\end{figure}
%%%%%% 
%%%%% 
\begin{figure}[H]
  \centering
  % \psfrag{2}{\tiny Two Step}
  \includegraphics[width=0.5\textwidth,  trim={0 0 0 0},clip]{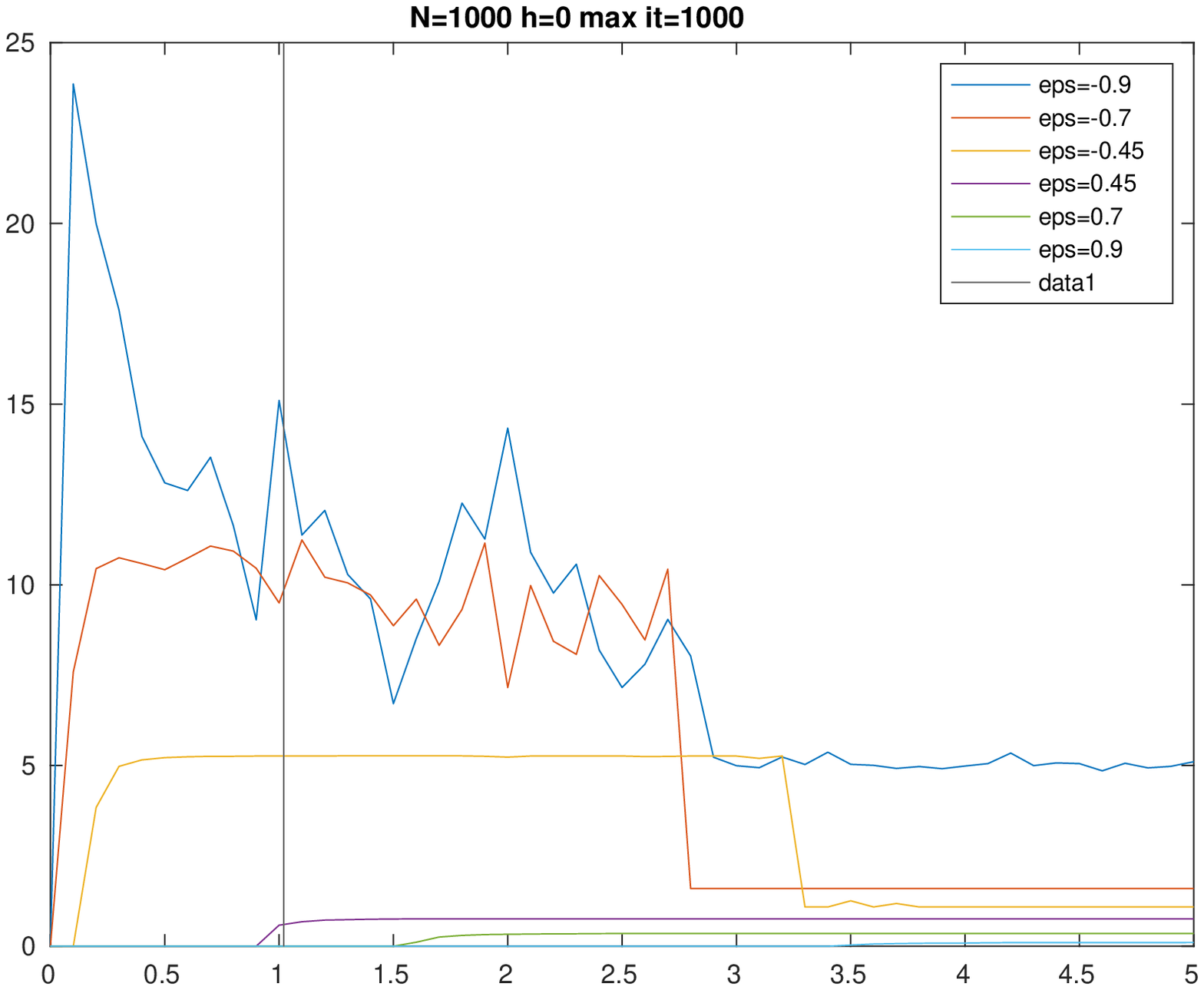}
  \caption{{\tiny} \space  1000 iterations of the $\varepsilon$-Banach algorithm for $N=1000, \beta\in [0, 5]$ and $h=0$. We simulated one uniformly in the corners of the hypercube distributed start value for $\varepsilon = -0.9, -0.7, -0.45, 0.45, 0.7, 0.9$. We plotted the {\it maximum absolute error} for the given parameters and used the same start value for all $\beta$.}
\end{figure}
\begin{figure}[H]
  \centering
  % \psfrag{2}{\tiny Two Step}
  \includegraphics[width=0.5\textwidth,  trim={0 0 0 0},clip]{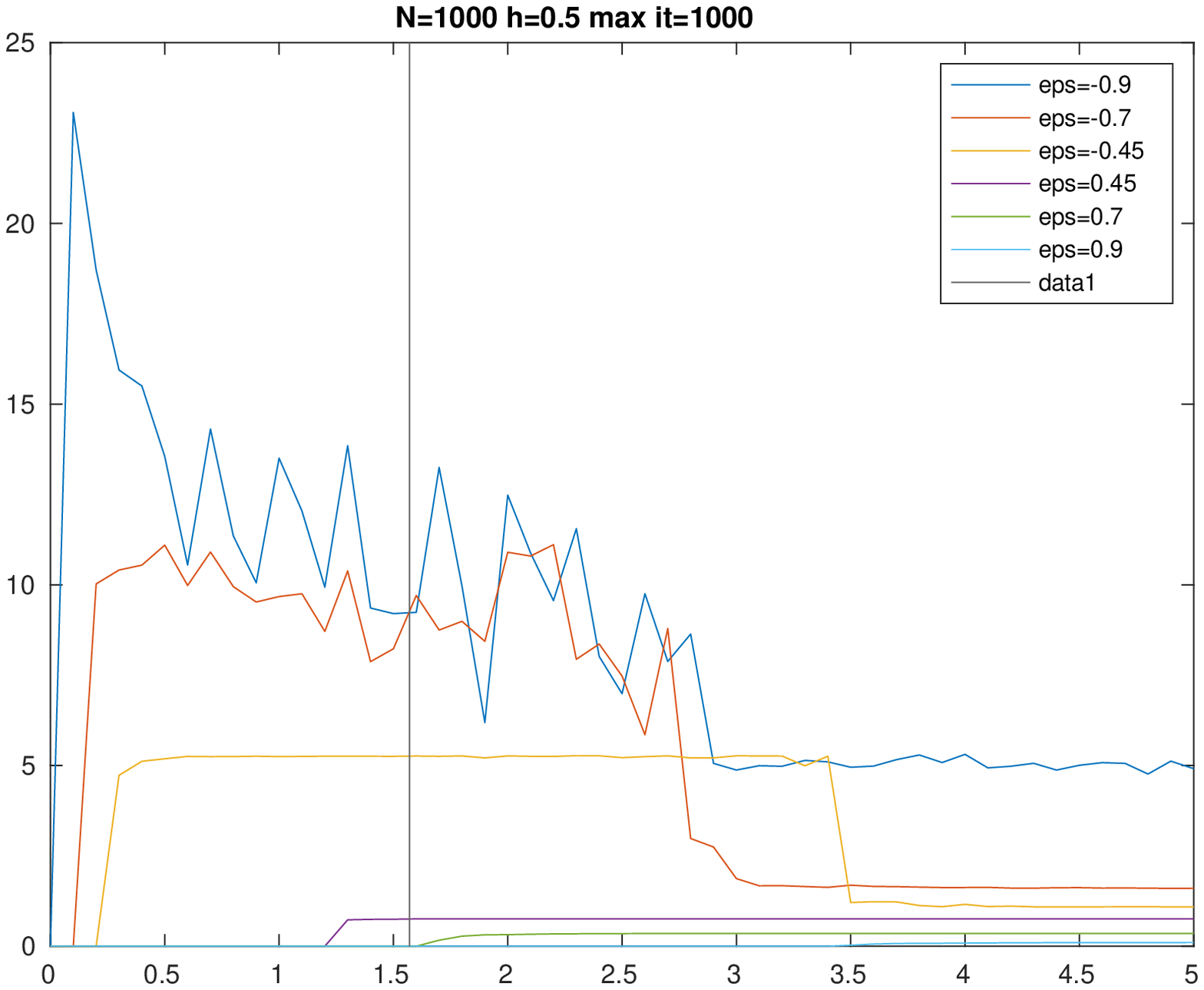}
  \caption{{\tiny} \space  1000 iterations of the $\varepsilon$-Banach algorithm for $N=1000, \beta\in [0, 5]$ and $h=0.5$. We simulated one uniformly in the corners of the hypercube distributed start value for $\varepsilon = -0.9, -0.7, -0.45, 0.45, 0.7, 0.9$. We plotted the {\it maximum absolute error} for the given parameters and used the same start value for all $\beta$.}
\end{figure}
\begin{figure}[H]
  \centering
  % \psfrag{2}{\tiny Two Step}
  \includegraphics[width=0.5\textwidth,  trim={0 0 0 0},clip]{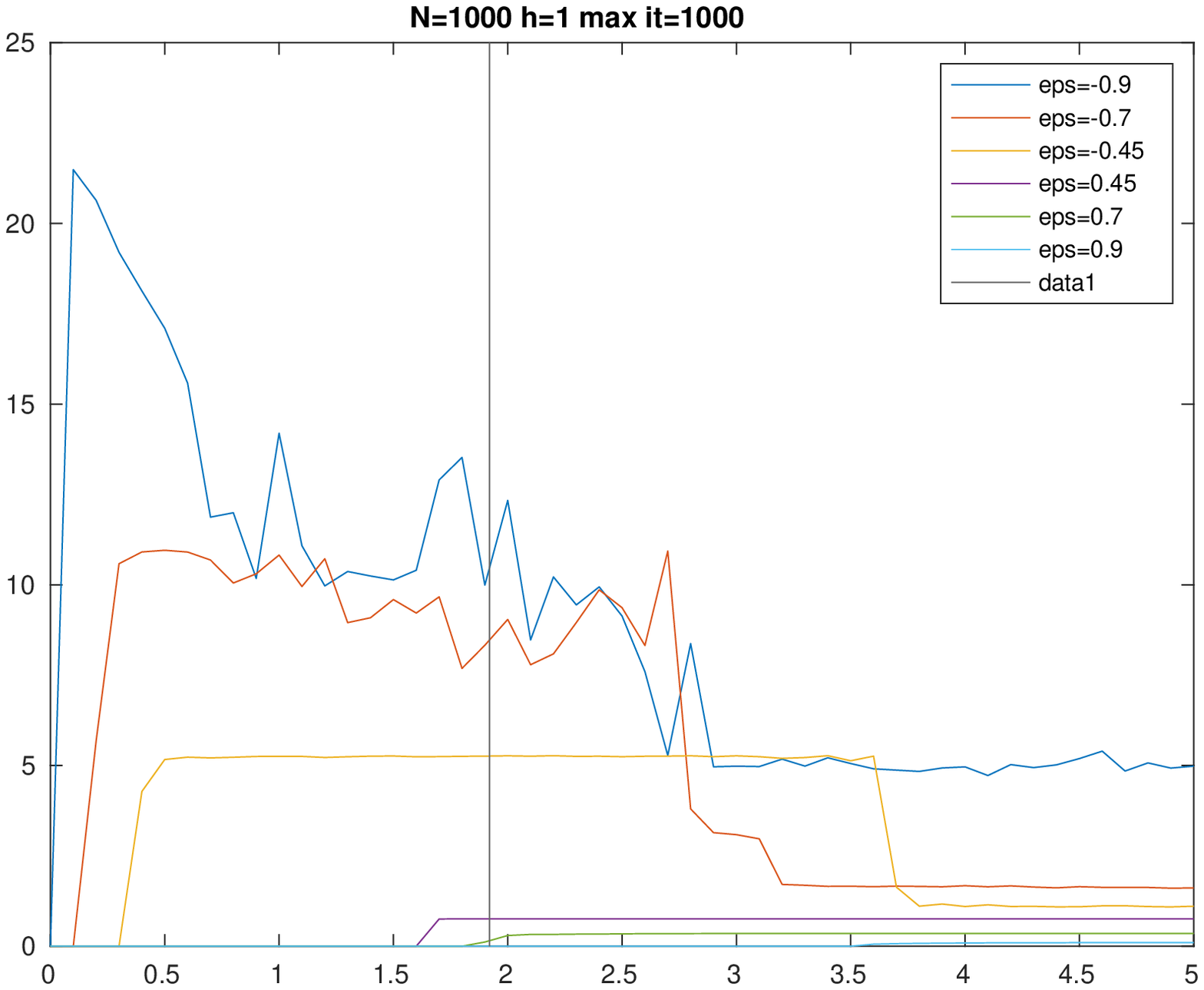}
  \caption{{\tiny} \space  1000 iterations of the $\varepsilon$-Banach algorithm for $N=1000, \beta\in [0, 5]$ and $h=1$. We simulated one uniformly in the corners of the hypercube distributed start value for $\varepsilon = -0.9, -0.7, -0.45, 0.45, 0.7, 0.9$. We plotted the {\it maximum absolute error} for the given parameters and used the same start value for all $\beta$.}
\end{figure}
\FloatBarrier\clearpage
\subsection{For $\beta =3 $ and $h=0.5$ dependence on the start values for  uniformly  in the hypercube chosen ones}
\subsubsection{$N=10$,  $\beta =3 $ and $h=0.5$}\leavevmode
For all of the following figures five we used the same start values:
\begin{figure}[H]
  \centering
  % \psfrag{2}{\tiny Two Step}
  \includegraphics[width=0.5\textwidth,  trim={0.45cm 0.4cm 0 0.4cm},clip]{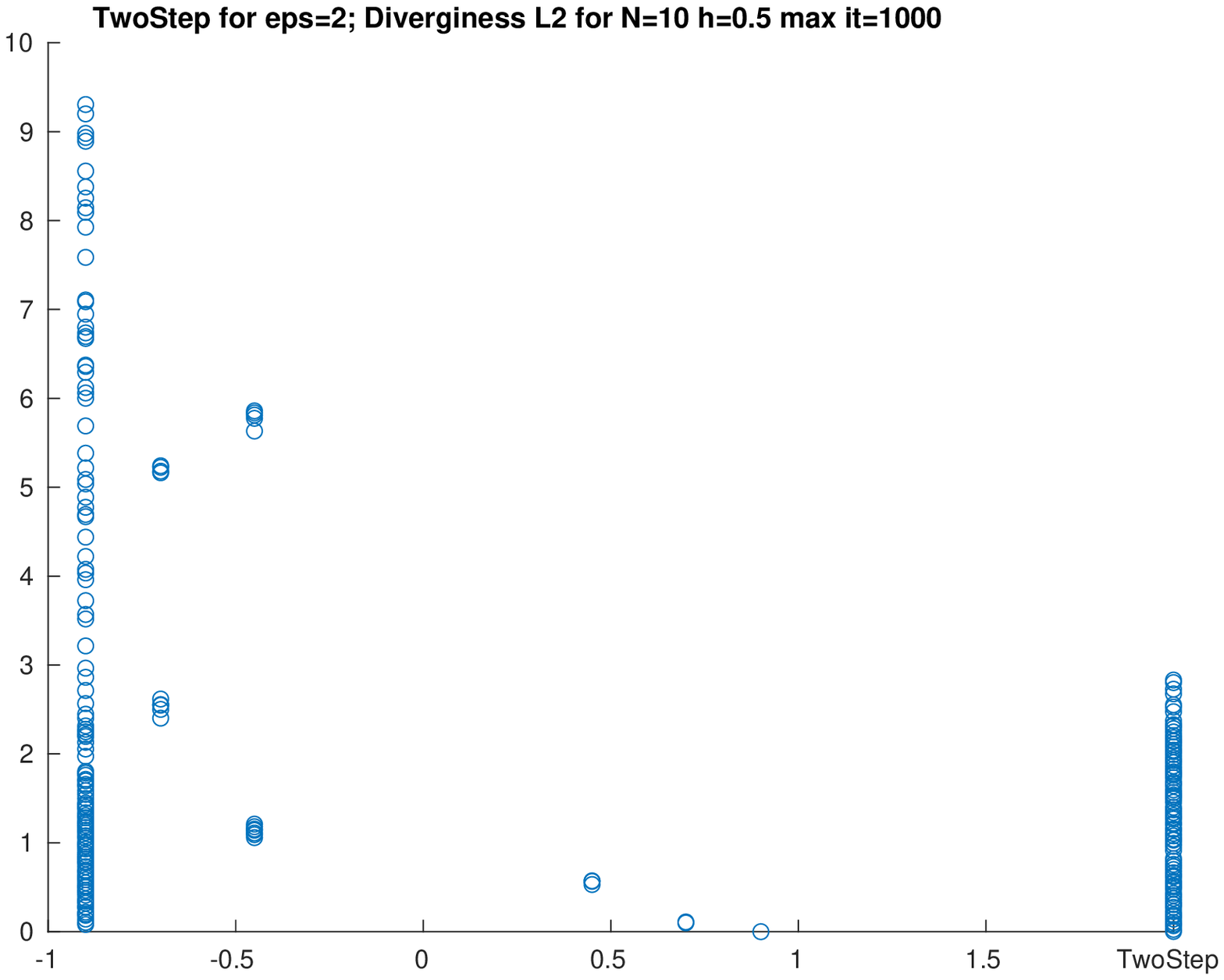}
  \caption{{\tiny} \space  1000 iterations of the $\varepsilon$-Banach and the Two-Step-Banach algorithm for $N=10, \beta=3$ and $h=0.5$. We simulated 1000 uniformly in the hypercube distributed start values for $\varepsilon = -0.9, -0.7, -0.45, 0.45, 0.7, 0.9$ and Two-Step-Banach. We plotted the {\it mean squared error} for the given parameters and used the same 1000 start values for all $\varepsilon$ and Two-Step-Banach.}
\end{figure}
\begin{figure}[H]
  \centering
  % \psfrag{2}{\tiny Two Step}
  \includegraphics[width=0.5\textwidth,  trim={0 0 0 0.4cm},clip]{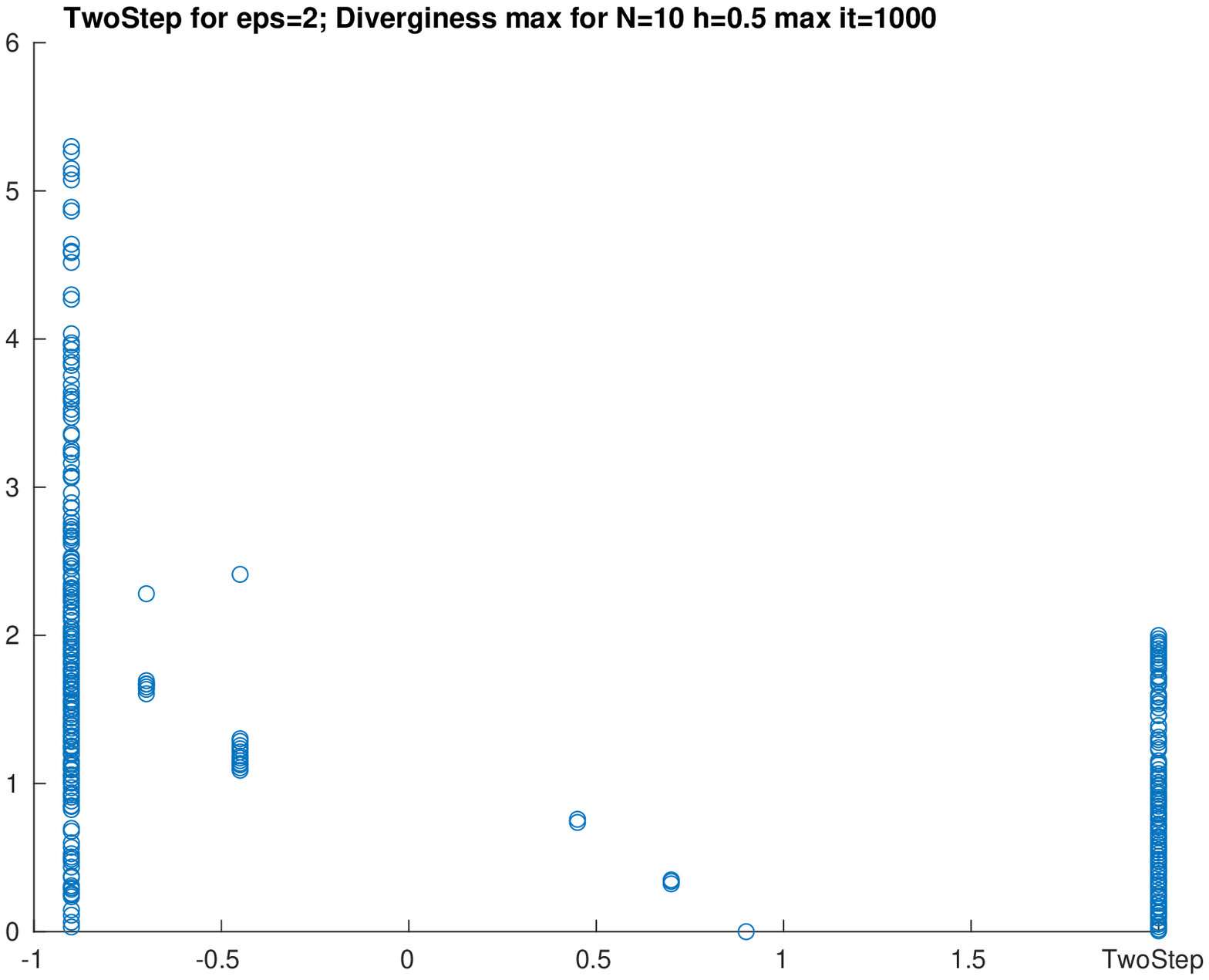}
  \caption{{\tiny} \space  1000 iterations of the $\varepsilon$-Banach and the Two-Step-Banach algorithm for $N=10, \beta=3$ and $h=0.5$. We simulated 1000 uniformly in the hypercube distributed start values for $\varepsilon = -0.9, -0.7, -0.45, 0.45, 0.7, 0.9$ and Two-Step-Banach. We plotted the {\it maximum absolute error} for the given parameters and used the same 1000 start values for all $\varepsilon$ and Two-Step-Banach.}
\end{figure}
\begin{figure}[H]
  \centering
  % \psfrag{2}{\tiny Two Step}
  \includegraphics[width=0.5\textwidth,  trim={0 0 0 0.4cm},clip]{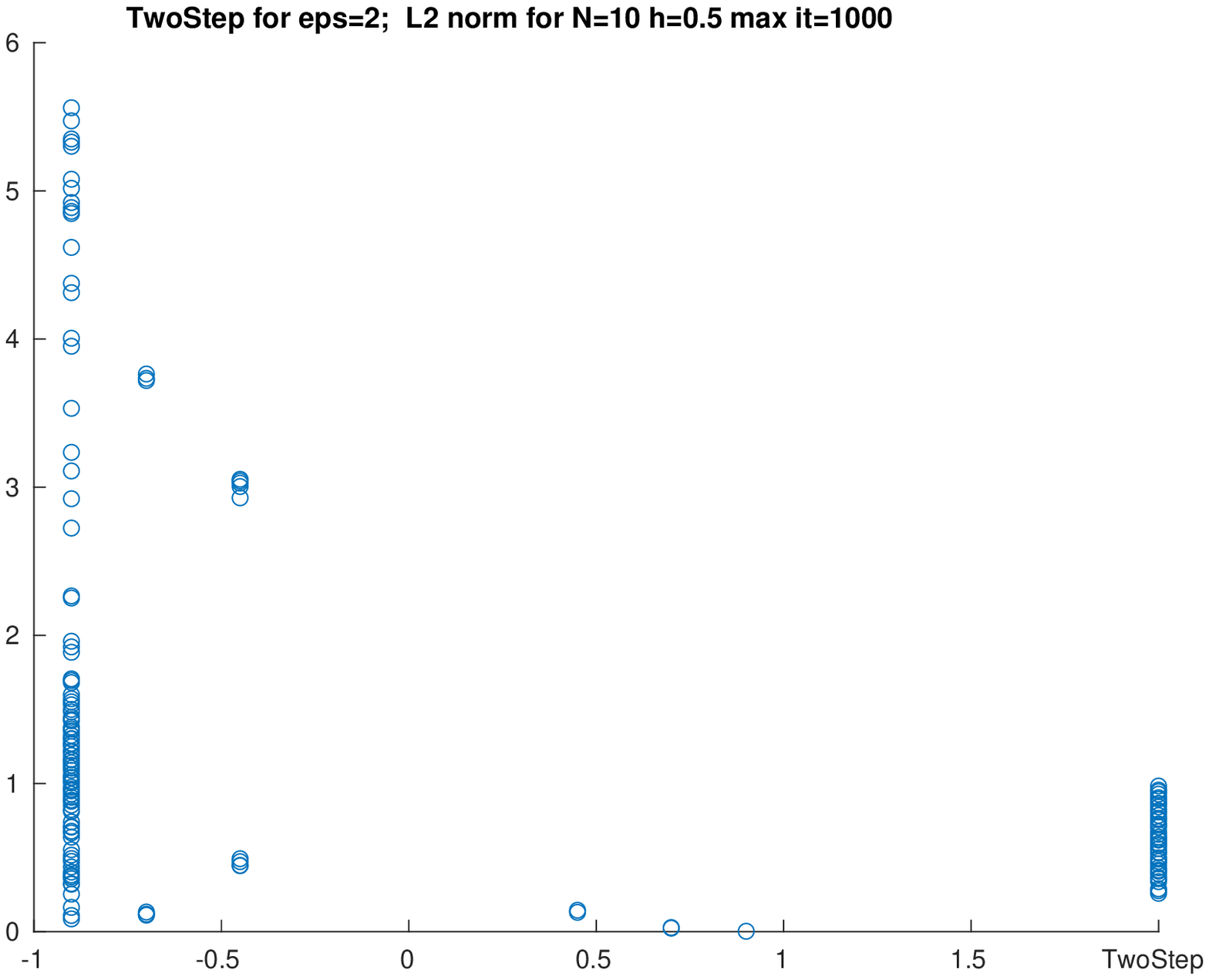}
  \caption{{\tiny} \space  1000 iterations of the $\varepsilon$-Banach and the Two-Step-Banach algorithm for $N=10, \beta=3$ and $h=0.5$. We simulated 1000 uniformly in the hypercube distributed start values for $\varepsilon = -0.9, -0.7, -0.45, 0.45, 0.7, 0.9$ and Two-Step-Banach. We plotted the {\it mean squared errors} for the given parameters and used the same 1000 start values for all $\varepsilon$ and Two-Step-Banach.}
\end{figure}
\begin{figure}[H]
  \centering
  % \psfrag{2}{\tiny Two Step}
  \includegraphics[width=0.5\textwidth,  trim={0 0 0 0.4cm},clip]{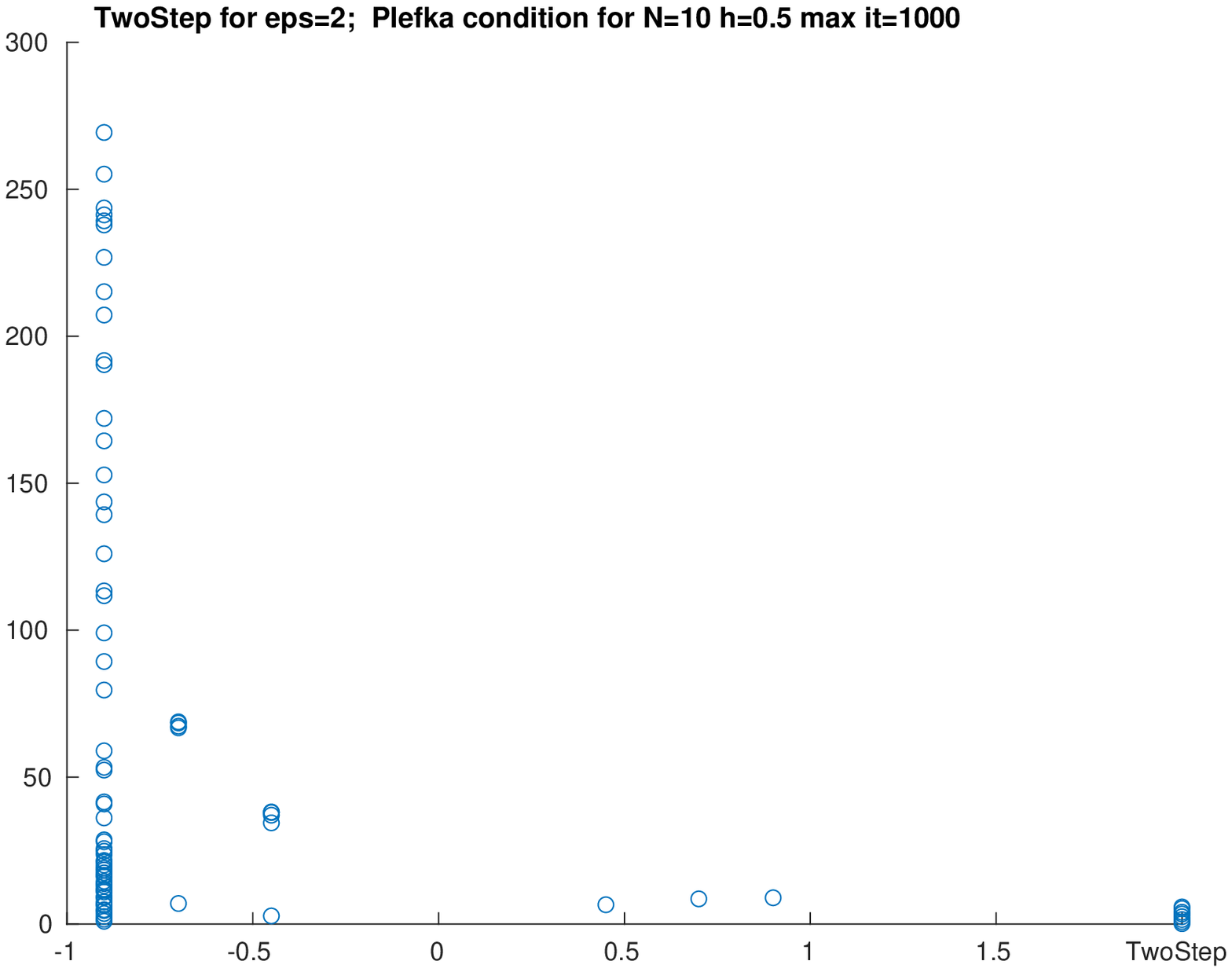}
  \caption{{\tiny} \space  1000 iterations of the $\varepsilon$-Banach and the Two-Step-Banach algorithm for $N=10, \beta=3$ and $h=0.5$. We simulated 1000 uniformly in the hypercube distributed start values for $\varepsilon = -0.9, -0.7, -0.45, 0.45, 0.7, 0.9$ and Two-Step-Banach. We plotted the Plefka values for the given parameters and used the same 1000 start values for all $\varepsilon$ and Two-Step-Banach.}
\end{figure}
\begin{figure}[H]
  \centering
  % \psfrag{2}{\tiny Two Step}
  \includegraphics[width=0.5\textwidth,  trim={0 0 0 0.4cm},clip]{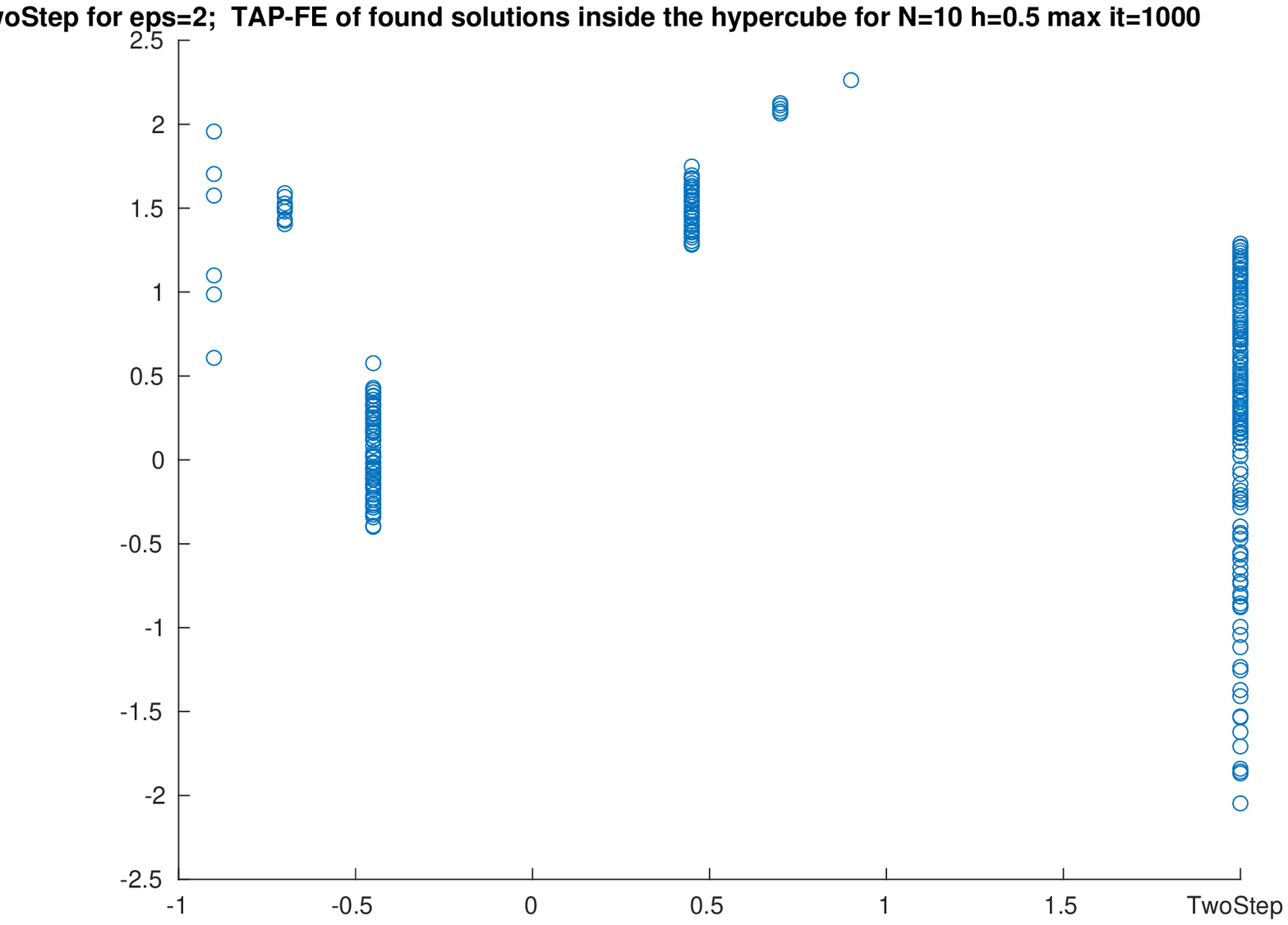}
  \caption{{\tiny} \space  1000 iterations of the $\varepsilon$-Banach and the Two-Step-Banach algorithm for $N=10, \beta=3$ and $h=0.5$. We simulated 1000 uniformly in the hypercube distributed start values for $\varepsilon = -0.9, -0.7, -0.45, 0.45, 0.7, 0.9$ and Two-Step-Banach. We plotted the TAP-free-energy values for all iterations ending inside the hypercube for the given parameters and used the same 1000 start values for all $\varepsilon$ and Two-Step-Banach.}
\end{figure}
\FloatBarrier\clearpage
\subsubsection{$N=25$,  $\beta =3 $ and $h=0.5$}\leavevmode
For all of the following figures five we used the same start values:
\begin{figure}[H]
  \centering
  % \psfrag{2}{\tiny Two Step}
  \includegraphics[width=0.5\textwidth,  trim={0 0 0 0.4cm},clip]{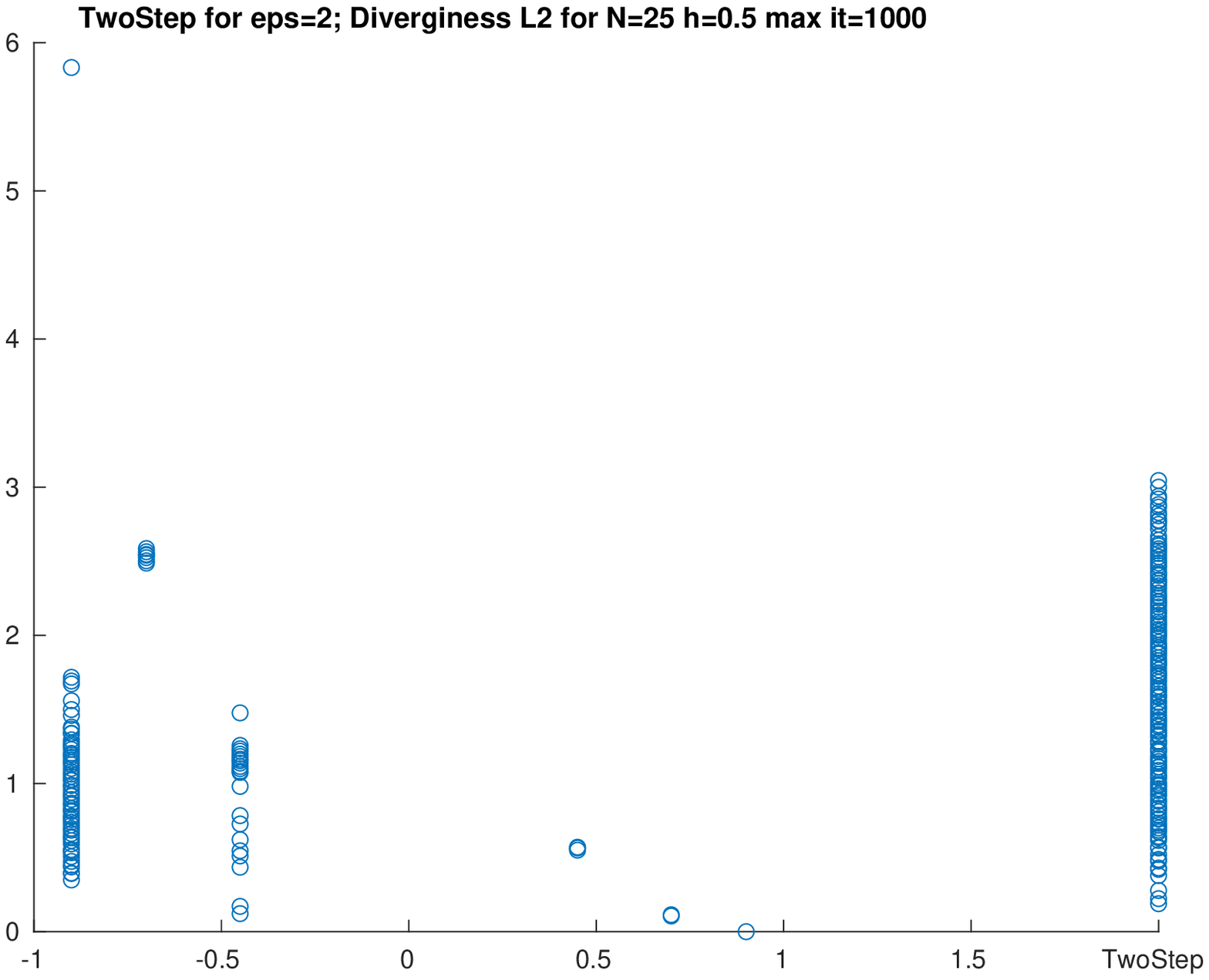}
  \caption{{\tiny} \space  1000 iterations of the $\varepsilon$-Banach and the Two-Step-Banach algorithm for $N=25, \beta=3$ and $h=0.5$. We simulated 1000 uniformly in the hypercube distributed start values for $\varepsilon = -0.9, -0.7, -0.45, 0.45, 0.7, 0.9$ and Two-Step-Banach. We plotted the {\it mean squared errors} for the given parameters and used the same 1000 start values for all $\varepsilon$ and Two-Step-Banach.}
\end{figure}
\begin{figure}[H]
  \centering
  % \psfrag{2}{\tiny Two Step}
  \includegraphics[width=0.5\textwidth,  trim={0 0 0 0.4cm},clip]{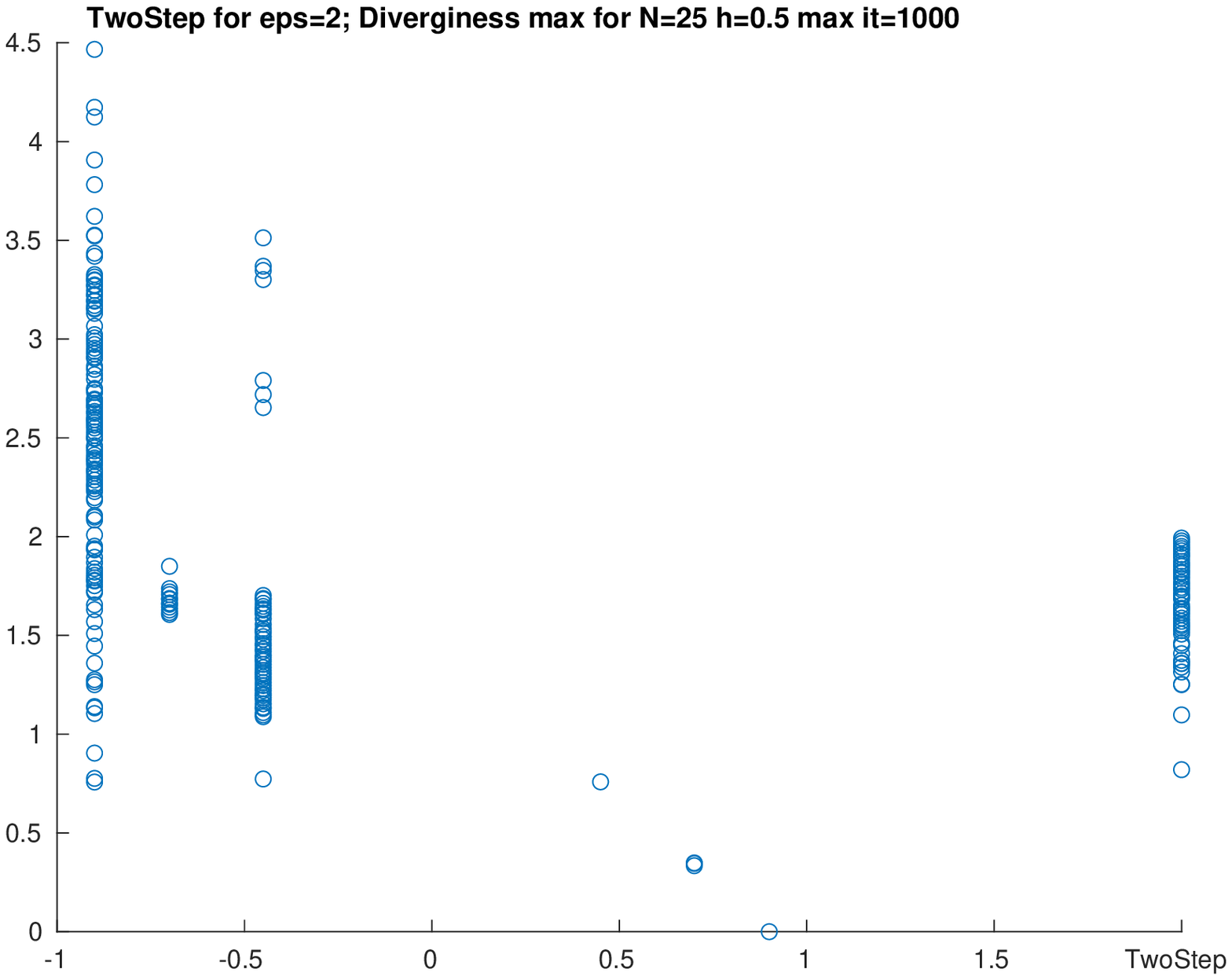}
  \caption{{\tiny} \space  1000 iterations of the $\varepsilon$-Banach and the Two-Step-Banach algorithm for $N=25, \beta=3$ and $h=0.5$. We simulated 1000 uniformly in the hypercube distributed start values for $\varepsilon = -0.9, -0.7, -0.45, 0.45, 0.7, 0.9$ and Two-Step-Banach. We plotted the {\it maximum absolute errors} for the given parameters and used the same 1000 start values for all $\varepsilon$ and Two-Step-Banach.}
\end{figure}
\begin{figure}[H]
  \centering
  % \psfrag{2}{\tiny Two Step}
  \includegraphics[width=0.5\textwidth,  trim={0 0 0 0.4cm},clip]{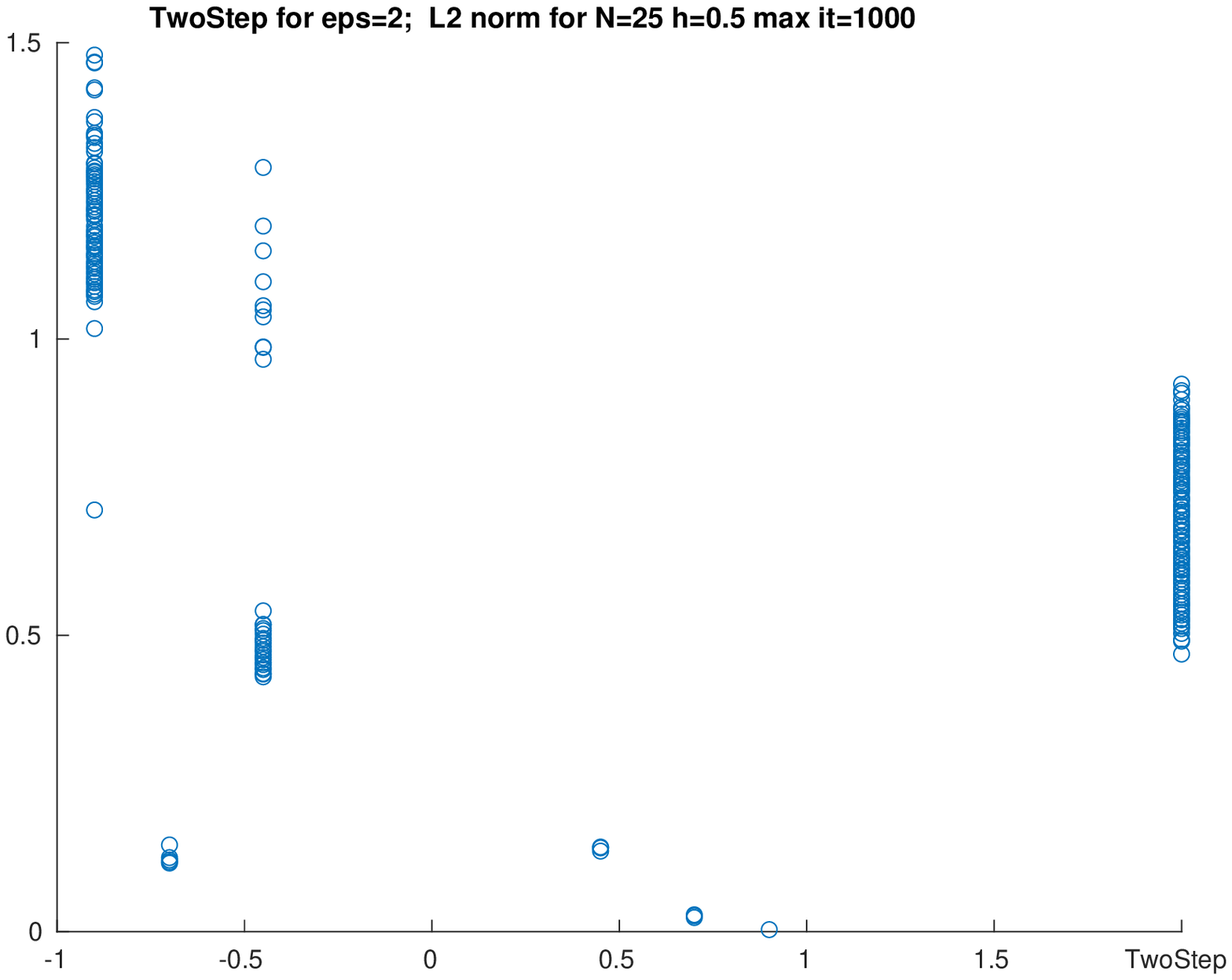}
  \caption{{\tiny} \space  1000 iterations of the $\varepsilon$-Banach and the Two-Step-Banach algorithm for $N=25, \beta=3$ and $h=0.5$. We simulated 1000 uniformly in the hypercube distributed start values for $\varepsilon = -0.9, -0.7, -0.45, 0.45, 0.7, 0.9$ and Two-Step-Banach. We plotted the {\it mean squared errors} for the given parameters and used the same 1000 start values for all $\varepsilon$ and Two-Step-Banach.}
\end{figure}
\begin{figure}[H]
  \centering
  % \psfrag{2}{\tiny Two Step}
  \includegraphics[width=0.5\textwidth,  trim={0 0 0 0.4cm},clip]{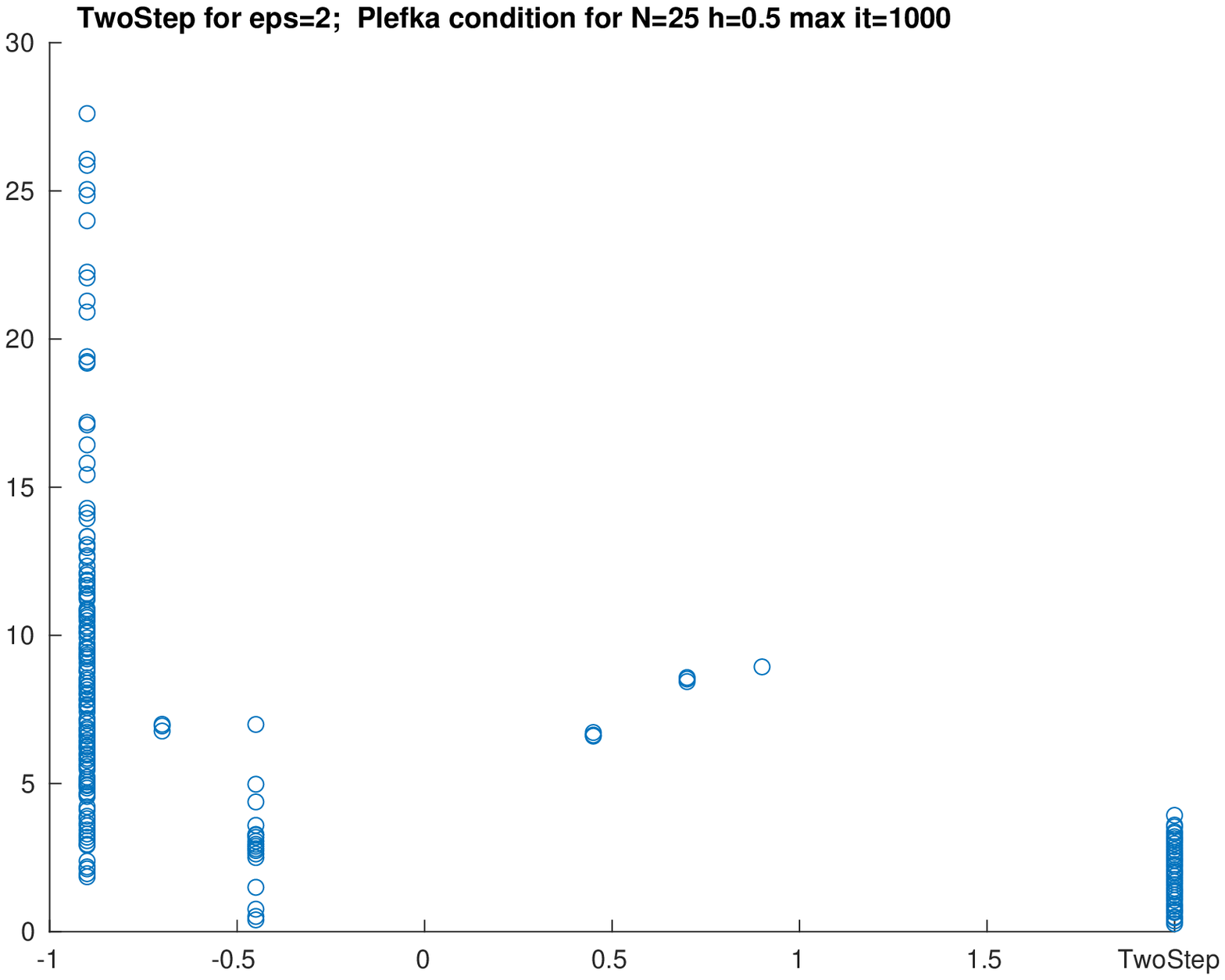}
  \caption{{\tiny} \space  1000 iterations of the $\varepsilon$-Banach and the Two-Step-Banach algorithm for $N=25, \beta=3$ and $h=0.5$. We simulated 1000 uniformly in the hypercube distributed start values for $\varepsilon = -0.9, -0.7, -0.45, 0.45, 0.7, 0.9$ and Two-Step-Banach. We plotted the Plefka values for the given parameters and used the same 1000 start values for all $\varepsilon$ and Two-Step-Banach.}
\end{figure}
\begin{figure}[H]
  \centering
  % \psfrag{2}{\tiny Two Step}
  \includegraphics[width=0.5\textwidth,  trim={0 0 0 0.4cm},clip]{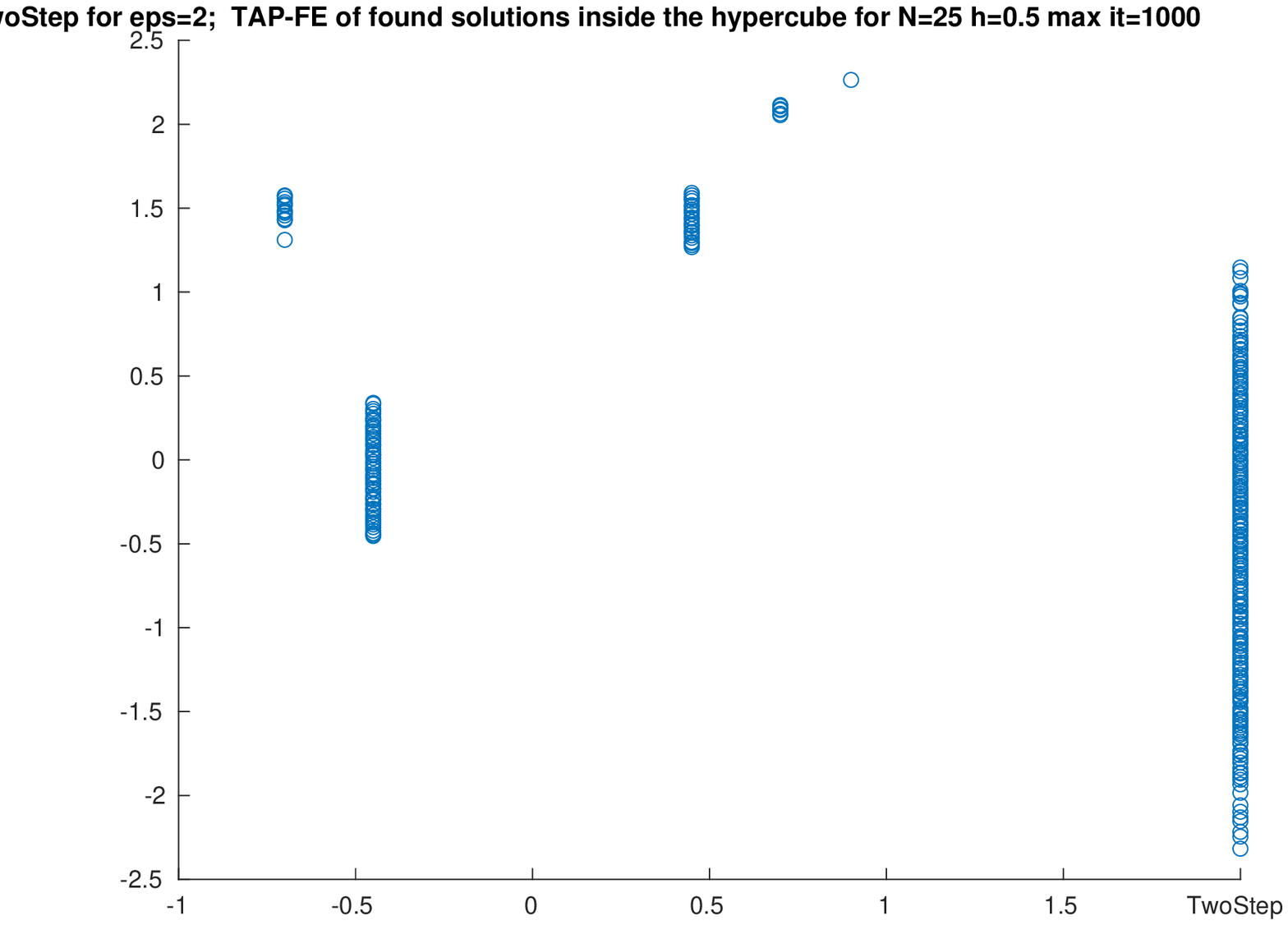}
  \caption{{\tiny} \space  1000 iterations of the $\varepsilon$-Banach and the Two-Step-Banach algorithm for $N=25, \beta=3$ and $h=0.5$. We simulated 1000 uniformly in the hypercube distributed start values for $\varepsilon = -0.9, -0.7, -0.45, 0.45, 0.7, 0.9$ and Two-Step-Banach. We plotted the TAP-free-energy values for all iterations ending inside the hypercube for the given parameters and used the same 1000 start values for all $\varepsilon$ and Two-Step-Banach.}
\end{figure}
\FloatBarrier\clearpage
\subsubsection{$N=100$,  $\beta =3 $ and $h=0.5$}\leavevmode
For all of the following figures five we used the same start values:
\begin{figure}[H]
  \centering
  % \psfrag{2}{\tiny Two Step}
  \includegraphics[width=0.5\textwidth,  trim={0 0 0 0.4cm},clip]{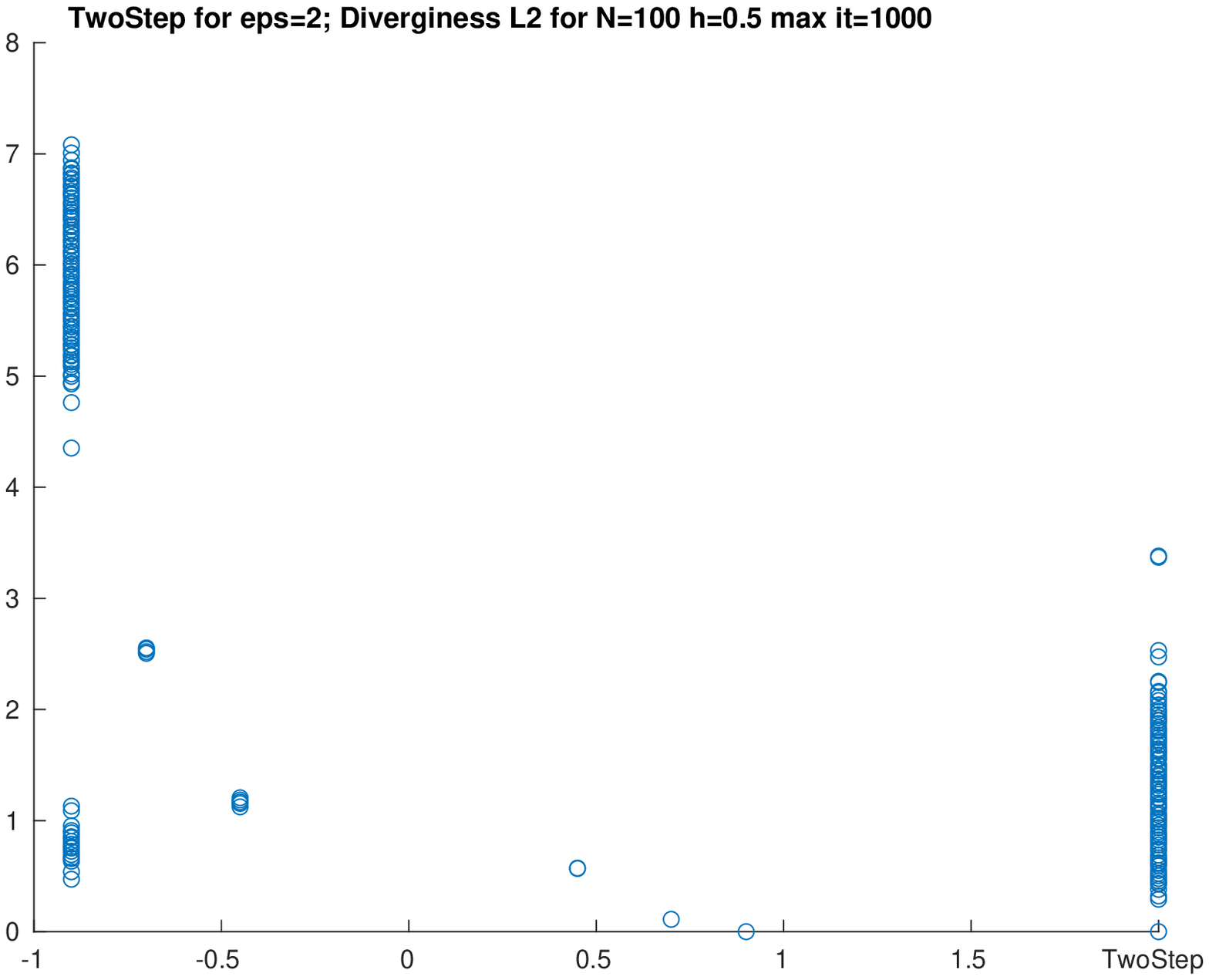}
  \caption{{\tiny} \space  1000 iterations of the $\varepsilon$-Banach and the Two-Step-Banach algorithm for $N=100, \beta=3$ and $h=0.5$. We simulated 1000 uniformly in the hypercube distributed start values for $\varepsilon = -0.9, -0.7, -0.45, 0.45, 0.7, 0.9$ and Two-Step-Banach. We plotted the {\it mean squared errors} for the given parameters and used the same 1000 start values for all $\varepsilon$ and Two-Step-Banach.}
\end{figure}
\begin{figure}[H]
  \centering
  % \psfrag{2}{\tiny Two Step}
  \includegraphics[width=0.5\textwidth,  trim={0 0 0 0.4cm},clip]{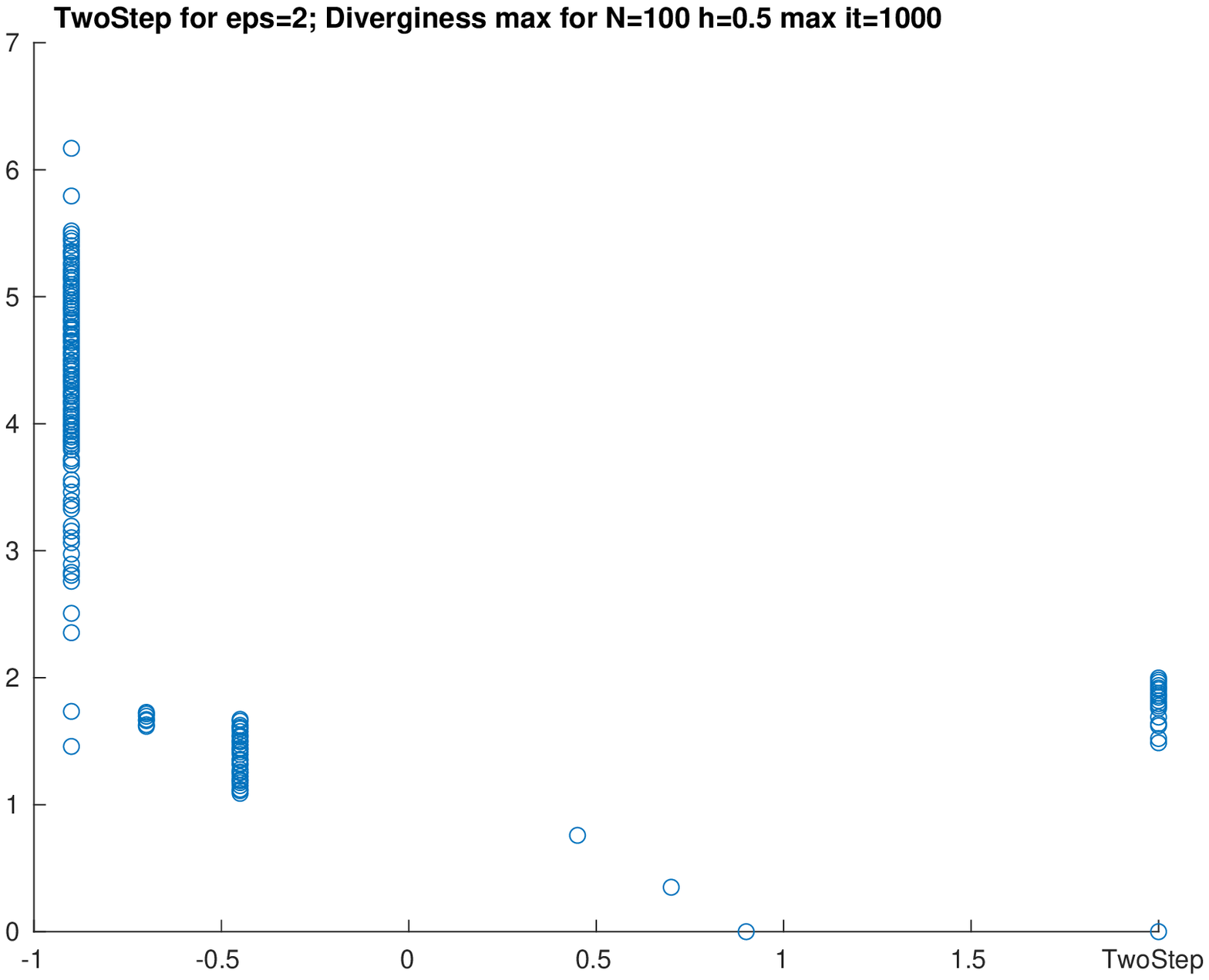}
  \caption{{\tiny} \space  1000 iterations of the $\varepsilon$-Banach and the Two-Step-Banach algorithm for $N=100, \beta=3$ and $h=0.5$. We simulated 1000 uniformly in the hypercube distributed start values for $\varepsilon = -0.9, -0.7, -0.45, 0.45, 0.7, 0.9$ and Two-Step-Banach. We plotted the {\it maximum absolute errors} for the given parameters and used the same 1000 start values for all $\varepsilon$ and Two-Step-Banach.}
\end{figure}
\begin{figure}[H]
  \centering
  % \psfrag{2}{\tiny Two Step}
  \includegraphics[width=0.5\textwidth,  trim={0 0 0 0.4cm},clip]{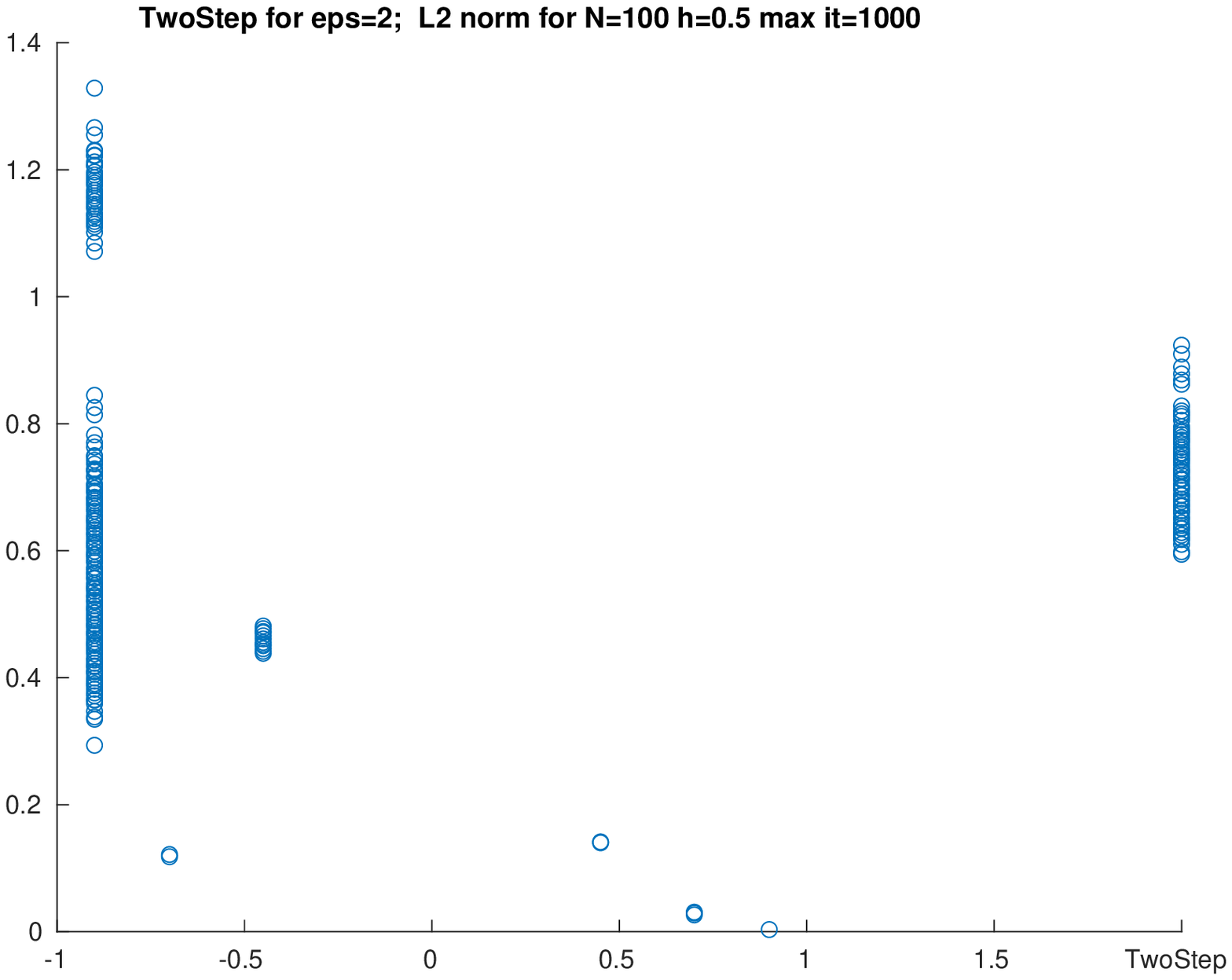}
  \caption{{\tiny} \space  1000 iterations of the $\varepsilon$-Banach and the Two-Step-Banach algorithm for $N=100, \beta=3$ and $h=0.5$. We simulated 1000 uniformly in the hypercube distributed start values for $\varepsilon = -0.9, -0.7, -0.45, 0.45, 0.7, 0.9$ and Two-Step-Banach. We plotted the {\it mean squared errors} for the given parameters and used the same 1000 start values for all $\varepsilon$ and Two-Step-Banach.}
\end{figure}
\begin{figure}[H]
  \centering
  % \psfrag{2}{\tiny Two Step}
  \includegraphics[width=0.5\textwidth,  trim={0 0 0 0.4cm},clip]{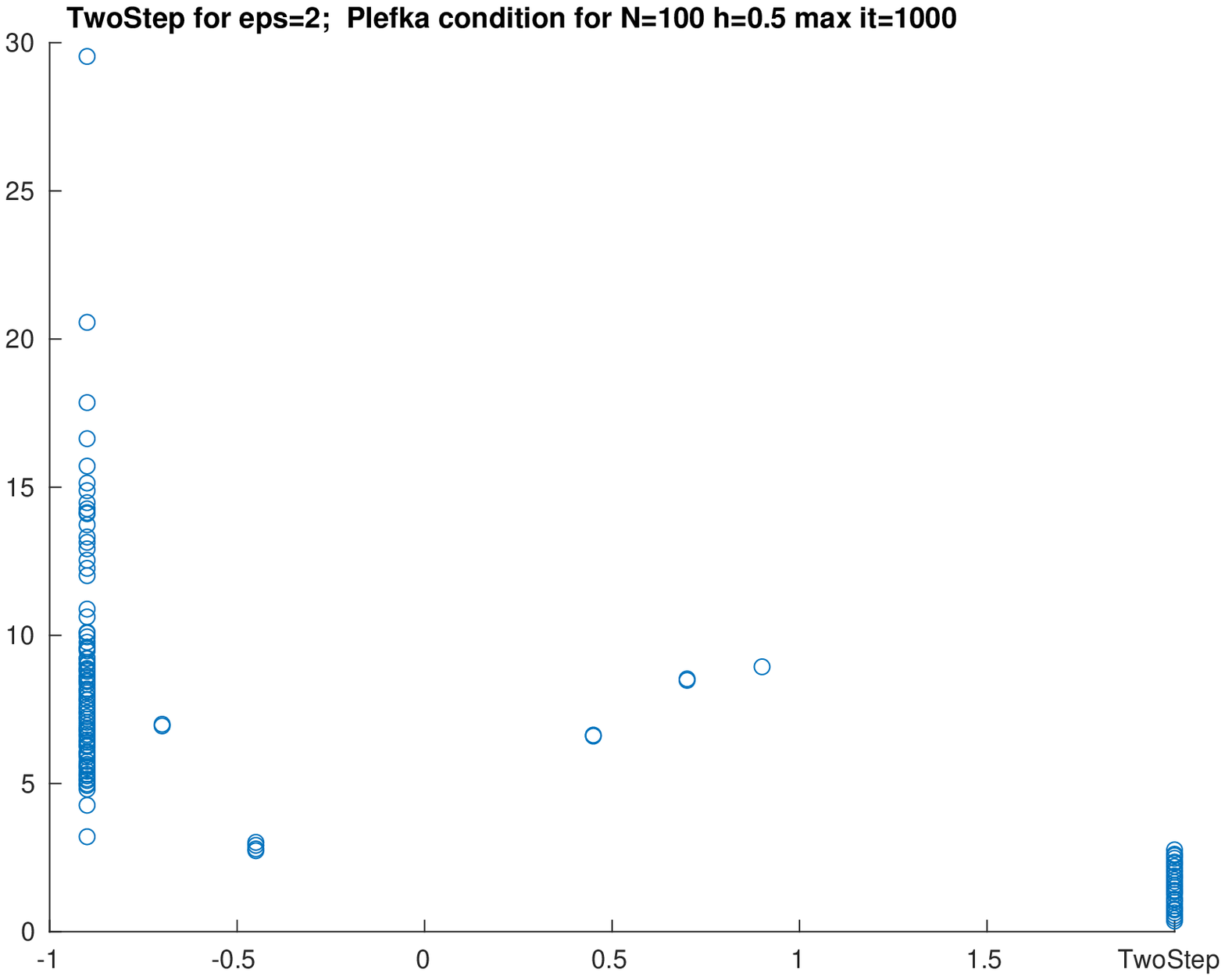}
  \caption{{\tiny} \space  1000 iterations of the $\varepsilon$-Banach and the Two-Step-Banach algorithm for $N=100, \beta=3$ and $h=0.5$. We simulated 1000 uniformly in the hypercube distributed start values for $\varepsilon = -0.9, -0.7, -0.45, 0.45, 0.7, 0.9$ and Two-Step-Banach. We plotted the Plefka values for the given parameters and used the same 1000 start values for all $\varepsilon$ and Two-Step-Banach.}
\end{figure}
\begin{figure}[H]
  \centering
  % \psfrag{2}{\tiny Two Step}
  \includegraphics[width=0.5\textwidth,  trim={0 0 0 0.4cm},clip]{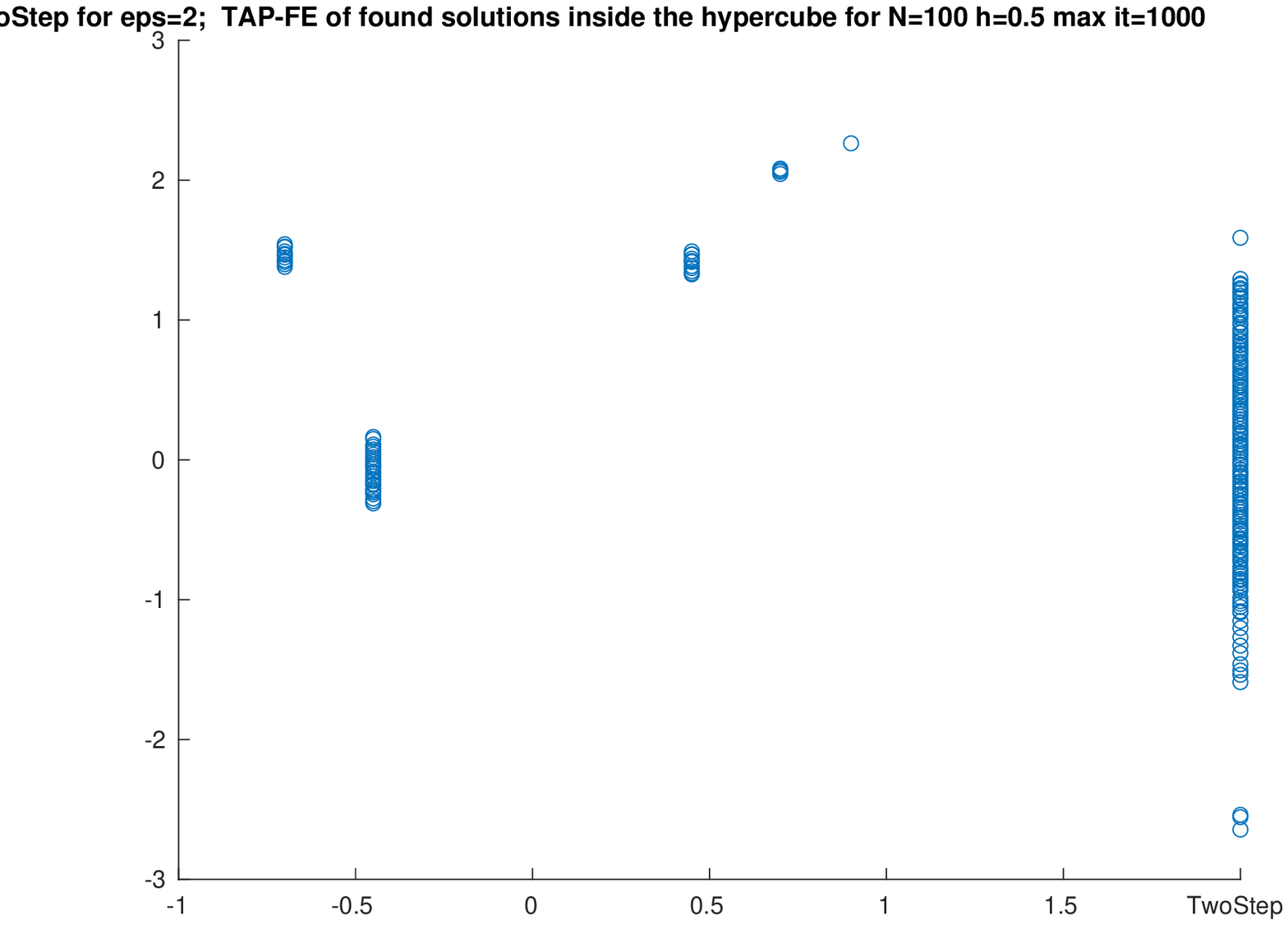}
  \caption{{\tiny} \space  1000 iterations of the $\varepsilon$-Banach and the Two-Step-Banach algorithm for $N=100, \beta=3$ and $h=0.5$. We simulated 1000 uniformly in the hypercube distributed start values for $\varepsilon = -0.9, -0.7, -0.45, 0.45, 0.7, 0.9$ and Two-Step-Banach. We plotted the TAP-free-energy values for all iterations ending inside the hypercube for the given parameters and used the same 1000 start values for all $\varepsilon$ and Two-Step-Banach.}
\end{figure}
\FloatBarrier\clearpage
\subsubsection{$N=500$,  $\beta =3 $ and $h=0.5$}\leavevmode
For all of the following figures five we used the same start values:
\begin{figure}[H]
  \centering
  % \psfrag{2}{\tiny Two Step}
  \includegraphics[width=0.5\textwidth,  trim={0 0 0 0.4cm},clip]{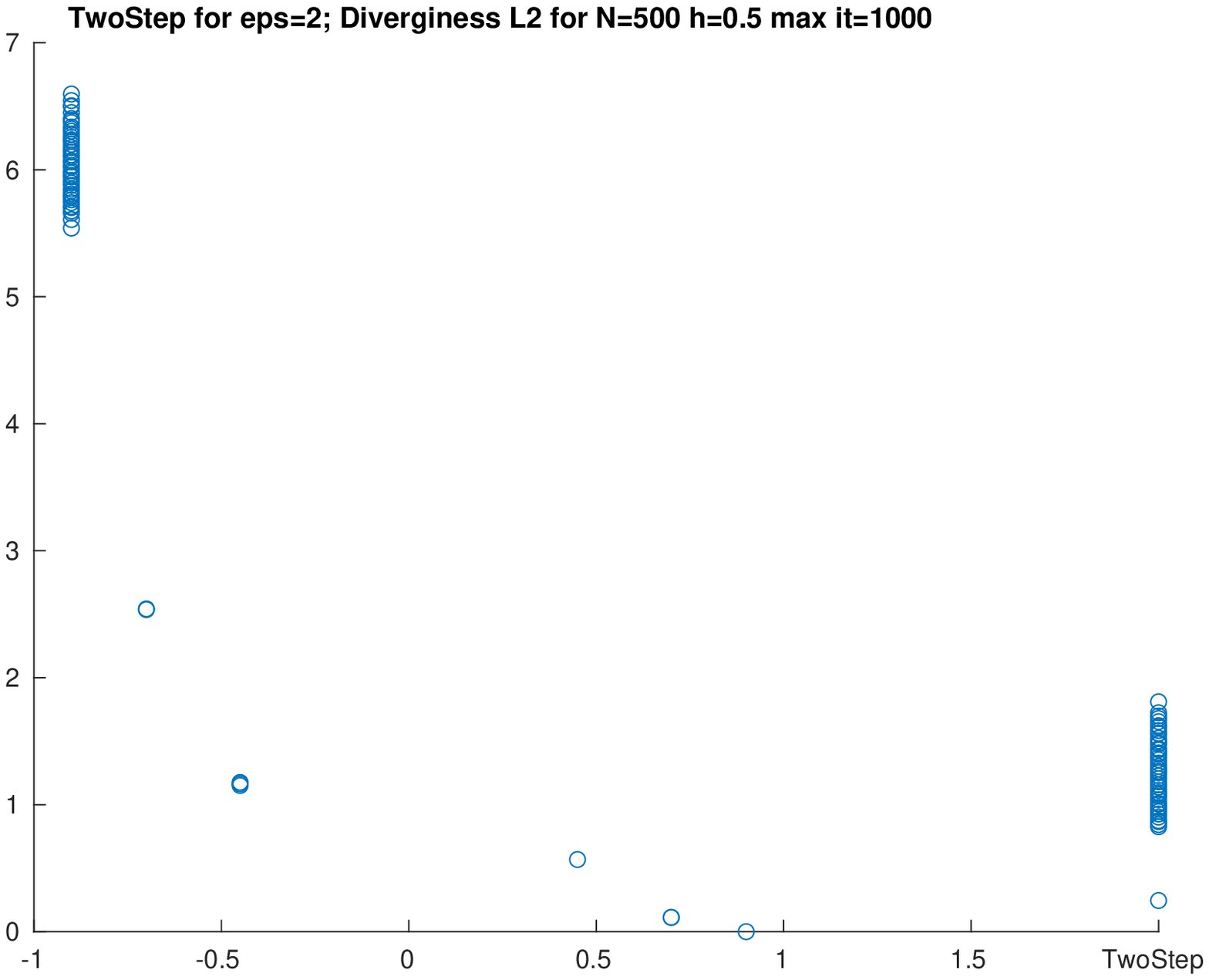}
  \caption{{\tiny} \space  1000 iterations of the $\varepsilon$-Banach and the Two-Step-Banach algorithm for $N=500, \beta=3$ and $h=0.5$. We simulated 1000 uniformly in the hypercube distributed start values for $\varepsilon = -0.9, -0.7, -0.45, 0.45, 0.7, 0.9$ and Two-Step-Banach. We plotted the {\it mean squared errors} for the given parameters and used the same 1000 start values for all $\varepsilon$ and Two-Step-Banach.}
\end{figure}
\begin{figure}[H]
  \centering
  % \psfrag{2}{\tiny Two Step}
  \includegraphics[width=0.5\textwidth,  trim={0 0 0 0.4cm},clip]{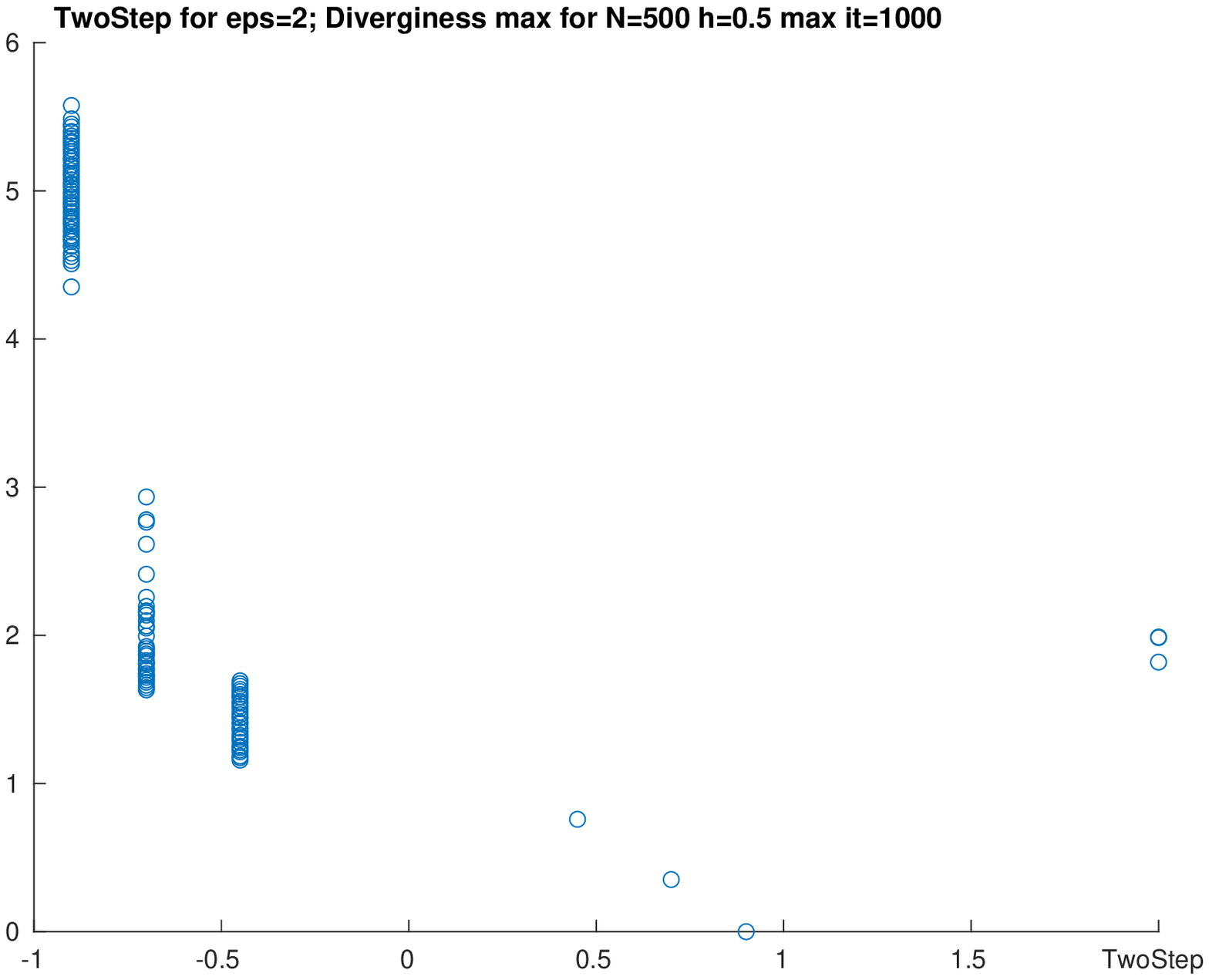}
  \caption{{\tiny} \space  1000 iterations of the $\varepsilon$-Banach and the Two-Step-Banach algorithm for $N=500, \beta=3$ and $h=0.5$. We simulated 1000 uniformly in the hypercube distributed start values for $\varepsilon = -0.9, -0.7, -0.45, 0.45, 0.7, 0.9$ and Two-Step-Banach. We plotted the {\it maximum absolute errors} for the given parameters and used the same 1000 start values for all $\varepsilon$ and Two-Step-Banach.}
\end{figure}
\begin{figure}[H]
  \centering
  % \psfrag{2}{\tiny Two Step}
  \includegraphics[width=0.5\textwidth,  trim={0 0 0 0.4cm},clip]{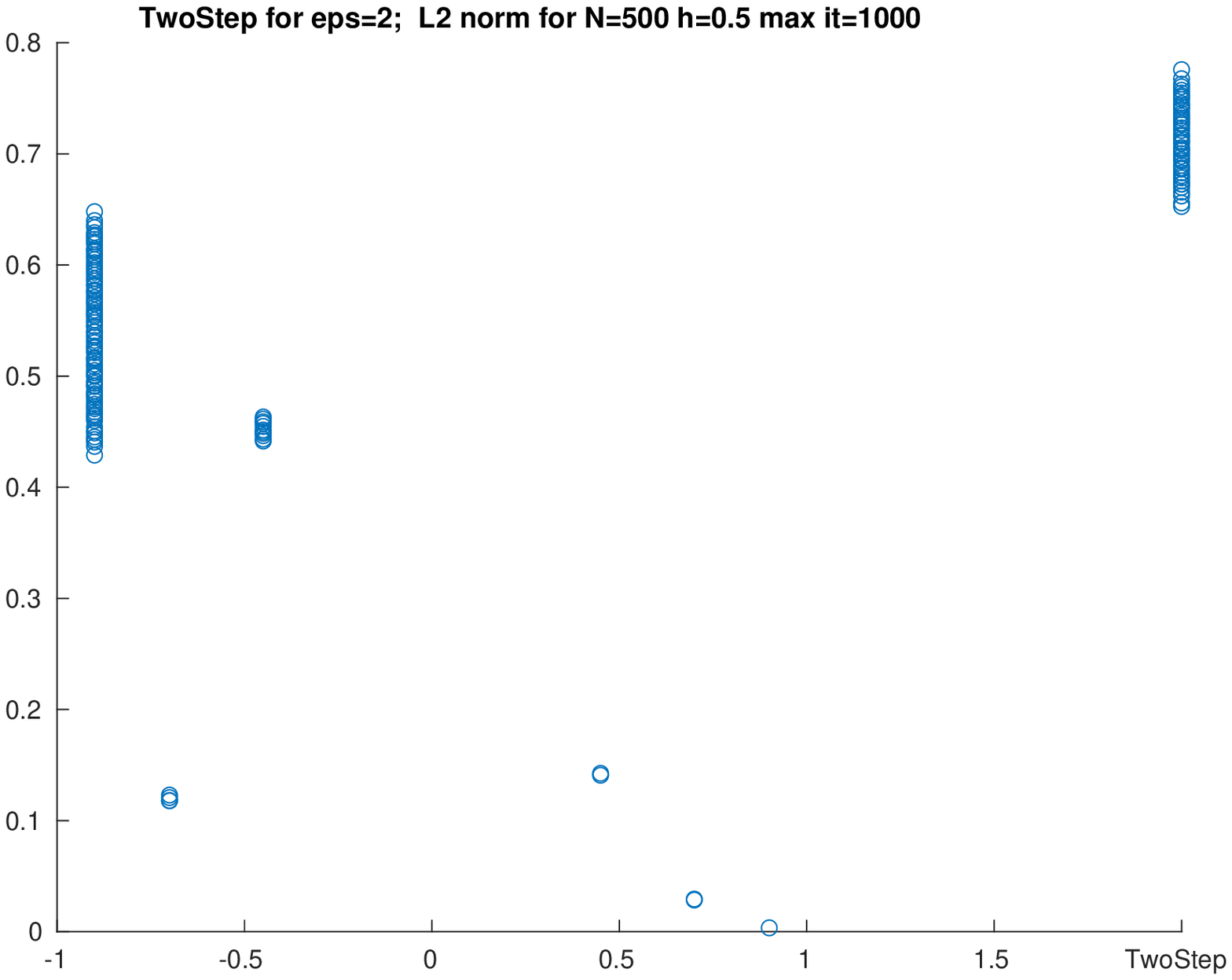}
  \caption{{\tiny} \space  1000 iterations of the $\varepsilon$-Banach and the Two-Step-Banach algorithm for $N=500, \beta=3$ and $h=0.5$. We simulated 1000 uniformly in the hypercube distributed start values for $\varepsilon = -0.9, -0.7, -0.45, 0.45, 0.7, 0.9$ and Two-Step-Banach. We plotted the {\it mean squared errors} for the given parameters and used the same 1000 start values for all $\varepsilon$ and Two-Step-Banach.}
\end{figure}
\begin{figure}[H]
  \centering
  % \psfrag{2}{\tiny Two Step}
  \includegraphics[width=0.5\textwidth,  trim={0 0 0 0.4cm},clip]{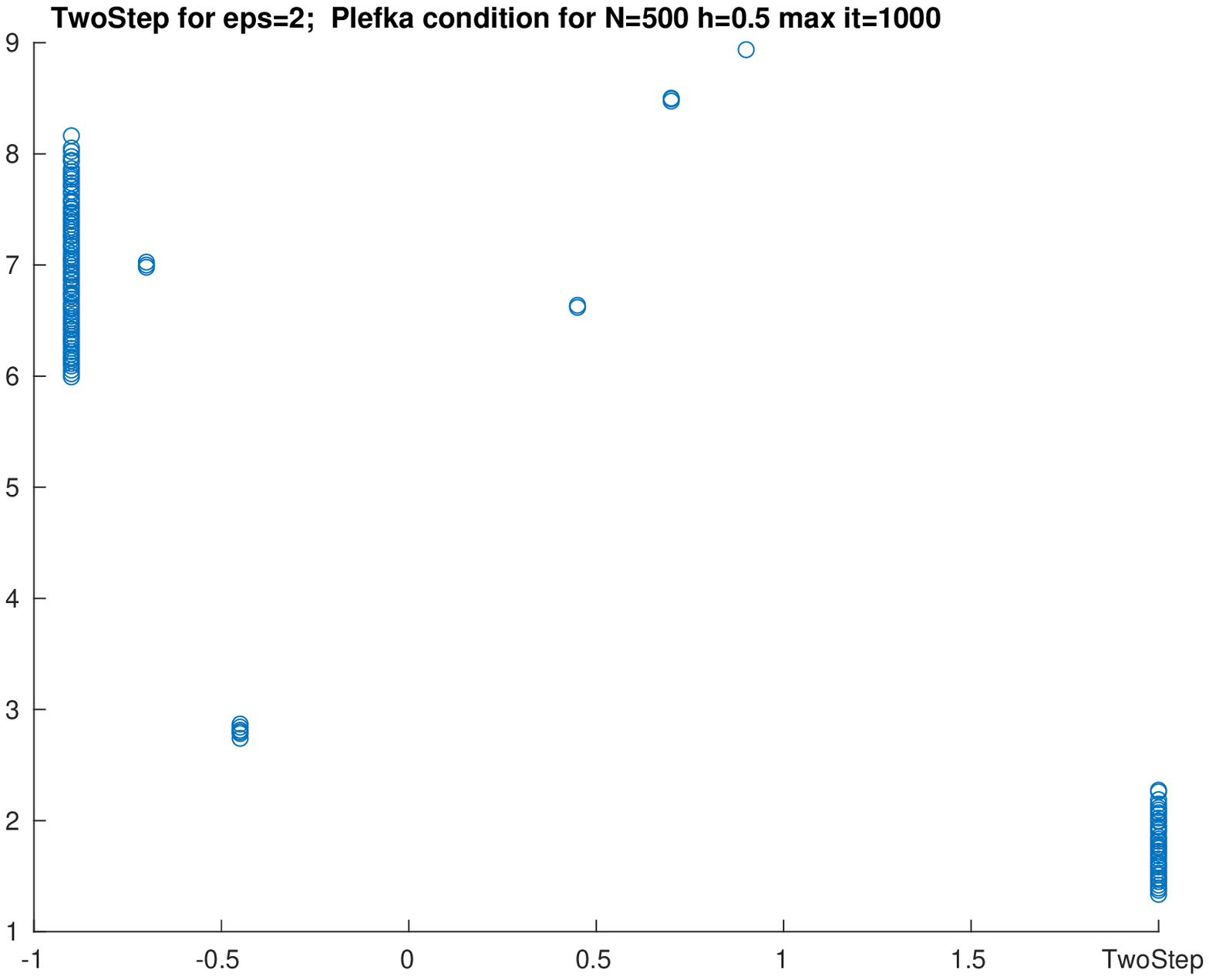}
  \caption{{\tiny} \space  1000 iterations of the $\varepsilon$-Banach and the Two-Step-Banach algorithm for $N=500, \beta=3$ and $h=0.5$. We simulated 1000 uniformly in the hypercube distributed start values for $\varepsilon = -0.9, -0.7, -0.45, 0.45, 0.7, 0.9$ and Two-Step-Banach. We plotted the Plefka values for the given parameters and used the same 1000 start values for all $\varepsilon$ and Two-Step-Banach.}
\end{figure}
\begin{figure}[H]
  \centering
  % \psfrag{2}{\tiny Two Step}
  \includegraphics[width=0.5\textwidth,  trim={0 0 0 0.4cm},clip]{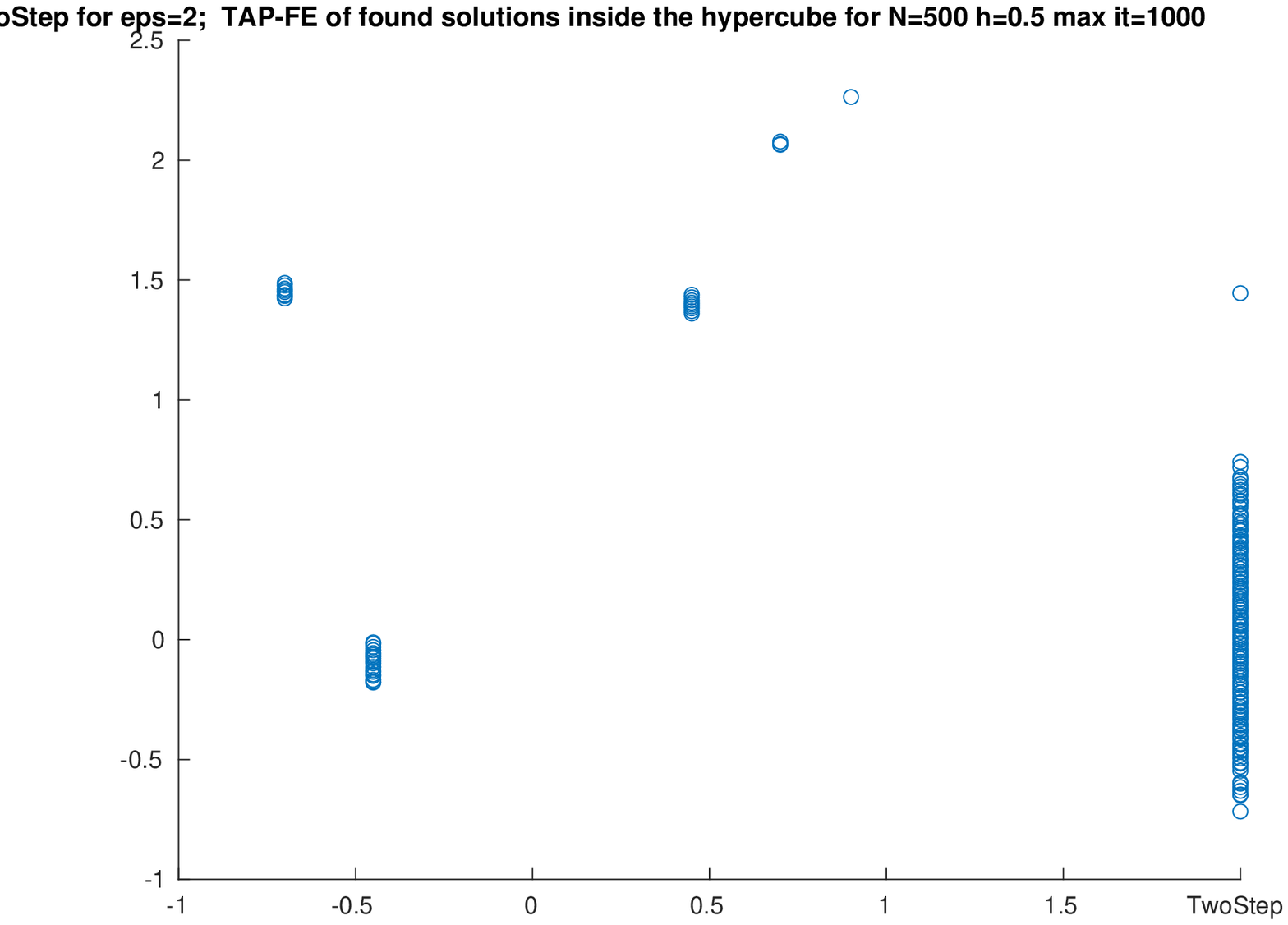}
  \caption{{\tiny} \space  1000 iterations of the $\varepsilon$-Banach and the Two-Step-Banach algorithm for $N=500, \beta=3$ and $h=0.5$. We simulated 1000 uniformly in the hypercube distributed start values for $\varepsilon = -0.9, -0.7, -0.45, 0.45, 0.7, 0.9$ and Two-Step-Banach. We plotted the TAP-free-energy values for all iterations ending inside the hypercube for the given parameters and used the same 1000 start values for all $\varepsilon$ and Two-Step-Banach.}
\end{figure}
\FloatBarrier\clearpage
\subsection{For $\beta =3 $ and $h=0.5$ dependence on the start values for  uniformly out of the corners of the hypercube chosen ones}
\subsubsection{$N=10$,  $\beta =3 $ and $h=0.5$}\leavevmode
For all of the following figures five we used the same start values:
\begin{figure}[H]
  \centering
  % \psfrag{2}{\tiny Two Step}
  \includegraphics[width=0.5\textwidth,  trim={0.45cm 0.4cm 0 0.4cm},clip]{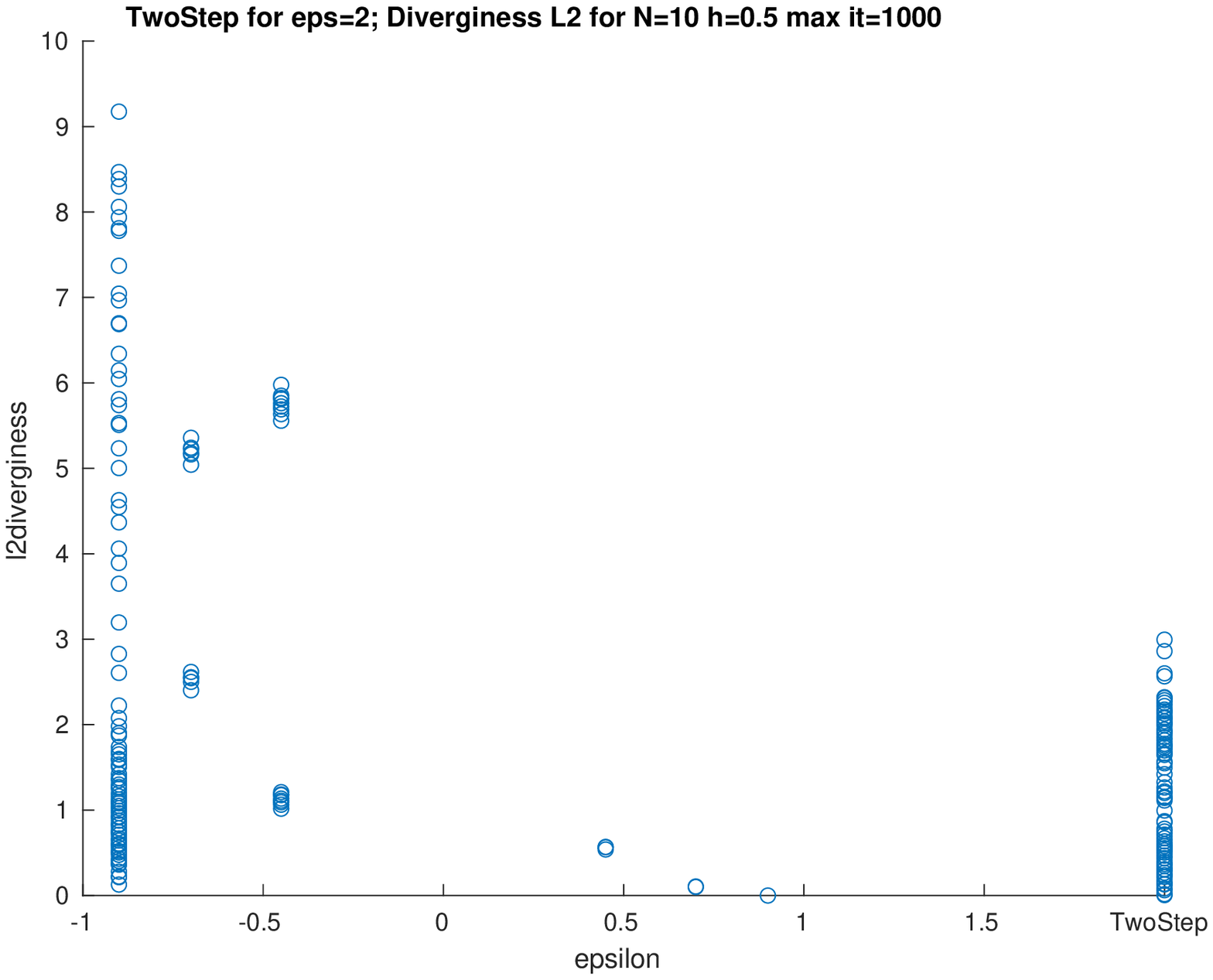}
  \caption{{\tiny} \space  1000 iterations of the $\varepsilon$-Banach and the Two-Step-Banach algorithm for $N=10, \beta=3$ and $h=0.5$. We simulated 1000 uniformly in the corners of the hypercube distributed start values for $\varepsilon = -0.9, -0.7, -0.45, 0.45, 0.7, 0.9$ and Two-Step-Banach. We plotted the {\it mean squared errors} for the given parameters and used the same 1000 start values for all $\varepsilon$ and Two-Step-Banach.}
\end{figure}
\begin{figure}[H]
  \centering
  % \psfrag{2}{\tiny Two Step}
  \includegraphics[width=0.5\textwidth,  trim={0.45cm 0.4cm 0 0.4cm},clip]{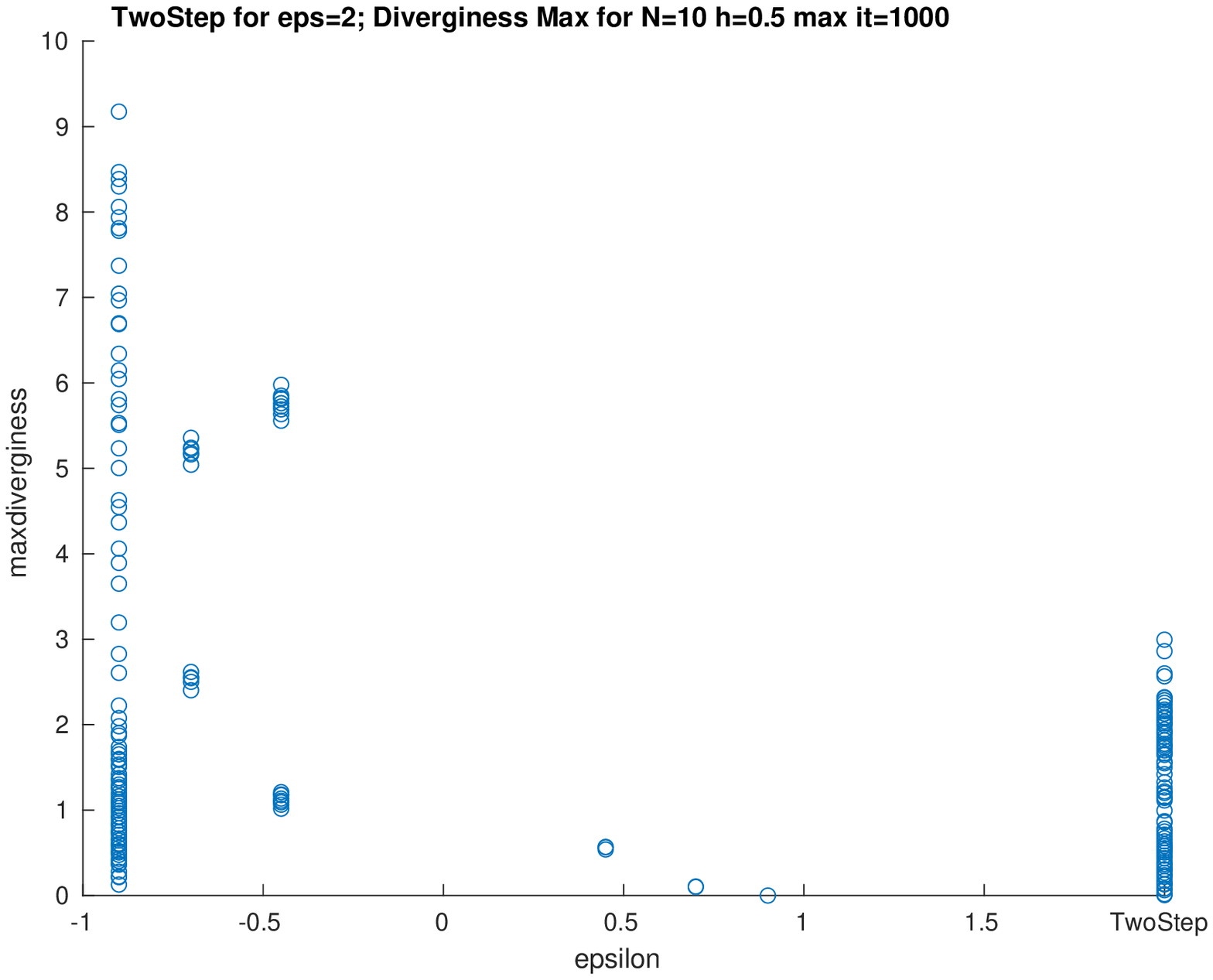}
  \caption{{\tiny} \space  1000 iterations of the $\varepsilon$-Banach and the Two-Step-Banach algorithm for $N=10, \beta=3$ and $h=0.5$. We simulated 1000 uniformly in the corners of the hypercube distributed start values for $\varepsilon = -0.9, -0.7, -0.45, 0.45, 0.7, 0.9$ and Two-Step-Banach. We plotted the {\it maximum absolute errors} for the given parameters and used the same 1000 start values for all $\varepsilon$ and Two-Step-Banach.}
\end{figure}
\begin{figure}[H]
  \centering
  % \psfrag{2}{\tiny Two Step}
  \includegraphics[width=0.5\textwidth,  trim={0.45cm 0.4cm 0 0.4cm},clip]{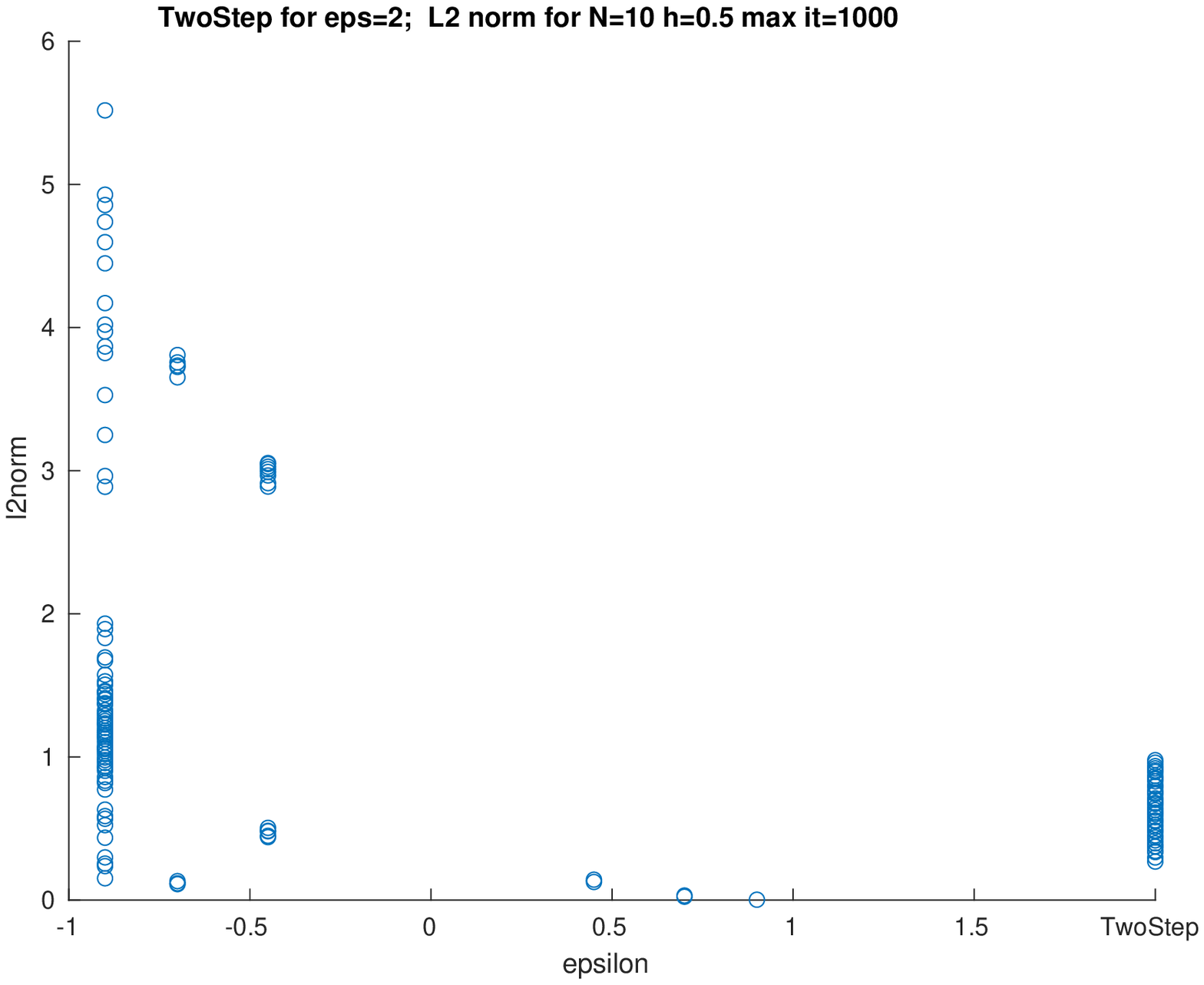}
  \caption{{\tiny} \space  1000 iterations of the $\varepsilon$-Banach and the Two-Step-Banach algorithm for $N=10, \beta=3$ and $h=0.5$. We simulated 1000 uniformly in the corners of the hypercube distributed start values for $\varepsilon = -0.9, -0.7, -0.45, 0.45, 0.7, 0.9$ and Two-Step-Banach. We plotted the {\it mean squared errors} for the given parameters and used the same 1000 start values for all $\varepsilon$ and Two-Step-Banach.}
\end{figure}
\begin{figure}[H]
  \centering
  % \psfrag{2}{\tiny Two Step}
  \includegraphics[width=0.5\textwidth,  trim={0.45cm 0.4cm 0 0},clip]{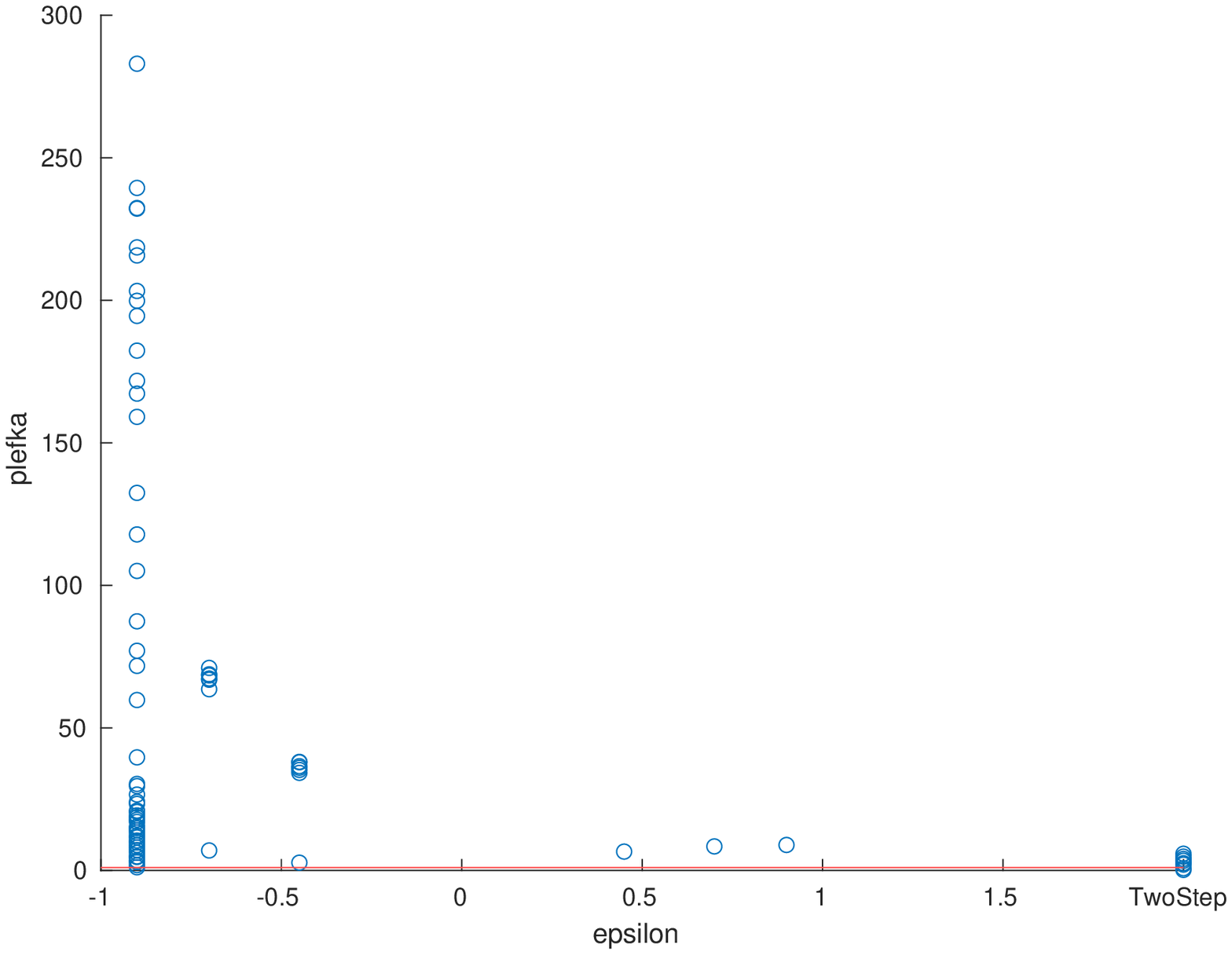}
  \caption{{\tiny} \space  1000 iterations of the $\varepsilon$-Banach and the Two-Step-Banach algorithm for $N=10, \beta=3$ and $h=0.5$. We simulated 1000 uniformly in the corners of the hypercube distributed start values for $\varepsilon = -0.9, -0.7, -0.45, 0.45, 0.7, 0.9$ and Two-Step-Banach. We plotted the Plefka values for the given parameters and used the same 1000 start values for all $\varepsilon$ and Two-Step-Banach.}
\end{figure}
\begin{figure}[H]
  \centering
  % \psfrag{2}{\tiny Two Step}
  \includegraphics[width=0.5\textwidth,  trim={0.45cm 0.4cm 0 0},clip]{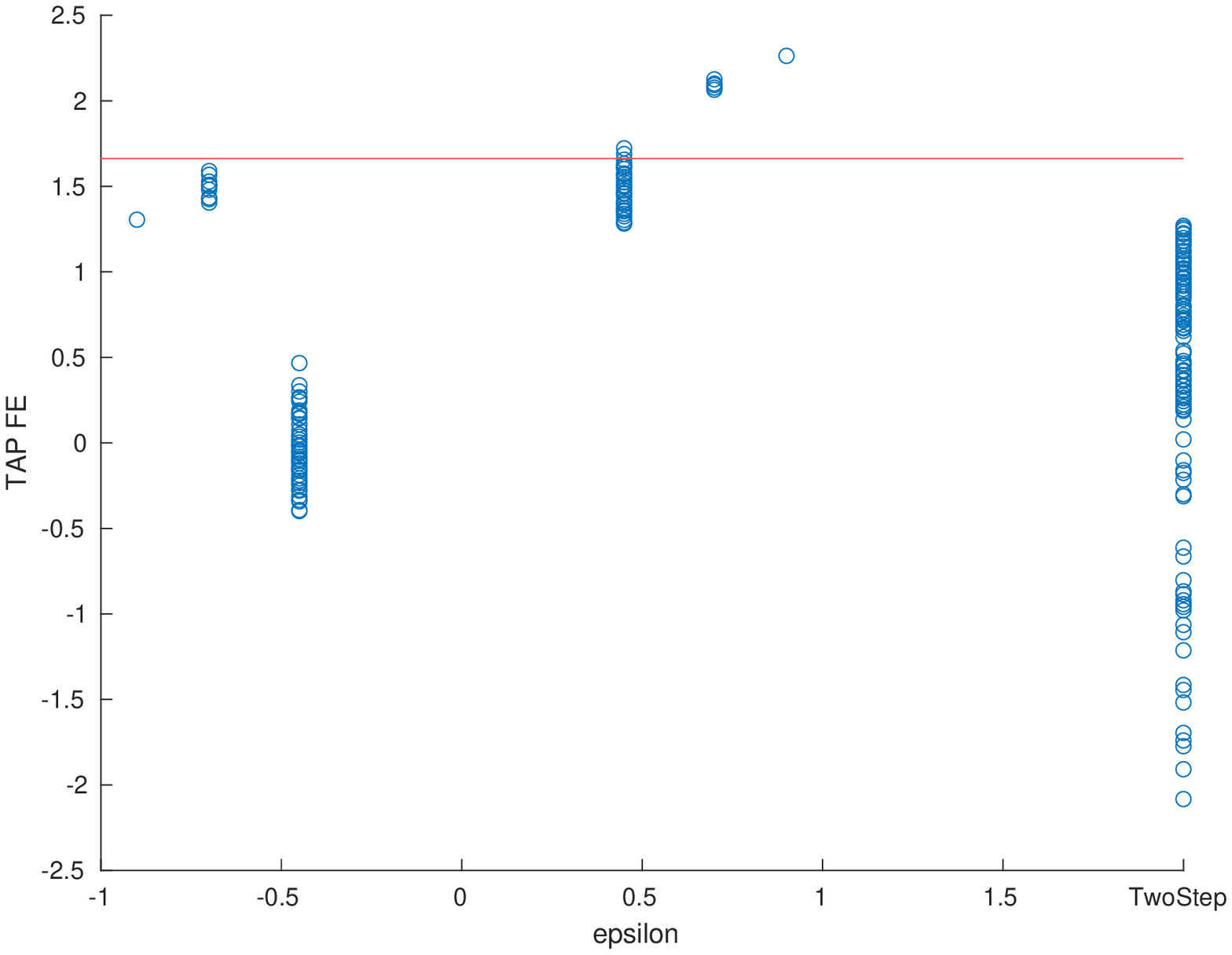}
  \caption{{\tiny} \space  1000 iterations of the $\varepsilon$-Banach and the Two-Step-Banach algorithm for $N=10, \beta=3$ and $h=0.5$. We simulated 1000 uniformly in the corners of the hypercube distributed start values for $\varepsilon = -0.9, -0.7, -0.45, 0.45, 0.7, 0.9$ and Two-Step-Banach. We plotted the TAP-free-energy values for all iterations ending inside the corners of the hypercube for the given parameters and used the same 1000 start values for all $\varepsilon$ and Two-Step-Banach.}
\end{figure}
\FloatBarrier\clearpage
\subsubsection{$N=25$,  $\beta =3 $ and $h=0.5$}\leavevmode
For all of the following figures five we used the same start values:
\begin{figure}[H]
  \centering
  % \psfrag{2}{\tiny Two Step}
  \includegraphics[width=0.5\textwidth,  trim={0.45cm 0.4cm 0 0.4cm},clip]{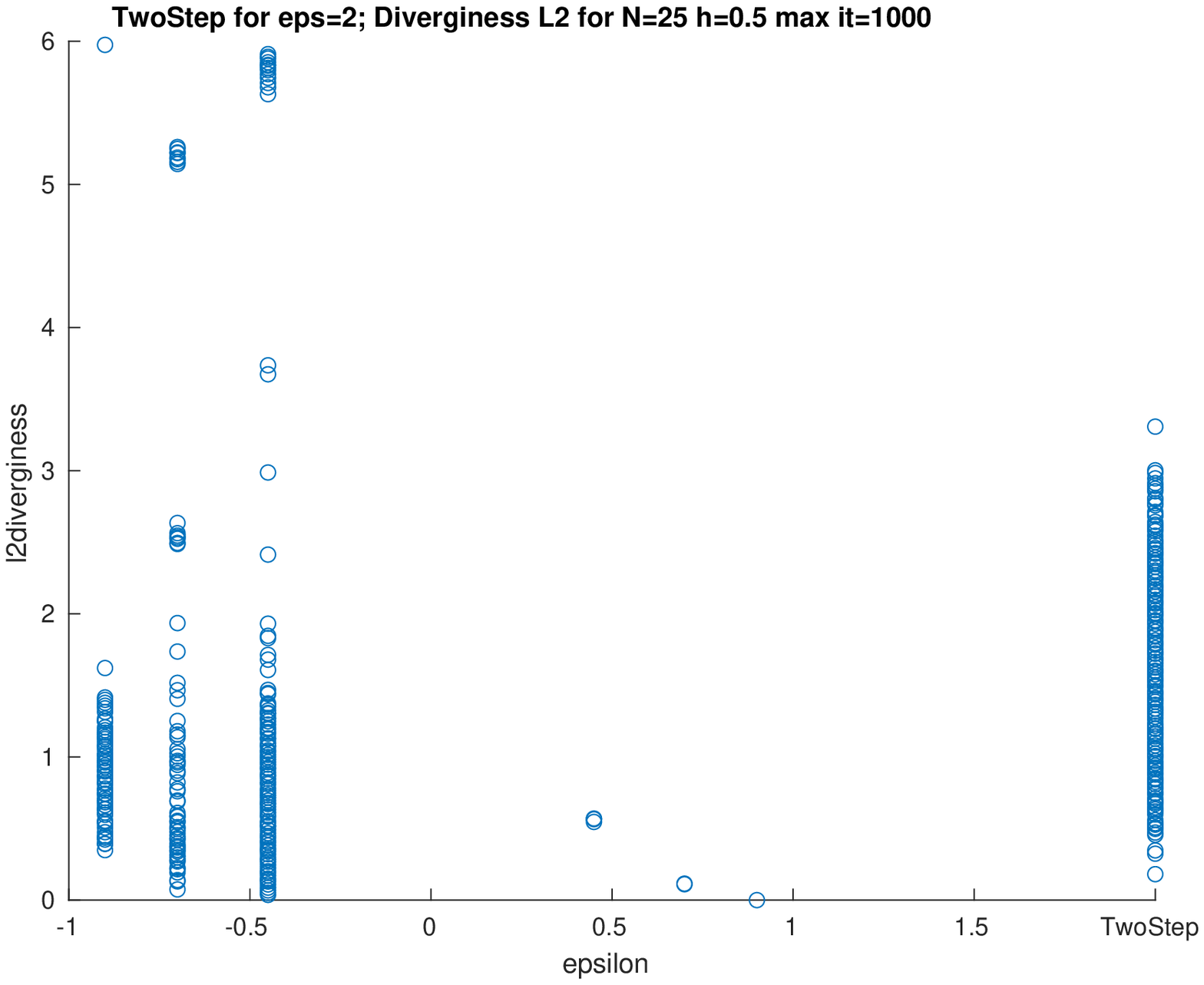}
  \caption{{\tiny} \space  1000 iterations of the $\varepsilon$-Banach and the Two-Step-Banach algorithm for $N=25, \beta=3$ and $h=0.5$. We simulated 1000 uniformly in the corners of the hypercube distributed start values for $\varepsilon = -0.9, -0.7, -0.45, 0.45, 0.7, 0.9$ and Two-Step-Banach. We plotted the {\it mean squared errors} for the given parameters and used the same 1000 start values for all $\varepsilon$ and Two-Step-Banach.}
\end{figure}
\begin{figure}[H]
  \centering
  % \psfrag{2}{\tiny Two Step}
  \includegraphics[width=0.5\textwidth,  trim={0.45cm 0.4cm 0 0.4cm},clip]{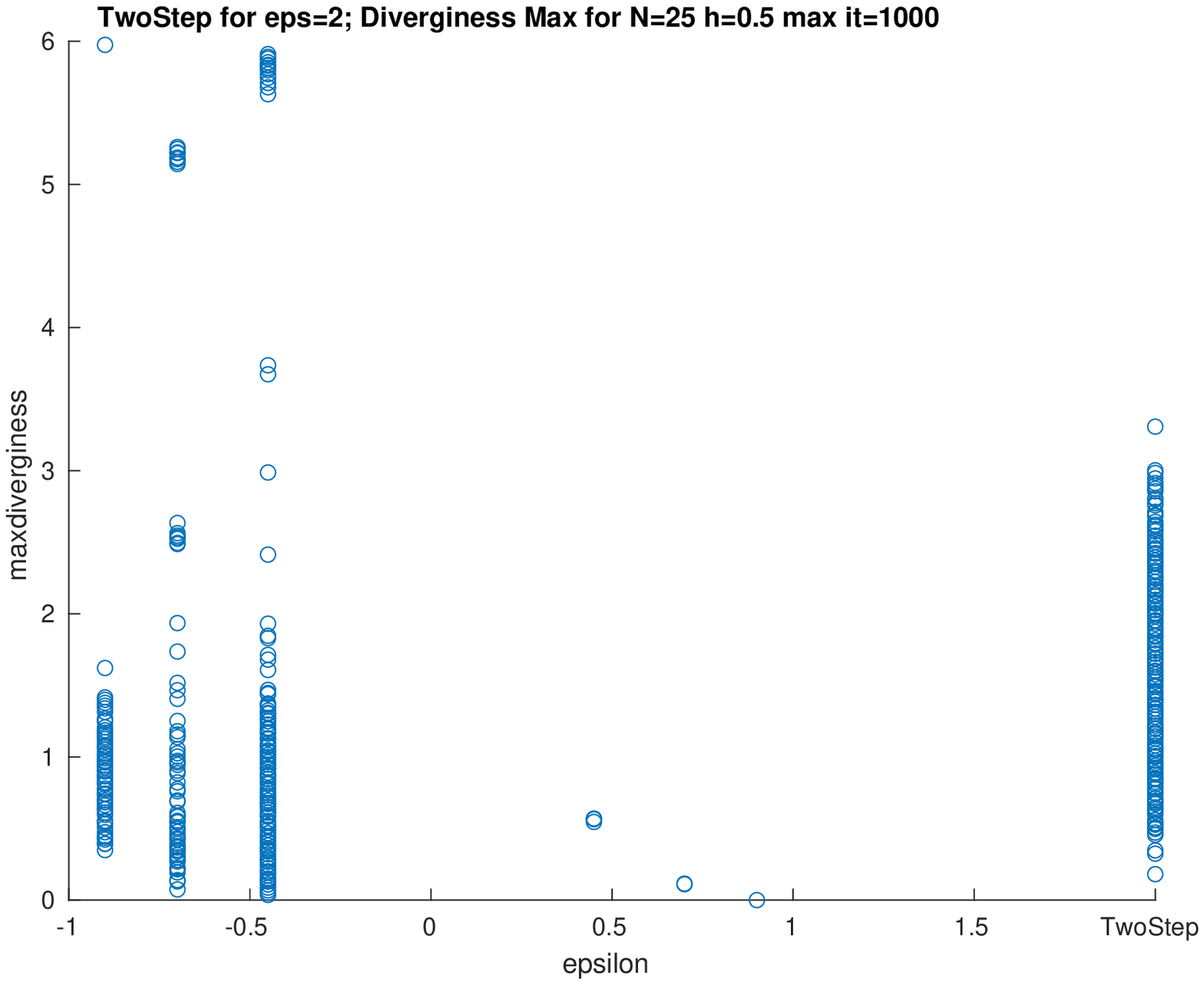}
  \caption{{\tiny} \space  1000 iterations of the $\varepsilon$-Banach and the Two-Step-Banach algorithm for $N=25, \beta=3$ and $h=0.5$. We simulated 1000 uniformly in the corners of the hypercube distributed start values for $\varepsilon = -0.9, -0.7, -0.45, 0.45, 0.7, 0.9$ and Two-Step-Banach. We plotted the {\it maximum absolute errors} for the given parameters and used the same 1000 start values for all $\varepsilon$ and Two-Step-Banach.}
\end{figure}
\begin{figure}[H]
  \centering
  % \psfrag{2}{\tiny Two Step}
  \includegraphics[width=0.5\textwidth,  trim={0.45cm 0.4cm 0 0.4cm},clip]{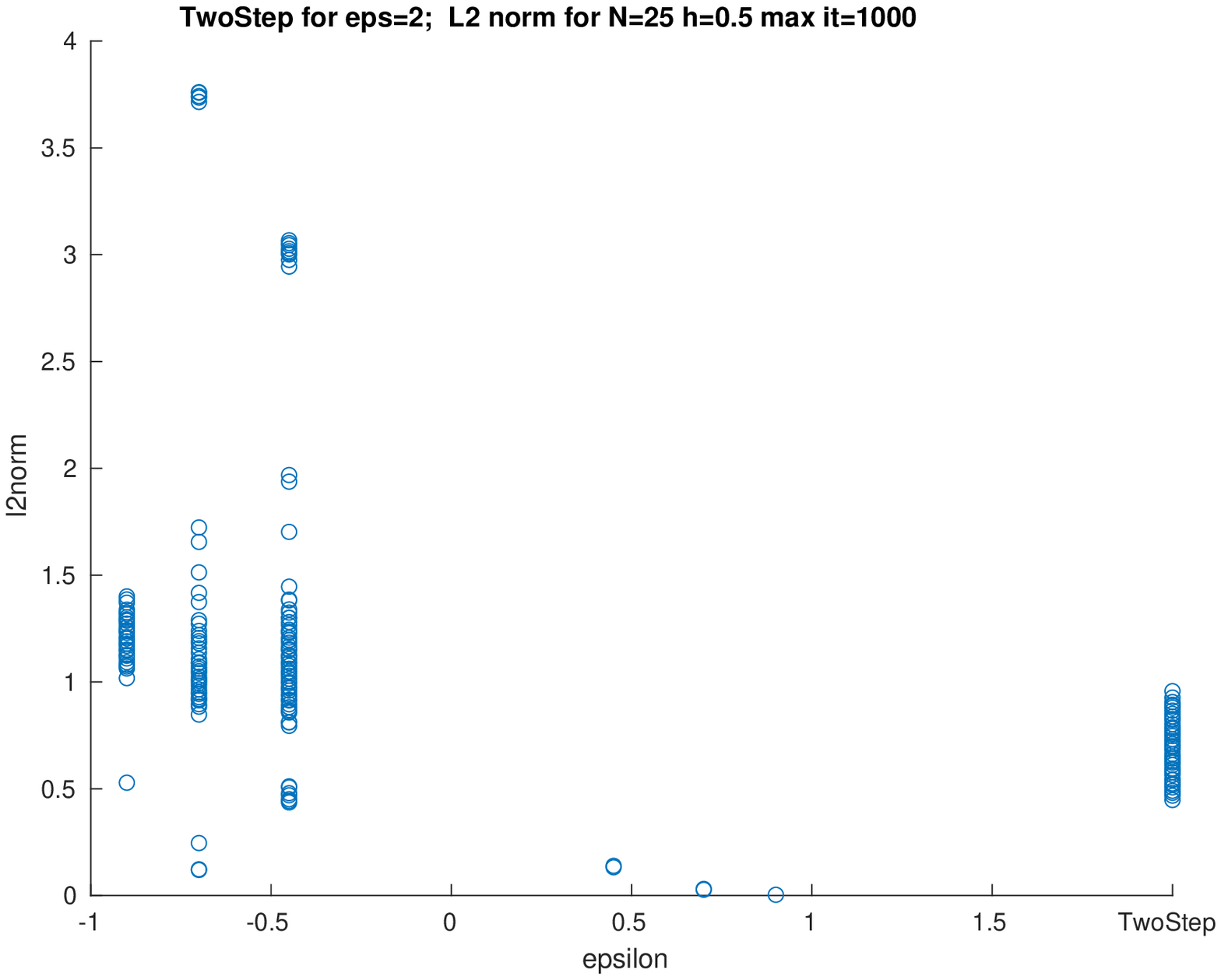}
  \caption{{\tiny} \space  1000 iterations of the $\varepsilon$-Banach and the Two-Step-Banach algorithm for $N=25, \beta=3$ and $h=0.5$. We simulated 1000 uniformly in the corners of the hypercube distributed start values for $\varepsilon = -0.9, -0.7, -0.45, 0.45, 0.7, 0.9$ and Two-Step-Banach. We plotted the {\it mean squared errors} for the given parameters and used the same 1000 start values for all $\varepsilon$ and Two-Step-Banach.}
\end{figure}
\begin{figure}[H]
  \centering
  % \psfrag{2}{\tiny Two Step}
  \includegraphics[width=0.5\textwidth,  trim={0.45cm 0.4cm 0 0},clip]{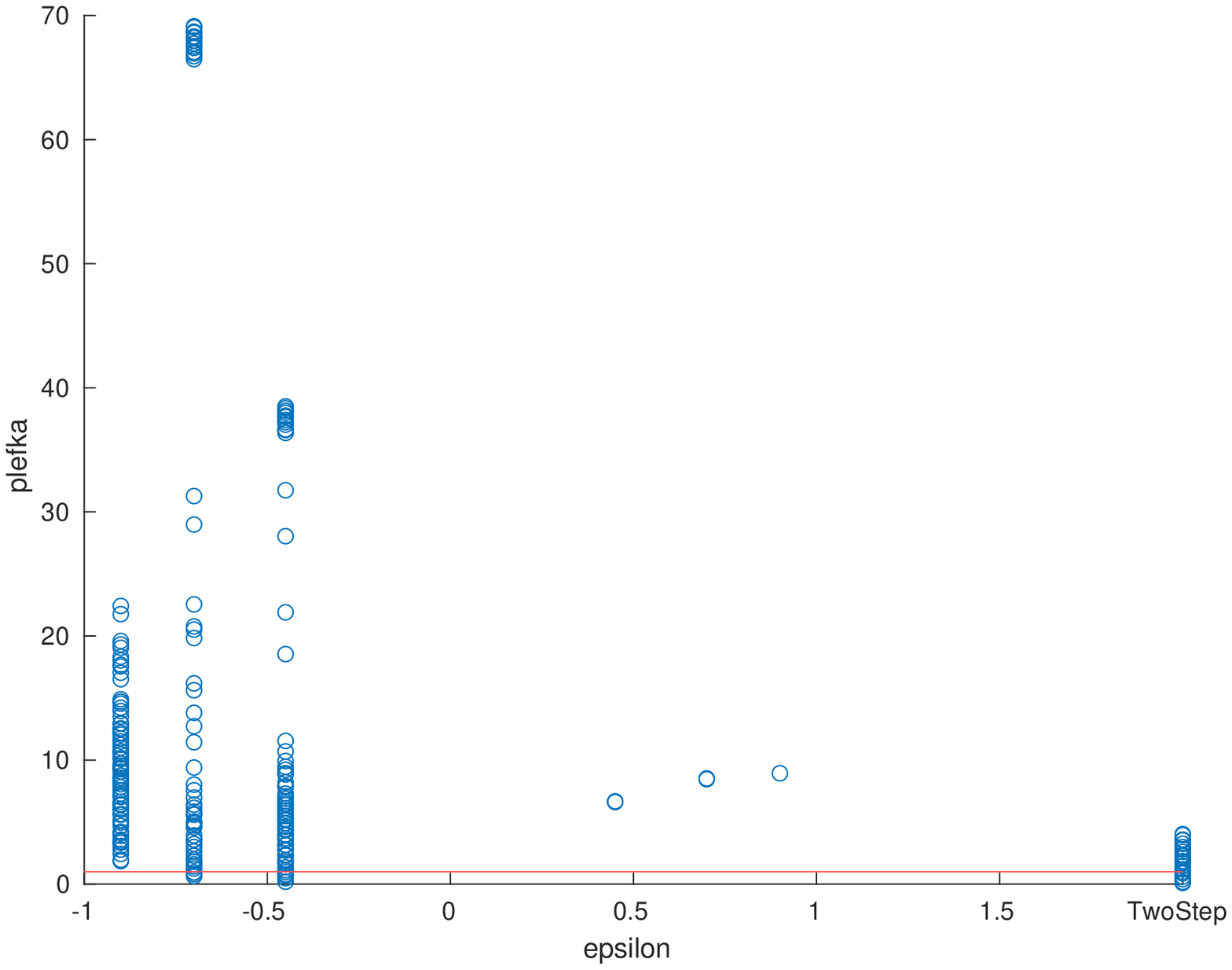}
  \caption{{\tiny} \space  1000 iterations of the $\varepsilon$-Banach and the Two-Step-Banach algorithm for $N=25, \beta=3$ and $h=0.5$. We simulated 1000 uniformly in the corners of the hypercube distributed start values for $\varepsilon = -0.9, -0.7, -0.45, 0.45, 0.7, 0.9$ and Two-Step-Banach. We plotted the Plefka values for the given parameters and used the same 1000 start values for all $\varepsilon$ and Two-Step-Banach.}
\end{figure}
\begin{figure}[H]
  \centering
  % \psfrag{2}{\tiny Two Step}
  \includegraphics[width=0.5\textwidth,  trim={0.45cm 0.4cm 0 0},clip]{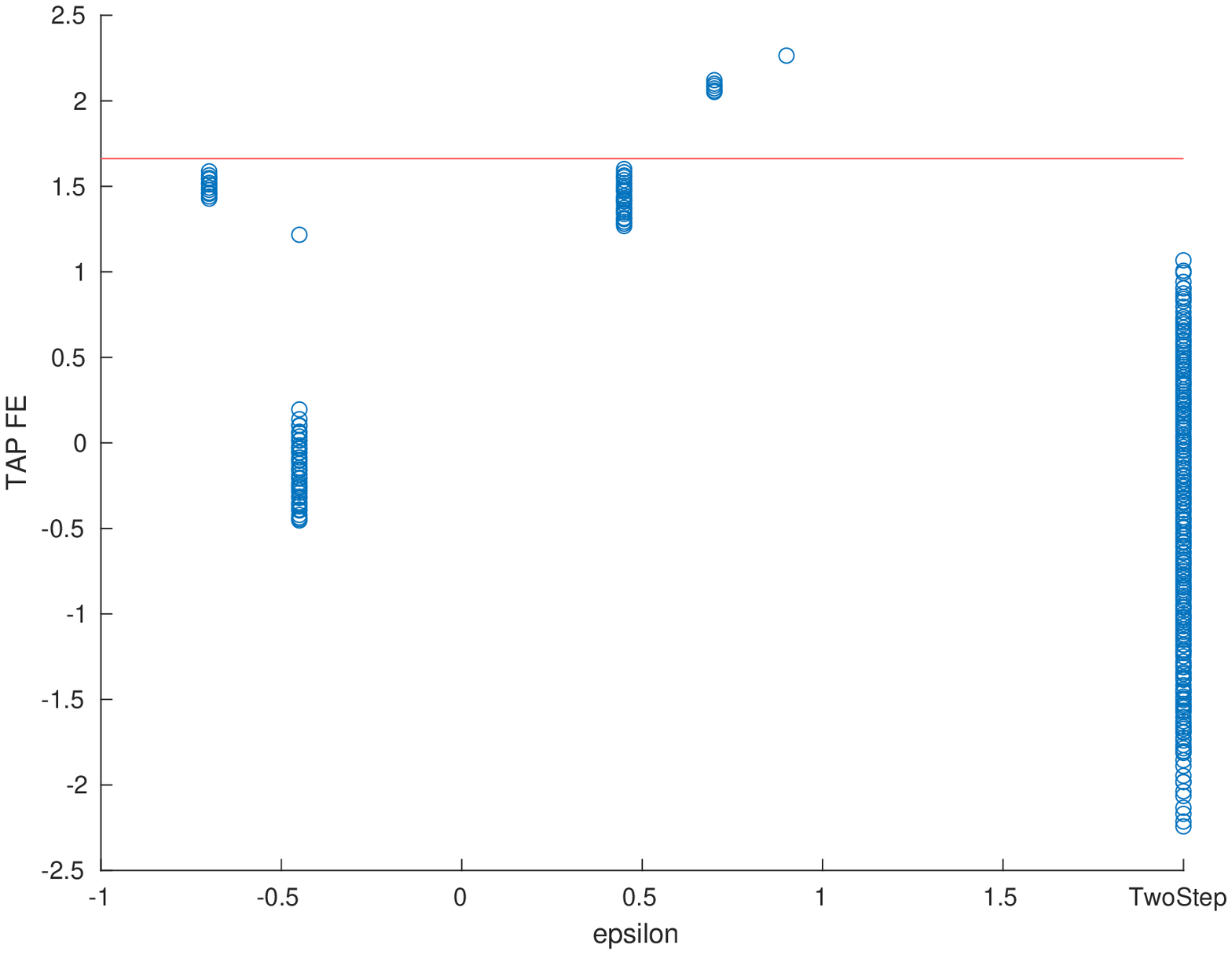}
  \caption{{\tiny} \space  1000 iterations of the $\varepsilon$-Banach and the Two-Step-Banach algorithm for $N=25, \beta=3$ and $h=0.5$. We simulated 1000 uniformly in the corners of the hypercube distributed start values for $\varepsilon = -0.9, -0.7, -0.45, 0.45, 0.7, 0.9$ and Two-Step-Banach. We plotted the TAP-free-energy values for all iterations ending inside the corners of the hypercube for the given parameters and used the same 1000 start values for all $\varepsilon$ and Two-Step-Banach.}
\end{figure}
\FloatBarrier\clearpage
\subsubsection{$N=100$,  $\beta =3 $ and $h=0.5$}\leavevmode
For all of the following figures five we used the same start values:
\begin{figure}[H]
  \centering
  % \psfrag{2}{\tiny Two Step}
  \includegraphics[width=0.5\textwidth,  trim={0.45cm 0.4cm 0 0.4cm},clip]{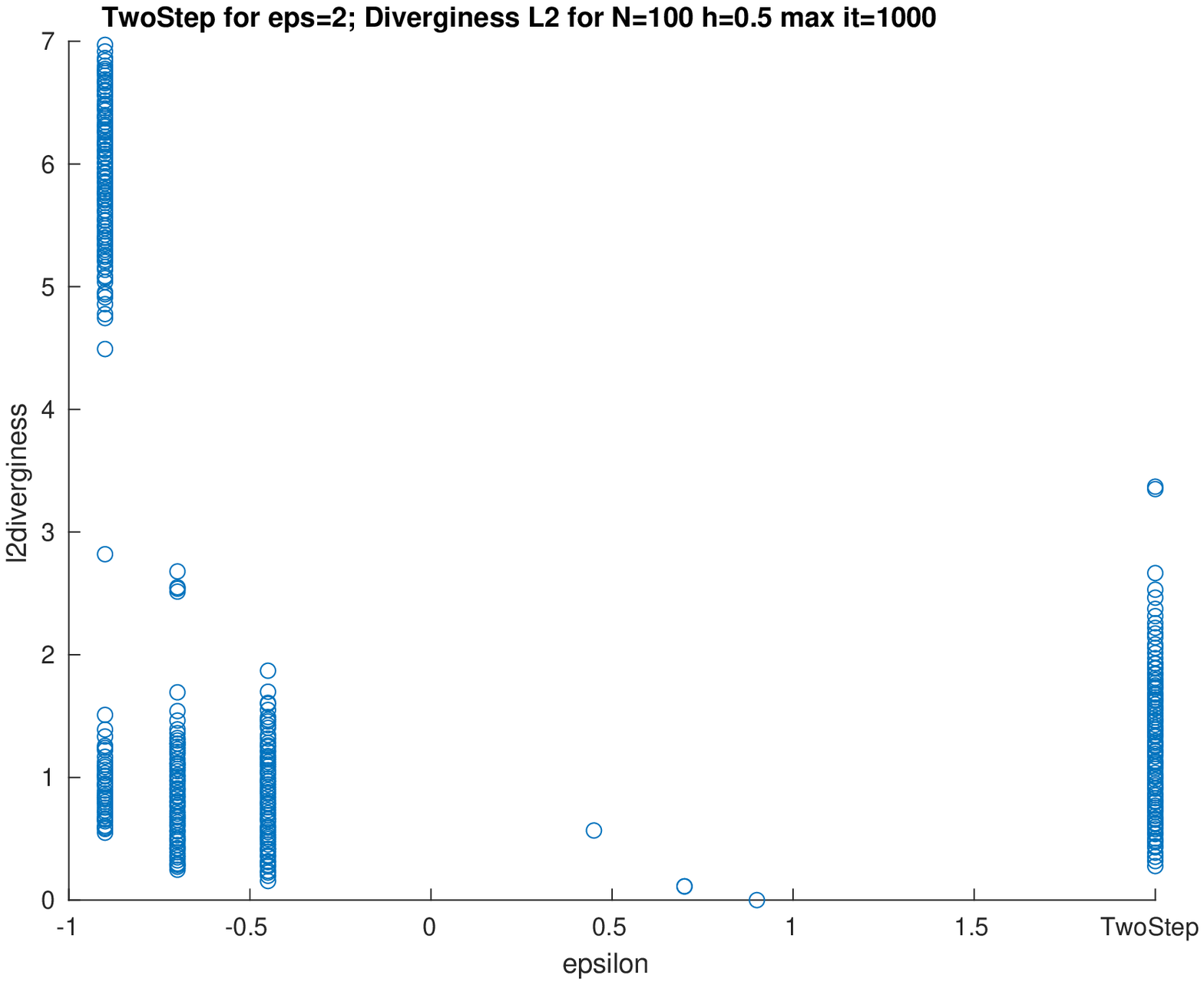}
  \caption{{\tiny} \space  1000 iterations of the $\varepsilon$-Banach and the Two-Step-Banach algorithm for $N=100, \beta=3$ and $h=0.5$. We simulated 1000 uniformly in the corners of the hypercube distributed start values for $\varepsilon = -0.9, -0.7, -0.45, 0.45, 0.7, 0.9$ and Two-Step-Banach. We plotted the {\it mean squared errors} for the given parameters and used the same 1000 start values for all $\varepsilon$ and Two-Step-Banach.}
\end{figure}
\begin{figure}[H]
  \centering
  % \psfrag{2}{\tiny Two Step}
  \includegraphics[width=0.5\textwidth,  trim={0.45cm 0.4cm 0 0.4cm},clip]{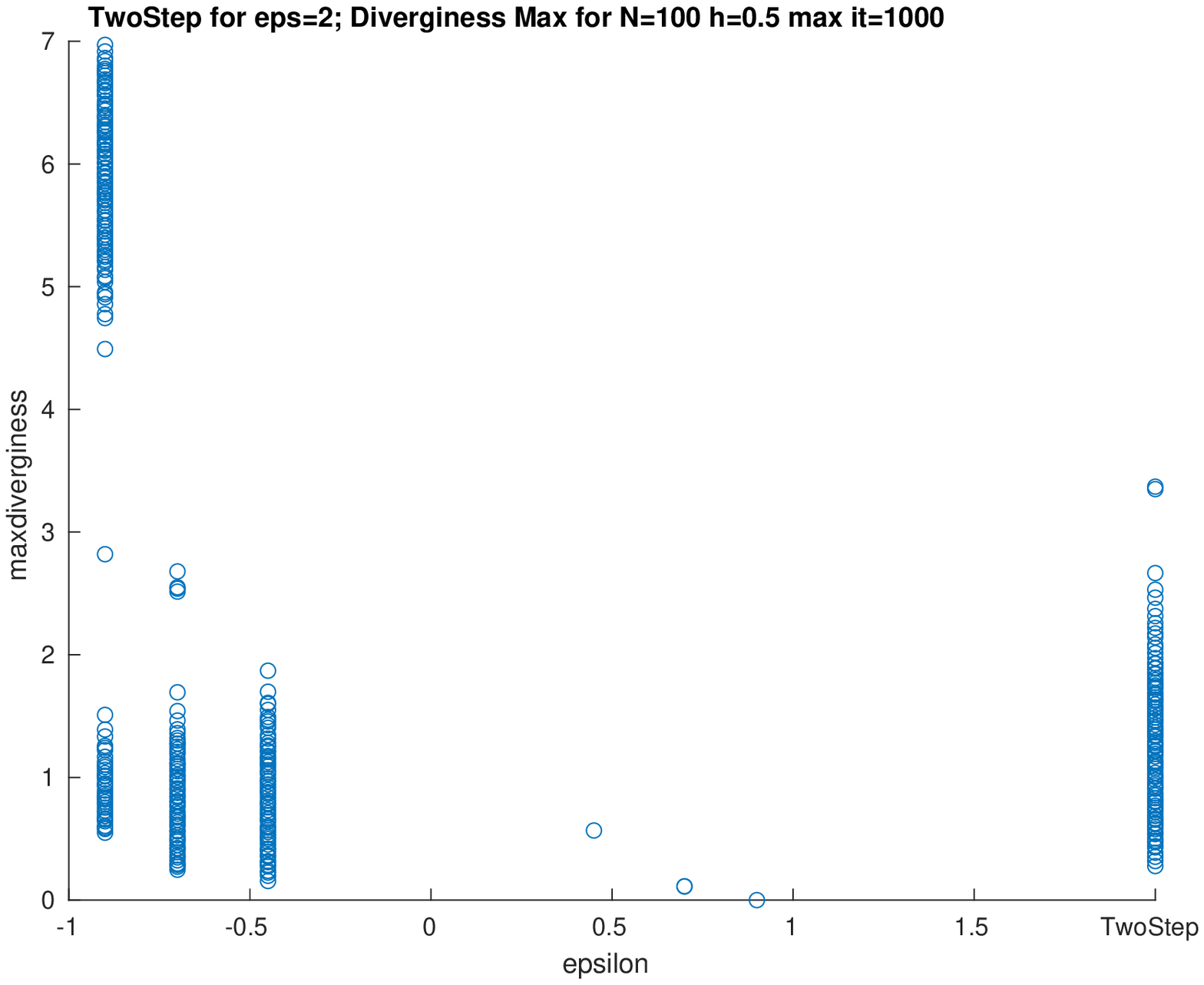}
  \caption{{\tiny} \space  1000 iterations of the $\varepsilon$-Banach and the Two-Step-Banach algorithm for $N=100, \beta=3$ and $h=0.5$. We simulated 1000 uniformly in the corners of the hypercube distributed start values for $\varepsilon = -0.9, -0.7, -0.45, 0.45, 0.7, 0.9$ and Two-Step-Banach. We plotted the {\it maximum absolute errors} for the given parameters and used the same 1000 start values for all $\varepsilon$ and Two-Step-Banach.}
\end{figure}
\begin{figure}[H]
  \centering
  % \psfrag{2}{\tiny Two Step}
  \includegraphics[width=0.5\textwidth,  trim={0.45cm 0.4cm 0 0.4cm},clip]{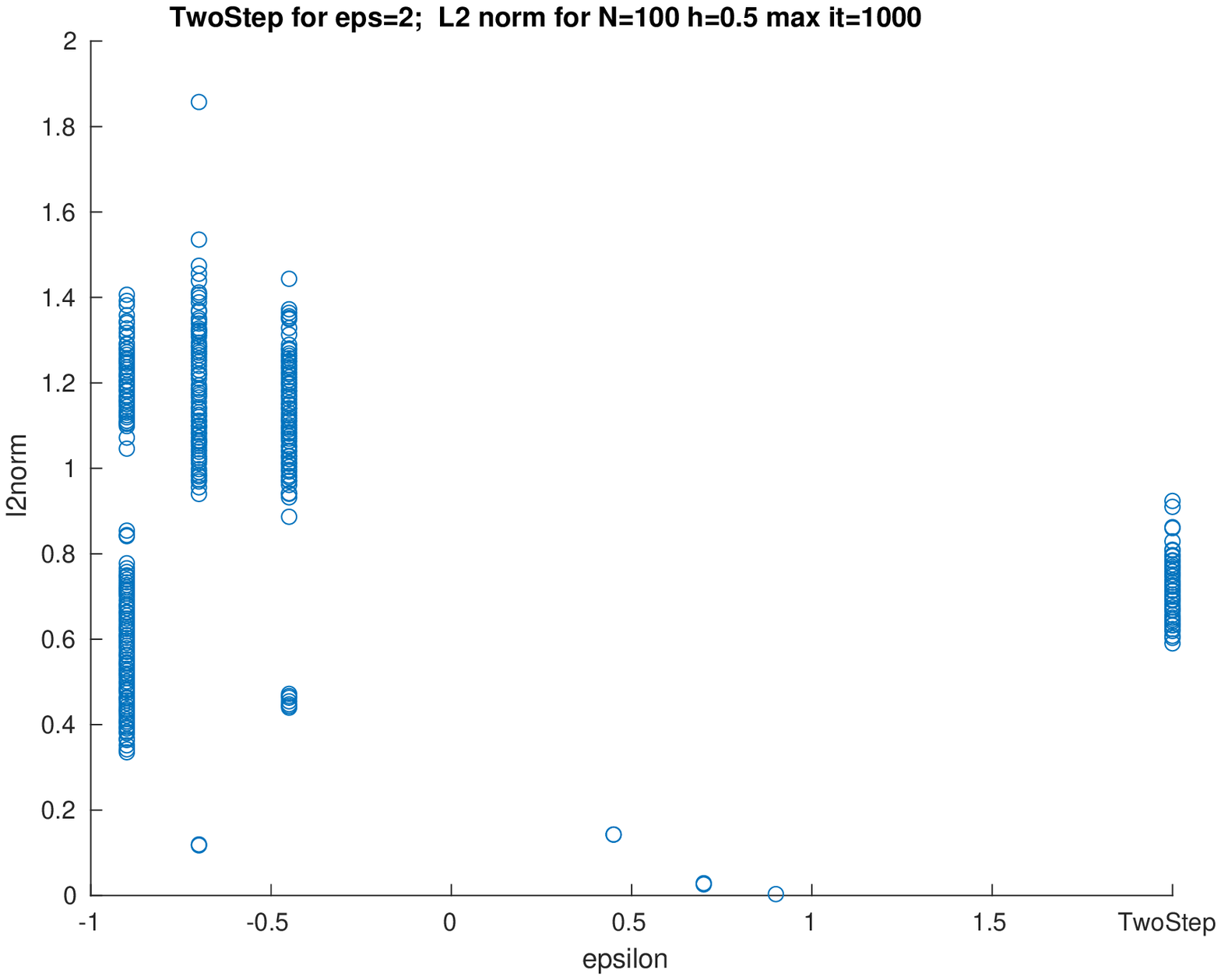}
  \caption{{\tiny} \space  1000 iterations of the $\varepsilon$-Banach and the Two-Step-Banach algorithm for $N=100, \beta=3$ and $h=0.5$. We simulated 1000 uniformly in the corners of the hypercube distributed start values for $\varepsilon = -0.9, -0.7, -0.45, 0.45, 0.7, 0.9$ and Two-Step-Banach. We plotted the {\it mean squared errorss} for the given parameters and used the same 1000 start values for all $\varepsilon$ and Two-Step-Banach.}
\end{figure}
\begin{figure}[H]
  \centering
  % \psfrag{2}{\tiny Two Step}
  \includegraphics[width=0.5\textwidth,  trim={0.45cm 0.4cm 0 0},clip]{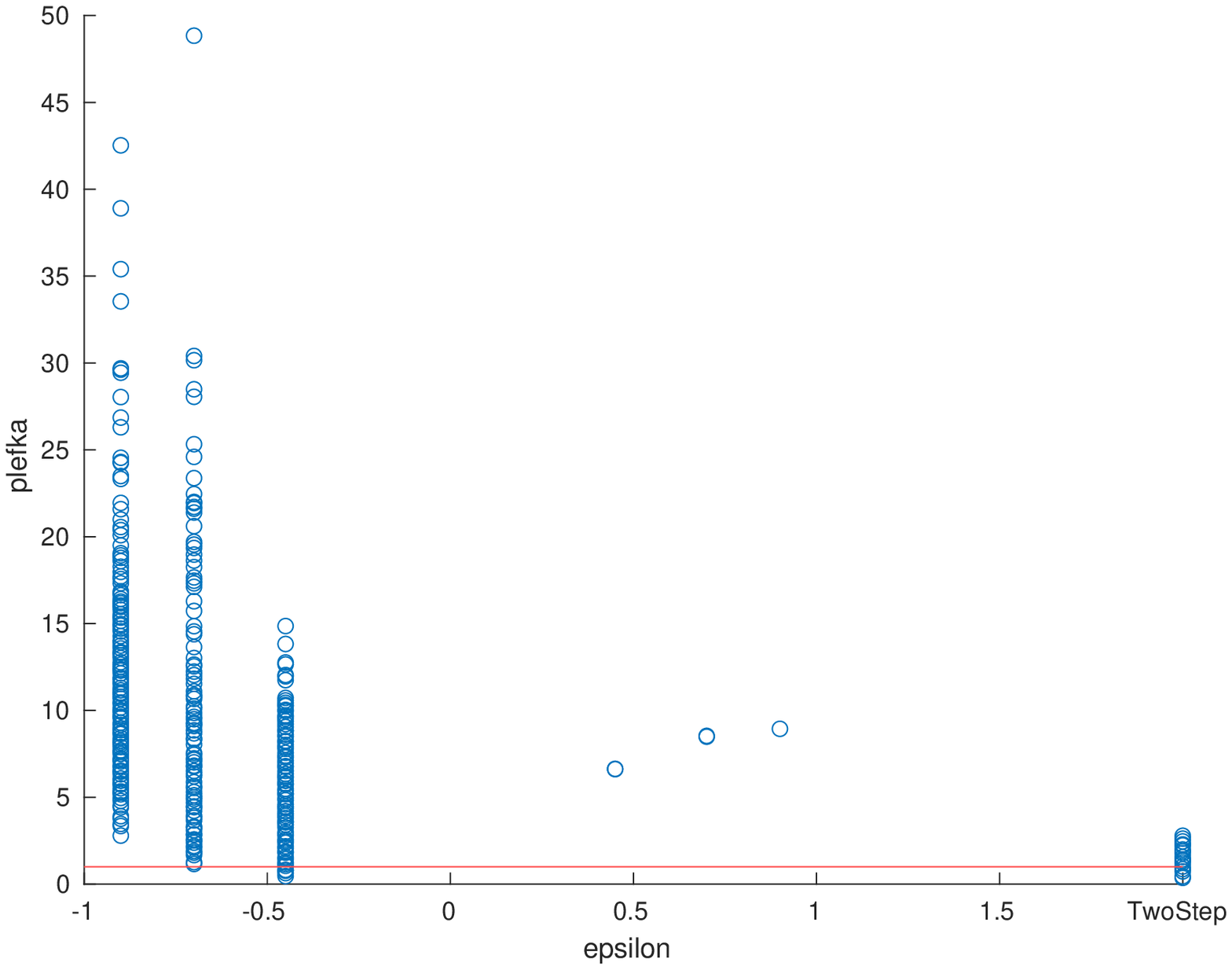}
  \caption{{\tiny} \space  1000 iterations of the $\varepsilon$-Banach and the Two-Step-Banach algorithm for $N=100, \beta=3$ and $h=0.5$. We simulated 1000 uniformly in the corners of the hypercube distributed start values for $\varepsilon = -0.9, -0.7, -0.45, 0.45, 0.7, 0.9$ and Two-Step-Banach. We plotted the Plefka values for the given parameters and used the same 1000 start values for all $\varepsilon$ and Two-Step-Banach.}
\end{figure}
\begin{figure}[H]
  \centering
  % \psfrag{2}{\tiny Two Step}
  \includegraphics[width=0.5\textwidth,  trim={0.45cm 0.4cm 0 0},clip]{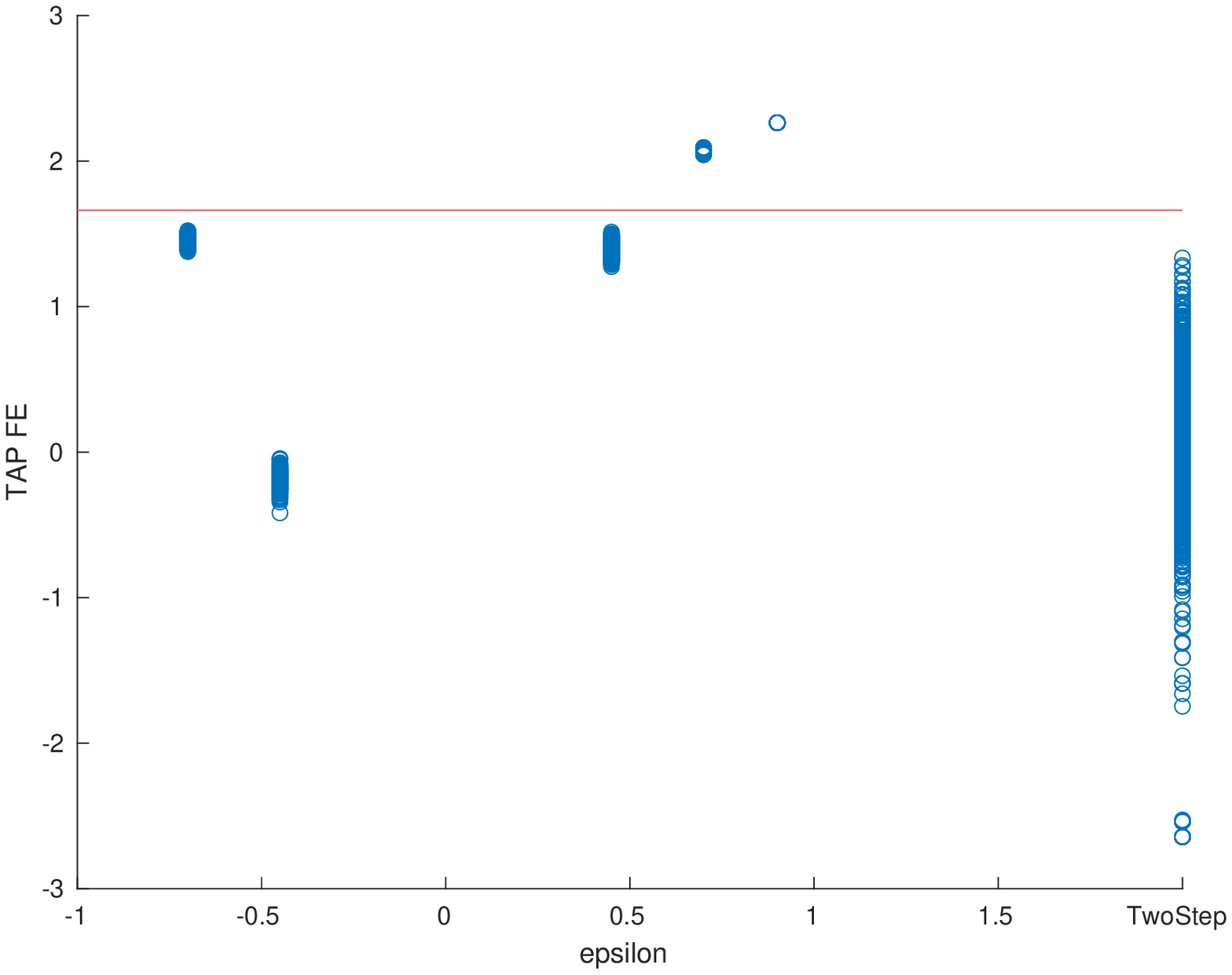}
  \caption{{\tiny} \space  1000 iterations of the $\varepsilon$-Banach and the Two-Step-Banach algorithm for $N=100, \beta=3$ and $h=0.5$. We simulated 1000 uniformly in the corners of the hypercube distributed start values for $\varepsilon = -0.9, -0.7, -0.45, 0.45, 0.7, 0.9$ and Two-Step-Banach. We plotted the TAP-free-energy values for all iterations ending inside the corners of the hypercube for the given parameters and used the same 1000 start values for all $\varepsilon$ and Two-Step-Banach.}
\end{figure}
\FloatBarrier\clearpage
\subsubsection{$N=500$,  $\beta =3 $ and $h=0.5$}\leavevmode
For all of the following figures five we used the same start values:
\begin{figure}[H]
  \centering
  % \psfrag{2}{\tiny Two Step}
  \includegraphics[width=0.5\textwidth,  trim={0.45cm 0.4cm 0 0.4cm},clip]{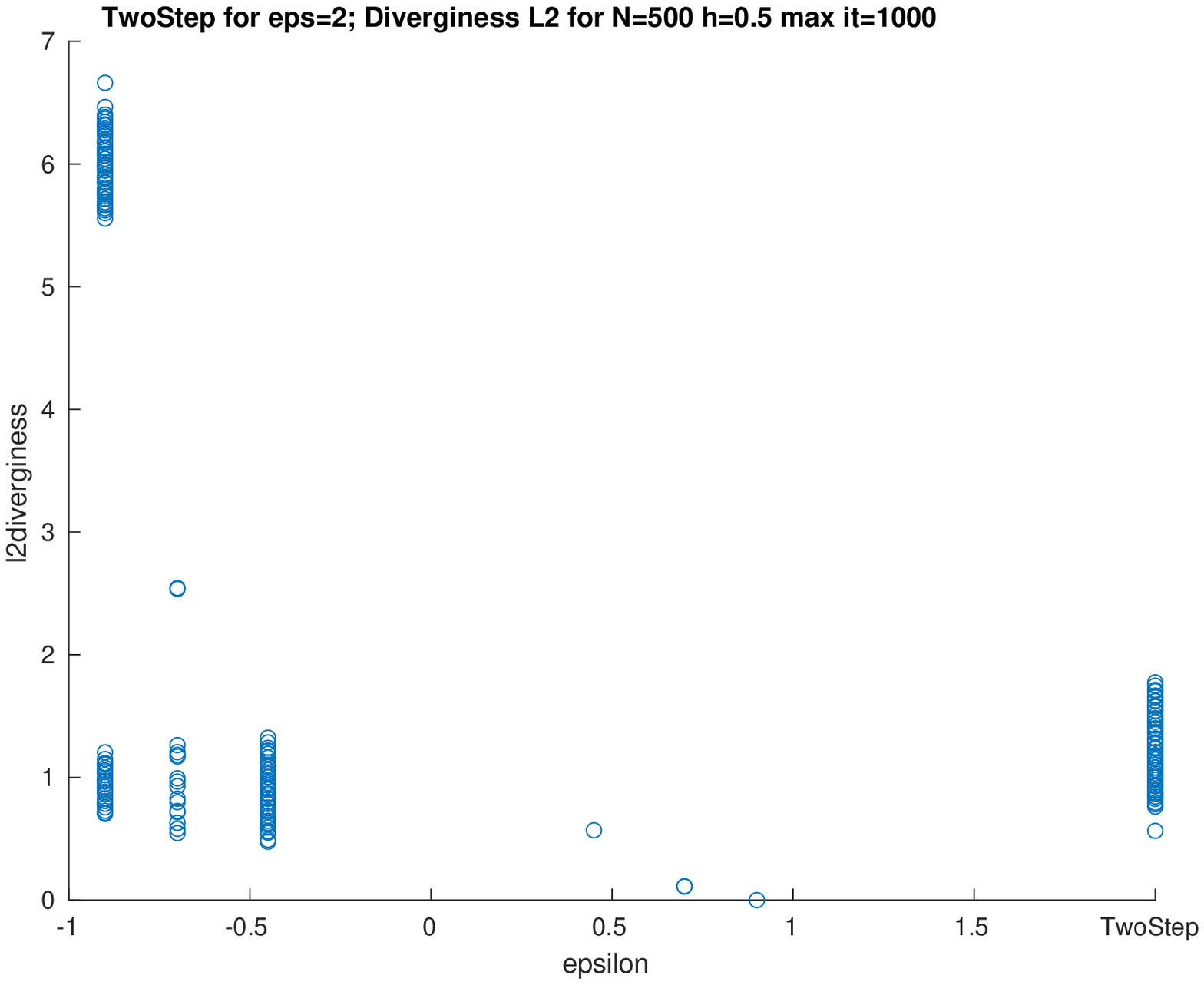}
  \caption{{\tiny} \space  1000 iterations of the $\varepsilon$-Banach and the Two-Step-Banach algorithm for $N=500, \beta=3$ and $h=0.5$. We simulated 1000 uniformly in the corners of the hypercube distributed start values for $\varepsilon = -0.9, -0.7, -0.45, 0.45, 0.7, 0.9$ and Two-Step-Banach. We plotted the {\it mean squared errors} for the given parameters and used the same 1000 start values for all $\varepsilon$ and Two-Step-Banach.}
\end{figure}
\begin{figure}[H]
  \centering
  % \psfrag{2}{\tiny Two Step}
  \includegraphics[width=0.5\textwidth,  trim={0.45cm 0.4cm 0 0.4cm},clip]{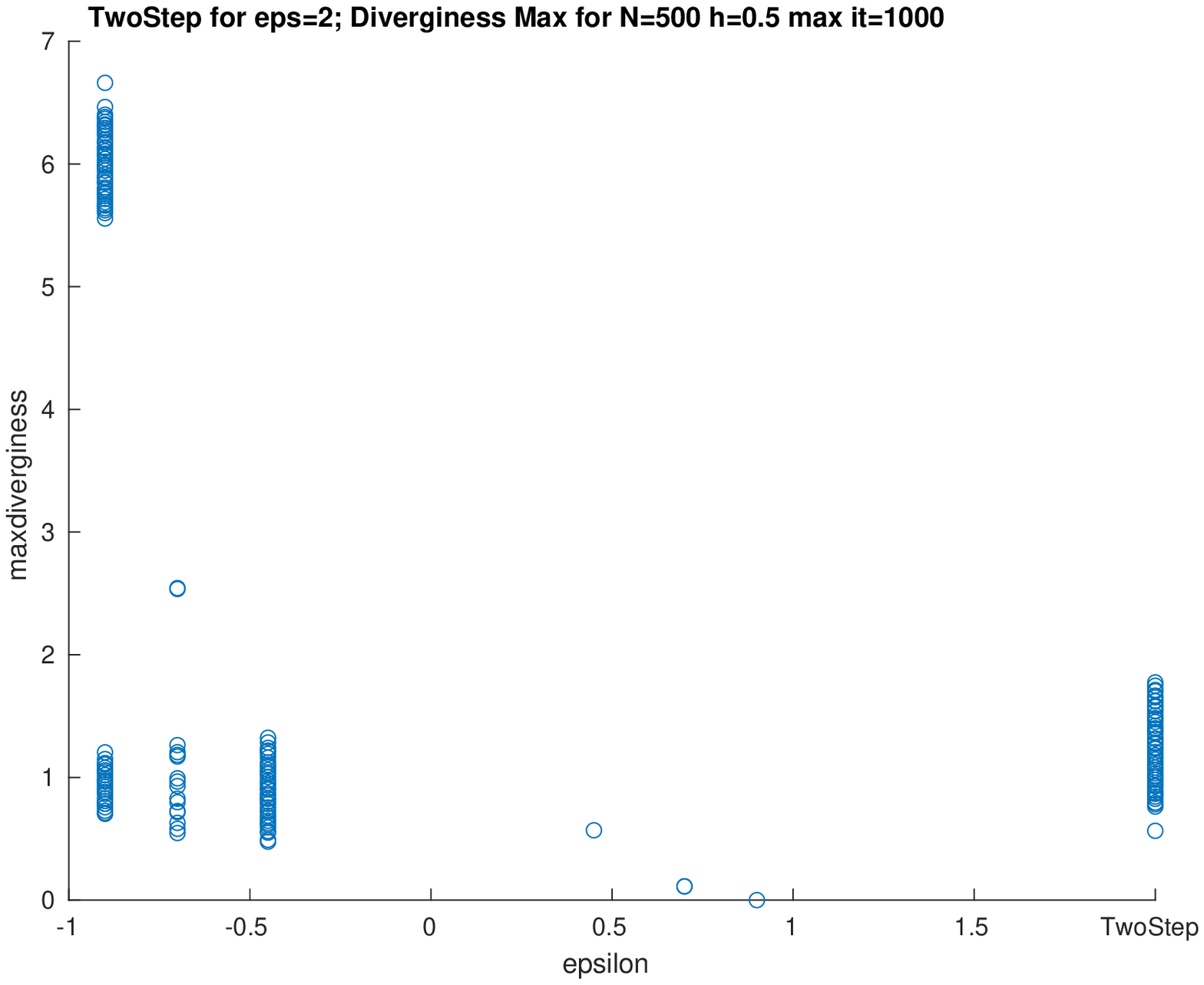}
  \caption{{\tiny} \space  1000 iterations of the $\varepsilon$-Banach and the Two-Step-Banach algorithm for $N=500, \beta=3$ and $h=0.5$. We simulated 1000 uniformly in the corners of the hypercube distributed start values for $\varepsilon = -0.9, -0.7, -0.45, 0.45, 0.7, 0.9$ and Two-Step-Banach. We plotted the {\it maximum absolute errors} for the given parameters and used the same 1000 start values for all $\varepsilon$ and Two-Step-Banach.}
\end{figure}
\begin{figure}[H]
  \centering
  % \psfrag{2}{\tiny Two Step}
  \includegraphics[width=0.5\textwidth,  trim={0.45cm 0.4cm 0 0.4cm},clip]{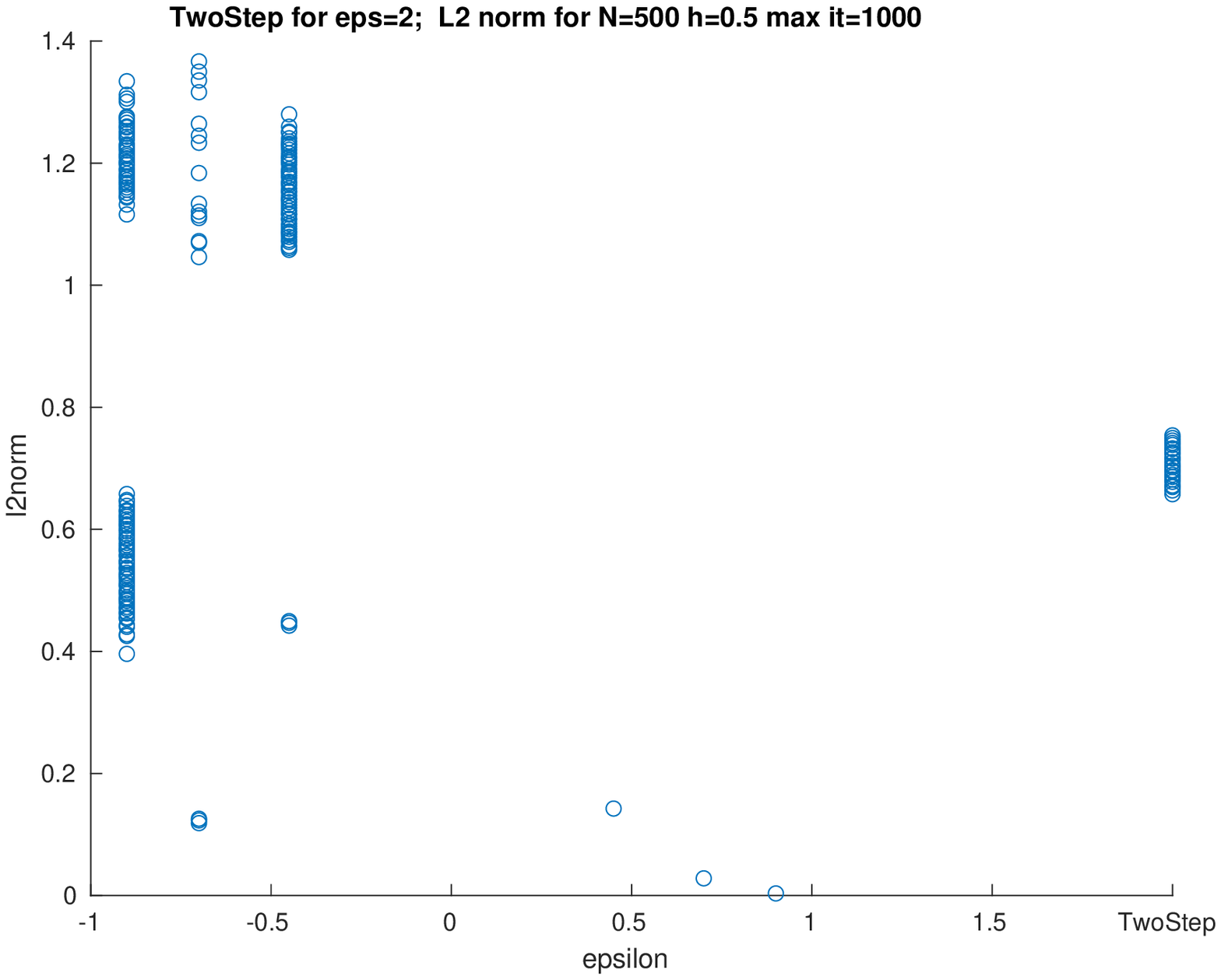}
  \caption{{\tiny} \space  1000 iterations of the $\varepsilon$-Banach and the Two-Step-Banach algorithm for $N=500, \beta=3$ and $h=0.5$. We simulated 1000 uniformly in the corners of the hypercube distributed start values for $\varepsilon = -0.9, -0.7, -0.45, 0.45, 0.7, 0.9$ and Two-Step-Banach. We plotted the {\it mean squared errors} for the given parameters and used the same 1000 start values for all $\varepsilon$ and Two-Step-Banach.}
\end{figure}
\begin{figure}[H]
  \centering
  % \psfrag{2}{\tiny Two Step}
  \includegraphics[width=0.5\textwidth,  trim={0.45cm 0.4cm 0 0},clip]{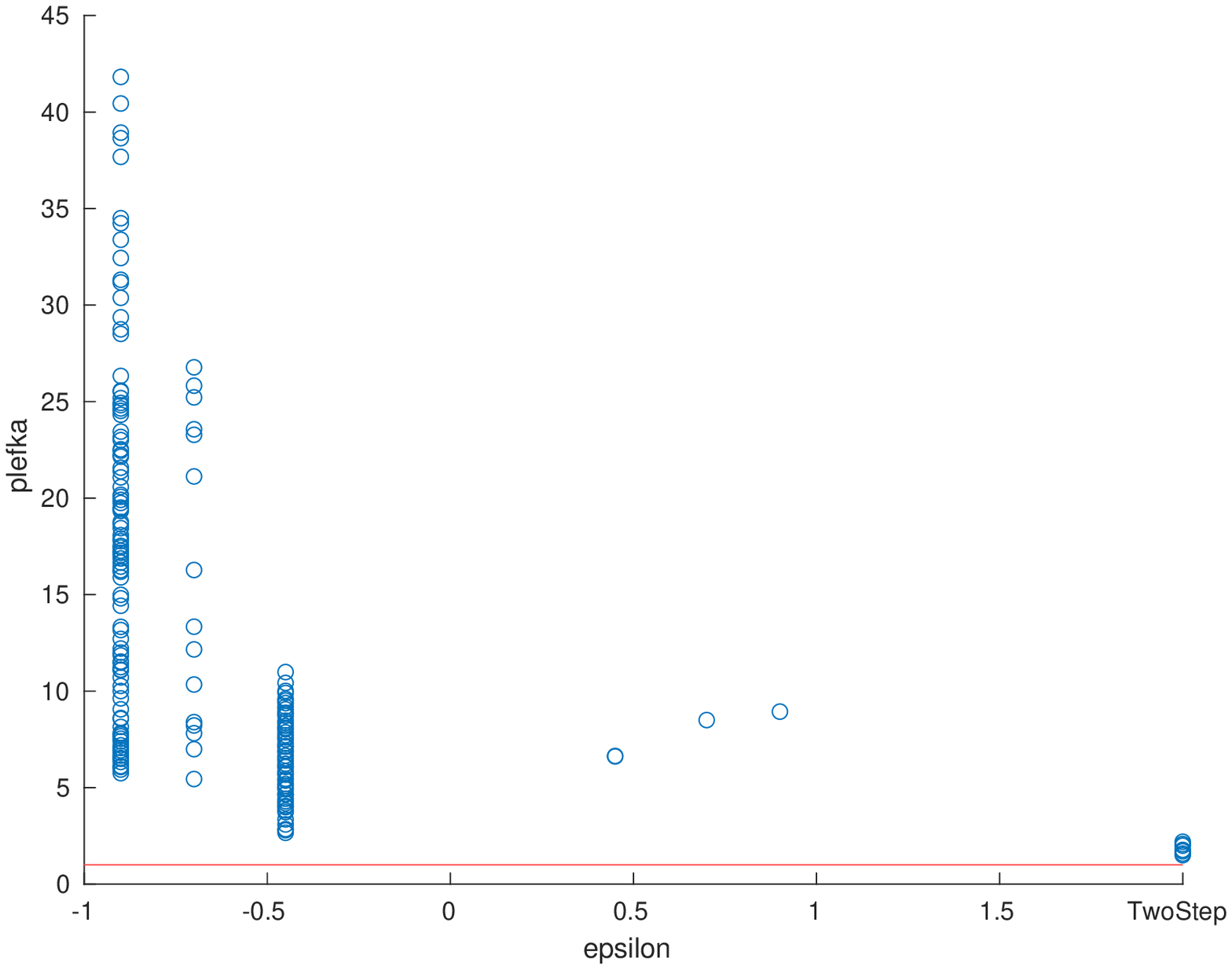}
  \caption{{\tiny} \space  1000 iterations of the $\varepsilon$-Banach and the Two-Step-Banach algorithm for $N=500, \beta=3$ and $h=0.5$. We simulated 1000 uniformly in the corners of the hypercube distributed start values for $\varepsilon = -0.9, -0.7, -0.45, 0.45, 0.7, 0.9$ and Two-Step-Banach. We plotted the Plefka values for the given parameters and used the same 1000 start values for all $\varepsilon$ and Two-Step-Banach.}
\end{figure}
\begin{figure}[H]
  \centering
  % \psfrag{2}{\tiny Two Step}
  \includegraphics[width=0.5\textwidth,  trim={0.45cm 0.4cm 0 0},clip]{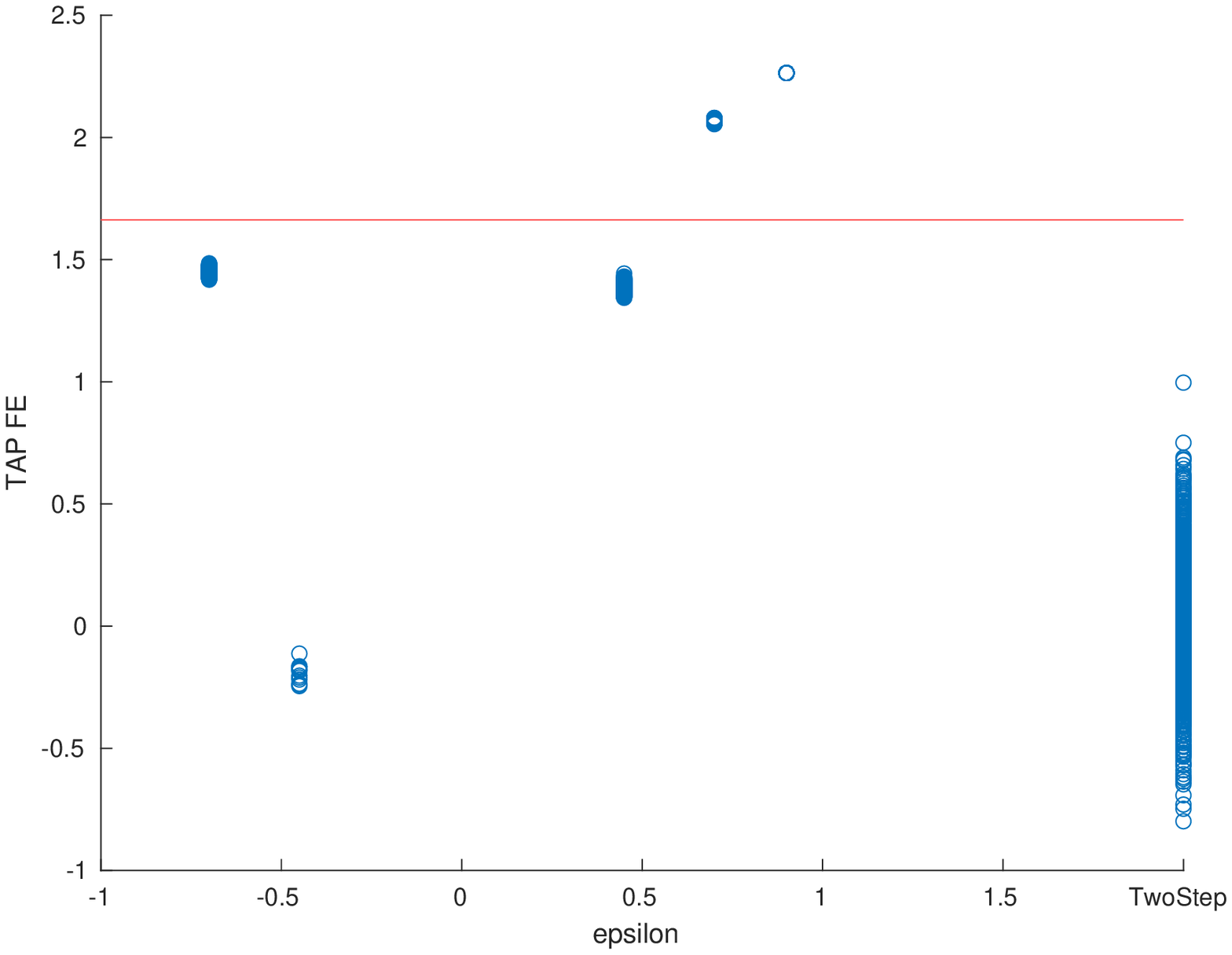}
  \caption{{\tiny} \space  1000 iterations of the $\varepsilon$-Banach and the Two-Step-Banach algorithm for $N=500, \beta=3$ and $h=0.5$. We simulated 1000 uniformly in the corners of the hypercube distributed start values for $\varepsilon = -0.9, -0.7, -0.45, 0.45, 0.7, 0.9$ and Two-Step-Banach. We plotted the TAP-free-energy values for all iterations ending inside the corners of the hypercube for the given parameters and used the same 1000 start values for all $\varepsilon$ and Two-Step-Banach.}
\end{figure}
\FloatBarrier\clearpage

\subsection{Short summary}
We see that the dependence on the start values is much stronger for negative $\varepsilon$ than for positive $\varepsilon$. For Banach Two Step everything seems to depend on the start values. It's quite interesting that for a given TAP free energy level there seems to be a range of optimal $\varepsilon$ and also a few values of $\varepsilon$ which won't give us solutions with the desired level of TAP free energy. (That seems true also for the Plefka value and the {\it mean squared error}.) It seems to be true that for negative $\varepsilon$ we have that for smaller $\varepsilon$ the dependence on the start value is much stronger. This is probably because of the fact that for negative $\varepsilon$ the iterations spend some time outside the hypercube. 
\FloatBarrier\clearpage
\section{Finding good solutions in high temperature: $\beta=3, h=0.5$}
We simulated 100.000 realizations of the disorder, if we talk about $n$ of them we talk about the first $n$ of them. We only accepted solutions satisfying the Plefka condition and being inside of the hypercube.
\FloatBarrier\clearpage
\subsection{100 realizations of disorder and 50 start values and 21 $\varepsilon$-values}\leavevmode
The following plots concern the same iterations with the same realizations, same start values.  We only accepted solutions satisfying the Plefka condition and being inside of the hypercube. 
\begin{figure}[H]
  \centering
  % \psfrag{2}{\tiny Two Step}
  \includegraphics[width=0.5\textwidth,  trim={0.45cm 0.4cm 0 0},clip]{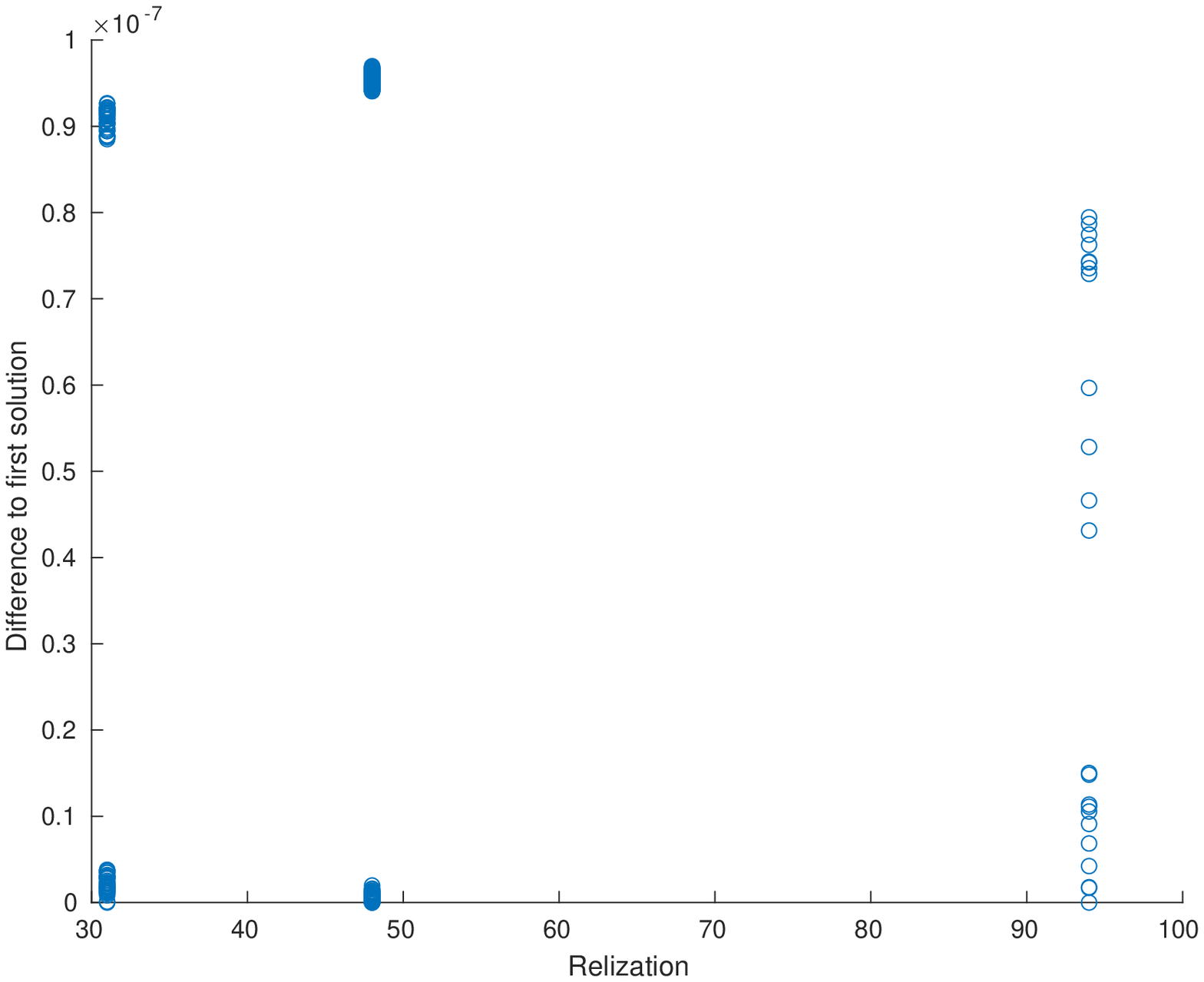}
  \caption{{\tiny} \space  For $N=25, \beta =3, h=0.5$ we simulated 100 realizations of the disorder and 50 start values used for all of them. For $\varepsilon = -0.705 - 0.001 \cdot k, k=0,\ldots, 20$ we made up to 1.000 iterations of the Banach-$\varepsilon$-algorithm where we stopped when we hit the target level of stability (measured with respect to the {\it maximum absolute error}) of $10^{-7}$.  For all realizations we plotted the {\it maximum absolute error} to the first solution of the realization we found.}
\end{figure}
\begin{figure}[H]
  \centering
  % \psfrag{2}{\tiny Two Step}
  \includegraphics[width=0.5\textwidth,  trim={0.45cm 0.4cm 0 0},clip]{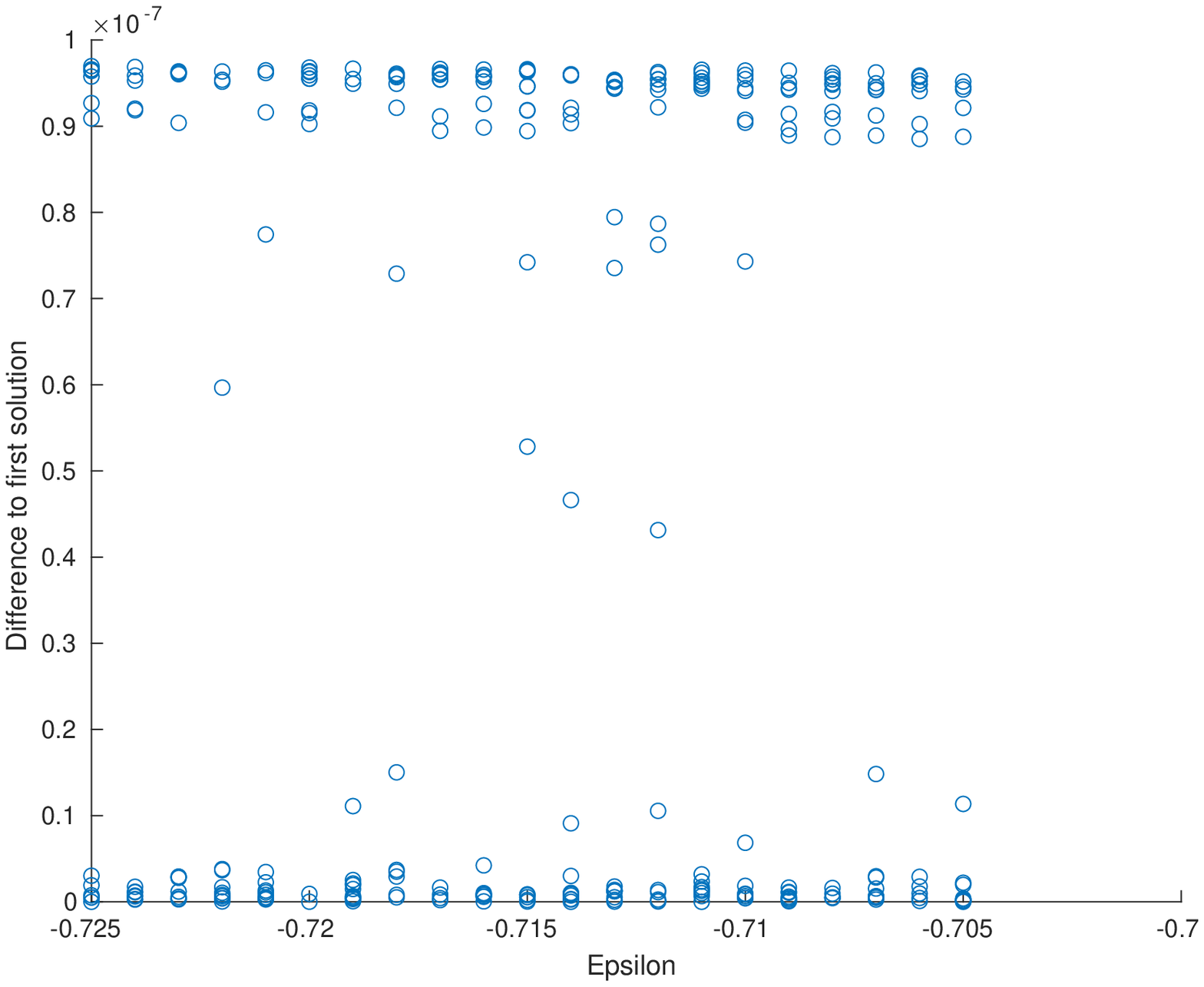}
  \caption{{\tiny} \space  For $N=25, \beta =3, h=0.5$ we simulated 100 realizations of the disorder and 50 start values used for all of them. For $\varepsilon = -0.705 - 0.001 \cdot k, k=0,\ldots, 20$ we made up to 1.000 iterations of the Banach-$\varepsilon$-algorithm where we stopped when we hit the target level of stability (measured with respect to the {\it maximum absolute error}) of $10^{-7}$.  For all $\varepsilon$ we plotted the {\it maximum absolute error} to the first solution of the realization we found.}
\end{figure}
\begin{figure}[H]
  \centering
  % \psfrag{2}{\tiny Two Step}
  \includegraphics[width=0.5\textwidth,  trim={0.45cm 0.4cm 0 0},clip]{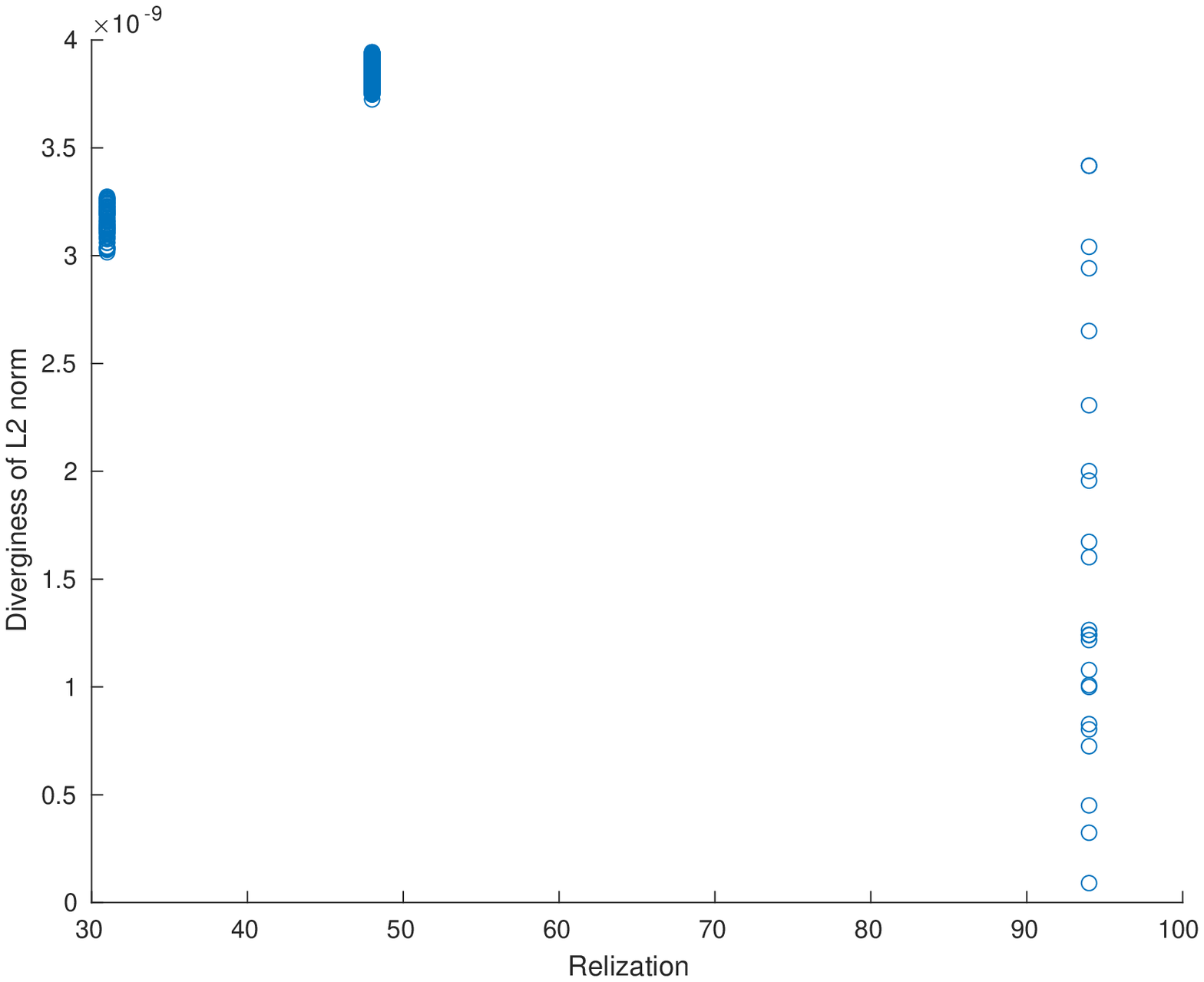}
  \caption{{\tiny} \space  For $N=25, \beta =3, h=0.5$ we simulated 100 realizations of the disorder and 50 start values used for all of them. For $\varepsilon = -0.705 - 0.001 \cdot k, k=0,\ldots, 20$ we made up to 1.000 iterations of the Banach-$\varepsilon$-algorithm where we stopped when we hit the target level of stability (measured with respect to the {\it maximum absolute error}) of $10^{-7}$.  For all realizations we plotted the difference of the {\it mean squared errors} of the last two iterates of the solutions found in the respective realization.}
\end{figure}
\begin{figure}[H]
  \centering
  % \psfrag{2}{\tiny Two Step}
  \includegraphics[width=0.5\textwidth,  trim={0.45cm 0.4cm 0 0},clip]{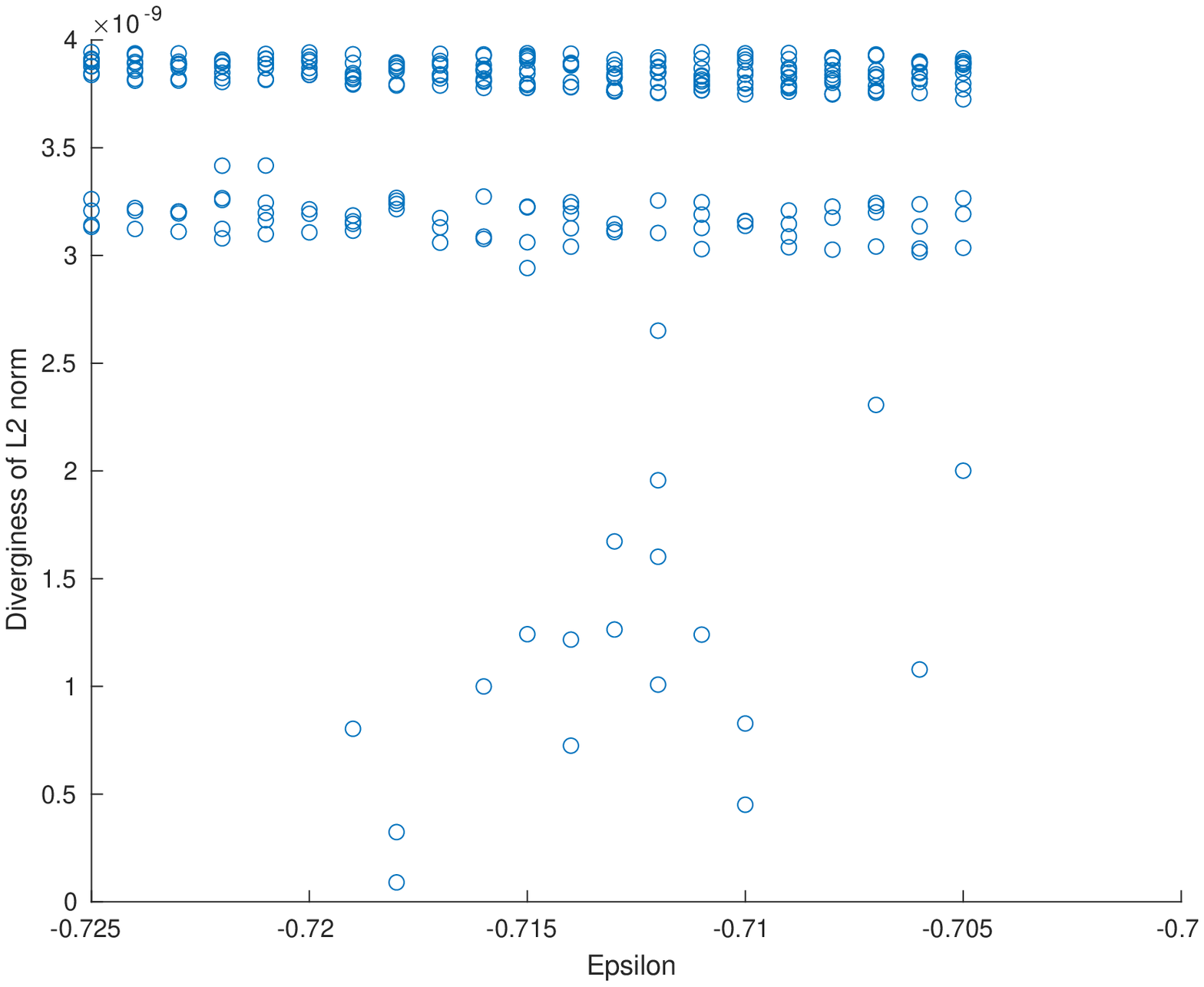}
  \caption{{\tiny} \space  For $N=25, \beta =3, h=0.5$ we simulated 100 realizations of the disorder and 50 start values used for all of them. For $\varepsilon = -0.705 - 0.001 \cdot k, k=0,\ldots, 20$ we made up to 1.000 iterations of the Banach-$\varepsilon$-algorithm where we stopped when we hit the target level of stability (measured with respect to the {\it maximum absolute error}) of $10^{-7}$.  For all $\varepsilon$ we plotted the difference of the {\it mean squared errors} of the last two iterates of the solutions found with the respective $\varepsilon$.}
\end{figure}
\begin{figure}[H]
  \centering
  % \psfrag{2}{\tiny Two Step}
  \includegraphics[width=0.5\textwidth,  trim={0.45cm 0.4cm 0 0},clip]{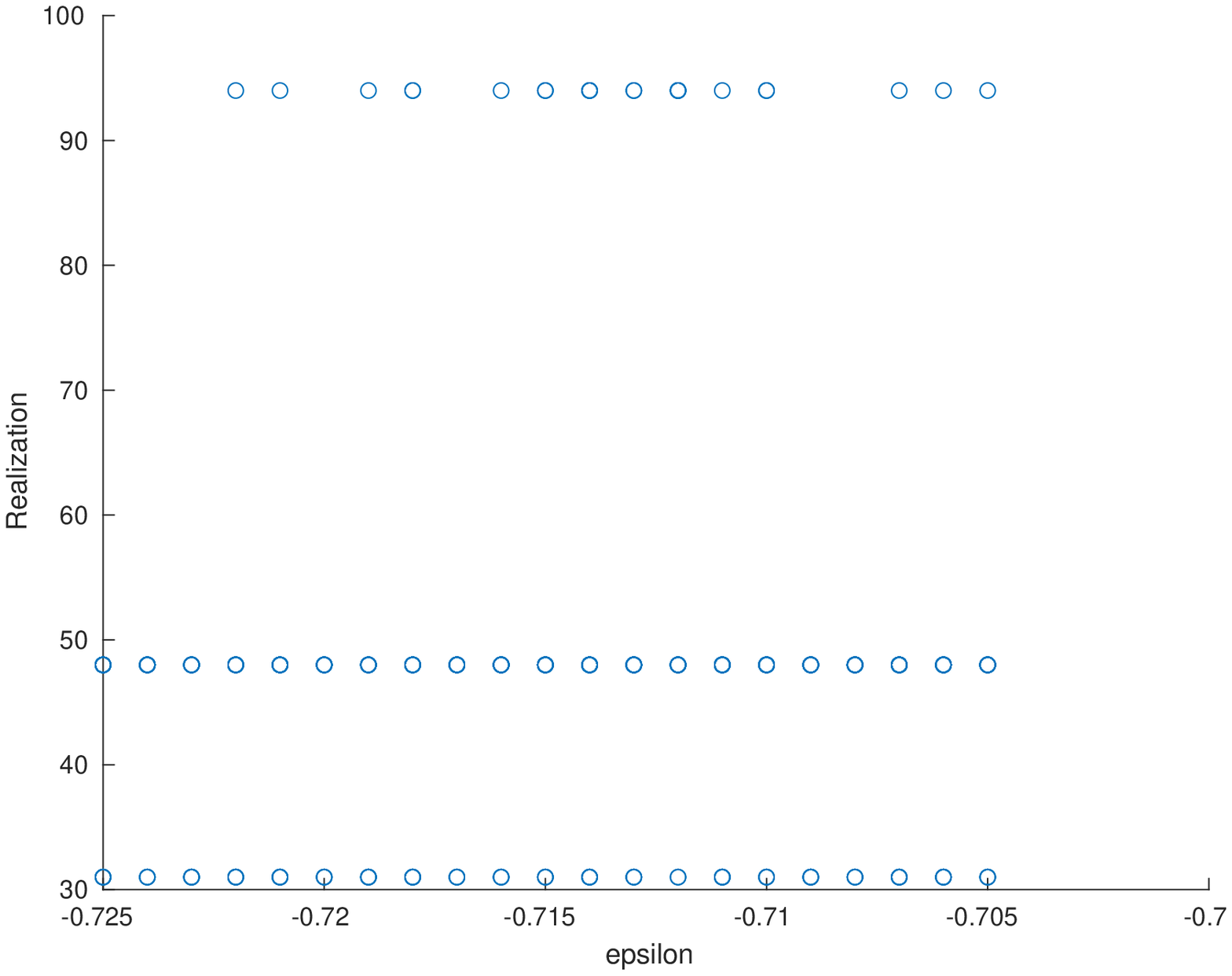}
  \caption{{\tiny} \space  For $N=25, \beta =3, h=0.5$ we simulated 100 realizations of the disorder and 50 start values used for all of them. For $\varepsilon = -0.705 - 0.001 \cdot k, k=0,\ldots, 20$ we made up to 1.000 iterations of the Banach-$\varepsilon$-algorithm where we stopped when we hit the target level of stability (measured with respect to the {\it maximum absolute error}) of $10^{-7}$.  For all solutions we found we plotted a point with $x$ coordinate specifying the used $\varepsilon$ and $y$ coordinate the concerning realization}
\end{figure}
\begin{figure}[H]
  \centering
  % \psfrag{2}{\tiny Two Step}
  \includegraphics[width=0.5\textwidth,  trim={0.45cm 0.4cm 0 0},clip]{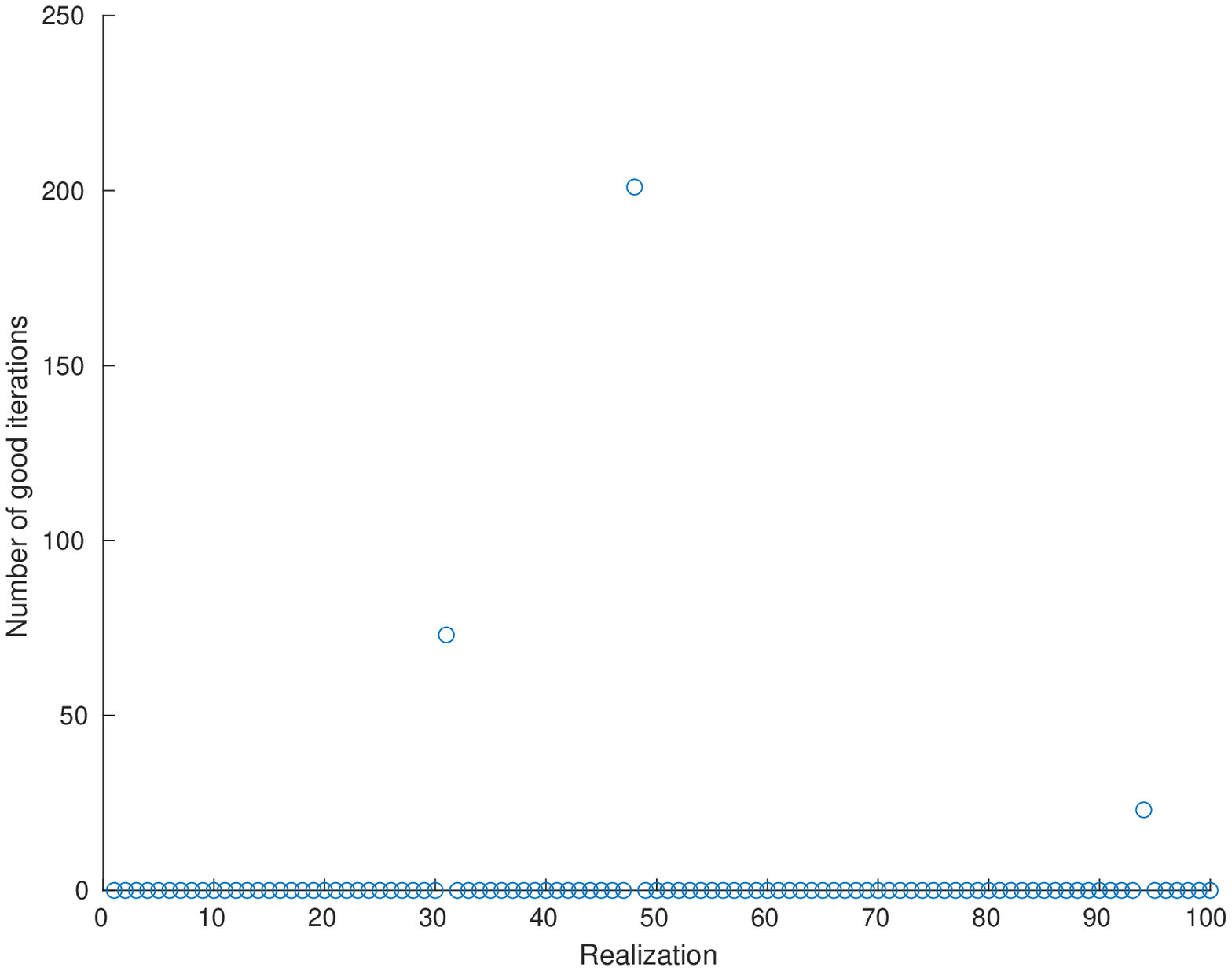}
  \caption{{\tiny} \space  For $N=25, \beta =3, h=0.5$ we simulated 100 realizations of the disorder and 50 start values used for all of them. For $\varepsilon = -0.705 - 0.001 \cdot k, k=0,\ldots, 20$ we made up to 1.000 iterations of the Banach-$\varepsilon$-algorithm where we stopped when we hit the target level of stability (measured with respect to the {\it maximum absolute error}) of $10^{-7}$.  For all realizations we plotted the number of good solutions. }
\end{figure}
\begin{figure}[H]
  \centering
  % \psfrag{2}{\tiny Two Step}
  \includegraphics[width=0.5\textwidth,  trim={0.45cm 0.4cm 0 0},clip]{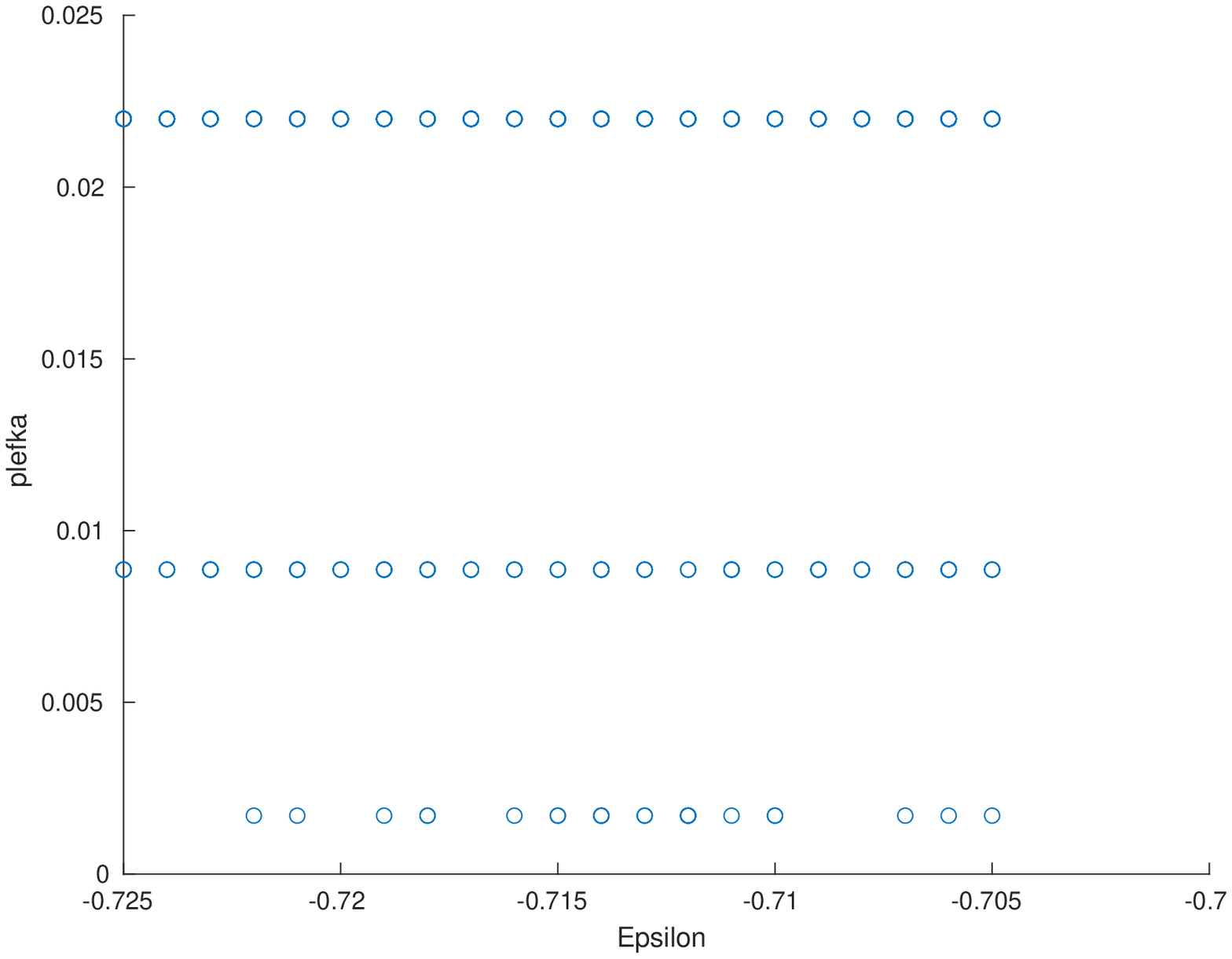}
  \caption{{\tiny} \space  For $N=25, \beta =3, h=0.5$ we simulated 100 realizations of the disorder and 50 start values used for all of them. For $\varepsilon = -0.705 - 0.001 \cdot k, k=0,\ldots, 20$ we made up to 1.000 iterations of the Banach-$\varepsilon$-algorithm where we stopped when we hit the target level of stability (measured with respect to the {\it maximum absolute error}) of $10^{-7}$.  For all $\varepsilon$ we plotted the Plefka values found with this $\varepsilon$. }
\end{figure}
\begin{figure}[H]
  \centering
  % \psfrag{2}{\tiny Two Step}
  \includegraphics[width=0.5\textwidth,  trim={0.45cm 0.4cm 0 0},clip]{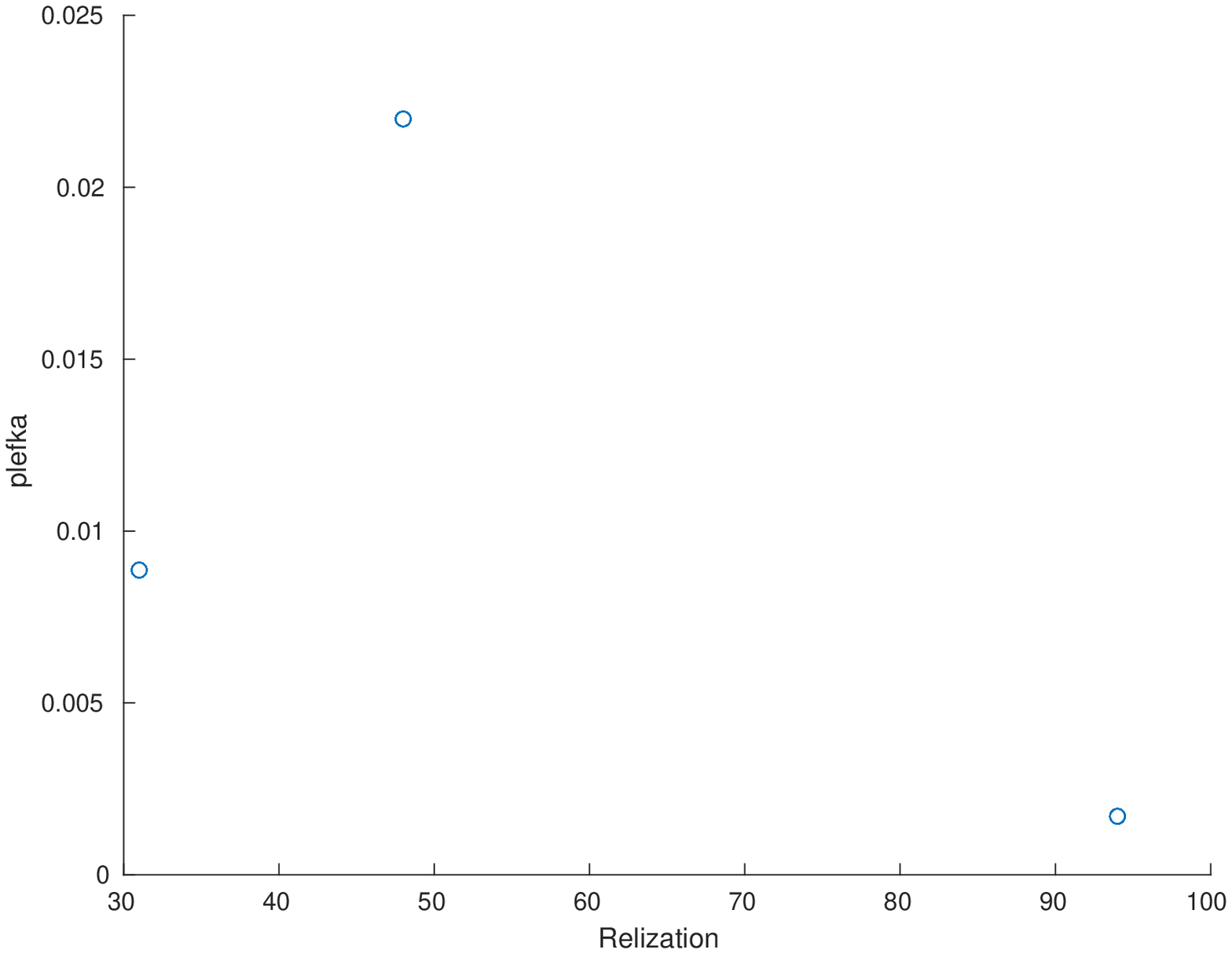}
  \caption{{\tiny} \space  For $N=25, \beta =3, h=0.5$ we simulated 100 realizations of the disorder and 50 start values used for all of them. For $\varepsilon = -0.705 - 0.001 \cdot k, k=0,\ldots, 20$ we made up to 1.000 iterations of the Banach-$\varepsilon$-algorithm where we stopped when we hit the target level of stability (measured with respect to the {\it maximum absolute error}) of $10^{-7}$.  For all realizations we plotted the Plefka values found in the concerning realization. }
\end{figure}
\begin{figure}[H]
  \centering
  % \psfrag{2}{\tiny Two Step}
  \includegraphics[width=0.5\textwidth,  trim={0.45cm 0.4cm 0 0},clip]{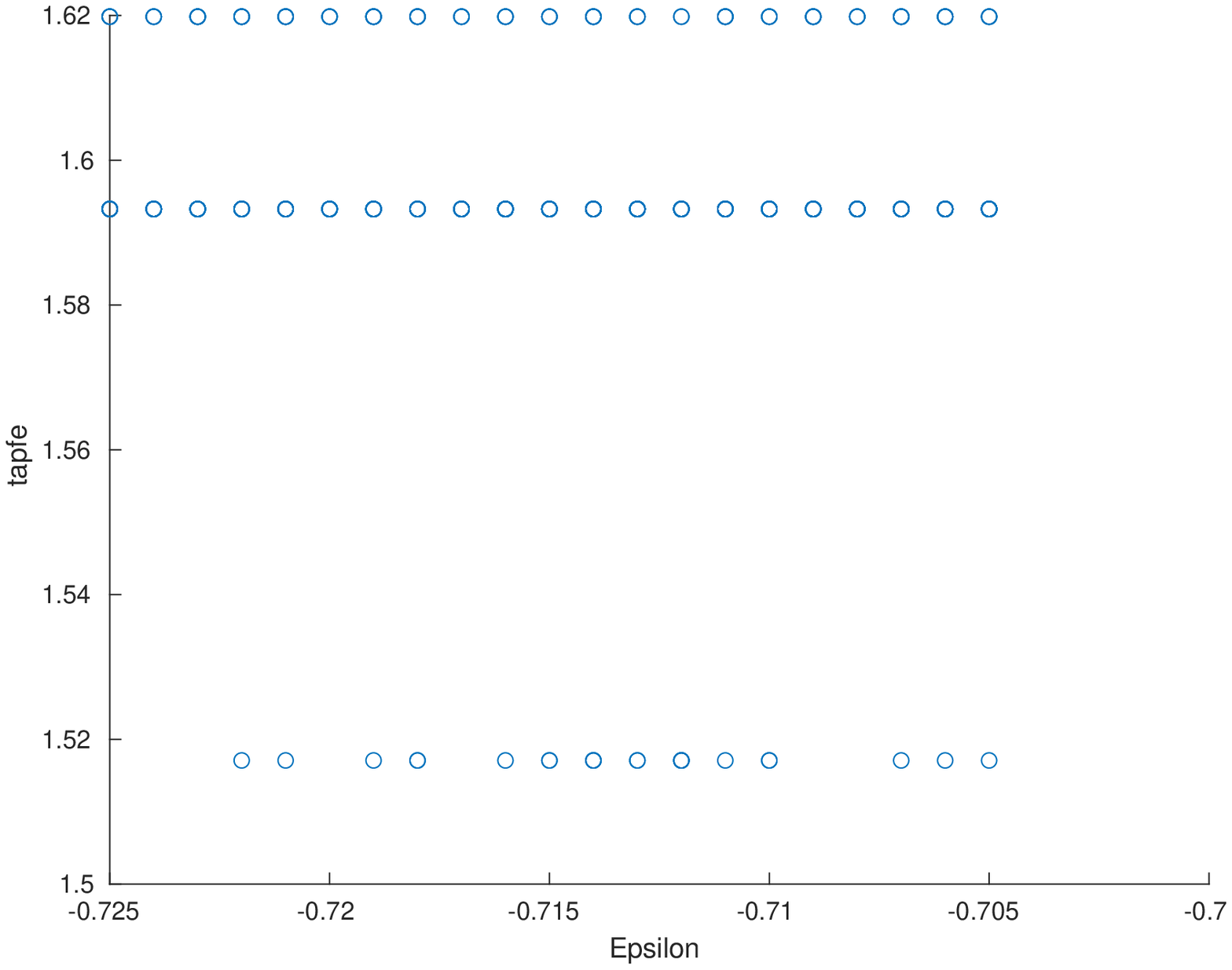}
  \caption{{\tiny} \space  For $N=25, \beta =3, h=0.5$ we simulated 100 realizations of the disorder and 50 start values used for all of them. For $\varepsilon = -0.705 - 0.001 \cdot k, k=0,\ldots, 20$ we made up to 1.000 iterations of the Banach-$\varepsilon$-algorithm where we stopped when we hit the target level of stability (measured with respect to the {\it maximum absolute error}) of $10^{-7}$.  For all $\varepsilon$ we plotted the TAP free energy values found with this $\varepsilon$. }
\end{figure}
\begin{figure}[H]
  \centering
  % \psfrag{2}{\tiny Two Step}
  \includegraphics[width=0.5\textwidth,  trim={0.45cm 0.4cm 0 0},clip]{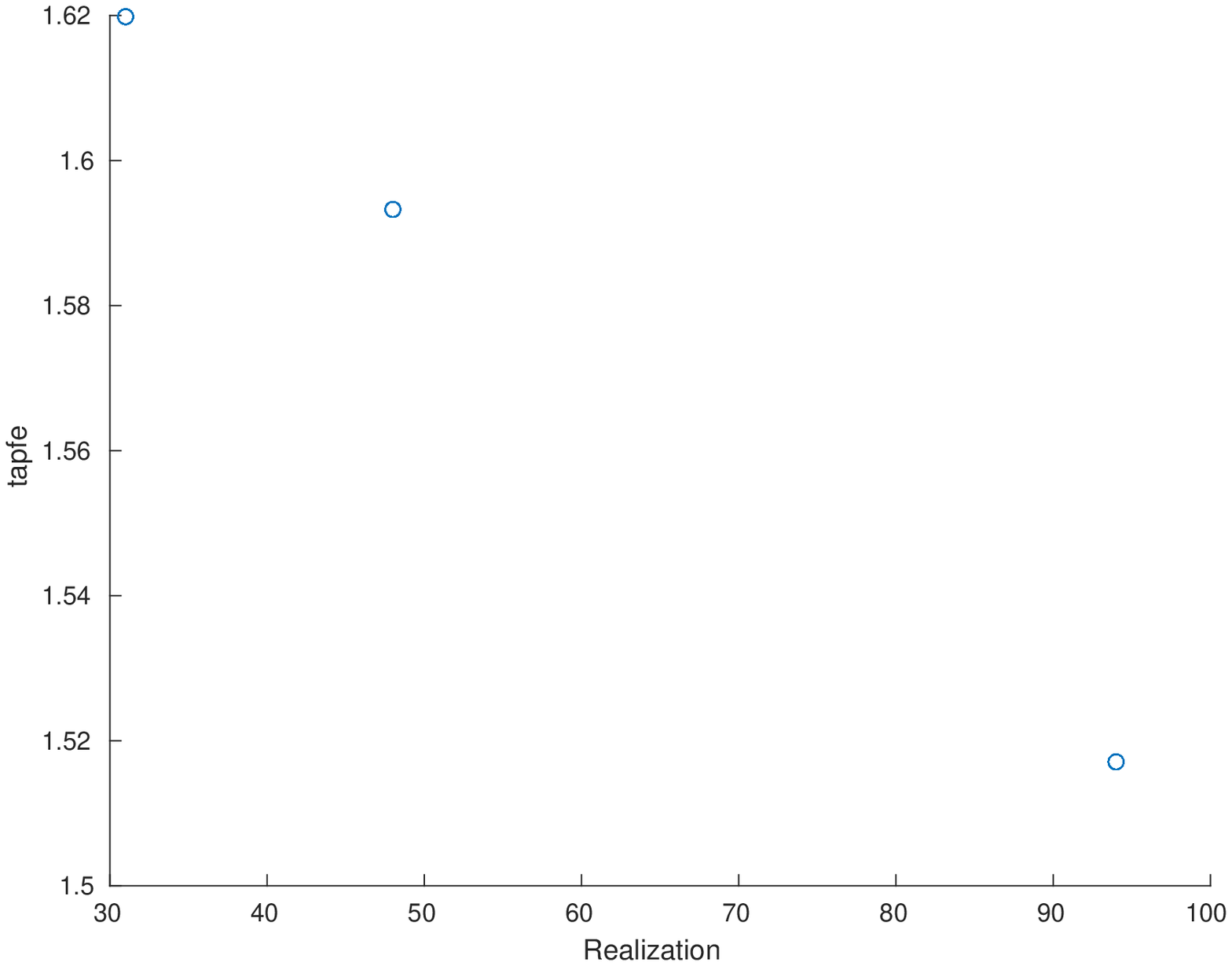}
  \caption{{\tiny} \space  For $N=25, \beta =3, h=0.5$ we simulated 100 realizations of the disorder and 50 start values used for all of them. For $\varepsilon = -0.705 - 0.001 \cdot k, k=0,\ldots, 20$ we made up to 1.000 iterations of the Banach-$\varepsilon$-algorithm where we stopped when we hit the target level of stability (measured with respect to the {\it maximum absolute error}) of $10^{-7}$.  For all realizations we plotted the TAP free energy values found for the concerning realization.}
\end{figure}
\begin{figure}[H]
  \centering
  % \psfrag{2}{\tiny Two Step}
  \includegraphics[width=0.5\textwidth,  trim={0cm 0.4cm 0 0},clip]{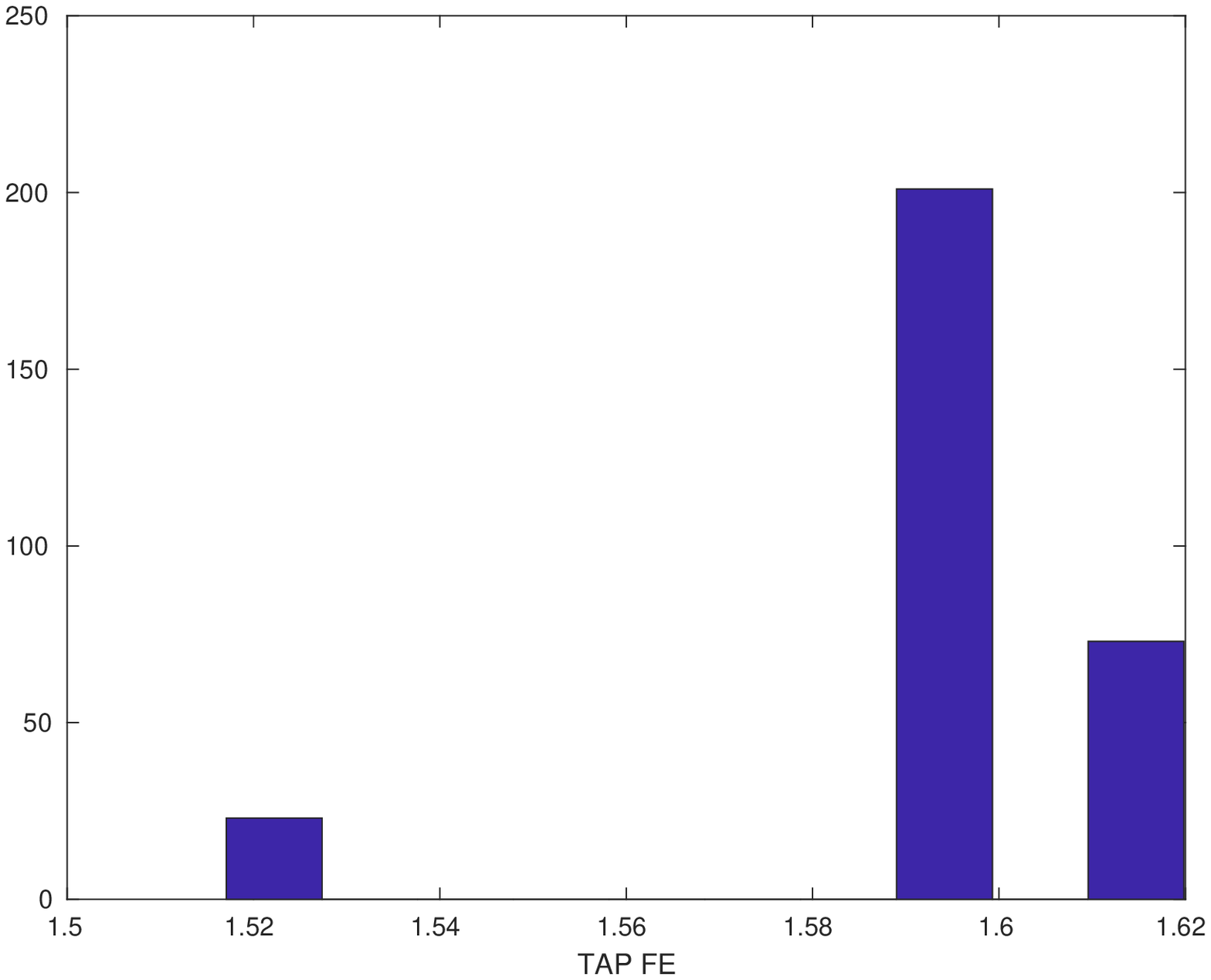}
  \caption{{\tiny} \space  For $N=25, \beta =3, h=0.5$ we simulated 100 realizations of the disorder and 50 start values used for all of them. For $\varepsilon = -0.705 - 0.001 \cdot k, k=0,\ldots, 20$ we made up to 1.000 iterations of the Banach-$\varepsilon$-algorithm where we stopped when we hit the target level of stability (measured with respect to the {\it maximum absolute error}) of $10^{-7}$.  Just a histogram of the TAP free energy values found.}
\end{figure}
\begin{figure}[H]
  \centering
  % \psfrag{2}{\tiny Two Step}
  \includegraphics[width=0.5\textwidth,  trim={0.45cm 0.4cm 0 0},clip]{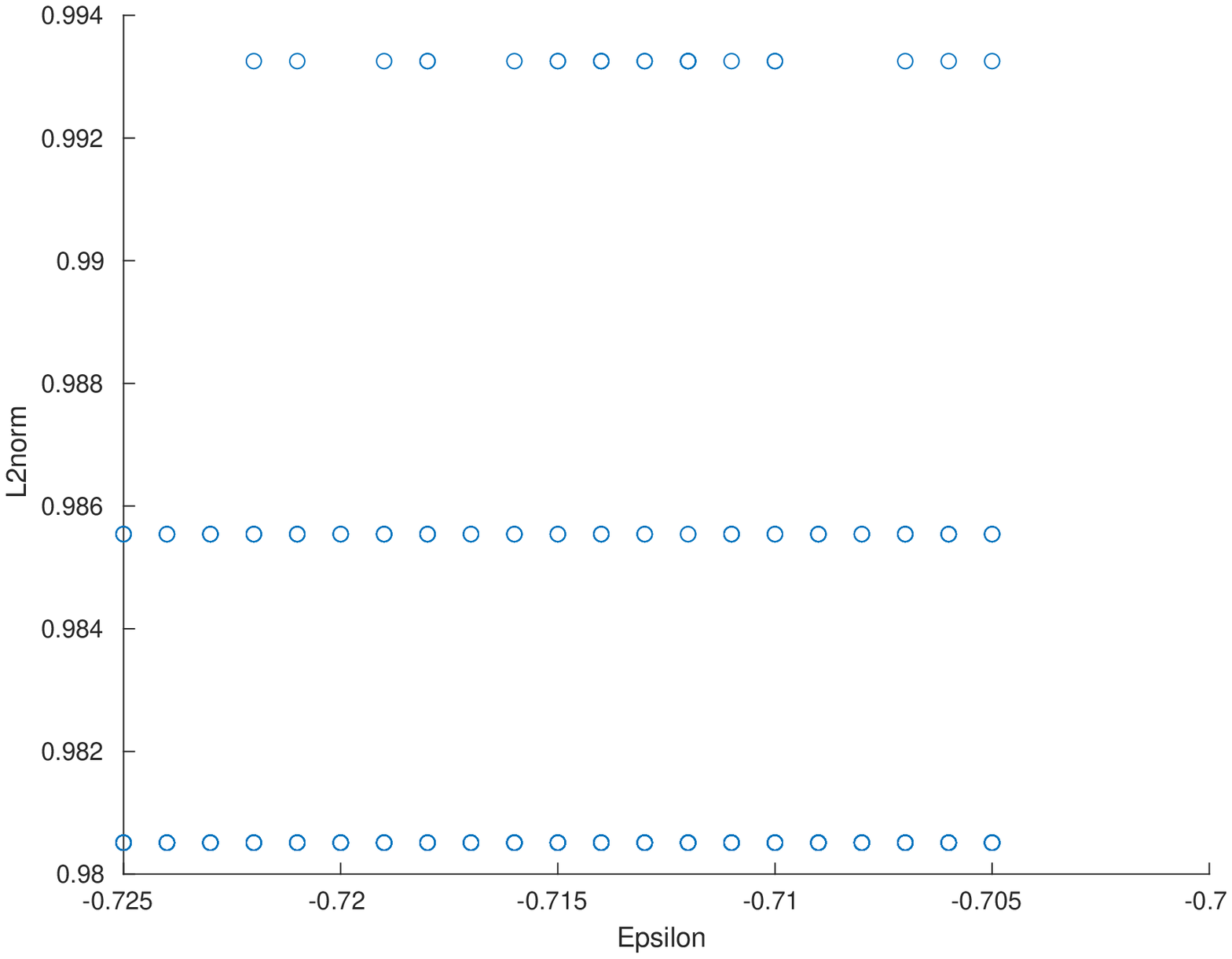}
  \caption{{\tiny} \space  For $N=25, \beta =3, h=0.5$ we simulated 100 realizations of the disorder and 50 start values used for all of them. For $\varepsilon = -0.705 - 0.001 \cdot k, k=0,\ldots, 20$ we made up to 1.000 iterations of the Banach-$\varepsilon$-algorithm where we stopped when we hit the target level of stability (measured with respect to the {\it maximum absolute error}) of $10^{-7}$.  For all $\varepsilon$ we plotted the {\it mean squared errors} found with this $\varepsilon$. }
\end{figure}
\begin{figure}[H]
  \centering
  % \psfrag{2}{\tiny Two Step}
  \includegraphics[width=0.5\textwidth,  trim={0.45cm 0.4cm 0 0},clip]{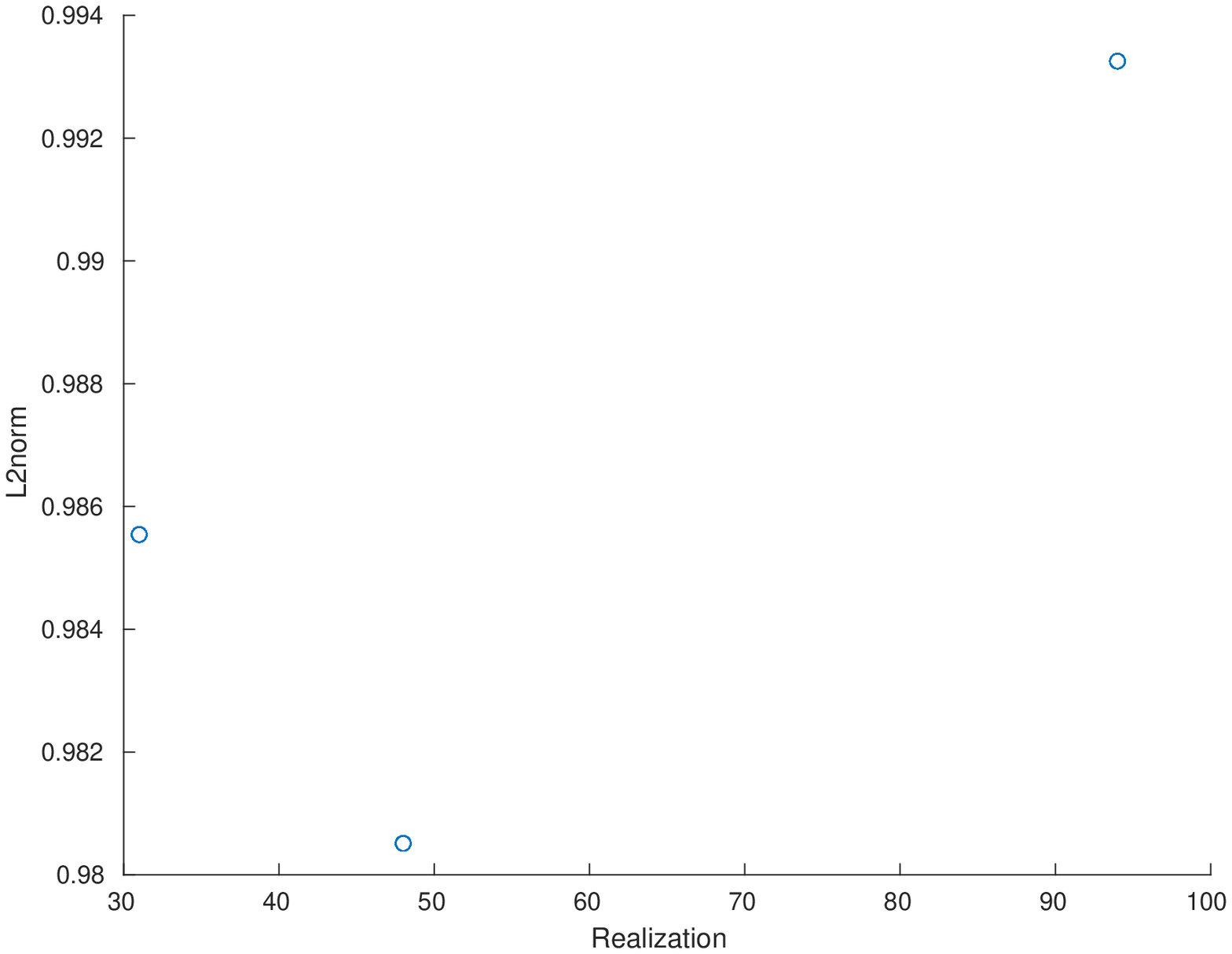}
  \caption{{\tiny} \space  For $N=25, \beta =3, h=0.5$ we simulated 100 realizations of the disorder and 50 start values used for all of them. For $\varepsilon = -0.705 - 0.001 \cdot k, k=0,\ldots, 20$ we made up to 1.000 iterations of the Banach-$\varepsilon$-algorithm where we stopped when we hit the target level of stability (measured with respect to the {\it maximum absolute error}) of $10^{-7}$.  For all realizations we plotted the {\it mean squared errors} found for the concerning realization.}
\end{figure}
\FloatBarrier\clearpage
\subsection{1000 realizations of disorder and 50 start values and 21 $\varepsilon$-values}\leavevmode
The following plots concern the same iterations with the same realizations, same start values.  We only accepted solutions satisfying the Plefka condition and being inside of the hypercube. 
\begin{figure}[H]
  \centering
  % \psfrag{2}{\tiny Two Step}
  \includegraphics[width=0.5\textwidth,  trim={0.45cm 0.4cm 0 0},clip]{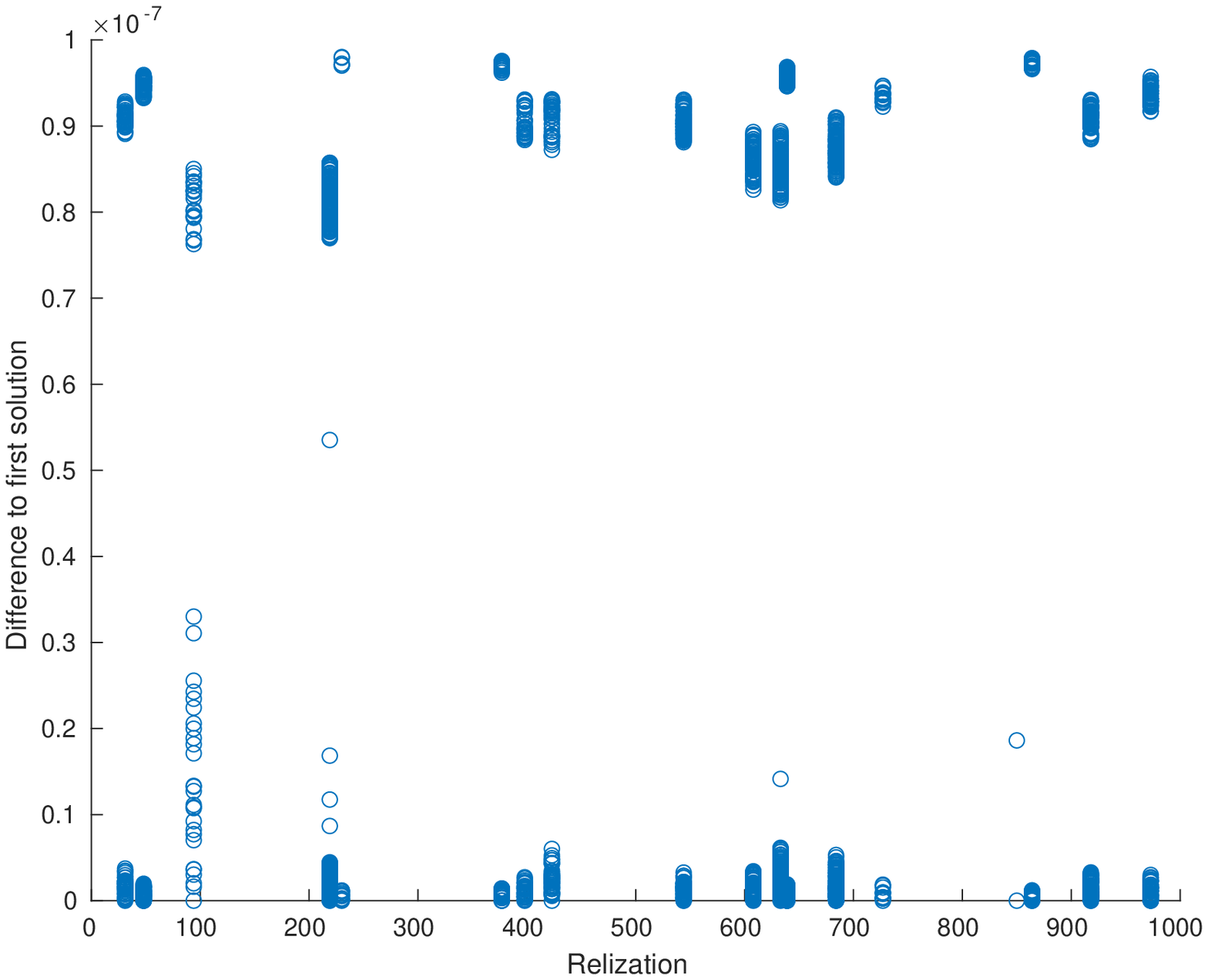}
  \caption{{\tiny} \space  For $N=25, \beta =3, h=0.5$ we simulated 1000 realizations of the disorder and 50 start values used for all of them. For $\varepsilon = -0.705 - 0.001 \cdot k, k=0,\ldots, 20$ we made up to 1.000 iterations of the Banach-$\varepsilon$-algorithm where we stopped when we hit the target level of stability (measured with respect to the {\it maximum absolute error}) of $10^{-7}$.  For all realizations we plotted the{\it maximum absolute error} to the first solution of the realization we found.}
\end{figure}
\begin{figure}[H]
  \centering
  % \psfrag{2}{\tiny Two Step}
  \includegraphics[width=0.5\textwidth,  trim={0.45cm 0.4cm 0 0},clip]{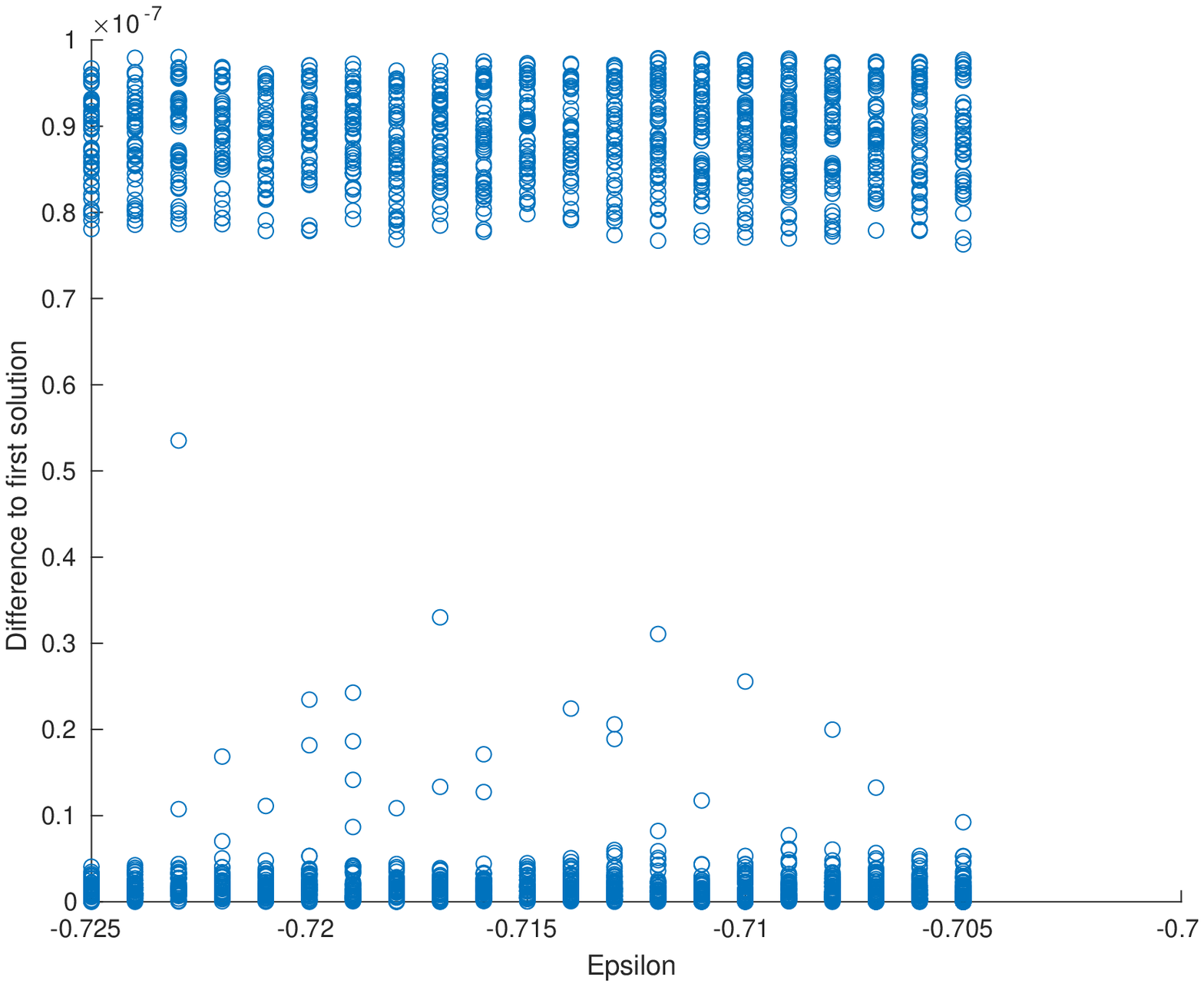}
  \caption{{\tiny} \space  For $N=25, \beta =3, h=0.5$ we simulated 1000 realizations of the disorder and 50 start values used for all of them. For $\varepsilon = -0.705 - 0.001 \cdot k, k=0,\ldots, 20$ we made up to 1.000 iterations of the Banach-$\varepsilon$-algorithm where we stopped when we hit the target level of stability (measured with respect to the {\it maximum absolute error}) of $10^{-7}$.  For all $\varepsilon$ we plotted the {\it maximum absolute error} to the first solution of the realization we found.}
\end{figure}
\begin{figure}[H]
  \centering
  % \psfrag{2}{\tiny Two Step}
  \includegraphics[width=0.5\textwidth,  trim={0.45cm 0.4cm 0 0},clip]{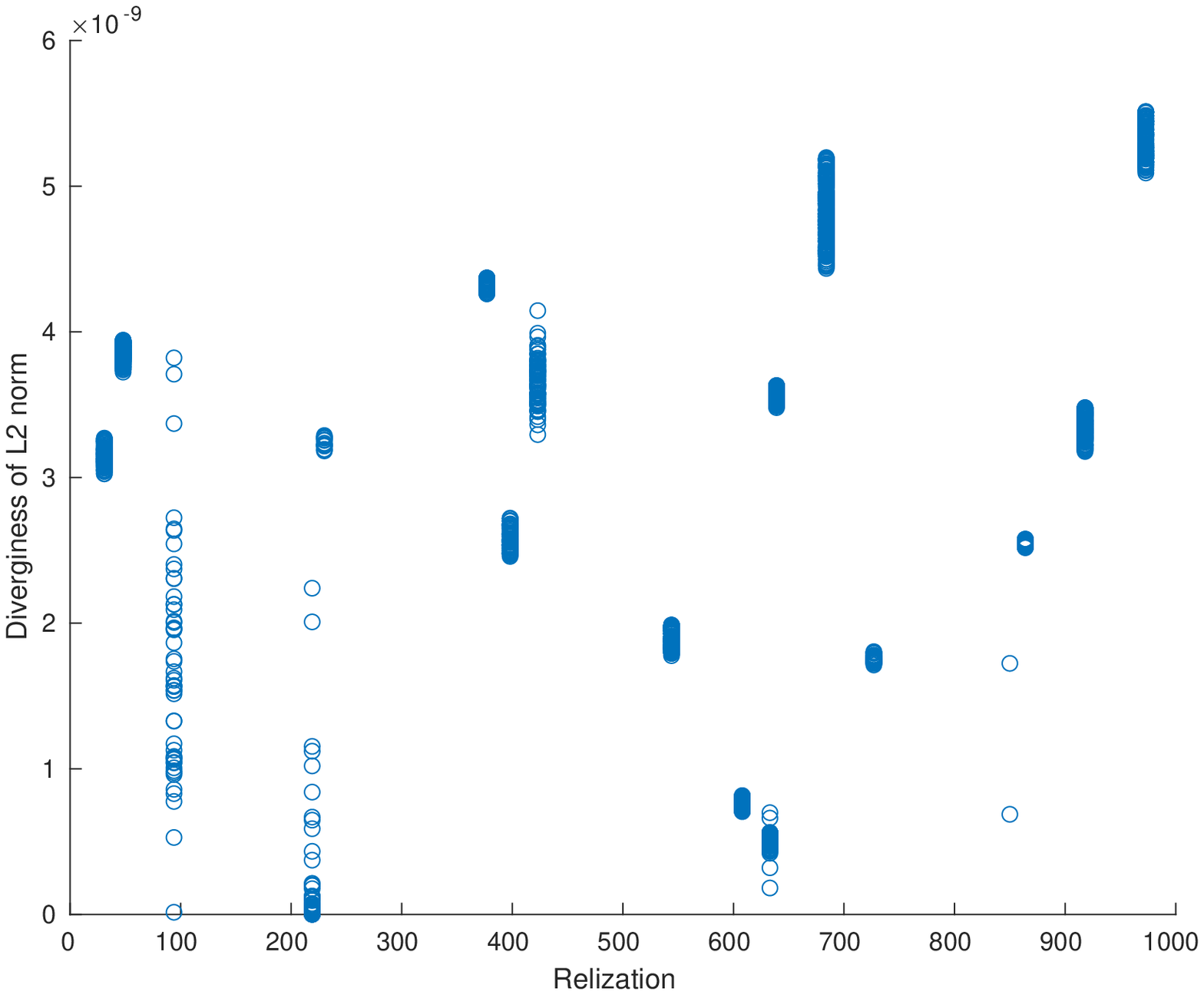}
  \caption{{\tiny} \space  For $N=25, \beta =3, h=0.5$ we simulated 1000 realizations of the disorder and 50 start values used for all of them. For $\varepsilon = -0.705 - 0.001 \cdot k, k=0,\ldots, 20$ we made up to 1.000 iterations of the Banach-$\varepsilon$-algorithm where we stopped when we hit the target level of stability (measured with respect to the {\it maximum absolute error}) of $10^{-7}$.  For all realizations we plotted the difference of the {\it mean squared errors} of the last two iterates of the solutions found in the respective realization.}
\end{figure}
\begin{figure}[H]
  \centering
  % \psfrag{2}{\tiny Two Step}
  \includegraphics[width=0.5\textwidth,  trim={0.45cm 0.4cm 0 0},clip]{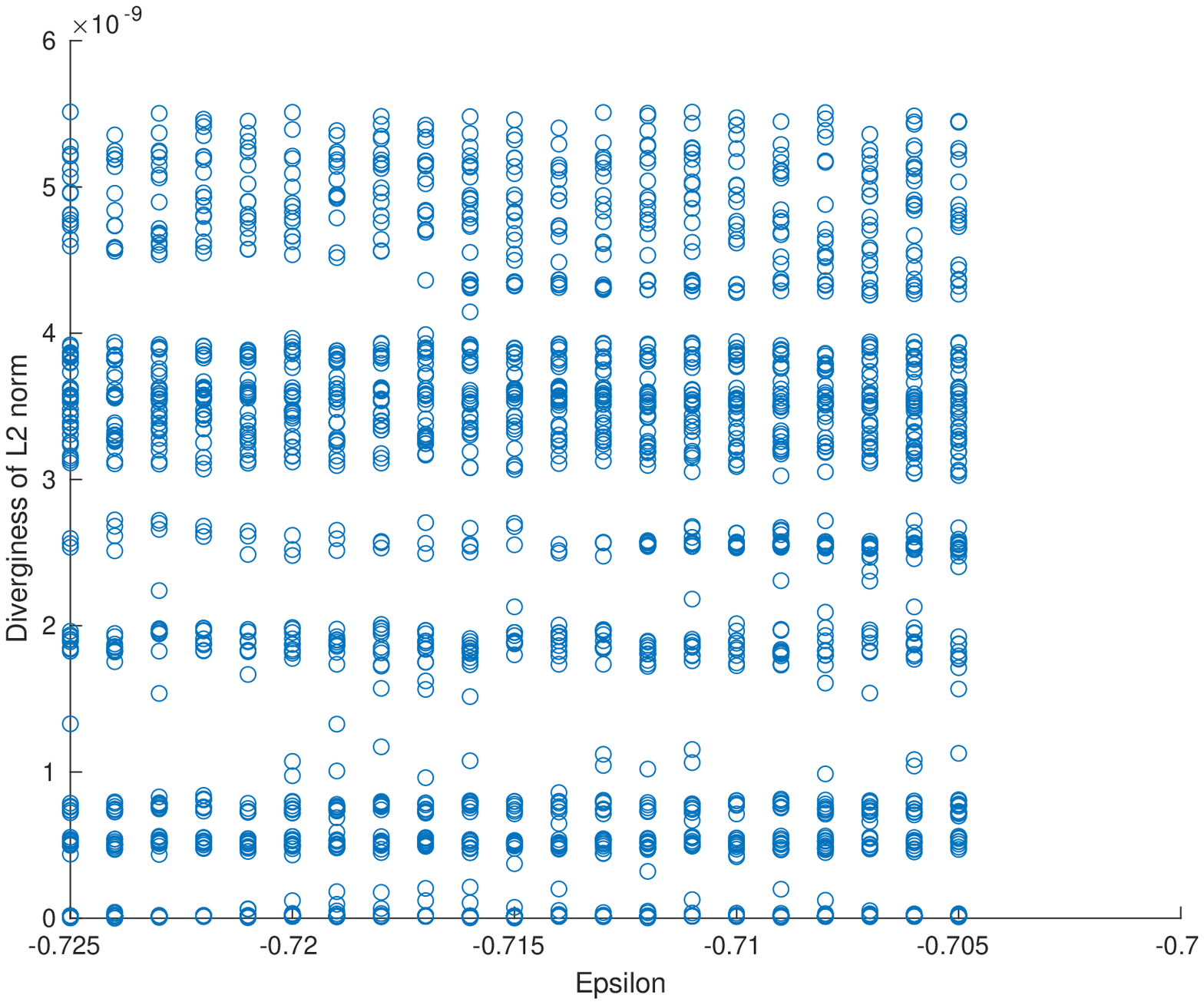}
  \caption{{\tiny} \space  For $N=25, \beta =3, h=0.5$ we simulated 1000 realizations of the disorder and 50 start values used for all of them. For $\varepsilon = -0.705 - 0.001 \cdot k, k=0,\ldots, 20$ we made up to 1.000 iterations of the Banach-$\varepsilon$-algorithm where we stopped when we hit the target level of stability (measured with respect to the {\it maximum absolute error}) of $10^{-7}$.  For all $\varepsilon$ we plotted the difference of the {\it mean squared errors} of the last two iterates of the solutions found with the respective $\varepsilon$.}
\end{figure}
\begin{figure}[H]
  \centering
  % \psfrag{2}{\tiny Two Step}
  \includegraphics[width=0.5\textwidth,  trim={0.45cm 0.4cm 0 0},clip]{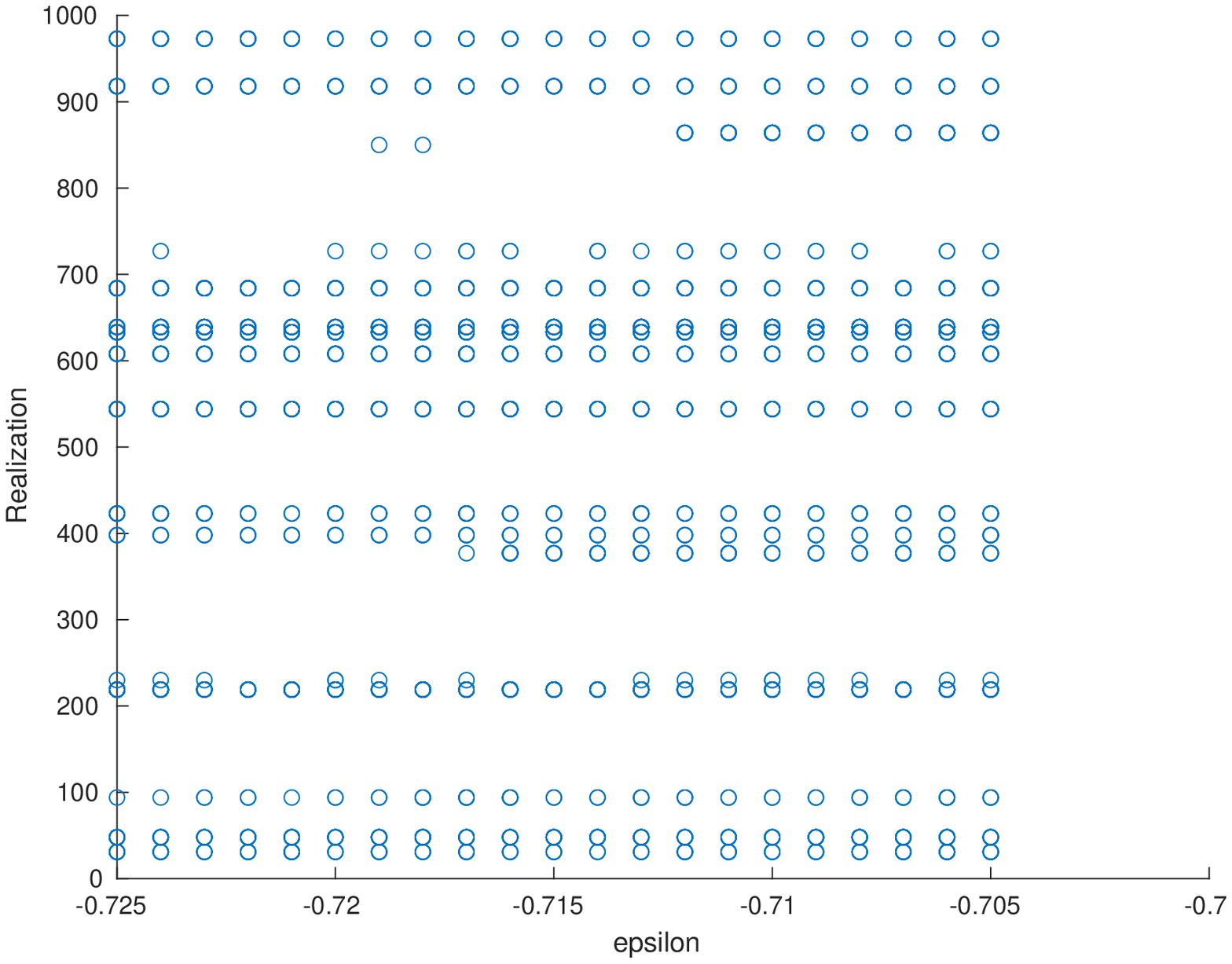}
  \caption{{\tiny} \space  For $N=25, \beta =3, h=0.5$ we simulated 1000 realizations of the disorder and 50 start values used for all of them. For $\varepsilon = -0.705 - 0.001 \cdot k, k=0,\ldots, 20$ we made up to 1.000 iterations of the Banach-$\varepsilon$-algorithm where we stopped when we hit the target level of stability (measured with respect to the {\it maximum absolute error}) of $10^{-7}$.  For all solutions we found we plotted a point with $x$ coordinate specifying the used $\varepsilon$ and $y$ coordinate the concerning realization}
\end{figure}
\begin{figure}[H]
  \centering
  % \psfrag{2}{\tiny Two Step}
  \includegraphics[width=0.5\textwidth,  trim={0.45cm 0.4cm 0 0},clip]{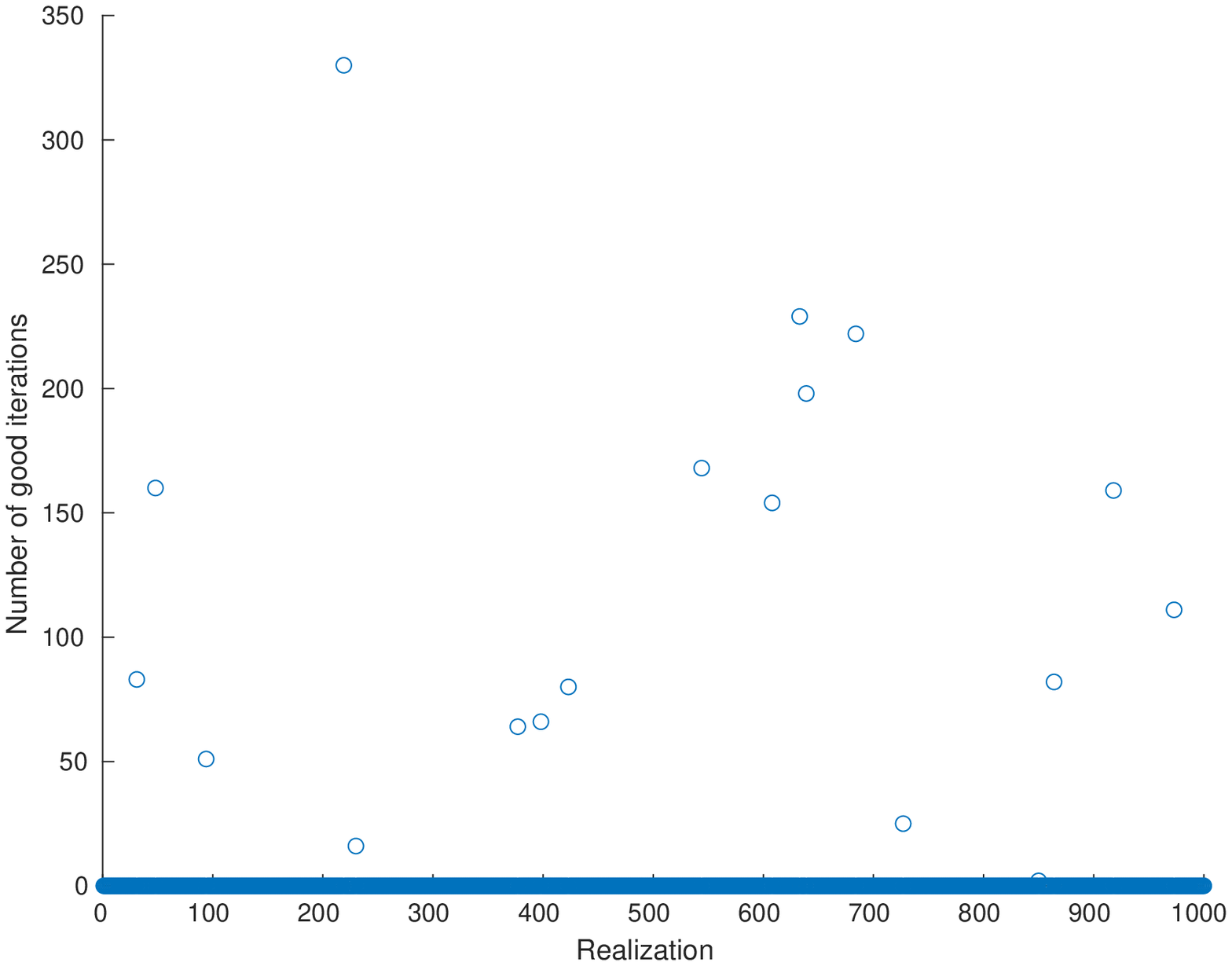}
  \caption{{\tiny} \space  For $N=25, \beta =3, h=0.5$ we simulated 1000 realizations of the disorder and 50 start values used for all of them. For $\varepsilon = -0.705 - 0.001 \cdot k, k=0,\ldots, 20$ we made up to 1.000 iterations of the Banach-$\varepsilon$-algorithm where we stopped when we hit the target level of stability (measured with respect to the {\it maximum absolute error}) of $10^{-7}$.  For all realizations we plotted the number of good solutions. }
\end{figure}
\begin{figure}[H]
  \centering
  % \psfrag{2}{\tiny Two Step}
  \includegraphics[width=0.5\textwidth,  trim={0.45cm 0.4cm 0 0},clip]{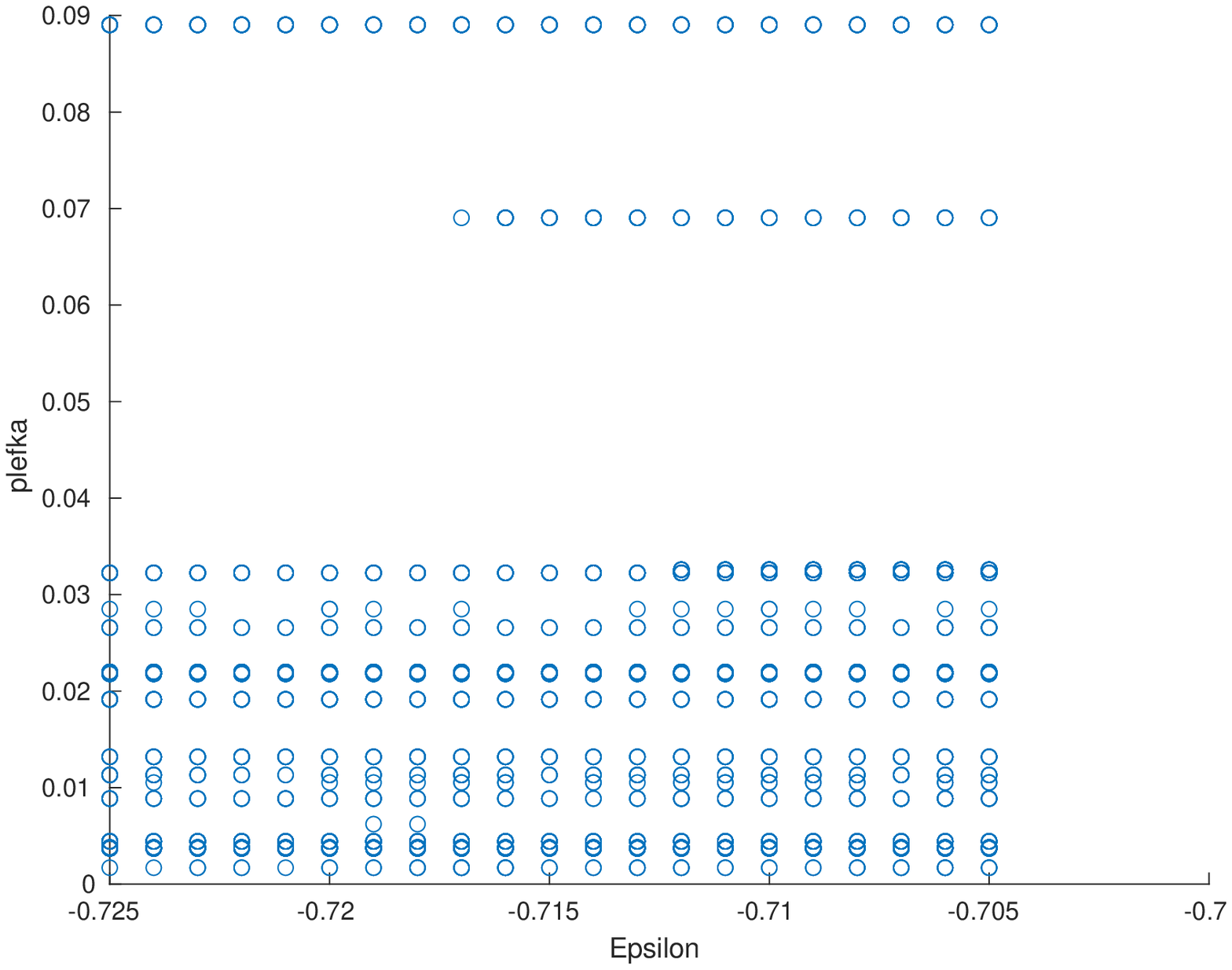}
  \caption{{\tiny} \space  For $N=25, \beta =3, h=0.5$ we simulated 1000 realizations of the disorder and 50 start values used for all of them. For $\varepsilon = -0.705 - 0.001 \cdot k, k=0,\ldots, 20$ we made up to 1.000 iterations of the Banach-$\varepsilon$-algorithm where we stopped when we hit the target level of stability (measured with respect to the {\it maximum absolute error}) of $10^{-7}$.  For all $\varepsilon$ we plotted the Plefka values found with this $\varepsilon$. }
\end{figure}
\begin{figure}[H]
  \centering
  % \psfrag{2}{\tiny Two Step}
  \includegraphics[width=0.5\textwidth,  trim={0.45cm 0.4cm 0 0},clip]{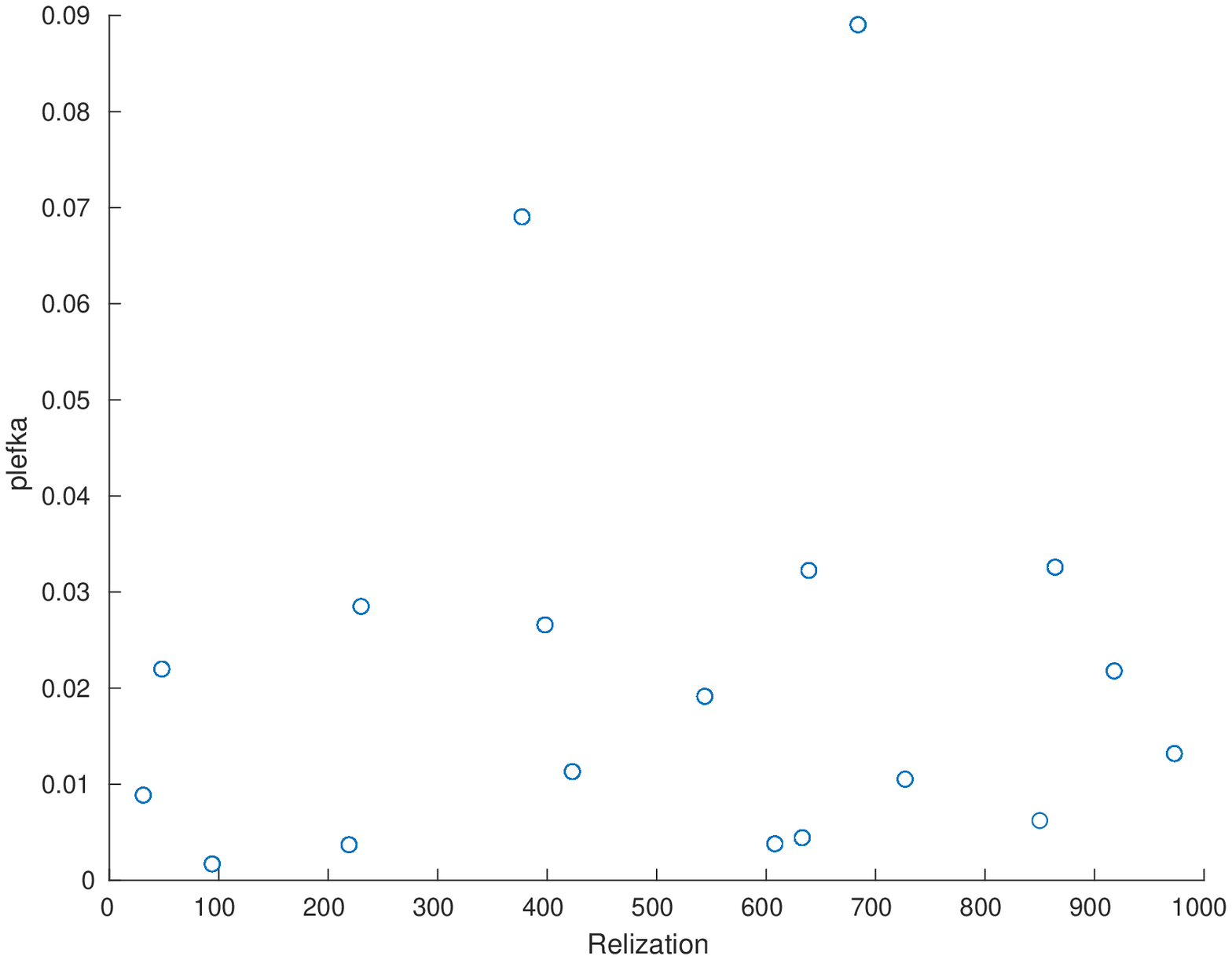}
  \caption{{\tiny} \space  For $N=25, \beta =3, h=0.5$ we simulated 1000 realizations of the disorder and 50 start values used for all of them. For $\varepsilon = -0.705 - 0.001 \cdot k, k=0,\ldots, 20$ we made up to 1.000 iterations of the Banach-$\varepsilon$-algorithm where we stopped when we hit the target level of stability (measured with respect to the {\it maximum absolute error}) of $10^{-7}$.  For all realizations we plotted the Plefka values found in the concerning realization. }
\end{figure}
\begin{figure}[H]
  \centering
  % \psfrag{2}{\tiny Two Step}
  \includegraphics[width=0.5\textwidth,  trim={0.45cm 0.4cm 0 0},clip]{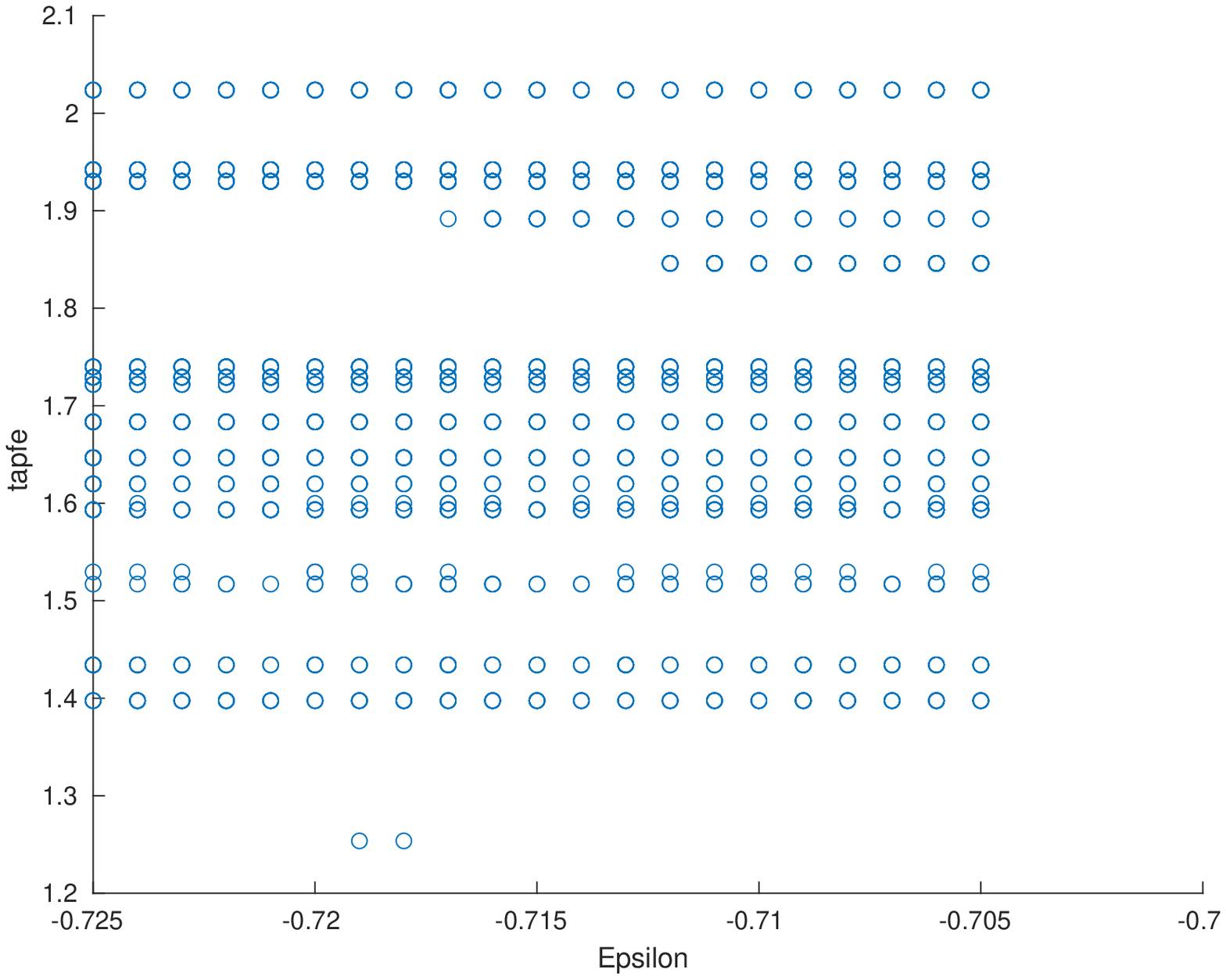}
  \caption{{\tiny} \space  For $N=25, \beta =3, h=0.5$ we simulated 1000 realizations of the disorder and 50 start values used for all of them. For $\varepsilon = -0.705 - 0.001 \cdot k, k=0,\ldots, 20$ we made up to 1.000 iterations of the Banach-$\varepsilon$-algorithm where we stopped when we hit the target level of stability (measured with respect to the {\it maximum absolute error}) of $10^{-7}$.  For all $\varepsilon$ we plotted the TAP free energy values found with this $\varepsilon$. }
\end{figure}
\begin{figure}[H]
  \centering
  % \psfrag{2}{\tiny Two Step}
  \includegraphics[width=0.5\textwidth,  trim={0.45cm 0.4cm 0 0},clip]{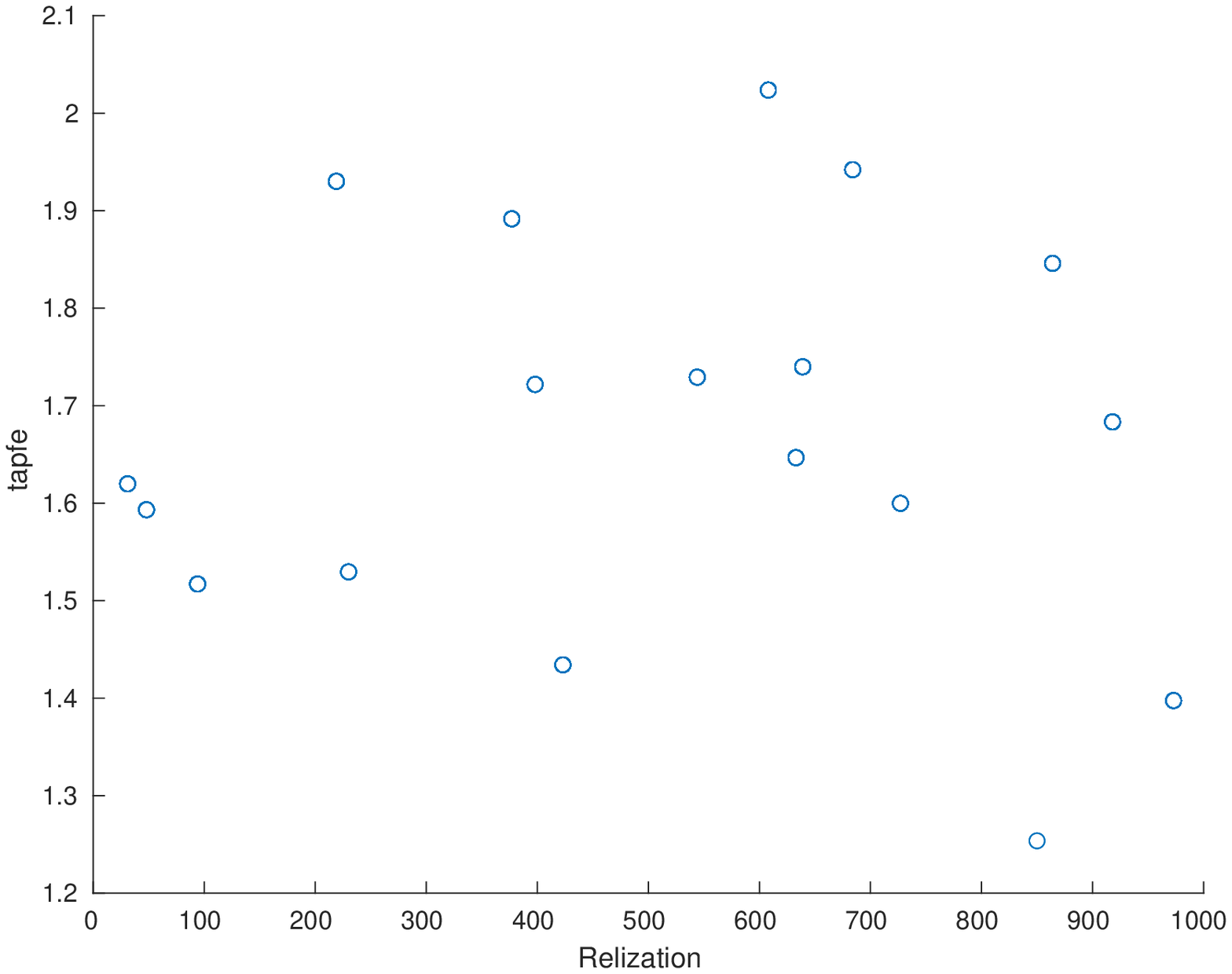}
  \caption{{\tiny} \space  For $N=25, \beta =3, h=0.5$ we simulated 1000 realizations of the disorder and 50 start values used for all of them. For $\varepsilon = -0.705 - 0.001 \cdot k, k=0,\ldots, 20$ we made up to 1.000 iterations of the Banach-$\varepsilon$-algorithm where we stopped when we hit the target level of stability (measured with respect to the {\it maximum absolute error}) of $10^{-7}$.  For all realizations we plotted the TAP free energy values found for the concerning realization.}
\end{figure}
\begin{figure}[H]
  \centering
  % \psfrag{2}{\tiny Two Step}
  \includegraphics[width=0.5\textwidth,  trim={0cm 0.4cm 0 0},clip]{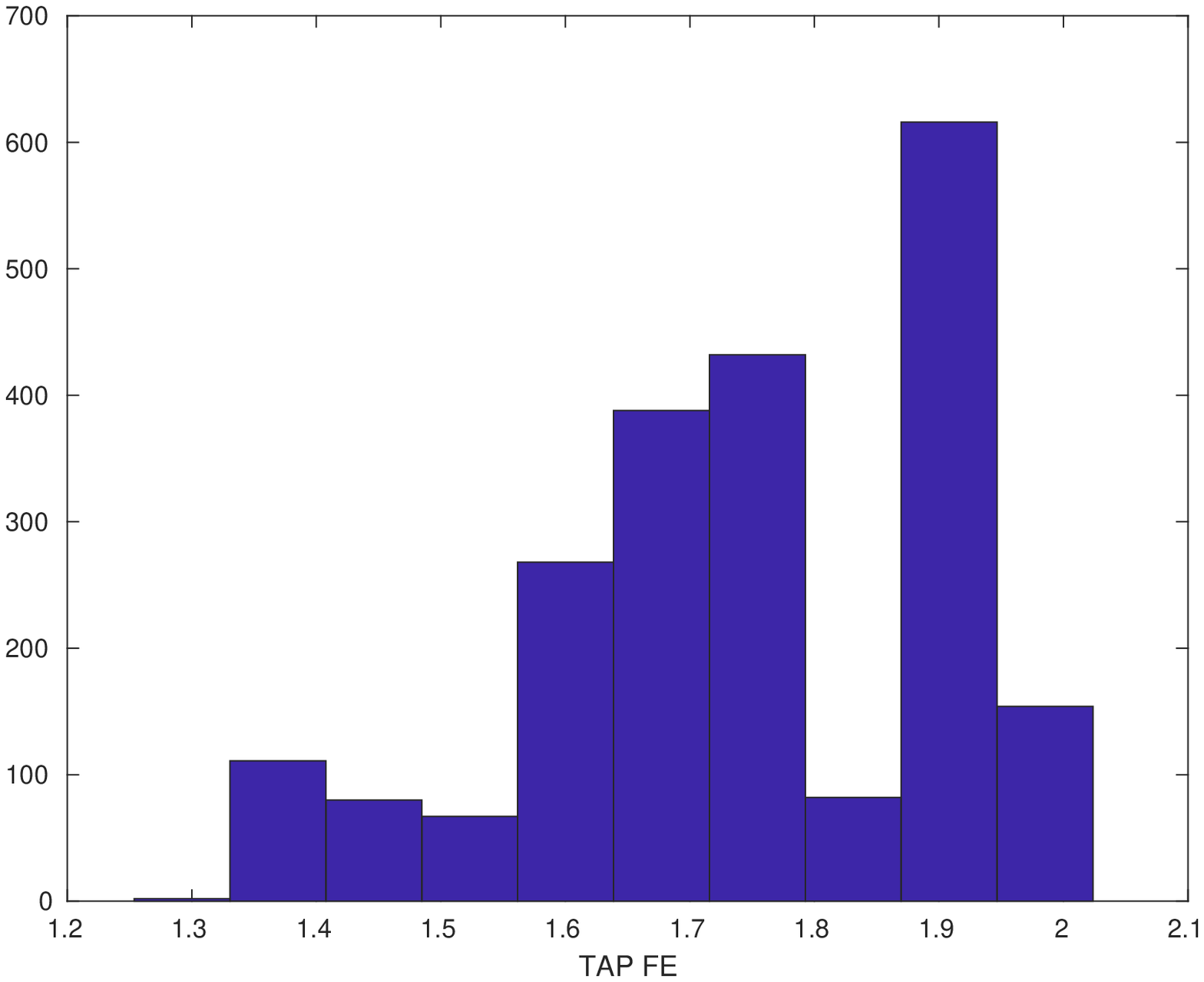}
  \caption{{\tiny} \space  For $N=25, \beta =3, h=0.5$ we simulated 1000 realizations of the disorder and 50 start values used for all of them. For $\varepsilon = -0.705 - 0.001 \cdot k, k=0,\ldots, 20$ we made up to 1.000 iterations of the Banach-$\varepsilon$-algorithm where we stopped when we hit the target level of stability (measured with respect to the {\it maximum absolute error}) of $10^{-7}$. Just a histogram of the TAP free energy values found.}
\end{figure}
\begin{figure}[H]
  \centering
  % \psfrag{2}{\tiny Two Step}
  \includegraphics[width=0.5\textwidth,  trim={0.45cm 0.4cm 0 0},clip]{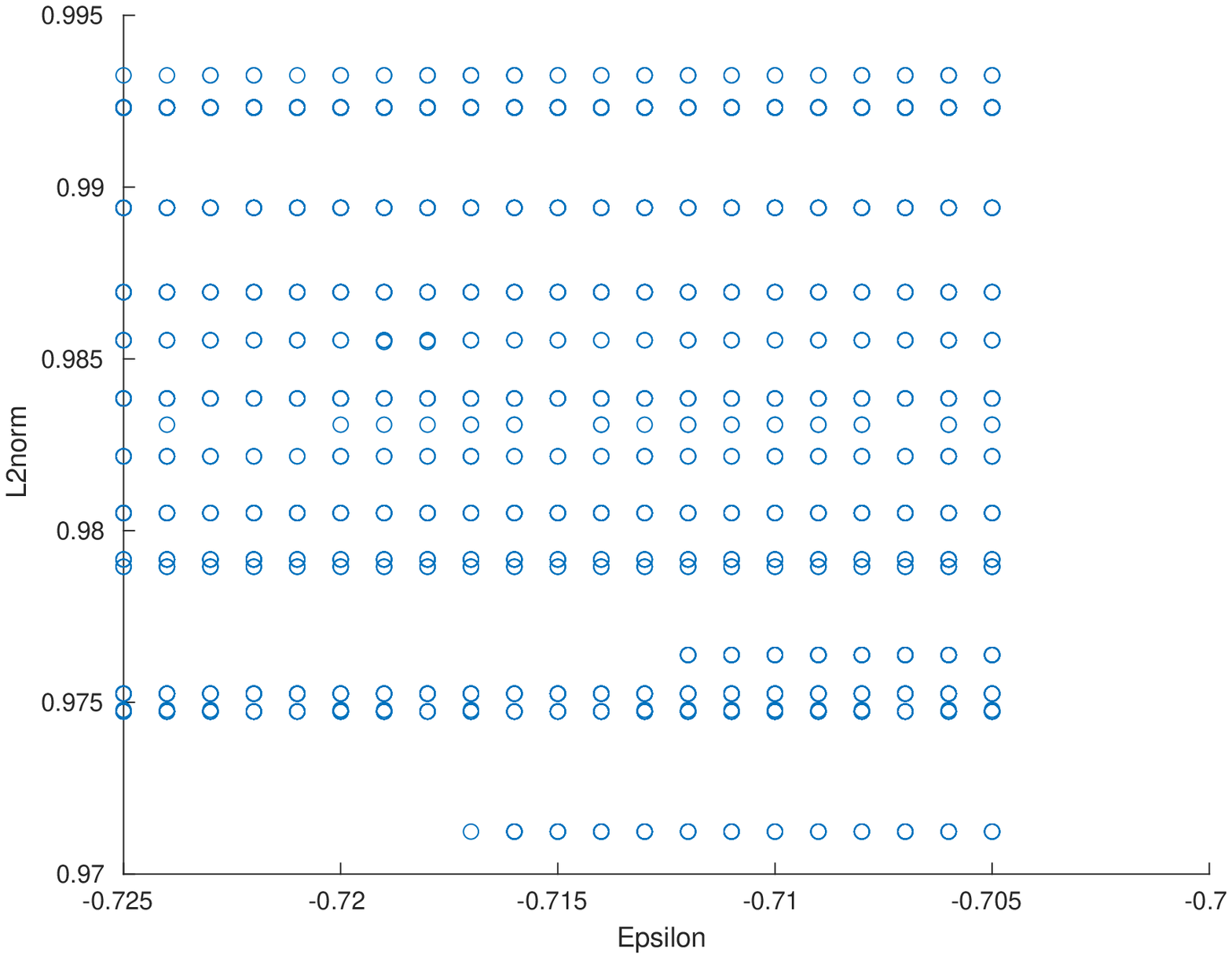}
  \caption{{\tiny} \space  For $N=25, \beta =3, h=0.5$ we simulated 1000 realizations of the disorder and 50 start values used for all of them. For $\varepsilon = -0.705 - 0.001 \cdot k, k=0,\ldots, 20$ we made up to 1.000 iterations of the Banach-$\varepsilon$-algorithm where we stopped when we hit the target level of stability (measured with respect to the {\it maximum absolute error}) of $10^{-7}$. For all $\varepsilon$ we plotted the {\it mean squared errors} found with this $\varepsilon$. }
\end{figure}
\begin{figure}[H]
  \centering
  % \psfrag{2}{\tiny Two Step}
  \includegraphics[width=0.5\textwidth,  trim={0.45cm 0.4cm 0 0},clip]{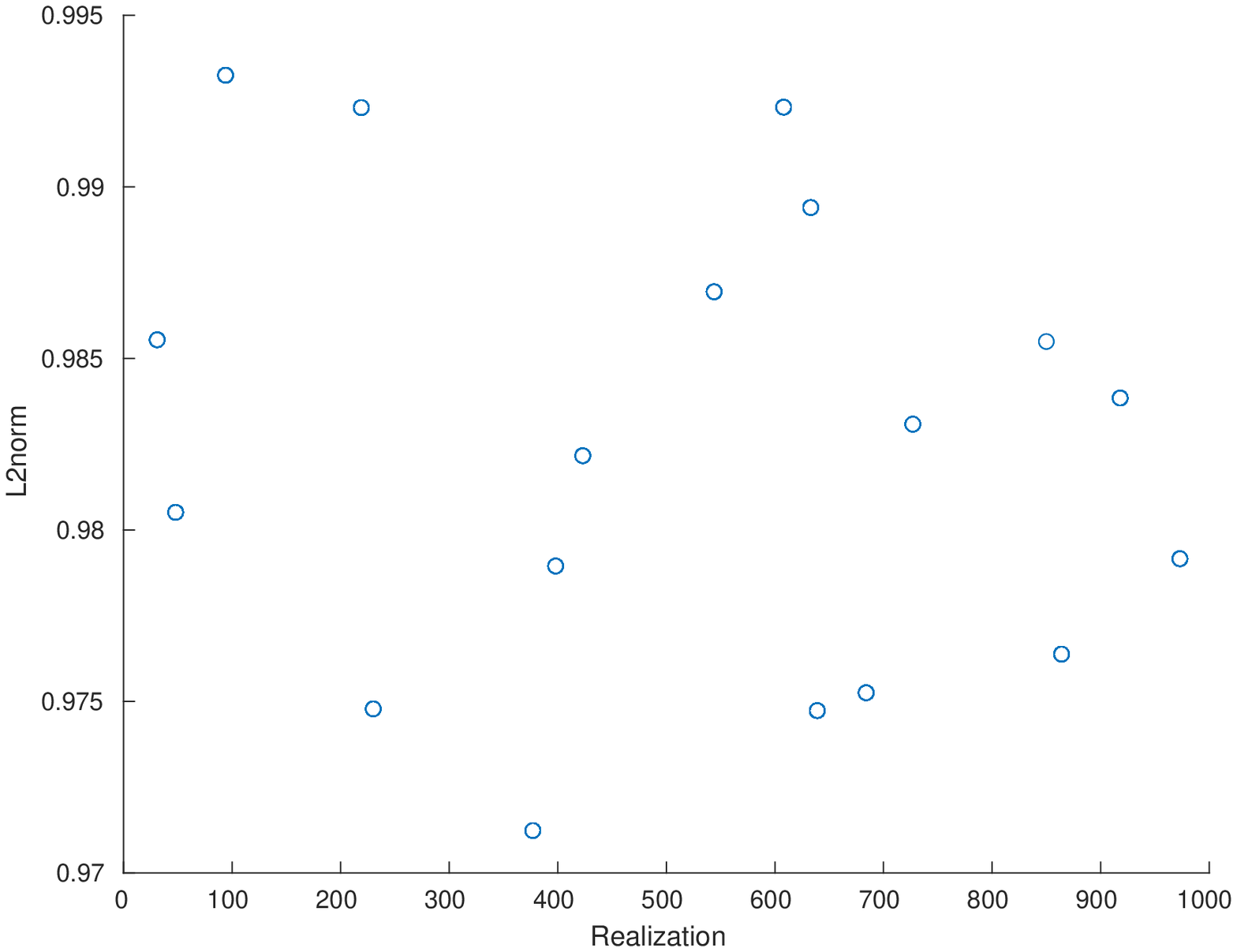}
  \caption{{\tiny} \space  For $N=25, \beta =3, h=0.5$ we simulated 1000 realizations of the disorder and 50 start values used for all of them. For $\varepsilon = -0.705 - 0.001 \cdot k, k=0,\ldots, 20$ we made up to 1.000 iterations of the Banach-$\varepsilon$-algorithm where we stopped when we hit the target level of stability (measured with respect to the {\it maximum absolute error}) of $10^{-7}$. For all realizations we plotted the {\it mean squared errors} found for the concerning realization.}
\end{figure}

\FloatBarrier\clearpage
\subsection{1000 realizations of disorder and 500 start values and 21 $\varepsilon$-values}\leavevmode
The following plots concern the same iterations with the same realizations, same start values.  We only accepted solutions satisfying the Plefka condition and being inside of the hypercube. 
\begin{figure}[H]
  \centering
  % \psfrag{2}{\tiny Two Step}
  \includegraphics[width=0.5\textwidth,  trim={0.45cm 0.4cm 0 0},clip]{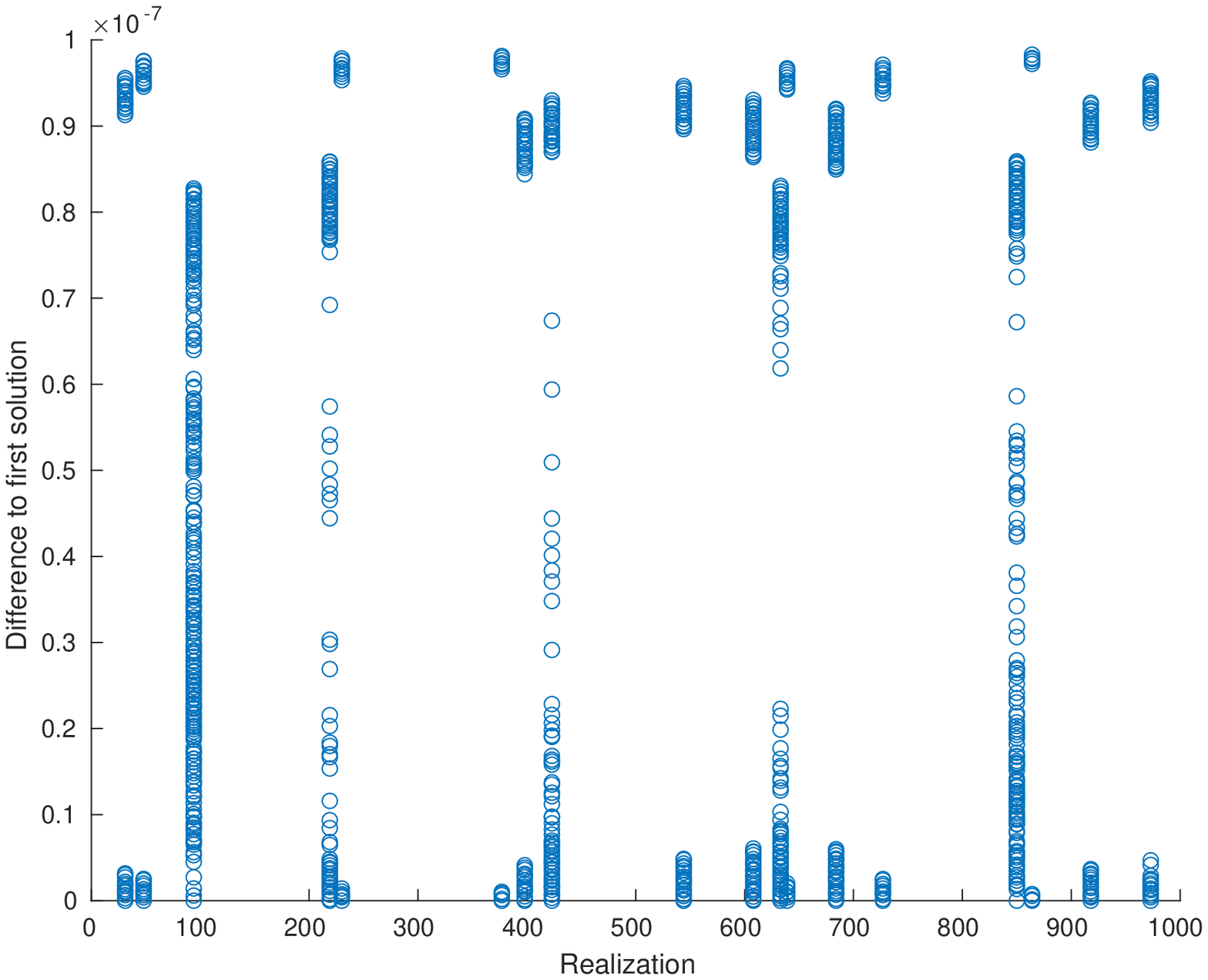}
  \caption{{\tiny} \space  For $N=25, \beta =3, h=0.5$ we simulated 1000 realizations of the disorder and 500 start values used for all of them. For $\varepsilon = -0.705 - 0.001 \cdot k, k=0,\ldots, 20$ we made up to 1.000 iterations of the Banach-$\varepsilon$-algorithm where we stopped when we hit the target level of stability (measured with respect to the {\it maximum absolute error}) of $10^{-7}$.  For all realizations we plotted the {\it maximum absolute error} to the first solution of the realization we found.}
\end{figure}
\begin{figure}[H]
  \centering
  % \psfrag{2}{\tiny Two Step}
  \includegraphics[width=0.5\textwidth,  trim={0.45cm 0.4cm 0 0},clip]{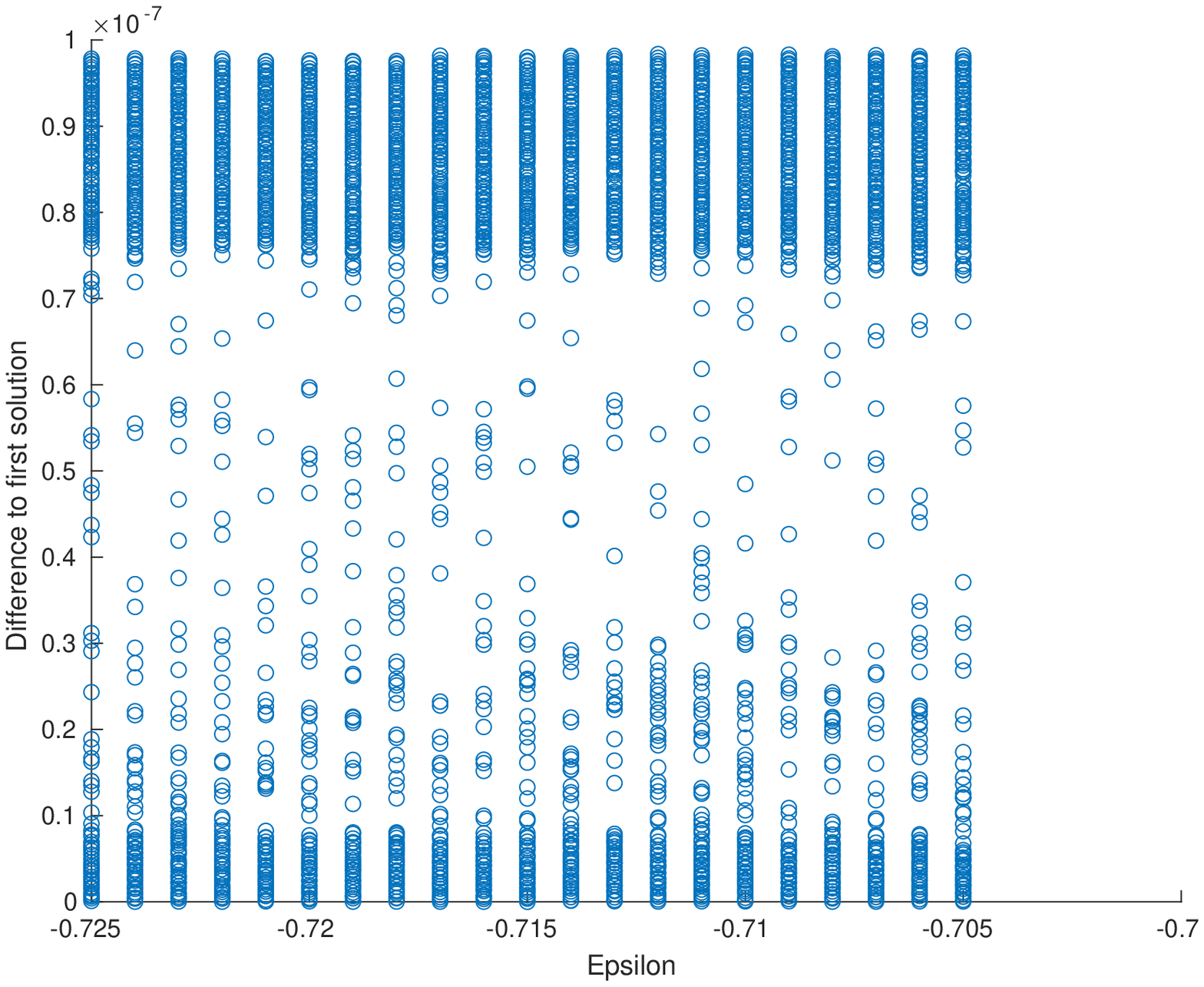}
  \caption{{\tiny} \space  For $N=25, \beta =3, h=0.5$ we simulated 1000 realizations of the disorder and 500 start values used for all of them. For $\varepsilon = -0.705 - 0.001 \cdot k, k=0,\ldots, 20$ we made up to 1.000 iterations of the Banach-$\varepsilon$-algorithm where we stopped when we hit the target level of stability (measured with respect to the {\it maximum absolute error}) of $10^{-7}$.  For all $\varepsilon$ we plotted the {\it maximum absolute error} to the first solution of the realization we found.}
\end{figure}
\begin{figure}[H]
  \centering
  % \psfrag{2}{\tiny Two Step}
  \includegraphics[width=0.5\textwidth,  trim={0.45cm 0.4cm 0 0},clip]{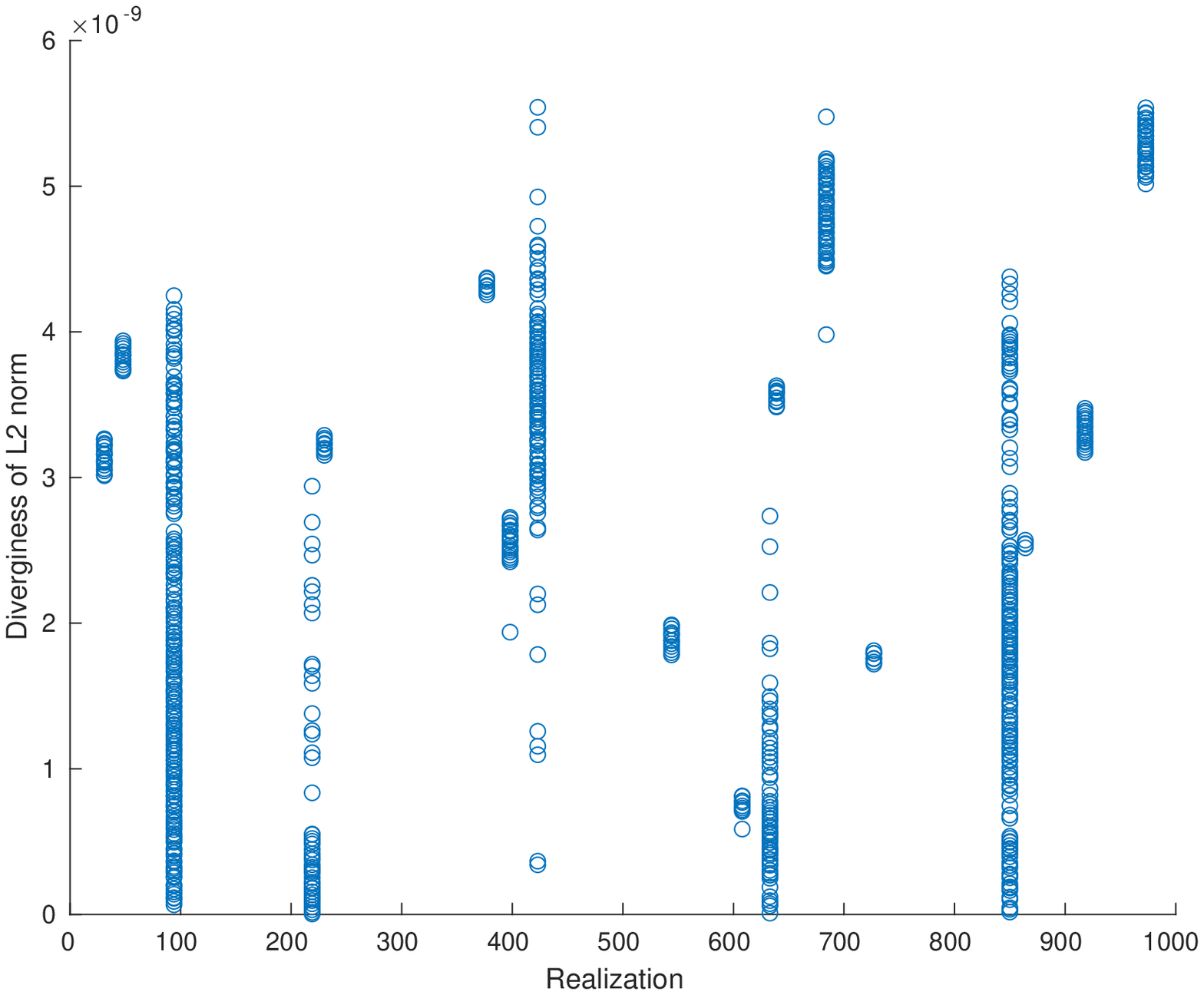}
  \caption{{\tiny} \space  For $N=25, \beta =3, h=0.5$ we simulated 1000 realizations of the disorder and 500 start values used for all of them. For $\varepsilon = -0.705 - 0.001 \cdot k, k=0,\ldots, 20$ we made up to 1.000 iterations of the Banach-$\varepsilon$-algorithm where we stopped when we hit the target level of stability (measured with respect to the {\it maximum absolute error}) of $10^{-7}$.  For all realizations we plotted the difference of the {\it mean squared errors} of the last two iterates of the solutions found in the respective realization.}
\end{figure}
\begin{figure}[H]
  \centering
  % \psfrag{2}{\tiny Two Step}
  \includegraphics[width=0.5\textwidth,  trim={0.45cm 0.4cm 0 0},clip]{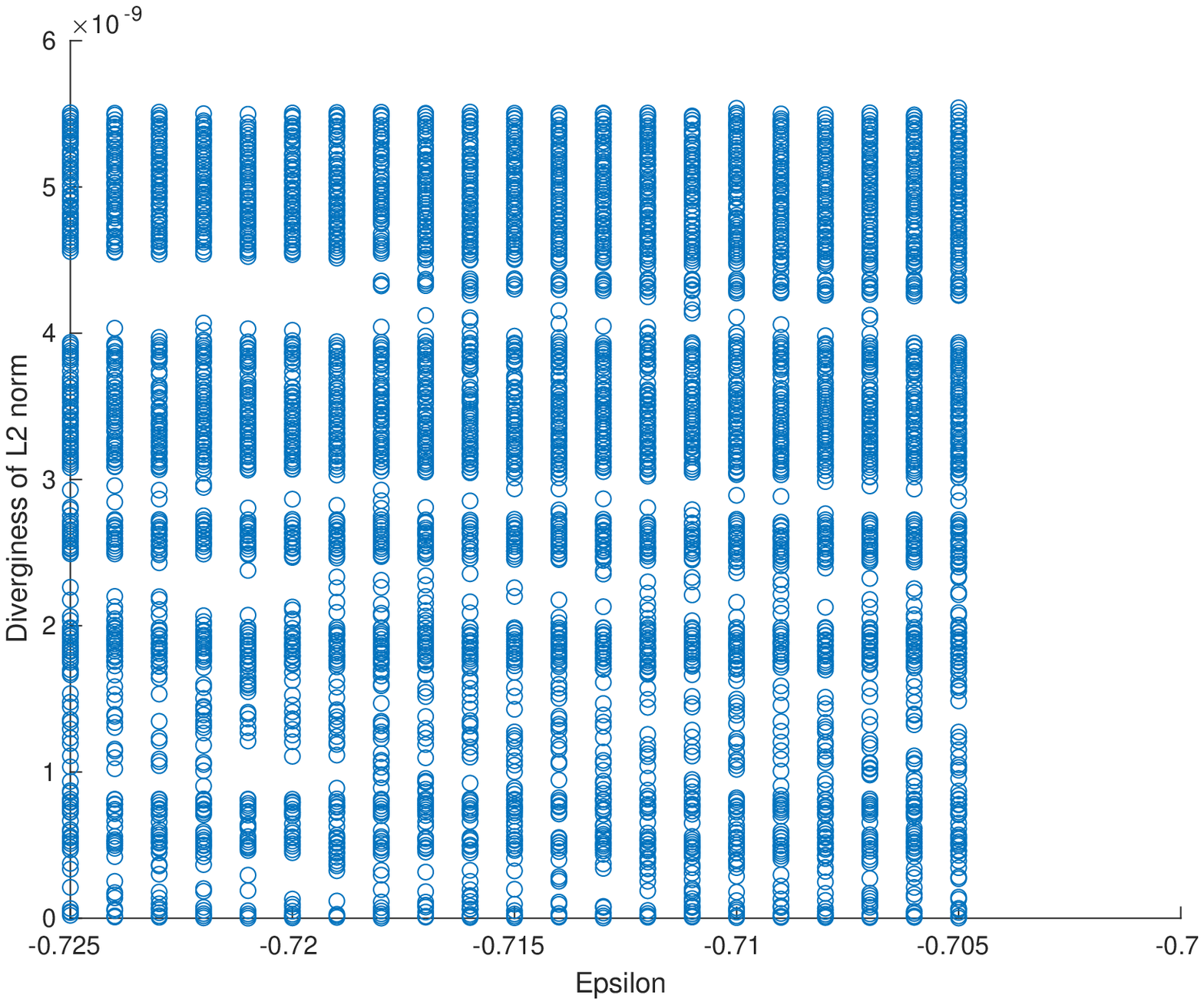}
  \caption{{\tiny} \space  For $N=25, \beta =3, h=0.5$ we simulated 1000 realizations of the disorder and 500 start values used for all of them. For $\varepsilon = -0.705 - 0.001 \cdot k, k=0,\ldots, 20$ we made up to 1.000 iterations of the Banach-$\varepsilon$-algorithm where we stopped when we hit the target level of stability (measured with respect to the {\it maximum absolute error}) of $10^{-7}$.  For all $\varepsilon$ we plotted the difference of the {\it mean squared errors} of the last two iterates of the solutions found with the respective $\varepsilon$.}
\end{figure}
\begin{figure}[H]
  \centering
  % \psfrag{2}{\tiny Two Step}
  \includegraphics[width=0.5\textwidth,  trim={0.45cm 0.4cm 0 0},clip]{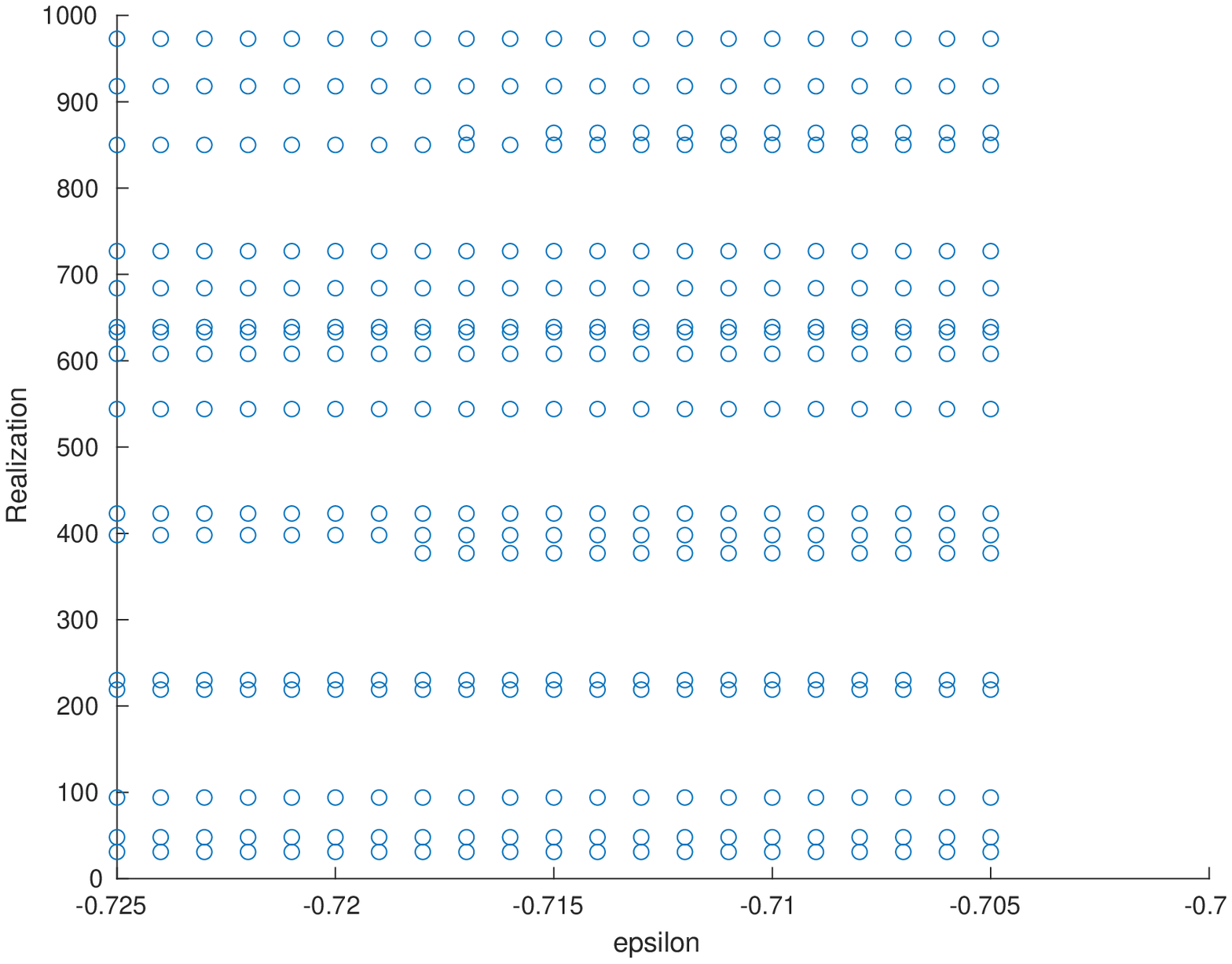}
  \caption{{\tiny} \space  For $N=25, \beta =3, h=0.5$ we simulated 1000 realizations of the disorder and 500 start values used for all of them. For $\varepsilon = -0.705 - 0.001 \cdot k, k=0,\ldots, 20$ we made up to 1.000 iterations of the Banach-$\varepsilon$-algorithm where we stopped when we hit the target level of stability (measured with respect to the {\it maximum absolute error}) of $10^{-7}$.  For all solutions we found we plotted a point with $x$ coordinate specifying the used $\varepsilon$ and $y$ coordinate the concerning realization}
\end{figure}
\begin{figure}[H]
  \centering
  % \psfrag{2}{\tiny Two Step}
  \includegraphics[width=0.5\textwidth,  trim={0.45cm 0.4cm 0 0},clip]{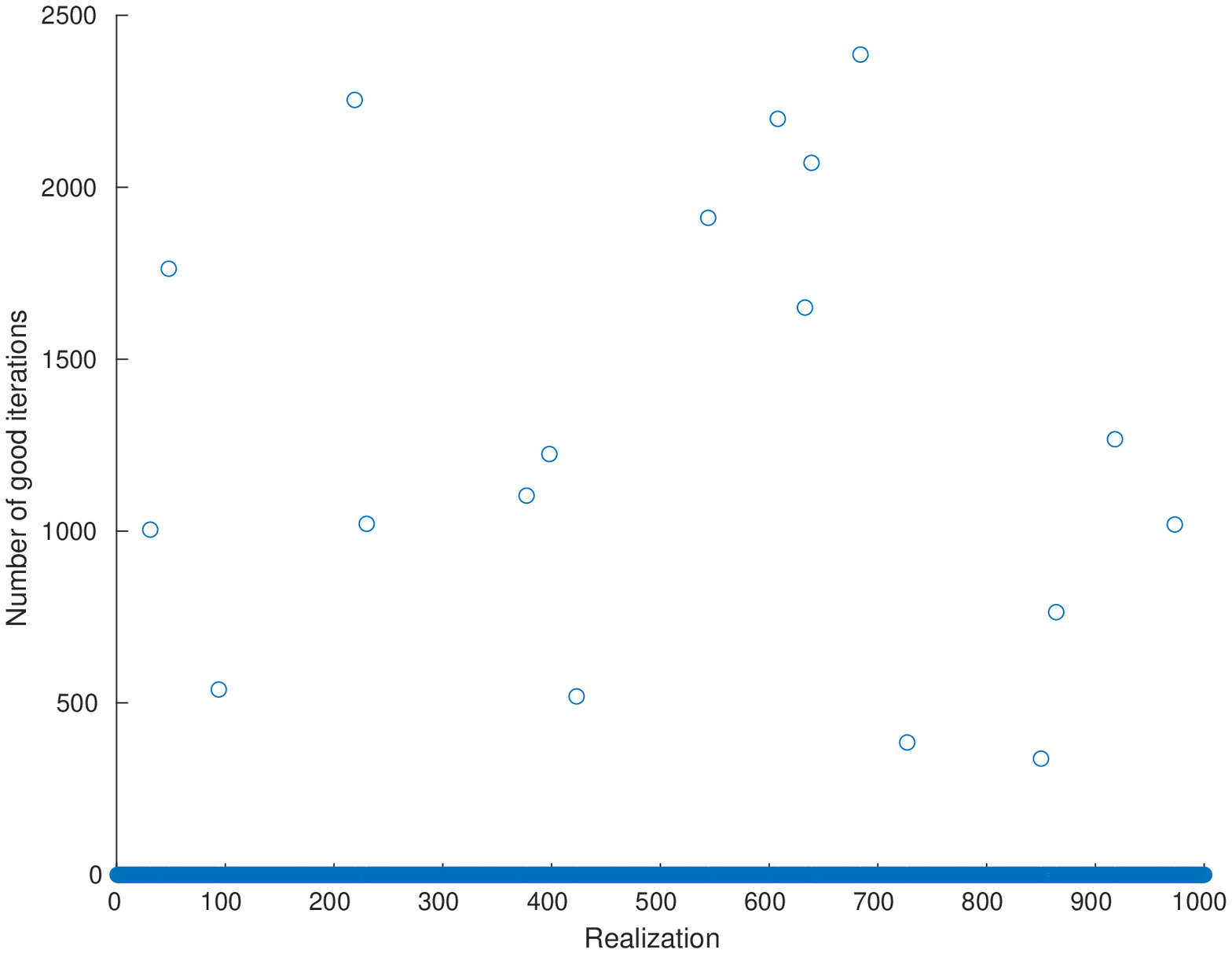}
  \caption{{\tiny} \space  For $N=25, \beta =3, h=0.5$ we simulated 1000 realizations of the disorder and 500 start values used for all of them. For $\varepsilon = -0.705 - 0.001 \cdot k, k=0,\ldots, 20$ we made up to 1.000 iterations of the Banach-$\varepsilon$-algorithm where we stopped when we hit the target level of stability (measured with respect to the {\it maximum absolute error}) of $10^{-7}$.  For all realizations we plotted the number of good solutions. }
\end{figure}
\begin{figure}[H]
  \centering
  % \psfrag{2}{\tiny Two Step}
  \includegraphics[width=0.5\textwidth,  trim={0.45cm 0.4cm 0 0},clip]{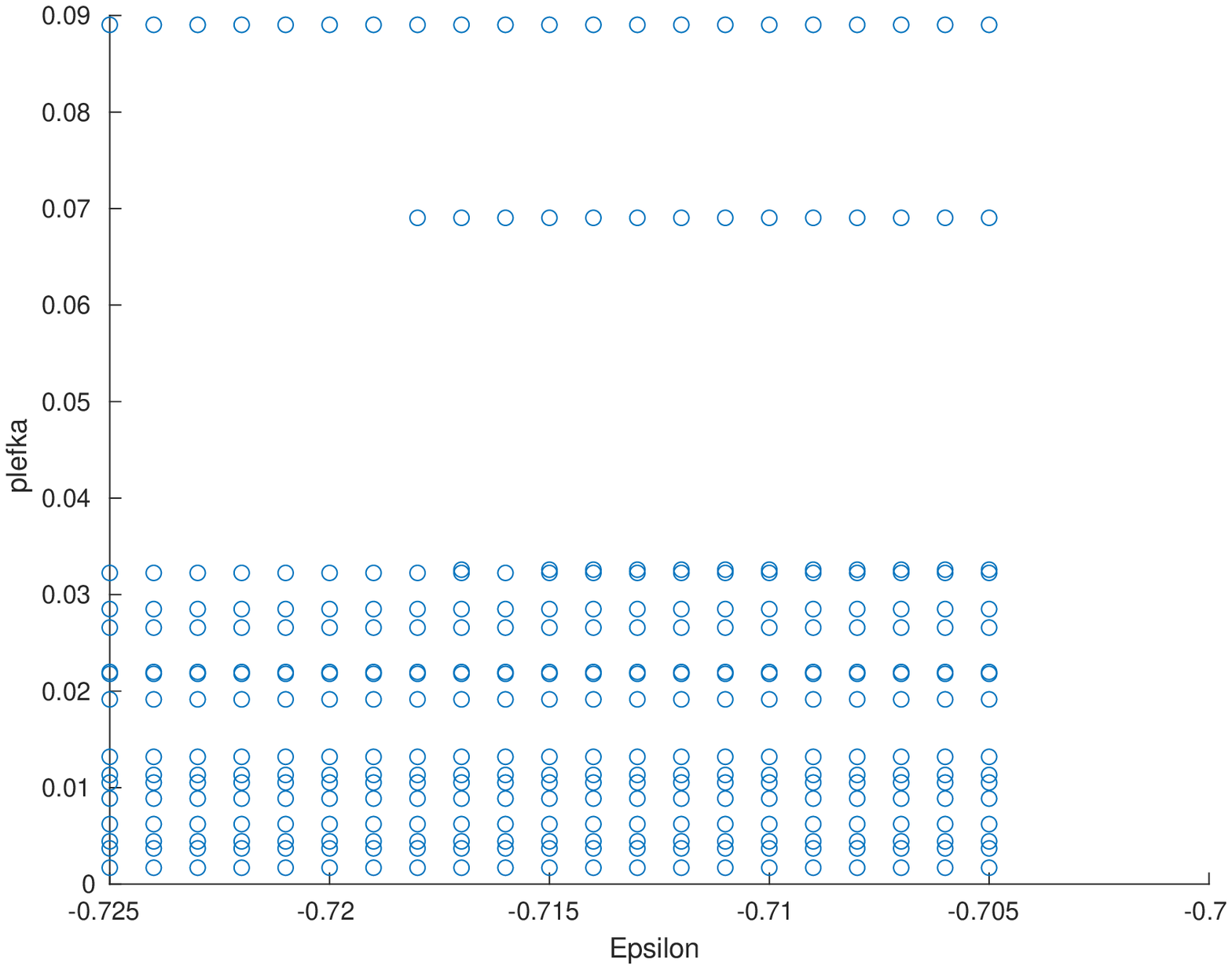}
  \caption{{\tiny} \space  For $N=25, \beta =3, h=0.5$ we simulated 1000 realizations of the disorder and 500 start values used for all of them. For $\varepsilon = -0.705 - 0.001 \cdot k, k=0,\ldots, 20$ we made up to 1.000 iterations of the Banach-$\varepsilon$-algorithm where we stopped when we hit the target level of stability (measured with respect to the {\it maximum absolute error}) of $10^{-7}$.  For all $\varepsilon$ we plotted the Plefka values found with this $\varepsilon$. }
\end{figure}
\begin{figure}[H]
  \centering
  % \psfrag{2}{\tiny Two Step}
  \includegraphics[width=0.5\textwidth,  trim={0.45cm 0.4cm 0 0},clip]{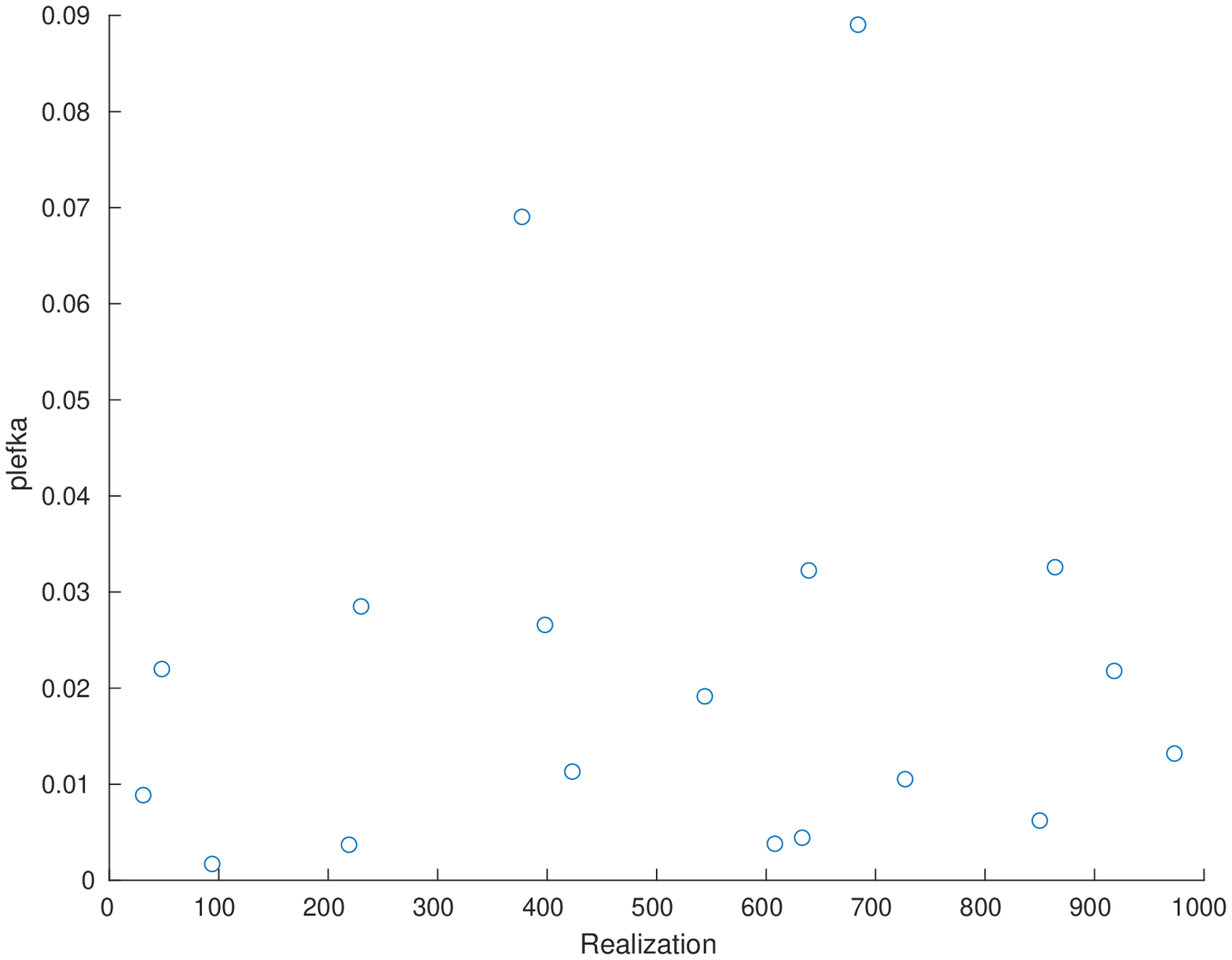}
  \caption{{\tiny} \space  For $N=25, \beta =3, h=0.5$ we simulated 1000 realizations of the disorder and 500 start values used for all of them. For $\varepsilon = -0.705 - 0.001 \cdot k, k=0,\ldots, 20$ we made up to 1.000 iterations of the Banach-$\varepsilon$-algorithm where we stopped when we hit the target level of stability (measured with respect to the {\it maximum absolute error}) of $10^{-7}$.  For all realizations we plotted the Plefka values found in the concerning realization. }
\end{figure}
\begin{figure}[H]
  \centering
  % \psfrag{2}{\tiny Two Step}
  \includegraphics[width=0.5\textwidth,  trim={0.45cm 0.4cm 0 0},clip]{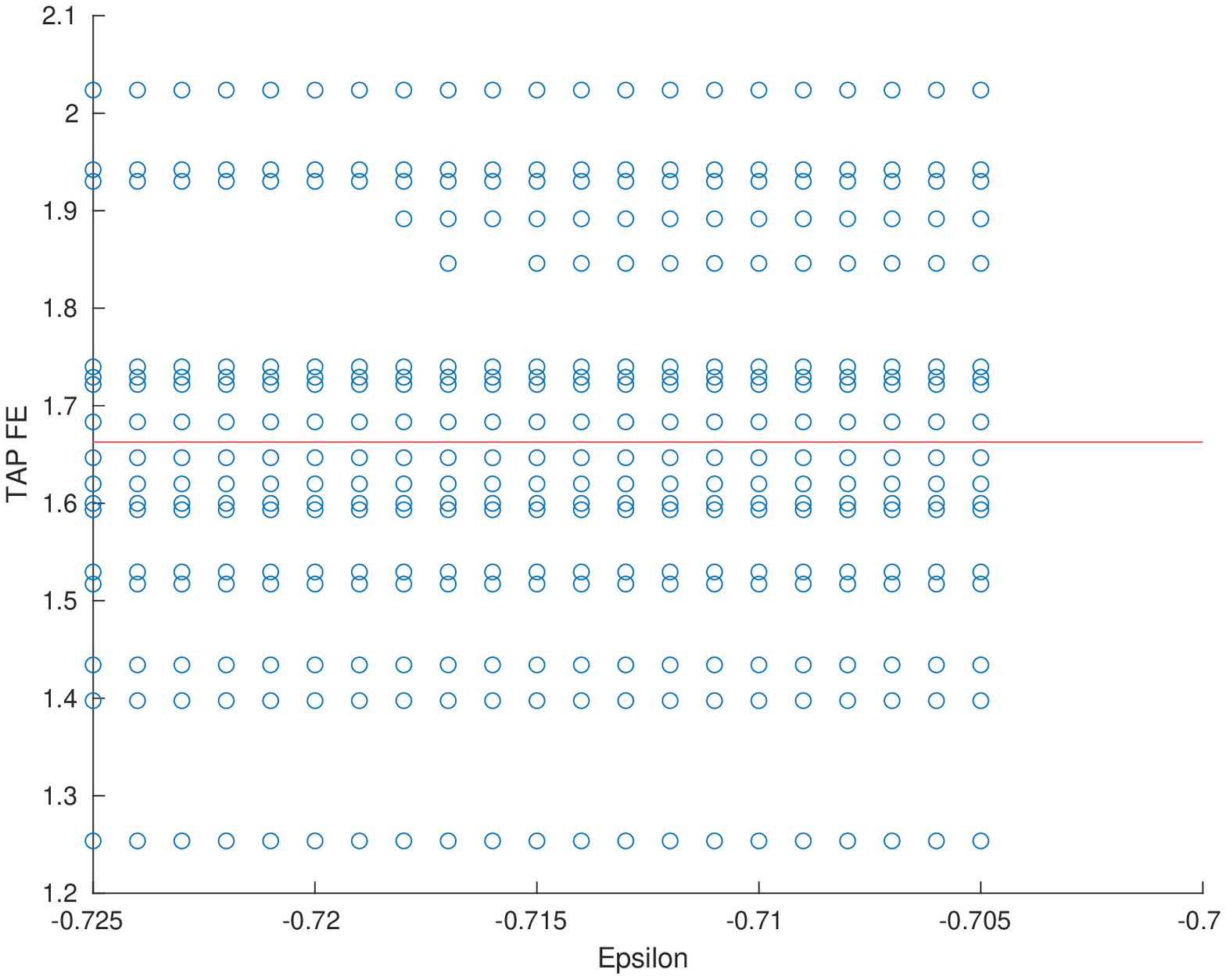}
  \caption{{\tiny} \space  For $N=25, \beta =3, h=0.5$ we simulated 1000 realizations of the disorder and 500 start values used for all of them. For $\varepsilon = -0.705 - 0.001 \cdot k, k=0,\ldots, 20$ we made up to 1.000 iterations of the Banach-$\varepsilon$-algorithm where we stopped when we hit the target level of stability (measured with respect to the {\it maximum absolute error}) of $10^{-7}$.  For all $\varepsilon$ we plotted the TAP free energy values found with this $\varepsilon$. }
\end{figure}
\begin{figure}[H]
  \centering
  % \psfrag{2}{\tiny Two Step}
  \includegraphics[width=0.5\textwidth,  trim={0.45cm 0.4cm 0 0},clip]{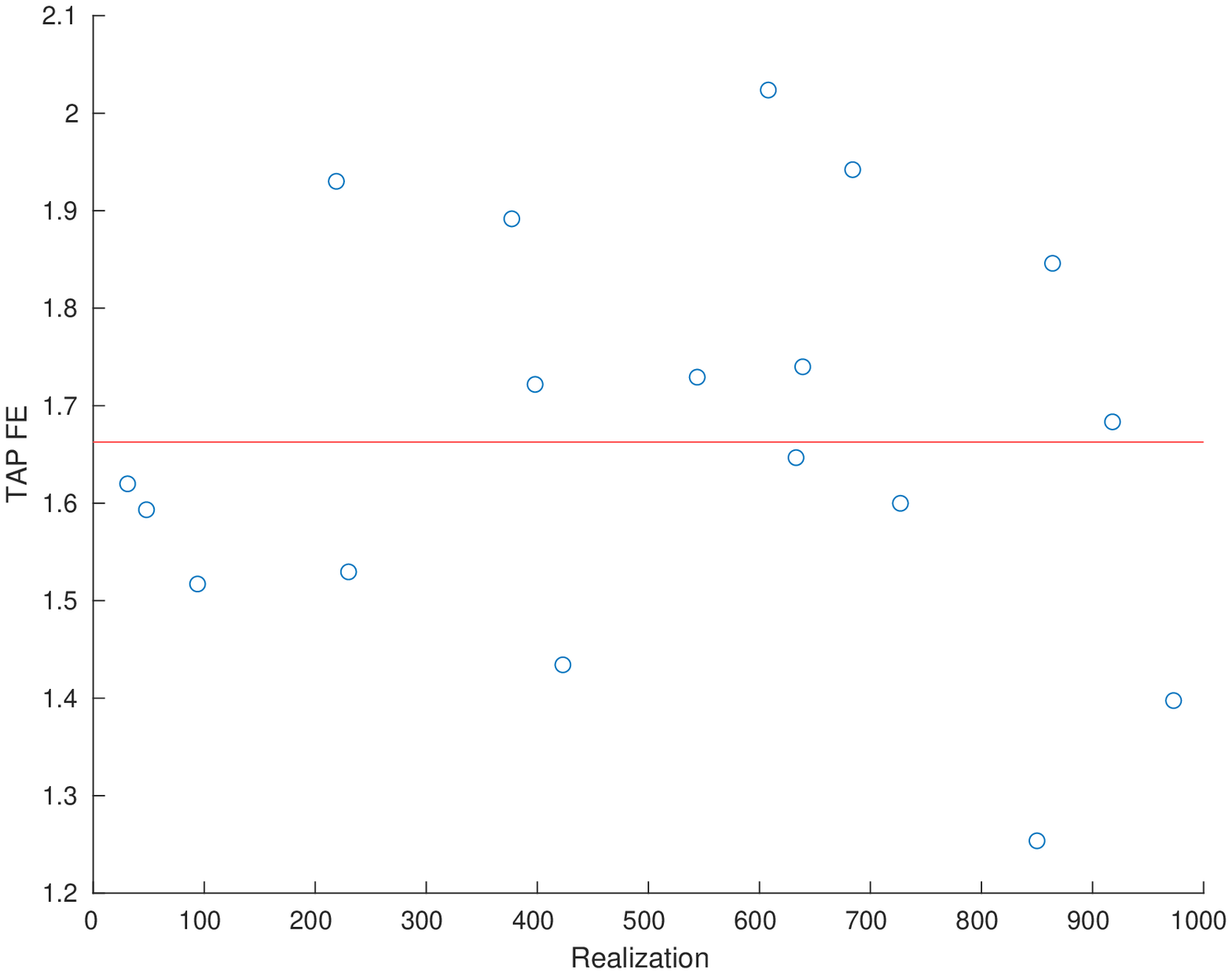}
  \caption{{\tiny} \space  For $N=25, \beta =3, h=0.5$ we simulated 1000 realizations of the disorder and 500 start values used for all of them. For $\varepsilon = -0.705 - 0.001 \cdot k, k=0,\ldots, 20$ we made up to 1.000 iterations of the Banach-$\varepsilon$-algorithm where we stopped when we hit the target level of stability (measured with respect to the {\it maximum absolute error}) of $10^{-7}$.  For all realizations we plotted the TAP free energy values found for the concerning realization.}
\end{figure}
\begin{figure}[H]
  \centering
  % \psfrag{2}{\tiny Two Step}
  \includegraphics[width=0.5\textwidth,  trim={0cm 0.4cm 0 0},clip]{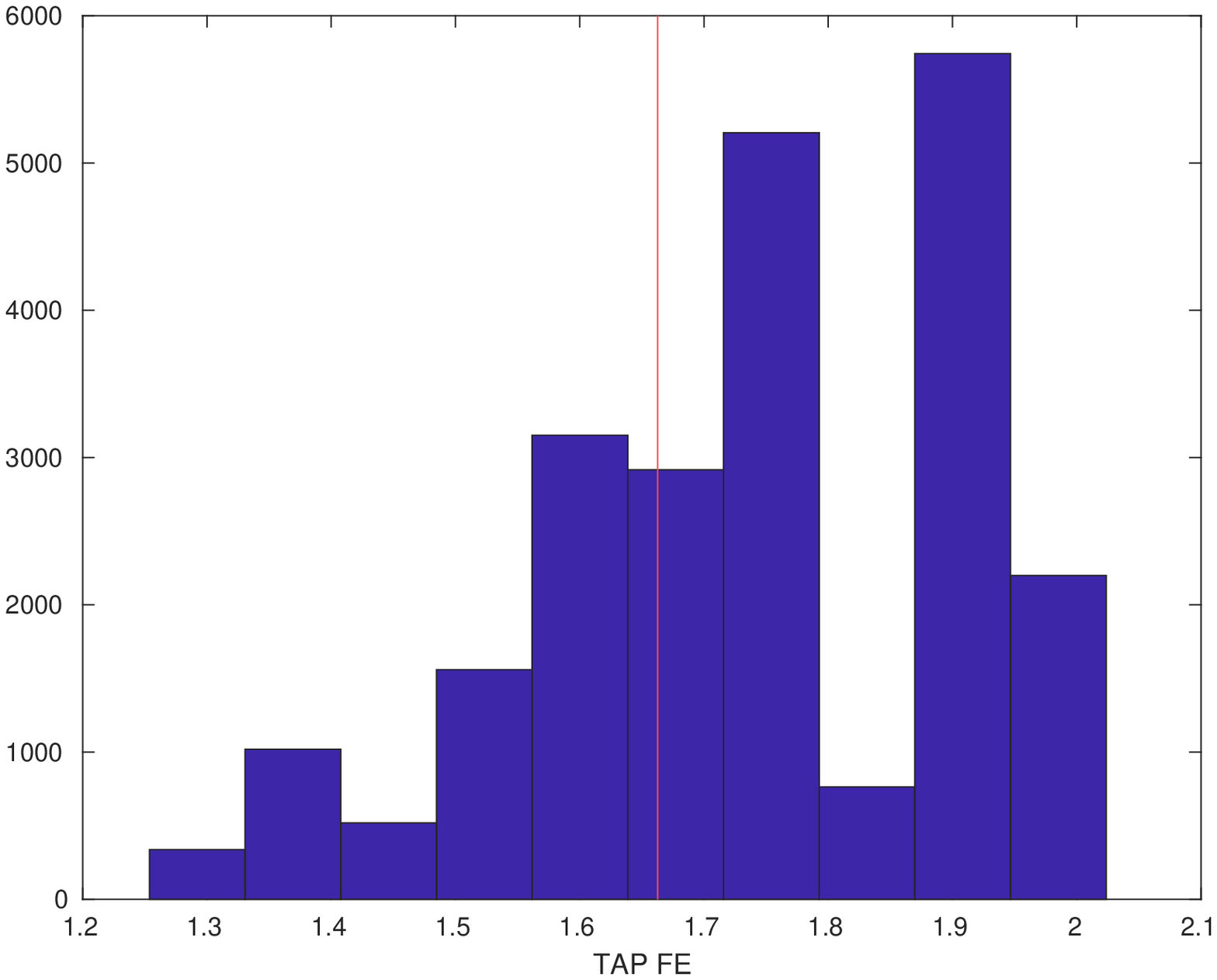}
  \caption{{\tiny} \space  For $N=25, \beta =3, h=0.5$ we simulated 1000 realizations of the disorder and 500 start values used for all of them. For $\varepsilon = -0.705 - 0.001 \cdot k, k=0,\ldots, 20$ we made up to 1.000 iterations of the Banach-$\varepsilon$-algorithm where we stopped when we hit the target level of stability (measured with respect to the {\it maximum absolute error}) of $10^{-7}$.  Just a histogram of the TAP free energy values found.}
\end{figure}
\begin{figure}[H]
  \centering
  % \psfrag{2}{\tiny Two Step}
  \includegraphics[width=0.5\textwidth,  trim={0.45cm 0.4cm 0 0},clip]{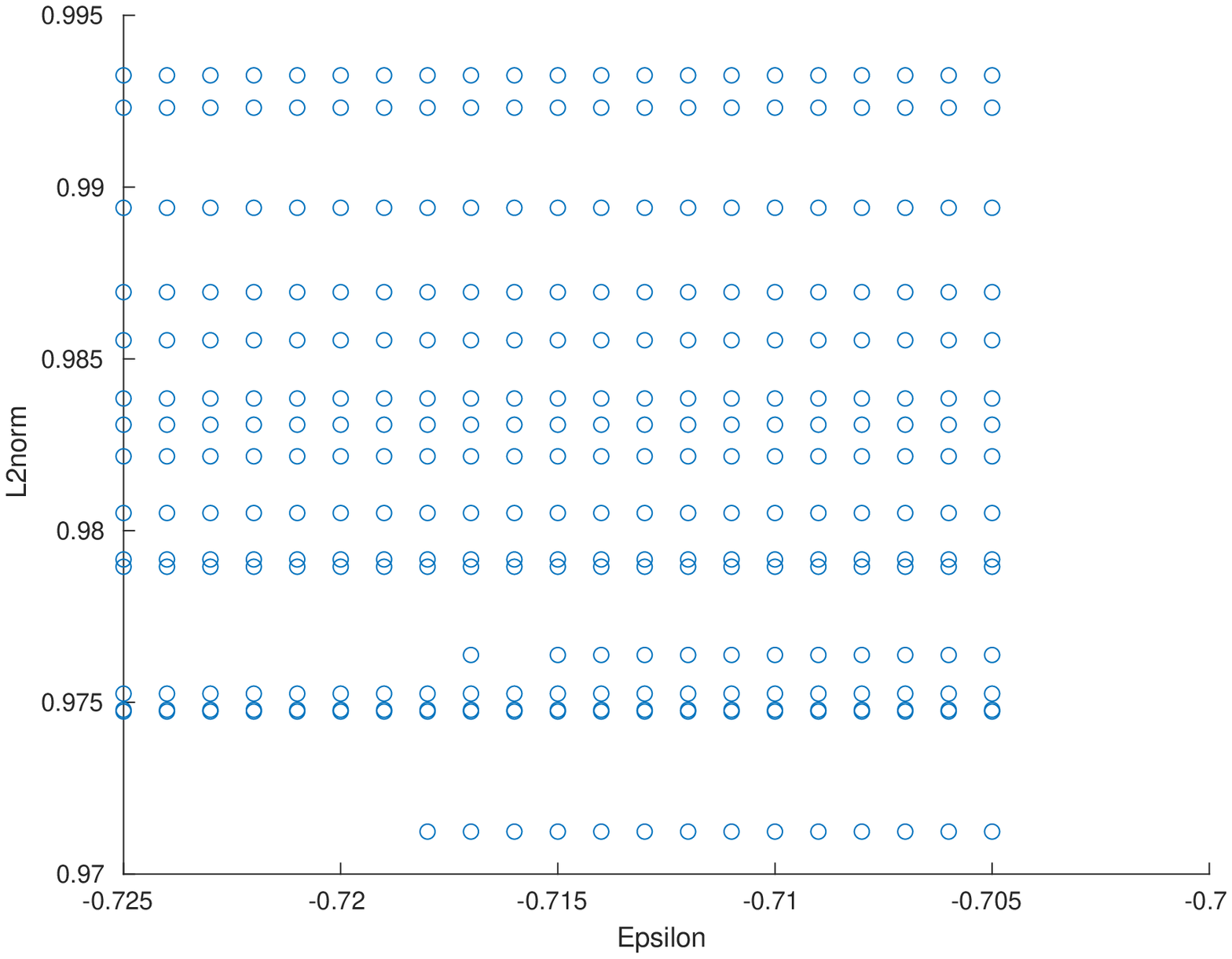}
  \caption{{\tiny} \space  For $N=25, \beta =3, h=0.5$ we simulated 1000 realizations of the disorder and 500 start values used for all of them. For $\varepsilon = -0.705 - 0.001 \cdot k, k=0,\ldots, 20$ we made up to 1.000 iterations of the Banach-$\varepsilon$-algorithm where we stopped when we hit the target level of stability (measured with respect to the {\it maximum absolute error}) of $10^{-7}$.  For all $\varepsilon$ we plotted the {\it mean squared errors} found with this $\varepsilon$. }
\end{figure}
\begin{figure}[H]
  \centering
  % \psfrag{2}{\tiny Two Step}
  \includegraphics[width=0.5\textwidth,  trim={0.45cm 0.4cm 0 0},clip]{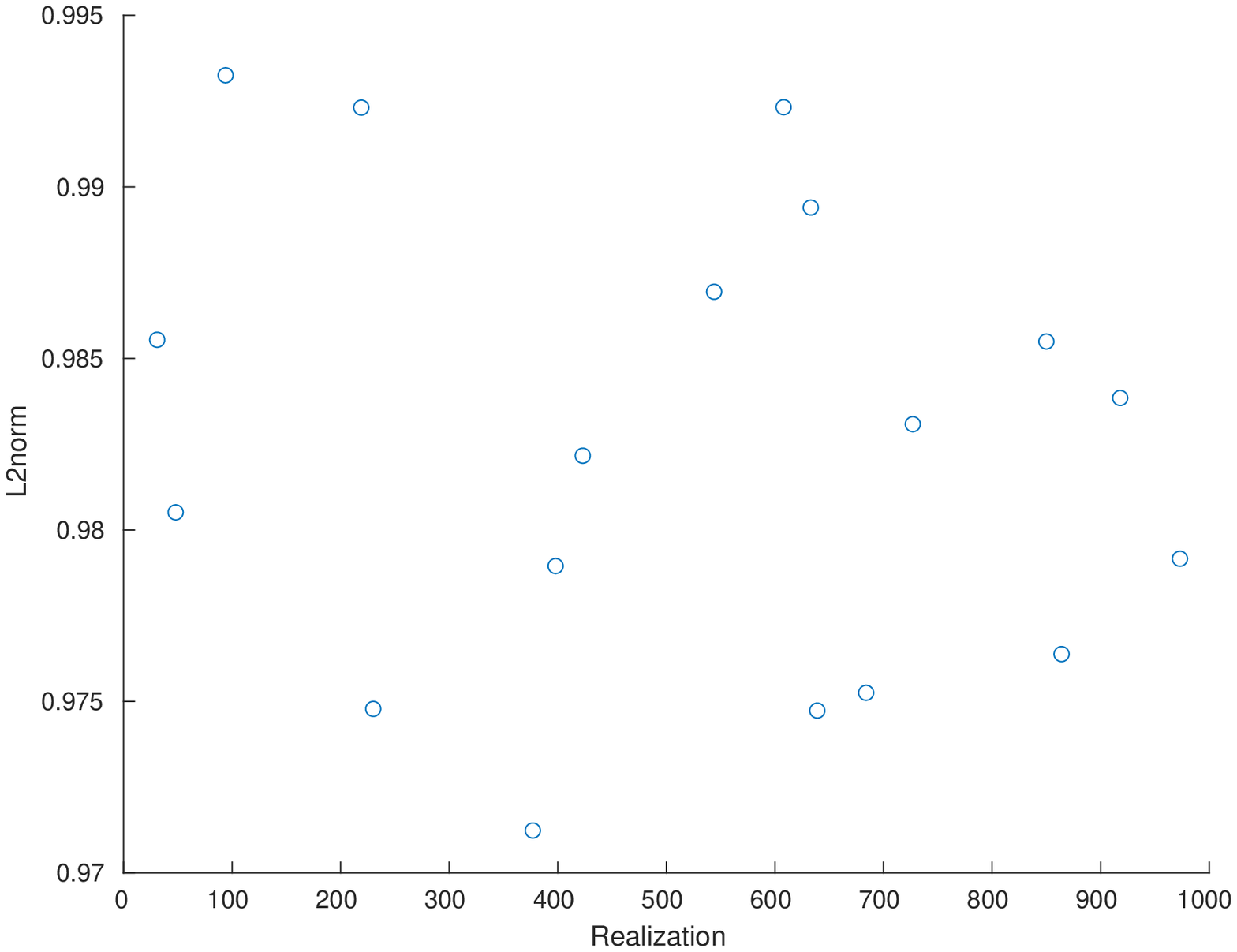}
  \caption{{\tiny} \space  For $N=25, \beta =3, h=0.5$ we simulated 1000 realizations of the disorder and 500 start values used for all of them. For $\varepsilon = -0.705 - 0.001 \cdot k, k=0,\ldots, 20$ we made up to 1.000 iterations of the Banach-$\varepsilon$-algorithm where we stopped when we hit the target level of stability (measured with respect to the {\it maximum absolute error}) of $10^{-7}$. For all realizations we plotted the {\it mean squared errors} found for the concerning realization.}
\end{figure}

\subsection{1000 realizations of disorder and 500 start values and 21 positive $\varepsilon$-values}\leavevmode
The following plots concern the same iterations with the same realizations, same start values.  We only accepted solutions satisfying the Plefka condition and being inside of the hypercube. 
\begin{figure}[H]
  \centering
  % \psfrag{2}{\tiny Two Step}
  \includegraphics[width=0.5\textwidth,  trim={0.45cm 0.4cm 0 0},clip]{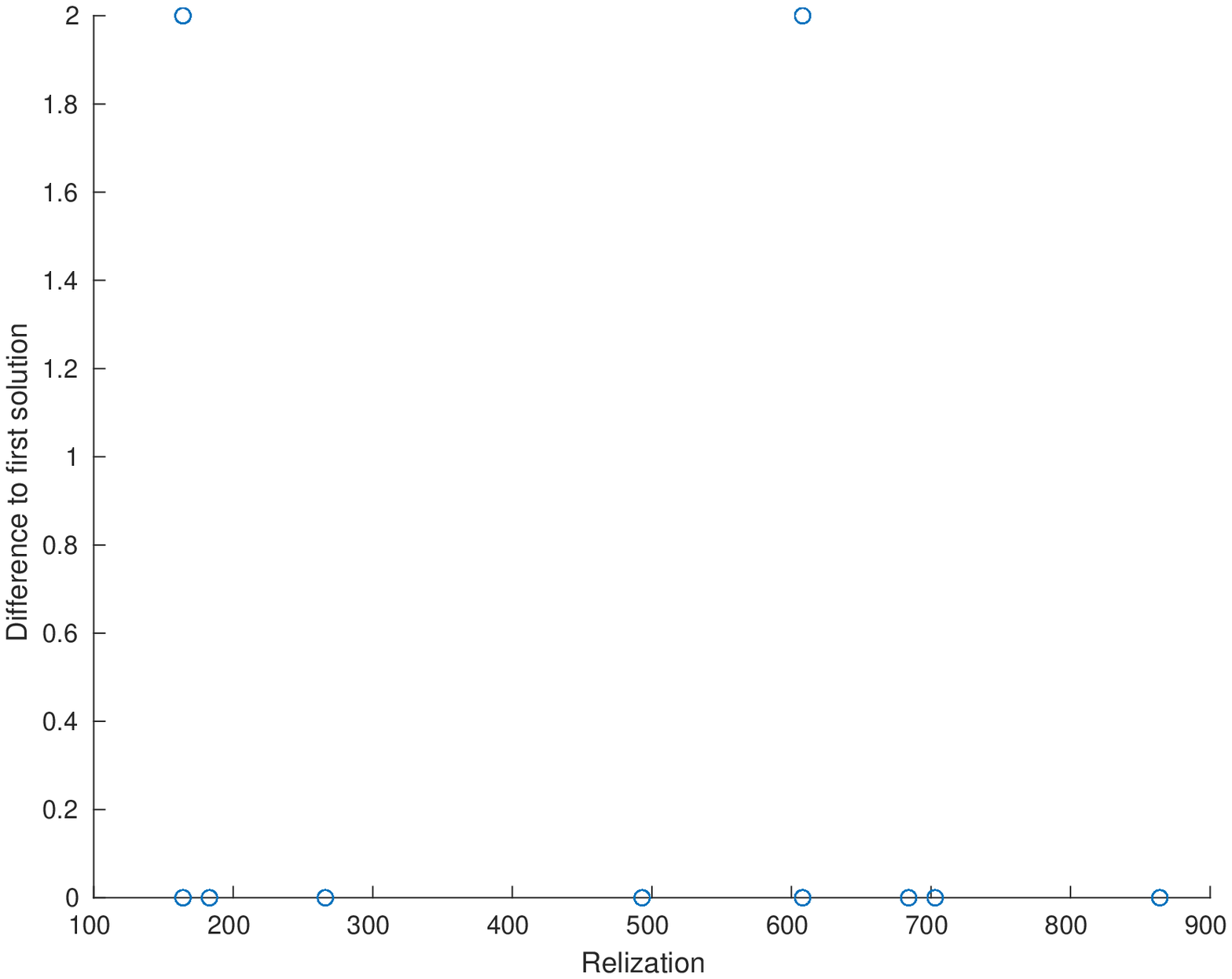}
  \caption{{\tiny} \space  For $N=25, \beta =3, h=0.5$ we simulated 1000 realizations of the disorder and 500 start values used for all of them. For $\varepsilon = 0.705 + 0.001 \cdot k, k=0,\ldots, 20$ we made up to 1.000 iterations of the Banach-$\varepsilon$-algorithm where we stopped when we hit the target level of stability (measured with respect to the {\it maximum absolute error}) of $10^{-7}$.  For all realizations we plotted the {\it maximum absolute error} to the first solution of the realization we found.}
\end{figure}
\begin{figure}[H]
  \centering
  % \psfrag{2}{\tiny Two Step}
  \includegraphics[width=0.5\textwidth,  trim={0.45cm 0.4cm 0 0},clip]{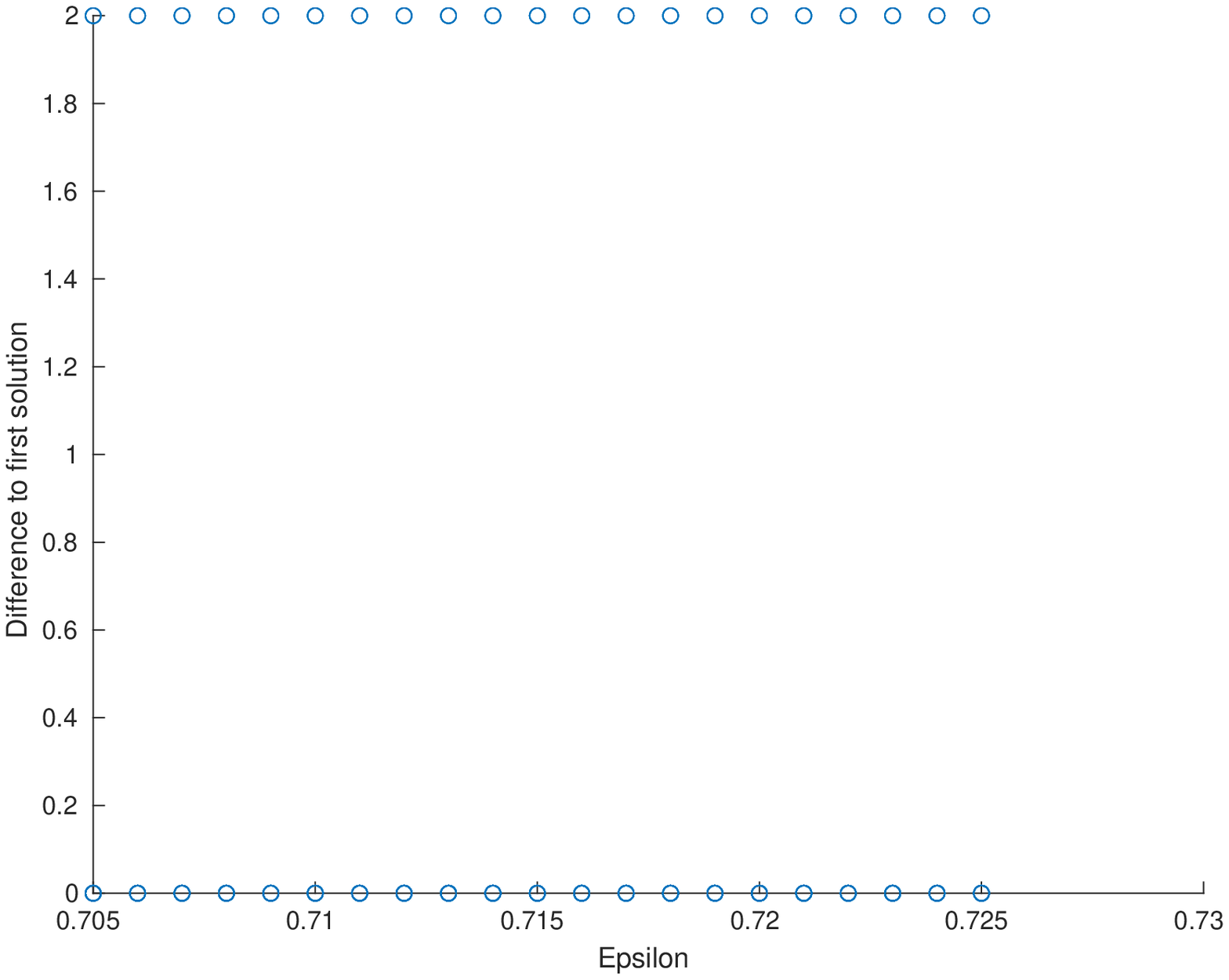}
  \caption{{\tiny} \space  For $N=25, \beta =3, h=0.5$ we simulated 1000 realizations of the disorder and 500 start values used for all of them. For $\varepsilon = 0.705 + 0.001 \cdot k, k=0,\ldots, 20$ we made up to 1.000 iterations of the Banach-$\varepsilon$-algorithm where we stopped when we hit the target level of stability (measured with respect to the {\it maximum absolute error}) of $10^{-7}$.  For all $\varepsilon$ we plotted the {\it maximum absolute error} to the first solution of the realization we found.}
\end{figure}
\begin{figure}[H]
  \centering
  % \psfrag{2}{\tiny Two Step}
  \includegraphics[width=0.5\textwidth,  trim={0.45cm 0.4cm 0 0},clip]{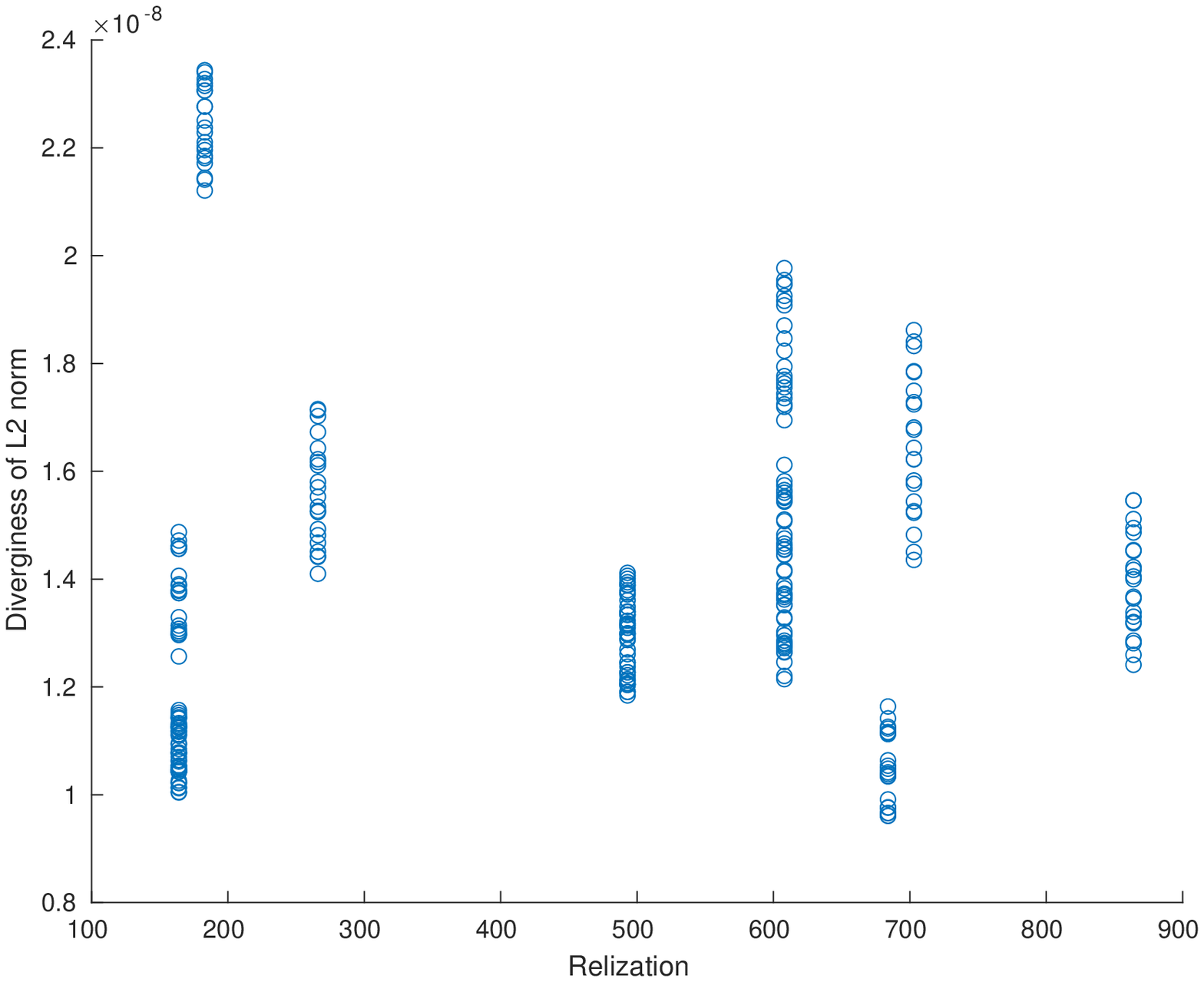}
  \caption{{\tiny} \space  For $N=25, \beta =3, h=0.5$ we simulated 1000 realizations of the disorder and 500 start values used for all of them. For $\varepsilon = 0.705 + 0.001 \cdot k, k=0,\ldots, 20$ we made up to 1.000 iterations of the Banach-$\varepsilon$-algorithm where we stopped when we hit the target level of stability (measured with respect to the {\it maximum absolute error}) of $10^{-7}$.  For all realizations we plotted the difference of the {\it mean squared errors} of the last two iterates of the solutions found in the respective realization.}
\end{figure}
\begin{figure}[H]
  \centering
  % \psfrag{2}{\tiny Two Step}
  \includegraphics[width=0.5\textwidth,  trim={0.45cm 0.4cm 0 0},clip]{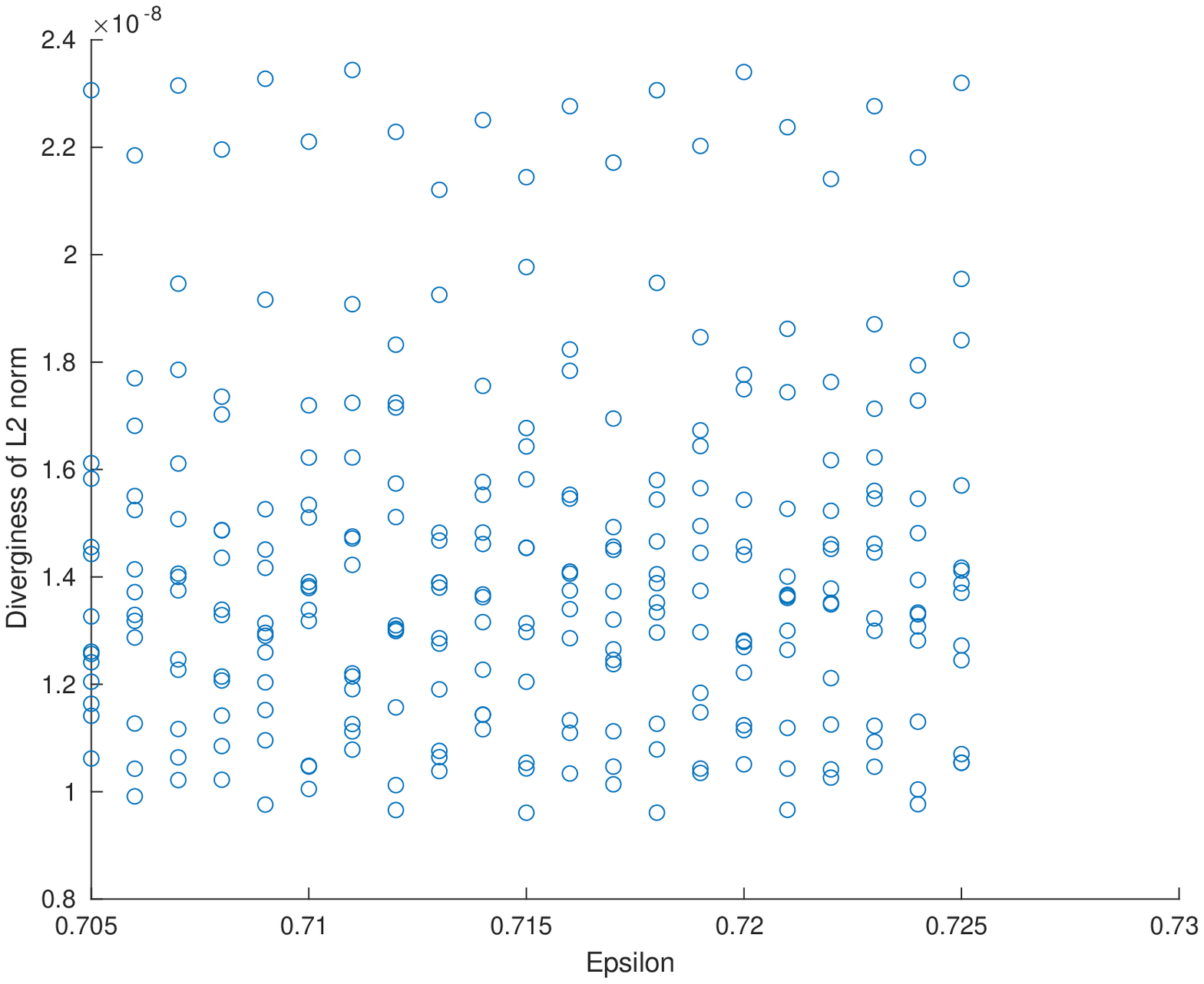}
  \caption{{\tiny} \space  For $N=25, \beta =3, h=0.5$ we simulated 1000 realizations of the disorder and 500 start values used for all of them. For $\varepsilon = 0.705 + 0.001 \cdot k, k=0,\ldots, 20$ we made up to 1.000 iterations of the Banach-$\varepsilon$-algorithm where we stopped when we hit the target level of stability (measured with respect to the {\it maximum absolute error}) of $10^{-7}$.  For all $\varepsilon$ we plotted the difference of the {\it mean squared errors} of the last two iterates of the solutions found with the respective $\varepsilon$.}
\end{figure}
\begin{figure}[H]
  \centering
  % \psfrag{2}{\tiny Two Step}
  \includegraphics[width=0.5\textwidth,  trim={0.45cm 0.4cm 0 0},clip]{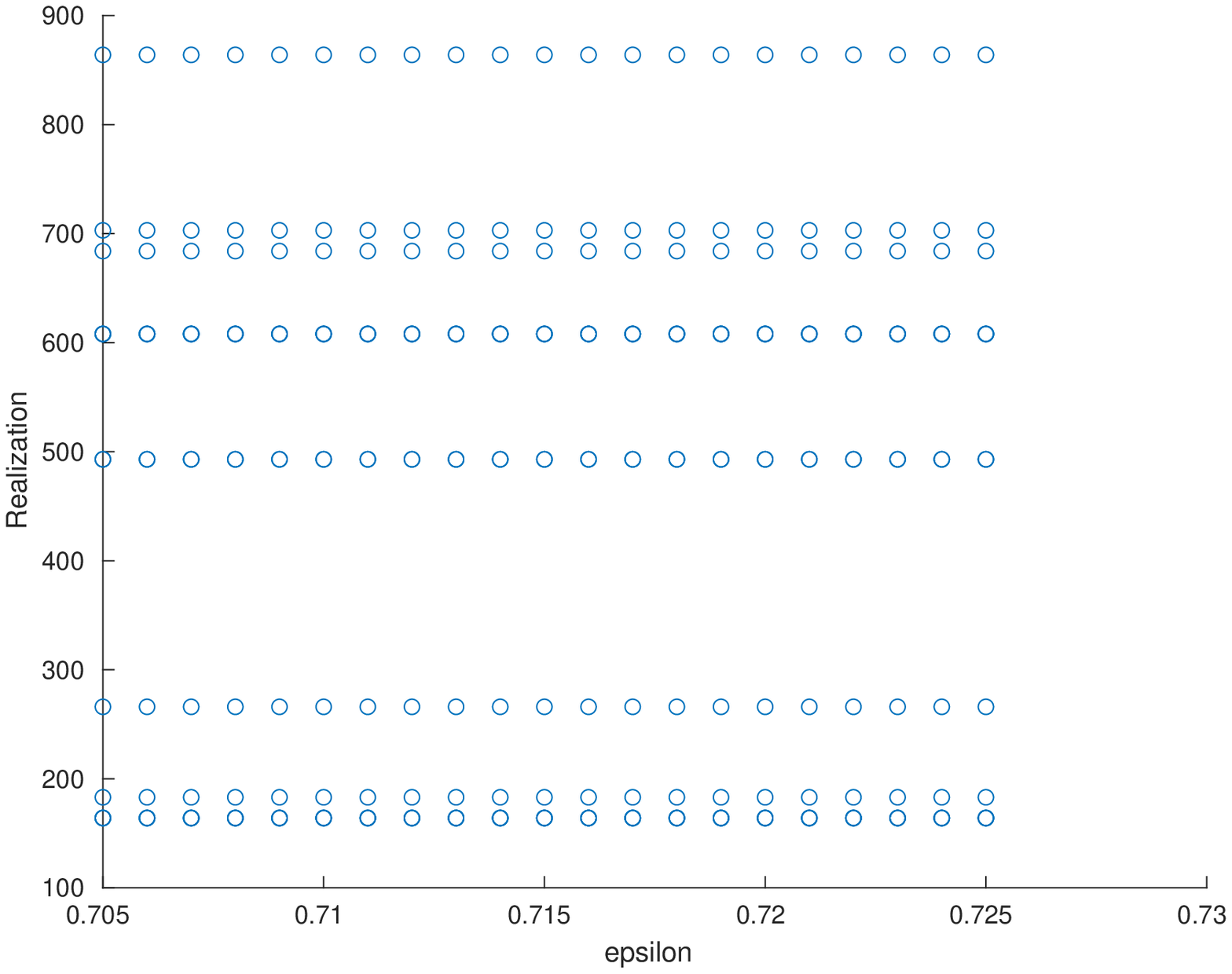}
  \caption{{\tiny} \space  For $N=25, \beta =3, h=0.5$ we simulated 1000 realizations of the disorder and 500 start values used for all of them. For $\varepsilon = 0.705 + 0.001 \cdot k, k=0,\ldots, 20$ we made up to 1.000 iterations of the Banach-$\varepsilon$-algorithm where we stopped when we hit the target level of stability (measured with respect to the {\it maximum absolute error}) of $10^{-7}$.  For all solutions we found we plotted a point with $x$ coordinate specifying the used $\varepsilon$ and $y$ coordinate the concerning realization}
\end{figure}
\begin{figure}[H]
  \centering
  % \psfrag{2}{\tiny Two Step}
  \includegraphics[width=0.5\textwidth,  trim={0.45cm 0.4cm 0 0},clip]{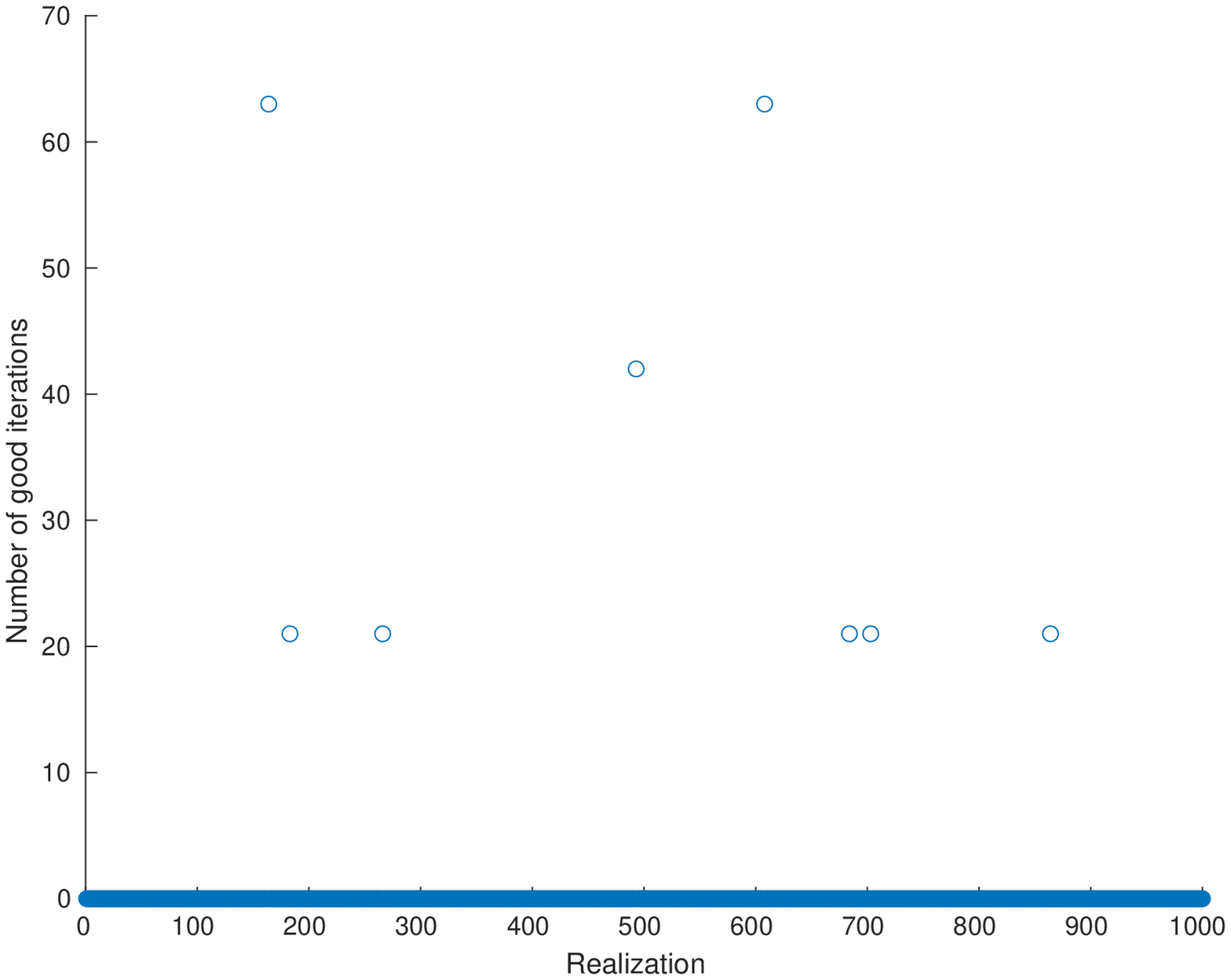}
  \caption{{\tiny} \space  For $N=25, \beta =3, h=0.5$ we simulated 1000 realizations of the disorder and 500 start values used for all of them. For $\varepsilon = 0.705 + 0.001 \cdot k, k=0,\ldots, 20$ we made up to 1.000 iterations of the Banach-$\varepsilon$-algorithm where we stopped when we hit the target level of stability (measured with respect to the {\it maximum absolute error}) of $10^{-7}$.  For all realizations we plotted the number of good solutions. }
\end{figure}
\begin{figure}[H]
  \centering
  % \psfrag{2}{\tiny Two Step}
  \includegraphics[width=0.5\textwidth,  trim={0.45cm 0.4cm 0 0},clip]{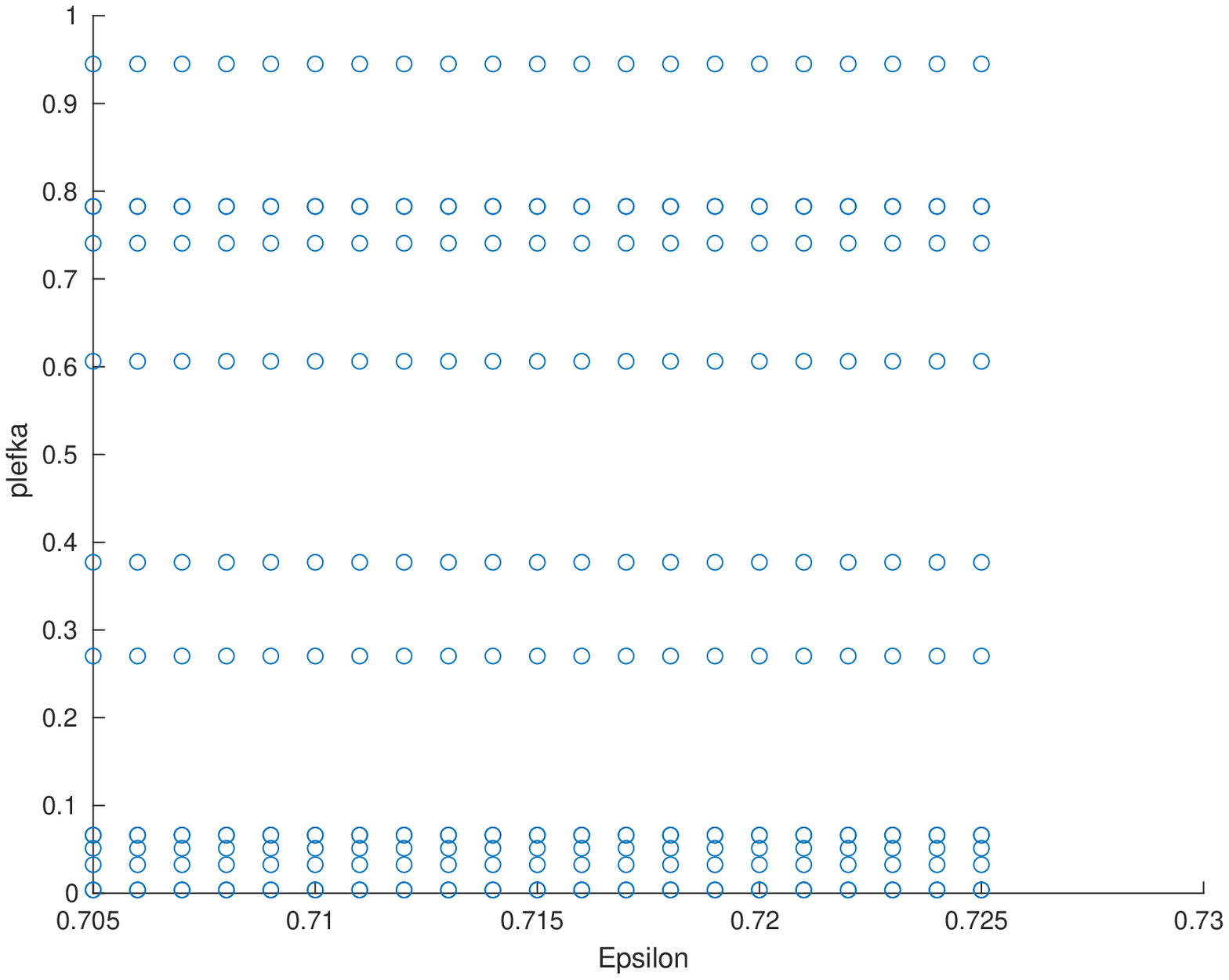}
  \caption{{\tiny} \space  For $N=25, \beta =3, h=0.5$ we simulated 1000 realizations of the disorder and 500 start values used for all of them. For $\varepsilon = 0.705 + 0.001 \cdot k, k=0,\ldots, 20$ we made up to 1.000 iterations of the Banach-$\varepsilon$-algorithm where we stopped when we hit the target level of stability (measured with respect to the {\it maximum absolute error}) of $10^{-7}$.  For all $\varepsilon$ we plotted the Plefka values found with this $\varepsilon$. }
\end{figure}
\begin{figure}[H]
  \centering
  % \psfrag{2}{\tiny Two Step}
  \includegraphics[width=0.5\textwidth,  trim={0.45cm 0.4cm 0 0},clip]{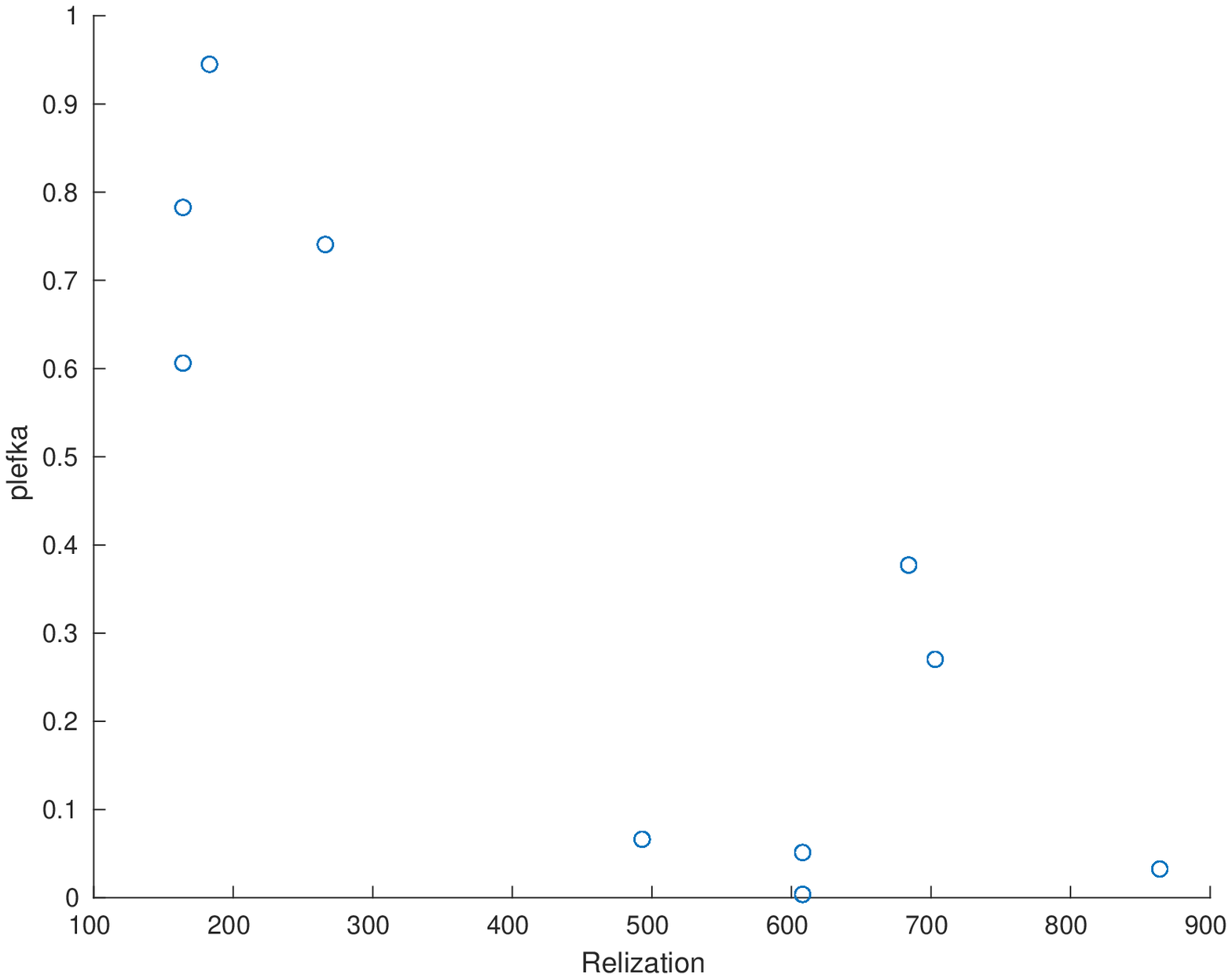}
  \caption{{\tiny} \space  For $N=25, \beta =3, h=0.5$ we simulated 1000 realizations of the disorder and 500 start values used for all of them. For $\varepsilon = 0.705 + 0.001 \cdot k, k=0,\ldots, 20$ we made up to 1.000 iterations of the Banach-$\varepsilon$-algorithm where we stopped when we hit the target level of stability (measured with respect to the {\it maximum absolute error}) of $10^{-7}$.  For all realizations we plotted the Plefka values found in the concerning realization. }
\end{figure}
\begin{figure}[H]
  \centering
  % \psfrag{2}{\tiny Two Step}
  \includegraphics[width=0.5\textwidth,  trim={0.45cm 0.4cm 0 0},clip]{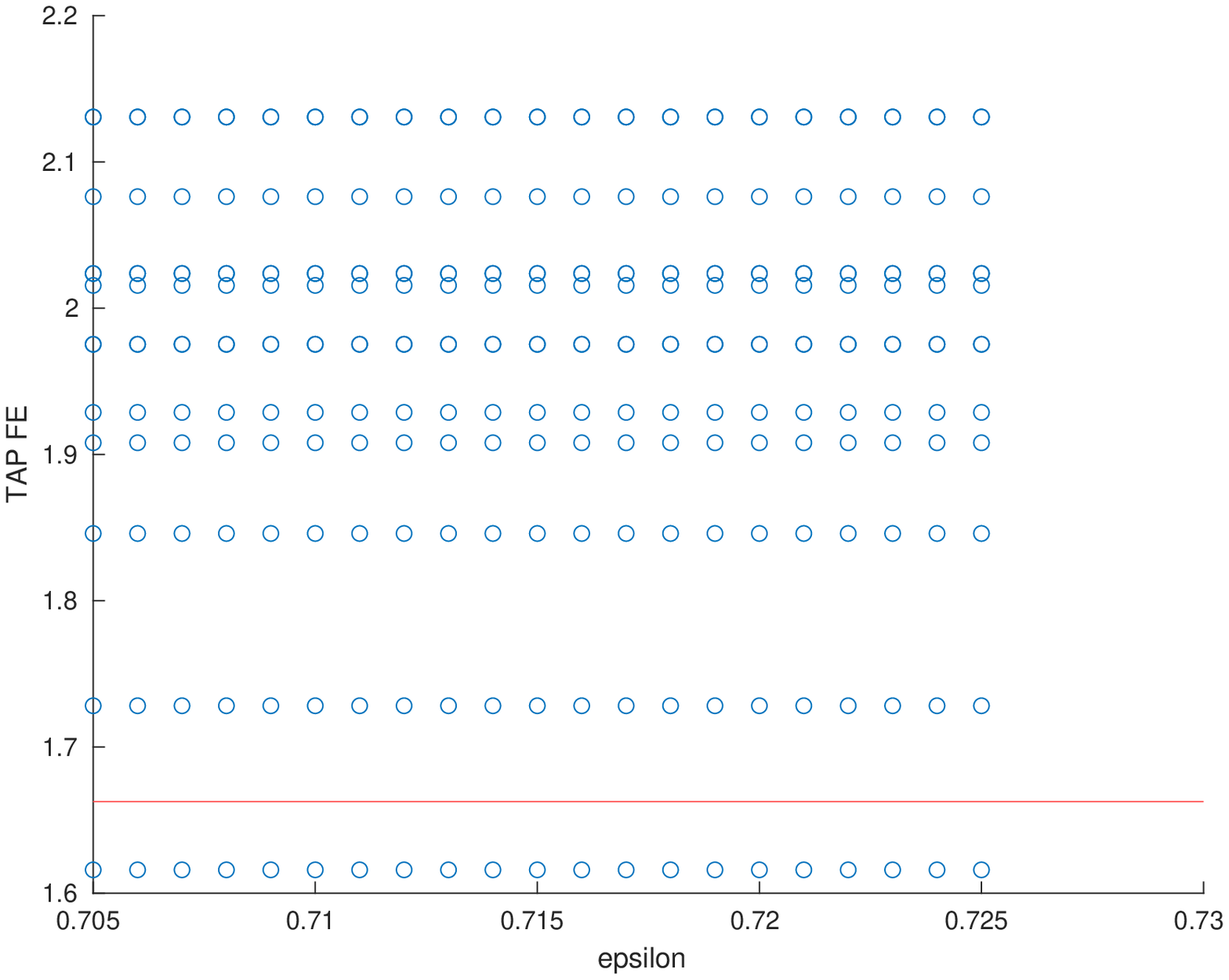}
  \caption{{\tiny} \space  For $N=25, \beta =3, h=0.5$ we simulated 1000 realizations of the disorder and 500 start values used for all of them. For $\varepsilon = 0.705 + 0.001 \cdot k, k=0,\ldots, 20$ we made up to 1.000 iterations of the Banach-$\varepsilon$-algorithm where we stopped when we hit the target level of stability (measured with respect to the {\it maximum absolute error}) of $10^{-7}$.  For all $\varepsilon$ we plotted the TAP free energy values found with this $\varepsilon$. }
\end{figure}
\begin{figure}[H]
  \centering
  % \psfrag{2}{\tiny Two Step}
  \includegraphics[width=0.5\textwidth,  trim={0.45cm 0.4cm 0 0},clip]{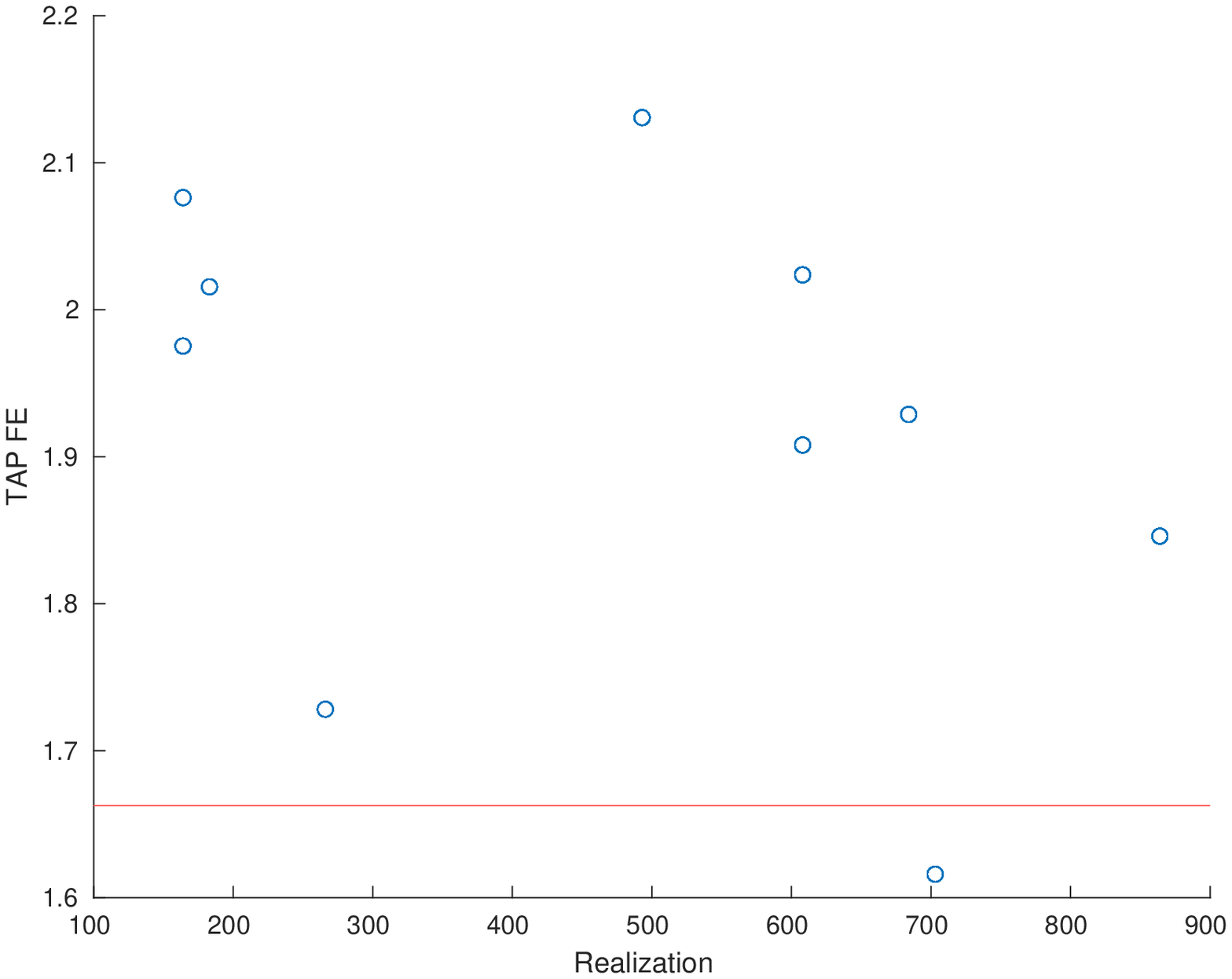}
  \caption{{\tiny} \space  For $N=25, \beta =3, h=0.5$ we simulated 1000 realizations of the disorder and 500 start values used for all of them. For $\varepsilon = 0.705 + 0.001 \cdot k, k=0,\ldots, 20$ we made up to 1.000 iterations of the Banach-$\varepsilon$-algorithm where we stopped when we hit the target level of stability (measured with respect to the {\it maximum absolute error}) of $10^{-7}$.  For all realizations we plotted the TAP free energy values found for the concerning realization.}
\end{figure}
\begin{figure}[H]
  \centering
  % \psfrag{2}{\tiny Two Step}
  \includegraphics[width=0.5\textwidth,  trim={0cm 0.4cm 0 0},clip]{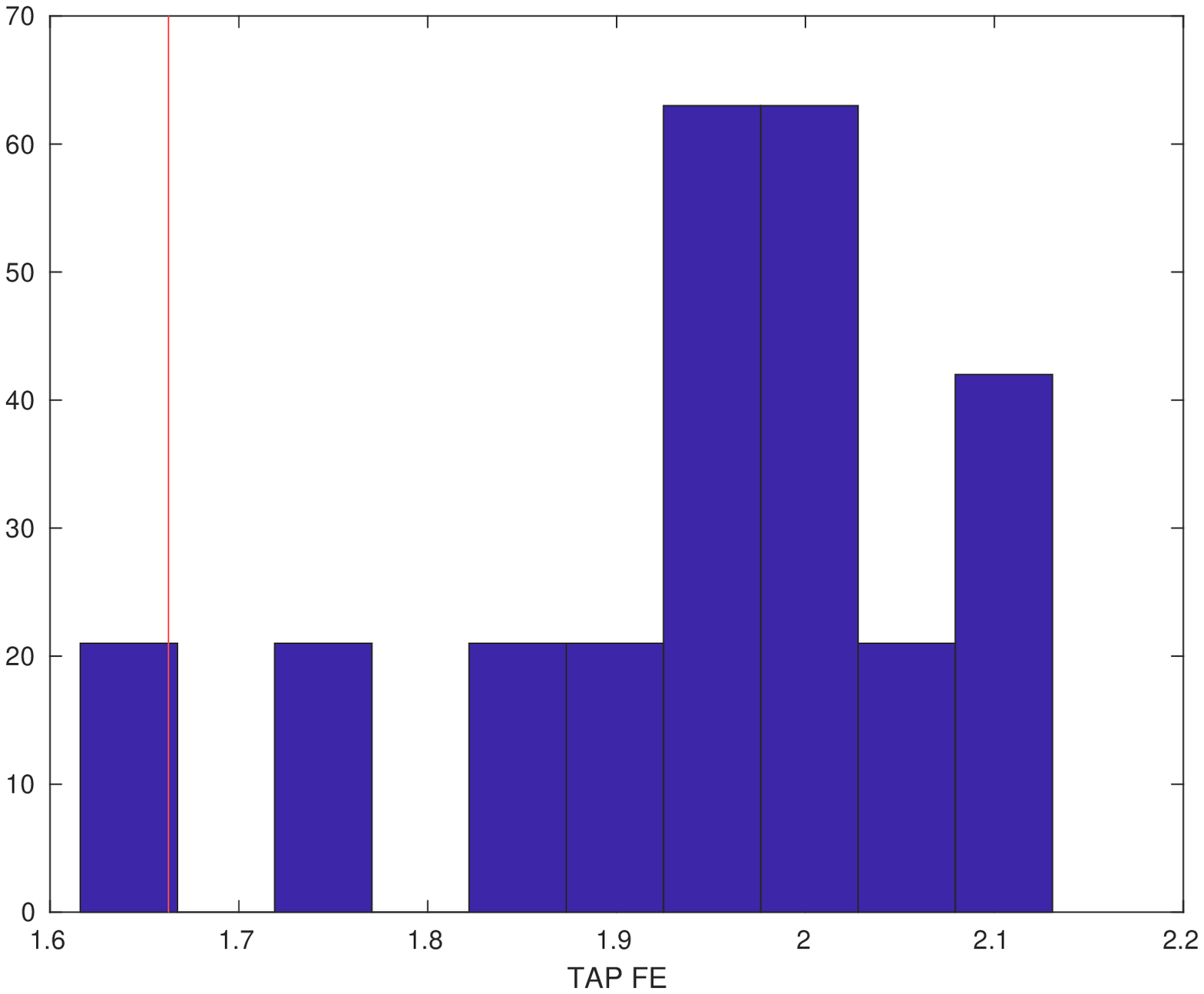}
  \caption{{\tiny} \space  For $N=25, \beta =3, h=0.5$ we simulated 1000 realizations of the disorder and 500 start values used for all of them. For $\varepsilon = 0.705 + 0.001 \cdot k, k=0,\ldots, 20$ we made up to 1.000 iterations of the Banach-$\varepsilon$-algorithm where we stopped when we hit the target level of stability (measured with respect to the {\it maximum absolute error}) of $10^{-7}$.  Just a histogram of the TAP free energy values found.}
\end{figure}
\begin{figure}[H]
  \centering
  % \psfrag{2}{\tiny Two Step}
  \includegraphics[width=0.5\textwidth,  trim={0.45cm 0.4cm 0 0},clip]{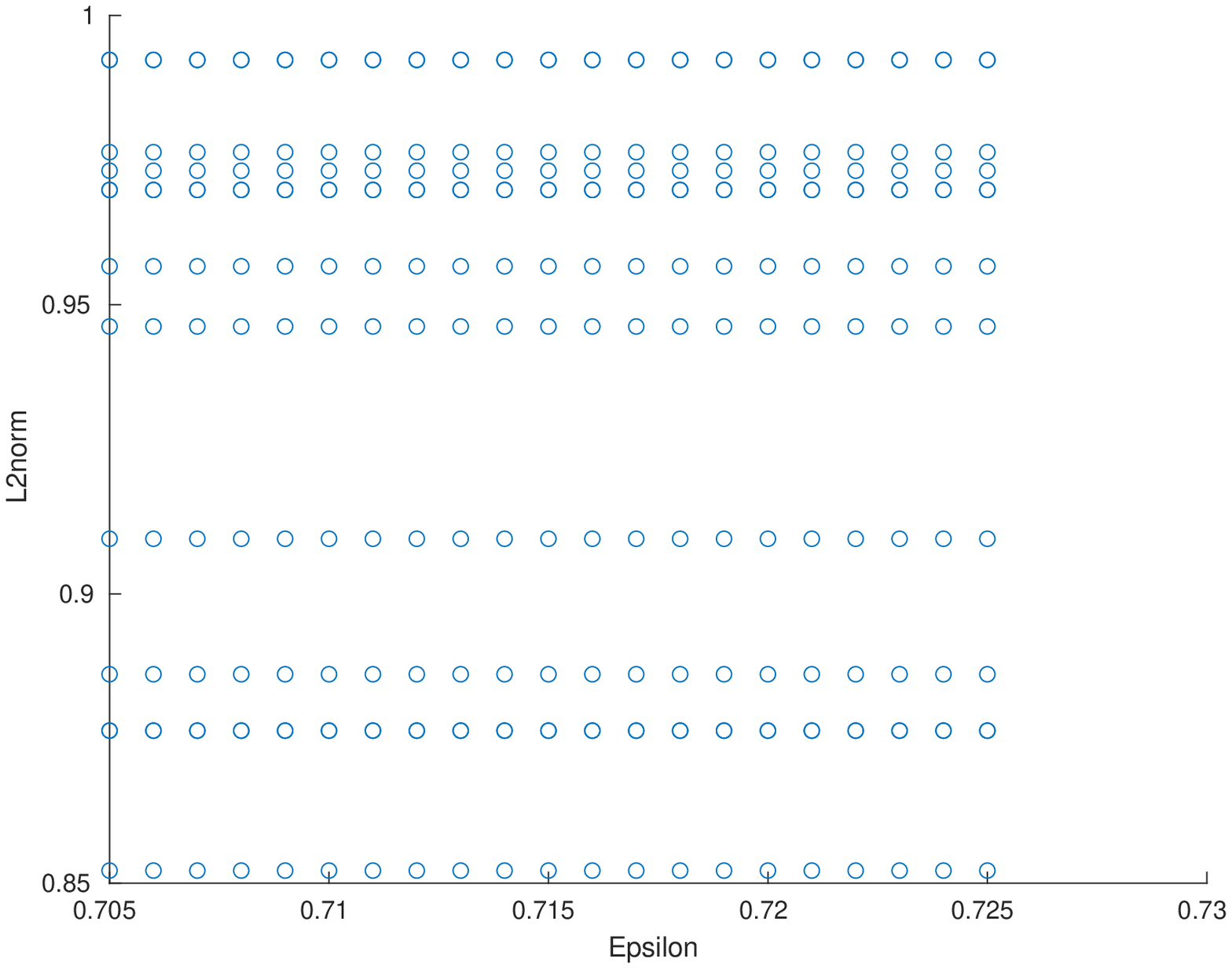}
  \caption{{\tiny} \space  For $N=25, \beta =3, h=0.5$ we simulated 1000 realizations of the disorder and 500 start values used for all of them. For $\varepsilon = 0.705 + 0.001 \cdot k, k=0,\ldots, 20$ we made up to 1.000 iterations of the Banach-$\varepsilon$-algorithm where we stopped when we hit the target level of stability (measured with respect to the {\it maximum absolute error}) of $10^{-7}$.  For all $\varepsilon$ we plotted the {\it mean squared errors} found with this $\varepsilon$. }
\end{figure}
\begin{figure}[H]
  \centering
  % \psfrag{2}{\tiny Two Step}
  \includegraphics[width=0.5\textwidth,  trim={0.45cm 0.4cm 0 0},clip]{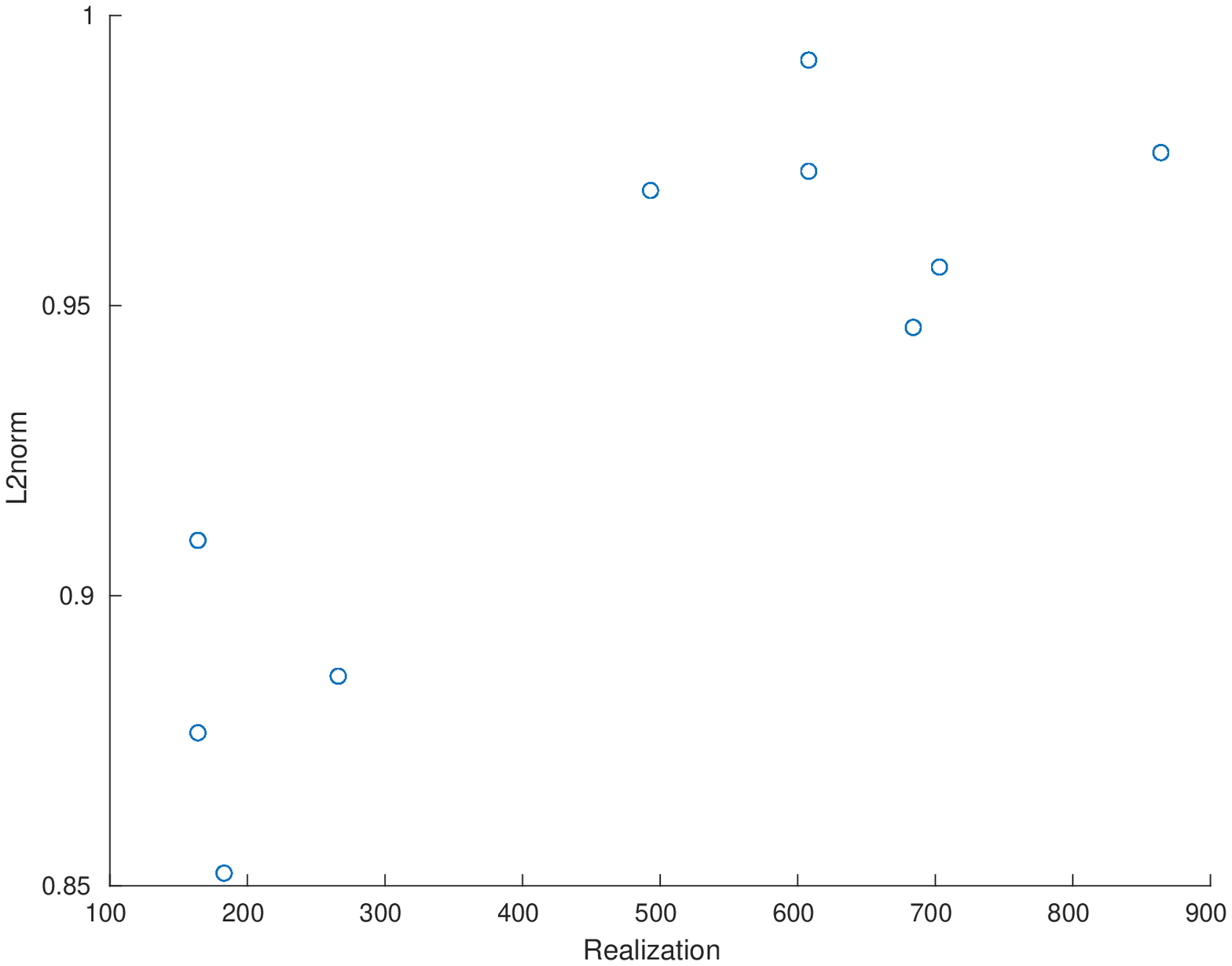}
  \caption{{\tiny} \space  For $N=25, \beta =3, h=0.5$ we simulated 1000 realizations of the disorder and 500 start values used for all of them. For $\varepsilon = 0.705 + 0.001 \cdot k, k=0,\ldots, 20$ we made up to 1.000 iterations of the Banach-$\varepsilon$-algorithm where we stopped when we hit the target level of stability (measured with respect to the {\it maximum absolute error}) of $10^{-7}$. For all realizations we plotted the {\it mean squared errors} found for the concerning realization.}
\end{figure}

\FloatBarrier\clearpage

\FloatBarrier\clearpage
\subsection{1000 realizations of disorder and 500 start values and 251 $\varepsilon$-values}\leavevmode
The following plots concern the same iterations with the same realizations, same start values.  We only accepted solutions satisfying the Plefka condition and being inside of the hypercube. 
\begin{figure}[H]
  \centering
  % \psfrag{2}{\tiny Two Step}
  \includegraphics[width=0.5\textwidth,  trim={0.45cm 0.4cm 0 0},clip]{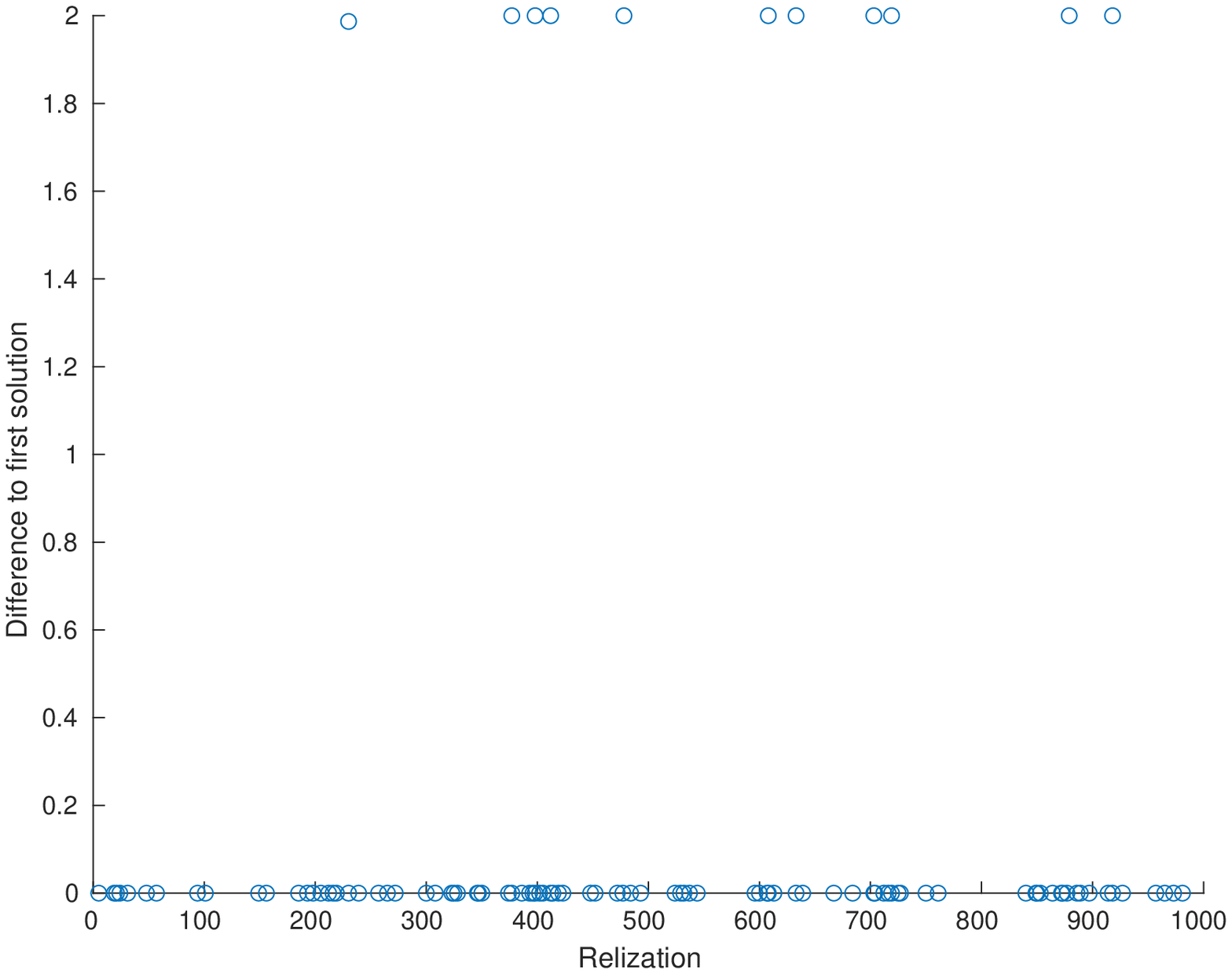}
  \caption{{\tiny} \space  For $N=25, \beta =3, h=0.5$ we simulated 1000 realizations of the disorder and 500 start values used for all of them. For $\varepsilon = -0.505 + 0.001 \cdot k, k=0,\ldots, 250$ we made up to 1.000 iterations of the Banach-$\varepsilon$-algorithm where we stopped when we hit the target level of stability (measured with respect to the {\it maximum absolute error}) of $10^{-7}$.  For all realizations we plotted the {\it maximum absolute error} to the first solution of the realization we found.}
\end{figure}
\begin{figure}[H]
  \centering
  % \psfrag{2}{\tiny Two Step}
  \includegraphics[width=0.5\textwidth,  trim={0.45cm 0.4cm 0 0},clip]{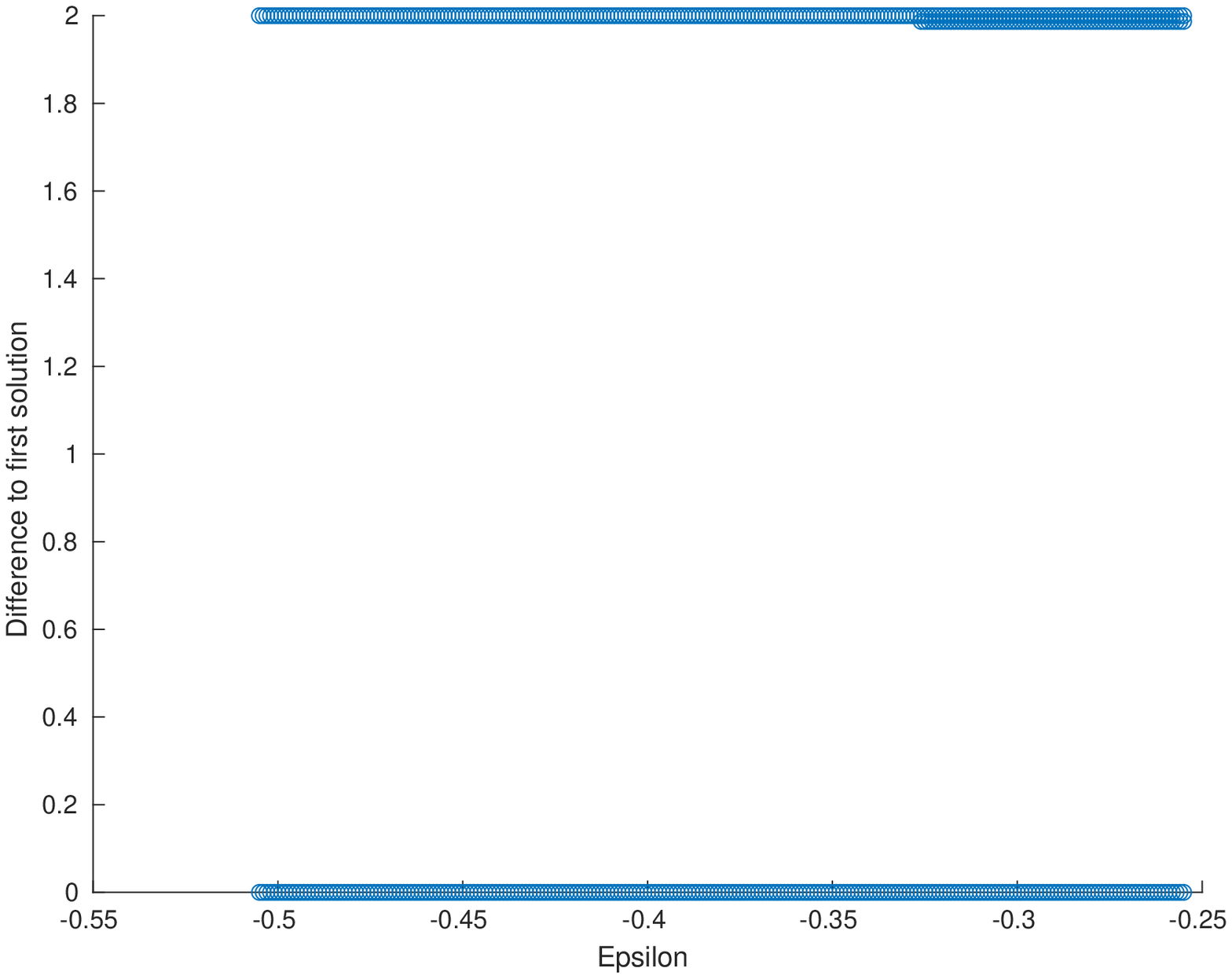}
  \caption{{\tiny} \space  For $N=25, \beta =3, h=0.5$ we simulated 1000 realizations of the disorder and 500 start values used for all of them. For $\varepsilon = -0.505 + 0.001 \cdot k, k=0,\ldots, 250$ we made up to 1.000 iterations of the Banach-$\varepsilon$-algorithm where we stopped when we hit the target level of stability (measured with respect to the {\it maximum absolute error}) of $10^{-7}$.  For all $\varepsilon$ we plotted the {\it maximum absolute error} to the first solution of the realization we found.}
\end{figure}
\begin{figure}[H]
  \centering
  % \psfrag{2}{\tiny Two Step}
  \includegraphics[width=0.5\textwidth,  trim={0.45cm 0.4cm 0 0},clip]{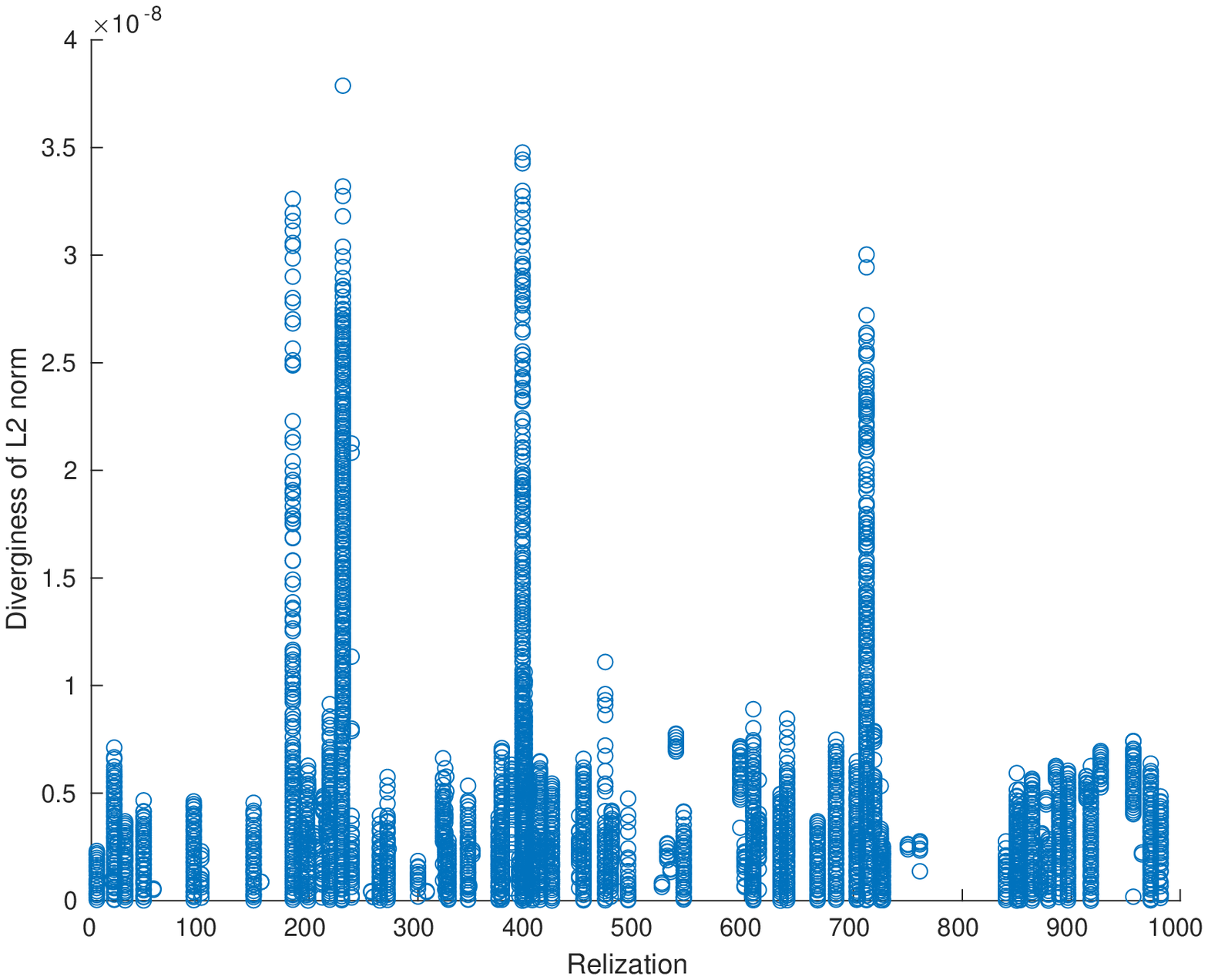}
  \caption{{\tiny} \space  For $N=25, \beta =3, h=0.5$ we simulated 1000 realizations of the disorder and 500 start values used for all of them. For $\varepsilon = -0.505 + 0.001 \cdot k, k=0,\ldots, 250$ we made up to 1.000 iterations of the Banach-$\varepsilon$-algorithm where we stopped when we hit the target level of stability (measured with respect to the {\it maximum absolute error}) of $10^{-7}$.  For all realizations we plotted the difference of the {\it mean squared errors} of the last two iterates of the solutions found in the respective realization.}
\end{figure}
\begin{figure}[H]
  \centering
  % \psfrag{2}{\tiny Two Step}
  \includegraphics[width=0.5\textwidth,  trim={0.45cm 0.4cm 0 0},clip]{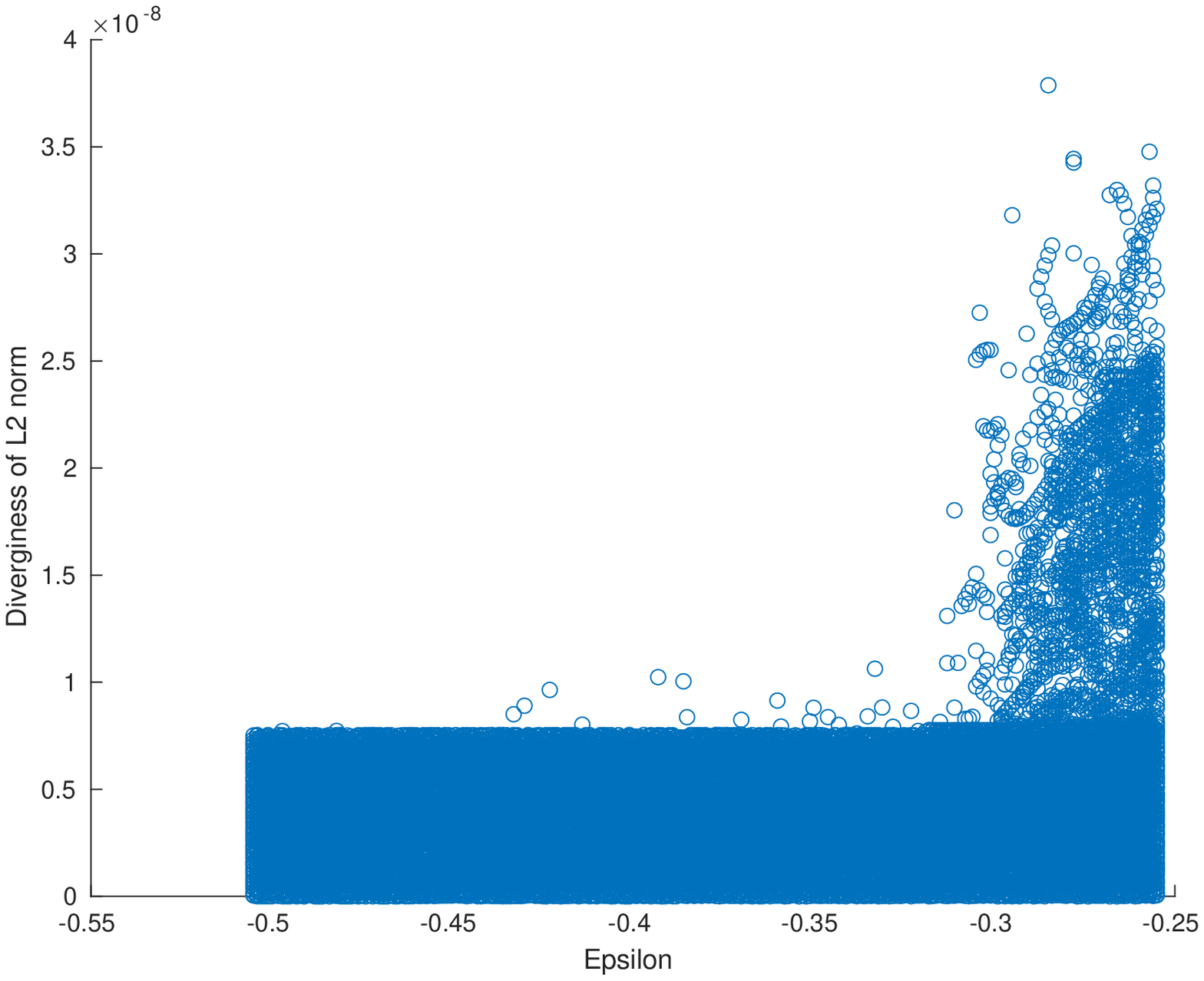}
  \caption{{\tiny} \space  For $N=25, \beta =3, h=0.5$ we simulated 1000 realizations of the disorder and 500 start values used for all of them. For $\varepsilon = -0.505 + 0.001 \cdot k, k=0,\ldots, 250$ we made up to 1.000 iterations of the Banach-$\varepsilon$-algorithm where we stopped when we hit the target level of stability (measured with respect to the {\it maximum absolute error}) of $10^{-7}$.  For all $\varepsilon$ we plotted the difference of the {\it mean squared errors} of the last two iterates of the solutions found with the respective $\varepsilon$.}
\end{figure}
\begin{figure}[H]
  \centering
  % \psfrag{2}{\tiny Two Step}
  \includegraphics[width=0.5\textwidth,  trim={0.45cm 0.4cm 0 0},clip]{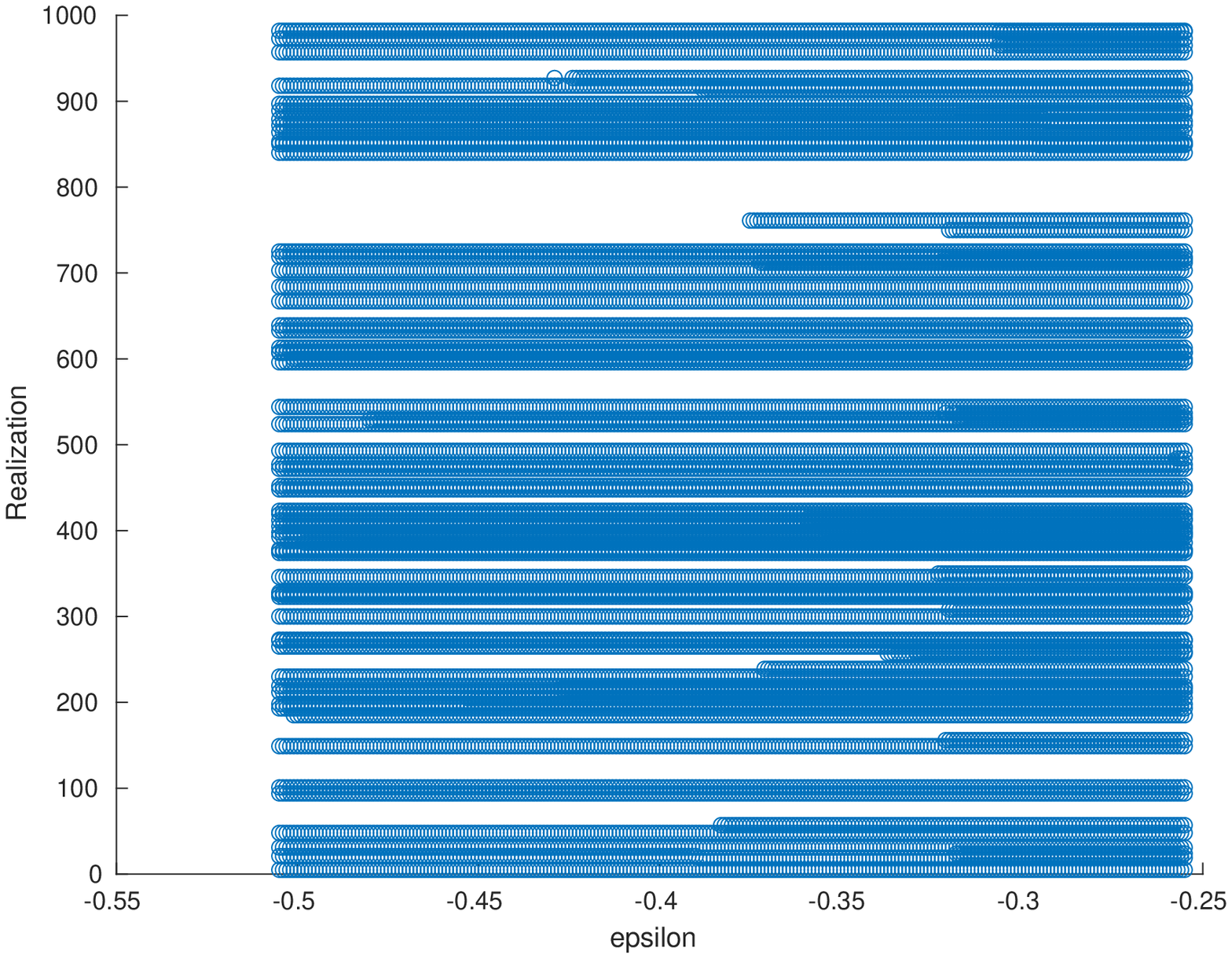}
  \caption{{\tiny} \space  For $N=25, \beta =3, h=0.5$ we simulated 1000 realizations of the disorder and 500 start values used for all of them. For $\varepsilon = -0.505 + 0.001 \cdot k, k=0,\ldots, 250$ we made up to 1.000 iterations of the Banach-$\varepsilon$-algorithm where we stopped when we hit the target level of stability (measured with respect to the {\it maximum absolute error}) of $10^{-7}$.  For all solutions we found we plotted a point with $x$ coordinate specifying the used $\varepsilon$ and $y$ coordinate the concerning realization}
\end{figure}
\begin{figure}[H]
  \centering
  % \psfrag{2}{\tiny Two Step}
  \includegraphics[width=0.5\textwidth,  trim={0.45cm 0.4cm 0 0},clip]{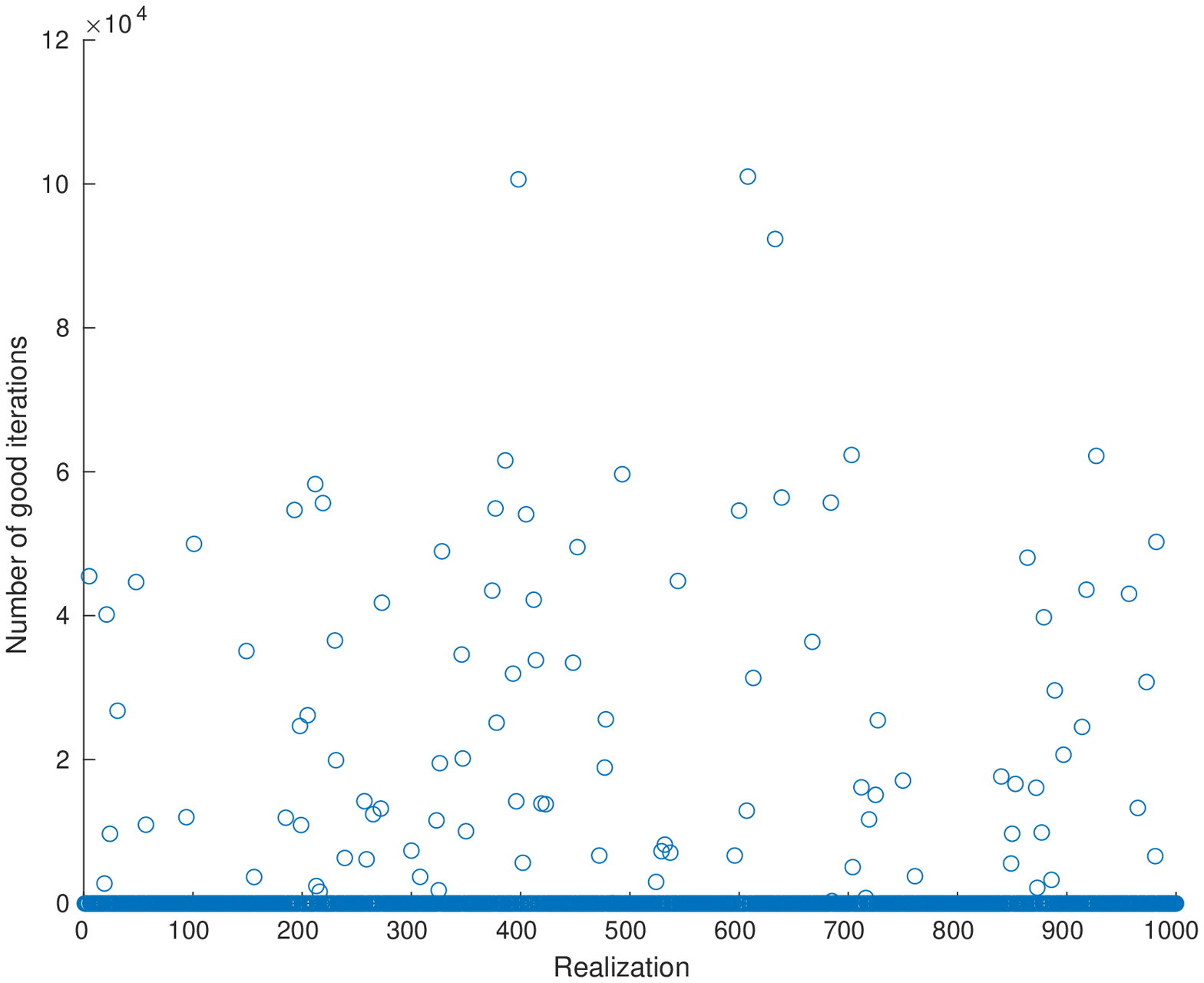}
  \caption{{\tiny} \space  For $N=25, \beta =3, h=0.5$ we simulated 1000 realizations of the disorder and 500 start values used for all of them. For $\varepsilon = -0.505 + 0.001 \cdot k, k=0,\ldots, 250$ we made up to 1.000 iterations of the Banach-$\varepsilon$-algorithm where we stopped when we hit the target level of stability (measured with respect to the {\it maximum absolute error}) of $10^{-7}$.  For all realizations we plotted the number of good solutions. }
\end{figure}
\begin{figure}[H]
  \centering
  % \psfrag{2}{\tiny Two Step}
  \includegraphics[width=0.5\textwidth,  trim={0.45cm 0.4cm 0 0},clip]{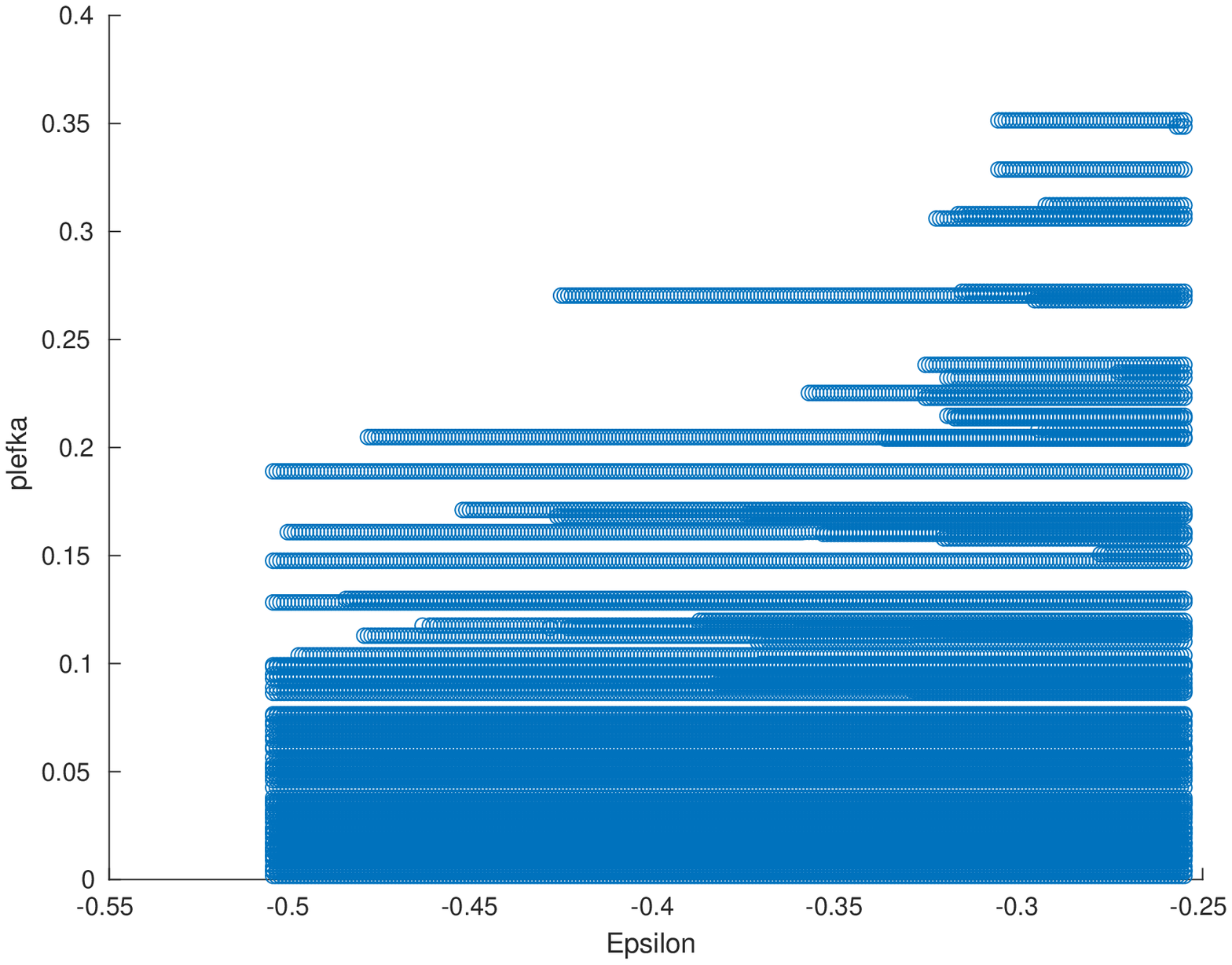}
  \caption{{\tiny} \space  For $N=25, \beta =3, h=0.5$ we simulated 1000 realizations of the disorder and 500 start values used for all of them. For $\varepsilon = -0.505 + 0.001 \cdot k, k=0,\ldots, 250$ we made up to 1.000 iterations of the Banach-$\varepsilon$-algorithm where we stopped when we hit the target level of stability (measured with respect to the {\it maximum absolute error}) of $10^{-7}$.  For all $\varepsilon$ we plotted the Plefka values found with this $\varepsilon$. }
\end{figure}
\begin{figure}[H]
  \centering
  % \psfrag{2}{\tiny Two Step}
  \includegraphics[width=0.5\textwidth,  trim={0.45cm 0.4cm 0 0},clip]{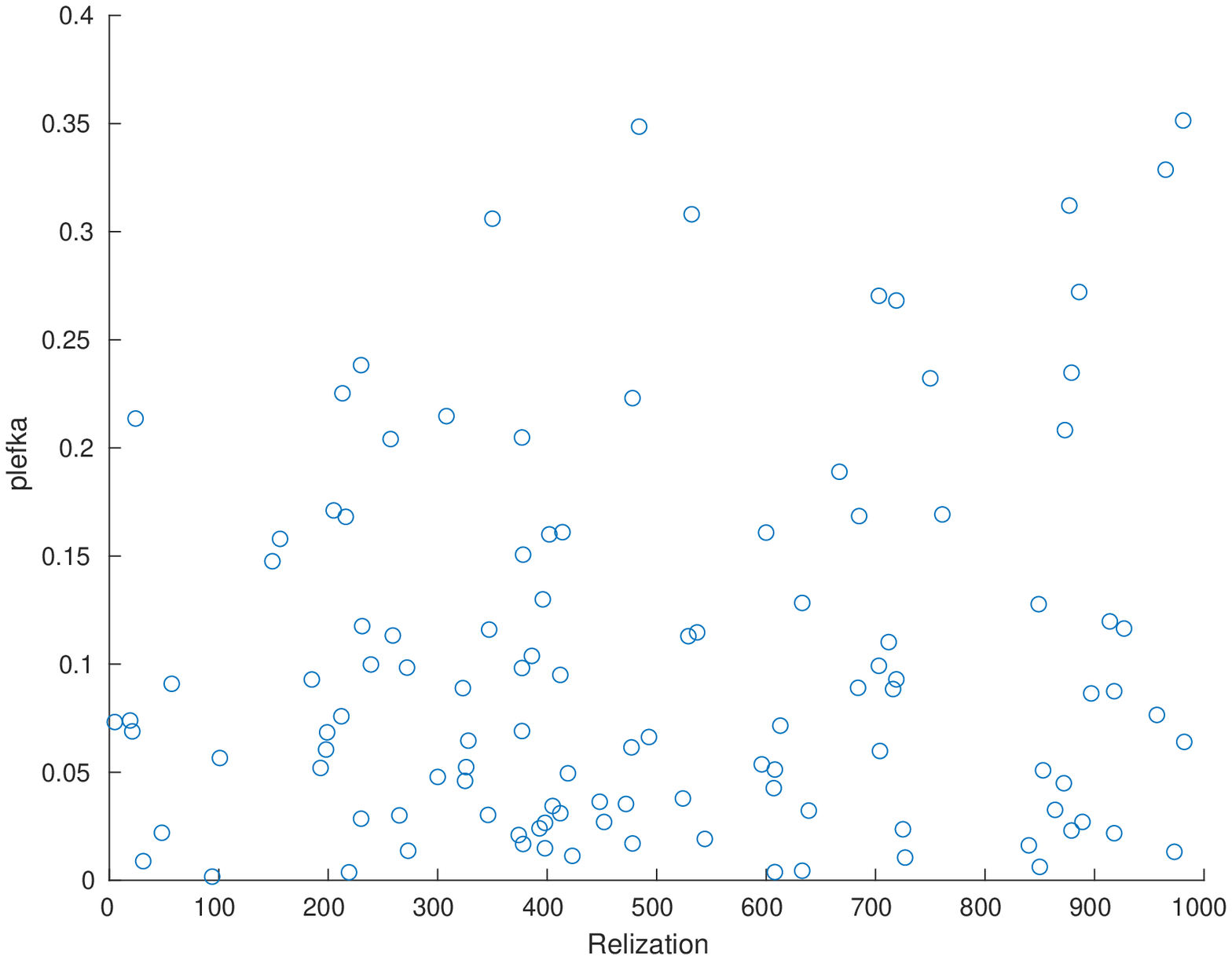}
  \caption{{\tiny} \space  For $N=25, \beta =3, h=0.5$ we simulated 1000 realizations of the disorder and 500 start values used for all of them. For $\varepsilon = -0.505 + 0.001 \cdot k, k=0,\ldots, 250$ we made up to 1.000 iterations of the Banach-$\varepsilon$-algorithm where we stopped when we hit the target level of stability (measured with respect to the {\it maximum absolute error}) of $10^{-7}$.  For all realizations we plotted the Plefka values found in the concerning realization. }
\end{figure}
\begin{figure}[H]
  \centering
  % \psfrag{2}{\tiny Two Step}
  \includegraphics[width=0.5\textwidth,  trim={0.45cm 0.4cm 0 0},clip]{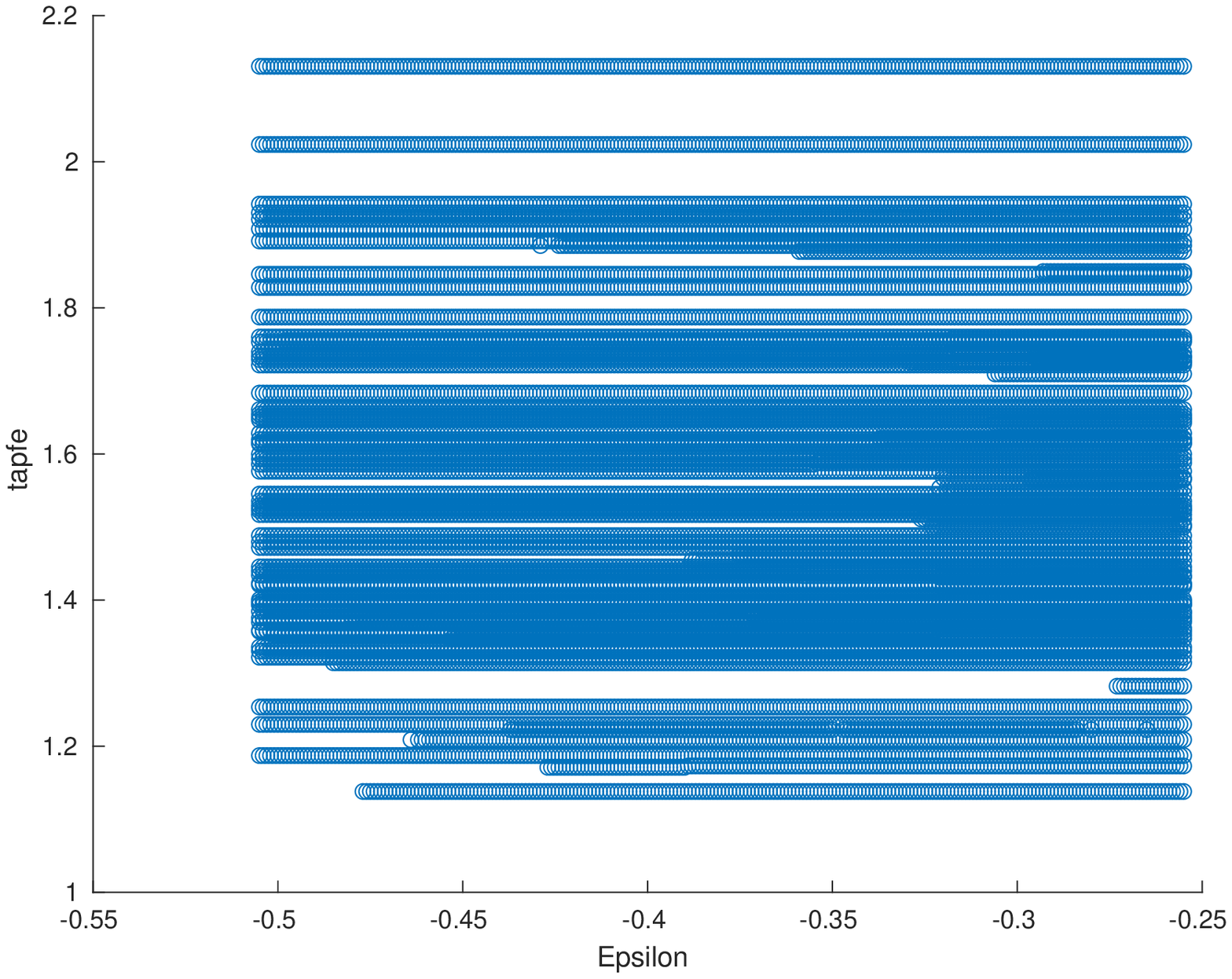}
  \caption{{\tiny} \space  For $N=25, \beta =3, h=0.5$ we simulated 1000 realizations of the disorder and 500 start values used for all of them. For $\varepsilon = -0.505 + 0.001 \cdot k, k=0,\ldots, 250$ we made up to 1.000 iterations of the Banach-$\varepsilon$-algorithm where we stopped when we hit the target level of stability (measured with respect to the {\it maximum absolute error}) of $10^{-7}$.  For all $\varepsilon$ we plotted the TAP free energy values found with this $\varepsilon$. }
\end{figure}
\begin{figure}[H]
  \centering
  % \psfrag{2}{\tiny Two Step}
  \includegraphics[width=0.5\textwidth,  trim={0.45cm 0.4cm 0 0},clip]{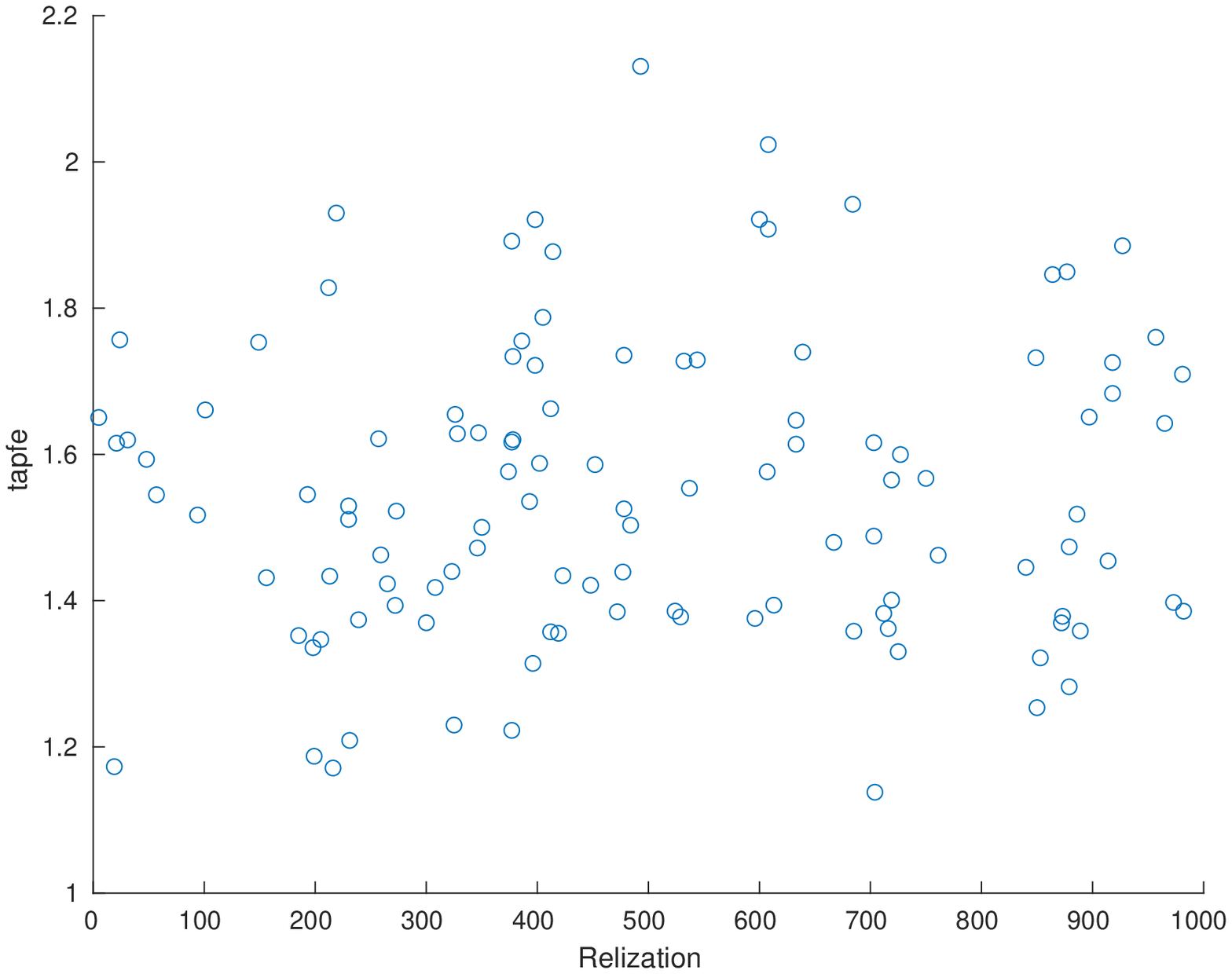}
  \caption{{\tiny} \space  For $N=25, \beta =3, h=0.5$ we simulated 1000 realizations of the disorder and 500 start values used for all of them. For $\varepsilon = -0.505 + 0.001 \cdot k, k=0,\ldots, 250$ we made up to 1.000 iterations of the Banach-$\varepsilon$-algorithm where we stopped when we hit the target level of stability (measured with respect to the {\it maximum absolute error}) of $10^{-7}$.  For all realizations we plotted the TAP free energy values found for the concerning realization.}
\end{figure}
\begin{figure}[H]
  \centering
  % \psfrag{2}{\tiny Two Step}
  \includegraphics[width=0.5\textwidth,  trim={0cm 0.4cm 0 0},clip]{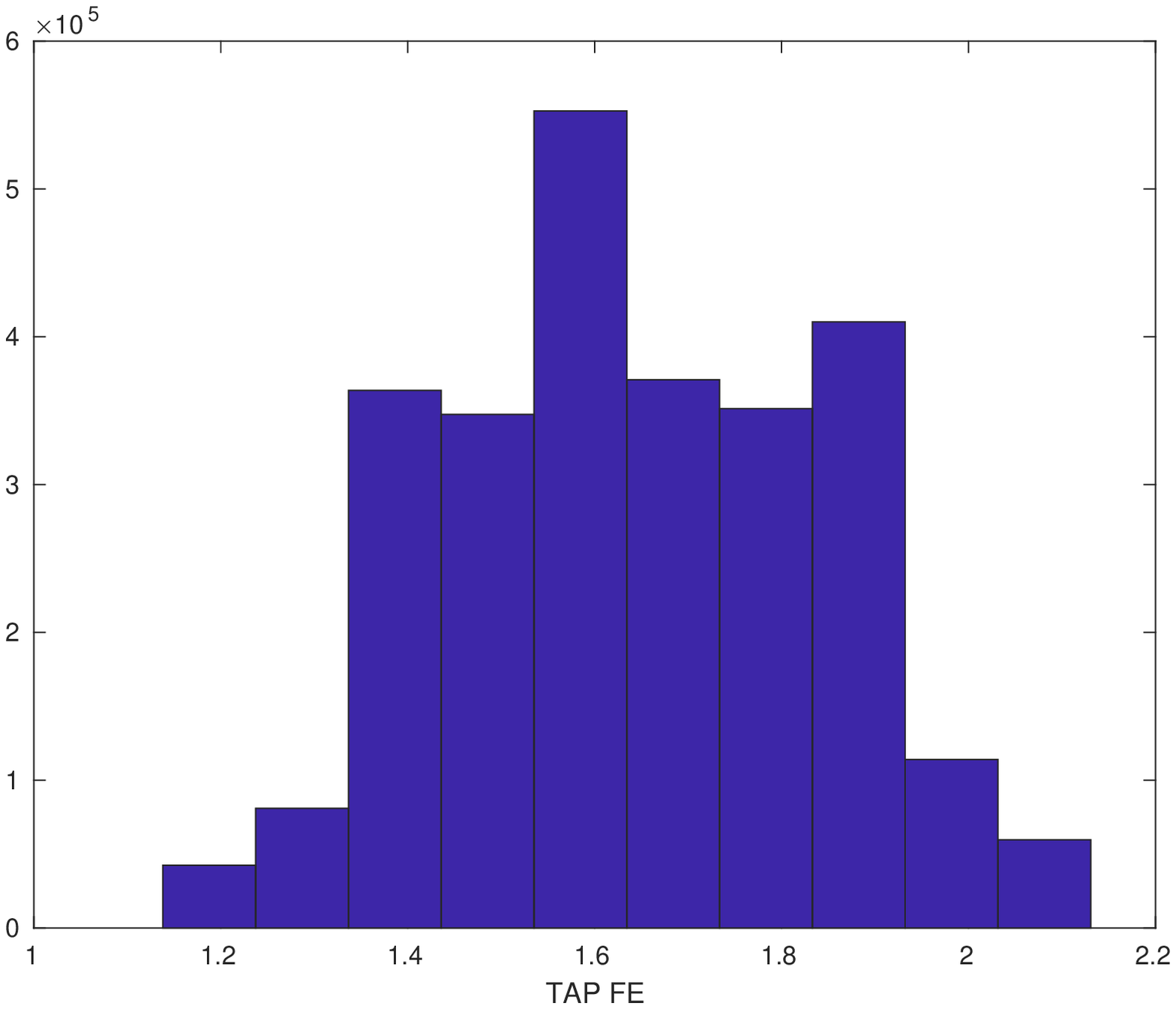}
  \caption{{\tiny} \space  For $N=25, \beta =3, h=0.5$ we simulated 1000 realizations of the disorder and 500 start values used for all of them. For $\varepsilon = -0.505 + 0.001 \cdot k, k=0,\ldots, 250$ we made up to 1.000 iterations of the Banach-$\varepsilon$-algorithm where we stopped when we hit the target level of stability (measured with respect to the {\it maximum absolute error}) of $10^{-7}$.  Just a histogram of the TAP free energy values found.}
\end{figure}
\begin{figure}[H]
  \centering
  % \psfrag{2}{\tiny Two Step}
  \includegraphics[width=0.5\textwidth,  trim={0.45cm 0.4cm 0 0},clip]{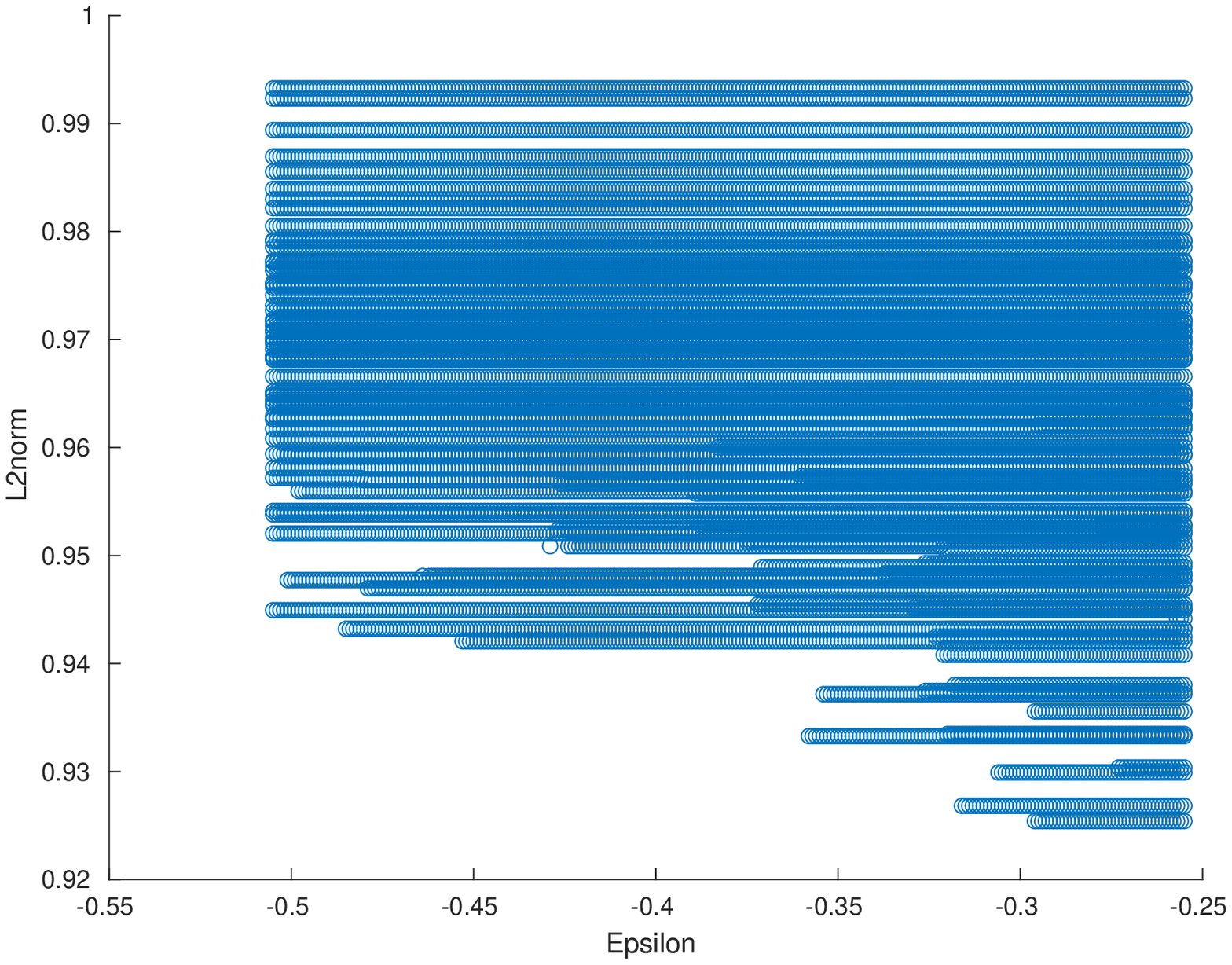}
  \caption{{\tiny} \space  For $N=25, \beta =3, h=0.5$ we simulated 1000 realizations of the disorder and 500 start values used for all of them. For $\varepsilon = -0.505 + 0.001 \cdot k, k=0,\ldots, 250$ we made up to 1.000 iterations of the Banach-$\varepsilon$-algorithm where we stopped when we hit the target level of stability (measured with respect to the {\it maximum absolute error}) of $10^{-7}$.  For all $\varepsilon$ we plotted the {\it mean squared errors} found with this $\varepsilon$. }
\end{figure}
\begin{figure}[H]
  \centering
  % \psfrag{2}{\tiny Two Step}
  \includegraphics[width=0.5\textwidth,  trim={0.45cm 0.4cm 0 0},clip]{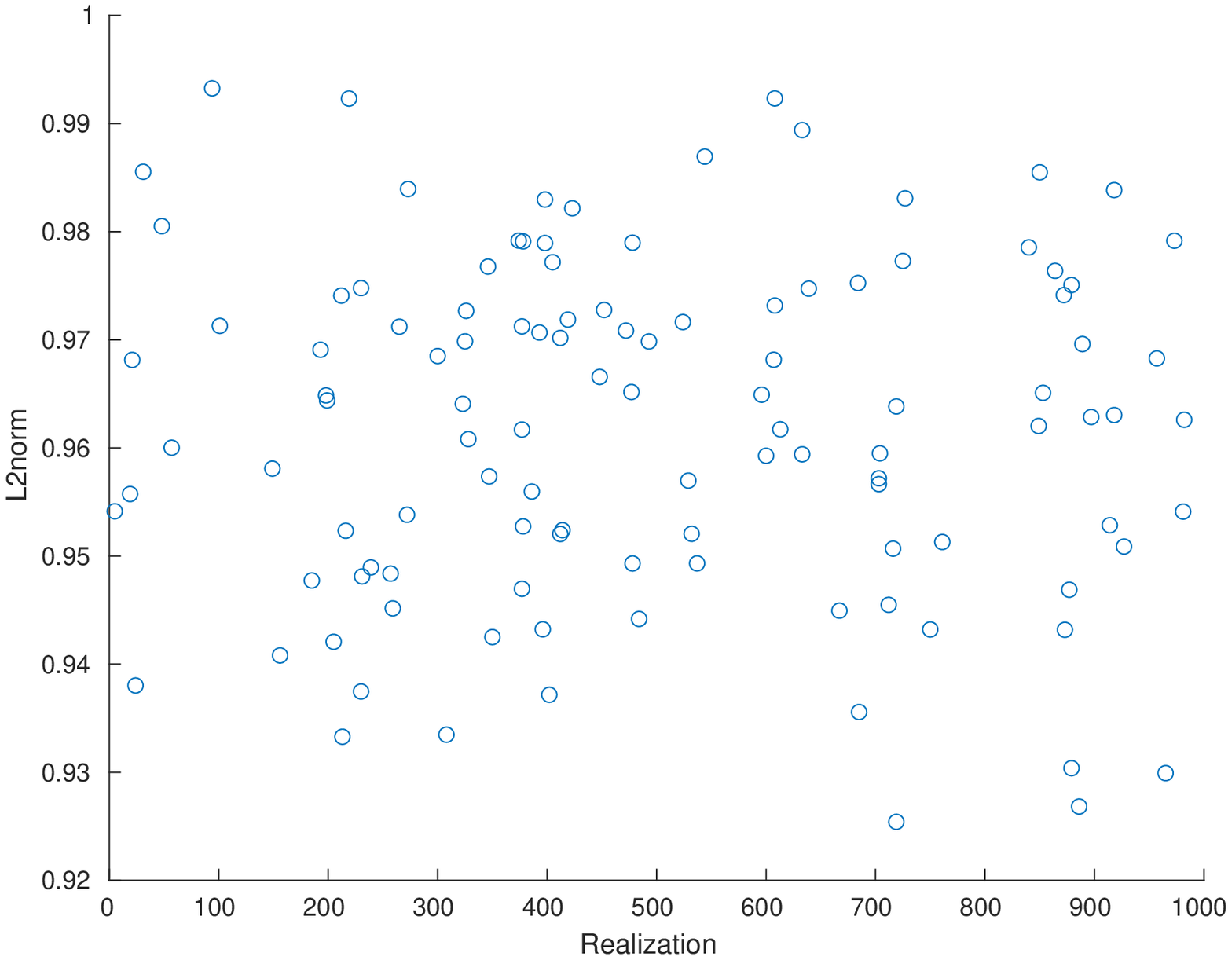}
  \caption{{\tiny} \space  For $N=25, \beta =3, h=0.5$ we simulated 1000 realizations of the disorder and 500 start values used for all of them. For $\varepsilon = -0.505 + 0.001 \cdot k, k=0,\ldots, 250$ we made up to 1.000 iterations of the Banach-$\varepsilon$-algorithm where we stopped when we hit the target level of stability (measured with respect to the {\it maximum absolute error}) of $10^{-7}$. For all realizations we plotted the {\it mean squared errors} found for the concerning realization.}
\end{figure}

\FloatBarrier\clearpage
\subsubsection{Different solutions found in 1000 realizations of disorder and with 500 start values and 251 $\varepsilon$-values}\leavevmode
We analyzed for each realization the number of differeent solutions found with the following algorithm:
\begin{algorithm}\begin{algorithmic}
  \State FilteredSolutions$ \gets$ \emph{empty array}
  \State FilteredSolutions[1]$ \gets$ FoundSolutions[$1$]
  \For {$k = 2,\ldots, $length(FoundSolutions)}
  \State $f \gets 0$
  \For {$j = 1,\ldots, $length(FilteredSolutions)}
  \If{supnorm(FoundSolutions[$k$]-FilteredSolutions[$j$]) $< 10^{-7}$}
  \State $f \gets 1$
  \EndIf
  \EndFor
  \If{$f=0$}
  \State FilteredSolutions[length(FilteredSolutions)$+1$] $ \gets$ FoundSolutions[$k$]
  \EndIf
  \EndFor
\end{algorithmic}\caption{Deciding how many solutions we found.}\label{algo:diff}\end{algorithm}
\begin{figure}[H]
  \centering
  % \psfrag{2}{\tiny Two Step}
  \includegraphics[width=0.5\textwidth,  trim={0.45cm 0.4cm 0 0},clip]{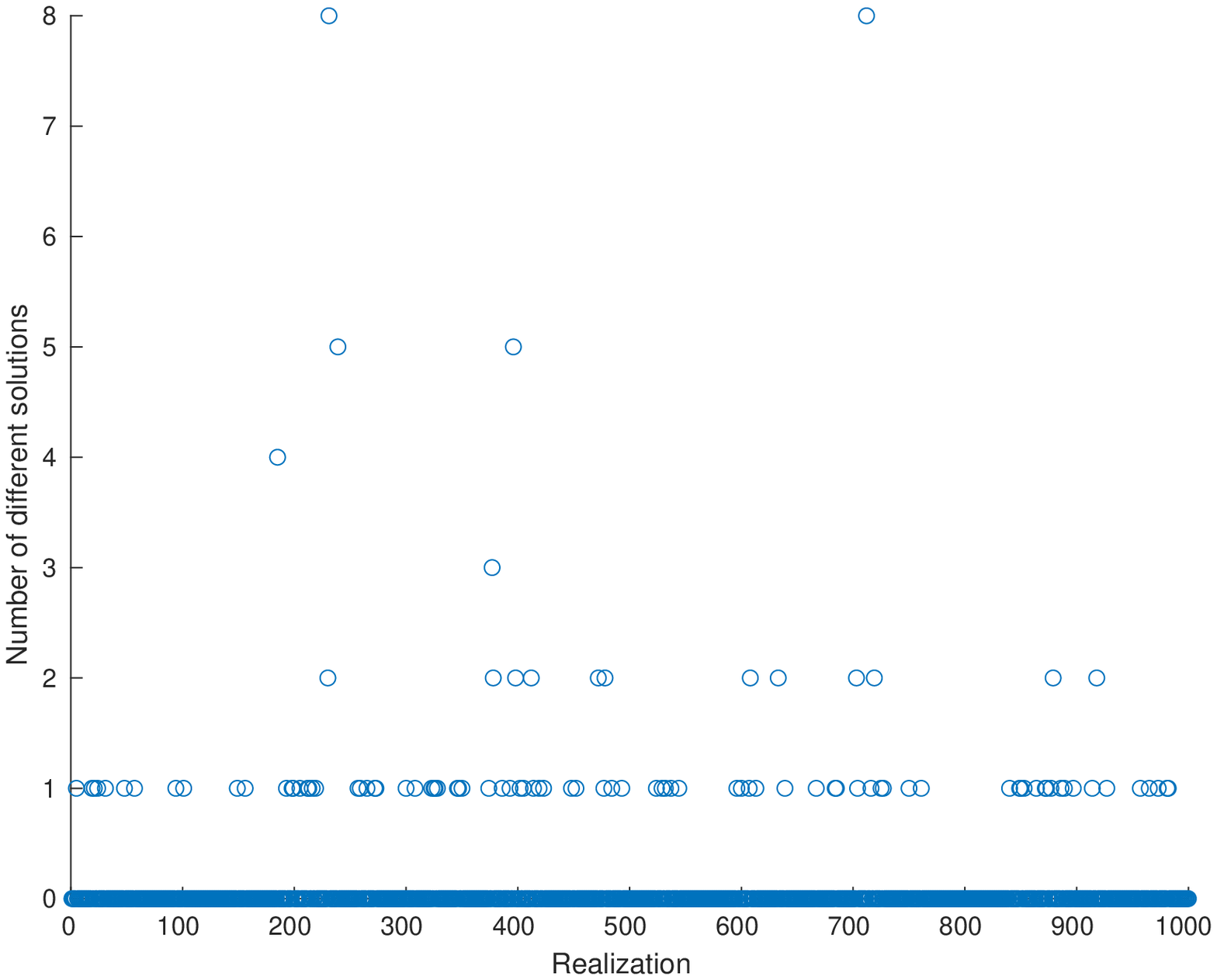}
  \caption{{\tiny} \space  For $N=25, \beta =3, h=0.5$ we simulated 1000 realizations of the disorder and 500 start values used for all of them. For $\varepsilon = -0.505 + 0.001 \cdot k, k=0,\ldots, 250$ we made up to 1.000 iterations of the Banach-$\varepsilon$-algorithm where we stopped when we hit the target level of stability (measured with respect to the {\it maximum absolute error}) of $10^{-7}$. For all realizations we plotted the number of different solutions found.}
\end{figure}
\begin{figure}[H]
  \centering
  % \psfrag{2}{\tiny Two Step}
  \includegraphics[width=0.5\textwidth,  trim={0.45cm 0.4cm 0 0},clip]{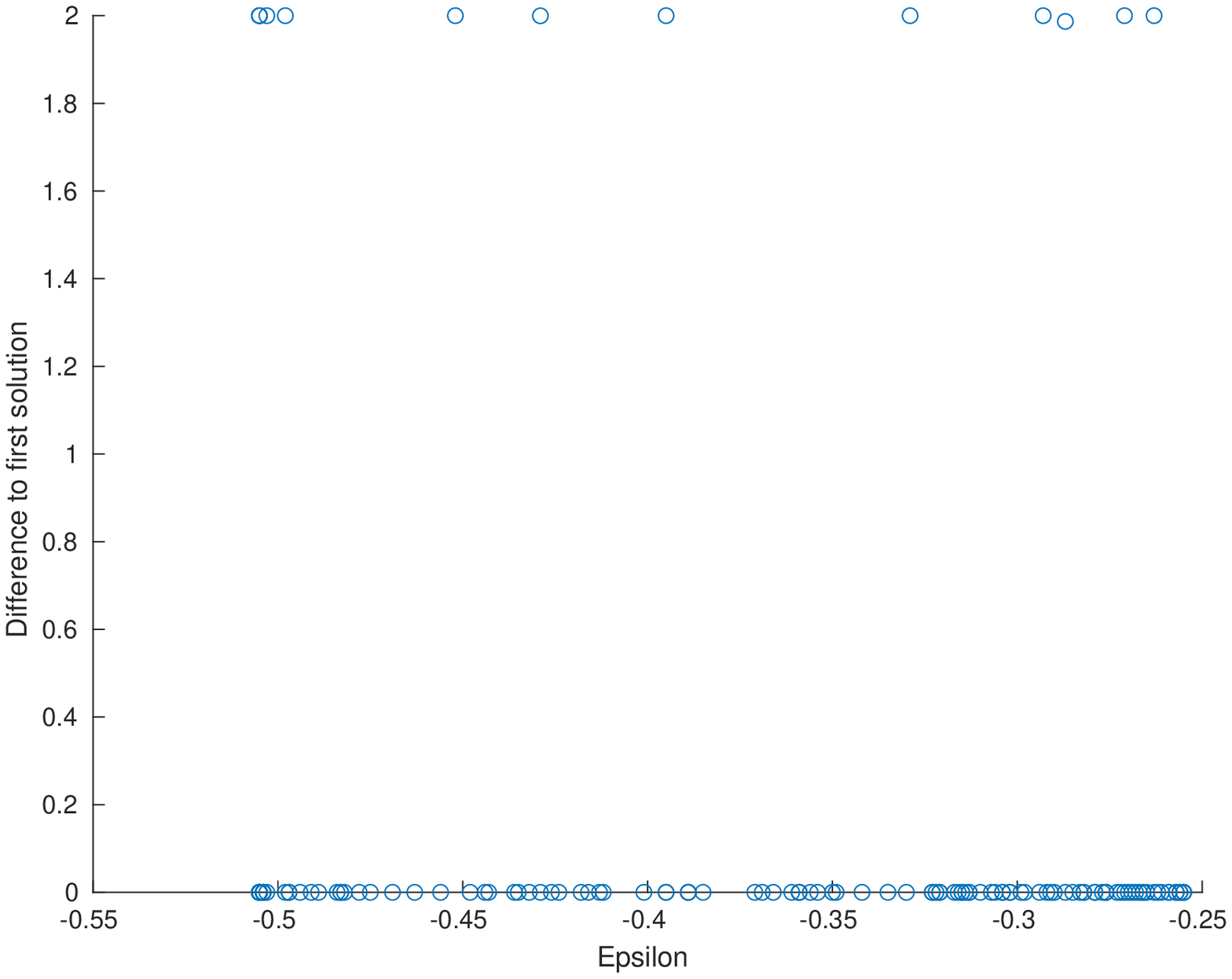}
  \caption{{\tiny} \space  For $N=25, \beta =3, h=0.5$ we simulated 1000 realizations of the disorder and 500 start values used for all of them. For $\varepsilon = -0.505 + 0.001 \cdot k, k=0,\ldots, 250$ we made up to 1.000 iterations of the Banach-$\varepsilon$-algorithm where we stopped when we hit the target level of stability (measured with respect to the {\it maximum absolute error}) of $10^{-7}$. For all realizations we plotted the {\it maximum absolute error} to the first solution found for all different solutions found.}
\end{figure}
\begin{figure}[H]
  \centering
  % \psfrag{2}{\tiny Two Step}
  \includegraphics[width=0.5\textwidth,  trim={0.45cm 0.4cm 0 0},clip]{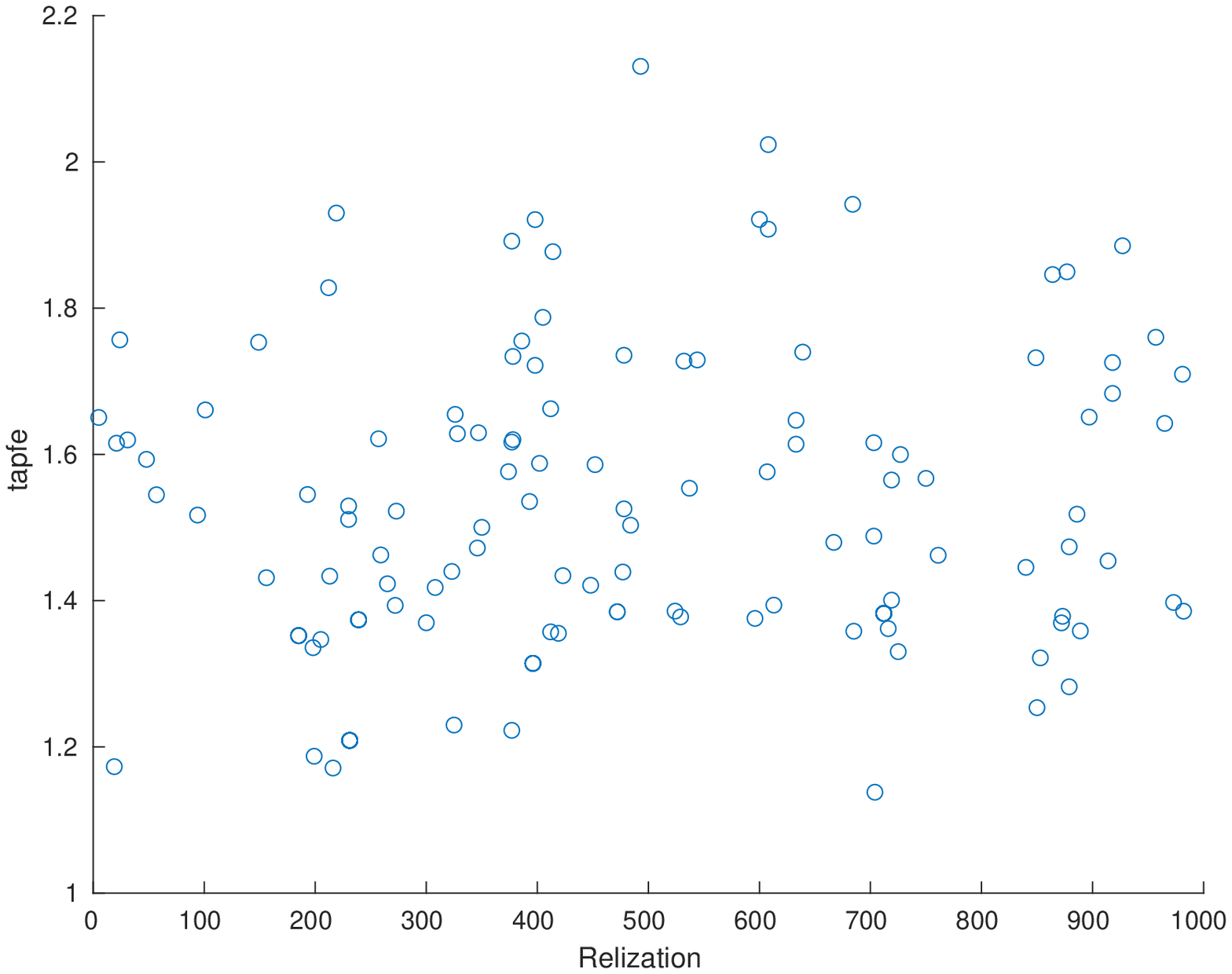}
  \caption{{\tiny} \space  For $N=25, \beta =3, h=0.5$ we simulated 1000 realizations of the disorder and 500 start values used for all of them. For $\varepsilon = -0.505 + 0.001 \cdot k, k=0,\ldots, 250$ we made up to 1.000 iterations of the Banach-$\varepsilon$-algorithm where we stopped when we hit the target level of stability (measured with respect to the {\it maximum absolute error}) of $10^{-7}$. For all realizations we plotted the TAP free energy of the found different solutions.}
\end{figure}
\begin{figure}[H]
  \centering
  % \psfrag{2}{\tiny Two Step}
  \includegraphics[width=0.5\textwidth,  trim={0.45cm 0.4cm 0 0},clip]{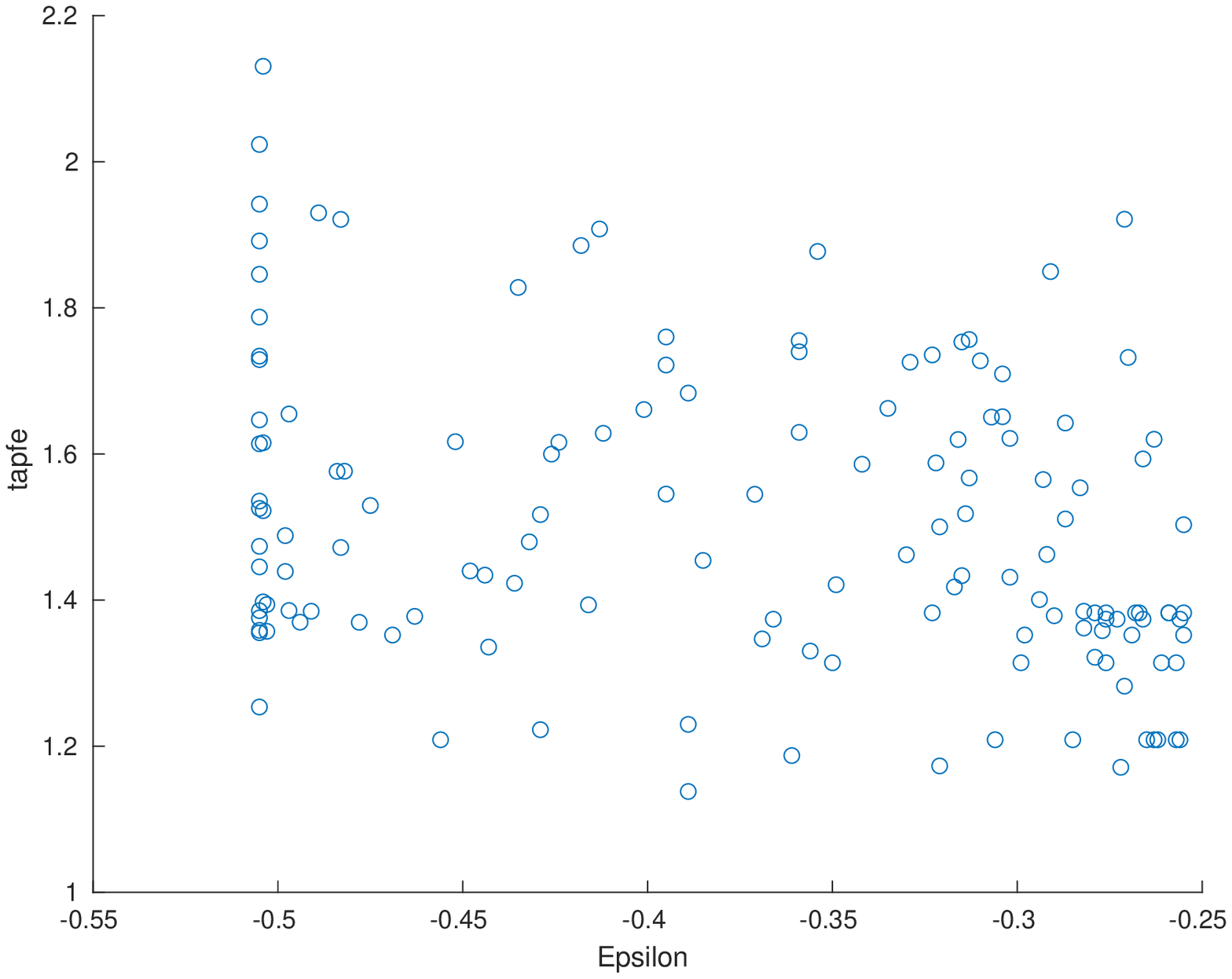}
  \caption{{\tiny} \space  For $N=25, \beta =3, h=0.5$ we simulated 1000 realizations of the disorder and 500 start values used for all of them. For $\varepsilon = -0.505 + 0.001 \cdot k, k=0,\ldots, 250$ we made up to 1.000 iterations of the Banach-$\varepsilon$-algorithm where we stopped when we hit the target level of stability (measured with respect to the {\it maximum absolute error}) of $10^{-7}$. For all $\varepsilon$ we plotted the TAP free energy of the found different solutions.}
\end{figure}
\begin{figure}[H]
  \centering
  % \psfrag{2}{\tiny Two Step}
  \includegraphics[width=0.5\textwidth,  trim={0.45cm 0.4cm 0 0},clip]{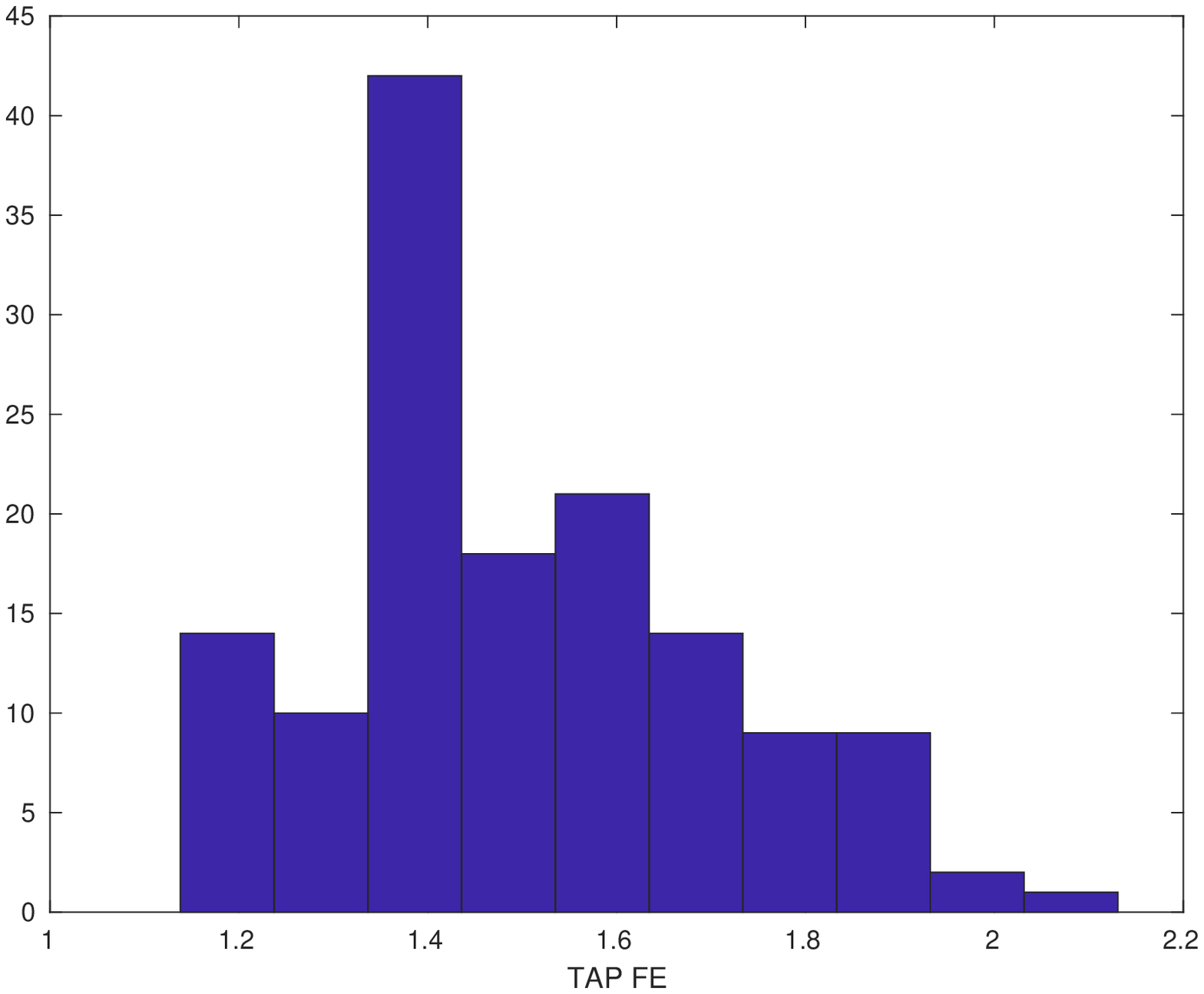}
  \caption{{\tiny} \space  For $N=25, \beta =3, h=0.5$ we simulated 1000 realizations of the disorder and 500 start values used for all of them. For $\varepsilon = -0.505 + 0.001 \cdot k, k=0,\ldots, 250$ we made up to 1.000 iterations of the Banach-$\varepsilon$-algorithm where we stopped when we hit the target level of stability (measured with respect to the {\it maximum absolute error}) of $10^{-7}$. We plotted a histogram of the TAP free energy values of the different solutions we found.}
\end{figure}\FloatBarrier
\newpage\null\thispagestyle{empty}